\documentclass[a4paper,11pt]{book}

\usepackage[T1]{fontenc}
\usepackage[utf8]{inputenc}
\usepackage{lmodern}

\usepackage{graphicx}
\usepackage{upgreek}
\usepackage{nicefrac,amssymb,amsthm,amsfonts,amsmath,epsfig,color,bm,yfonts,isomath}
\usepackage{algorithm,algorithmic}
\usepackage{booktabs}
\usepackage[mathscr]{eucal}
\usepackage{srcltx}
\usepackage{epsfig}
\usepackage{makeidx}
\usepackage{multicol}
\usepackage{upgreek}
\usepackage{tikz}
\usetikzlibrary{calc}
\usepackage{subcaption,comment}
\usepackage{slashbox}

\newcommand{\bA}{{\bm A}}
\newcommand{\bB}{{\bm B}}
\newcommand{\bD}{{\bm D}}
\newcommand{\bF}{{\bm F}}
\newcommand{\bI}{{\bm I}}
\newcommand{\bK}{{\bm K}}
\newcommand{\bM}{{\bm M}}

\newcommand{\bP}{{\bm P}}
\newcommand{\bQ}{{\bm Q}}
\newcommand{\bR}{{\bm R}}
\newcommand{\bS}{{\bm S}}
\newcommand{\bU}{{\bm U}}
\newcommand{\bW}{{\bm W}}
\newcommand{\bY}{{\bm Y}}
\newcommand{\bLambda}{{\bm \Lambda}}
\newcommand{\bPhi}{{\bm \Phi}}
\newcommand{\bPsi}{{\bm \Psi}}
\newcommand{\bSigma}{{\bm \Sigma}}

\newcommand{\bv}{{\boldsymbol v}}

\newcommand{\bx}{{\boldsymbol x}}
\newcommand{\bs}{{\boldsymbol s}}
\newcommand{\bn}{{\boldsymbol n}}
\newcommand{\bmu}{{\boldsymbol\mu}}
\newcommand{\bnu}{{\boldsymbol\nu}}

\newcommand{\tc}{{t_\circ}}
\newcommand{\te}{{t_\mathsf e}}
\newcommand{\yd}{{y^\mathsf d}}
\newcommand{\ydQ}{{y^\mathsf d_1}}
\newcommand{\ydT}{{y^\mathsf d_2}}
\newcommand{\yid}{{y^\mathsf{id}}}

\newcommand{\hydQ}{{\hat y^\mathsf d_1}}
\newcommand{\hydT}{{\hat y^\mathsf d_2}}
\newcommand{\un}{{u^\mathsf n}}
\newcommand{\udn}{{u^{\mathsf d\mathsf n}}}
\newcommand{\ubn}{{u^\mathsf n}}
\newcommand{\md}{{\mathsf m_\mathsf d}}
\newcommand{\mb}{{\mathsf m_\mathsf b}}
\newcommand{\X}{{\mathscr X}}
\newcommand{\Xad}{{\X_\mathsf{ad}}}
\newcommand{\Y}{{\mathscr Y}}
\newcommand{\U}{{\mathscr U}}
\newcommand{\Uad}{{\U_\mathsf{ad}}}
\newcommand{\Ha}{{\mathscr H}}

\newcommand{\tU}{{\tilde{\mathscr U}}}
\newcommand{\tUad}{{\tU_\mathsf{ad}}}
\newcommand{\mP}{{\mathscr P}}
\newcommand{\W}{{\mathscr W}}

\newcommand{\Z}{{\mathscr Z}}
\newcommand{\Zad}{{\Z_\mathsf{ad}}}
\newcommand{\Pf}{{\mathscr P_\mathsf f}}
\newcommand{\Ps}{{\mathscr P_\mathsf s}}
\newcommand{\Pfell}{{\mathscr P_\mathsf f^\ell}}
\newcommand{\Psell}{{\mathscr P_\mathsf s^\ell}}

\newcommand{\ya}{y_\mathsf{a}}
\newcommand{\yb}{y_\mathsf{b}}
\newcommand{\hya}{\hat{y}_\mathsf{a}}
\newcommand{\hyb}{\hat{y}_\mathsf{b}}
\newcommand{\yout}{y_\mathsf{out}}
\newcommand{\AaiU}{{\mathscr A_{ai}^\U}}
\newcommand{\AbiU}{{\mathscr A_{bi}^\U}}
\newcommand{\AaW}{{\mathscr A_{a}^\W}}
\newcommand{\AbW}{{\mathscr A_{b}^\W}}
\newcommand{\IiU}{{\mathscr I_i^\U}}
\newcommand{\IW}{{\mathscr I^\W}}
\newcommand{\uni}{{u_i^{\mathsf n}}}
\newcommand{\Xadeps}{{\X_\mathsf{ad}^{\varepsilon}}}
\newcommand{\Xadn}{{\X_\mathsf{ad}^{\varepsilon,n}}}

\newcommand{\ZAD}{{\mathfrak Z_\mathsf{ad}}}

\newcommand{\ua}{u_{\mathsf a}}
\newcommand{\ub}{u_{\mathsf b}}
\newcommand{\Jhat}{\hat{J}}

\DeclareMathOperator*{\argmax}{arg\hspace{0.1mm}max}
\DeclareMathOperator*{\esssup}{ess\hspace{0.1mm}sup}
\DeclareMathOperator*{\Span}{span}

\makeindex

\newtheorem{run}{Run}
\newtheorem{theorem}{Theorem}
\newtheorem{lemma}[theorem]{Lemma}
\newtheorem{corollary}[theorem]{Corollary}
\newtheorem{proposition}[theorem]{Proposition}
\newtheorem{definition}[theorem]{Definition}
\newtheorem{example}[theorem]{Example}
\newtheorem{remark}[theorem]{Remark}
\newtheorem{assumption}[theorem]{Assumption}
\newtheorem{problem}[theorem]{Problem}

\renewcommand{\thealgorithm}{\arabic{chapter}.\arabic{section}.\arabic{algorithm}}

\renewcommand{\theequation}{\arabic{chapter}.\arabic{section}.\arabic{equation}}
\renewcommand{\thecorollary}{\arabic{chapter}.\arabic{section}.\arabic{corollary}}
\renewcommand{\thetheorem}{\arabic{chapter}.\arabic{section}.\arabic{theorem}}
\renewcommand{\theproposition}{\arabic{chapter}.\arabic{section}.\arabic{proposition}}
\renewcommand{\theremark}{\arabic{chapter}.\arabic{section}.\arabic{remark}}
\renewcommand{\thedefinition}{\arabic{chapter}.\arabic{section}.\arabic{definition}}
\renewcommand{\thefigure}{\arabic{chapter}.\arabic{section}.\arabic{figure}}
\renewcommand{\thetable}{\arabic{chapter}.\arabic{section}.\arabic{figure}}
\renewcommand{\theexample}{\arabic{chapter}.\arabic{section}.\arabic{example}}
\renewcommand{\theassumption}{\arabic{chapter}.\arabic{section}.\arabic{assumption}}
\renewcommand{\theproblem}{\arabic{chapter}.\arabic{section}.\arabic{probleme}}

\renewenvironment{algorithm}[1][\relax]{\refstepcounter{algorithm}%
\addcontentsline{loa}{algorithm}%
    {\protect\numberline{Algorithm~\thealgorithm}{\ignorespaces#1}}%
\par\vspace{1\baselineskip}%
\expandafter\ifx#1\relax
\parindent0pt {\scshape\bfseries Algorithm~\thealgorithm.}\\
\else
\parindent0pt {\scshape\bfseries Algorithm~\thealgorithm.}\enspace{\bfseries#1.}\\
\fi}

\makeatletter
\renewcommand{\@chapapp}{}

\makeatother

\title{\Huge\textbf{POD Suboptimal Control}\\\textbf{of Evolution Problems}\\
\huge Theory and Applications}
\author{\textsc{S. Banholzer}, \textsc{D. Beermann}, \textsc{L. Mechelli}, \textsc{S. Volkwein}\thanks{Stefan Volkwein, Department of Mathematics and Statistics, Konstanz University, D-78457 Konstanz, Germany, e-mail: \texttt{Stefan.Volkwein@uni-konstanz.de}}}

\begin{document}


\maketitle


\tableofcontents


\mainmatter


\chapter{Introduction}
\label{SIAM-Book:Section1}

The optimization and control of systems governed by \index{Equation!partial differential, PDE}{\em partial differential equations (PDEs)} usually requires numerous evaluations of the forward problem or the optimality system. Despite the fact that many recent efforts have been made to limit or reduce the number of evaluations to 5-10, this cannot be achieved in all situations and even if this is possible, these evaluations may still require a formidable computational effort. For situations, where this effort is not acceptable, model-order reduction can be a means to significantly reduce the required computational resources. In this work we utilize the method of proper orthogonal decomposition (POD) in order to numerically solve optimization problems governed by PDEs. This reduced-order approach is based on projecting the dynamical system onto subspaces consisting of basis elements that contain characteristics of the expected solution. This stands in contrast to, e.g., finite element techniques, where the basis elements of the subspaces do not relate to the physical properties of the system that they approximate. 

\section{Reduced-order modeling for optimal control problems}
\label{SIAM-Book:Section1.1}

Optimal control problems for PDEs are often hard to tackle numerically because their discretization leads to very large-scale optimization problems. Therefore, different techniques of model reduction were developed to approximate these problems by smaller ones that are tractable with less effort.

The \index{Method!balanced truncation}{\em balanced truncation method} \cite{Ant05,SVR08,ZDG96} is one well studied model reduction technique for state-space systems. This method utilizes the solutions to two Lyapunov equations, the so-called controllability and observability Gramians. The balanced truncation method is based on transforming the state-space system into a balanced form so that these Gramians become diagonal and equal. Moreover, the states that are difficult to reach or to observe, are truncated. The advantage of this method is that it preserves the asymptotic stability in the reduced-order system. Furthermore, a-priori error bounds are available. Recently, the theory of balanced truncation model reduction was extended to descriptor systems; see, e.g., \cite{MS05} and \cite{HSS08}.

Recently the application of {\em reduced-order models} to linear, time varying and nonlinear systems, in particular to nonlinear control systems, has received an increasing amount of attention. The {\em reduced basis} (RB) method, as developed in \cite{GMNP07,PR05} and \cite{IR98}, is one such reduced-order method, where the basis elements correspond to the dynamics of expected control regimes. Let us refer to \cite{DH13,HO14,NRMQ13,OS13} for the successful use of the reduced basis method in PDE constrained optimization problems. Currently, {\em POD} is probably the mostly used and most successful model reduction technique for nonlinear optimal control problems, where the basis functions contain information from the solutions of the dynamical system at pre-specified time-instances, so-called snapshots; see, e.g., \cite{Cha00,HLBR12,Sir87,Vol12}. Due to a possible linear or almost-linear dependence on the time variable, the snapshots themselves are not appropriate as a basis. Hence a singular value decomposition is carried out and the leading generalized eigenfunctions are chosen as a basis, referred to as the POD basis. POD is successfully used in a variety of fields including fluid dynamics, coherent structures \cite{AH02,AFS00} and inverse problems \cite{BJWW00}. Moreover, in \cite{ABK01} POD is successfully applied to compute reduced-order controllers. The relationship between POD and balancing was considered in \cite{LMG02,Row05,WP02}. An error analysis for nonlinear dynamical systems in finite dimensions was carried out in \cite{RP02} and a missing point estimation in models described by POD was studied in \cite{AWWB08}. Let us mention that POD and the reduced basis method are successfully combined by variants of the POD greedy algorithm; see \cite{HO08} and \cite{Haa12}, for instance.

Reduced-order models are used in PDE-constrained optimization in various ways; see, e.g., \cite{HV05,SV10} for a survey. In optimal control problems, it is sometimes necessary to compute a feedback control law instead of a fixed optimal control. In the implementation of these feedback laws, models of reduced order can play an important and very useful role, see \cite{ABK01,KVX04,LV06,Rav00}. Another useful application is the use in optimization problems, where a PDE solver is part of the function evaluation. Obviously, thinking of a gradient evaluation or even a step-size rule in the optimization algorithm, an expensive function evaluation leads to an enormous amount of computing time. Here, the reduced-order model can replace the system given by a PDE in the objective function. It is quite common that a PDE can be replaced by a five- or ten-dimensional system of ordinary differential equations. This results computationally in a very fast method for optimization compared to the effort for the computation of a single solution of a PDE. There is a large amount of literature in engineering applications in this regard, we mention only the papers \cite{LT01,NAMTT03}. Recent applications can also be found in finance using the RB model \cite{Pir09} and the POD model \cite{SS13,Sch12} in the context of calibration for models in option pricing.

In the present work we use POD for deriving low order models of dynamical systems. These low order models then serve as surrogates for the dynamical system in the optimization process. We consider a linear-quadratic optimal control problem in an abstract setting and prove error estimates for the POD Galerkin approximations of the optimal control.

\section{PDE constrained optimal control problems}
\label{SIAM-Book:Section1.2}

To motivate the POD reduced-order modelling, we choose the following generic linear-quadratic optimal control problem: Suppose that the spatial $\Omega$ is a bounded spatial \index{Domain, spatial, $\Omega$}{\em domain} in $\mathbb R^\mathfrak n$ with $\mathfrak n\in\{1,2,3\}$. We write $\bx=(x_1,\ldots,x_\mathfrak n)$ for an element in $\Omega$. The (sufficiently smooth) boundary $\partial \Omega$ is denoted by $\Gamma$. For given $T>0$ we define the time-space cylinder $Q=(0,T)\times\Omega$ and the boundary set $\Sigma=(0,T)\times\Gamma$. All controls are assumed to have the structure
\begin{align*}
    u =(u^{\mathsf d},u^{\mathsf b}): \, (0,T) \to \mathbb R^\md \times \mathbb R^\mb
\end{align*}
with $\md,\,\mb\in\mathbb N$. The \index{Set!of admissible controls, $\Uad$}{\em set of admissible controls} has the simple form
\begin{align*}
    \Uad=\Big\{ &u:(0,T)\to\mathbb R^\md \times \mathbb R^\mb ~\big|~ u^{\mathsf d}_{\mathsf a i}(t) \le u^{\mathsf d}_i(t) \le u^{\mathsf d}_{\mathsf b i}(t),~i=1,...,\md, \\
	&\hspace{42.5mm}u^{\mathsf b}_{\mathsf a j}(t) \le u^{\mathsf b}_j(t) \le u^{\mathsf b}_{\mathsf b j}(t),~ j=1,...,\mb,t \in (0,T) \Big \},
\end{align*}
where the functions $u^{\mathsf d}_\mathsf a, u^{\mathsf d}_\mathsf b: (0,T) \to \mathbb R^\md$ and $u^{\mathsf b}_\mathsf a, u^{\mathsf b}_\mathsf b: (0,T) \to \mathbb R^\mb$ are lower and upper control bounds satisfying $u^{\mathsf d}_{\mathsf a i}(t) \le u^{\mathsf d}_{\mathsf b i}(t)$ and $u^{\mathsf b}_{\mathsf a j}(t) \le u^{\mathsf b}_{\mathsf b j}(t)$ for all $i=1,...,\md$, $j=1,...,\mb$ and $t \in (0,T)$. 
For any admissible control $u\in\Uad$, for a given (bounded) initial condition $y_\circ:\Omega\to\mathbb R$ and for inhomogeneities $f:Q\to\mathbb R$, $g: \Sigma\to\mathbb R$, the state variable $y=y(t,\bx)$ is the solution the the linear \index{Equation!heat}{\em heat equation} with convection which is a linear \index{Equation!partial differential, PDE}\index{Equation!partial differential!parabolic}{\em parabolic partial differential equation} (PDE)
\begin{subequations}
    \label{SIAM:EqMotPDE1}
    \begin{align}
        \label{SIAM:EqMotPDE1a}
        y_t(t,\bx)-\kappa \Delta y(t,\bx)+\bv(\bx) \cdot \nabla y(t,\bx)&=f(t,\bx)+\sum_{i=1}^\md u^{\mathsf d}_i(t)\chi_i(\bx)&&\text{for }(t,\bx)\in Q,\\
        \label{SIAM:EqMotPDE1b}
        \kappa \frac{\partial y}{\partial\bn}(t,\bs)+q(\bs)y(t,\bs)&=g(t,\bs)+\sum_{j=1}^\mb u^{\mathsf b}_j(t)\xi_j(\bs)&&\text{for }(t,\bs)\in\Sigma,\\
        \label{SIAM:EqMotPDE1c}
        y(0,\bx)&=y_\circ(\bx)&&\text{for }\bx\in\Omega.
    \end{align}
\end{subequations}
In \eqref{SIAM:EqMotPDE1}, the variables $\chi_i: \Omega \to \mathbb R$, $i=1,...,\md$, and $\xi_j: \Gamma \to \mathbb R$, $j=1,...,\mb$ are \index{Control shape functions}{\em control shape functions} which indicate the effect of each control function in $Q$ and $\Sigma$. Furthermore, the constant $\kappa>0$ and the vector-valued function $\bv=(v_i)_{1\le i\le\mathfrak n}:\Omega\to\mathbb R^\mathfrak n$ represent the diffusion and advection part of the equation while $q: \Gamma \to \mathbb R$ with $q\ge0$ in $Q$ is a boundary function. The dot between $\bv(\bx)$ and the gradient of $y$ in \eqref{SIAM:EqMotPDE1a} indicates the Euclidean inner product and is to be understood as 
\begin{align*}
    \bv(\bx)\cdot\nabla y(t,\bx)={\langle\bv(\bx),\nabla y(t,\bx)\rangle}_2=\sum_{j=1}^\mathfrak n v_j(\bx) \frac{\partial y}{\partial \bx_j}(t,\bx), \quad (t,\bx) \in Q.
\end{align*}
The term
\begin{align*}
    \frac{\partial y}{\partial \bn}(t,\bs)=\nabla y(t,\bs)\cdot\bn(\bx),\quad(t,\bs)\in\Sigma,
\end{align*}
stands for the (outward) {\em normal derivative} of $y$ and $\bn$ denotes the outward normal vector. Let for any $u \in \Uad$ there exists a unique solution $y = y(u)$ of \eqref{SIAM:EqMotPDE1}. This issue will be discussed in detail in Section~\ref{SIAM-Book:Section3.2}. Let us introduce the generic \index{Cost functional!quadratic}\index{Cost functional!tracking-type}{\em quadratic cost functional of tracking type} by
\begin{align*}
    J(y,u)&=\frac{\alpha_Q}{2}\int_0^T\int_\Omega\,\big|y(t,\bx)-y_Q(t,x)\big|^2\,\mathrm d\bx\mathrm dt+\frac{\alpha_\Omega}{2}\int_\Omega\,\big|y(T,\bx)-y_\Omega(\bx)\big|^2\,\mathrm d\bx\\
    &\quad+\frac{1}{2}\int_0^T\sum_{i=1}^\md \gamma^{\mathsf d}_i\,\big| u^{\mathsf d}_i(t)-u^{\mathsf d\mathsf n}_i(t) \big|^2+\sum_{j=1}^\mb\gamma^{\mathsf b}_j\,\big | u^{\mathsf b}_j(t)-u^{\mathsf b\mathsf n}_j(t) \big|^2\,\mathrm dt,
\end{align*}
where $\alpha_Q, \alpha_\Omega \ge 0$ are nonnegative \index{Cost functional!tracking weight, $\alpha_Q$, $\alpha_\Omega$}{\em tracking weights} for the state and $\gamma^\mathsf d_i$ ($1\le i\le\md)$, $\gamma^\mathsf b_j$ ($1\le j\le \mb)$ are positive \index{Cost functional!regularization parameter, $\gamma_i^\mathsf d$, $\gamma_i^\mathsf b$}{\em regularization parameters} for the controls. The functions $y_Q: Q \to \mathbb R$ and $y_\Omega: \Omega \to \mathbb R$ are given \index{Cost functional!desired states, $\ydQ$, $\ydT$}{\em desired states} and $\un=(\udn,\ubn):(0,T)\to\mathbb R^{m_\mathsf d}\times\mathbb R^{m_\mathsf b}$ is a chosen \index{Cost functional!nominal control, $\un$}{\em nominal} or \index{Cost functional!expected control, $\un$}{\em expected control}. Now the optimal control problem is of the following form:
\begin{equation}
    \label{OC}
    \tag{$\mathbf{OC}$}
    \min J(y,u)\quad\text{subject to (s.t.)}\quad u\in\Uad \text{ and }(y,u)\text{ satisfies \eqref{SIAM:EqMotPDE1}}.
\end{equation}
Notice that \eqref{OC} is a \index{Problem!quadratic programming!PDE-constrained}{\em PDE-constrained quadratic programming problem}; cf. \cite{BS12,HPUU09,IK08,Tro10}. Due to the fact that the state variable $y$ depends on the control variable $u$ via the state equations \eqref{SIAM:EqMotPDE1} we call \eqref{OC} an \index{Problem!optimal control}{\em optimal control problem}; cf. \cite{Tro10}. Furthermore, \eqref{OC} is said to be a \index{Problem!optimal control!linear-quadratic}{\em linear-quadratic optimal control problem}, because the cost functional is quadratic and the state equation linear.

Since \eqref{SIAM:EqMotPDE1} possesses a unique solution $y=y(u)$ for any $u\in\Uad$, we are able to define the {\em reduced cost functional}\index{Cost functional!reduced} $\hat J$ as
\begin{align*}
    \hat J(u)=J(y(u),u)\quad\text{for }u\in\Uad.
\end{align*}
Thus, \eqref{OC} can be equivalently written as the following \index{Problem!optimal control!reduced}{\em reduced optimal control problem}
\begin{equation}
    \label{OChat}
    \tag{$\mathbf{\widehat{OC}}$}
    \min\hat J(u)\quad\text{subject to (s.t.)}\quad u\in\Uad.
\end{equation}
If $\bar u$ is a local solution to \eqref{OChat}, then $(\bar y,\bar u)$ is a solution to \eqref{OC}, where $\bar y=y(\bar u)$ denotes the unique solution to \eqref{SIAM:EqMotPDE1} for $u=\bar u$.

\begin{remark}
    \rm The term
    \begin{align*}
        -\kappa\Delta y(t,\bx)+\bv(\bx)\cdot \nabla y(t,\bx),\quad(t,\bx)\in Q,
    \end{align*}
    in \eqref{SIAM:EqMotPDE1} can be replaced by the generalized elliptic differential operator
    \begin{align*}
        \big(\mathcal Ay\big)(t,\bx)=-\sum_{i=1}^\mathfrak n\bigg(\frac{\partial}{\partial x_i}\Big(\sum_{j=1}^\mathfrak na_{ij}(t,\bx)\frac{\partial y}{\partial x_j}(t,\bx)\Big)+v_i(t,\bx)\frac{\partial y}{\partial x_i}(t,\bx)\bigg)+c(t,\bx)y(t,\bx)
    \end{align*}
    for $(t,\bx)\in Q$ with bounded coefficient functions $a_{ij}$, $v_i$ and $c$. Moreover, more general objectives $J=J(y,u)$ can be studied; see \cite{HPUU09,Tro10}, for instance. To simplify the presentation we concentrate on the specific form \eqref{SIAM:EqMotPDE1}.\hfill$\blacksquare$
\end{remark}

\section{Guiding model problem for the numerical experiments}
\label{SIAM-Book:Section1.3}

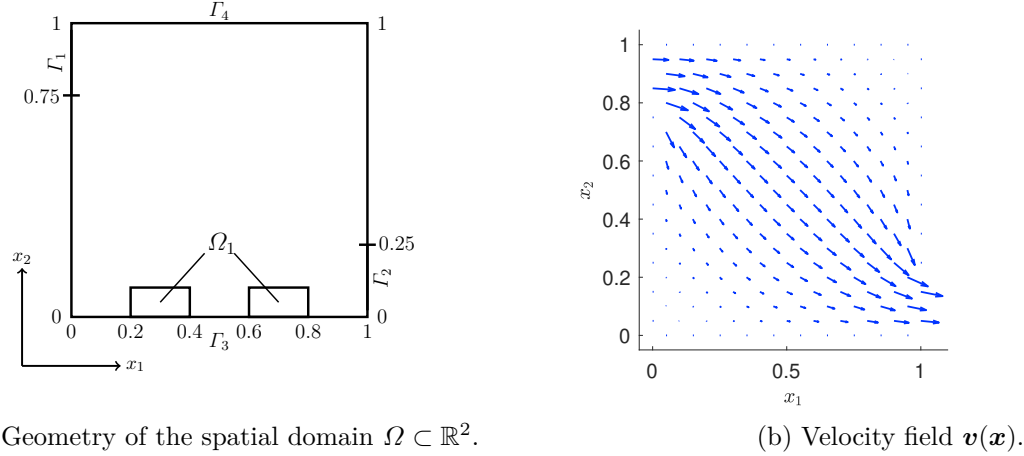
\begin{figure}
	\centering
	\begin{subfigure}{0.48\textwidth}	
	\centering
 	\resizebox{0.69\textwidth}{!}{
 	\begin{tikzpicture}
 	  \node [rectangle, draw, minimum size=6cm, line width=0.5mm] (Omega) at (0.5,0) {};
 	  \node [rectangle, draw, minimum width=1.2cm, minimum height=0.6cm, line width=0.5mm] (Omega) at (-0.7,-2.7) {};
		\node [rectangle, draw, minimum width=1.2cm, minimum height=0.6cm, line width=0.5mm] (Omega) at (1.7,-2.7) {};
		
 	  \node (upLeft) at (-2.5,3) {};
 	  \node (downRight) at (3.5,-3) {};
 	  
 	  \draw [-|,line width=0.5mm] (upLeft) -- ($(upLeft)-(0,1.5)$);
 	  \draw [-|,line width=0.5mm] (downRight) -- ($(downRight)+(0,1.5)$);
 	  \node at ($(upLeft)-(0.3,0)$) {\large 1};
 	  \node at ($(upLeft)-(0.6,1.5)$) {\large 0.75};
 	  \node at ($(upLeft)-(0.3,6)$) {\large 0};
 	  \node at ($(downRight)+(0.3,6)$) {\large 1};
 	  \node at ($(downRight)+(0.6,1.5)$) {\large 0.25};
 	  \node at ($(downRight)+(0.3,0)$) {\large 0};
 	  \node at ($(upLeft)-(0,6.3)$) {\large 0};
 	  \node at ($(upLeft)+(1.2,-6.3)$) {\large 0.2};
 	  \node at ($(upLeft)+(2.4,-6.3)$) {\large 0.4};
 	  \node at ($(upLeft)+(3.6,-6.3)$) {\large 0.6};
 	  \node at ($(upLeft)+(4.8,-6.3)$) {\large 0.8};
 	  \node at ($(upLeft)+(6.0,-6.3)$) {\large 1};
 	   	  
 	  \node [rotate=90] at ($(upLeft)-(0.3,0.8)$) {\large $\Gamma_1$};
		\node [rotate=90] at ($(downRight)+(0.3,0.8)$) {\large $\Gamma_2$};  	
 	  \node at (0.5,-3.5) {\large $\Gamma_3$};
 	  \node at (0.5,3.25) {\large $\Gamma_4$};
 	  
 	  \node at (0.55,-1.5) {\Large $\Omega_1$};
 	  \draw[-,line width=0.3mm] (0.2,-1.7) -- (-0.7,-2.7);
 	  \draw[-,line width=0.3mm] (0.8,-1.7) -- (1.7,-2.7);
 	  
 	  \draw [->,line width=0.4mm] (-3.5,-4) -- (-3.5,-2);
 	  \draw [->,line width=0.4mm] (-3.5,-4) -- (-1.5,-4);
 	  \node at (-3.5,-1.8) {$x_2$};
 	  \node at (-1.2,-4) {$x_1$};
 	    	  
  \end{tikzpicture}
  }
  \vspace{2.3ex}
  \caption{Geometry of the spatial domain $\Omega\subset\mathbb R^2$.}
  \label{fig:domainOmega}
	\end{subfigure}
	\hfill
	\begin{subfigure}{0.48\textwidth}	
	\includegraphics[height=50mm]{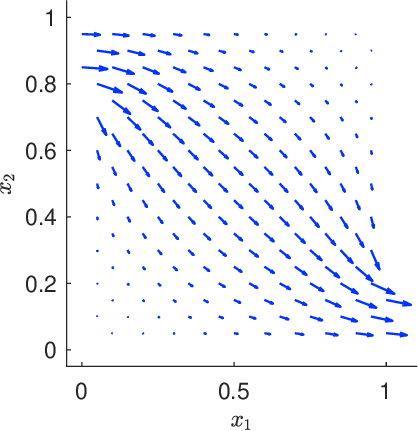}
	\caption{Velocity field $\bv(\bx)$.}
	\label{fig:navierStokesSolution}
	\end{subfigure}
	\caption{Numerical setup in Example \ref{ex:generalNumericSetup}.}
	\label{fig:generalNumericalSetup}
\end{figure}

Let us introduce a particular instance of \eqref{OC} that describes a heat evolution phenomenon and will be used repeatedly for numerical illustrations throughout this book. We choose the two-dimensional unit square $\Omega =(0,1)^2\subset\mathbb R^2$ as the spatial domain and a final time of $T=5$. The (thermal) diffusivity is given as $\kappa=0.5$. The advection field $\bv$ comes from the solution of an incompressible Navier-Stokes equation and can be interpreted as the air velocity between an inflow at $\{0\} \times [0.75,1]$ and an outflow at $\{1\} \times [0,0.25]$. The shape of the flow field is depicted alongside a sketch of the domain $\Omega$ in Figure~\ref{fig:generalNumericalSetup}.
We assume the existence of one distributed control $u^\mathsf d: (0,T) \to \mathbb R$ acting on the subdomain $\Omega_1=([0.2,0.4]\cup [0.6,0.8])\times[0,0.1]$ by a multiplication with a scaled indicator function $\chi: \Omega \to \{0,1\}$ which takes the value $0.1$ inside $\Omega_1$ and $0$ everywhere else. The inhomogeneity $f$ is set to $0$. On the boundary, we interpret the function $q$ as an insulation coefficient between the inside of $\Omega$ and the 'outside world' where we assume the temperature to be described by a function $y_\mathsf{out}(t) \in \mathbb R$ for $t \in (0,T)$. In absence of a boundary control element, this insulation effect can be described by the linear model
\begin{align*}
    \kappa\,\frac{\partial y}{\partial \bn}(t,\bs)+q(\bs)y(t,\bs)=q(\bs) y_\mathsf{out}(t),\quad (t,\bs) \in \Sigma.
\end{align*}
Notice that for every boundary point $\bs \in \Gamma$ and every time instance $t \in (0,T)$, the outward normal derivative $\tfrac{\partial y}{\partial \bn}(t,\bs)$ is proportional to the difference between the inside and outside temperature. In the notation of the general PDE \eqref{SIAM:EqMotPDE1}, this means that the boundary inhomogeneity is given by $g(t,\bs) = q(\bs) y_\mathsf{out}(t)$ for all $(t,\bs) \in \Sigma$. For the insulation function, we define the boundary subsegments $\Gamma_1 = \{0\} \times [0.75,1]$, $\Gamma_2 = \{1\} \times [0,0.25]$, $\Gamma_3 = [0,1] \times \{0\}$ and $\Gamma_4 = [0,1] \times \{1\}$, cf, Figure \ref{fig:domainOmega}. The function $q$ is then set to
\begin{align*}
    q(\bs)=\left\{
    \begin{aligned}
        &0.1&&\text{for }\bs\in\Gamma_1\cup\Gamma_2,\\
        &0 &&\text{for }\bs \in \Gamma_3 \cup \Gamma_4, \\
        &0.01&&\text{otherwise}.
    \end{aligned}
    \right .
\end{align*}
For the PDE model, this means that the upper and lower boundaries are assumed to be interior walls (with no interaction with the outside temperature) and that the insulation at the inflow and outflow segments $\Gamma_1$ and $\Gamma_2$ is comparably poorer than for the rest of these walls. Additionally, we add a single boundary control $u^\mathsf b:(0,T)\to\mathbb R$ entering at the segment $\Gamma_1$ with the scaled indicator function $\xi(\bs) = 0.1$ if $\bs \in \Gamma_1$ and $\xi(\bs)=0$ everywhere else. Moreover, the initial temperature is assumed to be constant at $y_\circ(\bx)=17$ for all $\bx\in\Omega$. 

For the cost function $J$, we would like to track the desired temperature $y_Q(t,x) = 17$ for all $(t,\bx) \in Q$. Additionally, we add a final-time tracking of the same temperature at $y_\Omega(\bx) = 17$ for all $\bx \in \Omega$ for the sake of numerical stability. Every deviation from the nominal controls $u^{\mathsf d,\circ}(t)=0$ and $u^{\mathsf b,\circ}(t)=0$ for $t \in (0,T)$ are penalized with the regularization parameters $\gamma^{\mathsf d} = 0.1$ and $\gamma^{\mathsf b} = 0.1$. In contrast, the state trackings are scaled by the parameters $\alpha_Q=1$ and $\alpha_\Omega=0.1$. 
	
Finally, we restrict the control variables to only take nonnegative values, which means that they can only be used to heat up, not cool down the system. Summarizing the linear-quadratic optimal control problem for our numerical experiments can now be described as follows:
\begin{align}
    \label{generalNumericalSetup_ocp}
    \begin{aligned}
        &\min J(y,u) = \frac{1}{2} \int_0^5 \int_\Omega(y(t,\bx)-17)^2 ~\mathrm d\bx\mathrm dt+ \frac{0.1}{2} \int_\Omega (y(T,\bx)-17)^2 ~\mathrm d\bx\\
        &\hspace{22mm}+ \frac{0.1}{2} \int_0^5 u^\mathsf d(t)^2+u^\mathsf b(t)^2\,\mathrm dt\\
        &\,\text{s.t. }0 \le u^{\mathsf d}(t), ~ 0 \le u^{\mathsf b}(t) \text{ for } t \in (0,5)\text{ and}\\
        &\hspace{4mm}
        \begin{aligned}
            &y_t(t,\bx) - 0.5 \Delta y(t,\bx) + \bv(\bx) \cdot \nabla y(t,\bx) = u^\mathsf d(t) \chi(\bx), &&(t,\bx) \in Q,\\
            &\frac{\partial y}{\partial \bn}(t,\bs) + q(\bs) y(t,\bs)= q(\bs) y_\mathsf{out}(t) + u^\mathsf b(t) \xi(\bs), &&(t,\bs)\in\Sigma,\\
            &y(0,\bx) = y_\circ(\bx), &&\bx \in \Omega.
        \end{aligned}
    \end{aligned}	
\end{align}
For the outside temperature, we assume a periodic function given by $y_\mathsf{out}(t) = 13 + 5 \cos (2 \pi t / 5)$.

\begin{run}
    \label{ex:generalNumericSetup}
    \rm Let us visualize the solution to the parabolic partial differential equation for three different control input functions, namely
    \begin{align*}
        u_j(t) = (u_j^{\mathsf d}(t), u_j^{\mathsf b}(t))^\top = \left \{ 
        \begin{aligned}
            &(0,0)^\top,&j=1, \\
            &(1.8,4.5)^\top,&j=2, \\
            &\big(1-\cos(2 \pi t/T)\big)(1.8,4.5)^\top,&j=3,
        \end{aligned}
        \right\}\quad t \in (0,T).
    \end{align*}
    These lead to solutions $y_j(t,\bx)$ for $(t,\bx) \in Q$, all $j=1,2,3$. The first control was chosen to see the solution without any control appearance. For the second control, we experimentally try to set both $u_2^{\mathsf d}(t)$ and $u_2^{\mathsf b}(t)$ to constant values which lead to a temperature evolution which is as close to $17$ as we can get, for which we have settled at the nominated values. Lastly, we observe that the outer temperature negatively affects the inner temperature with its periodical intensity. We try to counteract this by choosing the likewise periodical control $u_3^{\mathsf d}$ as denoted above. Figure \ref{fig:advectionsDiffusionSolutions_averages} shows us the time-average solutions 
    \begin{align*}
        \overline y_j(t) = \frac{1}{|\Omega|} \int_\Omega y_j(t,\bx)\,\mathrm d\bx, \quad t \in (0,T)\text{ and }j=1,2,3,
    \end{align*}
    where $|\Omega|$ stands for the (Lebesgue) measure of the bounded set $\Omega$.
    \begin{figure}
        \centering
        \includegraphics[height=50mm]{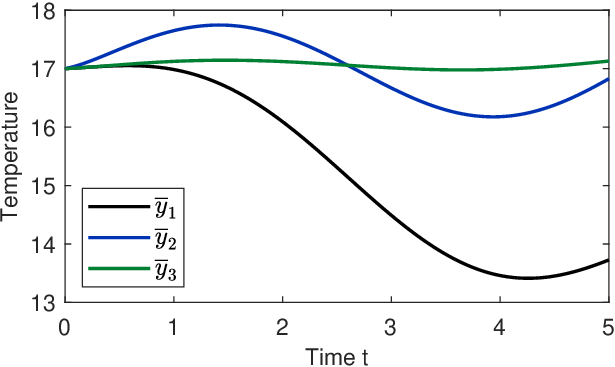}
        \caption{Run~\ref{ex:generalNumericSetup}. Time-average temperature evolution of the solutions to \eqref{SIAM:EqMotPDE1} for the fixed control inputs $u_1$, $u_2$, and $u_3$.}
        \label{fig:advectionsDiffusionSolutions_averages}
    \end{figure}
    As we can see, not applying any control effort as in $u_1$ results in a significant drop in temperature which follows the periodicity of the cosinal outside temperature $y_\mathsf{out}$. This effect is still visible for the constant control $u_2$ which results in an overshoot of the target at $17$ in the first half of the time horizon, and an undershoot in the second half. By setting the control $u_3$ to be likewise periodical but inverse in orientation, this effect is counteracted and we get a very close approximation of the desired temperature. In Figure \ref{fig:advectionsDiffusionSolutions_snapshots}, we can see so-called \emph{snapshots} of the solutions at the time instances $t=2$ and $t=4$.
    \begin{figure}
        \centering
        \includegraphics[width=\textwidth]{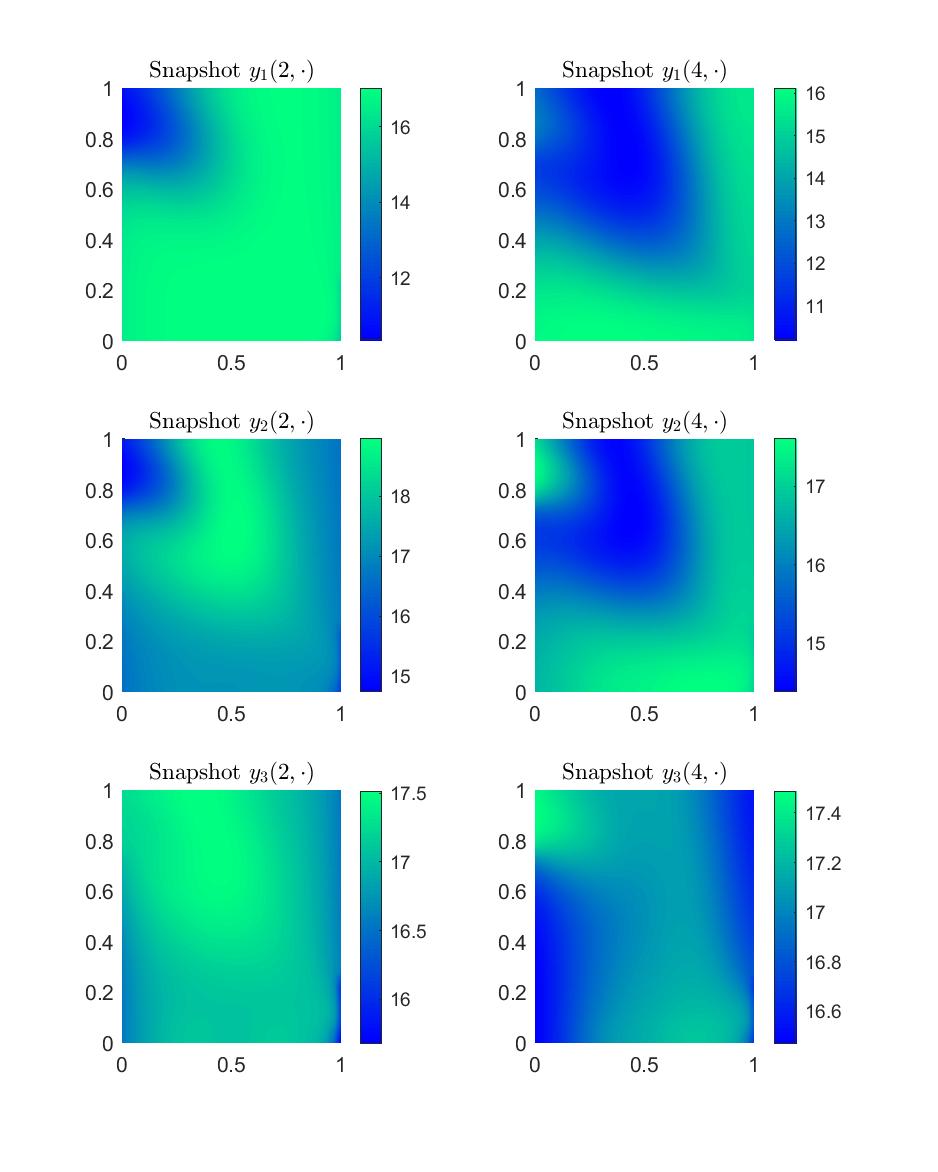}
        \caption{Run~\ref{ex:generalNumericSetup}. Solution snapshots of the solutions to \eqref{SIAM:EqMotPDE1} for the fixed control inputs $u_1$, $u_2$, and $u_3$.}
        \label{fig:advectionsDiffusionSolutions_snapshots}
    \end{figure}
    We can see that the advective cooldown effect from the outside temperature entering at $\Gamma_1$ is most prevalent for the uncontrolled case. For the time-constant control $u_2$, it is also still visible, but we can also recognise the heating effect taking over at $t=4$ in $\Gamma_1$ and $\Omega_1$. For the periodical control $u_3$, it can be seen that the advection effect is almost entirely counteracted. The snapshot at $t=4$ also identifies the regions close to the outside walls on the left and right as those which can not be heated as effectively with the current setup.\hfill$\Diamond$
\end{run}
	
We will later see controls which appear as a result from optimization algorithms applied to the optimal control problem \eqref{generalNumericalSetup_ocp}. It will interesting to see if they show a similar periodicity than our experimental candidate $u_3$. 
	
\section{Outline of the work}
\label{SIAM-Book:Section1.4}

The work is organized as follows. In Chapter~\ref{SIAM:Section-2} we introduce the POD method in finite and infinite-dimensional Hilbert spaces and discuss various applications. Chapter~\ref{SIAM-Book:Section3} is devoted to to POD-based Galerkin schemes for evolution problems. Mainly, we study linear problems taking different discretization methods into account. We provide a certified a-priori and a-posteriori error analysis. Furthermore, the numerical realizations are explained and illustrated by test examples. Quadratic programming problems governed by liner evolution problems are investigated in Chapter~\ref{SIAM-Book:Section4}. As in Chapter~\ref{SIAM-Book:Section3} we present a certified a-priori and a-posteriori error analysis. Moreover, we discuss basis update strategies. In Chapter~\ref{Advanced Topics in POD Suboptimal control} we give an outlook to further directions in reduced-order modeling in optimal control and optimization. More precisely, a nonlinear optimal control problem is studied. Moreover, state-constrained optimization problems are solved by a tailored combination of primal-dual active set methods and POD-aesed reduced-order modeling. Furthermore, POD Galerkin methods for multiobjective optimal control problems are investigated. Finally, some required theoretical foundations are summarized in the appendix.


\chapter{Proper Orthogonal Decomposition (POD)}
\label{SIAM:Section-2}

In this chapter we introduce proper orthogonal decomposition (POD) which is a method for model-order reduction. The chapter is organized as follows:

\begin{itemize}
    \item Section~\ref{SIAM:Section-2.1.1} introduces the discrete variant of the POD method. For a given, finite set of Hilbert space (data) vectors, we define and solve the problem of finding a low-dimensional POD subspace in which these vectors can be best approximated in a least-squares sense. We find that this subspace is spanned by the eigenvectors of a particular linear operator associated with the data vectors. The according eigenvalues are used to define an energy term which may be used as an indicator to find a suitable dimension of the POD subspace. We also consider some more concrete situations, which arises in numerical realizations, for the underlying problems such as unitary, Euclidean and finite-dimensional spaces as well as POD for finite-dimensional dynamical systems.
    \item In Section~\ref{SIAM:Section-2.1.2} we face the continuous variant of the POD method. In this case, there are infinitely many data vectors that depend on a variable (or parameter) $\bmu$ from a predetermined set. Analogous to Section~\ref{SIAM:Section-2.1.1} we derive a formulation to approximate all data vectors in a least-squares optimal way by a low-dimensional linear subspace. Moreover, the basis vectors of this space take the form of another linear operator which bears a close resemblance to the operator from the discrete variant of the POD method. 
    \item Section~\ref{Section:Perturbation analysis for the POD basis} concerns with the perturbation analysis for the POD basis. It deals with the case, where a continuous data trajectory as in Section~\ref{SIAM:Section-2.1.2} is discretized into a discrete set as in Section~\ref{SIAM:Section-2.1.1}. We compute the continuous and discrete POD bases to both trajectories and prove a stability result: If the discrete data trajectory approximates the continuous one in the limit, the discrete POD basis converges to the continuous one as well.
    \item In Section~\ref{Section:PODHilbert}, we discuss the specific choice when the underlying Hilbert space is part of a so-called Gelfand triple, a case which often occurs in evolution equations. It turns out that in this case it is possible to derive rate of convergence results which will be very useful in our a-priori error analysis carried out in Sections~\ref{SIAM-Book:Section3} and \ref{SIAM-Book:Section4}.
    \item For better readability all proofs from this chapter have been moved to Section~\ref{SIAM:Section-2.6}.
\end{itemize}

Throughout this chapter we suppose that -- unless stated otherwise -- $X$ is a separable complex Hilbert space endowed with the inner product $\langle \cdot\,,\cdot\rangle_X$ and induced norm $\|\cdot\|_X=\langle\cdot\,,\cdot\rangle_X^{1/2}$.

\section{The discrete variant of the POD method}
\label{SIAM:Section-2.1.1}
\setcounter{equation}{0}

Let us assume that we are given a snapshot space $\mathscr V^n$ spanned by a finite number of vectors which belong to a separable complex Hilbert space. Then $\mathscr V^n$ has a finite dimension $d \in \mathbb N$. The goal is to find an orthonormal set $V^\ell=\Span\{\psi_1,\ldots,\psi_\ell\}\subset\mathscr V^n$ with pairwise orthonormal vectors $\{\psi_i\}_{i=1}^\ell$ such that the following properties hold:
\begin{itemize}
    \item All elements in $\mathscr V^n$ can be approximated sufficiently accurate by a linear combination of the vectors $\psi_1,\ldots,\psi_\ell$.
    \item The dimension $\ell$ of $V^\ell$ is as small as possible.
\end{itemize}
In this section we show that the POD method offers a way to compute the $\{\psi_i\}_{i=1}^\ell$ in an optimal manner. It turns out that this approach is very much related to the eigendecomposition of a specific linear operator.

\subsection{Problem formulation and main result}
\label{SIAM:Section-2.1.1.1}

For fixed $n,K\in\mathbb N$ let the so-called {\em snapshots}\index{POD method!discrete variant!snapshots} $y_1^k,\ldots,y_n^k\in X$ be given for $1\le k\le{K}$. To avoid a trivial case we suppose that at least one of the $y_j^k$'s is non-zero. Then we introduce the finite dimensional, linear {\em snapshot space}\index{POD method!discrete variant!snapshot space}
\begin{equation}
    \label{SIAM:Eq-I.1.1.1}
    \mathscr V^n=\mathrm{span}\,\Big\{y_j^k\,|\,1\le j\le n\text{ and } 1\le k\le{K}\Big\}\subset X
\end{equation}
with finite dimension $d^n\in\{1,\ldots,nK\}$. In Section~\ref{Section:Perturbation analysis for the POD basis} we will consider the case, where the number $n$ tends to infinity. Therefore, we emphasize this dependence here by using the superscript $n$. We distinguish two cases:
\begin{itemize}
    \item [1)] The separable Hilbert space $X$ has finite dimension $m$: Then $X$ is isomorphic to $\mathbb C^m$; see, e.g., \cite[p.~47]{RS80}. We define the finite index set $\mathbb I=\{1,\ldots,m\}$. Clearly, we have $d^n\le\min(n{K},m)$. Especially in case of $X=\mathbb C^m$, the snapshots $y_j^k=(y_{ij}^k)_{1\le i\le m}$ are vectors in $\mathbb C^m$.
    \item [2)] $X$ is infinite-dimensional: Since $X$ is separable, each orthonormal basis of $X$ has countably many elements. In this case $X$ is isomorphic to the set $\ell_2$ of sequences $\{x_i\}_{i\in\mathbb N}$ of complex numbers which satisfy $\sum_{i=1}^\infty |x_i|^2<\infty$; see \cite[p.~47]{RS80}, for instance. The index set $\mathbb I$ is now the countable, but infinite set $\mathbb N$.
\end{itemize}
Recall that a \emph{complete} or \emph{maximal basis} $\{\psi_i\}_{i\in\mathbb I}$ in a separable Hilbert space $X$ satisfies for any $\varphi\in X$:
\begin{align*}
    {\langle\varphi,\psi_i\rangle}_X=0\text{ for all }i\in\mathbb I\quad\Longleftrightarrow\quad\varphi=0.
\end{align*}

In this section we are facing the following problem.

\begin{problem}[Discrete POD method]
    \label{ProblemPOD}
    The {\em discrete method of POD}\index{POD method!discrete variant} consists in choosing a complete orthonormal basis $\{\psi_i\}_{i\in \mathbb I}$ in $X$ such that for every $\ell\in\{1,\ldots,d^n\}$ the mean square error between the $n{K}$ elements $y_j^k$ and their corresponding $\ell$-th partial Fourier sum is minimized on average:
    \begin{equation}
        \tag{$\mathbf P^\ell_n$}
        \label{SIAM:Eq-I.1.1.2}
        \min\sum_{k=1}^{K} \omega_k^{K}\sum_{j=1}^n\alpha_j^n \Big\| y_j^k-\sum_{i=1}^\ell{\langle y_j^k,\psi_i\rangle}_X\,\psi_i\Big\|_X^2\text{ s.t. } \{\psi_i\}_{i=1}^\ell\subset X\text{ and }{\langle\psi_i,\psi_j\rangle}_X=\delta_{ij},~1 \le i,j \le \ell,
    \end{equation}
    where $\omega_k^{K}$ and $\alpha_j^n$ are positive weighting parameters for $k=1,\dots,{K}$ and $j=1,\dots,n$. Here, the symbol $\delta_{ij}$ denotes the Kronecker symbol satisfying $\delta_{ii}=1$ and $\delta_{ij}=0$ for $i\neq j$. 
\end{problem}

\begin{definition}
    An optimal solution $\{\hat\psi_i^n\}_{i=1}^\ell$ to Problem~{\em\ref{ProblemPOD}} is called a \index{POD method!discrete variant!POD basis of rank $\ell$}{\em POD basis of rank $\ell$}.
\end{definition}

\begin{remark}
    \rm To ease the presentation the weights $\alpha_j^n$ are chosen independently of $k\in\{1,\ldots,K\}$. Moreover, $n$ could depend on $k\in\{1,\ldots,{K}\}$ as well, with a generalization of the presented theory and algorithms being straightforward.\hfill$\blacksquare$
\end{remark}

\noindent
\textbf{Main result:} In Theorem~\ref{SIAM:Theorem-I.1.1.3} we prove that for every $\ell\in\{1,\ldots,d^n\}$ a solution $\{\hat\psi^n_i\}_{i=1}^\ell$ to \eqref{SIAM:Eq-I.1.1.2} is characterized by the eigenvalue problem
\begin{equation}
    \label{EigPODPro}
    \mathcal R^n\hat\psi_i^n=\hat\lambda_i^n\hat\psi_i^n\quad\text{for }1\le i\le\ell,
\end{equation}
where $\hat\lambda_1^n\ge\ldots\ge\hat\lambda_{d^n}^n>0$ denote the largest positive eigenvalues of the linear, bounded, non-negative and self-adjoint operator
\begin{equation}
    \label{OperatorR}
    \mathcal R^n: X \to X, \quad\mathcal R^n\psi=\sum_{k=1}^{K}\omega_k^{K}\sum_{j=1}^n\alpha_j^n\,{\langle\psi,y_j^k\rangle}_X y_j^k\quad\text{for }\psi\in X.
\end{equation}
Furthermore, we obtain the following approximation {\em error formula}\index{POD method!error formula}
\begin{align*}
    \sum_{k=1}^{K} \omega_k^{K}\sum_{j=1}^n\alpha_j^n \Big\| y_j^k-\sum_{i=1}^\ell{\langle y_j^k,\hat\psi_i^n\rangle}_X\,\hat\psi_i^n\Big\|_X^2=\sum_{i=\ell+1}^{d^n}\hat\lambda_i^n.
\end{align*}

\begin{remark}
    \rm The decay of the eigenvalues $\{\hat\lambda_i^n\}_{i=1}^{d^n}$ plays an essential role for a successful application of the POD method. In general, one has to utilize a complete orthonormal basis in $X$ to represent elements in the snapshot space $\mathscr V^n$ by their Fourier sum. This leads to a high-dimensional or even infinite-dimensional approximation scheme. Nevertheless, if the finite sum $\sum_{i=\ell+1}^{d^n}\hat\lambda_i^n$ is sufficiently small for a not too large $\ell$, elements in the subspace $\mathscr V^n$ can be approximated by a linear combination of the few basis elements $\{\hat\psi^n_i\}_{i=1}^\ell$. This offers the chance to reduce the number of terms in the Fourier series utilizing the POD basis of rank $\ell$, as shown in the following examples.\hfill$\blacksquare$
\end{remark}

\subsection{Three motivating examples}
\label{sec:threeExamples}

Before we focus on the mathematical theory we illustrate the POD method for three different examples. In the first one we show how to compress the information of a given function, which is depending on a time and a spatial variable, by using POD and apply this to the case of three functions being composed of sine and cosine functions.
 
\begin{example}[Three simple functions]
    \label{Example:cosExample}
    \rm We introduce a regular function $y:Q \to \mathbb R$, where $Q= (0,T)\times\Omega$ and $\Omega = (\bx_\mathsf a,\bx_\mathsf b) \subset \mathbb R$ are given with $T>0$ and $\bx_\mathsf a<\bx_\mathsf b$. These sets are equidistantly discretized into $n$ and $n_\bx$ points, respectively, as $0=t_1<\hdots < t_n = T$ and $\bx_\mathsf a=\bx_1<\ldots<\bx_{n_\bx}=\bx_\mathsf b$ with constant step sizes
    \begin{align*}
        \Delta t = \frac{T}{n-1} > 0, \quad h=\frac{\bx_\mathsf b - \bx_\mathsf a}{n_\bx-1}> 0.
    \end{align*}
    The function $y$ is approximated pointwise by introducting the $n$ vectors $y_j=(\mathrm y_{ij})_{1\le i\le n_\bx} \in \mathbb R^{n_\bx}$ with $\mathrm y_{ij}=y(t_j,\bx_i)$ for $i=1,...,n_\bx$ and $j=1,...,n$. In $X = \mathbb R^{n_\bx}$ we introduce a weighted inner product to represent a trapezoidal integration scheme over the $\bx$-variable :
    \begin{align*}
        {\langle\mathrm u,\mathrm v\rangle}_X= h\bigg(\frac{1}{2}\,\mathrm u_1\mathrm v_1+ \sum_{i=2}^{n_\bx-1}\mathrm u_i\mathrm v_i+ \frac{1}{2}\,\mathrm u_{n_\bx}\mathrm v_{n_\bx}\bigg)=\mathrm u^\top\bW\mathrm v\quad \text{for } \mathrm u,\mathrm v\in \mathbb R^{n_\bx}
    \end{align*}
    with the weighting matrix $\bW=h\,\mathrm{diag}\,(1/2,1,\ldots,1,1/2)\in\mathbb R^{n_\bx\times n_\bx}$. The associated norm is given as $\|\cdot\|_X=\langle\cdot\,,\cdot\rangle_X^{1/2}$. Similarly, we define the trapezoidal weights for the time variable as
    \begin{equation}
        \label{TrapW}
        \alpha_1^n = \frac{\Delta t}{2}, \quad \alpha_j^n = \Delta t \text{ for } j=2,...,n-1, \quad \alpha_n^n = \frac{\Delta t}{2}
    \end{equation}
    and consider the following instance of \eqref{SIAM:Eq-I.1.1.2}:
    \begin{equation}
        \label{cosExample_podProb}
        \min \sum_{j=1}^n \alpha_j^n\,\Big\| y_j-\sum_{i=1}^\ell{\langle y_j,\psi_i\rangle}_X\,\psi_i\Big\|_X^2\text{ s.t. }\{\psi_i\}_{i=1}^\ell\subset X\text{ and }{\langle\psi_i,\psi_j\rangle}_X=\delta_{ij},~1 \le i,j \le \ell.
    \end{equation}
    Problem~\eqref{cosExample_podProb} translates to finding $X$-orthonormal vectors $\psi_1,...,\psi_\ell \in X$ that approximate the discretized trajectory $\{y_1,...,y_n\} \subset X$ of the function $y$ as best as possible in a least-square-sense. We call $\{y_1^\ell,\hdots,y_n^\ell\} \subset \mathbb R^{n_\bx}$ given by $y_j^\ell = \sum_{i=1}^\ell \langle y_j,\psi_i \rangle_X \psi_i$ the $\ell$-th POD approximation of $y$. In our numerical example, we consider the three functions 
    \begin{align*}
        y_1: (0,2\pi) \times (0,2\pi) \to \mathbb R, \quad y_1(t,\bx)&=\cos(t)\cos(\bx), \\
        y_2: (0,2\pi) \times (0,2\pi) \to \mathbb R, \quad y_2(t,\bx)&=\cos(t+\bx), \\
        y_3: (0,1) \times (0,2\pi) \to \mathbb R, \quad y_3(t,\bx)&=\cos(t\bx).
    \end{align*}
    Each of these functions is discretized as described above and various POD approximations are computed by solving \eqref{cosExample_podProb} with methods that will be explained in detail in Section~\ref{SIAM:Section-2.1.1.2}. We denote the solutions by $\hat\psi_i^n$ for $i=1,...,\ell$.
    \begin{figure}
        \hspace{10mm}\includegraphics[height=50mm,width=0.9\textwidth]{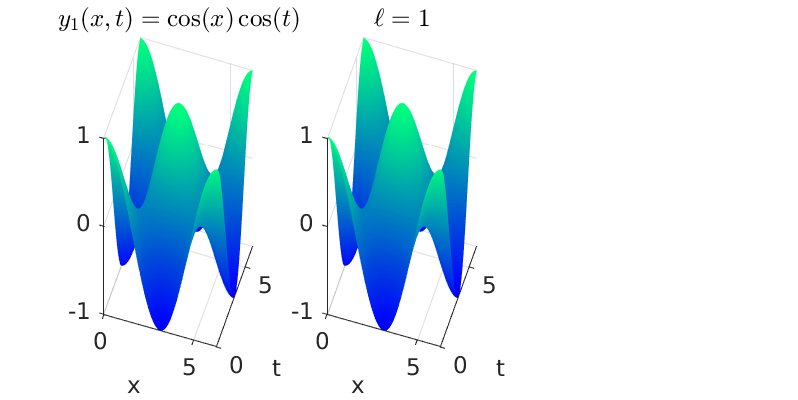}    
        \caption{Example~\ref{Example:cosExample}. Function $y_1(t,\bx)=\cos(t)\cos(\bx)$ (left plot); projection of $y_1$ onto the POD space $V^\ell=\Span\{\hat\psi_1^n,\ldots,\hat\psi_\ell^n\}$ for $\ell=1$ (right plot).}
        \label{fig:cosExample_podApproximation1}
    \end{figure}
    In Figure~\ref{fig:cosExample_podApproximation1} we observe that only one basis function $\hat\psi_1^n$ is enough to approximate the trajectory of $y_1$. This is due to the fact that the spatial and temporal variable are separated, meaning that the trajectory of $y_1$ is in fact one-dimensional, so each snapshot can be represented by a simple scaling of the $\cos(\bx)$ function. 
    \begin{figure}
        \hspace{5mm}\includegraphics[height=50mm,width=160mm]{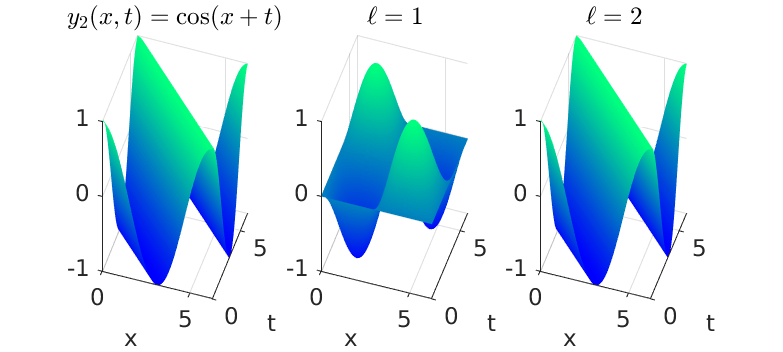}
        \caption{Example~\ref{Example:cosExample}. Function $y_2(t,\bx)=\cos(t+\bx)$ (left plot); projection of $y_2$ onto the POD space $V^\ell=\Span\{\hat\psi_1^n,\ldots,\hat\psi_\ell^n\}$ for $\ell=1$ (middle plot) and $\ell=2$ (right plot).}
        \label{fig:cosExample_podApproximation2}
    \end{figure}
    For $y_2$, two basis functions $\hat\psi_1^n$, $\hat\psi_2^n \in X$ are required to represent the trajectory; see Figure~\ref{fig:cosExample_podApproximation2}. The reason for this lies in the addition theorem for trigonometrical functions: 
    \begin{align*}
        \cos(t+\bx) = \cos(t)\cos(\bx) - \sin(t)\sin(\bx)
    \end{align*}
    The two POD bases functions $\hat\psi_1^n$ and $\hat\psi_2^n$ together span the same subspace as $\cos$ and $\sin$, the trajectory can be fully represented by its POD approximations for $\ell=2$.
    \begin{figure}
        \hspace{15mm}\includegraphics[height=50mm,width=140mm]{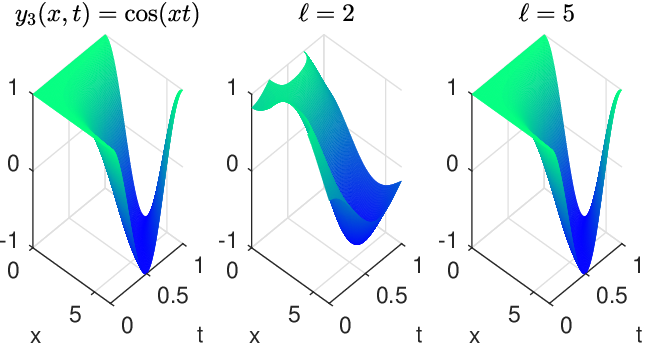}
        \caption{Example~\ref{Example:cosExample}. Function $y_3(t,\bx)=\cos(t\bx)$ (left plot); projection of $y_3$ onto the POD space $V^\ell=\Span\{\hat\psi_1^n,\ldots,\hat\psi_\ell^n\}$ for $\ell=2$ (middle plot) and $\ell=5$ (right plot).}
        \label{fig:cosExample_podApproximation3}
    \end{figure}
    In Figure~\ref{fig:cosExample_podApproximation3} we observe that the qualitative behaviour of $y_3$ is represented quite well by its POD approximation for $\ell=5$, while its approximation for $\ell=2$ lacks the sharp features of the original trajectory. There is no analytical property that allows us to decompose the term $\cos(t\bx)$ as it was done for $y_1$ and $y_2$. Hence, more basis functions are required to sufficiently capture the behaviour of $y_3$, as can also be seen in Figure~\ref{fig:cosExample_podApproximationError}.
    \begin{figure}
        \begin{center}
            \includegraphics[height=50mm]{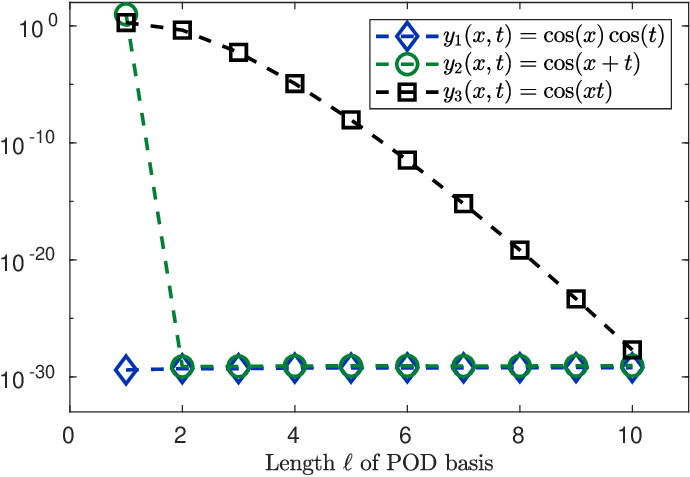}
        \end{center}
        \caption{Example~\ref{Example:cosExample}. Decay of the POD approximation error $\sum_{i=\ell+1}^{d^n}\hat\lambda_i^n$ for the three functions $y_1$, $y_2$ and $y_3$.}
        \label{fig:cosExample_podApproximationError}
    \end{figure}
    We refer to \cite[Section~3.3.3]{HLBR12}, where sufficient conditions are given so that the POD bases are Fourier modes.\hfill$\blacklozenge$
\end{example}

For the second example, we return to the parabolic evolution problem from Section~\ref{SIAM-Book:Section1.2} and present a typical situation how the discrete POD problem \eqref{SIAM:Eq-I.1.1.2} might arise in such a setting. 

\begin{example}[Evolution problem]
    \label{Example:PODParabolic}
    \rm Let us consider the parabolic problem introduced in Section~\ref{SIAM-Book:Section1.2} for ${K}$ different initial conditions $\{y_\circ^k\}_{k=1}^{K}$ In particular, we choose functions which have a doubled Gaussian temperature distribution added to the otherwise constant initial condition of $y_\circ(\bx)=17$ for all $\bx \in \Omega=(0,1)^2$:
    \begin{align*}
	   y_\circ^{i,j}(\bx) = 17 + \frac{2}{\sqrt{2 \pi} \sigma_j} \exp \bigg( - \frac{|\bx-p^{(i)}|_2^2}{2 \sigma_j^2} \bigg).
    \end{align*} 
    Here, $|\cdot|_2$ denotes the Euclidian norm, $p^{(i)} \in \mathbb R^2$ are the distribution peaks and $\sigma_i>0$ are different standard deviations. The latter are chosen as
    \begin{align*}
        p^{(i)}_1 = 0.15, ~p^{(i)}_2 \in\big\{0.1,0.3,0.5,0.7.0.9\big\}, \quad \sigma_j \in \big\{ 0.05,0.1,0.5,1\big\}.
    \end{align*}
    All possible combinations of $p^{(i)}$ and $\sigma_j$ give us $K=20$ different initial conditions. Since the initial peak temperature is always located on the left-hand side and the advective flow depicted in Figure~\ref{fig:navierStokesSolution} is oriented from the top left to the lower right, this peak will be transported to different levels of effectiveness by the advection effect. For the discretization of the time variable, we choose a temporal grid $0=t_1<t_2<\ldots<t_n=T$ with the step sizes $\delta t_j=t_j-t_{j-1}$ for $j=2,\ldots,n$. Let
    \begin{align*}
        \delta t=\min_{2\le j\le n}\delta t_j\quad\text{and}\quad\Delta t=\max_{2\le j\le n}\delta t_j
    \end{align*}
    be the minimal and maximal time step size, respectively. By $y_j^k:\Omega\to\mathbb R$, $1\le j\le n$ and $1\le k\le{K}$, we denote the solution to \eqref{SIAM:EqMotPDE1} at the time instance $t_j$ for initial condition $y^k_\circ$ and fixed control $u\in\Uad$. For $\alpha_j^n$ we take trapezoidal weights:
    \begin{align*}
        \alpha_1^n=\frac{\delta t_1}{2},\quad\alpha_j^n=\frac{\delta t_j+\delta t_{j-1}}{2}\text{ for }j=2,\ldots,n-1,\quad\alpha_n^n=\frac{\delta t_n}{2}.
    \end{align*}
    Then we can interpret the term
    \begin{align*}
        \sum_{j=1}^n\alpha_j^n\,\Big\| y_j^k-\sum_{i=1}^\ell{\langle y_j^k,\psi_i\rangle}_X\,\psi_i\Big\|_X^2
    \end{align*}
    as a trapezoidal approximation of the integral term
    \begin{align*}
        \int_0^T\Big\| y^k(t,\cdot)-\sum_{i=1}^\ell{\langle y_j^k(t,\cdot),\psi_i\rangle}_X\,\psi_i\Big\|_X^2\,\mathrm dt.
    \end{align*}
    The choice of the quadrature weights should be related to the (discretization) error of the approximations $y_j^k$. For example, the trapezoidal rule leads to an error of the size $\mathcal O(\Delta t^2)$, but if the discretization error is of higher order, then other quadrature rules, like e.g. the Simpson rule, should be utilized. A usual choice for the Hilbert space $X$ is the Lebesgue space $H=L^2(\Omega)=L^2(\Omega;\mathbb R)$; cf., Definition~\ref{Definition_L2_BanachValued}. However, we can also choose $X=V$ with the Sobolev space $V=H^1(\Omega)=H^1(\Omega;\mathbb R)$; cf., Definition~\ref{Definition_H1_BanachValued}. Since we do not want to weight the solutions for distinct initial conditions differently, we take $\omega_k^{K}=1$. Figure~\ref{fig:podInitialConditions_singVals} shows the decay of the first $25$ eigenvalues $\hat\lambda_i^n$ for the cases $X=V$ and $X=H$.
    \begin{figure}
        \centering
        \includegraphics[height=50mm]{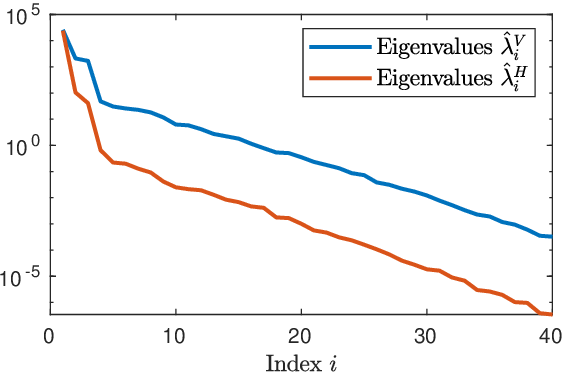}
        \caption{Example~\ref{Example:PODParabolic}. Eigenvalue decay with the choices $X=V=H^1(\Omega)$ and $X=H=L^2(\Omega)$.}
        \label{fig:podInitialConditions_singVals}
    \end{figure}
    Concentrating at this time only on the $V$-line, we can see that by the tenth eigenvalue, the order of magnitude has decreased from roughly $10^4$ to $10^1$. At the $20$-th eigenvalue, we have arrived at $10^{-1}$ and consideration of further eigenvectors is not necessary to capture the essential information contained in all $K=20$ snapshots trajectories. In Figure~\ref{fig:podInitialConditions_podVectors}, we can observe the first six orthonormal eigenvectors $\psi_1,\ldots,\psi_6 \in V = H^1(\Omega)$ for $X=V$.
    \begin{figure}
	   \centering
	   \includegraphics[width=110mm,height=120mm]{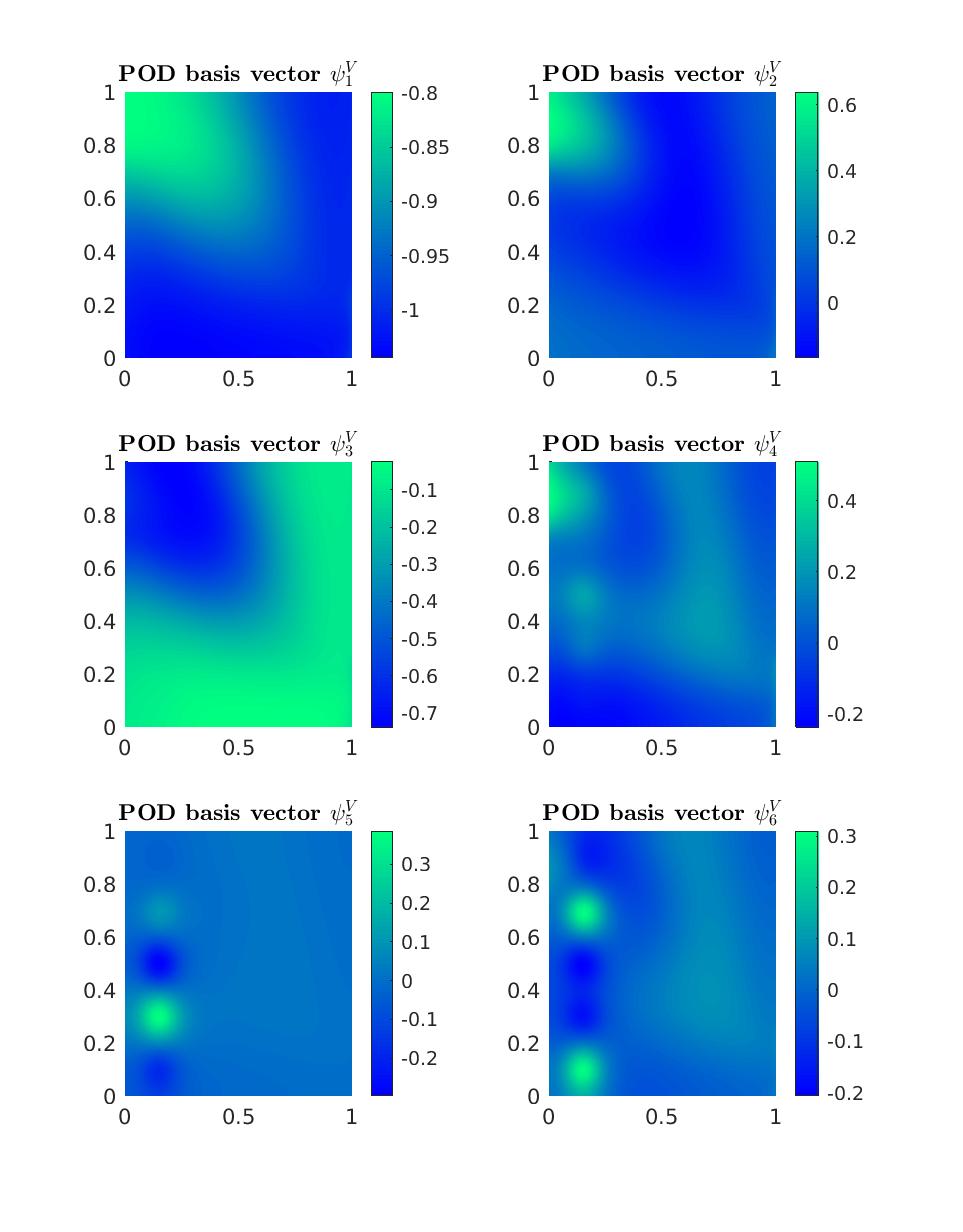}
	   \caption{Example~\ref{Example:PODParabolic}. First six POD basis vectors $\hat\psi_i^n$ with $X=V=H^1(\Omega)$.}
	   \label{fig:podInitialConditions_podVectors}
    \end{figure}
    We can see that the first three functions mainly capture the advection effect which is oriented from the upper left corner $\Gamma_1$ to the lower right corner $\Gamma_2$. After this, we can start to observe some circles on the left-hand sides of $\psi_5$ and $\psi_6$ which correspond to the intensity peaks from the initial conditions located between $p^{(1)} = (0.1,0.1)^\top$ and $p^{(5)} = (0.1,0.9)^\top$. The fact the ``advective'' POD vectors come before those belonging to the initial condition together with the fast decay in eigenvalues after $\ell=3$ tells us that the main challenge for the POD basis is not the varying initial conditions in this case, but the advection part. Note that this effect is only strengthened by the low outside temperature which may enter very effectively at the upper-left, poorly insulated inflow segment $\Gamma_1$. Finally, let us mention that one could also discuss different control inputs $u\in\Uad$ instead of, or in addition to, the various initial conditions.\hfill$\blacklozenge$
\end{example}

The third motivating example for the discrete POD method takes the form of a generic elliptic equation arising in the field of acoustics.

\begin{example}[Elliptic problem]
    \label{Example:PODElliptic}
    \rm The absorption behavior of trim parts inside a car is described by the frequency dependent impedance $Z(f)\in\mathbb C$. The acoustical simulations of car interior noise is governed by the well-known Helmholtz equation. A very important parameter in this model is the frequency-dependent impedance $Z(f)\in\mathbb C$ (see left plot of Figure~\ref{OmegaImpedance}), which the frequency-dependent absorption behavior of trim parts inside a car. Let a simplified interior of a car be given by the two-dimensional cross section $\Omega\subset\mathbb R^2$ plotted in the right plot of Figure~\ref{OmegaImpedance}; cf. \cite{LV14,Vol10,VH08}. 
    \begin{figure}
        \begin{center}
            \includegraphics[height=50mm]{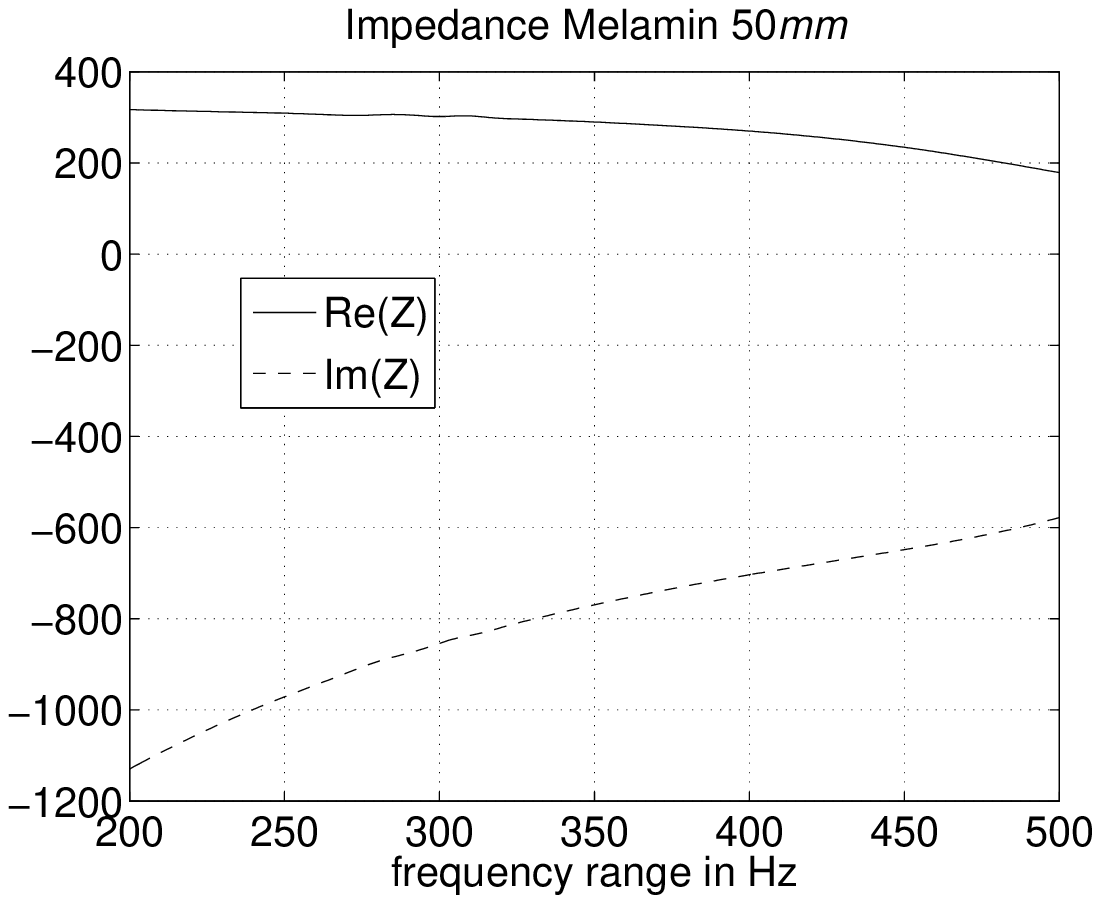}\hspace{10mm}
            \includegraphics[height=50mm]{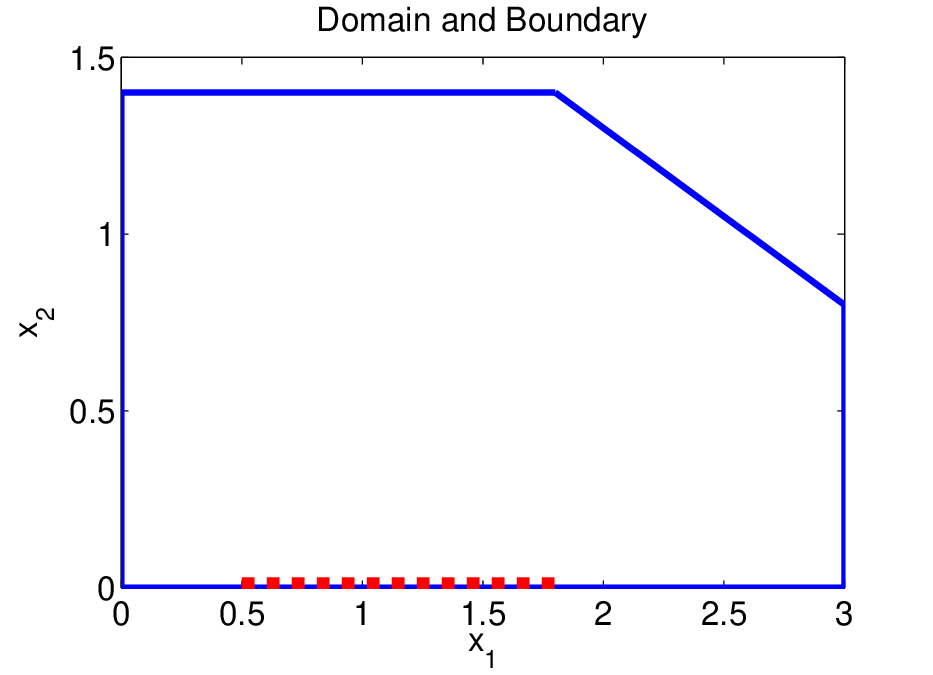}
        \end{center}
        \caption{Example~\ref{Example:PODElliptic}. Real and imaginary part of the impedance $Z(f)$ for the damping material Melamin {\em 50\,mm} over the frequency interval from {\em 200} to {\em 500\,Hz} (left plot); spatial domain $\Omega$ where $\Gamma_\mathrm R$ is shown by the dashed line ant the bottom (right plot); cf \cite{LV14,Vol10,VH08}.}
        \label{OmegaImpedance}
    \end{figure}
    The boundary $\Gamma=\partial\Omega$ is split into the disjoint parts $\Gamma_\mathrm R=\{(x_1,0)\in\mathbb R^2\,|\,0.5\le x_1\le 1,75\}$ (impedance part) and $\Gamma_\mathrm N=\Gamma\setminus\Gamma_\mathrm R$ (rigid body part). Furthermore, let the constants
    \begin{equation}
        \label{eq:constants}
        c=343.799\,\frac{\mathrm m}{\mathrm s},\quad k(f)=\frac{2\pi f}{c},\quad\varrho_\circ=1.19985\,\frac{\mathrm{kg}}{\mathrm m^3},\quad\omega(f)=2\pi f=ck(f)
    \end{equation}
    be given, where $k(f)$ and $\omega(f)$ are called the frequency-dependent {\em wave number}\index{Equation!Helmholtz!wave number} and {\em angular frequency}\index{Equation!Helmholtz!angular frequency}, respectively. In the right plot of Figure~\ref{OmegaImpedance} we see the complex-valued \index{Equation!Helmholtz!impedance}{\em impedance} for the fire-resistant damping material Melamin $50\,\mathrm{mm}$ in the frequency range from $200$ to $500\,\mathrm{Hz}$. The complex-valued \index{Equation!Helmholtz!sound pressure}{\em sound pressure} $p(f)=p(\cdot\,;f):\Omega\to\mathbb C$ is governed by the {\em Helmholtz equation}\index{Equation!Helmholtz}, which is the elliptic boundary value problem
    \begin{equation}
        \label{eq:Helmholtz-S2}
        \begin{aligned}
            -\Delta p(f)-k(f)^2p(f)&=q(f) && \text{in } \Omega,\\
            \frac{\jmath}{\varrho_\circ \omega(f)} \, \frac{\partial p(f)}{\partial\bn}&= 0 && \text{on } \Gamma_\mathrm N,\\
            \frac{\jmath}{\varrho_\circ \omega(f)} \, \frac{\partial p(f)}{\partial\bn}&= \frac{p(f)}{Z(f)}&& \text{on }  \Gamma_\mathrm R
        \end{aligned}
    \end{equation}
    for every frequency $f \in \mathscr F=[200,500]$. In \eqref{eq:Helmholtz-S2} the function $q(f)$ stands for the $f$-dependent source term modelling the excitation at the point $\bx_q= (0.28,1.21) \in \Omega$ with the frequency $f$:
    \begin{align*}
        q(\bx;f)= \frac{1}{5}\exp\left(\frac{\jmath\pi(f-200)}{50}\right) \exp\left(-50\,{|\bx-\bx_q|}_2^2\right)\quad\text{for }\bx\in\Omega,~f\in\mathscr F.
    \end{align*}
    The existence of a weak solution to \eqref{eq:Helmholtz-S2} follows from the Fredholm alternative; see, e.g., \cite[pp.~640-644]{Eva08}. For more details we refer to \cite[Theorem~2.2]{Vol10}.\hfill\\
    Next we apply the discrete POD method with ${K}=1$ and $\omega_1^K=1$. For $n=301$, we choose the equidistant frequency grid $f_j=199+j$ for $j=1,\ldots,n$. We utilize the (trapezoidal) weights
    \begin{align*}
        \alpha_1^n=\frac{1}{2},\quad\alpha_j^n=1\text{ for }j=2,\ldots,n-1,\quad\alpha_n^n=\frac{1}{2}.
    \end{align*}
    Let $p(f_j):\Omega\to\mathbb C$, $j=1,\ldots,n$, be the (approximate) weak solutions to \eqref{eq:Helmholtz-S2} for the frequency $f_j\in\mathscr F$. In Example~\ref{SIAM:Example-I.1.1.1} we will utilize piecewise linear finite elements to compute numerical solutions. Similarly to Example~\ref{Example:PODParabolic} a natural choice for the Hilbert space $X$ is the Lebesgue space $H=L^2(\Omega;\mathbb C)$. Of course, we could also choose $X=V$ with $V=H^1(\Omega;\mathbb C)$. The discrete POD method \eqref{SIAM:Eq-I.1.1.2} then takes the form
    \begin{equation}
        \label{SIAM:Eq:SoundPressurePOD}
        \left\{
        \begin{aligned}
            \min \quad &\sum_{j=1}^n\alpha_j^n \Big\| p(f_j)-\sum_{i=1}^\ell{\langle p(f_j),\psi_i\rangle}_H\,\psi_i\Big\|_H^2\\
            \text{s.t.} \quad &\{\psi_i\}_{i=1}^\ell\subset H\text{ and }{\langle\psi_i,\psi_j\rangle}_H=\delta_{ij},~1 \le i,j \le \ell.
        \end{aligned}
        \right.
    \end{equation}
    It corresponds to the determination of the subspace in which the solution of \eqref{eq:Helmholtz-S2} to each discrete frequency $f_j$ can approximated in the closest possible way.

    \noindent
    In Figure~\ref{ImpedanceEigen}, the decay of the normalized eigenvalues $\hat\lambda_i^n/\sum_{j=1}^{d^n}\hat\lambda_j^n\in[0,1]$ is presented for the two choices $X=H$ (left plot) and $X=V$ (right plot).
    \begin{figure}
        \begin{center}
            \includegraphics[height=50mm]{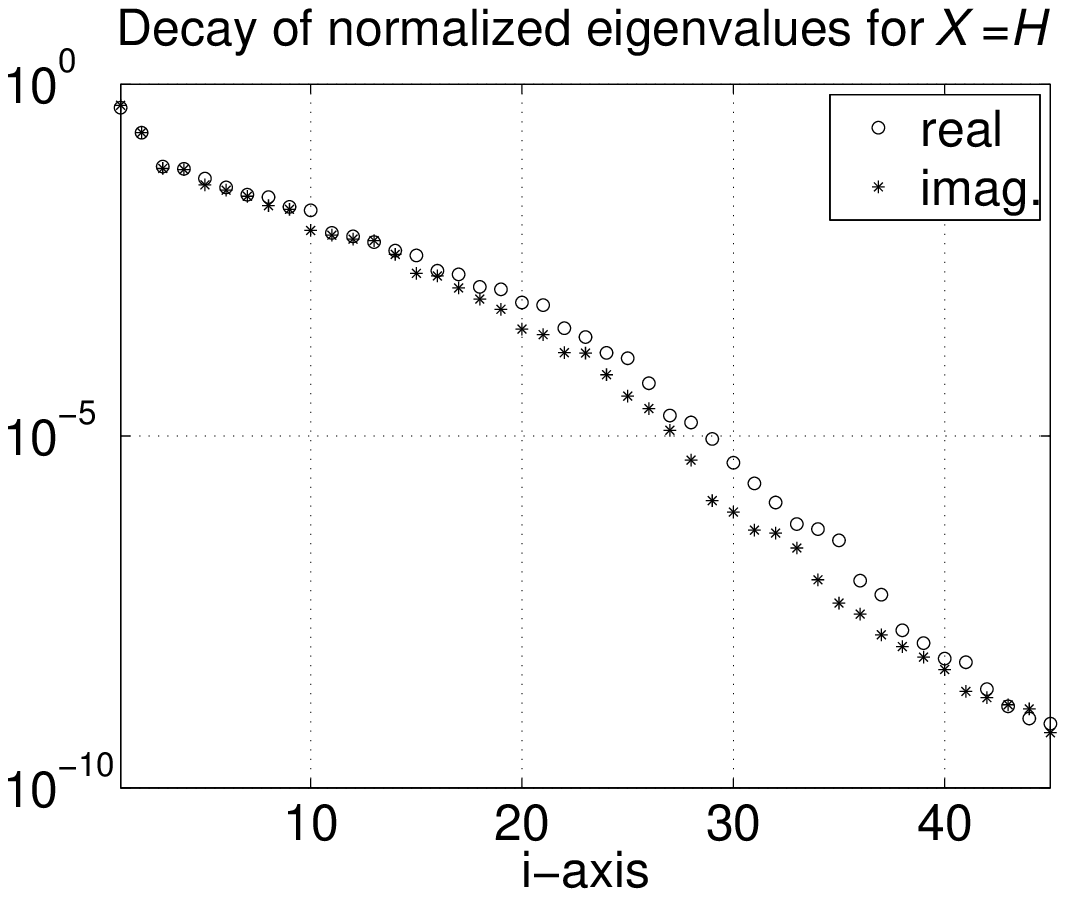}\hspace{10mm}\includegraphics[height=50mm]{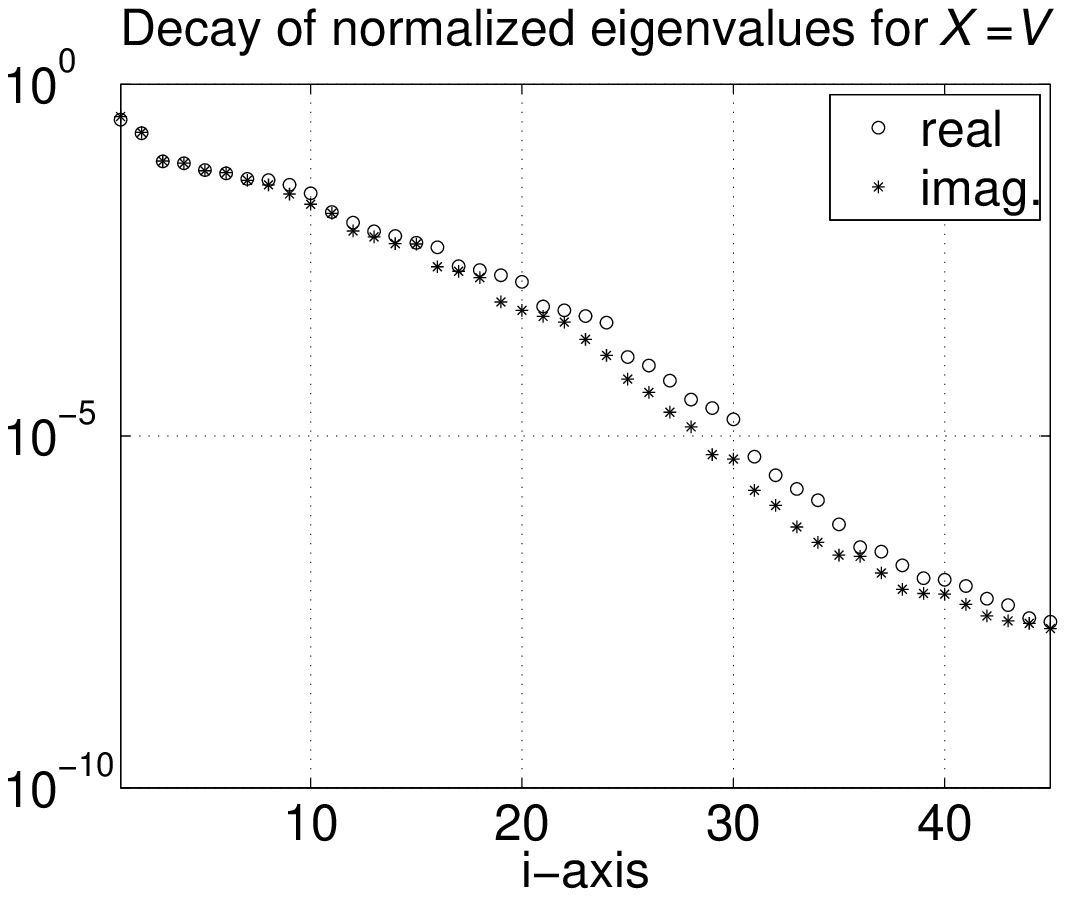}
        \end{center}
        \caption{Example~\ref{Example:PODElliptic}. Decay of the normalized eigenvalues $\hat\lambda_i^n/\sum_{j=1}^{d^n}\hat\lambda_j^n$ are presented for the two choices $X=H$ (left plot) and $X=V$ (right plot).}
        \label{ImpedanceEigen}
    \end{figure}
    We observe that the decay rates for the real and imaginary parts of the normalized eigenvalues behave similarly. Of course, the decay rate is slower for the stronger $V$-topology than for the weaker $H$-topology; compare Example~\ref{Example:PODParabolic} and Figure~\ref{fig:podInitialConditions_singVals}. We also refer to \cite{HDPV06} for more details.\hfill$\blacklozenge$
\end{example}

\begin{remark}
    \rm In Examples~\ref{Example:PODParabolic} and \ref{Example:PODElliptic} we have introduced two separable Hilbert spaces $H$ and $V$, which satisfy $V\subset H$. Furthermore, $V$ is continuously and densely embedded into $H$. This offers special properties for the POD method. We will address this issue in Section~\ref{Section:PODHilbert}.\hfill$\blacklozenge$
\end{remark}

\subsection{Derivation of the main result}

Recall the eigenvalue problem \eqref{EigPODPro} and the linear operator $\mathcal R^n$ introduced in \eqref{OperatorR}. In the next lemma we summarize properties of the operator $\mathcal R^n$ ensuring the existence of non-negative eigenvalues $\hat\lambda_i^n\ge\ldots\ge\hat\lambda_{d^n}^n>0$. A proof is given in Section~\ref{SIAM:Section-2.6.1}.

\begin{lemma}
    \label{SIAM:Lemma-I.1.1.1}
    Let $y_1^k,\ldots,y_n^k\in X$ be given snapshots for $1\le k\le{K}$. Define the linear operator $\mathcal R^n:X\to X$ as follows:
    \begin{equation}
        \label{SIAM:Eq-I.1.1.12}
        \mathcal R^n\psi=\sum_{k=1}^{K}\omega_k^{K}\sum_{j=1}^n\alpha_j^n\,{\langle \psi,y_j^k\rangle}_X\,y_j^k=\sum_{k=1}^{K}\omega_k^{K}\sum_{j=1}^n\alpha_j^n\,\overline{\langle y_j^k,\psi\rangle}_X\,y_j^k\quad\text{for }\psi \in X
    \end{equation}
    with positive weights $\omega_1^{K},\ldots,\omega_{K}^{K}$ and $\alpha_1^n,\ldots,\alpha_n^n$. Then $\mathcal R^n$ is a compact, self-adjoint and non-negative operator. 
\end{lemma}

\begin{remark}
    \rm Since $X$ is a separable complex Hilbert space and $\mathcal R^n:X\to X$ is a linear, compact, self-adjoint and non-negative operator, we can utilize the Riesz-Schauder and Hilbert-Schmidt theorems (Theorems~\ref{SIAM:Theorem-I.1.1.1} and \ref{SIAM:Theorem-I.1.1.2}, respectively) as well as Remark~\ref{SIAM:Remark-RealNonnegativeEigenvalue}: It follows that there exist a complete countable orthonormal basis $\{\hat\psi_i^n\}_{i\in\mathbb I}$ and a corresponding sequence of real eigenvalues $\{\hat\lambda_i^n\}_{i\in\mathbb I}$ satisfying
    \begin{equation}
        \label{SIAM:Eq-I.1.1.18}
        \mathcal R^n\hat\psi_i^n=\hat\lambda_i^n\hat\psi_i^n,\quad\hat\lambda_1^n\ge\ldots\ge\hat\lambda_{d^n}^n>\hat\lambda_{d^n+1}^n=\ldots =0.
    \end{equation}
    The spectrum of $\mathcal R^n$ is a pure point spectrum except for possibly $0$. Each non-zero eigenvalue of $\mathcal R^n$ has finite multiplicity and $0$ is the only possible accumulation point of the spectrum of $\mathcal R^n$.\hfill$\blacksquare$
\end{remark}

The next result is proved in Section~\ref{SIAM:Section-2.6.1}, too.

\begin{lemma}
    \label{Lem:PropRnOp}
    Let $y_1^k,\ldots,y_n^k\in X$ be given snapshots for $1\le k\le{K}$ and the linear operator $\mathcal R^n:X\to X$ be defined as in \eqref{SIAM:Eq-I.1.1.12}.
    \begin{enumerate}
        \item [\em 1)] We have
        \begin{equation}
            \label{SIAM:Eq-I.1.1.19}
            \sum_{k=1}^{K}\omega_k^{K}\sum_{j=1}^n\alpha_j^n\,\big|{\langle y_j^k,\hat\psi_i^n\rangle}_X\big|^2={\langle\mathcal R^n\hat\psi_i^n,\hat\psi_i^n\rangle}_X=\hat\lambda_i^n \quad\text{for any }i\in\mathbb I.
        \end{equation}
        In particular, it follows that
        \begin{equation}
            \label{SIAM:Eq-I.1.1.20}
            \sum_{k=1}^{K}\omega_k^{K}\sum_{j=1}^n\alpha_j^n\,\big|{\langle y_j^k,\hat\psi_i^n\rangle}_X\big|^2=0\quad\text{for all } i>d^n.
        \end{equation}
        \item [\em 2)] The identity
        \begin{equation}
            \label{SIAM:Eq-I.1.1.22}
            \sum_{k=1}^{K}\omega_k^{K}\sum_{j=1}^n\alpha_j^n\,{\|y_j^k\|}_X^2=\sum_{i=1}^{d^n}\hat\lambda_i^n.
        \end{equation}
        holds.
        \item [\em 3)] The set $\{\hat\psi_1^n,\ldots,\hat\psi_{d^n}^n\}$ forms an orthonormal basis of $\mathscr V^n$.
        \item [\em 4)] The objective of \eqref{SIAM:Eq-I.1.1.2} can be written as
        \begin{equation}
            \label{SIAM:Eq-I.1.1.23}
            \sum_{k=1}^{K}\omega_k^{K}\sum_{j=1}^n\alpha_j^n\Big\| y_j^k-\sum_{i=1}^\ell{\langle y_j^k,\psi_i\rangle}_X\,\psi_i\Big\|_X^2=\sum_{i=1}^{d^n}\hat\lambda_i^n-\sum_{k=1}^{K}\omega_k^{K}\sum_{j=1}^n\alpha_j^n\sum_{i=1}^\ell\big|{\langle y_j^k,\psi_i\rangle}_X\big|^2.
        \end{equation}
    \end{enumerate}
\end{lemma}

\begin{remark}
    \label{Remark:AltCost-1}
    \rm Due to Lemma~\ref{Lem:PropRnOp}-2) and -4) the minimization problem \eqref{SIAM:Eq-I.1.1.2} is equivalent to
    \begin{equation}
        \tag{$\mathbf{\hat P}^\ell_n$}
        \label{SIAM:Eq-I.1.1.4}
        \min -\sum_{k=1}^{K}\omega_k^{K}\sum_{j=1}^n\alpha_j^n\sum_{i=1}^\ell\big|{\langle y_j^k,\psi_i\rangle}_X\big|^2\text{ s.t. }\{\psi_i\}_{i=1}^\ell\subset X\text{ and }{\langle\psi_i,\psi_j\rangle}_X=\delta_{ij},~1 \le i,j \le \ell.
    \end{equation}
    \hfill$\blacksquare$
\end{remark}

Now we formulate the main result for \eqref{SIAM:Eq-I.1.1.2} and \eqref{SIAM:Eq-I.1.1.4}. For the proof we again refer the reader to Section~\ref{SIAM:Section-2.6.1}.

\begin{theorem}
    \label{SIAM:Theorem-I.1.1.3}
    Let $y_1^k,\ldots,y_n^k\in X$ for $1\le k\le{K}$ be given snapshots and $\mathcal R^n:X\to X$ be defined by \eqref{SIAM:Eq-I.1.1.12}. Suppose that $\{\hat\lambda_i^n\}_{i\in\mathbb I}$ and $\{\hat\psi_i^n\}_{i\in\mathbb I}$ denote the non-negative eigenvalues and associated orthonormal eigenfunctions of $\mathcal R^n$ satisfying \eqref{SIAM:Eq-I.1.1.18}. Then for every $\ell\in\{1,\ldots,d^n\}$, the first $\ell$ eigenfunctions $\{\hat\psi_i^n\}_{i=1}^\ell$ solve \eqref{SIAM:Eq-I.1.1.2} and \eqref{SIAM:Eq-I.1.1.4}. Moreover, the values of the costs evaluated at the optimal solution $\{\hat\psi_i^n\}_{i=1}^\ell$ satisfy
    \begin{equation}
        \label{SIAM:Eq-I.1.1.24}
        \sum_{k=1}^{K}\omega_k^{K}\sum_{j=1}^n\alpha_j^n\Big\|y_j^k-\sum_{i=1}^\ell{\langle y_j^k,\hat\psi_i^n\rangle}_X\,\hat\psi_i^n\Big\|_X^2=\sum_{i=\ell+1}^{d^n}\hat\lambda_i^n
    \end{equation}
    and
    \begin{equation}
        \label{SIAM:Eq-I.1.1.25}
            \sum_{k=1}^{K}\omega_k^{K}\sum_{j=1}^n\alpha_j^n\sum_{i=1}^\ell\big|{\langle y_j^k,\hat\psi_i^n\rangle}_X\big|^2=\sum_{i=1}^\ell\hat\lambda_i^n.
    \end{equation}
\end{theorem}

\begin{example}
    \label{Example:cosExample_eigenvalueDecay}
    \rm Let us return to Example~\ref{Example:cosExample}, where the trajectories of three time- and space-dependent functions $y_i$ ($i=1,2,3$) are captured using POD bases. We infer from \eqref{SIAM:Eq-I.1.1.24} and \eqref{SIAM:Eq-I.1.1.25} that the decay rate of the eigenvalues $\hat\lambda_1^n, \hat\lambda_2^n, \ldots$ can be regarded as an indicator of the quality of a POD approximation. 
    \begin{figure}
        \begin{center}
            \includegraphics[height=50mm]{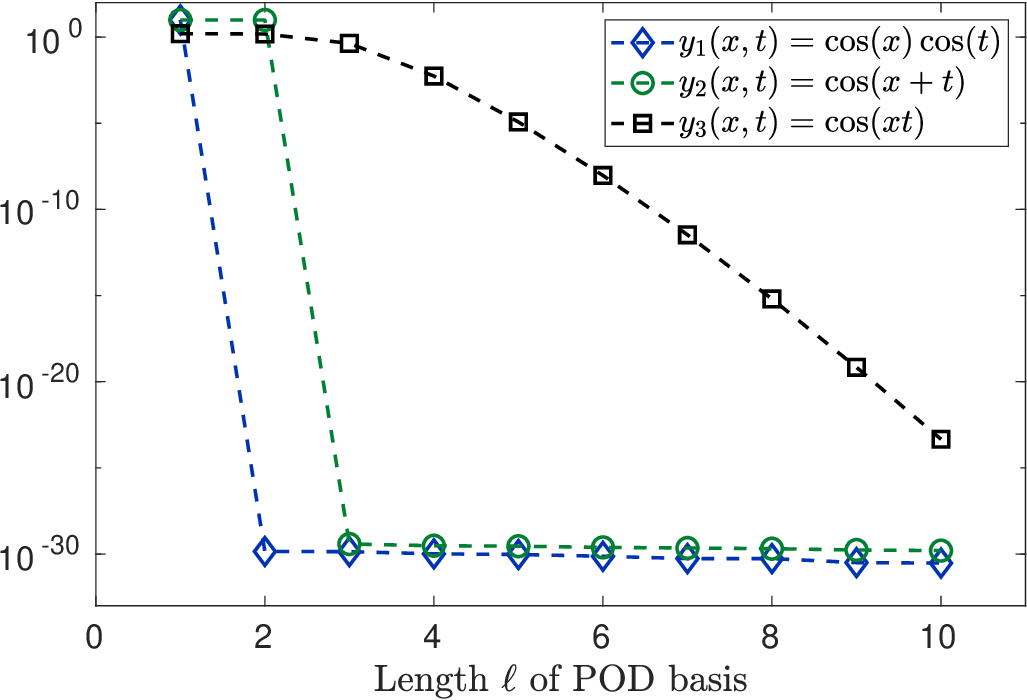}
        \end{center}
        \caption{Example~\ref{Example:cosExample_eigenvalueDecay}. Decay of the eigenvalues $\hat \lambda_1^n,\ldots,\hat \lambda_{10}^n$ for the three basic functions $y_1$, $y_2$ and $y_3$ from Example \ref{Example:cosExample}.}
        \label{fig:cosExample_eigenvalueDecay}
    \end{figure}
    In Figure~\ref{fig:cosExample_eigenvalueDecay} this decay is plotted over the index $i$. We can see that $\hat\lambda_2^n, \hat\lambda_3^n, \hdots$ are essentially zero for the first function $y_1(t,x) = \cos(t) \cos(x)$. This fits our theoretical and numerical observation from Example~\ref{Example:cosExample} that this trajectory is one-dimensional. The same is true for $y_2(t,x) = \cos(t+x)$ where $\hat\lambda_3^n, \hat\lambda_4^n, \hdots$ vanish in accordance with the two-dimensionality of the trajectory. Lastly, the eigenvalues of $y_3(t,x) = \cos(tx)$ decay much slower and more uniformly, indicating to us that there is no inherent low-dimensional subspace that will be able to fully capture this trajectory.\hfill$\blacklozenge$
\end{example}

\subsection{Relationship to singular value decomposition}
\label{SIAM:Section-RelationshipToSVD}

POD is strongly related to the {\em singular value decomposition}\index{Singular value decomposition, SVD}\index{POD method!discrete variant!SVD} (SVD) for linear operators. To show this relationship, let us endow the linear space $\mathbb C^{{K}\times n}$ of complex-valued $({K}\times n)$ matrices with the weighted inner product
\begin{align*}
    {\langle\bPhi,\bPsi\rangle}_{\mathbb C^{{K}\times n}}=\sum_{k=1}^{K}\omega_k^{K}\sum_{j=1}^n\alpha_j^n\Phi_{kj}\overline\Psi_{kj}\quad\text{for }\bPhi=\big(\big(\Phi_{kj}\big)\big),\,\bPsi=\big(\big(\Psi_{kj}\big)\big)\in\mathbb C^{{K}\times n},
\end{align*}
where the weights $\omega_1^{K},\ldots,\omega_{K}^{K}$ and $\alpha_1^n,\ldots,\alpha_n^n$ are assumed to be positive. Furthermore, we utilize the induced norm
\begin{align*}
    {\|\bPhi\|}_{\mathbb C^{{K}\times n}}={\langle\bPhi,\bPhi\rangle}_{\mathbb C^{{K}\times n}}^{1/2}=\bigg(\sum_{k=1}^{K}\omega_k^{K}\sum_{j=1}^n\alpha_j^n\,|\Phi_{kj}|^2\bigg)^{1/2}\text{ for }\bPhi=\big(\big(\Phi_{kj}\big)\big)\in\mathbb C^{{K}\times n}.
\end{align*}
Recall that for fixed $n,{K}\in\mathbb N$ we have introduced the snapshots $y_1^k,\ldots,y_n^k\in X$ for $1\le k\le{K}$ at the beginning of Section~\ref{SIAM:Section-2.1.1.1}. Now we define the linear operator $\mathcal Y^n:\mathbb C^{{K}\times n}\to X$ as
\begin{equation}
    \label{Eq:OpY}
    \mathcal Y^n\bPhi=\sum_{k=1}^{K}\omega_k^{K}\sum_{j=1}^n\alpha_j^n\Phi_{kj}y_j^k\quad\text{for }\bPhi=\big(\big(\Phi_{kj}\big)\big)\in\mathbb C^{{K}\times n}.
\end{equation}

Properties of $\mathcal Y^n$ are presented in the next lemma which is proved in Section~\ref{SIAM:Section-2.6.1}.

\begin{lemma}
    \label{Lem:OpY}
    \begin{enumerate}
        \item [\em 1)] The operator $\mathcal Y^n$ is bounded and compact.
        \item [\em 2)] Its adjoint operator $\mathcal Y^{n,\star}=(\mathcal Y^n)^\star:X\to\mathbb C^{{K}\times n}$ of $\mathcal Y^n$ is given by
        \begin{align*}
            (\mathcal Y^{n,\star}\psi)_{kj} = \langle\psi,y^k_j\rangle_X \quad \text{for } 1\le k\le{K} \text{ and } 1\le j\le n
        \end{align*}
        for all $\psi \in X$.
        \item [\em 3)] The operator $\mathcal Y^n\mathcal Y^{n,\star}:X\to X$ is linear, bounded, compact, self-adjoint and non-negative. In particular, it holds $\mathcal Y^n\mathcal Y^{n,\star} = \mathcal R^n$, where the operator $\mathcal R^n$ has been introduced in \eqref{SIAM:Eq-I.1.1.12}.
    \end{enumerate}
\end{lemma}

\begin{remark}
    \label{SIAM:Remark_YY_adj_R}
    \rm In particular, by introducing the operator $\mathcal Y^n$ we have proved Lemma~\ref{SIAM:Lemma-I.1.1.1} in a different way.\hfill$\blacksquare$
\end{remark}

Suppose that $\{\hat\lambda_i^n\}_{i\in\mathbb I}$ and $\{\hat\psi_i^n\}_{i\in\mathbb I}$ denote the non-negative eigenvalues and associated orthonormal eigenfunctions of $\mathcal R^n$ given by Theorem~\ref {SIAM:Theorem-I.1.1.3}. Then we define
\begin{equation}
    \label{Eq:SVD-1}
    \hat\sigma_i^n=\sqrt{\hat\lambda_i^n}>0\quad\text{and}\quad \hat\bPhi^n_i=\frac{1}{\hat\sigma_i^n}\,\mathcal Y^{n,\star}\hat\psi_i^n\in\mathbb C^{{K}\times n}\quad \text{for }i=1,\ldots,d^n.
\end{equation}
The SVD of the linear operator $\mathcal Y^n$ is considered in the next lemma, which is proved in Section~\ref{SIAM:Section-2.6.1}.

\begin{lemma}
    \label{Lem:SVD_Y}
    The elements $\{\hat\sigma_i^n \}_{i=1}^{d^n}$, $\{ \hat\bPhi^n_i \}_{i=1}^{d^n}\subset\mathbb C^{K\times n}$ and $\{ \hat\psi_i^n \}_{i=1}^{d^n}\subset X$ form the SVD of the operator $\mathcal Y^n$, i.e. $\{ \hat\bPhi^n_i \}_{i=1}^{d^n}$ and $\{ \hat\psi_i^n \}_{i=1}^{d^n}$ are respectively orthonormal in $\mathbb C^{{K}\times n}$ and $X$, and it holds
    \begin{equation}
        \label{Eq:SVD_total}
        \mathcal Y^n\hat\bPhi^n_i=\hat\sigma_i^n\hat\psi_i^n\in X\quad\text{and}\quad\mathcal Y^{n,\star}\hat\psi_i^n=\hat\sigma_i^n\hat\bPhi_i^n\in\mathbb C^{{K}\times n}\quad\text{for }i=1,\ldots,d^n.
    \end{equation}
    Moreover, we have the representation
    \begin{equation}
        \label{Eq:Lem:SVD_Y}
        y_j^k=\sum_{i=1}^{d^n}{\langle y_j^k,\hat\psi_i^n\rangle}_X\hat\psi_i^n=\sum_{i=1}^{d^n}\overline{\big(\mathcal Y^{n,\star}\hat\psi_i^n\big)}_{kj}\hat\psi_i^n=\sum_{i=1}^{d^n}\hat\sigma_i^n\overline{\big(\hat\Phi_i^n\big)}_{kj}\hat\psi_i^n
    \end{equation}
    for $1\le k\le{K}$ and $1\le j\le n$.
\end{lemma}

\begin{remark}
    \label{SIAM:Remark_Operator_K}
    \rm Let the liner operator $\mathcal K^n$ be given by $\mathcal K^n=\mathcal Y^{n,\star}\mathcal Y^n:\mathbb C^{{K}\times n}\to\mathbb C^{{K}\times n}$. Then, for any $\bPhi=((\Phi_{\nu l}))\in\mathbb C^{{K}\times n}$ the components of the matrix $\mathcal K^n \bPhi\in\mathbb C^{{K}\times n}$ are given as
    \begin{align*}
        \big(\mathcal K^n\bPhi\big)_{kj}=\Bigg(\mathcal Y^{n,\star}\bigg(\sum_{\nu=1}^{K}\omega_\nu^{K}\sum_{l=1}^n\alpha_l^n\Phi_{\nu l}y_l^\nu\bigg)\Bigg)_{kj}=\sum_{\nu=1}^{K}\omega_\nu^{K}\sum_{l=1}^n\alpha_l^n\Phi_{\nu l}{\langle y^\nu_l,y_j^k\rangle}_X
    \end{align*}
    for $1\le k\le{K}$ and $1\le j\le n$. It follows from Lemma~\ref{Lem:OpY} that $\mathcal K^n$ is compact, self-adjoint and non-negative. In particular, \eqref{Eq:SVD_total} implies that
    \begin{align*}
        \mathcal K^n\hat\bPhi_i^n=\mathcal Y^{n,\star}\mathcal Y^n\hat\bPhi_i^n=\hat\sigma_i^n\mathcal Y^{n,\star}\hat\psi_i^n=\hat\lambda_i^n\hat\bPhi_i^n\quad\text{for }i=1,\ldots,d^n.
    \end{align*}
    Consequently, the matrices $\{\hat\bPhi_i^n\}_{i=1}^{d^n}$ are orthonormal eigenvectors of $\mathcal K^n$ associated with the positive eigenvalues $\hat\lambda_1^n\ge\ldots\ge\hat\lambda_{d^n}^n>0$. \hfill$\blacksquare$
\end{remark} 

Next we use the SVD of the operator $\mathcal Y^n$ to develop an elementwise error estimate for the POD in the sense of \eqref{SIAM:Eq-I.1.1.24}. For the proof we refer the reader to Section~\ref{SIAM:Section-2.6.1}.

\begin{lemma}
    \label{Lem:PointwiseErrEstPOD}
    Let $(k,j)\in\{1,\ldots,{K}\}\times\{1,\ldots,n\}$ be arbitrary. Then it holds
    \begin{equation}
        \label{Eq:SVD2norm}
        \omega_k^{K}\alpha_j^n\,\Big\|y_j^k-\sum_{i=1}^\ell{\langle y_j^k,\hat\psi_i^n\rangle}_X\,\hat\psi_i^n\Big\|_X^2\le\hat\lambda_{\ell+1}^n.
    \end{equation}
\end{lemma} 

\subsection{Computational aspects for POD}

For the application of POD to concrete problems the choice of $\ell$ is certainly of central importance. It appears that no general a-priori rules are available. Rather the choice of $\ell$ is based on heuristic considerations combined with observing the ratio of the \index{POD method!discrete variant!relative total energy}{\em relative total energy} contained in the snapshots $y_1^k,\ldots,y_n^k$, $1\le k\le{K}$, which is expressed by
\begin{align*}
    \mathcal E(\ell) = \frac{\sum_{i=1}^\ell\hat\lambda_i^n}{\sum_{i=1}^{d^n}\hat\lambda_i^n}\in[0,1].
\end{align*}
Utilizing \eqref{SIAM:Eq-I.1.1.18} we have
\begin{equation}
    \label{SIAM:Eq-I.1.1.39b}
    \mathcal E(\ell)=\frac{\sum_{i=1}^\ell\hat\lambda_i^n}{\sum_{k=1}^{K}\omega_k^{K}\sum_{j=1}^n\alpha_j^n\,{\| y_j^k\|}_X^2},
\end{equation}
i.e., the computation of the eigenvalues $\{\hat\lambda_i\}_{i=\ell+1}^{d^n}$ is not necessary. This is utilized in numerical implementations when iterative eigenvalue solvers like the Lanczos method; see \cite[Chapter~10]{Ant05}, are applied. In these kind of solvers not all eigenvalues are computed, but only the first $\ell$ largest one. Let us summarize the computation of a POD basis of rank $\ell$ in Algorithm~\ref{SIAM:Algorithm-I.1.1.1}.

\begin{remark}
    \rm We will address the solution of the eigenvalue problems in lines 2 and 8 in the Sections~\ref{SIAM:Section-2.1.1.2}-\ref{SIAM:Section-2.1.1.4}.\hfill$\blacksquare$
\end{remark}

\bigskip
\hrule
\vspace{-3.5mm}
\begin{algorithm}[(Discrete POD in Hilbert spaces)]
    \label{SIAM:Algorithm-I.1.1.1}
    \vspace{-3mm}
    \hrule
    \vspace{0.5mm}
	\begin{algorithmic}[1]
        \REQUIRE Snapshots $\{y_j^k\}_{j=1}^n\subset X$ for $1\le k\le{K}$ and rank $\ell \le d^n$ or tolerance $0<\varepsilon_{\mathsf{POD}}\ll1$;
        \STATE Compute $\Lambda^n=\sum_{k=1}^{K}\omega_k^{K}\sum_{j=1}^n\alpha_j^n\,\| y_j^k\|_X^2$;
        \IF{$\ell$ is given}
            \STATE Compute pairs $\{(\hat\lambda_i^n,\hat\psi_i^n)\}_{i=1}^\ell$ by solving \eqref{SIAM:Eq-I.1.1.18} with $\hat\lambda_1^n\ge\ldots\ge\hat\lambda^n_\ell$ and $\langle\hat\psi_i^n,\hat\psi_j^n\rangle_X=\delta_{ij}$ for $i,j=1,\ldots,\ell$;
            \STATE Set $\mathcal E(\ell)=\sum_{i=1}^\ell\hat\lambda_i^n/\Lambda^n$;
        \ELSIF{$\varepsilon_\mathsf{POD}$ is given}
            \STATE Set $i=0$ and $\mathcal E(0)=0$;
            \WHILE{$\mathcal E(i)\le1-\varepsilon_{\mathsf{POD}}$ {\bf and} $i<d^n$}
                \STATE Set $i=i+1$;
                \STATE Compute a pair $(\hat\lambda_i^n,\hat\psi_i^n)$ by solving $\mathcal R^n\hat\psi_i^n=\hat\lambda_i^n\hat\psi_i^n$ subject to  $\hat\lambda_i^n\ge\hat\lambda_j^n$ for $j=i+1,\ldots,d^n$ and $\langle\hat\psi_i^n,\hat\psi_j^n\rangle_X=\delta_{ij}$ for $j=1,\ldots,i-1$;
                \STATE Define $\mathcal E(i)=\mathcal E(i-1)+\hat\lambda_i^n/\Lambda^n$;
            \ENDWHILE
            \STATE Set $i=\ell$;
        \ENDIF
        \RETURN POD basis $\{\hat\psi_i^n\}_{i=1}^\ell$, eigenvalues $\{\hat\lambda_i^n\}_{i=1}^\ell$ and ratio $\mathcal E(\ell)$;
    \end{algorithmic}
    \hrule
\end{algorithm}

\subsection{POD-greedy method for multiple snapshots}
\label{SIAM:Section-PODGreedy}

Let us mention the following approach for multiple snapshots which is related to the {\em POD-greedy algorithm}\index{POD method!discrete variant!greedy} for parametrized evolution problems; cf. \cite{Gre12,Haa12,HO08}, for instance. Suppose again that we are given ${K}$ different snapshot trajectories $\{y_j^k\}_{j=1}^n\subset X$ for $k=1,\ldots,{K}$. Starting with $\nu=1$, we choose a $k_\nu=k_1\in\{1,\ldots,K\}$ and consider the \index{POD method!discrete variant!snapshot space}{\em snapshot space} defined by the $k_\nu$-th trajectory 
\begin{equation}
    \label{Greedy-1}
    \mathscr V^n_\nu=\mathrm{span}\,\big\{y_{\nu j}\,\big|\,1\le j\le n\big\}\subset X,\quad y_{\nu j}=y_j^{k_\nu}\text{ for }j=1,\ldots,n
\end{equation}
with the finite dimension $d^n_\nu=\dim\mathscr V^n_\nu\le n$. We define the linear operator
\begin{equation}
    \label{Greedy-2}
    \mathcal R^n_\nu\psi=\sum_{j=1}^n\alpha_j^n\,{\langle y_{\nu j},\psi_i\rangle}_X\,y_{\nu j}\quad\text{for }\psi\in X,
\end{equation}
which is compact, self-adjoint and non-negative by Lemma~\ref{SIAM:Lemma-I.1.1.1}. Suppose that $\{\hat\lambda_{\nu i}^n\}_{i\in\mathbb I}$ and $\{\hat\psi_{\nu i}^n\}_{i\in\mathbb I}$ are the non-negative eigenvalues and associated orthonormal eigenfunctions satisfying
\begin{equation}
    \label{Greedy-3}
    \mathcal R^n_\nu\hat\psi_{\nu i}^n=\hat\lambda_{\nu i}^n\hat\psi_{\nu i}^n,\quad\hat\lambda_{\nu 1}^n\ge\ldots\ge\hat\lambda_{\nu d_\nu^n}^n>\hat\lambda_{\nu d_\nu^n+1}^n=\ldots=0.
\end{equation}
For a given $\ell_\nu\le d^n_\nu$, the eigenfunctions $\{\hat\psi_{\nu i}^n\}_{i=1}^{\ell_\nu}$ solve
\begin{align*}
    \sum_{j=1}^n\alpha_j^n \Big\| y_{\nu j}-\sum_{i=1}^{\ell_\nu}{\langle y_{\nu j},\psi_i\rangle}_X\,\psi_i\Big\|_X^2\text{ s.t. } \{\psi_i\}_{i=1}^{\ell_\nu}\subset X\text{ and }{\langle\psi_i,\psi_j\rangle}_X=\delta_{ij},~1 \le i,j \le \ell_\nu,
\end{align*}
i.e., $\{\hat\psi^n_{\nu i}\}_{i=1}^{\ell_\nu}$ is a POD basis of rank $\ell_\nu$ for the snapshot space $\mathscr V^n_\nu$. In the next step we test if this POD basis leads also to a good POD basis for the other snapshot ensembles $\{y_j^k\}_{j=1}^n\subset X$ for $k\neq k_1$. If this is not the case, we add a POD basis of the snapshot ensemble $\{y_j^{k_2}\}_{j=1}^n\subset X$, $k_2\in\{1,\ldots,k\}\setminus\{k_1\}$, who leads to the worst approximation quality to the already computed POD bases functions. The method in summarized in Algorithm~\ref{SIAMBook:Alg_PODgreedy}.

\bigskip
\hrule
\vspace{-3.5mm}
\begin{algorithm}[(POD greedy method)]
    \label{SIAMBook:Alg_PODgreedy}
    \vspace{-3mm}
    \hrule
    \vspace{0.5mm}
    \begin{algorithmic}[1]
        \REQUIRE Snapshots $\{y_j^k\}_{j=1}^n\subset X$ for $1\le k\le{K}$ and rank $\ell \le d^n$ or tolerance $0<\varepsilon_{\mathsf{POD}}\ll1$;
        \STATE Set $\nu=1$, $\ell=0$ and $\mathscr N=\emptyset$;
        \REPEAT
            \STATE Compute the projection errors
            \begin{align*}
                \mathsf{err}_\nu^\ell(k)=\sum_{j=1}^n\alpha_j^n \Big\| y_j^k-\sum_{i=1}^\ell{\langle y_j^k,\hat\psi_i^n\rangle}_X\,\hat\psi_i^n\Big\|_X^2\quad\text{for }k\in\big\{1,\ldots,{K}\big\}\setminus\mathscr N;
            \end{align*}
            \STATE Determine the index $k_\nu\in\argmax\{\mathsf{err}_\nu^\ell(k)\,|\,k\in\{1,\ldots,{K}\}\setminus\mathscr N\}$;
            \STATE Set $\mathscr N=\mathscr N\cup\{k_\nu\}$;
            \STATE Choose the snapshots $y_{\nu j}=y_j^{k_\nu}-\sum_{i=1}^\ell{\langle y_j^{k_\nu},\hat\psi_i^n\rangle}_X\,\hat\psi_i^n$ for $j=1,\ldots,n$;
            \STATE Solve \eqref{Greedy-3} to get a POD basis $\{\hat\psi^n_{\nu i}\}_{i=1}^{\ell_\nu}$ satisfying $\mathsf{err}_\nu^\ell(k_\nu)\le\varepsilon_\mathsf{POD}$;
            \STATE Set $\hat\psi_{\ell+i}^n=\hat\psi^n_{\nu i}$ for $i=1,\ldots,\ell_\nu$ and put $\ell=\ell+\ell_\nu$;
        \UNTIL{$\max\big\{\mathsf{err}_\nu^\ell(k)\,|\,k=1,\ldots,{K}\big\}\le\varepsilon_\mathsf{POD}$ \textbf{or} $\nu={K}$}
        \RETURN POD basis $\{\hat\psi_i^n\}_{i=1}^\ell$ and eigenvalues $\{\hat\lambda_i^n\}_{i=1}^\ell$;
    \end{algorithmic}
    \hrule
\end{algorithm}

\begin{remark}
    \rm
    \begin{enumerate}
        \item [1)] In steps 3 and 6 we set
        \begin{align*}
            \sum_{i=1}^\ell{\langle y_j^{k_\nu},\hat\psi_i^n\rangle}_X\,\hat\psi_i^n=0,\quad k_\nu\in\big\{1,\ldots,{K}\big\}
        \end{align*}
        for $\ell=0$. In particular, we have 
        \begin{align*}
            k_1\in\argmax\bigg\{\sum_{j=1}^n\alpha_j^n \big\| y_j^k\big\|_X^2\,\Big|\,k\in\{1,\ldots,{K}\}\bigg\}
        \end{align*}
        for $\ell=0$.
        \item [2)] By construction we get an orthonormal set $\{\hat\psi_i\}_{i=1}^\ell$. However, the \index{Method!Gram-Schmidt}{\em Gram-Schmidt orthogonalization}\index{Method!Gram-Schmidt} \cite[pp.~56-61]{TB97} should be employed to ensure this property numerically.\hfill$\blacksquare$
    \end{enumerate}
\end{remark}

\subsection{POD in the unitary spaces}
\label{SIAM:Section-2.1.1.2}

Let us discuss the POD method for the specific case $X=\mathbb C^m$ with $m\in\mathbb N$ and a single snapshot ensemble, meaning that we set ${K}=1$ and $\omega_1^{K}=1$. Then we have $n$ snapshot vectors $y_1,\ldots,y_n$ and introduce the rectangular matrix $\bY=[y_1\,|\ldots|\,y_n]\in\mathbb C^{m\times n}$ with rank $d^n\le\min(m,n)$. Choosing $\alpha_j^n=1$ for $1\le j\le n$, Problem \eqref{SIAM:Eq-I.1.1.2} has the form
\begin{equation}
    \label{SIAM:Eq-I.1.1.41}
    \min\sum_{j=1}^n\Big\| y_j-\sum_{i=1}^\ell{\langle y_j,\psi_i\rangle}_{\mathbb C^m}\,\psi_i\Big\|_{\mathbb C^m}^2\text{ s.t. } \{\psi_i\}_{i=1}^\ell\subset \mathbb C^m\text{ and }{\langle\psi_i,\psi_j\rangle}_{\mathbb C^m}=\delta_{ij},~1 \le i,j \le \ell.
\end{equation}
In \eqref{Eq:OpY} we have introduced the operator $\mathcal Y^n$, which is now of the simplified form
\begin{align*}
    \mathcal Y^n\phi=\sum_{j=1}^n\phi_jy_j=\bY\phi\quad\text{for }\phi\in\mathbb C^n,
\end{align*}
i.e., $\mathcal Y^n=\bY$ holds true. This implies $\mathcal Y^{n,\star}=\bY^\mathsf H$, where ``$\mathsf H$'' stands for the Hermitian on $\bY$. Furthermore, $\mathcal R^n=\bY\bY^\mathsf H$ holds true. We infer that \eqref{SIAM:Eq-I.1.1.18} leads to the symmetric $m\times m$ eigenvalue problem
\begin{equation}
    \label{SIAM:Eq-I.1.1.37}
    \bY\bY^\mathsf H\hat\psi_i^n=\hat\lambda_i^n\hat\psi_i^n,\quad\hat\lambda_1^n\ge\ldots\ge\hat\lambda_{d^n}^n>\hat\lambda_{d^n+1}^n=\ldots =\hat\lambda_m^n=0.
\end{equation}
We have discussed the \index{Singular value decomposition, SVD}\index{POD method!discrete variant!SVD}{\em singular value decomposition} (SVD) for linear operators in Section~\ref{SIAM:Section-RelationshipToSVD}. Since $\mathcal Y^n$ is just the matrix $\bY$, the SVD for linear operators coincides with the well-known SVD for matrices \cite[pp.~25-30]{TB97}: There exist real numbers $\hat\sigma_1^n\ge\ldots\ge\hat\sigma_{d^n}^n>0$ and unitary matrices $\bPsi \in \mathbb C^{m \times m}$ with column vectors $\{\hat\psi_i^n\}_{i=1}^m\subset\mathbb R^m$ and $\bPhi\in\mathbb C^{n \times n}$ with column vectors $\{\hat\phi_i^n\}_{i=1}^n\subset\mathbb R^n$ such that
\begin{equation}
    \label{SIAM:Eq-I.1.1.38}
    \bPsi^\mathsf H \bY \bPhi = \left(
    \begin{array}{cc}
        \bD~ & 0\\
        0~ & 0
    \end{array}
    \right)=: \bSigma \in \mathbb R^{m \times n},
\end{equation}
where $\bD=\mathrm{diag}\,(\hat\sigma_1^n,\ldots,\hat\sigma_{d^n}^n) \in \mathbb R^{d^n \times d^n}$ and the zeros in \eqref{SIAM:Eq-I.1.1.38} denote matrices of appropriate dimensions. Moreover the vectors $\{\hat\psi_i^n\}_{i=1}^{d^n}$ and $\{\hat\phi_i^n\}_{i=1}^{d^n}$ satisfy
\begin{equation}
    \label{SIAM:Eq-I.1.1.39}
    \bY\hat\phi_i^n=\hat\sigma_i^n\hat\psi_i^n\quad\text{and}\quad\bY^\mathsf H\hat\psi_i^n=\hat\sigma_i^n\hat\phi_i^n \quad \text{for } i=1,\ldots,d^n.
\end{equation}
They are eigenvectors of $\bY\bY^\mathsf H$ and $\bY^\mathsf H\bY$, respectively, with eigenvalues $\hat\lambda_i^n=(\hat\sigma_i^n)^2>0$, $i=1,\ldots,d^n$. The vectors $\{\hat\psi_i^n\}_{i=d^n+1}^m$ and $\{\hat\phi_i^n\}_{i=d^n+1}^n$ (if $d^n<m$ respectively $d^n<n$) are eigenvectors of $\bY\bY^\mathsf H$ and $\bY^\mathsf H\bY$ with eigenvalue $0$. Consequently, in the case $n<m$ one can determine the POD basis of rank $\ell$ as follows: Compute the eigenvectors $\hat\phi_1^n,\ldots,\hat\phi_\ell^n \in \mathbb C^n$ by solving the symmetric $n \times n$ eigenvalue problem
\begin{subequations}
    \label{MethodsOfSnapshots-1}
    \begin{equation}
        \label{MethodsOfSnapshots-1a}
        \bY^\mathsf H\bY\hat\phi_i^n=\hat\lambda_i^n\hat\phi_i^n \quad \text{for } i=1,\ldots,\ell
    \end{equation}
    and set, by \eqref{SIAM:Eq-I.1.1.39},
    \begin{equation}
        \label{MethodsOfSnapshots-1b}
        \hat\psi_i^n=\frac{1}{\hat\sigma_i^n}\,\bY\hat\phi_i^n \quad \text{for } i=1,\ldots,\ell.
    \end{equation}
\end{subequations}
Recall that by Theorem~\ref{SIAM:Theorem-I.1.1.3} we have
\begin{equation}
    \label{SIAM:Eq-I.1.1.42}
    \sum_{j=1}^n\Big\| y_j-\sum_{i=1}^\ell{\langle y_j,\hat\psi_i^n\rangle}_{\mathbb C^m}\,\hat\psi_i\Big\|_{\mathbb C^m}^2=\sum_{i=\ell+1}^{d^n}\hat\lambda_i=\sum_{i=\ell+1}^{d^n}(\hat\sigma_i^n)^2.
\end{equation}
From \eqref{Eq:SVD2norm} we deduce that
\begin{align*}
    \Big\|y_j-\sum_{i=1}^\ell{\langle y_j,\hat\psi_i^n\rangle}_{\mathbb C^m}\,\hat\psi_i^n\Big\|_{\mathbb C^m}\le\hat\sigma_{\ell+1}^n\quad\text{for every }j\in\{1,\ldots,n\}.
\end{align*} 
Thus, the $(\ell+1)$-th singular value gives an upper bound for the approximation error for each snapshot vector $y_j$, $j=1,\ldots,n$.

\begin{remark}
    \label{Remark:Sirovich}
    \rm
    \begin{enumerate}
    \item [1)] For historical reasons the method of determining the POD basis by solving \eqref{MethodsOfSnapshots-1} is sometimes called the \index{POD method!discrete variant!methods of snapshots}{\em method of snapshots}; see \cite{Sir87}. For $i\in\{1,\ldots,\ell\}$ let us denote by $(\hat\phi_i)_j$, $1\le j\le n$ the $j$-th component of the eigenvector $\hat\phi_i\in\mathbb C^n$. Using $\bY=[y_1\,|\ldots|\,y_n]$ and \eqref{MethodsOfSnapshots-1b} we find that the $i$-th POD basis is given by the snapshots $\{y_j\}_{j=1}^n$ as follows:
    \begin{align*}
        \hat\psi_i^n=\frac{1}{\hat\sigma_i^n}\,\bY\hat\phi_i^n=\sum_{j=1}^n\bigg(\frac{1}{\hat\sigma_i^n}(\hat\phi_i)_j\bigg)y_j \quad \text{for } i=1,\ldots,\ell.
    \end{align*}
    \item [2)] If $m<n$ holds, we compute the POD basis by solving the $m \times m$ eigenvalue problem \eqref{SIAM:Eq-I.1.1.37}. If the matrix $\bY$ is badly scaled, we should avoid to build the matrix product $\bY\bY^\mathsf H$ (or $\bY^\mathsf H\bY$). In this case the SVD turns out to be more stable for the numerical computation of the POD basis of rank $\ell$.\hfill$\blacksquare$
    \end{enumerate}
\end{remark}

Note that
\begin{align*}
    \sum_{j=1}^n{\|y_j\|}_{\mathbb C^m}^2={\|\bY\|}_F^2 \quad \text{for } j=1,\ldots,n
\end{align*}
holds, where the Frobenius norm $\|\cdot\|_F$ has been recalled in Section~\ref{App:Basics}. Thus, we infer from \eqref{SIAM:Eq-I.1.1.39b}
\begin{align*}
    \mathcal E(\ell)=\frac{\sum_{i=1}^\ell\hat\lambda_i^n}{\sum_{j=1}^n{\| y_j\|}_{\mathbb C^m}^2} = \frac{\sum_{i=1}^\ell\hat\lambda_i^n}{{\|\bY\|}_F^2}.
\end{align*}
Let $\{\hat\psi_i^n\}_{i=1}^m\subset\mathbb C^m$ be a solution to \eqref{SIAM:Eq-I.1.1.37} with associated eigenvalues $\{\hat\lambda_i^n\}_{i=1}^m$. Since $\{\hat\psi_1^n,...,\hat\psi_{d^n}^n\}$ form an orthonormal basis of $\mathscr V^n$, we have
\begin{equation}
    \label{SIAM:Eq-I.1.1.43}
    y_j=\sum_{i=1}^{d^n}{\langle y_j,\hat\psi_i^n\rangle}_{\mathbb C^m}\,\hat\psi_i^n\quad\text{for any }j\in\{1,\ldots,n\}.
\end{equation}
We define the matrices $\hat\bPsi^{d^n}=[\hat\psi_1^n\,|\ldots,|\hat\psi_{d^n}^n]\in\mathbb C^{m\times d^n}$ and $\hat\bB^{d^n}$ by $\hat B^{d^n}_{ij}=\langle y_j,\hat\psi_i^n\rangle_{\mathbb C^m}$ for $1 \le i \le d^n$, $1\le j \le n$. Then we can write \eqref{SIAM:Eq-I.1.1.43} as $\bY=\hat\bPsi^{d^n}\hat\bB^{d^n}$. From \eqref{SIAM:Eq-I.1.1.42} it follows that
\begin{align*}
    {\|\bY-\hat\bPsi^\ell\hat\bB^\ell\|}_F^2&=\sum_{i=1}^m \sum_{j=1}^n \Big|Y_{ij}-\sum_{k=1}^\ell\hat \Psi_{ik}^\ell\hat B_{kj}\Big|^2=\sum_{j=1}^n\sum_{i=1}^m\Big|Y_{ij}-\sum_{k=1}^\ell {\langle y_j,\hat\psi_i^n\rangle}_{\mathbb C^m}\hat \Psi_{ik}^\ell \Big|^2\\
    &= \sum_{j=1}^n\Big\|y_j-\sum_{k=1}^\ell{\langle y_j,\hat\psi_i^n\rangle}_{\mathbb C^m}\hat\psi_k\Big\|_{\mathbb C^m}^2=\sum_{i=\ell+1}^{d^n}(\hat\sigma_i^n)^2,
\end{align*}
where $\hat\bPsi^\ell$ contains the first $\ell$ columns of $\hat\bPsi^{d^n}$ and $\hat\bB^\ell$ the first $\ell$ rows of $\hat\bB^{d^n}$.

The next results follows from \eqref{SIAM:Eq-I.1.1.41} and \eqref{SIAM:Eq-I.1.1.37}. It states that for every $\ell \le d^n\le m$ the approximation of the columns of $\bY$ by the first $\ell$ singular vectors $\{\hat\psi_i\}_{i=1}^\ell$ is optimal with respect to the Frobenius norm among all rank $\ell$ approximations to the columns of $\bY$. For a proof we refer to Section~\ref{SIAM:Section-2.6.1}.

\begin{corollary}
    \label{SIAM:Corollary-I.1.1.1}
    Let all hypotheses of Theorem~{\rm \ref{SIAM:Theorem-I.1.1.3}} hold.
    \begin{enumerate}
        \item [\em 1)] {\em Optimality of the POD method}\index{POD method!optimality}: Suppose that $\bPsi^{d^n}=[\psi_1\,|\ldots|\,\psi_{d^n}] \in \mathbb C^{m \times d^n}$ denotes a matrix with $d^n$ pairwise orthonormal vectors $\psi_i\in\mathbb C^m$ and that the expansion of the columns of $\bY$ in the basis $\{\psi_i\}_{i=1}^{d^n}$ is given by
        \begin{align*}
            \bY=\bPsi^{d^n}\bB^{d^n}, \quad \text{where} \quad B^{d^n}_{ij}={\langle y_j,\psi_i\rangle}_{\mathbb C^m} \quad\text{for } 1 \le i \le d^n, 1\le j \le n.
        \end{align*}
        Then for every $\ell \in \{1,\ldots,d^n\}$ we have
        \begin{equation}
            \label{SIAM:Eq-I.1.1.47}
            {\|\bY-\hat\bPsi^\ell\hat\bB^\ell\|}_F \le {\|\bY-\bPsi^\ell\bB^\ell\|}_F.
        \end{equation}
        In \eqref{SIAM:Eq-I.1.1.47} the matrix $\bPsi^\ell$ contains the first $\ell$ columns of $\bPsi^{d^n}$ and $\bB^\ell$ are the first $\ell$ rows of $\bB^{d^n}$. Similarly, $\hat \bPsi^\ell$ and $\hat\bB^\ell$ are defined.
        \item [\em 2)] {\em Uncorrelated POD coefficients}: We have
        \begin{align*}
            \sum_{j=1}^n\overline{\langle y_j,\hat\psi_i^n\rangle}_{\mathbb C^m}{\langle y_j,\hat\psi_k^n\rangle}_{\mathbb C^m}=\hat\lambda_i^n\delta_{ik} \quad \text{for } 1 \le i,k \le \ell.
        \end{align*}
    \end{enumerate}
\end{corollary}

In Algorithm~\ref{SIAM:Algorithm-I.1.1.2} we summarize the steps for the computation of a POD basis of rank $\ell$ for the specific case $X=\mathbb C^m$.

\bigskip
\hrule
\vspace{-3.5mm}
\begin{algorithm}[(POD method in $\boldsymbol{\mathbb C^m}$)]
    \label{SIAM:Algorithm-I.1.1.2}
    \vspace{-3mm}
    \hrule
    \vspace{0.5mm}
    \begin{algorithmic}[1]
        \REQUIRE Snapshots $\{y_j\}_{j=1}^n \subset \mathbb C^m$, POD rank $\ell \le d^n$ and {\tt flag} for the solver;
        \STATE Set $\bY=[y_1,\ldots,y_n] \in \mathbb C^{m \times n}$;
        \IF{{\tt flag} = 0}
            \STATE Compute singular value decomposition $[\bPsi,\bSigma,\bPhi]=\mathrm{svd}\,(\bY)$;
            \STATE Set $\hat\psi_i^n=\Psi_{\cdot,i}\in \mathbb C^m$ and $\hat\lambda_i^n=\bSigma_{ii}^2$ for $i=1,\ldots,\ell$;
        \ELSIF{{\tt flag} = 1}
            \STATE Determine $\mathcal R^n=\bY\bY^\mathsf H \in \mathbb C^{m \times m}$;
            \STATE Compute eigenvalue decomposition $[\bPsi,\bLambda]=\mathrm{eig}\,(\mathcal R^n)$;
            \STATE Set $\hat\psi_i^n=\bPsi_{\cdot,i} \in \mathbb C^m$ and $\hat\lambda_i^n=\Lambda_{ii}$ for $i=1,\ldots,\ell$;
        \ELSIF{{\tt flag} = 2}
            \STATE Determine $\mathcal K^n=\bY^\mathsf H\bY \in \mathbb C^{n \times n}$;
            \STATE Compute eigenvalue decomposition $[\bPhi,\bLambda]=\mathrm{eig}\,(\mathcal K^n)$;
            \STATE Set $\hat\lambda_i^n=\Lambda_{ii}$ and $\hat\psi_i^n=\bY\bPhi_{\cdot,i}/\Lambda_{ii}^{1/2}\in \mathbb C^m$ for $i=1,\ldots,\ell$;
        \ENDIF
        \STATE Compute $\mathcal E(\ell)=\sum_{i=1}^\ell\hat\lambda_i^n/{\|\bY\|}_F^2$;
        \RETURN POD basis $\{\hat\psi_i^n\}_{i=1}^\ell$, eigenvalues $\{\hat\lambda_i^n\}_{i=1}^\ell$ and ratio $\mathcal E(\ell)$;
    \end{algorithmic}
    \hrule
\end{algorithm}

\subsection{POD in Euclidean spaces}
\label{SIAM:Section-2.1.1.2a}

Let us briefly discuss the special case $X=\mathbb R^m$, where we keep all other choices of Section~\ref{SIAM:Section-2.1.1.2}. We proceed in the same way, with the exception to consider the transpose ``$\top$'' instead of the hermitian conjugate ``$\mathsf{H}$''. More precisely, \eqref{SIAM:Eq-I.1.1.37} reads
\begin{equation}
    \label{SIAM:Eq-I.1.1.54}
    \bY\bY^\top\hat\psi_i^n=\hat\lambda_i^n\hat\psi_i^n,\quad\hat\lambda_1^n\ge\ldots\ge\hat\lambda_{d^n}^n>\hat\lambda_{d^n+1}^n=\ldots =\hat\lambda_m^n=0
\end{equation}
for the snapshot matrix $\bY=((Y_{ij}))\in\mathbb R^{m\times n}$ with $\bY=[y_1\,|\ldots|\,y_m]$. Let us explain for the special case $X=\mathbb R^m$ how we can derive \eqref{SIAM:Eq-I.1.1.54} by techniques of non-linear programming. To simplify the presentation here, we consider only the case $\ell=1$ and suppose that $\hat\lambda_1^n>\hat\lambda_2^n$ holds. In that case, \eqref{SIAM:Eq-I.1.1.4} is equivalent to the maximization problem
\begin{equation}
    \label{SIAM:Eq-I.1.1.55}
    \tag{$\mathbf{\hat P}^1_n$}
    \max\sum_{j=1}^n{\langle y_j,\psi\rangle}_2^2\quad \text{s.t.}\quad\psi=\big(\uppsi_i\big)_{1\le i\le m}\in\mathbb R^m\text{ and }{|\psi|}_2=1.
\end{equation}
Problem \eqref{SIAM:Eq-I.1.1.55} is a constrained optimization problem. Therefore, we define the \index{Lagrange functional}{\em Lagrange functional} $\mathcal L:\mathbb R^m \times \mathbb R \to \mathbb R$ associated with \eqref{SIAM:Eq-I.1.1.55} by
i.e.,
\begin{align*}
    \mathcal L(\psi,\mu)=-\sum_{j=1}^n {\langle y_j,\psi\rangle}_2^2+\mu e(\psi)\quad \text{for } (\psi,\mu) \in \mathbb R^m \times \mathbb R,
\end{align*}
where $e:\mathbb R^m\to\mathbb R$ is given by $e(\psi)=| \psi |^2_2-1$ for $\psi\in\mathbb R^m$. Suppose that $\hat\psi^n=(\hat\uppsi_i^n)_{1\le i\le m}\in \mathbb R^m$ is a solution to \eqref{SIAM:Eq-I.1.1.55}. Then $\nabla e(\hat\psi^n)=2(\hat\psi^n)^\top\neq0$ holds. Thus, $\hat\psi^n$ is.a \emph{regular point} for the constrained function $e(\cdot)$ and \index{Optimality conditions!first-order!necessary}{\em first-order necessary optimality conditions}\index{Optimality conditions!first-order necessary} \cite[Chapter~12]{NW06} for \eqref{SIAM:Eq-I.1.1.55} are given as follows: there exists an associated (unique) Lagrange multiplier $\hat\mu^n\in\mathbb R$ satisfying
\begin{align*}
    \nabla \mathcal L(\hat\psi^n,\hat\mu^n)\stackrel{!}{=}0\quad\text{in }\mathbb R^m\times \mathbb R.
\end{align*}
We compute the gradient of $\mathcal L$ with respect to $\psi$. For any component $i\in\{1,\ldots,m\}$ we get
\begin{align*}
    &\frac{\partial \mathcal L}{\partial \psi_i}\,(\hat\psi^n,\hat\mu^n)=\frac{\partial}{\partial\psi_i}\,\Bigg(-\sum_{j=1}^n \bigg|\sum_{k=1}^mY_{kj}\hat\uppsi_k^n\bigg|^2+\bigg(\sum_{k=1}^m(\hat\uppsi_k^n)^2-1\bigg)\hat\mu^n\Bigg)=2\sum_{j=1}^n\bigg(-\sum_{k=1}^mY_{kj}\hat\uppsi_k^n\bigg)Y_{ij}+2\hat\mu^n\hat\uppsi_i^n\\
    &\quad=2\sum_{k=1}^m\bigg(-\sum_{j=1}^nY_{ij}(\bY)^\top_{jk}\hat\uppsi_k^n\bigg) +2\hat\mu^n\hat\uppsi_i^n=2 \sum_{k=1}^m\big(-\bY\bY^\top)_{ik}\hat\uppsi_k^n\big) +2\hat\mu^n\hat\uppsi_i^n.
\end{align*}
Thus,
\begin{equation}
    \label{SIAM:Eq-I.1.1.59}
    \nabla_\psi \mathcal L(\hat\psi^n,\hat\mu^n)=-2\big(\bY\bY^\top\hat\psi^n-\hat\mu^n\hat\psi^n\big)\stackrel{!}{=}0\quad\text{in }\mathbb R^m.
\end{equation}
For $\hat\mu^n=\hat\lambda^n$ equation \eqref{SIAM:Eq-I.1.1.59} yields the eigenvalue problem \eqref{SIAM:Eq-I.1.1.54}. From
\begin{align*}
    \frac{\partial \mathcal L}{\partial \lambda}(\hat\psi^n,\hat\lambda^n)\stackrel{!}{=}0\quad\text{in }\mathbb R
\end{align*}
we infer the constraint $|\hat\psi^n|_2=1$. However, it is not clear yet which eigenvalue-eigenvector pair $(\hat{\lambda}_i^n,\hat{\psi}_i^n)$ is a solution to \eqref{SIAM:Eq-I.1.1.55}. Therefore, let us turn to \index{Optimality conditions!second-order}{\em second-order optimality conditions} \cite[Chapter~12]{NW06}. Note that $\nabla_{\psi\psi} \mathcal L(\psi,\mu)=-2(\bY\bY^\top-\mu\bI)\in\mathbb R^{m \times m}$ holds for any $(\psi,\mu)\in\mathbb R^m\times\mathbb R$. Let $\psi \in \mathbb R^m$ be chosen arbitrarily. Let $\{\hat\psi_i^n\}_{i=1}^m\subset\mathbb R^m$ denote $m$ orthonormal eigenvectors of $\bY\bY^\top$ satisfying \eqref{SIAM:Eq-I.1.1.54}. Then we can write $\psi$ in the form
\begin{align*}
    \psi=\sum_{i=1}^m{\langle\psi,\hat\psi_i^n\rangle}_2\,\hat\psi_i^n.
\end{align*}
At $(\hat\psi^n,\hat\lambda^n)=(\hat\psi^n_1,\hat\lambda^n_1)$ we conclude from $\hat\lambda^n_1\ge\hat\lambda_2^n\ge\ldots\ge\hat\lambda_m^n\ge 0$ that
\begin{align}
    \label{SIAM:Eq-I.1.1.43-A}
    & {\langle\psi,\nabla_{\psi\psi}\mathcal L(\hat\psi^n,\hat\lambda^n)\psi\rangle}_2=-2\,{\langle\psi,\big(\bY\bY^\top-\hat\lambda^n\bI\big)\psi\rangle}_2 \notag\\
    &=-2\sum_{i=1}^m \sum_{j=1}^m {\langle\psi,\hat\psi_i^n\rangle}_2{\langle \psi,\hat\psi_j^n \rangle}_2{\langle\hat\psi_i^n,\big(\bY\bY^\top-\hat\lambda^n\bI\big)\hat\psi_j^n\rangle}_2\notag\\
    &=-2\sum_{i=1}^m \sum_{j=1}^m \big(\hat\lambda_j^n-\hat\lambda^n\big)\,{\langle\psi,\hat\psi_i^n\rangle}_2{\langle \psi,\hat\psi_j^n\rangle}_2{\langle\hat\psi_i^n,\hat\psi_j^n\rangle}_2 \notag\\
    & = 2\sum_{i=1}^m \big(\hat\lambda^n-\hat\lambda_i^n\big){\langle \psi,\hat\psi_i^n\rangle}_2^2\ge 0. 
\end{align}
Thus, $(\hat\psi^n,\hat\lambda^n)$ satisfies the \index{Optimality conditions!second-order!necessary} {\em second-order necessary optimality conditions} for a minimum. Additionally, we get that every eigenvalue-eigenvector pair $(\hat{\lambda}_i^n,\hat{\psi}_i^n)$ with $\lambda_i^n < \lambda_1^n$ does not fulfill the necessary second-order condition.

The \index{Optimality conditions!second-order!sufficient}{\em second-order sufficient optimality conditions} would require $\nabla_{\psi \psi} \mathcal L(\hat\psi^n,\hat \lambda^n)$ to be positive definite on the tangential hyperplane
\begin{align*}
    \ker \big(\nabla e(\hat\psi^n)^\top\big) = \big \{ \psi \in \mathbb R^m \,\big|\, {\langle\hat\psi^n,\psi\rangle}_2= 0 \big \} = \{ \hat\psi^n \}^\perp = \mathrm{span}\,\big\{ \hat{\psi}_2^n,\ldots,\hat{\psi}_m^n \big\}.
\end{align*}
Inserting $\hat{\psi}_i^n$ into \eqref{SIAM:Eq-I.1.1.43-A} for an arbitrary $i \in \{2,\ldots,m \}$, we get the condition $\hat\lambda^n > \hat\lambda_i^n$. Therefore, the sufficient second-order conditions hold true if $\hat\lambda_1^n > \hat\lambda_2^n$. Otherwise, $\psi = \hat\psi_2^n$ is a non-zero vector from the tangential hyperplane that satisfies $\langle\hat\psi_2 ^n,\nabla_{\psi\psi}\mathcal L(\hat\psi^n,\hat\lambda^n)\hat\psi_2^n\rangle_{\mathbb R^m} = 0$, so that $\hat\psi^n$ is only a saddle point. In this case, $\hat\lambda_1^n=\hat\lambda_2^n$ implies that $\hat\psi_2^n$ is also a global solution to \eqref{SIAM:Eq-I.1.1.55}, as can be seen from \eqref{SIAM:Eq-I.1.1.19}. To prove that $\hat\psi^n$ actually solves \eqref{SIAM:Eq-I.1.1.55} one has to proceed as in the proof of Theorem~\ref{SIAM:Theorem-I.1.1.3}. The observation that a multiplicity of solutions to \eqref{SIAM:Eq-I.1.1.55} coincides with second-order sufficient optimality conditions not holding true carries over to the case $\ell>1$: If $\hat\lambda_\ell^n > \hat\lambda_{\ell+1}^n$, then it can be shown that second-order sufficient optimality conditions hold for $\hat\psi_1^n,\hdots,\hat\psi_\ell^n$.

\subsection{POD with weighted inner product}
\label{SIAM:Section-2.1.1.3}

Let us return to the unitary case of Section \ref{SIAM:Section-2.1.1.2}, but this time with the inclusion of the weights $\alpha_1^n,\ldots,\alpha_n^n$ from the original problem \eqref{SIAM:Eq-I.1.1.2} and the inner product in $X=\mathbb C^m$ being given by a positive definite and Hermitian \index{Matrix!weighting, $\bW$}{\em weighting matrix} $\bW\in\mathbb C^{m \times m}$:
\begin{equation}
    \label{SIAM:Eq-I.1.1.61}
    {\langle \psi,\phi\rangle}_\bW=\psi^\top \overline{\bW\phi}={\langle \psi,\bW\phi\rangle}_\mathbb C^m\quad\text{for }\psi,\phi\in\mathbb C^m
\end{equation}
The associated norm is defined as $|\cdot|_\bW=\langle\cdot\,,\cdot\rangle_\bW^{1/2}$. From $\bW=\bW^\mathsf H$ we infer that
\begin{align*}
    {\langle \psi,\psi\rangle}_\bW=\psi^\top \overline{\bW\psi}=\big(\psi^\top \overline{\bW\psi}\big)^\top=\overline\psi^\top\bW^\mathsf H\psi=\overline\psi^\top\bW\psi=\overline{\psi^\top\overline{\bW\psi}}=\overline{{\langle \psi,\psi\rangle}}_\bW
\end{align*}
for every $\psi \in \mathbb C^m$. Consequently, $\langle \psi,\psi\rangle_\bW\in\mathbb R$ holds for every $\psi \in \mathbb C^m$.

We choose ${K}=1$ and $\omega_1^{K}=1$. Then problem \eqref{SIAM:Eq-I.1.1.2} has the form
\begin{align*}
    \min\sum_{j=1}^n\alpha_j^n\,\Big| y_j-\sum_{i=1}^\ell{\langle y_j,\psi_i\rangle}_\bW\,\psi_i\Big|_\bW^2\text{ s.t. } \{\psi_i\}_{i=1}^\ell\subset \mathbb C^m\text{ and }{\langle\psi_i,\psi_j\rangle}_\bW=\delta_{ij},~1 \le i,j \le \ell,
\end{align*}
where the associated induced norm is given as usual by $|\cdot|_\bW=\langle\cdot\,,\cdot\rangle_\bW^{1/2}$. Let us now motivate the inner product by the three examples from Section \ref{sec:threeExamples}. In all cases, we are dealing with discretized representations of infinite-dimensional functions, for which the inner product matrix $\bW$ is supposed to represent the finite-dimensional version of the underlying function space's inner product. 
\begin{example}
    \label{SIAM:Example-I.1.1.1}
    \rm
    \begin{enumerate}
        \item [1)] In Example~\ref{cosExample_podProb}, each data vector $y_j \in \mathbb R^{n_{\bx}}$ was representative of the function $y(t_j,\cdot): (\bx_\mathsf a, \bx_\mathsf b) = \Omega \to \mathbb R$. A natural inner product space for these types of functions is given by the space $H=L^2(\Omega)$ of square integrable functions; cf. Definition \ref{Definition_L2_BanachValued}: 
        \begin{align*}
            {\langle f,g \rangle}_H= \int_{\bx_\mathsf a}^{\bx_\mathsf b} f(\bx) g(\bx)\,\mathrm d\bx\quad\text{for }f,g\in H.
        \end{align*}
        We want to approximate the inner product in $H$ by a quadrature rule. This can be done by introducing an appropriate weighting matrix $\bW \in \mathbb R^{n_{\bx}\times n_{\bx}}$ and a weighted inner product in $\mathbb R^{n_\bx}$. A natural way to do this has been already introduced in Example \ref{cosExample_podProb} itself by use of trapezoidal integration: it holds for all $f,g\in H$
        \begin{align*}
            {\langle f,g \rangle}_H\approx h \bigg( \frac{\mathrm f_1^h\mathrm g_1}{2} + \sum_{i=2}^{n_{\bx}-1} \mathrm f_i^h\mathrm g_i^h + \frac{\mathrm f_{n_\bx}^h\mathrm g_{n_\bx}^h}{2} \bigg) = (\mathrm f^h)^\top \bW\mathrm g^h,
        \end{align*}
        where the vector $\mathrm f^h \in \mathbb R^{n_{\bx}}$ is the discretization of $f$ given by the $n_\bx$ components
        \begin{align*}
            \mathrm f^h_i=\left\{
            \begin{aligned}
                &\frac{2}{h}\int_{\bx_1}^{\bx_1+h/2}f(\bx)\,\mathrm d\bx&&\text{for }i=1,\\
                &\frac{1}{h}\int_{\bx_i-h/2}^{\bx_i+h/2}f(\bx)\,\mathrm d\bx&&\text{for }i=2,\ldots,n_\bx-1,\\
                &\frac{2}{h}\int_{\bx_{n_\bx}-h/2}^{\bx_{n_\bx}}f(\bx)\,\mathrm d\bx&&\text{for }i=n_\bx.
            \end{aligned}
            \right.
        \end{align*}
        Analogously, the vector $\mathrm g^h=(\mathrm g_i^h) \in \mathbb R^{n_{\bx}}$ is defined. Then, the matrix $\bW$ is diagonal and given by $\bW = \text{diag}(h/2,h,\ldots,h,h/2)$. 
        \item [2)] For the parabolic problem in Example \ref{Example:PODParabolic}, it was already argued that two natural choices for the underlying infinite-dimensional function space are given by $H = L^2(\Omega)$ and $V = H^1(\Omega)$. We apply a standard piecewise linear \index{Method!finite element, FE}{\em finite element (FE) discretization} with $m \in \mathbb N$ degrees of freedom; cf., e.g., \cite{BS08,Tho97}. Let $\{\varphi_i\}_{i=1}^m \subset H^1(\Omega)$ denote the FE piecewise ansatz functions where $\varphi_i$ takes the value one on the $i$-th node of the grid and zero on all other nodes. Then the discretized FE space is given by the linear hull
        \begin{align*}
            X^h=\bigg\{\varphi\in V\,\Big|\,\varphi(\bx)=\sum_{i=1}^m \mathrm v_i \varphi_i(\bx)\text{ for }\mathrm v=(\mathrm v_i)\in \mathbb R^m\text{ and }\bx \in \Omega\bigg\}.
        \end{align*}
        We introduce the \index{Matrix!mass!FE, $\bM$}{\em mass matrix} $\bM = ((M_{ij})) \in \mathbb R^{m \times m}$ with the elements $M_{ij} = \langle \varphi_j, \varphi_i \rangle_{L^2(\Omega)}$ for $1 \le i,j \le m$, and the \index{Matrix!stiffness!FE, $\bS$}{\em stiffness matrix} $\bS = ((S_{ij})) \in \mathbb R^{m \times m}$ with the elements $S_{ij} = \langle \varphi_i, \varphi_j \rangle_V$. Note that both $\bM$ and $\bS$ are symmetric and positive definite. Furthermore, the $H$- and $V$-inner product of two FE function 
        \begin{align*}
            v(\bx)=\sum_{i=1}^m \mathrm v_i \varphi_i(\bx),\quad w(\bx)=\sum_{i=1}^m \mathrm w_i \varphi_i(\bx),\quad\bx\in\Omega,~\mathrm v_i,\mathrm w_i\in\mathbb R
        \end{align*}
        can be described by
        \begin{align*}
            {\langle v,w \rangle}_H&=\sum_{i=1}^m\sum_{j=1}^m \mathrm v_i  M_{ij}\mathrm w_j=\mathrm v^\top\bM\mathrm w=:{\langle\mathrm v,\mathrm w\rangle}_\bM,\\
            {\langle v,w\rangle}_V&=\sum_{i=1}^m\sum_{j=1}^m \mathrm v_i  S_{ij}\mathrm w_j=\mathrm v^\top\bS\mathrm w=:{\langle\mathrm v,\mathrm w\rangle}_\bS
        \end{align*}
        with the coefficient vectors $\mathrm v=(\mathrm v_i),\,\mathrm w=(\mathrm w_i)\in\mathbb R^m$. Both $\bM$ and $\bS$ are possible choices for the inner product matrix $\bW$.
        \item [3)] For the elliptic problem in Example \ref{Example:PODElliptic}, we may also utilize the matrices $\bM$ and $\bS$ from part 2). The only difference is that the underlying function space is now given by $L^2(\Omega;\mathbb C)$, so we have to allow complex vectors $\mathrm v, \mathrm w \in \mathbb C^m$ as FE coefficients. 
    \end{enumerate}

    \medskip\noindent
    Both for the parabolic and the elliptic problem considered in 2) and 3), we utilize a Galerkin discretization to project the infinite-dimensional PDEs to the FE space, which will be explained in detail in Section \ref{SIAM-Book:Section3.4}. This gives us ansatz fuctions $y^h(t)$ for all times $t \in (0,T)$ in the parabolic case and $p^h(f) \in X^h$ for all frequencies $f \in \mathscr F$ in the elliptic case. Let $\mathrm y(t) \in \mathbb R^m$ and $\mathrm p(f) \in \mathbb C^m$ denote their coefficient vectors. With the given discretizations of $(0,T)$ and $\mathscr F$ in Examples \ref{Example:PODParabolic} and \ref{Example:PODElliptic}, this leads to data trajectories $\mathrm y(t_j) \in \mathbb R^m$ and $\mathrm p(f_j) \in \mathbb C^m$ for $j=1,...,n$. The discrete versions of the associated POD problems are then given by
    \begin{align*}
        \min \sum_{j=1}^n \alpha_j^n\,\Big| \mathrm y(t_j) - \sum_{i=1}^\ell{\langle \mathrm y(t_j), \mathrm v_i \rangle}_{\bW}\, \mathrm v_i \Big|_{\bW}^2\text{ s.t. }\{\mathrm v_i\}_{i=1}^\ell \subset \mathbb R^m \text{ and } {\langle \mathrm v_i, \mathrm v_j \rangle}_{\bW} = \delta_{ij}, ~1 \le i,j \le \ell,
    \end{align*}
    and
    \begin{align*}
        \min\sum_{j=1}^n \alpha_j^n\,\Big| \mathrm p(f_j) - \sum_{i=1}^\ell {\langle \mathrm p(f_j), \mathrm v_i \rangle}_{\bW}\,\mathrm v_i \Big|_{\bW}^2\text{ s.t. }\{\mathrm v_i\}_{i=1}^\ell \subset \mathbb C^m \text{ and } {\langle \mathrm v_i, \mathrm v_j \rangle}_{\bW} = \delta_{ij}, ~1 \le i,j \le \ell.
    \end{align*}
    In both cases, $\bW$ denotes either $\bM$ or $\bS$ and $\mathrm v_i$ is the coordinate vector of the $i$-th POD basis vector $\psi_i \in X_h$, which lives in the FE space $X^h$. \hfill
\end{example}

As in Section~\ref{SIAM:Section-2.1.1.2} we introduce the snapshot matrix $\bY=[y_1\,|\ldots\,|\,y_n]\in\mathbb C^{m\times n}$ with rank $d^n\le\min(m,n)$. Moreover, we define the diagonal matrix $\bD=\mathrm{diag}\,(\alpha_1^n,$ $\ldots,\alpha_n^n)\in\mathbb R^{n\times n}$. We find that
\begin{align*}
    \big(\mathcal R^n \psi\big)_i&=\Big(\sum_{j=1}^n\alpha_j^n\,{\langle \psi,y_j\rangle}_W\,y_j\Big)_i=\sum_{j=1}^n\sum_{l=1}^m \sum_{\nu=1}^m \alpha_j^n\psi_l W_{l\nu}\overline Y_{\nu j}Y_{ij}\\
    &=\big(\bY\bD\bY^\mathsf H \bW\psi\big)_i\quad \text{for all } \psi\in\mathbb C^m \text{ and } 1\le i \le m.
\end{align*}
Consequently, \eqref{SIAM:Eq-I.1.1.18} leads to the eigenvalue problem
\begin{equation}
    \label{SIAM:Eq-I.1.1.70}
    \bY\bD\bY^\mathsf H\bW\hat\psi_i^n=\hat\lambda_i^n\hat\psi_i^n,\quad\hat\lambda_1^n\ge\ldots\ge \hat\lambda_{d^n}^n>\hat\lambda_{d^n+1}^n=\ldots =\hat\lambda_m^n=0.
\end{equation}
Due to the fact that $\bW$ is a Hermitian positive-definite matrix, it can be decomposed as $\bW=\bQ\bB\bQ^\mathsf H$, where the diagonal matrix $\bB=\mathrm{diag}\,(\beta_1,\ldots,\beta_m)\in\mathbb R^{m\times m}$ contains the eigenvalues $\beta_1\ge\ldots\ge\beta_m>0$ of $\bW$ and $\bQ\in \mathbb C^{m \times m}$ is a unitary matrix, i.e. $\bQ^{-1}=\bQ^\mathsf H$ holds. We have
\begin{align*}
    \bW^r=\bQ\bB^r\bQ^\mathsf H\quad\text{for }r\in\mathbb R,
\end{align*}
with $\bB^r=\mathrm{diag}\,(\beta_1^r,\ldots,\beta_m^r)$. Note that $(\bW^r)^{-1}=\bW^{-r}$ and $\bW^{r+s}=\bW^r \bW^s$ for $r,\,s \in \mathbb R$. Moreover,
\begin{align*}
    {\langle \psi,\phi \rangle}_\bW&=\psi^\top\overline{\bW\phi}=\psi^\top\overline{\bQ\bB\bQ^\mathsf H\phi}=\psi^\top\overline{\bQ}\bB^{1/2}\bQ^\top  \overline{\bQ\bB^{1/2}\bQ^\mathsf H\phi}\\
    &=\big(\bQ\bB^{1/2}\bQ^\mathsf H\psi\big)^\top\overline{\bQ\bB^{1/2}\bQ^\mathsf H\phi}=\big(\bW^{1/2}\psi\big)\overline{\bW^{1/2}\phi}={\langle \bW^{1/2} \psi,\bW^{1/2}\phi\rangle}_{\mathbb C^m} 
\end{align*}
and $| \psi |_W=|\bW^{1/2} \psi |_2$ for $\psi,\phi \in \mathbb C^m$. Note that $\bD^{1/2}=\mathrm{diag}\,(\alpha_1^{1/2},\ldots,\alpha_n^{1/2})$ holds. Inserting $\psi_i^n=\bW^{1/2} \hat\psi_i^n$ in \eqref{SIAM:Eq-I.1.1.70}, multiplying \eqref{SIAM:Eq-I.1.1.70} by $\bW^{1/2}$ from the left and setting $\tilde{\bY}=\bW^{1/2}\bY\bD^{1/2}$ yield the symmetric $m\times m$ eigenvalue problem
\begin{align*}
    \tilde{\bY}\tilde{\bY}^\mathsf H\psi_i^n=\hat\lambda_i^n\psi_i^n \quad \text{for } 1\le i\le\ell.
\end{align*}
Note that
\begin{equation}
    \label{SIAM:Eq-I.1.1.74}
    \tilde{\bY}^\mathsf H\tilde{\bY}=\bD^{1/2} \bY^\mathsf H \bW \bY \bD^{1/2} \in \mathbb C^{n \times n}.
\end{equation}
Thus, the POD basis $\{\hat\psi_i^n\}_{i=1}^\ell$ of rank $\ell$ can also be computed by the methods of snapshots as follows: First solve the symmetric $n \times n$ eigenvalue problem
\begin{align*}
    \tilde{\bY}^\mathsf H\tilde{\bY}\phi_i^n=\hat\lambda_i^n\phi_i^n,~1\le i\le\ell\quad\text{and} \quad {\langle \phi_i^n,\phi_j^n \rangle}_2=\delta_{ij} \quad \text{for } 1 \le i,j \le \ell.
\end{align*}
Then we set (by using the SVD of $\tilde{\bY}$)
\begin{equation}
    \label{SIAM:Eq-I.1.1.76}
    \hat\psi_i^n=\bW^{-1/2} \psi_i^n=\frac{1}{\hat\sigma_i^n}\,\bW^{-1/2}\tilde{\bY}\phi_i^n=\frac{1}{\hat\sigma_i^n}\,\bY\bD^{1/2}\phi_i^n \quad \text{for } 1\le i\le \ell.
\end{equation}
We find that
\begin{align*}
    {\langle \hat\psi_i^n,\hat\psi_j^n\rangle}_\bW&=(\hat\psi_i^n)^\top \overline{\bW\hat\psi_j^n}=\frac{1}{\hat\sigma_i^n\hat\sigma_j^n}\,(\phi_i^n)^\top \overline{\bD^{1/2} \bY^\mathsf H \bW \bY \bD^{1/2}\phi_j^n}\\
    &=\frac{1}{\hat\sigma_i^n\hat\sigma_j^n}\,{\langle\phi_i^n,\tilde{\bY}^\mathsf H\tilde{\bY}\phi_j^n\rangle}_2=\delta_{ij} \quad \text{for } 1 \le i,j \le \ell.
\end{align*}
Thus, the POD basis $\{\hat\psi_i^n\}_{i=1}^\ell$ of rank $\ell$ is orthonormal in $\mathbb C^m$ with respect to the inner product $\langle\cdot \, , \cdot \rangle_\bW$. We observe from \eqref{SIAM:Eq-I.1.1.74} and \eqref{SIAM:Eq-I.1.1.76} that the computation of $\bW^{1/2}$ and $\bW^{-1/2}$ is not required for this last approach. For applications like Example~\ref{SIAM:Example-I.1.1.1}, where $\bW$ is not just a diagonal matrix, the method of snapshots turns out to be more attractive with respect to the computational costs even if $n>m$ holds. It holds
\begin{equation}
    \label{SIAM:Eq-I.1.1.78}
    \sum_{j=1}^n\alpha_j^n\,{|y_j|}_\bW^2={\|\tilde{\bY}\|}_F^2 \quad \text{for } j=1,\ldots,n.
\end{equation}
Thus, we infer from \eqref{SIAM:Eq-I.1.1.39b} that the \index{POD method!discrete variant!relative total energy}{\em relative total energy} is given as
\begin{align*}
    \mathcal E(\ell)=\frac{\sum_{i=1}^\ell\hat\lambda_i^n}{{\|\tilde{\bY}\|}_F^2}.
\end{align*}
Let us summarize all computational steps in Algorithm~\ref{SIAM:Algorithm-I.1.1.3}.

\bigskip
\hrule
\vspace{-3.5mm}
\begin{algorithm}[(POD in $\mathbb C^m$ with a weighted inner product)]
    \label{SIAM:Algorithm-I.1.1.3}
    \vspace{-3mm}
    \hrule
    \vspace{0.5mm}
    \begin{algorithmic}[1]
        \REQUIRE Snapshots $\{y_j\}_{j=1}^n \subset \mathbb C^m$, positive weights $\{\alpha_j^n\}_{j=1}^n$, POD rank $\ell \le d^n$, Hermitian positive-definite matrix $\bW \in \mathbb C^{m \times m}$ and {\tt flag} for the eigenvalue solver;
        \STATE Set $\bY=[y_1,\ldots,y_n] \in \mathbb C^{m \times n}$ and $\bD=\mathrm{diag}\,(\alpha_1^n,\ldots,\alpha_n^n)\in\mathbb R^{n\times n}$;
        \IF{{\tt flag} = 0}
            \STATE Determine $\tilde{\bY}=\bW^{1/2}\bY\bD^{1/2} \in \mathbb C^{m \times n}$;
            \STATE Compute singular value decomposition $[\bPsi,\bSigma,\bPhi]=\mathrm{svd} \,(\tilde{\bY})$;
            \STATE Set $\hat\psi_i^n=\bW^{-1/2}\bPsi_{\cdot,i}\in \mathbb C^m$ and $\hat\lambda_i^n=\Sigma_{ii}^2$ for $i=1,\ldots,\ell$;
            \STATE Compute $\mathcal E(\ell)=\sum_{i=1}^\ell \lambda_i/{\|\tilde{\bY}\|}_F^2$;
        \ELSIF{{\tt flag} = 1}
            \STATE Determine $\bR=\bW^{1/2}\bY\bD\bY^\mathsf H\bW^{1/2} \in \mathbb C^{m \times m}$;
            \STATE Compute eigenvalue decomposition $[\bPsi,\bLambda]=\mathrm{eig}\, (\bR)$;
            \STATE Set $\hat\psi_i^n=\bW^{-1/2}\bPsi_{\cdot,i}\in \mathbb C^m$ and $\hat\lambda_i^n=\Lambda_{ii}$ for $i=1,\ldots,\ell$;
            \STATE Compute $\mathcal E(\ell)=\sum_{i=1}^\ell \lambda_i/{\|\tilde{\bY}\|}_F^2$;
        \ELSIF{{\tt flag} = 2}
            \STATE Determine $\bK=\bD^{1/2}\bY\bW\bY^\mathsf H\bD^{1/2} \in \mathbb C^{n \times n}$;
            \STATE Compute eigenvalue decomposition $[\bPhi,\bLambda]=\mathrm{eig}\, (\bK)$;
            \STATE Set $\hat\lambda_i^n=\Lambda_{ii}$ and $\hat\psi_i^n=\bY\bD^{1/2}\bPhi_{\cdot,i}/\Lambda_{ii}^{1/2}\in \mathbb C^m$ for $i=1,\ldots,\ell$;
            \STATE Compute $\mathcal E(\ell)=\sum_{i=1}^\ell \lambda_i/{\|\tilde{\bY}\|}_F^2$;
        \ENDIF
        \RETURN POD basis $\{\hat\psi_i^n\}_{i=1}^\ell$, eigenvalues $\{\hat\lambda_i^n\}_{i=1}^\ell$ and ratio $\mathcal E(\ell)$;
    \end{algorithmic}
    \hrule
\end{algorithm}

\begin{example}
    \label{Example:PODParabolic_timeComparison}
    \rm We want to apply all three strategies of Algorithm~\ref{SIAM:Algorithm-I.1.1.3} to Example~\ref{Example:PODParabolic}. In addition to the temporal discretization mentioned in that example, the infinite-dimensional state equation \eqref{SIAM:EqMotPDE1} is discretized using linear finite elements, thus replacing the state space $V = H^1(\Omega)$ by a linear space $V^m \subset V$ of dimension $m \in \mathbb N$. This discretization will be carried out in detail in Section~\ref{SIAM-Book:Section3.4}, so we do not go into more detail right now. The important part is that we have different accuracy levels indicated by the numbers $n$ for the temporal discretization and $m$ for the spatial discretization. Since the three strategies in Algorithm~\ref{SIAM:Algorithm-I.1.1.3} solve different kinds of decomposition problems, we would like to compare their computational efficiency for three different sizes of $m$ and $n$. In particular, we choose $n \in \{100,200,400\}$ and $m \in \{765,1319,2999\}$ and apply all three strategies to the problem using $K=20$ different trajectories, exactly as in Problem~\ref{Example:PODParabolic}. We fix the basis length at $\ell=5$. Table~\ref{tab:podParabolic_timeComparison} shows the results.
	\begin{table}[ht]
	    \centering
	    \begin{tabular}{|c|*{3}{|c}|*{3}{|c}|*{3}{|c}|}\hline
	    	{} & \multicolumn{3}{c||}{\texttt{flag=0}} & \multicolumn{3}{c||}{\texttt{flag=1}} & \multicolumn{3}{c|}{\texttt{flag=2}} \\
		    \hline \backslashbox{$m$}{$n$} & \makebox[1.5em]{100} & \makebox[1.8em]{200} & \makebox[1.8em]{400} & \makebox[1.8em]{100} & \makebox[1.8em]{200} & \makebox[1.8em]{400} & \makebox[1.8em]{100} & \makebox[1.8em]{200} & \makebox[1.8em]{400} \\
		    \hline \makebox[3em]{\phantom{1}765} & 0.2 & \phantom{1}0.5 & \phantom{1}1.5 & \phantom{1}0.3 & \phantom{1}0.3 & \phantom{1}0.3 & 2.4 & 16.4 & 117.6 \\
		    \hline \makebox[3em]{1319} & 0.7 & \phantom{1}1.2 & \phantom{1}3.2 & \phantom{1}1.0 & \phantom{1}1.1 & \phantom{1}1.3 & 2.4 & 16.6 & 116.0 \\
		    \hline \makebox[3em]{2999} & 4.0 & 11.2 & 13.0 & 10.2 & 10.3 & 10.8 & 2.6 & 17.3 & 118.5 \\\hline
	   \end{tabular}
	   \caption{Example~\ref{Example:PODParabolic_timeComparison}. CPU-times for the POD basis computation by the three different strategies from Algorithm \ref{SIAM:Algorithm-I.1.1.3}. All times are in seconds, with computations carried out sequentally on a system with 3.5GHz CPU and 16GB RAM.}
	   \label{tab:podParabolic_timeComparison}
    \end{table}
    We can immediately tell that computation times for \texttt{flag=1} depend mostly on $m$ while those for \texttt{flag=2} depend on $n$. This makes sense since we mostly have to solve eigenvalue problems of size $(nK) \times (nK)$ and $m \times m$ in both respectively. This underlines the importance of using the method of snapshots ($\texttt{flag=2}$) only for comparably small numbers of $nK$. Alternatively, the singular value decomposition in the case of \texttt{flag=0} offers a compromise in terms of efficiency, never giving excessively high computation times but being beaten for large quotients of $(nK)/m$ by \texttt{flag=1} and for small quotients of $(nK)/m$ by \texttt{flag=2}. \hfill$\blacklozenge$
\end{example}

\subsection{POD for finite-dimensional dynamical systems}
\label{SIAM:Section-2.1.1.4}

In Example \ref{SIAM:Example-I.1.1.1}-b), we have already encountered a finite-dimensional dynamical system which arose from a semi-discretization of a parabolic PDE. In this section, we will study how POD can be applied to these types of systems. We will focus on linear systems where the linear operator and the inhomogeneity may additionally depend on a parameter from a compact set, which will result in $K>1$ data trajectories for the POD problem. 

Let us choose $X=\mathbb C^m$ endowed with the weighted inner product \eqref{SIAM:Eq-I.1.1.61}. Assume that $\mathscr N\subset\mathbb R^\mathfrak p$ is a compact set and $\mathfrak p\in\mathbb N$. For $T>0$ we consider the parametrized linear initial value problem
\begin{equation}
    \label{SIAM:Eq-I.1.1.87}
    \dot y(t)=\bA(\bnu)y(t)+f(t;\bnu)\quad \text{for } t \in (0,T],\quad y(0)=y_\circ,
\end{equation}
where $y_\circ \in \mathbb C^m$ is a given initial condition and $\bA(\bnu) \in \mathbb C^{m \times m}$ is a coefficient matrix, which depends on the parameter $\bnu\in\mathscr D$. Moreover, the parameter-dependent mapping $f(\cdot\,;\bnu):[0,T] \times \mathbb C^m \to \mathbb C^m$ is continuous for any parameter $\bnu\in\mathscr N$. It is well-known that for any $\bnu$ there exists a unique (classical) solution $y=y(\bnu) \in C^1(0,T;\mathbb C^m) \cap C([0,T]; \mathbb C^m)$ to \eqref{SIAM:Eq-I.1.1.87}. Recall that the solution $y$ to \eqref{SIAM:Eq-I.1.1.87} is given by the integral representation
\begin{align*}
    y(t;\bnu)=e^{\bA(\bnu)t} y_\circ+\int_0^t e^{\bA(\bnu)(t-s)} f(s;\bnu) \, \mathrm ds\quad\text{for }\bnu\in\mathscr N
\end{align*}
with $e^{\bA(\bnu) t}=\sum_{n=0}^\infty t^n\bA(\bnu)^n/n!$ is the exponential of matrix $\bA(\bnu)$; cf. \cite{Paz83}.

\begin{example}
    \rm For the spatial domain $\Omega=(0,1)\subset\mathbb R$ we define the space-time cylinder $Q=(0,T)\times\Omega$. The parameter set $\mathscr D$ is given by the three-dimensional rectangular domain
    \begin{align*}
        \mathscr N=\big[\upnu_{\mathsf a 1},\upnu_{\mathsf b 1}\big]\times\big[\upnu_{\mathsf a 2},\upnu_{\mathsf b 2}\big]\times\big[\upnu_{\mathsf a 3},\upnu_{\mathsf b 3}\big]=:\big[\bnu_{\mathsf a},\bnu_{\mathsf b}\big]\subset\mathbb R^3
    \end{align*}
    with $\bnu_{\mathsf a}=(\upnu_{\mathsf a i})_{1\le i\le3},\bnu_{\mathsf b}=(\upnu_{\mathsf b i})_{1\le i\le3}\in\mathscr N$ and $0\le\upnu_{\mathsf a i}\le\upnu_{\mathsf b i}$ for $i=1,2,3$. For any $\bnu=(\upnu_i)_{1\le i\le3}\in\mathscr N$ let us consider the following one-dimensional advection-diffusion equation:
    \begin{equation}
        \label{heat-ex}
        \begin{aligned}
            &y_t(t,\bx)=\upnu_1 c(\bx)y_{\bx\bx}(t,\bx)+\upnu_2 v(\bx)y_\bx(t,\bx)+\upnu_3 f(t,\bx),&&(t,\bx)\in Q,\\
            &y(t,0)=y(t,1)=0,&&t \in (0,T],\\
            &y(0,\bx)=y_\circ(\bx),&&\bx\in\Omega,
        \end{aligned}
    \end{equation}
    where $c,v,y_\circ$ are given functions in $C(\overline\Omega)$ with $c(\bx)\geq c_0>0$ for all $\bx\in\Omega$ and a given inhomogeneity $f\in C(\overline Q)$. To solve \eqref{heat-ex} numerically we apply a standard finite difference approximation for the spatial variable $x$; cf., e.g., \cite{Str04}. To simplify the presentation we suppose the spatial grid is equidistant, i.e., $x_i=ih$ for $0\le i\le m+1$ with a step size $h=1/(m+1)$. Setting $\mathrm c_i=c(\bx_i)$, $\mathrm v_i=v(\bx_i)$, $\mathrm y_{\circ i}=y_\circ(\bx_i)$ and $\mathrm f_i(t)= f(t,\bx_i)$ for $1 \le i\le m$, $y\in\mathbb R$ and $\bnu\in\mathscr N$,  we consider instead of \eqref{heat-ex} the problem
    \begin{subequations}
        \label{heat-ex-d}
        \begin{align}
            \label{heat-ex-da}
            y_t(t,\bx_i)&=\upnu_1\mathrm c_i y_{xx}(t,\bx_i)+\upnu_2 \mathrm v_iy_\bx(t,\bx_i)+\upnu_3 \mathrm f_i(t),&&\hspace{-2.5mm}t\in(0,T],1\le i\le m,\\
            \label{heat-ex-db}
            y(t,\bx_0)&=y(t,\bx_{m+1})=0,&&\hspace{-2.5mm}t \in (0,T],\\
            \label{heat-ex-dc}
            y(0,\bx_i)&=\mathrm y_{\circ i},&&\hspace{-2.5mm}1\le i\le m.
        \end{align}
    \end{subequations}
    The first and second partial derivatives of $y$ with respect to $x$ are discretized by centered difference quotients of second-order. Then we have (cf. \cite{Str04})
    \begin{equation}
        \label{SIAM:Eq_CentrDiff}
        \begin{aligned}
            y_{\bx\bx}(t,\bx_i)&=\frac{y(t,\bx_{i+1})-2y(t,\bx_i)+y(t,\bx_{i-1})}{h^2}+\mathcal O\big(h^2\big),\\
            y_\bx(t,\bx_i)&=\frac{y(t,\bx_{i+1})-y(t,\bx_{i-1})}{2h}+\mathcal O\big(h^2\big)
        \end{aligned}
    \end{equation}
    for $i\in\{1,\ldots,m\}$ and $t\in[0,T]$ provided the fourth partial derivative of $v$ with respect to $x$ exists for any $t\in[0,T]$ and is continuous on $\overline Q$. Let us denote by $\mathrm y_i:[0,T]\to\mathbb R$ the numerical approximation for $y(\cdot\,,\bx_i)$ for $i=1,\ldots,m$. From \eqref{heat-ex-db} we get $\mathrm y_0=\mathrm y_{m+1}=0$ on $[0,T]$ so that we obtain the following ordinary differential equations for the time-dependent functions $\mathrm y_i$:
    \begin{subequations}
        \label{heat-ex-2}
        \begin{equation}
            \label{heat-ex-2a}
            \left\{
            \begin{aligned}
                \dot{\mathrm y}_1(t)&=a_{11}(\bnu)\mathrm y_1(t)+a_{12}(\bnu)\mathrm y_2(t)+\upnu_3\mathrm  f_1(t),\\
                \dot{\mathrm y}_i(t)&=a_{i,i-1}(\bnu)\mathrm y_{i-1}(t)+a_{ii}(\bnu)\mathrm y_i(t)+a_{i,i+1}(\bnu)\mathrm y_{i+1}(t)+\upnu_3 \mathrm f_i(t),\quad i=2,\ldots,m-1,\\
                \dot{\mathrm y}_m(t)&=a_{m,m-1}(\bnu)\mathrm y_{m-1}(t)+a_{mm}(\bnu)\mathrm y_m(t)+\upnu_3 \mathrm f_m(t)
            \end{aligned}
            \right.
        \end{equation}
        with the coefficients
        \begin{align*}
            a_{i,i-1}(\bnu)=\frac{\upnu_1 \mathrm c_i}{h^2}-\frac{\upnu_2 \mathrm v_i}{2h},\quad a_{ii}(\bnu)=-\frac{2\upnu_1 \mathrm c_i}{h^2},\quad a_{i,i+1}(\bnu)=\frac{\upnu_1 \mathrm c_i}{h^2}+\frac{\upnu_2 \mathrm v_i}{2h}.
        \end{align*}
        From \eqref{heat-ex-dc} we infer the initial condition
        \begin{equation}
            \mathrm y_i(0)=\mathrm y_{\circ i}\quad\text{for }i=1,\ldots,m.
        \end{equation}
    \end{subequations}
    Introducing the matrix $\bA(\bnu)=((a_{ij}(\bnu)))\in\mathbb R^{m\times m}$ and the vectors
    \begin{align*}
        \mathrm y(t)=\left(
        \begin{array}{c}
            \mathrm y_1(t)\\
            \vdots\\
            \mathrm y_m(t)
        \end{array}
        \right),\quad\mathrm f(t)=\left(
        \begin{array}{c}
            \mathrm f_1(t)\\
            \vdots\\
            \mathrm f_m(t)
        \end{array}
        \right)\mbox{ for } t \in [0,T],\quad\mathrm y_\circ=\left(
        \begin{array}{c}
            \mathrm y_{\circ 1}\\
            \vdots\\
            \mathrm y_{\circ m}
        \end{array}
        \right)\in\mathbb R^m
    \end{align*}
    we can express \eqref{heat-ex-2} in the form
    \begin{equation}
        \label{heat-ex-3}
        \dot{\mathrm y}(t)=\bA(\bnu)\mathrm y(t)+\upnu_3 \mathrm f(t) \quad\text{for } t \in (0,T],\quad \mathrm y(0)=\mathrm y_\circ.
    \end{equation}
    The vector $\mathrm y(t)$, $t \in [0,T]$, represents a function in $\Omega$ evaluated at $m$ grid points. Therefore, we should supply $\mathbb R^m$ by a weighted inner product representing a discretized inner product in an appropriate function space. Here we choose the inner product introduced in Example~\ref{SIAM:Example-I.1.1.1}.\hfill$\blacklozenge$
\end{example}

Let $0=t_1<t_2<\ldots<t_n=T$ be a given time grid in the interval $[0,T]$, where the step sizes are given as $\delta t_j=t_j-t_{j-1}$ for $j=2,\ldots,n$. We suppose that we know the solution $y=y(t_j;\bnu)$ to \eqref{SIAM:Eq-I.1.1.87} at the given time instances $t_j$, $j \in \{1,\ldots,n\}$ and for the parameters $\bnu_k\in\mathscr N$, $k\in\{1,\ldots,{K}\}$. Our goal is to determine a POD basis of rank $\ell \le \min\{m,n\}$ for the snapshot subspace
\begin{align*}
    \mathscr V^n=\mathrm{span}\,\big\{y_j^k=y(t_j;\bnu_k)\,\big|\,1\le j\le n\text{ and }1\le k\le K\big\}\subset\mathbb C^m.
\end{align*}
For $1\le k\le K$ we find that
\begin{align*}
    \sum_{j=1}^n\alpha_j^n\,\Big|y_j^k-\sum_{i=1}^\ell{\langle y_j^k,\psi_i\rangle}_\bW\,\psi_i\Big|_\bW^2=\sum_{j=1}^n\alpha_j^n\,\Big|y(t_j;\bnu^k)-\sum_{i=1}^\ell{\langle y(t_j;\bnu^k),\psi_i\rangle}_\bW\,\psi_i\Big|_\bW^2
\end{align*}
which can be interpreted as a quadrature rule for the integral
\begin{align*}
    \int_0^T \Big|y(t;\bnu^k)-\sum_{i=1}^\ell{\langle y(t;\bnu^k),\psi_i\rangle}_\bW\,\psi_i\Big|_\bW^2\,\mathrm dt,
\end{align*}
provided that $\{\alpha_j^n\}_{j=1}^n$ are chosen, e.g., as trapezoidal weights
\begin{equation}
    \label{SIAM:Eq-I.1.1.90}
    \alpha_1^n=\frac{\delta t_1}{2},\quad\alpha_j=\frac{\delta t_{j-1}+\delta t_j}{2}\text{ for }j=2,\ldots,n-1,\quad\alpha_n=\frac{\delta t_m}{2}.
\end{equation} 
It is worth to mention that a choice of trapezoidal weights leads to the same order of the approximation error as the discretization of the first and second derivatives by second-order centered difference quotients in \eqref{SIAM:Eq_CentrDiff}. This issue will be addressed in a more detail in Section~\ref{Section:Perturbation analysis for the POD basis}. Of course, other choices are also possible, e.g., Simpson weights. Now, \eqref{SIAM:Eq-I.1.1.2} reads
\begin{equation}
    \label{SIAM:EqPODFinDynSys}
    \left\{
    \begin{aligned}
        &\min\sum_{k=1}^K \omega_k^{K}\sum_{j=1}^n\alpha_j^n\, \Big|y(t_j;\bnu^k)-\sum_{i=1}^\ell{\langle y(t_j;\bnu^k),\psi_i\rangle}_\bW\,\psi_i\Big|_\bW^2\\
        &\hspace{0.5mm}\text{s.t. } \{\psi_i\}_{i=1}^\ell\subset\mathbb C^m\text{ and }{\langle\psi_i,\psi_j\rangle}_\bW=\delta_{ij},~1 \le i,j \le \ell.
    \end{aligned}
    \right.
\end{equation}
Next we investigate \eqref{SIAM:EqPODFinDynSys}. As in Section~\ref{SIAM:Section-2.1.1.3} let $\bD=\mathrm{diag}\,(\alpha_1^n,\ldots,\alpha_n^n)\in\mathbb R^{n\times n}$. We introduce the matrices $\bY^k=[y_1^k\,|\ldots\,|\,y_n^{K}]\in\mathbb C^{m\times n}$ for $k=1,\ldots,{K}$ and set $\bY=[\bY^1\,|\ldots\,|\,\bY^k]\in\mathbb C^{m\times (n{K})}$ with $\mathrm{rank}\,\bY=d^n\le\min(m,n{K})$. Then we find
\begin{align*}
    \mathcal R^n\psi=\sum_{k=1}^{K}\omega_k^{K}\sum_{j=1}^n\alpha_j^n\,{\langle\psi,y_j^k\rangle}_\bW\,y_j^k=\sum_{k=1}^{K} \omega_k^{K} \bY^k\bD(Y^k)^\mathsf H\bW\psi=\bY \tilde{\bD}\bY^\mathsf H\bW\psi
\end{align*}
for $\psi\in\mathbb C^m$ with the diagonal block matrix $\tilde{\bD}=\mathrm{diag}(\omega_1^{K} \bD,\ldots,\omega_{K}^{K} \bD)\in\mathbb R^{(n{K})\times(n{K})}$. Hence, \eqref{SIAM:Eq-I.1.1.18} corresponds to the eigenvalue problem
\begin{equation}
    \label{SIAM:Eq-I.1.1.81}
    \bY\tilde{\bD}Y^\mathsf H \bW\hat\psi_i^n=\hat\lambda_i^n\hat\psi_i^n,\quad\hat\lambda_1^n\ge\ldots\ge\hat\lambda_{d^n}^n>\hat\lambda_{d^n+1}^n=\ldots =\hat\lambda_m^n=0.
\end{equation}
Setting $\psi_i^n=\bW^{1/2} \hat\psi_i^n$ in \eqref{SIAM:Eq-I.1.1.81} and multiplying by $\bW^{1/2}$ from the left yield
\begin{equation}
    \label{SIAM:Eq-I.1.1.82}
    \bW^{1/2} \bY \tilde{\bD} \bY^\mathsf H W^{1/2}\psi_i^n=\hat\lambda_i^n\psi_i^n.
\end{equation}
Let $\tilde{\bY}=\bW^{1/2}\bY \tilde{\bD}^{1/2} \in \mathbb C^{m \times (n{K})}$. Using $\bW^\mathsf H=\bW$ as well as $\tilde{\bD}^\top=\tilde{\bD}$ we infer from \eqref{SIAM:Eq-I.1.1.82} that the POD basis $\{\hat\psi_i^n\}_{i=1}^\ell$ of rank $\ell$ is given by the symmetric $m \times m$ eigenvalue problem
\begin{align*}
    \tilde{\bY} \tilde{\bY}^\mathsf H \psi_i^n=\hat\lambda_i^n\psi_i^n \quad \text{for } 1 \le i \le \ell, \quad\text{and} \quad {\langle\psi_i^n,\psi_j^n \rangle}_{\mathbb C^m}=\delta_{ij} \quad\text{for } 1 \le i,j \le \ell,
\end{align*}
and $\hat\psi_i^n=\bW^{-1/2}\psi_i^n$. Note that
\begin{align*}
    \tilde{\bY}^\mathsf H \tilde{\bY}=\tilde{\bD}^{1/2}\bY^\mathsf H\bW\bY \tilde{\bD}^{1/2} \in \mathbb C^{(n{K})\times (n{K})}.
\end{align*}
Thus, the POD basis of rank $\ell$ can also be computed by the methods of snapshots as follows: First solve the symmetric $(n{K})\times (n{K})$ eigenvalue problem
\begin{align*}
    \tilde{\bY}^\mathsf H \tilde{\bY}\phi_i^n=\hat\lambda_i^n\phi_i^n \quad \text{for }1 \le i \le \ell \quad \mbox{and}\quad {\langle\phi_i^n,\phi_j^n \rangle}_2=\delta_{ij} \quad \text{for }1 \le i,j \le \ell.
\end{align*}
Then we set (by SVD)
\begin{align*}
    \tilde\psi_i^n=\bW^{-1/2}\psi_i^n=\frac{1}{\hat\sigma_i^n}\, \bW^{-1/2} \tilde{\bY}\phi_i^n=\frac{1}{\hat\sigma_i^n}\,\bY \tilde{\bD}^{1/2}\phi_i^n \quad \text{for }1 \le i \le \ell.
\end{align*}
In real applications the solution of \eqref{SIAM:Eq-I.1.1.87} has to be approximated by utilizing a time integration method. Here, we apply the {\em Crank-Nicolson scheme}\index{Method!Crank-Nicolson} \cite[Chapter~6.3]{Str04} and compute approximations $y_j=y_j(\bnu)\in\mathbb C^m$, $j=1,\ldots,n$, of the solution to \eqref{SIAM:Eq-I.1.1.87} at the time instances $\{t_j\}_{j=1}^n$ for $\bnu\in\mathscr N$:
\begin{equation}
    \label{SIAM:Eq-thetaScheme}
    \begin{aligned}
        &y_1=y_\circ,\\
        &\text{\bf for }j=2\text{ \bf to }n\\
        &\quad\frac{y_j-y_{j-1}}{\delta t_j}=\frac{1}{2}\,\Big(A(\bnu)\big(y_{j-1}+y_j\big)+f(t_{j-1};\bnu)+f(t_j;\bnu)\Big).
    \end{aligned}
\end{equation}
From \eqref{SIAM:Eq-thetaScheme} we compute approximations
\begin{align*}
    \mathbb C^m\ni y_j^k=y_j(\bnu_k)\approx y(t_j;\bnu_k)\in\mathbb C^m,
\end{align*}
which serve as our snapshots in the computation of the POD basis; see Algorithm~\ref{SIAM:Algorithm-I.1.1.4}. Notice that we simplify our presentation by choosing the same temporal grid for every parameter $\bnu$. In practise, we usually have to choose different, in particular adaptive time grids. Let us refer to the recent work \cite{Gra19,GHV21}.

\bigskip
\hrule
\vspace{-3.5mm}
\begin{algorithm}[(Discrete POD for dynamical systems)]
    \label{SIAM:Algorithm-I.1.1.4}
    \vspace{-3mm}
    \hrule
    \vspace{0.5mm}
    \begin{algorithmic}[1]
        \REQUIRE Positive weights $\{\omega_k^{K}\}_{k=1}^{K}$ and $\{\alpha_j^n\}_{j=1}^n$, POD rank $\ell \le \min(m,n{K})$, Hermitian positive-definite matrix $\bW \in \mathbb C^{m \times m}$, parameter grid $\{\bnu_k\}_{k=1}^{K}\subset\mathscr N$ and {\tt flag} for the eigenvalue solver;
        \FOR{$k=1$ {\bf to} ${K}$}
            \STATE Solve \eqref{SIAM:Eq-thetaScheme} for $y_j^k=y_j(\bnu_k)\in\mathbb C^m$, $j=1,\ldots,n$;
            \STATE Set $\bY^k=[y_1^k\,|\ldots|\,y_n^k]\in\mathbb C^{m\times n}$;
        \ENDFOR
        \STATE Define $\bY=[\bY^1\,|\ldots|\,\bY^{K}] \in \mathbb C^{m \times (n{K})}$;
        \STATE Set $\bD=\mathrm{diag}\,(\alpha_1^n,\ldots,\alpha_n^n)\in\mathbb R^{n\times n}$ and $\tilde{\bD}=\mathrm{diag}\,(\omega_1^{K}\bD,\ldots,\omega_{K}^{K}\bD)\in\mathbb R^{(n{K})\times(n{K})}$;
        \IF{{\tt flag} = 0}
            \STATE Determine $\tilde{\bY}=\bW^{1/2}\bY\tilde{\bD}^{1/2} \in \mathbb C^{m \times (n{K})}$;
            \STATE Compute singular value decomposition $[\bPsi,\bSigma,\bPhi]=\mathrm{svd} \,(\tilde{\bY})$;
            \STATE Set $\hat\psi_i^n=\bW^{-1/2}\bPsi_{\cdot,i}\in \mathbb C^m$ and $\hat\lambda_i^n=\Sigma_{ii}^2$ for $i=1,\ldots,\ell$;
            \STATE Compute $\mathcal E(\ell)=\sum_{i=1}^\ell \lambda_i/{\|\tilde{\bY}\|}_F^2$;
        \ELSIF{{\tt flag} = 1}
            \STATE Determine $\bR=\bW^{1/2}\bY\tilde{\bD}\bY^\mathsf H\bW^{1/2} \in \mathbb C^{m \times m}$;
            \STATE Compute eigenvalue decomposition $[\bPsi,\bLambda]=\mathrm{eig}\, (\bR)$;
            \STATE Set $\hat\psi_i^n=\bW^{-1/2}\bPsi_{\cdot,i}\in \mathbb C^m$ and $\hat\lambda_i^n=\Lambda_{ii}$ for $i=1,\ldots,\ell$;
            \STATE Compute $\mathcal E(\ell)=\sum_{i=1}^\ell \lambda_i/{\|\tilde{\bY}\|}_F^2$;
        \ELSIF{{\tt flag} = 2}
            \STATE Determine $\bK=\tilde{\bD}^{1/2}\bY^\mathsf H\bW\bY\tilde{\bD}^{1/2} \in \mathbb C^{(n{K})\times (n{K})}$;
            \STATE Compute eigenvalue decomposition $[\bPhi,\bLambda]=\mathrm{eig}\, (\bK)$;
            \STATE Set $\hat\lambda_i^n=\Lambda_{ii}$ and $\hat\psi_i^n=\Y\tilde{\bD}^{1/2}\bPhi_{\cdot,i}/\Lambda_{ii}^{1/2}\in \mathbb C^m$ for $i=1,\ldots,\ell$;
            \STATE Compute $\mathcal E(\ell)=\sum_{i=1}^\ell \lambda_i/{\|\tilde{\bY}\|}_F^2$;
        \ENDIF
        \RETURN POD basis $\{\hat\psi_i^n\}_{i=1}^\ell$, eigenvalues $\{\hat\lambda_i^n\}_{i=1}^\ell$ and ratio $\mathcal E(\ell)$;
    \end{algorithmic}
    \hrule
\end{algorithm}

\section{The continuous variant of the POD method}
\label{SIAM:Section-2.1.2}
\setcounter{equation}{0}
\setcounter{theorem}{0}
\setcounter{figure}{0}
\setcounter{run}{0}

In Section \ref{SIAM:Section-2.1.1} we have considered the case of a finite number of snapshots $y_1^k,\ldots,y_n^k\in X$ for $1\le k\le{K}$, leading to the discrete variant of the POD method. As we have seen in Section~\ref{SIAM:Section-2.1.1.4} the discrete set of snapshots is in many cases an approximation of a snaphot continuum and the weights $\{\alpha_j^n\}_{j=1}^n$ are given as quadrature weights to approximate an integral inner product in this case. \\
So let us suppose that we are given a parameter set $\mathscr D$, which is a closed and bounded subset of $\mathbb R^\mathfrak p$ with $\mathfrak p\in\mathbb N$, and trajectories $y^k\in C(\mathscr D;X)$, $1\le k\le{K}$. So far, for the discrete variant of the POD method, we have assumed a finite parameter grid
\begin{align*}
    \mathscr D^n= \{ \bmu_1,\ldots,\bmu_n \} \subset \mathscr D
\end{align*}
and an approximation of the snapshots by $y_j^k=y^k(\bmu_j)\in X$ or $y_j^k\approx y^k(\bmu_j)\in X$. Then the snapshot subspace $\mathscr V^n$ of dimension $d^n$ introduced in \eqref{SIAM:Eq-I.1.1.1} depends on the chosen grid points $\{\bmu_j\}_{j=1}^n$. Consequently, the operator $\mathcal R^n$ defined in \eqref{SIAM:Eq-I.1.1.12} and thus the POD basis $\{\hat\psi_i^n\}_{i=1}^\ell$ of rank $\ell$ as well as the corresponding eigenvalues $\{\hat\lambda_i^n\}_{i=1}^\ell$ depend also on the grid $\mathscr D^n$ (which has already been indicated by the superindex $n$). In this sense, enlarging the number $n$ of snapshots $\{y_j^k\}_{j=1}^n$ and weights $\{\alpha_j^n\}_{j=1}^n$ corresponds to a finer approximation of the continuous snapshots and the respective integrals. \\
We now study the POD method when taking the continuous snapshots $y^k\in C(\mathscr D;X)$, $1\le k\le{K}$, instead of the discretized ones, leading to the continuous variant of the POD method. In detail, the discrete trajectory set $\mathcal V^n$ from \eqref{SIAM:Eq-I.1.1.1} is now replaced by a continuous $X$-valued set. \\ 
In this section we show that the theory from the discrete variant can be adapted to the continuous case and that the POD basis can be obtained analogously to the discrete case. \\
Later in Section 2.4 we will then investigate how the grid points $\{\bmu_j\}_{j=1}^n$ and the weights $\{\alpha_j^n\}_{j=1}^n$ have to be chosen to obtain the convergence of the discrete POD basis $\{\hat\psi_i^n\}_{i=1}^\ell$ to the continuous POD basis $\{\hat\psi_i\}_{i=1}^\ell$ as $n \to \infty$, i.e., when the parameter grid is being increasingly refined.

\subsection{Problem formulation and main result}
\label{SIAM:Section-2.1.2.1}

As mentioned above, we suppose that $\mathscr D$ is a closed and bounded subset of $\mathbb R^\mathfrak p$ with $\mathfrak p\in\mathbb N$. In Section~\ref{SIAM:Section-2.1.2.2} we will study the special case $\mathscr D=[0,T]$ in more detail. Suppose that we have trajectories or \index{POD method!continuous variant!snapshots}{\em snapshots} $y^k\in L^2(\mathscr D;X)$, $1\le k\le{K}$. Let us now introduce the \index{POD method!continuous variant!snapshot space}{\em snapshot space} by 
\begin{align*}
    \mathscr V=\Bigg\{\sum_{k=1}^K\int_\mathscr Dy^k(\bmu)\phi^k(\bmu)\,\mathrm d\bmu\,\Big|\,\phi^k\in L^2(\mathscr D)\text{ for } 1\le k\le K\Bigg\}\subset X
\end{align*}
with the dimension
\begin{align*}
    d=\left\{
    \begin{aligned}
        &\dim\mathscr V&&\text{if }\dim\mathscr V\text{ is finite},\\
        &\infty&&\text{otherwise}.
    \end{aligned}
    \right.
\end{align*}
Instead of Problem~\ref{ProblemPOD} we consider the following problem in this section:

\begin{problem}[Continuous POD method]
    \label{ProblemPODCont}
    For any $\ell\le d$, the \index{POD method!continuous variant}{\em continuous method of POD} consists in determining a POD basis of rank $\ell$ which minimizes the mean square error between the trajectories $y^k$ and the corresponding $\ell$-th partial Fourier sums on average in the closed, bounded set $\mathscr D$:
    \begin{equation}
        \tag{$\mathbf P^\ell$}
        \label{SIAM:Eq-I.1.2.3}
        \left\{
        \begin{aligned}
            &\min\sum_{k=1}^{K} \omega_k^{K}\int_\mathscr D \Big\| y^k(\bmu)-\sum_{i=1}^\ell {\langle y^k(\bmu),\psi_i\rangle}_X\,\psi_i\Big\|_X^2\,\mathrm d\bmu\\
            &\hspace{0.5mm}\text{s.t. } \{\psi_i\}_{i=1}^\ell\subset X\text{ and }{\langle\psi_i,\psi_j\rangle}_X=\delta_{ij},~1 \le i,j \le \ell.
        \end{aligned}
        \right.
    \end{equation}
\end{problem}

\begin{definition}
    An optimal solution $\{\hat\psi_i\}_{i=1}^\ell$ to Problem~\emph{\ref{ProblemPODCont}} is called a \index{POD method!continuous variant!POD basis of rank $\ell$}{\em POD basis of rank $\ell$}\index{POD method!POD basis}.
\end{definition}

\begin{remark}
    \rm Analogous to \eqref{SIAM:Eq-I.1.1.4} we can -- instead of \eqref{SIAM:Eq-I.1.2.3} -- consider the problem
    \begin{equation}
        \tag{$\mathbf{\hat P}^\ell$}
        \label{SIAM:Eq-I.1.2.4}
        \left\{
        \begin{aligned}
            &\max\sum_{k=1}^{K}\omega_k^{K}\int_\mathscr D \sum_{i=1}^\ell \big|{\langle y^k(\bmu),\psi_i\rangle}_X\big|^2\,\mathrm d\bmu\\
            &\hspace{0.5mm}\text{s.t. } \{\psi_i\}_{i=1}^\ell\subset X\text{ and }{\langle\psi_i,\psi_j\rangle}_X=\delta_{ij},~1 \le i,j \le \ell;
        \end{aligned}
        \right.
    \end{equation}
    cf. Remark~\ref{Remark:AltCost-1}. A solution to \eqref{SIAM:Eq-I.1.2.3} solves also \eqref{SIAM:Eq-I.1.2.4} and vice versa.\hfill$\blacksquare$
\end{remark}

\noindent
{\bf Main result:} In Theorem~\ref{Theorem2.2.1} we will show that for every $\ell\in\{1,\ldots,d\}$ a solution $\{\hat\psi_i\}_{i=1}^\ell$ to \eqref{SIAM:Eq-I.1.2.3} and to \eqref{SIAM:Eq-I.1.2.4} is characterized by the eigenvalue problem
\begin{equation}
    \label{EigContPODMeth}
    \mathcal R\hat\psi_i=\hat\lambda_i\hat\psi_i,\quad1\le i\le\ell,
\end{equation}
where $\hat\lambda_1\ge\ldots\ge\hat\lambda_d>0$ and the linear integral operator $\mathcal R:X\to X$ is given as
\begin{equation}
    \label{SIAM:Eq-I.1.2.5}
    \mathcal R\psi=\sum_{k=1}^{K}\omega_k^{K}\int_\mathscr D {\langle \psi,y^k(\bmu)\rangle}_X\,y^k(\bmu)\,\mathrm d\bmu\quad\text{for }\psi\in X.
\end{equation}

\subsection{Derivation of the main result}
\label{SIAM:Section-2.1.2.1a}

To proof our main result we define the linear operator $\mathcal Y:L^2(\mathscr D;\mathbb C^{K}) \to X$ by
\begin{equation}
    \label{Eq2.2.6}
    \mathcal Y \phi=\sum_{k=1}^{K}\omega_k^{K}\int_\mathscr D\phi^k(\bmu)y^k(\bmu)\,\mathrm d\bmu\quad\text{for }\phi=(\phi^1,\ldots,\phi^{K})\in L^2(\mathscr D;\mathbb C^K).
\end{equation}
Notice that $\mathrm{ran}\,(\mathcal Y)=\mathscr V$ holds true. The following result is proved in Section~\ref{SIAM:Section-2.6.2}.

\begin{lemma}
    \label{Lemma2.2.0}
    Let $y^k\in L^2(\mathscr D;X)$, $1\le k\le K$, be given snapshot trajectories.
    \begin{enumerate}
        \item [\em 1)] The operator $\mathcal Y$ is bounded.
        \item [\em 2)] The adjoint $\mathcal Y^\star:X\to L^2(\mathscr D;X)$ is given by
        \begin{equation}
            \label{Eq2.2.10}
            (\mathcal Y^\star \psi)(\bmu)=\left(
            \begin{array}{c}
                {\langle\psi,y^1(\bmu)\rangle}_X\\
                \vdots\\
                {\langle\psi,y^{K}(\bmu)\rangle}_X\\
            \end{array}
            \right)\quad\text{for }\psi\in X\text{ and }\bmu\in\mathscr D\text{ a.e.,}
        \end{equation}
        where \index{a.e., almost everywhere}{\em a.e.} stands for {\em almost everywhere} and the (weighted) inner product in $L^2(\mathscr D;\mathbb C^K)$ is given as
        \begin{align*}
            {\langle \phi,\tilde\phi\rangle}_{L^2(\mathscr D;\mathbb C^K)}=\sum_{k=1}^{K}\omega_k^{K}\,\int_\mathscr D\phi^k(\bmu)\tilde\phi^k(\bmu)\,\mathrm d\bmu\quad\text{for }\phi,\tilde\phi\in L^2(\mathscr D;\mathbb C^K).
        \end{align*}
        Furthermore, $\mathcal Y^\star$ is compact.
    \end{enumerate}
\end{lemma}

From Lemma~\ref{Lemma2.2.0} we infer the next result, which is proved in Section~\ref{SIAM:Section-2.6.2}.

\begin{corollary}
    \label{Lemma2.2.1}
    Let $y^k\in L^2(\mathscr D;X)$, $1\le k\le K$, be given snapshot trajectories.
    \begin{enumerate}
        \item [\em 1)] The operator $\mathcal R$ introduced in \eqref{SIAM:Eq-I.1.2.5} is is given as $\mathcal Y\mathcal Y^\star$. Furthermore, $\mathcal R$ is compact, self-adjoint and non-negative.
        \item [\em 2)] For every $\bmu\in\mathscr D$ the linear operator $\mathcal K=\mathcal Y^\star\mathcal Y:L^2(\mathscr D;\mathbb C^{K}) \to L^2(\mathscr D;\mathbb C^{K})$ is given as
        \begin{align*}
            \big(\mathcal K\phi\big)(\bmu) =\left(
            \begin{array}{c}
                \sum\limits_{k=1}^{K} \int_\mathscr D {\langle y^k(\bnu),y^1(\bmu)\rangle}_X\phi^k(\bnu)\,\mathrm d\bnu\\
                \vdots\\
                \sum\limits_{k=1}^{K} \int_\mathscr D {\langle y^k(\bnu),y^{K}(\bmu)\rangle}_X\phi^k(\bnu)\,\mathrm d\bnu
            \end{array}
            \right) \quad \text{for } \phi\in L^2(\mathscr D;\mathbb C^{K}).
        \end{align*}
        Moreover, $\mathcal K$ is compact.
    \end{enumerate}
\end{corollary}

\begin{remark}
    \rm In \eqref{Eq:OpY} we have introduced the linear and compact operator $\mathcal Y^n:\mathbb C^{{K}\times n}\to X$, which can be interpreted as an approximation of the operator $\mathcal Y$ defined in \eqref{Eq2.2.6} provided
    \begin{align*}
        \int_\mathscr D\phi^k(\bmu)y^k(\bmu)\,\mathrm d\bmu\approx\sum_{j=1}^n\alpha_j^n\Phi_{kj}y^k_j\quad\text{for }1\le k\le {K}.
    \end{align*}
    Since $\Phi^k(\bmu_j)$ as well as $y^k(\bmu_j)$ might not be well-defined for a given $\Phi^k\in L^2(\mathscr D;\mathbb C)$ and $y^k\in L^2(\mathscr D;X)$, we set
    \begin{align*}
        \Phi_{kj}=\frac{1}{|\mathscr M_j|}\int_{\mathscr M_j}\phi^k(\bmu)\,\mathrm d\bmu\quad\text{and}\quad y^k_j=\frac{1}{|\mathscr M_j|}\int_{\mathscr M_j}y^k(\bmu)\,\mathrm d\bmu,
    \end{align*}
    where $\mathscr M_1,\ldots,\mathscr{M}_n \subset\mathscr D$ are open sets with non-zero measure $|\mathscr M_j|$ for all $j=1,\ldots,n$, such that $\cup_{j=1}^n \mathscr M_j = \mathscr D$ and $\bmu_j\in \mathscr M_j$. We will address this issue in Section~\ref{Section:Perturbation analysis for the POD basis} in more detail.\hfill$\blacksquare$
\end{remark}

Recall the definition of the index set $\mathbb I$ in Section~\ref{SIAM:Section-2.1.1.1}. In the next theorem we formulate how the solution to \eqref{SIAM:Eq-I.1.2.3} and \eqref{SIAM:Eq-I.1.2.4} can be found. The proof is given in Section~\ref{SIAM:Section-2.6.2}.

\begin{theorem}
    \label{Theorem2.2.1}
    Let $y^k\in L^2(\mathscr D;X)$ be given trajectories for $1\le k\le{K}$. Suppose that the linear operator $\mathcal R$ is defined by \eqref{SIAM:Eq-I.1.2.5}. Then there exist non-negative eigenvalues $\{\hat\lambda_i\}_{i\in\mathbb I}$ and associated orthonormal eigenfunctions $\{\hat\psi_i\}_{i\in\mathbb I}$ satisfying
    \begin{equation}
        \label{Eq2.2.16}
        \mathcal R\hat\psi_i=\hat\lambda_i\hat\psi_i,\hspace{2mm}\left\{
        \begin{aligned}
            &\hat\lambda_1\ge\ldots\ge\hat\lambda_d>\hat\lambda_{d+1}=\ldots =0&&\text{if }d<\infty,\\
            &\hat\lambda_1\ge\hat\lambda_2\ge\ldots\text{ and }\lim_{i\to\infty}\hat\lambda_i=0&&\text{if }d=\infty.
        \end{aligned}
        \right.
    \end{equation}
    For every $\ell\in\mathbb N$ (with $\ell\le d$ if $d$ is finite) the first $\ell$ eigenfunctions $\{\hat\psi_i\}_{i=1}^\ell$ solve \eqref{SIAM:Eq-I.1.2.3} and \eqref{SIAM:Eq-I.1.2.4}. Moreover, the value of the objectives evaluated at the optimal solution $\{\hat\psi_i\}_{i=1}^\ell$ satisfies
    \begin{equation}
        \label{Eq2.2.17}
        \sum_{k=1}^{K}\omega_k^{K}\int_\mathscr D\Big\|y^k(\bmu)-\sum_{i=1}^\ell{\langle y^k(\bmu),\hat\psi_i\rangle}_X\,\hat\psi_i\Big\|_X^2\,\mathrm d\bmu=\sum_{i>\ell}\hat\lambda_i
    \end{equation}
    and
    \begin{equation}
        \label{Eq2.2.18}
        \sum_{k=1}^{K}\omega_k^{K}\int_\mathscr D\sum_{i=1}^\ell\big|{\langle y^k(\bmu),\hat\psi_i\rangle}_X\big|^2\,\mathrm d\bmu=\sum_{i=1}^\ell\hat\lambda_i,
    \end{equation}
    respectively.
\end{theorem}

We obtain also a result similar to Lemma~\ref{Lem:PropRnOp}-1). Its proof is given in Section~\ref{SIAM:Section-2.6.2}.

\begin{lemma}
    \label{Lem:PropRnOp-b}
    Let $y^k\in L^2(\mathscr D;X)$ be given trajectories for $1\le k\le{K}$. Suppose that the linear operator $\mathcal R$ is defined by \eqref{SIAM:Eq-I.1.2.5}. Moreover, let $\{\hat\lambda_i\}_{i\in\mathbb I}$ and $\{\hat\psi_i\}_{i\in\mathbb I}$ be the non-negative eigenvalues and associated orthonormal eigenfunctions, respectively, satisfying \eqref{Eq2.2.16}.
    \begin{enumerate}
        \item [\em 1)] Then we have
        \begin{equation}
            \label{SIAM:Eq-I.1.1.19-b}
            \sum_{k=1}^{K}\omega_k^{K}\int_\mathscr D\big|{\langle y^k(\bmu),\hat\psi_i\rangle}_X\big|^2\,\mathrm d\bmu={\langle\mathcal R\hat\psi_i,\hat\psi_i\rangle}_X=\hat\lambda_i\quad\text{for any }i\in\mathbb I.
        \end{equation}
        In particular, if $d<\infty$ holds it follows that
        \begin{equation}
            \label{SIAM:Eq-I.1.1.20-2}
            \sum_{k=1}^{K}\omega_k^{K}\int_\mathscr D\big|{\langle y^k(\bmu),\hat\psi_i\rangle}_X\big|^2\,\mathrm d\bmu=0\quad\text{for all } i>d.
        \end{equation}
        \item [\em 2)] The identity
        \begin{equation}
            \label{Eq2.2.19}
            \sum_{k=1}^{K}\int_\mathscr D{\|y^k(\bmu)\|}_X^2\,\mathrm d\bmu=\left\{
            \begin{aligned}
                &\sum_{i=1}^d\hat\lambda_i&&\text{if }\mathrm{rank}(\mathcal R)=d<\infty,\\
                &\sum_{i=1}^\infty\hat\lambda_i&&\text{otherwise}
            \end{aligned}
            \right.
        \end{equation}
        holds.
    \end{enumerate}
\end{lemma}

\begin{remark}[Singular value decomposition]
    \label{Remark2.2.1a}
    \rm Suppose that $y^k\in L^2(\mathscr D;X)$ holds. We have introduced the linear operator $\mathcal Y:L^2(\mathscr D;\mathbb R^{K})\to X$ in \eqref{Eq2.2.6}. Recall that the linear operators $\mathcal R=\mathcal Y\mathcal Y^*$ and $\mathcal K=\mathcal Y^*\mathcal Y$ are compact, non-negative and self-adjoint. Due to Theorem~\ref{Theorem2.2.1} there exist non-negative eigenvalues $\{\hat\lambda_i\}_{i\in\mathbb I}$ and associated orthonormal eigenfunctions $\{\hat\psi_i\}_{i\in\mathbb I}$ satisfying \eqref{Eq2.2.16}. Furthermore, there is a sequence $\{\hat\phi_i\}_{i\in\mathbb I}$ of orthonormal eigenfunctions such that
    \begin{align*}
        \mathcal K\hat\phi_i=\hat\lambda_i\hat\phi_i \quad\text{for } i \in \mathbb I.
    \end{align*}
    We set $\mathbb R^+_0=\{s\in\mathbb R\,|\,s\ge 0\}$ and $\hat\sigma_i=\hat\lambda_i^{1/2}$ for $i \in \mathbb I$. The sequence $\{\hat\sigma_i,\hat\phi_i,\hat\psi_i\}_{i\in\mathbb I}$ in $\mathbb R^+_0\times L^2(\mathscr D;\mathbb R^{K})\times X$ can be interpreted as a {\em singular value decomposition}\index{POD method!continuous variant!SVD} of the mapping $\mathcal Y$. In fact, we have
    \begin{align*}
        \mathcal Y\hat\phi_i=\hat\sigma_i\hat\psi_i\quad \text{and} \quad\mathcal Y^\star\hat\psi_i=\hat\sigma_i\hat\phi_i\quad \text{for } i\in\mathbb I.
    \end{align*}
    Since $\hat\sigma_i>0$ holds for $1=1\ldots,d$, we have $\hat\phi_i=\mathcal Y^\star\hat\psi_i/\hat\sigma_i$ for $i=1,\ldots,d$. For the specific case $X=\mathbb C^m$, $K=1$, $\omega_1^K=1$ and one trajectory $x\in L^2(\mathscr D;\mathbb C^m)$ we find
    \begin{align*}
        \hat\psi_i=\frac{1}{\hat\sigma_i}\,\mathcal Y\hat\phi_i=\int_\mathscr D\bigg(\frac{1}{\hat\sigma_i}\,\hat\phi_i(\bmu)\bigg)y(\bmu)\,\mathrm d\bmu\quad\text{for }i=1,\ldots,\ell,
    \end{align*}
    which can be interpreted as the continuous variant of the \index{POD method!continuous variant!methods of snapshots}{\em method of snapshots}; cf. Remark~\ref{Remark:Sirovich}.\hfill$\blacksquare$
\end{remark}

\subsection{POD for finite-dimensional dynamical systems}
\label{SIAM:Section-2.1.2.2}

In this section we continue the discussion of Section~\ref{SIAM:Section-2.1.1.4} by introducing a continuous version of the POD method for the given initial-value problem \eqref{SIAM:Eq-I.1.1.87}. Let $X=\mathbb C^m$ be endowed with the weighted inner product \eqref{SIAM:Eq-I.1.1.61}. In the context of Section~\ref{SIAM:Section-2.1.2.1} we choose $\mathscr D=[0,T]$. We write $t$ instead of $\bmu$. Then \eqref{SIAM:Eq-I.1.2.3} has the form
\begin{equation}
    \label{SIAM:EqPODFinDynSysCont}
    \left\{
    \begin{aligned}
        &\min \sum_{k=1}^K\omega_k^K\int_0^T\Big\|y(t;\bnu_k)-\sum_{i=1}^\ell{\langle y(t;\bnu_k),\psi_i\rangle}_\bW\,\psi_i\Big\|_\bW^2\,\mathrm dt\\
        &\hspace{0.5mm}\text{s.t. } \{\psi_i\}_{i=1}^\ell\subset\mathbb C^m\text{ and }{\langle\psi_i,\psi_j\rangle}_\bW=\delta_{ij},~1 \le i,j \le \ell.
    \end{aligned}
    \right.
\end{equation}
We introduce the matrix $\bR=((R_{ij}))\in\mathbb C^{m\times m}$ with the elements
\begin{align*}
    R_{ij}=\sum_{k=1}^{K} \omega_k^{K}\int_0^T \mathrm y^k_i(t;\bnu_k)\overline{\mathrm y^k_j(t;\bnu_k)}\,\mathrm dt \quad\text{for }1\le i,j\le m,
\end{align*}
where $\mathrm y^k_i(t;\bnu_k)$ stands for the $i$-th component of the vector $y^k(t;\bnu_k)\in\mathbb C^m$. From
\begin{align*}
    \overline R_{ij}=\sum_{k=1}^{K} \omega_k^{K}\int_0^T \overline{\mathrm y^k_i(t;\bnu_k)}\mathrm y^k_j(t;\bnu_k)\,\mathrm dt=R_{ji} \quad \text{for } 1\le i,j\le m,
\end{align*}
it follows that $\bR^\mathsf H=\bR$ holds true. Let $\psi\in\mathbb C^m$ be arbitrarily chosen and $i\in\{1,\ldots,\ell\}$. Then we have
\begin{align*}
    \big(\mathcal R\psi\big)_i&=\sum_{k=1}^{K} \omega_k^{K}\int_0^T{\langle \psi,y(t;\bnu_k)\rangle}_\bW\mathrm y^k_i(t;\bnu_k)\,\mathrm dt\\
    &=\sum_{k=1}^{K} \omega_k^{K}\sum_{j=1}^m\sum_{l=1}^m\psi_lW_{lj}\int_0^T \mathrm y^k_i(t;\bnu_k)\overline{\mathrm y_j(t;\bnu_k)}\,\mathrm dt=\sum_{j=1}^m\sum_{l=1}^mR_{ij}W_{jl}\psi_l\\
    &=\big(\bR\bW\psi\big)_i.
\end{align*}
Consequently, the eigenvalue problem \eqref{Eq2.2.16} reads
\begin{equation}
    \label{Eq2.2.16EigVal}
    \bR\bW\hat\psi_i=\hat\lambda_i\hat\psi_i,\hspace{2mm}\hat\lambda_1\ge\ldots\ge\hat\lambda_d>\hat\lambda_{d+1}=\ldots =0.
\end{equation}
Setting $\psi_i=\bW^{1/2}\hat\psi_i$ and multiplying \eqref{Eq2.2.16EigVal} by $\bW^{1/2}$ from the left we find
\begin{equation}
    \label{Eq2.2.16EigValSym}
    \tilde{\bR}\psi_i=\hat\lambda_i\psi_i,\hspace{2mm}\hat\lambda_1\ge\ldots\ge\hat\lambda_d>\hat\lambda_{d+1}=\ldots =0,
\end{equation}
with $\tilde{\bR}=\bW^{1/2}\bR\bW^{1/2}$. Since $\tilde{\bR}^\mathsf H=\tilde{\bR}$ is valid, \eqref{Eq2.2.16EigValSym} is a Hermitian eigenvalue problem in $\mathbb C^m$. 

\section{Perturbation analysis for the POD basis}
\label{Section:Perturbation analysis for the POD basis}
\setcounter{equation}{0}
\setcounter{theorem}{0}
\setcounter{figure}{0}
\setcounter{run}{0}

So far, we have introduced two variants of the POD method: In Section \ref{SIAM:Section-2.1.1}, the snapshots $\{y_i^k\}_{i=1}^n$ are given over a discrete index set $\{1,...,n\}$ for $1 \le k \le {K}$, which is why the corresponding problem \eqref{SIAM:Eq-I.1.1.2} is also called the \textit{discrete POD method}; cf. Problem~\ref{ProblemPOD}. Meanwhile, a continuous snapshot trajectory $\{y^k(\bmu) ~|~ \bmu \in \mathscr D\}$ with $\mathscr D \subset \mathbb R^{\mathfrak p}$ occurred in Section \ref{SIAM:Section-2.1.2} along with the corresponding \textit{continuous POD method}; cf. Problem~\ref{ProblemPODCont}. It was shown that the solutions to these problems are given by orthonormal eigenvectors $\{\hat\psi_1^n,\hdots,\hat\psi_\ell^n\}$ and $\{\hat\psi_1,\hdots,\hat\psi_\ell\}$ to the $\ell$ largest eigenvalues of the operators $\mathcal R^n, \mathcal R: X \to X$ introduced in \eqref{SIAM:Eq-I.1.1.12} and \eqref{SIAM:Eq-I.1.2.5}, respectively. In this section, we investigate what happens when the discrete snapshots represent an approximation of the continuous trajectory that gets more and more precise as $n \to \infty$. This notion was already introduced at the beginning of Section \ref{SIAM:Section-2.1.2} by discretizing the parameter space $\mathscr D$ into $n$ sample points $\{\bmu_1^n,\hdots,\bmu_n^n \} \subset \mathscr D$ and setting
\begin{align*}
    y_j^k = y^k(\bmu_j) \quad \text{for } 1 \le j \le n\text{ and }1 \le k \le {K}.
\end{align*}
Just like this discrete trajectory represents the continuous one, we will choose the positive weights $\alpha_1^n,\hdots,\alpha_n^n$ in a way that the sum in the discrete problem \eqref{SIAM:Eq-I.1.1.2} provides an approximation for the integral in the continuous problem \eqref{SIAM:Eq-I.1.2.3}.  Under certain conditions, this will result in a form of convergence of the discrete solution vectors $\{\hat\psi_i^n\}_{i=1}^\ell$ to the continuous ones $\{\hat\psi_i\}_{i=1}^\ell$ and of the discrete eigenvalues $\{\hat\lambda_i^n\}_{i=1}^\ell$ to the continuous ones $\{\hat\lambda_i\}_{i=1}^\ell$ as $n$ tends to $\infty$. 

\subsection{General case in complex Hilbert spaces}
\label{Section:PertGenCase}

In general the spectrum $\sigma(\mathcal T)$ of an operator $\mathcal T\in\mathscr L(X)$ will not depend continuously on $\mathcal T$. This is an essential difference to the finite dimensional case, where, e.g., the eigenvalues of a quadratic matrix $A=((A_{ij})) \in \mathbb R^{n \times n}$ depend continuously on the entries $A_{ij}$. However, if the operator contains isolated spectral points, a similar behavior can be proven for those points in the infinite-dimensional case, see Section~\ref{SIAM:Section_SpectralTheory}. 

We will now apply the general perturbation results from Section~\ref{SIAM:Section_SpectralTheory} to the concrete case of the operators $\mathcal R^n$ are $\mathcal R$ for which we will assume that convergence in the operator norm holds.

\begin{theorem}
    \label{Theorem2.3.1}
    Suppose that the sample points $\{\bmu_j^n\}_{j=1}^n$ and weighting parameters $\{\alpha_j^n\}_{j=1}^n$ are chosen in a way such that for every $F \in C(\mathscr D;X)$, it holds
    \begin{equation}
        \label{eq:perturbationIntegration}
        \sum_{j=1}^n \alpha_j^n F(\bmu_j^n) \to \int_{\mathscr D} F(\bmu) ~\mathrm d\bmu \quad  \text{as } n \to \infty.
    \end{equation}
    Further, assume that for $1 \le k \le {K}$, the trajectory is given by $y^k \in C(\mathscr D;X)$ and its discrete counterpart is chosen as $y_j^k := y^k(\bmu_j)$ for $1 \le j \le n$. Let the operators $\mathcal R^n, \mathcal R \in \mathscr L(X)$ defined in \eqref{EigPODPro} and \eqref{SIAM:Eq-I.1.2.5} satisfy
    \begin{equation}
        \label{Eq2.3.1}
        {\|\mathcal R^n-\mathcal R\|}_{\mathscr L(X)} \to 0 \quad \text{as } n \to \infty.
    \end{equation}
    Let $\{(\hat\lambda_i^n,\hat\psi_i^n)\}_{i \in \mathbb I}$ and $\{(\hat\lambda_i,\hat\psi_i)\}_{i \in \mathbb I}$ denote the sequence of decreasing eigenvalues and of eigenvectors of $\mathcal R^n$ and $\mathcal R$ satisfying \eqref{SIAM:Eq-I.1.1.18} and \eqref{Eq2.2.16}, respectively. Suppose that $\ell \in \{1,...,d\}$ is finite. Then we have $\hat\lambda_i^n \to \hat\lambda_i$ as $n \to \infty$ for every $1 \le i \le \ell$. Furthermore, if $\ell$ satisfies $\hat\lambda_\ell > \hat\lambda_{\ell+1}$, we have
    \begin{subequations}
        \begin{align}
            \label{eq:perturbationLambda}
            \lim_{n \to \infty}\sum_{i=\ell+1}^{d^n} \hat\lambda_i^n &= \sum_{i=\ell+1}^d \hat\lambda_i, \\
            \label{eq:perturbationPsi}
            \lim\limits_{n \to \infty} \sum_{i=1}^\ell \langle \varphi,\hat\psi_i^n \rangle_X \hat\psi_i^n &= \sum_{i=1}^\ell \langle \varphi, \hat\psi_i \rangle_X \hat\psi_i \quad\text{for } \varphi \in X.
        \end{align}
    \end{subequations}
\end{theorem}

\begin{example}
    \label{ex:podMeshRefinemenent}
    \rm We want to illustrate the validity of Theorem~\ref{Section:PertGenCase} by continuing our general Example~\ref{ex:generalNumericSetup}. The control value is set to $u(t) = (0,0)^\top$ for all $t \in (0,T)$. The parameter domain $\mathscr D$ coincides with the time axis $(0,T) \subset \mathbb R$. We discretize the infinite-dimensional space $V = H^1(\Omega)$ by the \index{Method!finite element, FE}piecewise linear FE method (cf., e.g., \cite{BS08,Tho97} and Section~\ref{SIAM-Book:Section3.4}) using $m = 1319$ ansatz functions. Therefore, the Hilbert space $X$ is given by $\mathbb R^m$ with the inner product induced by the stiffness matrix $\bS \in \mathbb R^{m \times m}$ from Example~\ref{SIAM:Example-I.1.1.1} approximating the inner product in $V$. In the context of \eqref{SIAM:Eq-I.1.2.3}, we have $K=1$, $\omega_1^1 = 1$ and $y^1 \in L^2(0,T;\mathbb R^m)$ being the weak solution to the FE-discretized PDE \eqref{SIAM:EqMotPDE1} with still-continuous time axis. We ensure the integral approximation \eqref{eq:perturbationIntegration} by choosing the usual equidistant trapezoidal weights given in \eqref{TrapW}. The step size $\Delta t$ for the time discretization is inversely proportional to $n$ and we choose the values 
    \begin{equation}
        \label{eq:podMeshRefinemenent_n}
        n \in \{10,50,100,200,400,1.000,4000,10000\} =: \mathrm N.
    \end{equation}
    For the time integration of the FE-discretized PDE \eqref{SIAM:EqMotPDE1} we apply an implicit Euler scheme with $n+1$ time instances $0=t_0<t_1<\ldots<t_n=T$. The discretization will be described in detail in Section~\emph{\ref{SIAM-Book:Section3.5}}. The discrete solutions are denoted by $\{y^n_j\}_{j=1}^n \subset \mathbb R^m$ for all $n$ in \eqref{eq:podMeshRefinemenent_n}, and are also the data trajectories for the discrete POD problem \eqref{SIAM:Eq-I.1.1.2}. Thus having set up both the discrete and the continuous POD problem, we would like to solve both problems and compare the continuous eigenvalues $\hat\lambda_i$ to the discrete eigenvalues $\hat\lambda_i^n$ for increasing $n$. However, it was already mentioned that the continuous problem can not be solved numerically, so we use the discrete problem with $n=10000$ as an approximation and introduce for the first eight eigenvalues
    \begin{equation}
        \label{eq:podMeshRefinement_mu}
        \mu_i^n:=\big| \hat\lambda_i^n-\hat\lambda_i^{10000} \big| \qquad \text{for } n \in \mathrm N,~i=1,...,8.
    \end{equation}
    \begin{figure}
        \centering
        \includegraphics[width=125mm,height=50mm]{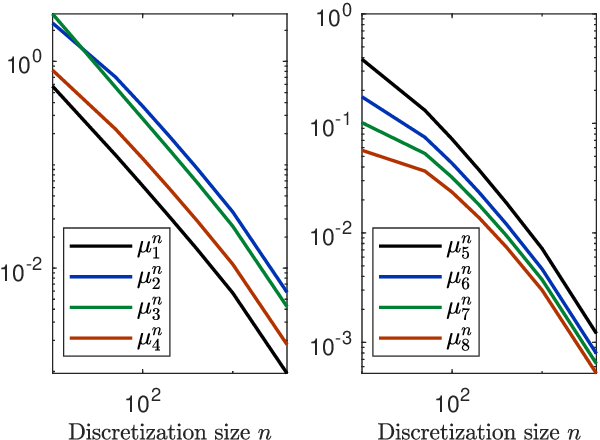}
        \caption{Example~\ref{ex:podMeshRefinemenent}. The values $\bmu_i^n$ from \eqref{eq:podMeshRefinement_mu} indicating the convergence of POD eigenvalues $\hat\lambda_i^n$ as $n \to \infty$.}
        \label{fig:podMeshRefinement_time}
    \end{figure}
    The results are displayed in Figure~\ref{fig:podMeshRefinement_time}, where we can observe a clearly monotonous decay with $n$ in the ``error'' terms $\mu_i^n$, suggesting that these discrete eigenvalues truly approximate the analytical eigenvalues as $n\to\infty$, as it was suggested by Theorem~\ref{Section:PertGenCase}.\hfill$\blacklozenge$
\end{example}

\begin{example}
    \label{ex:podMeshRefinemenent-2}
    \rm We can observe a similar asymptotic convergence as in Example~\ref{ex:podMeshRefinemenent} if we fix the time discretization at $n = 400$ and choose different dimensions $m$ of the FE spaces
    \begin{equation}
        \label{eq:podMeshRefinemenet_m}
        m \in \{139, 203, 263, 359, 494, 765, 1319, 2999 \} =: \mathrm M.
    \end{equation}
    This gives us discrete solutions $\{y_j^m\}_{j=1}^n\subset X^m$, to each of which we solve the discrete POD problem \eqref{SIAM:Eq-I.1.1.2} and get eigenvectors $\hat\lambda_i^m$, this time indicated by the dependence on $m$. Figure~\ref{fig:podMeshRefinement_space} shows the decay of the error in these eigenvalues, which is computed, as in \eqref{eq:podMeshRefinement_mu}, by
    \begin{equation}
        \label{eq:podMeshRefinement_nu}
        \nu_i^m = \big| \hat\lambda_i^m-\hat\lambda_i^{2999} \big| \qquad \text{for } m \in \mathrm M, ~i=1,...,8.
    \end{equation}
    \begin{figure}
        \centering
        \includegraphics[width=125mm,height=50mm]{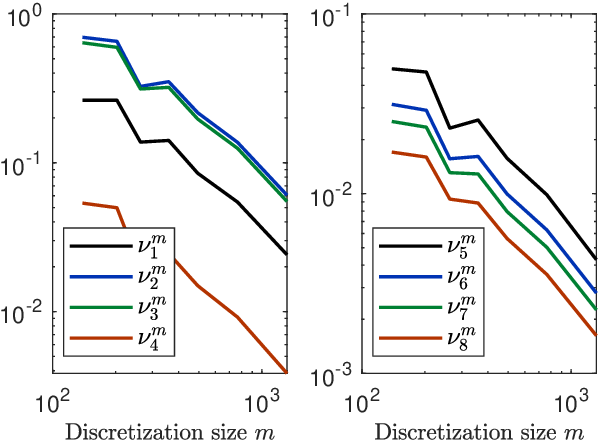}
        \caption{Example~\ref{ex:podMeshRefinemenent-2}. The values $\nu_i^m$ from \eqref{eq:podMeshRefinement_nu} indicating the convergence of POD eigenvalues $\hat\lambda_i^m$ as $m \to \infty$.}
        \label{fig:podMeshRefinement_space}
    \end{figure}
    A similar decay as for the time-refinement case can be seen in Figure~\ref{fig:podMeshRefinement_space}, indicating that by taking the limit $m \to \infty$, we are approaching the eigenvalues for the time-discretized weak solution to the PDE \eqref{SIAM:EqMotPDE1} which is still continuous in the space-dimension. The reason for this lies in the fact that the subspace $X^m$ is isometrically isomorphic to a subspace $V^m$ of $V = H^1(\Omega)$ (spanned by the FE basis vectors) which becomes an increasingly good approximation of it as $m \to \infty$. Taking both error decays as $n,m \to \infty$ together suggests that we are approaching the solution to the fully-continuous POD problem; where the time axis is still continuous and $X$ is given by $V$, exceedingly well.\hfill$\blacklozenge$
\end{example}

\subsection{The discrete POD method revisited}
\label{Section:DicPODrev}

Theorem~\ref{Theorem2.3.1} relies essentially on the convergence of $\mathcal R^n$ to $\mathcal R$ in the operator norm; cf. \eqref{Eq2.3.1}. In this section we give a sufficient condition to ensure \eqref{Eq2.3.1}. To simplify the presentation we consider the case ${K}=1$ and apply a (one-dimensional) trapezoidal approximation of the integral operator $\mathcal R$ defined in \eqref{SIAM:Eq-I.1.2.5}. As in Section~\ref{SIAM:Section-2.1.2.1} let the parameter domain $\mathscr D\subset\mathbb R$ be given by the interval $\mathscr D=[\bmu^a,\bmu^b]$. For chosen $n\in\mathbb N$ we define a grid on $[\bmu_\mathsf a,\bmu_\mathsf b]$ as follows:
\begin{align*}
    \bmu_\mathsf a=\bmu_1<\ldots<\bmu_n=\bmu_\mathsf b,\quad\delta\bmu_j=\bmu_j-\bmu_{j-1}\quad\text{for }2\le j\le n.
\end{align*}
We set $\Delta\bmu=\max_{2\le j\le n}\delta\bmu_j$ and choose the trapezoidal weights
\begin{align*}
    \alpha_1^n=\frac{\delta\bmu_2}{2},\quad \alpha_j^n=\frac{\delta\bmu_j+\delta\bmu_{j-1}}{2} \quad \text{for } 2\le j\le n-1,\quad \alpha_n^n=\frac{\delta\bmu_n}{2}.
\end{align*}
With the above choices of $\{\bmu_j\}_{j=1}^n$ and $\{\alpha_j^n\}_{j=1}^n$, we can view $\mathcal R^n$ as a trapezoidal approximation of $\mathcal R$.

\begin{proposition}
    \label{Pro2.3.1}
    Suppose that the snapshot trajectories $y^k$ belong to $H^1(\mathscr D;X)$ for $1\le k\le {K}$ and we have $\Delta \bmu \to 0$ as $n \to \infty$. Then \eqref{Eq2.3.1} holds true.
\end{proposition}

\begin{remark}
    \label{Remark2.3.1}
    \rm
    \begin{enumerate}
        \item [1)] Theorem~\ref{Theorem2.3.1} and Proposition~\ref{Pro2.3.1} give an answer to the two questions posed at the beginning of Section~\ref{SIAM:Section-2.1.2}: The grid points $\{\bmu_j\}_{j=1}^n$ and the associated positive weights $\{\alpha_j^n\}_{j=1}^n$ should be chosen such that $\mathcal R^n$ is a quadrature approximation of $\mathcal R$ and $\|\mathcal R^n-\mathcal R\|_{\mathscr L(X)}$ is small (for reasonable $n$). In Proposition~\ref{Pro2.3.1} it was shown that this can for example be obtained if the weights $\{\alpha_j^n\}_{j=1}^n$ are chosen as trapezoidal weights and the maximal grid size converges to 0. Clearly, other choices for the weights $\{\alpha_j^n\}_{j=1}^n$ are also possible provided \eqref{Eq2.3.1} is guaranteed. For instance, we can choose the Simpson weights. A different strategy is applied in \cite{KV10,LV14}, where the time instances $\{\bmu_j\}_{j=1}^n$ are chosen by an optimization approach. 
        \item [2)] In a more general case, $\mathscr D \subset \mathbb R^{\mathfrak p}$ could be an arbitrary closed, bounded and connected domain with $\mathfrak p > 1$. In this case, the simple one-dimensional discretization of $\mathscr D$ has to be replaced by a more complicated finite difference or FE grid and according weights to approximate the integration over $\mathscr D$. The above proof of convergence $\mathcal R^n \to \mathcal R$ in the operator norm then has to be adapted for this more general case.\hfill$\blacksquare$
    \end{enumerate}
\end{remark}

\subsection{Perturbation analysis in $\bm{\mathbb C^m}$}
\label{Section:PertEuclCase}

Let us assume that $X=\mathbb C^m$. To simplify the presentation we assume $K=1$ and $\omega_1^{K}=1$. The snapshots $y_j=(\mathrm y_{ij})_{1\le i\le m}\in\mathbb C^m$, $1\le j\le n$, are again supposed to be approximations of $y(\bmu)=(\mathrm y_i(\bmu))_{1\le i\le m}\in\mathbb C^m$ at the grid points $\bmu=\bmu_j\in\mathscr D$, $1\le j\le n$. Again we set $\bY=[y_1\,|\ldots|\,y_n]\in\mathbb C^{m\times n}$. In Section~\ref{SIAM:Section-2.1.1.2} we have shown that $\mathcal R^n=\bY\bY^\mathsf H$ holds. Next we derive a matrix representation of the operator $\mathcal R$. For $\mathrm v=(\mathrm v_i)_{1\le i\le m}\in\mathbb C^m$ and $1\le i\le m$ we have
\begin{align*}
    (\mathcal R\mathrm v)_i&=\int_\mathscr D {\langle\mathrm v,y(\bmu)\rangle}_{\mathbb C^m}\,\mathrm y_i(\bmu)\,\mathrm d\bmu=\int_\mathscr D \sum_{j=1}^m \mathrm y_i(\bmu)\overline{\mathrm y_j(\bmu)}\mathrm v_j\,\mathrm d\bmu\\
    &=\sum_{j=1}^m\Big(\int_\mathscr D  \mathrm y_i(\bmu)\overline{\mathrm y_j(\bmu)} \,\mathrm d\bmu\Big)\mathrm v_j.
\end{align*}
Therefore, the matrix $\mathcal R\in\mathbb C^{m\times m}$ is given by the elements
\begin{align*}
    \mathcal R_{ij}=\int_\mathscr D \mathrm y_i(\bmu)\overline{\mathrm y_j(\bmu)} \,\mathrm d\bmu\quad\text{for }1\le i,j\le m.
\end{align*}
Notice that the matrix $\mathcal R$ can be interpreted as a correlation matrix. Since both $\mathcal R^n$ and $\mathcal R$ are Hermitian matrices, we can apply well-known results from perturbation theory to estimate the differences between associated eigenvalues; see, e.g., \cite[Chapter~1]{Cha11}.

\begin{proposition}
    Suppose that the first $\ell\le d$ eigenvalues of $\mathcal R$ have algebraic and geometric multiplicity one, i.e., $\hat\lambda_1>\hat\lambda_2>\ldots>\hat\lambda_\ell$. If $\varepsilon=\|\mathcal R-\mathcal R^n\|_2$ is small enough, we have
    \begin{align*}
        \big|\hat\lambda_i-\hat\lambda_i^n\big|=\mathcal O(\varepsilon)\quad\text{and}\quad\mathrm{dist}\,\left(\hat\psi_i^n,\mathrm{ker}\big(\mathcal R-\hat\lambda_i\bI\big)\right)=\mathcal O(\varepsilon)\quad\text{for }1\le i\le\ell
    \end{align*}
    with the minimal distance function
    \begin{align*}
        \mathrm{dist}\,\left(\hat\psi_i^n,\mathrm{ker}\big(\mathcal R-\hat\lambda_i\bI\big)\right)=\min\left\{{\|\hat\psi_i^n-v\|}_2\,\big|\,v\in\mathrm{ker}\big(\mathcal R-\hat\lambda_i\bI\big)\right\}
    \end{align*}
    and $\bI\in\mathbb C^{m\times m}$ stands for the \index{Matrix!identity, $\bI$}{\em identity matrix}.
\end{proposition}

\begin{remark}
    \rm For the more general case that the algebraic and geometric multiplicity of the eigenvalues are different from one, we refer the reader to \cite[Chapter~2]{Wil65} and \cite[Chapter~2]{Kat80}.\hfill$\blacksquare$
\end{remark}

\section{POD for embedded Hilbert spaces}
\label{Section:PODHilbert}
\setcounter{equation}{0}
\setcounter{theorem}{0}
\setcounter{figure}{0}
\setcounter{run}{0}

In this section we discuss a specific choice for the Hilbert space $X$ which often occurs in conjunction with reduced-order modelling of evolution problems, cf. Chapter \ref{SIAM-Book:Section3}. It turns out that in this case it is possible to derive rate of convergence results which will be very useful in our a-priori error analysis carried out in Sections~\ref{SIAM-Book:Section3} and \ref{SIAM-Book:Section4}. The presented results are essentially based on \cite{Sin14}. We make the following basic assumption throughout this section.
\begin{assumption}
    \label{AssHV}
    Let $H$ and $V$ be two separable Hilbert spaces with $V\subset H$. We assume that $V$ is dense in $H$, i.e., for any $\varphi\in H$ there exists a sequence $\{\varphi_\nu\}_{\nu\in\mathbb N}$ in $V$ such that $\|\varphi_\nu-\varphi\|_H\to0$ as $\nu\to\infty$. Moreover, $V$ is continuously embedded in $H$.
\end{assumption}

\begin{remark}
    \label{rem:gelfandTiple}
    \rm It follows from Assumption~\ref{AssHV} that there exists a constant $c_V>0$ satisfying
    \begin{equation}
        \label{Poincare}
        {\|\varphi\|}_H\le c_V\,{\|\varphi\|}_V\quad\text{for all }\varphi\in V.
    \end{equation}
    Throughout the book we refer to \eqref{Poincare} with the name \index{Inequality!embedding}{\em embedding inequality}. Identifying the Hilbert space $H$ with its dual space $H'$, we get the so-called \index{Gelfand triple}{\em Gelfand triple} $V\hookrightarrow H\simeq H'\hookrightarrow V'$.\hfill$\blacksquare$
\end{remark}

\begin{example}
    \rm Let $\Omega \subset \mathbb{R}^m$ be a bounded domain.
    \begin{enumerate}
        \item [1)] We have previously defined the Hilbert spaces $H=L^2(\Omega)$ and $V=H^1(\Omega)$ in Examples~\ref{Example:PODParabolic} and \ref{Example:PODElliptic} with their usual inner products
        \begin{align*}
            {\langle\varphi,\phi\rangle}_H&=\int_\Omega\varphi(\bx)\phi(\bx)\,\mathrm d\bx&&\text{for }\varphi,\phi\in H,\\
            {\langle\varphi,\phi\rangle}_V&=\int_\Omega\varphi(\bx)\phi(\bx)+\nabla\varphi(\bx)\cdot\nabla\phi(\bx)\,\mathrm d\bx&&\text{for }\varphi,\phi\in V
        \end{align*}
        and the induced norms $\|\cdot\|_H=\langle\cdot\,,\cdot\rangle_H^{1/2}$ and $\|\cdot\|_V=\langle\cdot\,,\cdot\rangle_V^{1/2}$, respectively. It is well-known that these spaces are separable and that $V$ is densely as well as continuously embedded in $H$; see, e.g., \cite{Eva08}. In this case the embedding constant is given by $c_V = 1$.
        \item [2)] Another example for Hilbert spaces which form a Gelfand triple is given by $H=L^2(\Omega)$ and $V=H^1_0(\Omega)$, where the inner product in $V$ is defined as
        \begin{align*}
            {\langle\varphi,\phi\rangle}_V=\int_\Omega\nabla\varphi(\bx)\cdot\nabla\phi(\bx)\,\mathrm d\bx\quad\text{for }\varphi,\phi\in V
        \end{align*}
        with the induced norm $\|\cdot\|_V=\langle\cdot\,,\cdot\rangle_V^{1/2}$. In this case the embedding $V \hookrightarrow H$ can be obtained by the Poincar\'{e} inequality and the embedding constant $c_V$ is hence given by the Poincar\'{e} constant; see \cite[p.~265]{Eva08}, for instance.\hfill$\blacklozenge$
    \end{enumerate}
    \hfill
\end{example}

\subsection{Continuous POD version}
\label{Section:ContPODHilbert}

Now \index{POD method!continuous variant}we study \eqref{SIAM:Eq-I.1.2.3} for two different choices of the separable Hilbert space $X$, namely for $X=H$ and $X=V$. Suppose that we are given snapshots $y^k\in L^2(\mathscr D;V)$, $1\le k\le {K}$. Then the snapshot space $\mathscr V$ is a subset of $V$, and thus in particular of $H$. We consider
\begin{equation}
    \tag{$\mathbf P^\ell_H$}
    \label{SIAM:PellH}
    \left\{
    \begin{aligned}
        &\min\sum_{k=1}^{K} \omega_k^{K}\int_\mathscr D \Big\| y^k(\bmu)-\sum_{i=1}^\ell {\langle y^k(\bmu),\psi_i\rangle}_H\,\psi_i\Big\|_H^2\,\mathrm d\bmu\\
        &\hspace{0.5mm}\text{s.t. } \{\psi_i\}_{i=1}^\ell\subset H\text{ and }{\langle\psi_i,\psi_j\rangle}_H=\delta_{ij},~1 \le i,j \le \ell,
    \end{aligned}
    \right.
\end{equation}
and
\begin{equation}
    \tag{$\mathbf P^\ell_V$}
    \label{SIAM:PellV}
    \left\{
    \begin{aligned}
        &\min\sum_{k=1}^{K} \omega_k^{K}\int_\mathscr D \Big\| y^k(\bmu)-\sum_{i=1}^\ell {\langle y^k(\bmu),\psi_i\rangle}_V\,\psi_i\Big\|_V^2\,\mathrm d\bmu\\
        &\hspace{0.5mm}\text{s.t. } \{\psi_i\}_{i=1}^\ell\subset V\text{ and }{\langle\psi_i,\psi_j\rangle}_V=\delta_{ij},~1 \le i,j \le \ell.
    \end{aligned}
    \right.
\end{equation}

\begin{example}
    \label{ExamplePertAnal-1}
    \rm In Example~\ref{Example:PODParabolic}, we have already considered the parabolic example with the choices $X=V$ and $X=H$. For convenience, we consider only the continuous problem \eqref{SIAM:Eq-I.1.2.3} in our notation. Of course, all numerical results are generated by fully discrete methods and can only serve as an illustration of this infinite-dimensional case. Figure~\ref{fig:podInitialConditions_singVals} already shows the difference between the eigenvalues $\hat\lambda_i^V$ and $\hat\lambda_i^H$. It can be seen the $\hat\lambda_i^V$ is always larger than $\hat\lambda_i^H$. The reason for this lies in the fact that we have 
    \begin{align*}
        {\|\psi\|}_V^2 = {\|\psi\|}_H^2 + \int_\Omega |\nabla \psi(\bx)|^2_2 ~\mathrm d\bx \ge{\|\psi\|}_H^2 \qquad \text{for all }\psi \in V.
    \end{align*}
    This means that the $V$-basis has to approximate both the values of $y^k(t,\cdot)$ and their gradients $\nabla y^k(t,\cdot)$, while the $H$-basis only has to do the former, which is clearly the easier task. By the growing discrepancy between $\hat\lambda_i^V$ and $\hat\lambda_i^H$ as $i$ increases, we can observe that this difference in the approximation demand increases as well. In addition, note that $\{\hat\psi_i^H\}_{i=1}^\ell$ ``only'' has to be $H$-orthonormal whereas $\{\hat\psi_i^V\}_{i=1}^\ell$ is supposed to be $V$-orthonormal. In terms of functions, this especially implies that each $\hat\psi_i^V$ is limited in terms of sharp edges since high values of the gradient have to be avoided due to $\|\hat\psi_i^V\|_V=1$. 
    \begin{figure}
        \centering
        \includegraphics[width=110mm,height=120mm]{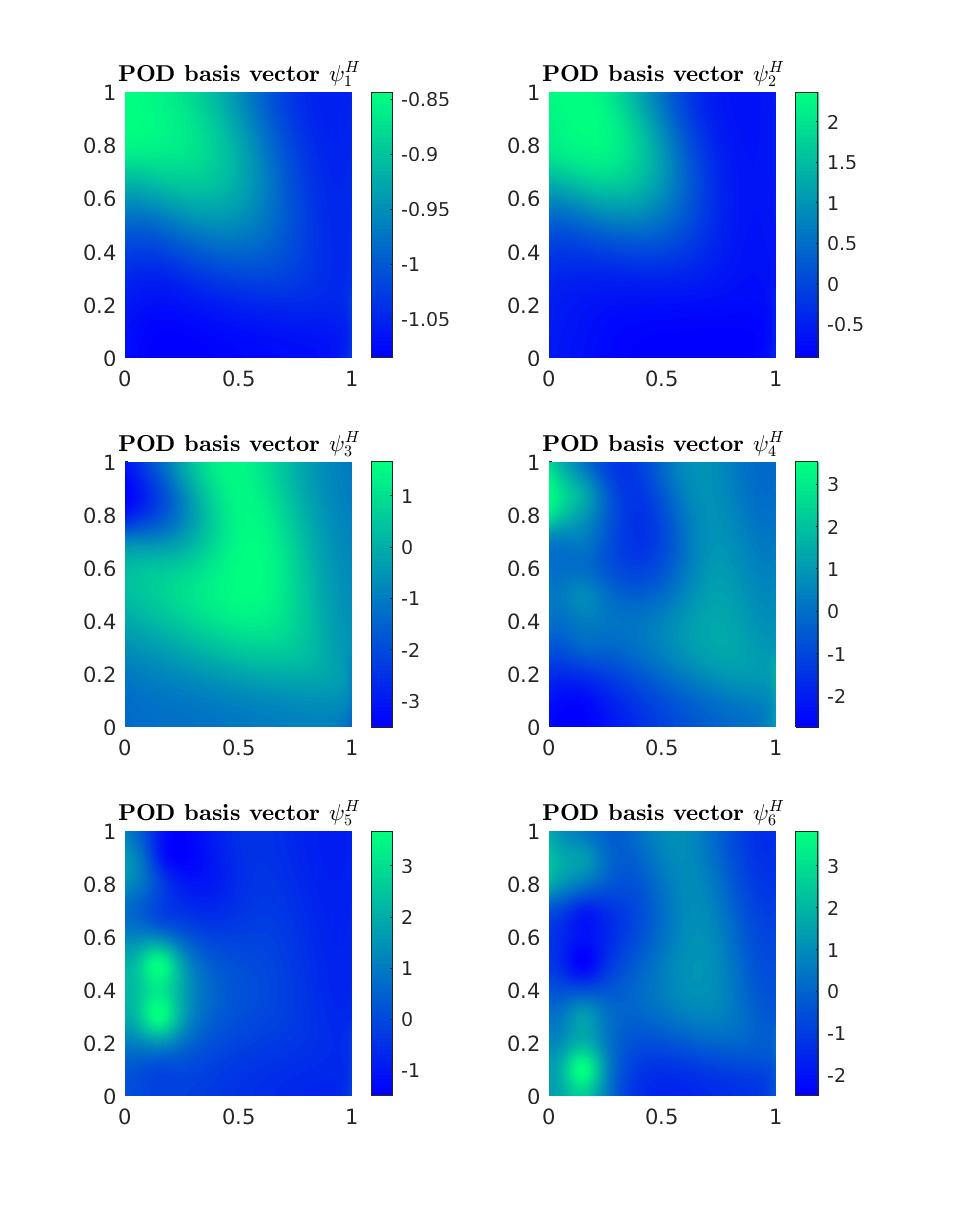}
        \caption{Example~\ref{ExamplePertAnal-1}. First six POD basis vectors $\hat\psi_i^H$.}
        \label{fig:podInitialConditions_podVectors_H}
    \end{figure}
    In Figure \ref{fig:podInitialConditions_podVectors_H}, we see the first six POD vectors for the choice $X=H$, which can be compared to Figure \ref{fig:podInitialConditions_podVectors} with the choice $X=V$. We can immediately validate the above observation since every basis vector $\hat\psi_i^H$ has a much wider range than its $V$-counterpart $\hat\psi_i^V$. 
	
    \noindent
    In conclusion, the choice $X=H$ makes the problem much easier since the basis is less limited in terms of gradient jumps and the demanded approximation quality is lower than for $X=V$. In practise, the users must make their choice based on the problem at hand. For example, if it is essential to include gradient information in the basis, the $V$-approach is imperative. We can also consider mixed cases, where, e.g., we require a $V$-orthonormal POD basis, but look for only $H$-approximation of the trajectories, which leads to a problem of the type 
    \begin{align*}
        \left\{
        \begin{aligned}
            &\min\sum_{k=1}^K \omega_k \int_0^T \bigg\| y^k(t)-\sum_{i=1}^\ell{\langle y^k(t),\psi_i \rangle}_H \psi_i\bigg\|_H^2\,\mathrm dt\\
            &\hspace{1mm}\text{s.t. }\{\psi_i\}_{i=1}^\ell\subset V\text{ and }{\langle \psi_i, \psi_j \rangle}_V = \delta_{ij} \quad \text{for } i,j=1,...,\ell.
        \end{aligned}
        \right.
    \end{align*}
    Now, the first-order conditions are given as a generalized eigenvalue problem.\hfill$\blacklozenge$
\end{example}

\subsubsection{Optimal solutions to \eqref{SIAM:PellH} and \eqref{SIAM:PellV}}

The solutions to \eqref{SIAM:PellH} and \eqref{SIAM:PellV} are given by Theorem~\ref{Theorem2.2.1}. Therefore, let us define the two linear operators 
\begin{equation}
    \label{OperatorsR}
    \begin{aligned}
        \mathcal R_H\psi&=\sum_{k=1}^{K}\omega_k^{K}\int_\mathscr D {\langle\psi,y^k(\bmu) \rangle}_H\,y^k(\bmu)\,\mathrm d\bmu&&\text{for } \psi \in H,\\
        \mathcal R_V\psi&=\sum_{k=1}^{K}\omega_k^{K}\int_\mathscr D {\langle\psi,y^k(\bmu) \rangle}_V\,y^k(\bmu)\,\mathrm d\bmu&&\text{for } \psi \in V.
    \end{aligned}
\end{equation}
We set
\begin{align*}
    d^H=\left\{
    \begin{aligned}
        &\mathrm{rank}\big(\mathrm{ran}(\mathcal R_H)\big)&&\text{if }\mathrm{rank}\big(\mathrm{ran}(\mathcal R_H)\big)\text{ is finite},\\
        &\infty&&\text{otherwise}
    \end{aligned}
    \right.
\end{align*}
and
\begin{align*}
    d^V=\left\{
    \begin{aligned}
        &\mathrm{rank}\big(\mathrm{ran}(\mathcal R_V)\big)&&\text{if }\mathrm{rank}\big(\mathrm{ran}(\mathcal R_V)\big)\text{ is finite},\\
        &\infty&&\text{otherwise}.
    \end{aligned}
    \right.
\end{align*}
Since $y^k\in L^2(0,T;V)$ holds, we have $\mathcal R_H\psi\in V$ for all $\psi\in H$ and $\mathcal R_V \psi \in V$ for all $\psi \in V$. Nevertheless, $\mathcal R_H$ is considered by definition as an operator mapping from $H$ to $H$. By Lemma~\ref{Lemma2.2.1} both $\mathcal R_H$ and $\mathcal R_V$ are compact, self-adjoint and non-negative. Thus, it follows from Theorem~\ref{Theorem2.2.1} that there exist non-negative eigenvalues $\{\hat\lambda_i^H\}_{i\in\mathbb I}$ and associated orthonormal eigenfunctions $\{\hat\psi_i^H\}_{i\in\mathbb I} \subset H$ satisfying
\begin{equation}
    \label{Eq:PODEigH}
    \mathcal R_H\hat\psi_i^H=\hat\lambda_i^H\hat\psi_i^H,\hspace{2mm}\left\{
    \begin{aligned}
        &\hat\lambda_1^H\ge\ldots\ge\hat\lambda_{d^H}^H>\hat\lambda_{d^H+1}^H=\ldots =0&&\text{if }d^H<\infty,\\
        &\hat\lambda_1^H\ge\hat\lambda_2^H\ge\ldots\text{ and }\lim_{i\to\infty}\hat\lambda_i^H=0&&\text{if }d^H=\infty.
    \end{aligned}
    \right.
\end{equation}
with $d^H\in\mathbb N \cup \{ \infty \}$. For every $\ell\in\{1,\ldots,d^H\}$ the first $\ell$ eigenfunctions $\{\hat\psi_i^H\}_{i=1}^\ell$ solve \eqref{SIAM:PellH}. When the POD basis $\{\hat\psi_i^H\}_{i=1}^\ell$ of rank $\ell$ is computed, we set
\begin{align*}
    H^\ell=\mathrm{span}\,\{\hat\psi_1^H,\ldots,\hat\psi_\ell^H\}.
\end{align*}
Analogously, Theorem~\ref{Theorem2.2.1} implies the existence of non-negative eigenvalues $\{\hat\lambda_i^V\}_{i\in\mathbb I}$ and associated orthonormal eigenfunctions $\{\hat\psi_i^V\}_{i\in\mathbb I} \subset V$ satisfying
\begin{equation}
    \label{Eq:PODEigV}
    \mathcal R_V\hat\psi_i^V=\hat\lambda_i^V\hat\psi_i^V,\hspace{2mm}\left\{
    \begin{aligned}
        &\hat\lambda_1^V\ge\ldots\ge\hat\lambda_{d^V}^V>\hat\lambda_{d^V+1}^V=\ldots =0&&\text{if }d^V<\infty,\\
        &\hat\lambda_1^V\ge\hat\lambda_2^V\ge\ldots\text{ and }\lim_{i\to\infty}\hat\lambda_i^V=0&&\text{if }d^V=\infty.
    \end{aligned}
    \right.
\end{equation}
with $d^V\in\mathbb N \cup \{ \infty \}$. The first $\ell$ eigenfunctions $\{\hat\psi_i^V\}_{i=1}^\ell$ solve \eqref{SIAM:PellV} for every $\ell\in\{1,\ldots,d^V\}$. We define the $\ell$-dimensional subspace
\begin{align*}
    V^\ell=\mathrm{span}\,\{\hat\psi_1^V,\ldots,\hat\psi_\ell^V\}\subset V.
\end{align*}
Moreover, we infer from \eqref{Eq2.2.17} the approximation formulas
\begin{subequations}
    \label{Eq:PODErrForm}
    \begin{align}
        \label{Eq:PODErrForm-H}
        \sum_{k=1}^{K} \omega_k^{K}\int_\mathscr D\Big\| y^k(\bmu)-\sum_{i=1}^\ell {\langle y^k(\bmu),\hat\psi_i^H\rangle}_H\,\hat\psi_i^H\Big\|_H^2\,\mathrm d\bmu=\sum_{i>\ell}\hat\lambda_i^H,\\
        \label{Eq:PODErrForm-V}
        \sum_{k=1}^{K} \omega_k^{K}\int_\mathscr D\Big\| y^k(\bmu)-\sum_{i=1}^\ell {\langle y^k(\bmu),\hat\psi_i^V\rangle}_V\,\hat\psi_i^V\Big\|_V^2\,\mathrm d\bmu=\sum_{i>\ell}\hat\lambda_i^V.
    \end{align}
\end{subequations}

\begin{remark}
    \label{SIAM:Remark_PODBasis_LinIndep}
    \rm As $\{\hat\psi_i^H\}_{i=1}^\ell$ and $\{\hat\psi_i^V\}_{i=1}^\ell$ are orthonormal in $H$ and $V$, respectively, we derive that $((\langle \hat\psi_i^H,\hat\psi_j^H \rangle_H)) = \bI$ and $((\langle \hat\psi_i^V,\hat\psi_j^V \rangle_V)) = \bI$, respectively. Additionally, as both $\{\hat\psi_i^H\}_{i=1}^\ell$ and $\{\hat\psi_i^V\}_{i=1}^\ell$ are linearly independent in $V$ and $H$, respectively, it is easy to see that the matrices $((\langle \hat\psi_i^H,\hat\psi_j^H \rangle_V))$ and $((\langle \hat\psi_i^V,\hat\psi_j^V \rangle_H))$ are positive definite, respectively.\hfill$\blacksquare$
\end{remark}

The relationship between the eigenvalues $\hat\lambda_i^H$ and $\hat\lambda_i^V$ is investigated in the next lemma, which is taken from \cite[Lemma 3.1]{Sin14}. We also refer the reader to \cite[Lemma~3.4]{GV17}.

\begin{lemma}
    \label{SIAM:Lemma3.2.1}
    Let Assumption~{\em\ref{AssHV}} be satisfied. Suppose that the snapshots fulfill $y^k\in L^2(0,T;V)$ for $k=1,\ldots,{K}$. Then:
    \begin{enumerate}
        \item [\rm 1)] For all $i\in\{1,\ldots,d_H\}$ we have $\hat\psi^H_i\in V$.
        \item [\rm 2)] $\lambda_i^V=0$ for all $i > d$ with some $d\in\mathbb N$ if and only if $\lambda_i^H=0$ for all $i > d$, i.e., we have $d^H=d^V$ provided the rank of $\mathcal R_V$ or of $\mathcal R_H$ is finite.
        \item [\rm 3)] $\lambda_i^V>0$ for all $i\in\mathbb I$ if and only if $\lambda_i^H>0$ for all $i\in\mathbb I$.
    \end{enumerate}
\end{lemma}

\begin{remark}
    \label{Remark:DimD}
    \rm We infer from Lemma~\ref{SIAM:Lemma3.2.1}-1) it follows that the space $H^\ell$ is an $\ell$-dimensional subspace of $V$. Furthermore, Lemma~\ref{SIAM:Lemma3.2.1}-1) and 2) imply that only one of the following cases can occur:
    \begin{enumerate}
        \item [1)] There exists a $d\in\mathbb N$ such that $d=d_H=d_V$ and $\lambda_i^H=\lambda_i^V=0$ for $i>d$.
        \item [2)] For all $i\in\mathbb N$ we have $\lambda_i^H>0$ as well as $\lambda_i^V>0$.\hfill$\blacksquare$
    \end{enumerate}
\end{remark}

\subsubsection{POD projection operators}

Let us introduce the following four projection operators
\begin{subequations}
    \label{SIAM:Eq3.2.12}
    \begin{align}
        \label{SIAM:Eq3.2.12a}
        &\mathcal P^\ell_H:H\to H^\ell,&&v^\ell=\mathcal P^\ell_H\varphi&&\text{solves }\min_{w^\ell\in H^\ell}{\|\varphi-w^\ell\|}_H\text{ for }\varphi\in H,\\
        \label{SIAM:Eq3.2.12b}
        &\mathcal P^\ell_V:V\to H^\ell,&&v^\ell=\mathcal P^\ell_V\varphi&&\text{solves }\min_{w^\ell\in H^\ell}{\|\varphi-w^\ell\|}_V\text{ for }\varphi\in V,\\
        \label{SIAM:Eq3.2.12c}
        &\mathcal Q^\ell_H:H\to V^\ell,&&u^\ell=\mathcal Q^\ell_H\varphi&&\text{solves } \min_{w^\ell\in V^\ell}{\|\varphi-w^\ell\|}_H\text{ for }\varphi\in H,\\
        \label{SIAM:Eq3.2.12d}
        &\mathcal Q^\ell_V:V\to V^\ell,&&u^\ell=\mathcal Q^\ell_V\varphi&&\text{solves } \min_{w^\ell\in V^\ell}{\|\varphi-w^\ell\|}_V\text{ for }\varphi\in V.
    \end{align}
\end{subequations}
These operators will play an essential role in our error analysis carried out in Sections~\ref{SIAM-Book:Section3} and \ref{SIAM-Book:Section4}.
 
In the following lemma we state how the projection operators can be computed. Furthermore, we list essential properties of these operators. For the proof we refer the reader to Section~\ref{SIAM:Section-2.6.4}.
 
\begin{proposition}
    \label{Lemma:ProjOper}
    Let Assumption~{\em\ref{AssHV}} hold. Suppose that the snapshots fulfill $y^k\in L^2(0,T;V)$ for $k=1,\ldots,{K}$. The four projections in \eqref{SIAM:Eq3.2.12} are given as follows:
    \begin{equation}
        \label{Eq:Prolongations-12}
        \begin{aligned}
            \mathcal P^\ell_H\varphi&=\sum_{i=1}^\ell {\langle\varphi,\hat\psi_i^H\rangle}_H\,\hat\psi_i^H&&\text{for }\varphi\in H,&
            \mathcal P^\ell_V\varphi&=\sum_{i=1}^\ell \mathrm v_i\hat\psi_i^H&&\text{for }\varphi\in V,\\
            \mathcal Q^\ell_V\varphi&=\sum_{i=1}^\ell {\langle\varphi,\hat\psi_i^V\rangle}_V\,\hat\psi_i^V&&\text{for }\varphi\in V,&\mathcal Q^\ell_H\varphi&=\sum_{i=1}^\ell \mathrm u_i\hat\psi_i^V&&\text{for }\varphi\in H,
        \end{aligned}
    \end{equation}
    where the vectors $\mathrm v=(\mathrm v_1,\ldots,\mathrm v_\ell)^\top,\,\mathrm u=(\mathrm  u_1,\ldots,\mathrm u_\ell)^\top\in\mathbb R^\ell$ solve the following linear systems
    \begin{subequations}
        \label{SIAM:Eq3.2.14}
        \begin{align}
            \label{SIAM:Eq3.2.14b}
            \sum_{j=1}^\ell{\langle\hat\psi_j^H,\hat\psi_i^H\rangle}_V\,\mathrm v_j={\langle\varphi,\hat\psi_i^H\rangle}_V&&\text{for }1\le i\le \ell,\\
            \label{SIAM:Eq3.2.14a}
            \sum_{j=1}^\ell{\langle\hat\psi_j^V,\hat\psi_i^V\rangle}_H\,\mathrm u_j={\langle\varphi,\hat\psi_i^V\rangle}_H&&\text{for }1\le i\le \ell,
        \end{align}
    \end{subequations}
    respectively. All four projections are linear and bounded. Furthermore, $\mathcal P_H^\ell,\,\mathcal Q_H^\ell$ are orthogonal in $H$ and $\mathcal P_V^\ell,\,\mathcal Q_V^\ell$ are orthogonal in  $V$.
\end{proposition}

\begin{remark}
    \label{Rem:ConvergenceRate}
    \rm Let Assumption~\ref{AssHV} hold. Then we infer from Proposition~\ref{Lemma:ProjOper} that the approximation formulas \eqref{Eq:PODErrForm} can be expressed as
    \begin{subequations}
        \label{ErrRate}
        \begin{align}
            \label{ErrRateH}
            \sum_{k=1}^{K} \omega_k^{K}\int_\mathscr D\Big\| y^k(\bmu)-\mathcal P_H^\ell y^k(\bmu)\Big\|_H^2\,\mathrm d\bmu=\sum_{i>\ell}\hat\lambda_i^H,\\
            \label{ErrRate-V}
            \sum_{k=1}^{K} \omega_k^{K}\int_\mathscr D\Big\| y^k(\bmu)-\mathcal Q_V^\ell y^k(\bmu)\Big\|_V^2\,\mathrm d\bmu=\sum_{i>\ell}\hat\lambda_i^V.
        \end{align}
    \end{subequations}
    Hence, we do not only have $\|y^k-\mathcal P_H^\ell y^k\|^2_{L^2(\mathscr D;H)}\to0$ and $\|y^k-\mathcal Q_V^\ell y^k\|^2_{L^2(\mathscr D;V)}\to0$ as $\ell\to\infty$ for $k=1,\ldots,K$, but the rate of convergence is given as $\sum_{i>\ell}\hat\lambda_i^H$ and $\sum_{i>\ell}\hat\lambda_i^V$, respectively. Especially, the decay rate of the eigenvalues $\{\hat\lambda_i^H\}_{i\in\mathbb I}$ or $\{\hat\lambda_i^V\}_{i\in\mathbb I}$ is essential for the speed of convergence.\hfill$\blacksquare$
\end{remark}

\subsubsection{Singular value decomposition}

In Remark~\ref{Remark2.2.1a} we have introduced the \index{Singular value decomposition, SVD}{\em singular value decomposition} of the operator $\mathcal Y$ defined by \eqref{Eq2.2.6}. By $\mathcal Y_H$ and $\mathcal Y_V$ we denote the operator $\mathcal Y$ considered as a mapping in $\mathscr L(L^2(\mathscr D,\mathbb C^{K}),H)$ and $\mathscr L(L^2(\mathscr D,\mathbb C^{K}),V)$, respectively. Then the associated adjoints $\mathcal Y_H^\star$ and $\mathcal Y_V^\star$ satisfy
\begin{align*}
    \big(\mathcal Y^\star_H\psi)(\bmu)&=\left(
    \begin{array}{c}
        {\langle\psi,y^k(\bmu)\rangle}_H\\
        \vdots\\
        {\langle\psi,y^k(\bmu)\rangle}_H
    \end{array}
    \right)&&\text{for }\psi\in H\text{ and }\bmu\in\mathscr D,\\
    \big(\mathcal Y^\star_V\psi)(\bmu)&=\left(
    \begin{array}{c}
        {\langle\psi,y^k(\bmu)\rangle}_V\\
        \vdots\\
        {\langle\psi,y^k(\bmu)\rangle}_V
    \end{array}
    \right)&&\text{for }\psi\in V\text{ and }\bmu\in\mathscr D.
\end{align*}
The operators $\mathcal K_H=\mathcal Y^{H,\star}\mathcal Y^H,\mathcal K_V=\mathcal Y^{V,\star}\mathcal Y^V:L^2(\mathscr D;\mathbb C^{K}) \to L^2(\mathscr D;\mathbb C^{K})$ are given by
\begin{align*}
    \big(\mathcal K_H\phi\big)(\bmu) =\left(
    \begin{array}{c}
        \sum\limits_{k=1}^{K} \int_\mathscr D{\langle y^k(\bnu),y^1(\bmu)\rangle}_H\phi^k(\bnu)\,\mathrm d\bnu\\
        \vdots\\
        \sum\limits_{k=1}^{K} \int_\mathscr D{\langle y^k(\bnu),y^{K}(\bmu)\rangle}_H\phi^k(\bnu)\,\mathrm d\bnu
    \end{array}
    \right)\quad\text{for }\phi\in L^2(\mathscr D;\mathbb C^{K}),
\end{align*}
and
\begin{align*}
    \big(\mathcal K_V\phi\big)(\bmu) =\left(
    \begin{array}{c}
        \sum\limits_{k=1}^{K} \int_\mathscr D{\langle y^k(\bnu),y^1(\bmu)\rangle}_V\phi^k(\bnu)\,\mathrm d\bnu\\
        \vdots\\
        \sum\limits_{k=1}^{K} \int_\mathscr D{\langle y^k(\bnu),y^{K}(\bmu)\rangle}_V\phi^k(\bnu)\,\mathrm d\bnu
    \end{array}
    \right)\quad\text{for }\phi\in L^2(\mathscr D;\mathbb C^{K}),
\end{align*}
respectively. Due to \index{POD method!continuous variant!SVD}singular value decomposition, there exists a sequence $\{\hat\phi_i^H\}_{i\in\mathbb I}\subset L^2(\mathscr D;\mathbb C^{K})$ satisfying
\begin{align*}
    \mathcal K_H\hat\phi_i^H=\hat\lambda_i^H\hat\phi_i^H,\quad \mathcal R_H\hat\psi_i^H=\hat\lambda_i^H\hat\psi_i^H, \quad \mathcal Y_H\hat\phi_i^H=\hat\sigma_i^H\hat\psi_i^H,\quad\mathcal Y^\star_H\hat\psi_i^H=\hat\sigma_i^H\hat\phi_i^H
\end{align*}
for $i\in\mathbb I$ and $\hat\sigma_i^H=(\hat\lambda_i^H)^{1/2}$. Analogously, there is a sequence $\{\hat\phi_i^V\}_{i\in\mathbb I}\subset L^2(\mathscr D;\mathbb C^{K})$ satisfying
\begin{align*}
    \mathcal K_V\hat\phi_i^V=\hat\lambda_i^V\hat\phi_i^V,\quad \mathcal R_V\hat\psi_i^V=\hat\lambda_i^V\hat\psi_i^H,\quad \mathcal Y_V\hat\phi_i^V=\hat\sigma_i^V\hat\psi_i^V,\quad\mathcal Y^\star_V\psi_i^V=\hat\sigma_i^V\hat\phi_i^V
\end{align*}
for $i\in\mathbb I$ and $\hat\sigma_i^V=(\hat\lambda_i^V)^{1/2}$.

\subsubsection{A-priori error analysis}

In the next proposition we summarize convergence properties of the projection operators introduced in \eqref{SIAM:Eq3.2.12}. For a proof we refer the reader to \cite[Proposition~5.5]{Sin14}.

\begin{proposition}
    \label{Proposition:ProjConvV}
    Suppose that Assumption~{\em\ref{AssHV}} holds and the snapshots fulfill $y^k\in L^2(0,T;V)$ for $k=1,\ldots,{K}$. Then $\lim_{\ell\to\infty}\|\mathcal Q_V^\ell\varphi-\varphi\|_V=0$ for every $\varphi\in V$. If $\hat\lambda_i^H>0$ for all $i\in\mathbb N$, we have
    \begin{enumerate}
        \item [\em 1)] $\lim_{\ell\to\infty}\|\mathcal P_H^\ell\varphi-\varphi\|_H=0$ for every $\varphi\in\mathrm{ran}\,(\mathcal Y_H)$,
        \item [\em 2)] $\lim_{\ell\to\infty}\|\mathcal P_V^\ell\varphi-\varphi\|_V=0$ for every $\varphi\in V$.
    \end{enumerate}
\end{proposition}

Convergence properties in the $V$-topology for the snapshots $\{y^k\}_{k=1}^K\subset L^2(\mathscr D,V)$ are presented in the following theorem. The proof is given in \cite[Theorems~5.2 and 5.3]{Sin14}.

\begin{theorem}
    \label{Theorem:ProjConvV}
    Suppose that Assumption~{\em\ref{AssHV}} and $y^k\in L^2(0,T;V)$ for $k=1,\ldots,K$ hold. 
    \begin{enumerate}
        \item [\em 1)] If $\hat\lambda_i^H>0$ for all $i\in\mathbb N$, it follows that
        \begin{align*}
            \lim_{\ell\to\infty}{\|\mathcal P_H^\ell y^k-y^k\|}_{L^2(\mathscr D;V)}=\lim_{\ell\to\infty}{\|\mathcal P_V^\ell y^k-y^k\|}_{L^2(\mathscr D;V)}=0\quad\text{for }1\le k\le K.
        \end{align*}
        \item [\em 2)] If $\hat\lambda_i^V>0$ for all $i\in\mathbb N$ and $\|\mathcal Q_H^\ell\|_{\mathscr L(V)}$ is bounded independently of $\ell$, we have
        \begin{align*}
            \lim_{\ell\to\infty}{\|\mathcal Q_H^\ell y^k-y^k\|}_{L^2(\mathscr D;V)}=0\quad\text{for }1\le k\le K.
        \end{align*}
        \item [\em 3)] It holds
        \begin{align*}
            \lim_{\ell\to\infty}{\|\mathcal Q_V^\ell y^k-y^k\|}_{L^2(\mathscr D;V)}=0\quad\text{for }1\le k\le K.
        \end{align*}
    \end{enumerate}
\end{theorem}

\begin{remark}
    \rm The assumption that $\|\mathcal Q^\ell_H\|_{\mathscr L(V)}$ is bounded independently of $\ell$ was also supposed in \cite{CGS12}.\hfill$\blacksquare$
\end{remark}

Recall the error formulas \eqref{Eq:PODErrForm}. Moreover, we have further error formulas for the $V$-topology, which will play an essential role in the error analysis carried out in Chapters~\ref{SIAM-Book:Section3} and \ref{SIAM-Book:Section4}. Proofs can be found in \cite[Theorems~5.2, 5.3 and 5.4]{Sin14}.

\begin{theorem}
    \label{Prop:VTopology}
    Suppose that Assumption~{\em\ref{AssHV}} and $y^k\in L^2(0,T;V)$ for $k=1,\ldots,K$ hold. For every $\ell\ge1$ we have
    \begin{subequations}
        \label{Prop:Rate-V}
        \begin{align}
            \label{RateHH}
            \sum_{k=1}^{K}\omega_k^{K}\int_\mathscr D\big\|y^k(\bmu)-\mathcal P^\ell_Hy^k(\bmu)\big\|_H^2\,\mathrm d\bmu&=\sum_{i>\ell}\hat\lambda_i^H<\infty,\\\label{RatePH-V}
            \sum_{k=1}^{K}\omega_k^{K}\int_\mathscr D\big\|y^k(\bmu)-\mathcal P^\ell_Hy^k(\bmu)\big\|_V^2\,\mathrm d\bmu&=\sum_{i>\ell}^d\hat\lambda_i^H\,{\|\hat\psi_i^H\|}_V^2<\infty,\\
            \label{RatePV-V}
            \hspace{-1mm}\sum_{k=1}^{K}\omega_k^{K}\int_\mathscr D\big\|y^k(\bmu)-\mathcal P^\ell_Vy^k(\bmu)\big\|_H^2\,\mathrm d\bmu&=\sum_{i>\ell}\hat\lambda_i^H\big\|\hat\psi_i^H-\mathcal P^\ell_V\hat\psi_i^H\big\|_V^2<\infty,\\
            \label{RateQH-V}
            \hspace{-5mm}\sum_{k=1}^{K}\omega_k^{K}\int_\mathscr D\big\|y^k(\bmu)-\mathcal Q^\ell_Hy^k(\bmu)\big\|_V^2\,\mathrm d\bmu&=\sum_{i=\ell+1}^d\hat\lambda_i^V\,{\|\hat\psi_i^V-\mathcal Q_H^\ell\hat\psi_i^V\|}_V^2<\infty,\\
            \label{RateQV-V}
            \sum_{k=1}^{K}\omega_k^{K}\int_\mathscr D\big\|y^k(\bmu)-\mathcal Q^\ell_Vy^k(\bmu)\big\|_V^2\,\mathrm d\bmu&=\sum_{i>\ell}\hat\lambda_i^V<\infty.
        \end{align}
    \end{subequations}
\end{theorem}

\begin{remark}
    \rm Instead of only having a decay rate of the eigenvalues as in the formulas \eqref{RateHH} and \eqref{RateQV-V}, there are additional damping factors in the a-priori error formulas \eqref{RatePH-V}, \eqref{RatePV-V} and \eqref{RateQH-V}.\hfill$\blacksquare$
\end{remark}

\subsection{Discrete POD version}
\label{Section:DiscPODHilbert}

The presented results in Section~\ref{Section:ContPODHilbert} can also be applied for the \index{POD method!discrete variant}discrete variant of the POD method introduced in Section~\ref{SIAM:Section-2.1.1}. Again the separable Hilbert space $X$ stands either for $H$ or $V$. We suppose that the \index{POD method!discrete variant!snapshots}{\em snapshots} $y_j^k\in V$ are given for $1\le j\le n$ and $1\le k\le{K}$. Then the \index{POD method!discrete variant!snapshot space}{\em snapshot space} is defined as
\begin{align*}
    \mathscr V^n=\mathrm{span}\,\big\{y_j^k\,\big|\,1\le j\le n\text{ and }1\le k\le{K}\big\}\subset V\hookrightarrow H.
\end{align*}
Of course, now the snapshot set $\mathscr V^n$ has finite dimension. Due to Remark~\ref{Remark:DimD}-1) we have $d^n=\dim\mathscr V^n\le nK$ for $X=H$ and $X=V$. We consider the minimization problems
\begin{equation}
    \tag{$\mathbf P^{n\ell}_H$}
    \label{SIAM:PellHDisc}
    \left\{
    \begin{aligned}
        &\min\sum_{k=1}^{K} \omega_k^{K}\sum_{j=1}^n\alpha_j^n\,\Big\| y^k_j-\sum_{i=1}^\ell {\langle y^k_j,\psi_i\rangle}_H\,\psi_i\Big\|_H^2\\
        &\hspace{0.5mm}\text{s.t. } \{\psi_i\}_{i=1}^\ell\subset H\text{ and }{\langle\psi_i,\psi_j\rangle}_H=\delta_{ij},~1 \le i,j \le \ell,
    \end{aligned}
    \right.
\end{equation}
and
\begin{equation}
    \tag{$\mathbf P^{n\ell}_V$}
    \label{SIAM:PellVDisc}
    \left\{
    \begin{aligned}
        &\min\sum_{k=1}^{K} \omega_k^{K}\sum_{j=1}^n\alpha_j^n\,\Big\| y^k_j-\sum_{i=1}^\ell {\langle   y^k_j,\psi_i\rangle}_V\,\psi_i\Big\|_V^2\\
        &\hspace{0.5mm}\text{s.t. } \{\psi_i\}_{i=1}^\ell\subset V\text{ and }{\langle\psi_i,\psi_j\rangle}_V=\delta_{ij},~1 \le i,j \le \ell.
    \end{aligned}
    \right.
\end{equation}

\subsubsection{Optimal solutions to \eqref{SIAM:PellHDisc} and \eqref{SIAM:PellVDisc}}

The solutions to \eqref{SIAM:PellHDisc} and \eqref{SIAM:PellVDisc} are given by Theorem~\ref{SIAM:Theorem-I.1.1.3}. Therefore, let us define the two linear operators
\begin{align*}
    \mathcal R_H^n\psi=\sum_{k=1}^{K}\omega_k^K\sum_{j=1}^n\alpha_j^n\,{\langle y^k_j,\psi_i\rangle}_Hy_j^k\text{ for }\psi\in H,\quad\mathcal R_V^n\psi=\sum_{k=1}^{K}\omega_k^K\sum_{j=1}^n\alpha_j^n\,{\langle y^k_j,\psi_i\rangle}_Vy_j^k\text{ for }\psi\in V.
\end{align*}
Clearly, $\mathcal R_H^n$ and $\mathcal R_V^n$ are finite rank operators having the same rank $d^n\le nK$. Moreover, it follows from Lemma~\ref{SIAM:Lemma-I.1.1.1} that both $\mathcal R_H^n$ and $\mathcal R_V^n$ are compact, self-adjoint and non-negative. Let $\{(\hat\lambda_i^{nH},\hat\psi_i^{nH})\}_{i\in\mathbb I}$ and $\{(\hat\lambda_i^{nV},\hat\psi_i^{nV})\}_{i\in\mathbb I}$ denote the eigenvalue-eigenfunction pairs satisfying
\begin{align*}
    \mathcal R_H^n\hat\psi_i^{nH}=\hat\lambda_i^{nH}\hat\psi_i^{nH},\quad\hat\lambda_1^{nH}\ge\ldots\ge\hat\lambda_{d^n}^{nH}>\hat\lambda_{d^n+1}^{nH}=\ldots=0
\end{align*}
and
\begin{align*}
    \mathcal R_V^n\hat\psi_i^{nV}=\hat\lambda_i^{nV},n\hat\psi_i^{nV},\quad\hat\lambda_1^{nV}\ge\ldots\ge\hat\lambda_{d^n}^{nV}>\hat\lambda_{d^n+1}^{nV}=\ldots=0.
\end{align*}
For any $\ell\in\{1,\ldots,d^n\}$ the first $\ell$ eigenfunctions $\{\hat\psi_i^{nH}\}_{i=1}^\ell$ solve \eqref{SIAM:PellHDisc}. Further, $\{\hat\psi_i^{nV}\}_{i=1}^\ell$ is a solution to \eqref{SIAM:PellVDisc} for any $\ell\in\{1,\ldots,d^n\}$. Since $\mathcal R_H^n$ and $\mathcal R_V^n$ are finite rank operators, Lemma~\ref{SIAM:Lemma3.2.1} implies the next result.

\begin{lemma}
    \label{SIAM:Lemma3.2.1Disc}
    Suppose that the snapshots $y^k_j$ belong to $V$ for all $1\le j\le n$ and $1 \le k\le{K}$. Then $\hat\psi_i^{nH}\in V$ for all $i\in\{1,\ldots,d^n\}$.
\end{lemma}

\subsubsection{POD projection operators}

We introduce the two subspaces
\begin{align*}
    H^{n\ell}=\mathrm{span}\,\big\{\hat\psi_1^{nH},\ldots,\hat\psi_\ell^{nH}\big\}\subset V\quad\text{and}\quad V^{n\ell}=\mathrm{span}\,\big\{\hat\psi_1^{nV},\ldots,\hat\psi_\ell^{nV}\big\}\subset V.
\end{align*}
As in \eqref{SIAM:Eq3.2.12} we define the projection operators
\begin{equation}
    \label{Eq:ProjOperatorsDisc}
    \begin{aligned}
        &\mathcal P^{n\ell}_H:H\to H^{n,\ell},&&v^\ell=\mathcal P^{n\ell}_H\varphi&\text{solves }\min_{w^\ell\in H^{n\ell}}{\|\varphi-w^\ell\|}_H\text{ for }\varphi\in H,\\
        &\mathcal P^{n\ell}_V:V\to H^{n,\ell},&&u^\ell=\mathcal P^{n\ell}_V\varphi&\text{ solves } \min_{w^\ell\in H^{n\ell}}{\|\varphi-w^\ell\|}_V\text{ for }\varphi\in V,\\
        &\mathcal Q^{n\ell}_H:H\to V^{n,\ell},&&v^\ell=\mathcal Q^{n\ell}_H\varphi&\text{solves }\min_{w^\ell\in V^{n\ell}}{\|\varphi-w^\ell\|}_H\text{ for }\varphi\in H,\\
        &\mathcal Q^{n\ell}_V:V\to V^{n,\ell},&&u^\ell=\mathcal Q^{n\ell}_V\varphi&\text{solves } \min_{w^\ell\in V^{n\ell}}{\|\varphi-w^\ell\|}_V\text{ for }\varphi\in V.
    \end{aligned}
\end{equation}

The next result can be shown similarly as Proposotion~\ref{Lemma:ProjOper}.

\begin{proposition}
    \label{Lemma:ProjOperDisc}
    Let Assumption~{\em\ref{AssHV}} be satisfied. Suppose that the snapshots fulfill $y^k_j\in V$ for $k=1,\ldots,{K}$, $1\le j\le n$. The four projections in \eqref{Eq:ProjOperatorsDisc} are given as follows:
    \begin{equation}
        \label{Eq:Prolongations-12Disc}
        \begin{aligned}
            \mathcal P^{n\ell}_H\varphi&=\sum_{i=1}^\ell {\langle\varphi,\hat\psi_i^{nH}\rangle}_H\,\hat\psi_i^{nH}&&\hspace{-3mm}\text{for }\varphi\in H,&
            \mathcal P^{n\ell}_V\varphi&=\sum_{i=1}^\ell \mathrm v_i\hat\psi_i^{nH}&&\hspace{-3mm}\text{for }\varphi\in V,\\
            \mathcal Q^{n\ell}_V\varphi&=\sum_{i=1}^\ell {\langle\varphi,\hat\psi_i^{nV}\rangle}_V\,\hat\psi_i^{nV}&&\hspace{-3mm}\text{for }\varphi\in V,&\mathcal Q^{n\ell}_H\varphi&=\sum_{i=1}^\ell \mathrm u_i\hat\psi_i^{nV}&&\hspace{-3mm}\text{for }\varphi\in H,
        \end{aligned}
    \end{equation}
    where the vectors $\mathrm v=(\mathrm v_1,\ldots,\mathrm v_\ell)^\top,\,\mathrm u=(\mathrm u_1,\ldots,\mathrm u_\ell)^\top\in\mathbb R^\ell$ solve the following linear systems
    \begin{subequations}
        \label{SIAM:Eq3.2.14Disc}
        \begin{align}
            \label{SIAM:Eq3.2.14bDisc}
            \sum_{j=1}^\ell{\langle\hat\psi_j^{nH},\hat\psi_i^{nH}\rangle}_V\,\mathrm v_j={\langle\varphi,\hat\psi_i^{nH}\rangle}_V&&\text{for }1\le i\le \ell,\\
            \label{SIAM:Eq3.2.14aDisc}
            \sum_{j=1}^\ell{\langle\hat\psi_j^{nV},\hat\psi_i^{nV}\rangle}_H\,\mathrm u_j={\langle\varphi,\hat\psi_i^{nV}\rangle}_H&&\text{for }1\le i\le \ell,
        \end{align}
    \end{subequations}
    respectively. All four projections are linear and bounded. Furthermore, $\mathcal P_H^{n\ell},\,\mathcal Q_H^{n\ell}$ are orthogonal in $H$ and $\mathcal P_V^{n\ell},\,\mathcal Q_V^{n\ell}$ are orthogonal in  $V$.
\end{proposition}

\subsubsection{A-priori error analysis}

We also transfer Theorem~\ref{Prop:VTopology} to the discrete case.

\begin{theorem}
    \label{Prop:VTopologyDisc}
    Suppose that Assumption~{\em\ref{AssHV}} and $y^k_j\in V$ for $k=1,\ldots,{K}$, $1\le j\le n$, hold. We have for all $1\le\ell\le d^n$
    \begin{subequations}
        \label{Prop:Rate-VDisc}
        \begin{equation}
            \label{RatePH-VDisc}
            \sum_{k=1}^{K}\omega_k^{K}\sum_{j=1}^n\alpha_j^n\,\big\|y^k_j-\mathcal P^{n\ell}_Hy^k_j\big\|_V^2=\sum_{i=\ell+1}^{d^n}\hat\lambda_i^{nH}\,{\|\hat\psi_i^{nH}\|}_V^2
        \end{equation}
        For all $1\le\ell\le d^n$ it follows that
        \begin{equation}
            \label{RatePV-VDisc}
            \sum_{k=1}^{K}\omega_k^{K}\sum_{j=1}^n\alpha_j^n\,\big\|y^k_j-\mathcal P^{n\ell}_Vy^k_j\big\|_H^2=\sum_{i=\ell+1}^{d^n}\hat\lambda_i^{nH}\big\|\hat\psi_i^{nH}-\mathcal P^{n\ell}_V\hat\psi_i^{nH}\big\|_V^2.
        \end{equation}
        We further have
        \begin{equation}
            \label{RateQH-VDisc}
            \sum_{k=1}^{K}\omega_k^{K}\sum_{j=1}^n\alpha_j^n\,\big\|y^k_j-\mathcal Q^{n\ell}_Hy^k_j\big\|_V^2=\sum_{i=\ell+1}^{d^n}\hat\lambda_i^{nV}\,{\|\hat\psi_i^{nV}-\mathcal Q_H^{n\ell}\hat\psi_i^{nV}\|}_V^2
        \end{equation}
        for all $1\le\ell\le d^n$. Finally, we have for all $1\le\ell\le d^n$
        \begin{equation}
            \label{RateQV-VDisc}
            \sum_{k=1}^{K}\omega_k^{K}\sum_{j=1}^n\alpha_j^n\,\big\|y^k_j-\mathcal      Q^{n\ell}_Vy^k_j\big\|_V^2=\sum_{i=\ell+1}^{d^n}\hat\lambda_i^{nV}.
        \end{equation}
    \end{subequations}
\end{theorem}

\begin{example}
    \label{ex:embeddedEqualities}
    \rm We continue with Example \ref{Example:PODParabolic} to illustrate the validity of the equalities in \eqref{Prop:Rate-VDisc}. For this purpose, we set $d_{\mathsf{max}} = 50$ and compute for $\ell = 1,\hdots,d_{\mathsf{max}}$ the values 
    \begin{align*}
    	\Delta_V(\mathcal P_H^{n\ell}) &= \Big| \sum_{k=1}^{K}\omega_k^{K}\sum_{j=1}^n\alpha_j^n\,\big\|y^k_j-\mathcal P^{n\ell}_Hy^k_j\big\|_V^2-\sum_{i=\ell+1}^{d_{\mathsf{max}}}\hat\lambda_i^{nH}\,{\|\hat\psi_i^{nH}\|}_V^2\Big|,\\
    	\Delta_H(\mathcal P_V^{n\ell}) &= \Big| \sum_{k=1}^{K}\omega_k^{K}\sum_{j=1}^n\alpha_j^n\,\big\|y^k_j-\mathcal P^{n\ell}_Vy^k_j\big\|_H^2-\sum_{i=\ell+1}^{d_{\mathsf{max}}}\hat\lambda_i^{nH}\big\|\hat\psi_i^{nH}-\mathcal P^{n\ell}_V\hat\psi_i^{nH}\big\|_V^2 \Big|,\\
    	\Delta_V(\mathcal Q^{n\ell}_H) &= \Big| \sum_{k=1}^{K}\omega_k^{K}\sum_{j=1}^n\alpha_j^n\,\big\|y^k_j-\mathcal Q^{n\ell}_Hy^k_j\big\|_V^2-\sum_{i=\ell+1}^{d_{\mathsf{max}}}\hat\lambda_i^{nV}\,{\|\hat\psi_i^{nV}-\mathcal Q_H^{n\ell}\hat\psi_i^{nV}\|}_V^2 \Big|,\\
    	\Delta_V(\mathcal Q^{n\ell}_V) &= \Big| \sum_{k=1}^{K}\omega_k^{K}\sum_{j=1}^n\alpha_j^n\,\big\|y^k_j-\mathcal Q^{n\ell}_Vy^k_j\big\|_V^2=\sum_{i=\ell+1}^{d_{\mathsf{max}}}\hat\lambda_i^{nV} \Big|.
    \end{align*}
    According to Proposition \ref{Prop:VTopologyDisc}, all these terms should vanish for ${d_{\mathsf{max}}}=d^n$.\hfill$\blacklozenge$
\end{example}

\section{Proofs of Section~\ref{SIAM:Section-2}}
\label{SIAM:Section-2.6}
\setcounter{equation}{0}
\setcounter{theorem}{0}
\setcounter{figure}{0}
\setcounter{run}{0}

\subsection{Proofs of Section~\ref{SIAM:Section-2.1.1}}
\label{SIAM:Section-2.6.1}

\noindent{\bf\em Proof of Lemma~{\em\ref{SIAM:Lemma-I.1.1.1}}.} Since the inner product $\langle\cdot\,,\cdot\rangle_X$ is linear in the first argument, it follows that $\mathcal R^n$ is a linear operator. From
\begin{align*}
    {\|\mathcal R^n\psi\|}_X\le\sum_{k=1}^{K}\omega_k^{K}\sum_{j=1}^n\alpha_j^n \big|{\langle\psi,y_j^k\rangle}_X\big|\,{\|y_j^k\|}_X\quad\text{for }\psi\in X
\end{align*}
and the Cauchy-Schwarz inequality (see Lemma~\ref{App:CSIneq}) we conclude that $\mathcal R^n$ is bounded, i.e., $\mathcal R^n\in\mathscr L(X)$ holds true. Since $\mathcal R^n\psi \in\mathscr V^n$ holds for all $\psi\in X$, the range of $\mathcal R^n$ is finite dimensional. Thus, $\mathcal R^n$ is a finite rank operator which is compact; see \cite[p.~199]{RS80}. Next we show that $\mathcal R^n$ is self-adjoint. For that purpose we choose arbitrary elements $\psi,\tilde\psi\in X$ and consider
\begin{align*}
    {\langle \mathcal R^n\psi,\tilde\psi\rangle}_X&=\sum_{k=1}^{K}\omega_k^{K}\sum_{j=1}^n\alpha_j^n\,{\langle\psi,y_j^k\rangle}_X\,{\langle y_j^k,\tilde\psi\rangle}_X=\sum_{k=1}^{K}\omega_k^{K}\sum_{j=1}^n{\langle\psi,\alpha_j^n\,\overline{\langle y_j^k,\tilde\psi\rangle}_Xy_j^k\rangle}_X\\
    &=\bigg\langle\psi,\sum_{k=1}^{K}\omega_k^{K}\sum_{j=1}^n\alpha_j^n\,{\langle\tilde\psi,y_j^k\rangle}_Xy_j^k\bigg\rangle_X={\langle\psi,\mathcal R^n\tilde\psi\rangle}_X,
\end{align*}
where we have used that the weights $\alpha_j^n$ are real-valued for $1\leq j \leq n$. Thus, $\mathcal R^n$ is self-adjoint. For any $\psi \in X$ we derive
\begin{equation}
    \label{SIAM:Eq-I.1.1.16}
    \begin{aligned}
        {\langle \mathcal R^n\psi,\psi\rangle}_X&=\sum_{k=1}^{K}\omega_k^{K}\sum_{j=1}^n\alpha_j^n\,{\langle\psi,y_j^k\rangle}_X\,{\langle y_j^k,\psi\rangle}_X=\sum_{k=1}^{K}\omega_k^{K}\sum_{j=1}^n\alpha_j^n\,\big|{\langle y_j^k,\psi\rangle}_X\big|^2\ge 0.
    \end{aligned}
\end{equation}
Hence, $\mathcal R^n$ is also non-negative.\hfill$\Box$

\bigskip\noindent{\bf\em Proof of Lemma~{\em\ref{Lem:PropRnOp}}.}
\begin{enumerate}
    \item [1)] We derive \eqref{SIAM:Eq-I.1.1.19} directly from \eqref{SIAM:Eq-I.1.1.16}, \eqref{SIAM:Eq-I.1.1.18} and $\|\hat\psi_i^n\|_X=1$. Since $\hat\lambda_i^n=0$ holds for $i>d^n$, we find \eqref{SIAM:Eq-I.1.1.20}.
    \item [2)] Since $\{\hat\psi_i^n\}_{i\in\mathbb I}$ is a complete orthonormal basis and $\|y_j^k\|_X<\infty$ holds for $1\le k\le{K}$, $1\le j\le n$, we have
    \begin{equation}
        \label{SIAM:Eq-I.1.1.6}
        {\|y_j^k\|}_X^2=\bigg\langle\sum_{i\in\mathbb I}{\langle y_j^k,\psi_i\rangle}_X\,\psi_i,\sum_{\nu\in\mathbb I}{\langle y_j^k,\psi_\nu\rangle}_X\,\psi_\nu\bigg\rangle_X=\sum_{i\in\mathbb I}\big|{\langle y_j^k,\psi_i\rangle}_X\big|^2
    \end{equation}
    for $1\le k\le{K}$ and $1\le j\le n$. Thus, we derive from \eqref{SIAM:Eq-I.1.1.6}, \eqref{SIAM:Eq-I.1.1.19} and \eqref{SIAM:Eq-I.1.1.20} that
    \begin{align*}
        &\sum_{k=1}^{K}\omega_k^{K}\sum_{j=1}^n\alpha_j^n\,{\|y_j^k\|}_X^2=\sum_{k=1}^{K}\omega_k^{K}\sum_{j=1}^n\alpha_j^n\sum_{i\in\mathbb I}\big|{\langle y_j^k,\hat\psi_i^n\rangle}_X\big|^2\\
        &\quad=\sum_{i\in\mathbb I}\sum_{k=1}^{K}\omega_k^{K}\sum_{j=1}^n\alpha_j^n\,\big|{\langle y_j^k,\hat\psi_i^n\rangle}_X\big|^2=\sum_{i\in\mathbb I}\hat\lambda_i^n=\sum_{i=1}^{d^n}\hat\lambda_i^n.
    \end{align*}
    \item [3)] Note that \eqref{SIAM:Eq-I.1.1.18} implies that $\hat\psi_i^n \in \mathscr V^n$ for $i \le d^n$ and from \eqref{SIAM:Eq-I.1.1.20} we can conclude that $\hat\psi_i^n$, $i>d^n$, belongs to the orthogonal complement $(\mathscr V^n)^\perp$.
    \item [4)] By \eqref{SIAM:Eq-I.1.1.6} the (probably infinite) sum $\sum_{i\in\mathbb I}\hat\lambda_i^n$ is bounded. Notice that
    \begin{equation}
        \label{SIAM:Eq-I.1.1.3}
        \begin{aligned}
            \Big\| y_j^k-\sum_{i=1}^\ell{\langle y_j^k,\psi_i\rangle}_X\,\psi_i\Big\|_X^2&=\Big\langle y_j^k-\sum_{i=1}^\ell{\langle y_j^k,\psi_i\rangle}_X\,\psi_i,y_j^k-\sum_{l=1}^\ell{\langle y_j^k,\psi_l\rangle}_X\,\psi_l\Big\rangle_X\\
            &\quad={\|y_j^k\|}_X^2-\sum_{i=1}^\ell\left({\langle y_j^k,\psi_i\rangle}_X{\langle \psi_i,y_j^k\rangle}_X+\overline{{\langle y_j^k,\psi_i\rangle}}_X{\langle y_j^k,\psi_i\rangle}_X\right)\\
            &\qquad+\sum_{i=1}^\ell\sum_{l=1}^\ell{\langle y_j^k,\psi_i\rangle}_X\overline{{\langle y_j^k,\psi_l\rangle}}_X{\langle\psi_i,\psi_l\rangle}_X\\
            &\quad={\|y_j^k\|}_X^2-\sum_{i=1}^\ell\big|{\langle y_j^k,\psi_i\rangle}_X\big|^2
        \end{aligned}
    \end{equation}
    holds for any set $\{\psi_i\}_{i=1}^\ell\subset X$ satisfying $\langle\psi_i,\psi_l \rangle_X=\delta_{il}$ for $1 \leq i,l \leq \ell$. Combining \eqref{SIAM:Eq-I.1.1.22} and \eqref{SIAM:Eq-I.1.1.3}, the objective of \eqref{SIAM:Eq-I.1.1.2} can be written as stated in \eqref{SIAM:Eq-I.1.1.23}.\hfill$\Box$
\end{enumerate}

\bigskip\noindent{\bf\em Proof of Theoren~{\em\ref{SIAM:Theorem-I.1.1.3}}.} First, let $\{\psi_1,\hdots,\psi_\ell\}\subset X$ be an arbitrary orthonormal set in $X$. In terms of the orthonormal basis $\{\hat\psi_i^n\}_{i \in \mathbb I}$ we get
\begin{align*}
    \psi_i = \sum_{\mu \in \mathbb I} \langle \psi_i, \hat\psi_\mu^n \rangle_X \hat\psi_\mu^n \quad \text{for } i=1,\ldots,\ell.
\end{align*}
Inserting these identities into the objective function of \eqref{SIAM:Eq-I.1.1.4} yields
\begin{align*}
    \mathscr F(\psi_1,\hdots,\psi_\ell) :&= \sum_{k=1}^{K} \omega_k^{K} \sum_{j=1}^n \alpha_j^n \sum_{i=1}^\ell \left |{\langle y_j^k,\psi_i \rangle}_X \right |^2 \\
	&= \sum_{k=1}^{K} \omega_k^{K} \sum_{j=1}^n \alpha_j^n \sum_{i=1}^\ell \sum_{\mu,\nu \in \mathbb I}{\langle y_j^k, \hat\psi_\mu^n \rangle}_X \overline{{\langle y_j^k, \hat\psi_\nu^n \rangle}_X} {\langle \psi_i,\hat\psi_\nu^n \rangle}_X \overline{{\langle \psi_i, \hat\psi_\mu^n \rangle}_X}  \\
	&= \sum_{\mu,\nu \in \mathbb I} \bigg \langle \sum_{k=1}^{K} \omega_k^{K} \sum_{j=1}^n \alpha_j^n\,{\langle \hat\psi_\nu^n, y_j^k \rangle}_X y_j^k, \hat\psi_\mu^n \bigg \rangle_X \sum_{i=1}^\ell{\langle \psi_i, \hat\psi_\nu^n \rangle}_X \overline{{\langle \psi_i, \hat\psi_\mu^n \rangle}_X} \\
	&=\sum_{\mu,\nu \in \mathbb I}{\langle \mathcal R^n \hat\psi_\nu^n, \hat\psi_\mu^n \rangle}_X \sum_{i=1}^\ell{\langle \psi_i, \hat\psi_\nu^n \rangle}_X \overline{{\langle \psi_i, \hat\psi_\mu^n \rangle}_X}.
\end{align*}
We use the fact that the vector family $\{\hat\psi_\nu^n\}_{\nu \in \mathbb I}$ is a series of orthonormal eigenvectors of $\mathcal R^n$ to the eigenvalues $\{\hat\lambda^n_\nu\}_{\nu \in \mathbb I}$ and further deduce 
\begin{align}
	\label{eq:POD_proof1}
	\begin{aligned}
		\mathscr F(\psi_1,\hdots,\psi_\ell) &= \sum_{\mu,\nu \in \mathbb I} \hat\lambda_\nu^n \langle \hat\psi_\nu^n, \hat\psi_\mu^n \rangle_X \sum_{i=1}^\ell \langle \psi_i, \hat\psi_\nu^n \rangle_X \overline{ \langle \psi_i, \hat\psi_\mu^n \rangle}_X= \sum_{\mu \in \mathbb I} \hat\lambda_\mu^n \sum_{i=1}^\ell \big | \langle \psi_i,\hat\psi_\mu^n \rangle_X \big |^2.
	\end{aligned} 
\end{align}
Next, we add a zero term by utilizing $1 = \| \psi_i \|_X^2 = \sum_{\mu \in \mathbb I} \big | \langle \psi_i,\hat\psi_\mu^n \rangle_X \big |^2$:
\begin{align*}
	\mathscr F(\psi_1,\hdots,\psi_\ell) &= \sum_{i=1}^\ell \Big( \sum_{\mu \in \mathbb I} \hat\lambda_\mu^n \big | \langle \psi_i,\hat\psi_\mu^n \rangle_X \big |^2 + \hat\lambda_\ell^n - \hat\lambda_\ell^n \sum_{\mu \in \mathbb I} \big | \langle \psi_i, \hat\psi_\mu^n \rangle_X \big |^2 \Big) \\
	&= \sum_{i=1}^\ell \Big( \hat\lambda_\ell^n + \sum_{\mu \in \mathbb I} (\hat\lambda_\mu^n - \hat\lambda_\ell^n) \big | \langle \psi_i,\hat\psi_\mu^n \rangle_X \big |^2 \Big)\le \sum_{i=1}^\ell \Big( \hat\lambda_\ell^n + \sum_{\mu=1}^\ell (\hat\lambda_\mu^n - \hat\lambda_\ell^n) \big | \langle \psi_i, \hat\psi_\mu^n \rangle_X \big |^2 \Big),
\end{align*}
where we have used that $\hat\lambda_\mu^n \le \hat\lambda_\ell^n$ for $\mu \ge \ell$. Rearranging the sum, we finally end up with
\begin{align}
	\mathscr F(\psi_1,\hdots,\psi_\ell) &= \sum_{\mu=1}^\ell \Big( \hat\lambda_\ell^n + (\hat\lambda_\mu^n - \hat\lambda_\ell^n) \sum_{i=1}^\ell \big | \langle \psi_i,\hat\psi_\mu^n \rangle_X \big |^2 \Big) \notag \\
	&\le \sum_{\mu=1}^\ell \Big( \hat\lambda_\ell^n + (\hat\lambda_\mu^n - \hat\lambda_\ell^n) \sum_{i \in \mathbb I} \big | \langle \psi_i,\hat\psi_\mu^n \rangle_X \big |^2 \Big) \notag \\
	&= \sum_{\mu=1}^\ell \Big( \hat\lambda_\ell^n + (\hat\lambda_\mu^n - \hat\lambda_\ell^n) \Big) = \sum_{\mu=1}^\ell \hat\lambda_\mu^n. \label{eq:POD_proof2}
\end{align}
We have now reached an upper bound to the cost function $\mathscr F$ for any orthonormal set $\{\psi_1,\hdots,\psi_\ell\}$. If we particularly insert the set $\{\hat\psi_1^n,\hdots,\hat\psi_\ell^n\}$ into \eqref{eq:POD_proof1} the objective value is
\begin{align*}
	\mathscr F(\hat\psi_1^n,\hdots,\hat \psi_\ell^n) &= \sum_{\mu \in \mathbb I} \hat\lambda_\mu^n \sum_{i=1}^\ell \big | \langle \hat\psi_i^n, \hat\psi_\mu^n \rangle_X \big |^2 = \sum_{\mu=1}^\ell \hat\lambda_\mu^n
\end{align*}
which coincides with the upper bound \eqref{eq:POD_proof2} and proves \eqref{SIAM:Eq-I.1.1.25}. Thus,  $\{\hat\psi_1^n,\hdots,\hat\psi_\ell^n\}$ is a global solution to \eqref{SIAM:Eq-I.1.1.4} and therefore also to \eqref{SIAM:Eq-I.1.1.2}. Finally, we have to prove \eqref{SIAM:Eq-I.1.1.24}. In order to do this, we use \eqref{SIAM:Eq-I.1.1.23} and \eqref{SIAM:Eq-I.1.1.25} and get
\begin{align*}
    \sum_{k=1}^{K}\omega_k^{K}\sum_{j=1}^n\alpha_j^n\Big\|y_j^k-\sum_{i=1}^\ell{\langle y_j^k,\hat\psi_i^n\rangle}_X\,\hat\psi_i^n\Big\|_X^2 & = \sum_{i=1}^{d^n}\hat\lambda_i^n-\sum_{k=1}^{K}\omega_k^{K}\sum_{j=1}^n\alpha_j^n\sum_{i=1}^\ell\big|{\langle y_j^k,\psi_i\rangle}_X\big|^2 \\
	& = \sum_{i=1}^{d^n}\hat\lambda_i^n - \sum_{i=1}^\ell \hat\lambda_i^n = \sum_{i=\ell+1}^{d^n}\hat\lambda_i^n
\end{align*} 
which proves \eqref{SIAM:Eq-I.1.1.24}.\hfill$\Box$

\begin{remark}
    \label{SIAM:Remark-I.1.1.2}
    \rm Theorem~\ref{SIAM:Theorem-I.1.1.3} can also be proved by using the theory of non-linear programming; see \cite{HLBR12,Vol01}, for instance. For that purpose we introduce the quadratic cost functional
    \begin{equation*}
        J(\psi_1,\ldots,\psi_\ell)=-\sum_{k=1}^{K}\omega_k^{K}\sum_{j=1}^n\alpha_j^n\sum_{i=1}^\ell\big|{\langle y_j^k,\psi_i\rangle}_X\big|^2\quad\text{for }\{\psi_i\}_{i=1}^\ell\subset X.
    \end{equation*}
    Moreover, let
    \begin{equation*}
        e(\psi_1,\ldots,\psi_\ell)=\big(\big({\langle\psi_i,\psi_j\rangle}_X-\delta_{ij}\big)\big)\in\mathbb C^{\ell\times\ell}\quad\text{for }\{\psi_i\}_{i=1}^\ell\subset X.
    \end{equation*}
    Note that $e(\psi_1,\ldots,\psi_\ell)=e(\psi_1,\ldots,\psi_\ell)^\mathsf H$ holds, i.e., $e(\psi_1,\ldots,\psi_\ell)$ is a Hermitian matrix. Now \eqref{SIAM:Eq-I.1.1.4} can be expressed as the equality constrained optimization problem
    \begin{equation}
        \label{SIAM:Eq-Lagra}
        \min J(\psi_1,\ldots,\psi_\ell)\quad \text{subject to}\quad\{\psi_i\}_{i=1}^\ell\subset X\text{ and }e(\psi_1,\ldots,\psi_\ell)=0.
    \end{equation}
    Endowing the matrix space $\mathbb C^{\ell\times\ell}$ by the inner product
    \begin{align*}
        {\langle A,B\rangle}_{\mathbb C^{\ell\times\ell}}=\sum_{i=1}^\ell\sum_{j=1}^\ell a_{ij}\overline b_{ij}\quad\text{for }A=\big(\big(a_{ij}\big)\big),\,B=\big(\big(b_{ij}\big)\big)\in\mathbb C^{\ell\times\ell},
    \end{align*}
    the Lagrange functional associated with \eqref{SIAM:Eq-Lagra} is given as
    \begin{align*}
        \mathcal L(\psi_1,\ldots,\psi_\ell,\Lambda)=J(\psi_1,\ldots,\psi_\ell)+{\langle e(\psi_1,\ldots,\psi_\ell),\Lambda\rangle}_{\mathbb C^{\ell\times\ell}}
    \end{align*}
    for $\{\psi_i\}_{i=1}^\ell\subset X\text{ and }$ and $\Lambda=((\lambda_{ij}))\in\mathbb C^{\ell\times\ell}$. Applying a Lagrangian framework \cite[Chapter~12]{NW06} it turns out that \eqref{SIAM:Eq-I.1.1.18} are first-order necessary optimality conditions for \eqref{SIAM:Eq-I.1.1.4}. Especially, $\lambda_{ii}=\hat\lambda_i^n$ and $\lambda_{ij}=0$ hold for $1\le i,j\le \ell$. For the special case $X=\mathbb R^m$ endowed with the Euclidean inner product, for ${K}=1$ and $\ell=1$ first- and second-order optimality conditions are discussed in Section~\ref{SIAM:Section-2.1.1.2a}.\hfill$\blacksquare$
\end{remark}

\noindent{\bf\em Proof of Lemma~{\em \ref{Lem:OpY}}.} 
\begin{enumerate}
    \item [1)] First we proof that $\mathcal Y^n$ is bounded. Using the Cauchy-Schwarz inequality and \eqref{SIAM:Eq-I.1.1.22} we infer that
    \begin{align*}
        {\|\mathcal Y^n\Phi\|}_X&\le\sum_{k=1}^{K}\omega_k^{K}\sum_{j=1}^n\alpha_j^n\,|\Phi_{kj}|{\|y_j^k\|}_X\le\bigg(\sum_{k=1}^{K}\omega_k^{K}\sum_{j=1}^n\alpha_j^n\,\big|\Phi_{kj}\big|^2\bigg)^{1/2}\bigg(\sum_{k=1}^{K}\omega_k^{K}\sum_{j=1}^n\alpha_j^n\,{\|y_j^k\|}_X^2\bigg)^{1/2}\\
        &=\bigg(\sum_{i=1}^{d^n}\hat\lambda_i^n\bigg)^{1/2}{\|\Phi\|}_{\mathbb C^{{K}\times n}}\quad\text{for every }\Phi=\big(\big(\Phi_{kj}\big)\big)\in\mathbb C^{{K}\times n}.
    \end{align*}
    Consequently, $\mathcal Y^n$ is bounded. The image space of $\mathcal Y^n$ is given by $\mathscr V^n$, where $\mathscr V^n$ has been defined in \eqref{SIAM:Eq-I.1.1.1}. Thus, $\mathcal Y^n$ is a finite rank operator and therefore compact; see \cite[p.~199]{RS80}, for example.
    \item [2)] Let $\psi \in X$ be arbitrary. Then it holds
    \begin{align*}
        {\langle\mathcal Y^{n,\star}\psi,\Phi\rangle}_{\mathbb C^{{K}\times n}}&={\langle\psi,\mathcal Y^n\Phi\rangle}_X=\sum_{k=1}^{K}\omega_k^{K}\sum_{j=1}^n\alpha_j^n\overline\Phi_{kj}\,{\langle\psi,y^k_j\rangle}_X= {\langle (({\langle\psi,y^k_j\rangle}_X)) , \Phi \rangle}_{\mathbb C^{{K}\times n}},
    \end{align*}
    which gives the claim.
    \item [3)] Clearly, the composition of two linear operators is linear, too. By Lemma~\ref{Lem:OpY} the operator $\mathcal Y^n$ is bounded, so that its adjoint $\mathcal Y^{n,\star}$ is bounded as well; see, e.g., \cite[p.~186]{RS80}. Hence, $\mathcal Y^n\mathcal Y^{n,\star}$ is also bounded. Since $\mathcal Y^n$ is compact (see Lemma~\ref{Lem:OpY}) and $\mathcal Y^{n,\star}$ is bounded, $\mathcal Y^n\mathcal Y^{n,\star}$ is compact; see, e.g., \cite[p.~200]{RS80}. From $(\mathcal Y^n\mathcal Y^{n,\star})^\star=\mathcal Y^n\mathcal Y^{n,\star}$ we derive that $\mathcal Y^n\mathcal Y^{n,\star}$ is self-adjoint. As we have
    \begin{align*}
        {\langle \mathcal Y^n\mathcal Y^{n,\star}\psi,\psi\rangle}_X={\langle \mathcal Y^{n,\star}\psi,\mathcal Y^{n,\star}\psi\rangle}_{\mathbb C^{{K}\times n}}={\|\mathcal Y^{n,\star}\psi\|}_{\mathbb C^{{K}\times n}}\ge0\quad\text{for every }\psi\in X,
    \end{align*}
    the operator $\mathcal Y^n\mathcal Y^{n,\star}$ is even non-negative. Finally, we get for any $\psi\in X$
    \begin{align*}
        \mathcal Y^n\mathcal Y^{n,\star}\psi=\mathcal Y^n\Big(\big({\langle\psi,y^k_j\rangle}_X\big)_{kj}\Big)=\sum_{k=1}^{K}\omega_k^{K}\sum_{j=1}^n\alpha_j^n\,{\langle\psi,y^k_j\rangle}_Xy_j^k=\mathcal R^n\psi.
    \end{align*}
    Thus, $\mathcal R^n=\mathcal Y^n\mathcal Y^{n,\star}$ holds. \hfill$\Box$
\end{enumerate}

\noindent{\bf\em Proof of Lemma~{\em \ref{Lem:SVD_Y}}.} We already know that $\{ \hat\psi_i^n \}_{i=1}^{d^n}$ are orthonormal in $X$. Utilizing \eqref{Eq:SVD-1}, $\mathcal R^n=\mathcal Y^n\mathcal Y^{n,\star}$ and \eqref{SIAM:Eq-I.1.1.18} we find
\begin{align*}
    {\langle\hat\bPhi_i^n,\hat\bPhi_j^n\rangle}_{\mathbb C^{{K}\times n}}&=\frac{1}{\hat\sigma_i^n\hat\sigma_j^n}\,{\langle\mathcal Y^{n,\star}\hat\psi_i^n,\mathcal Y^{n,\star}\hat\psi_j^n\rangle}_{\mathbb C^{{K}\times n}}=\frac{1}{\hat\sigma_i^n\hat\sigma_j^n}\,{\langle\hat\psi_i^n,\mathcal R^n\hat\psi_j^n\rangle}_X=\frac{\hat\lambda_j^n}{\hat\sigma_i^n\hat\sigma_j^n}\,{\langle\hat\psi_i^n,\hat\psi_j^n\rangle}_X\\
    &=\frac{\hat\lambda_j^n\delta_{ij}}{\hat\sigma_i^n\hat\sigma_j^n}\quad\text{for }i,j=1,\ldots,d^n.
\end{align*}
Thus, the elements $\{\hat\Phi_i\}_{i=1}^{d^n}$ are orthonormal in $\mathbb C^{{K}\times n}$. Again by \eqref{Eq:SVD-1}, $\mathcal R^n=\mathcal Y^n\mathcal Y^{n,\star}$ and \eqref{SIAM:Eq-I.1.1.18} we obtain
\begin{equation*}
    \mathcal Y^n\hat\bPhi^n_i=\frac{1}{\hat\sigma_i^n}\mathcal Y^n\mathcal Y^{n,\star}\hat\psi_i^n=\frac{1}{\hat\sigma_i^n}\mathcal R^n\hat\psi_i^n=\hat\sigma_i^n\hat\psi_i^n\quad\text{for }i=1,\ldots,d^n.
\end{equation*}
Finally, \eqref{Eq:SVD_total} follows directly from the definition of the elements $\{\hat\bPhi_i\}_{i=1}^{d^n}$ in \eqref{Eq:SVD-1}. From \eqref{Eq:SVD_total} and with the fact that $\{\hat\psi_1^n,...,\hat\psi_{d^n}^n\}$ forms an orthonormal basis of $\mathscr V^n$, we derive directly \eqref{Eq:Lem:SVD_Y}.\hfill$\Box$

\bigskip\noindent{\bf\em Proof of Lemma~{\em\ref{Lem:PointwiseErrEstPOD}}.} Let $(k,j)\in\{1,\ldots,{K}\}\times\{1,\ldots,n\}$ be chosen arbitrarily. We define $\bPhi=((\Phi_{\nu l}))\in\mathbb C^{K\times n}$ by
\begin{align*}
    \Phi_{\nu l}=\left\{
    \begin{aligned}
        &\frac{1}{\omega_k^{K}\alpha_j^n}&&\text{if }(\nu,l)=(k,j),\\
        &0&&\text{otherwise}
    \end{aligned}
    \right.
\end{align*} 
for $1\le \nu\le K$ and $1\le l\le n$. Since $\{\hat\bPhi_i^n\}_{i=1}^{{K} n}$ forms an orthonormal basis of $\mathbb C^{K\times n}$, we have the representation
\begin{align*}
    \bPhi=\sum_{i=1}^{{K} n}{\langle\bPhi,\hat\bPhi_i^n\rangle}_{\mathbb C^{{K}\times n}}\hat\bPhi_i^n.
\end{align*}
Therefore, we obtain
\begin{align*}
    \frac{1}{\omega_k^{K}\alpha_j^n}=\Phi_{kj}&= \left( \sum_{i=1}^{{K} n}{\langle\bPhi,\hat\bPhi_i^n\rangle}_{\mathbb C^{{K}\times n}}\hat\bPhi_i^n \right)_{kj}=\sum_{i=1}^{{K} n}\bigg(\sum_{\nu=1}^{K}\omega_\nu^{K}\sum_{l=1}^n\alpha_l^n\Phi_{\nu l}\overline{(\hat\bPhi_i^n)}_{\nu l}\bigg)\big(\hat\bPhi^n_i\big)_{kj}\\
    &=\sum_{i=1}^{{K} n}\bigg(\sum_{\nu=1}^{K}\omega_\nu^{K}\sum_{l=1}^n\alpha_l^n\,\frac{\delta_{k\nu}\delta_{jl}}{\omega_k^{K}\alpha_j^n}\,\overline{(\hat\bPhi_i^n)}_{\nu l}\bigg)\big(\hat\bPhi^n_i\big)_{kj}=\sum_{i=1}^{{K} n}\big|\big(\hat\bPhi^n_i\big)_{kj}\big|^2.
\end{align*}
Consequently,
\begin{align*}
    \omega_k^{K}\alpha_j^n\sum_{i=1}^{{K} n}\big|\big(\hat\bPhi^n_i\big)_{kj}\big|^2=1
\end{align*}
for every $(k,j)\in\{1,\ldots,{K}\}\times\{1,\ldots,n\}$. This implies
\begin{align*}
    \Big\|y_j^k-\sum_{i=1}^\ell{\langle y_j^k,\hat\psi_i^n\rangle}_X\,\hat\psi_i^n\Big\|_X^2&=\Big\|\sum_{i=\ell+1}^{d^n}\hat\sigma_i^n\overline{\big(\hat\bPhi_i^n\big)}_{kj}\,\hat\psi_i^n\Big\|_X^2=\sum_{i=\ell+1}^{d^n}(\hat\sigma_i^n)^2\big|\big(\hat\bPhi_i^n\big)_{kj}\big|^2\le(\hat\sigma_{\ell+1}^n)^2\sum_{i=\ell+1}^{d^n}\big|\big(\hat\bPhi_i^n\big)_{kj}\big|^2\\
    &\le(\hat\sigma_{\ell+1}^n)^2\sum_{i=1}^{{K} n}\big|\big(\hat\bPhi_i^n\big)_{kj}\big|^2=\frac{(\hat\sigma_{\ell+1}^n)^2}{\omega_k^{K}\alpha_j^n},
\end{align*}
which implies \eqref{Eq:SVD2norm}.\hfill$\Box$

\bigskip\noindent{\bf\em Proof of Corollary~{\em \ref{SIAM:Corollary-I.1.1.1}}.}
\begin{enumerate}
    \item [1)] {\em Optimality of the POD basis}: We have
    \begin{align*}
        {\|Y-\hat\Psi^\ell\hat B^\ell\|}_F^2 & = \sum_{j=1}^n \sum_{k=1}^m \Big| y_{kj} - \sum_{i=1}^\ell {\langle y_j,\hat\psi_i\rangle}_{\mathbb C^m} \hat\psi_{ki} \Big|^2= \sum_{j=1}^n \| y_j - \sum_{i=1}^\ell {\langle y_j,\hat\psi_i\rangle}_{\mathbb C^m} \hat\psi_i \|_{\mathbb{C}^m}^2,
    \end{align*}
    and with the same argument
    \begin{align*}
        {\|Y-\Psi^\ell B^\ell\|}_F^2 & = \sum_{j=1}^n \| y_j - \sum_{i=1}^\ell {\langle y_j,\psi_i\rangle}_{\mathbb C^m} \psi_i \|_{\mathbb{C}^m}^2.
    \end{align*}
    Since $\{\hat\psi_i^n\}_{i=1}^\ell$ solves \eqref{SIAM:Eq-I.1.1.41}, we immediately get
    \begin{align*}
        {\|Y-\hat\Psi^\ell\hat B^\ell\|}_F \leq {\|Y-\Psi^\ell B^\ell\|}_F,
    \end{align*}
    which was the claim.
    \item [2)] The claim follows from \eqref{SIAM:Eq-I.1.1.37} and $\langle\hat\psi_i^n,\hat\psi_k^n\rangle_X=\delta_{ik}$ for $1 \le i,k \le \ell$:
    \begin{align*}
        \sum_{j=1}^n\overline{\langle y_j,\hat\psi_i^n\rangle}_{\mathbb C^m}{\langle y_j,\hat\psi_k^n\rangle}_{\mathbb C^m}&=\bigg\langle\sum_{j=1}^n{\langle\hat\psi_i^n,y_j\rangle}_{\mathbb C^m}\,y_j,\hat\psi_k^n\bigg\rangle_{\mathbb C^m}={\langle YY^\mathsf H\hat\psi_i^n,\hat\psi_k^n\big\rangle}_{\mathbb C^m}={\langle\hat\lambda_i^n\hat\psi_i^n,\hat\psi_k^n\big\rangle}_{\mathbb C^m}\\
        &=\hat\lambda_i^n\delta_{ik} \quad \text{for } 1\le i,k\le \ell.
    \end{align*}
    \hfill\endproof
\end{enumerate}

\subsection{Proofs of Section~\ref{SIAM:Section-2.1.2}}
\label{SIAM:Section-2.6.2}

\noindent{\bf\em Proof of Lemma~{\em\ref{Lemma2.2.0}}.} 
\begin{enumerate}
    \item [1)] Utilizing the Cauchy-Schwarz inequality and $y^k\in L^2(\mathscr D;X)$ for $1 \le k\le {K}$ we infer that
    \begin{align*}
        {\|\mathcal Y \phi\|}_X&\le\sum_{k=1}^{K}\omega_k^{K}\int_\mathscr D\big|\phi^k(\bmu)\big| {\|y^k(\bmu)\|}_X\,\mathrm d\bmu\le \sum_{k=1}^{K}\omega_k^{K}\,{\|\phi^k\|}_{L^2(\mathscr D;\mathbb{C})}{\|y^k\|}_{L^2(\mathscr D;X)}\\
        &\le\bigg(\sum_{k=1}^{K} \,{\|\phi^k\|}_{L^2(\mathscr D;\mathbb C)}^2\bigg)^{1/2}\bigg(\sum_{k=1}^{K}\omega_k^{K}\,{\|  y^k(\bmu)\|}_X^2\bigg)^{1/2}=C_\mathcal Y\,{\|\phi\|}_{L^2(\mathscr D;\mathbb C^{K})}
    \end{align*}
    for any $\phi\in L^2(\mathscr D;\mathbb C^{K})$, where we set $C_\mathcal Y=(\sum_{k=1}^{K}\omega_k^{K}\,{\|y^k(\bmu)\|}_X^2)^{1/2}<\infty$. Hence, $\mathcal Y$ is bounded.
    \item [2)] The Hilbert space adjoint $\mathcal Y^\star:X \to L^2(\mathscr D;\mathbb C^{K})$ of $\mathcal Y$ satisfies
    \begin{align*}
        {\langle \mathcal Y^\star\psi,\phi\rangle}_{L^2(\mathscr D;\mathbb C^{K})}={\langle\psi,\mathcal Y\phi\rangle}_X\quad\text{for }\psi\in X\text{ and }\phi\in L^2(\mathscr D;\mathbb C^{K}).
    \end{align*}
    Since we derive
    \begin{align*}
        {\langle \mathcal Y^\star\psi,\phi\rangle}_{L^2(\mathscr D;\mathbb C^{K})}&={\langle\psi,\mathcal Y\phi\rangle}_X=\bigg\langle\psi,\sum_{k=1}^{K}\omega^K_k\int_\mathscr D\phi^k(\bmu)y^k(\bmu)\,\mathrm d\bmu\bigg\rangle_X\\
        &=\sum_{k=1}^{K}\omega^K_k\int_\mathscr D{\langle\psi,y^k(\bmu)\rangle}_X\overline{\phi^k(\bmu)}\,\mathrm d\bmu=\left\langle{\big({\langle\psi,y^k(\cdot)\rangle}_X\big)}_{1\le k\le {K}},\phi\right\rangle_{L^2(\mathscr D;\mathbb C^{K})}
    \end{align*}
    for $\psi\in X$ and $\phi\in L^2(\mathscr D;\mathbb C^{K})$, the adjoint operator is given by \eqref{Eq2.2.10}. Next we prove that $\mathcal Y^\star$ is compact. Let $\{\chi_n\}_{n\in\mathbb N}\subset X$ be a sequence converging weakly to an element $\chi\in X$, i.e.,
    \begin{align*}
        \lim_{n\to\infty}{\langle\chi_n,\psi\rangle}_X={\langle\chi,\psi\rangle}_X\quad\text{for all }\psi\in X.
    \end{align*}
    This implies that
    \begin{align*}
        \lim_{n\to\infty}(\mathcal Y^\star \chi_n)(t)= \lim_{n\to\infty} \left(
        \begin{array}{c}
            {\langle\chi_n,y^1(\bmu)\rangle}_X\\
            \vdots\\
            {\langle\chi_n,y^{K}(\bmu)\rangle}_X\\
        \end{array}
        \right)=\left(
        \begin{array}{c}
            {\langle\chi,y^1(\bmu)\rangle}_X\\
            \vdots\\
            {\langle\chi,y^{K}(\bmu)\rangle}_X\\
        \end{array}
        \right)=\big(\mathcal Y^\star\chi\big)(t)
    \end{align*}
    for $\bmu\in\mathscr D$ a.e. Thus, the sequence $\{\mathcal Y^\star\chi_n\}_{n\in\mathbb N}$ converges weakly to $\mathcal Y^\star\chi$ in $L^2(\mathscr D;\mathbb C^{K})$. \hfill$\Box$
\end{enumerate}

\medskip\noindent{\bf\em Proof of Corollary~{\em\ref{Lemma2.2.1}}.}
\begin{enumerate}
    \item [1)] From
    \begin{align*}
        \big(\mathcal Y\mathcal Y^\star\big)\psi=\mathcal Y\left(
        \begin{array}{c}
            {\langle\psi,y^1(\cdot)\rangle}_X\\
            \vdots\\
            {\langle\psi,y^{K}(\cdot)\rangle}_X\\
        \end{array}
        \right)=\sum_{k=1}^{K}\omega_k^{K}\int_\mathscr D{\langle\psi,y^k(\bmu)\rangle}_Xy^k(\bmu)\,\mathrm d\bmu\quad\text{for }\psi\in X
    \end{align*}
    and \eqref{SIAM:Eq-I.1.2.5} we infer that $\mathcal R=\mathcal Y\mathcal Y^\star$ holds. Due to Lemma~\ref{Lemma2.2.0} the operator $\mathcal Y$ is bounded and its adjoint $\mathcal Y^\star$ even compact. Consequently, $\mathcal R=\mathcal Y\mathcal Y^\star$ is compact. We have $\mathcal R^\star=(\mathcal Y\mathcal Y^\star)^\star=\mathcal R$, i.e., the operator $\mathcal R$ is self-adjoint. From
    \begin{align*}
        {\langle \mathcal R \psi,\psi \rangle}_X&=\bigg\langle\sum_{k=1}^{K}\omega_k^{K}\int_\mathscr D{\langle \psi,y^k(\bmu)\rangle}_X \, y^k(\bmu)\,\mathrm d\bmu,\psi\bigg\rangle_X\\
        &= \sum_{k=1}^{K}\omega_k^{K}\int_\mathscr D\big|{\langle\psi,y^k(\bmu)\rangle}_X\big|^2\,\mathrm d\bmu\ge 0\quad\text{for all }\psi \in X
    \end{align*}
    we infer that $\mathcal R$ is non-negative.
    \item [2)] It follows from the proof of Lemma~\ref{Lemma2.2.1} that $\mathcal K=\mathcal Y^\star\mathcal Y$ is compact. However, the compactness of $\mathcal K$ can also be shown as follows: Notice that the kernel function
    \begin{align*}
        r_{ik}: \mathscr D \times \mathscr D \to \mathbb C, \quad r_{ik}(\bnu,\bmu) = {\langle y^k(\bnu),y^i(\bmu)\rangle}_X \quad \text{for } 1\le i,k\le{K},
    \end{align*}
    belongs to $L^2(\mathscr D \times \mathscr D; \mathbb C)$. Then it follows from \cite[pp. 197 and 277]{Yos95} that the linear integral operator $\mathcal K_{ik}:L^2(\mathscr D;\mathbb C) \to L^2(\mathscr D;\mathbb C)$ defined by
    \begin{align*}
        (\mathcal K_{ik} \phi)(\bmu)=\int_\mathscr Dr_{ik}(\bnu,\bmu)\phi(\bnu)\,\mathrm d\bnu \quad \text{for } \phi\in L^2(\mathscr D;\mathbb C),
    \end{align*}
    is a compact operator. This implies that the operator $\sum_{k=1}^{K}\mathcal K_{ik}$ is compact for $1\le i\le {K}$ as well.\hfill$\Box$
\end{enumerate}

\medskip\noindent{\bf\em Proof of Theorem~{\em \ref{Theorem2.2.1}}.} The existence of sequences $\{\hat\lambda_i\}_{i\in\mathbb I}$ of eigenvalues and $\{\hat\psi_i\}_{i\in\mathbb I}$ of associated eigenfunctions satisfying \eqref{Eq2.2.16} follows from Corrolary~\ref{Lemma2.2.1}-1), Theorem~\ref{SIAM:Theorem-I.1.1.1} and Theorem~\ref{SIAM:Theorem-I.1.1.2}. Analogous to the proof of Theorem~\ref{SIAM:Theorem-I.1.1.3} one can show that $\{\hat\psi_i\}_{i=1}^\ell$ solves \eqref{SIAM:Eq-I.1.2.3} as well as \eqref{SIAM:Eq-I.1.2.4} and that \eqref{Eq2.2.17} and \eqref{Eq2.2.18}, respectively, are valid.\hfill$\Box$
\endproof

\bigskip\noindent{\bf\em Proof of Lemma~{\em\ref{Lemma2.2.0}}.} 
\begin{enumerate}
    \item [1)] In fact,
    \begin{align*}
        \mathcal R\hat\psi_i=\sum_{k=1}^{K}\omega_k^K\int_\mathscr D{\langle\hat\psi_i,y^k(\bmu)\rangle}_X\,y^k(\bmu)\,\mathrm d\bmu\quad\text{for every }i\in\mathbb I.
    \end{align*}
    Taking the inner product with $\hat\psi_i$ and using \eqref{Eq2.2.16} we get
    \begin{align*}
        \sum_{k=1}^{K}\omega_k^K\int_\mathscr D\big|{\langle\hat\psi_i,y^k(\bmu)\rangle}_X\big|^2\,\mathrm d\bmu={\langle\mathcal R\hat\psi_i,\hat\psi_i\rangle}_X=\hat\lambda_i,
    \end{align*}
    which gives \eqref{SIAM:Eq-I.1.1.19-b}. Obviously, \eqref{SIAM:Eq-I.1.1.19} and $\hat\lambda_i=0$ for $i>d$ imply \eqref{SIAM:Eq-I.1.1.20-2}.
    \item [2)] Expanding each $y^k(\bmu) \in X$ in terms of $\{\hat\psi_i\}_{i\in\mathbb I}$ for each $1 \le k\le{K}$ we have
    \begin{align*}
        y^k(\bmu)=\sum_{i=1}^d{\langle y^k(\bmu),\hat\psi_i\rangle}_X\,\hat\psi_i
    \end{align*}
    and hence
    \begin{align*}
        \sum_{k=1}^{K}\int_\mathscr D{\|y^k(\bmu)\|}_X^2\,\mathrm d\bmu=\sum_{k=1}^{K}\sum_{i=1}^d\int_\mathscr D\big|{\langle y^k(\bmu),\hat\psi_i\rangle}_X\big|^2\,\mathrm d\bmu=\sum_{i=1}^d\hat\lambda_i,
    \end{align*}
    which is \eqref{Eq2.2.19}.\hfill$\Box$
\end{enumerate}

\subsection{Proofs of Section~\ref{Section:Perturbation analysis for the POD basis}}
\label{SIAM:Section-2.6.3}

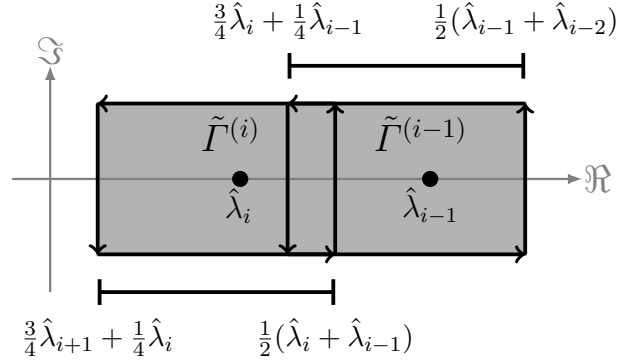
\begin{figure}
    \centering
    \begin{tikzpicture}
    \coordinate (reAxisMin) at (-3,0);
    \coordinate (reAxisMax) at (4.5,0);
    \coordinate (imAxisMin) at (-2.5,-1.5);
    \coordinate (imAxisMax) at (-2.5,1.5);
    \coordinate (aLeft) at (-1.875,0);
    \coordinate (aRight) at (1.25,0);
    \coordinate (pLeft) at (0.624,0);
    \coordinate (pRight) at (3.75,0);
    \coordinate (aLambda) at (0,0);
    \coordinate (pLambda) at (2.5,0);
    \filldraw [thick, fill=gray!60, even odd rule] ($(aLeft)-(0,1)$)  coordinate (GeneralStart) -- ($(aRight)-(0,1)$) -- ($(aRight)+(0,1)$) -- ($(aLeft)+(0,1)$) -- cycle;
    \filldraw [thick, fill=gray!60, even odd rule] ($(pLeft)-(0,1)$)  coordinate (GeneralStart) -- ($(pRight)-(0,1)$) -- ($(pRight)+(0,1)$) -- ($(pLeft)+(0,1)$) -- cycle;
    \draw [line width=0.3mm, gray,-latex] (reAxisMin) -- (reAxisMax);
    \draw [line width=0.3mm, gray,-latex] (imAxisMin) -- (imAxisMax);
    \draw[->,line width=0.5mm] ($(aLeft)-(0,1)$) -- ($(aRight)-(0,1)$);
    \draw[->,line width=0.5mm] ($(aRight)-(0,1)$) -- ($(aRight)+(0,1)$);
    \draw[->,line width=0.5mm] ($(aRight)+(0,1)$) -- ($(aLeft)+(0,1)$);
    \draw[->,line width=0.5mm] ($(aLeft)+(0,1)$) -- ($(aLeft)-(0,1)$);
    \draw[->,line width=0.5mm] ($(pLeft)-(0,1)$) -- ($(pRight)-(0,1)$);
    \draw[->,line width=0.5mm] ($(pRight)-(0,1)$) -- ($(pRight)+(0,1)$);
    \draw[->,line width=0.5mm] ($(pRight)+(0,1)$) -- ($(pLeft)+(0,1)$);
    \draw[->,line width=0.5mm] ($(pLeft)+(0,1)$) -- ($(pLeft)-(0,1)$);
    \draw[|-|,line width=0.5mm] ($(aLeft)-(0,1.5)$) -- ($(aRight)-(0,1.5)$);
    \draw[|-|,line width=0.5mm] ($(pLeft)+(0,1.5)$) -- ($(pRight)+(0,1.5)$);
    \node[draw,circle,inner sep=2pt,fill] at (aLambda) {};
    \node[draw,circle,inner sep=2pt,fill] at (pLambda) {};
    \node at ($(aLambda)+(-0.02,-0.35)$) {\large $\hat\lambda_i$};
    \node at ($(pLambda)+(0.0,-0.35)$) {\large $\hat\lambda_{i-1}$};
    \node at ($(aLambda)+(-0.125,0.6)$) {\Large \textbf{$\tilde\Gamma^{(i)}$}};
    \node at ($(pLambda)+(-0.125,0.6)$) {\Large \textbf{$\tilde\Gamma^{(i-1)}$}};
    \node at ($(aLeft)-(0,2.1)$) {$\tfrac{3}{4} \hat\lambda_{i+1} + \tfrac{1}{4}\hat\lambda_i$};
    \node at ($(aRight)-(0,2.1)$) {$\tfrac{1}{2} (\hat\lambda_{i} + \hat\lambda_{i-1})$};
    \node at ($(pLeft)+(0,2.1)$) {$\tfrac{3}{4} \hat\lambda_{i} + \tfrac{1}{4}\hat\lambda_{i-1}$};
    \node at ($(pRight)+(0,2.1)$) {$\tfrac{1}{2} (\hat\lambda_{i-1} + \hat\lambda_{i-2})$};
    \node[gray] at ($(reAxisMax)+(0.2,0.0)$) {\Large $\Re$};
    \node[gray] at ($(imAxisMax)+(0.0,0.2)$) {\Large $\Im$};
    \end{tikzpicture}
    \caption{Choice of the curves $\tilde\Gamma^{(i)}$ and $\tilde\Gamma^{(i-1)}$ in the proof of Theorem \ref{Theorem2.3.1}.}
    \label{fig:curveChoice}
\end{figure}

\noindent{\bf\em Proof of Theorem~{\em\ref{Theorem2.3.1}}.} We employ a proof by induction over $i$, starting at $i=1$. Our first step is to isolate the eigenvalue $\hat\lambda_1>0$ from the rest of the spectrum by a closed curve. Since all occuring eigenvalues are real, we will always choose the separating curves as the boundaries of rectangles in the complex plane
\begin{align*}
    R_{[a,b]} = [a,b] + i[-1,1] = \left \{ z \in \mathbb C ~|~ \Re(z) \in [a,b], \Im(z) \in [-1,1] \right \} \quad \text{for } a,b \in \mathbb R
\end{align*}
and denote the curve by $\Gamma_{[a,b]} = \partial R_{[a,b]}$. For a first curve, let us choose $\Gamma^{(1)} := \Gamma_{[\hat\lambda_1-\varepsilon,\hat\lambda_1+\varepsilon]}$. By the ordering of eigenvalues, we have
\begin{align*}
    \hat\lambda_1 = \hdots = \hat\lambda_r > \hat\lambda_{r+1},
\end{align*}
where $r \in \mathbb N$ is the multiplicity of $\hat\lambda_1$. Therefore, if $\varepsilon$ is chosen small enough, $\hat\lambda_1$ is the only spectral point of $\mathcal R$ in $\Gamma^{(1)}$. As $\mathcal R^n \to \mathcal R$ as $n \to \infty$, we can apply Theorem \ref{thm:continuityOfEigenvalues} to $\mathcal T = \mathcal R$, $\mathcal S = \mathcal R^n$, $\Sigma = \{\hat\lambda_1\}$ and $\Sigma' = \{\hat\lambda_i\}_{i > r}$ for large enough $n$. This means that for such $n$, the spectrum of $\mathcal R^n$ is also split in two parts $\Sigma(\mathcal R^n)$ inside of $\Gamma^{(1)}$ and $\Sigma'(\mathcal R^n)$ outside of $\Gamma^{(1)}$. In the according decompositions $X = M(\mathcal R) \oplus M'(\mathcal R)$ and $X = M(\mathcal R^n) \oplus M'(\mathcal R^n)$ from Theorem \ref{thm:spectralDecomposition}, property 2) states that $M(\mathcal R)$ and $M(\mathcal R^n)$ are isomorphic, so in particular of same dimension. Since $\mathcal R$ is compact and self-adjoint, we may employ Corollary \ref{cor:eigenvalueDecomposition} which states that $M(\mathcal R)$ is given by the eigenspace to $\hat\lambda_1$. The operator $\mathcal R^n$ is likewise compact and self-adjoint, so by Remark~\ref{rem:eigenvalueDecomposition} $\Sigma(\mathcal R^n)$ consists of at most $r$ different eigenvalues $\hat\lambda^n_{j_1^n},\hdots,\hat\lambda^n_{j_r^n}$ of $\mathcal R^n$ contained in $\Gamma^{(1)}$ and the $r$-dimensional space $M(\mathcal R^n)$ consists of eigenvectors to these eigenvalues. If we now send $\varepsilon \to 0$, it follows
\begin{align*}
    \hat\lambda_{j_1^n}^n,\ldots,\hat\lambda_{j_r^n}^n \to \hat\lambda_1\quad\text{as }n \to \infty.
\end{align*}
We would like to show that these eigenvalues can only be given by $\hat\lambda_1^n,\hdots,\hat\lambda_r^n$. To do this, we consider the fixed curve 
\begin{align*}
    \tilde\Gamma^{(1)} = \Gamma_{[ \nicefrac{3\hat\lambda_2}{4}  + \nicefrac{\hat\lambda_1}{4} , \hat\lambda_1+1]}.
\end{align*}
Obviously, $\hat\lambda_1$ is still the only eigenvalue of $\mathcal R$ enclosed by $\tilde\Gamma^{(1)}$ and $\tilde\Gamma^{(1)}$ encloses $\Gamma^{(1)}$ for small enough $\varepsilon$. If $n$ is large enough, Theorem \ref{thm:continuityOfEigenvalues} again states that there are at most $r$ different eigenvalues of $\mathcal R^n$ inside $\tilde\Gamma^{(1)}$ which can only be given by $\hat\lambda_{j_1^n}^n,\hdots,\hat\lambda_{j_r^n}^n$. Also, for large enough $n$, $\mathcal R^n$ has no more eigenvalues larger than $\hat\lambda_1+1$ because of $\| \mathcal R^n-\mathcal R \|_{\mathscr L(X)} \to 0$. So the $r$ largest eigenvalues (counting with multiplicity) of $\mathcal R^n$ are enclosed by $\tilde\Gamma^{(1)}$. In other words, $\hat\lambda_{j_1^n}^n,\hdots,\hat\lambda_{j_r^n}^n$ are given by $\hat\lambda_1^n,\hdots,\hat\lambda_r^n$ for large enough $n$. This proves $\hat\lambda_1^n,\hdots,\hat\lambda_r^n \to \hat\lambda_1$. Also, $\hat\psi_1,\hdots,\hat\psi_r$ is an orthonormal basis of $M(\mathcal R)$ and $\hat\psi_1^n,\hdots,\hat\psi_r^n$ is an orthonormal basis of $M(\mathcal R^n)$ if $n$ is large enough. Therefore, the projections $P(\mathcal R)$ and $P(\mathcal R^n)$ onto $M(\mathcal R)$ and $M(\mathcal R^n)$ from Theorem \ref{thm:continuityOfEigenvalues}-3) are given by
\begin{align*}
    P(\mathcal R)\psi = \sum_{j=1}^r \langle \psi,\hat\psi_j \rangle_X \hat\psi_j, \quad P(\mathcal R^n)\psi = \sum_{j=1}^r \langle \psi,\hat\psi_j^n \rangle_X \hat\psi_j^n \quad \text{for } \psi \in X.
\end{align*}
Because of this the convergence result \eqref{eq:perturbationPsi} for $\ell=r$ follows directly from the fact that $P(\mathcal R^n) \to P(\mathcal R)$ by Theorem \ref{thm:continuityOfEigenvalues}. This concludes the induction base. \\
For the induction step, let us assume that we are given a certain index $i$ which shall satisfy $\hat\lambda_i < \hat\lambda_{i-1}$ without loss of generality. Like in the induction base, we have
\begin{align*}
    \hat\lambda_i = \hdots = \hat\lambda_{i+r-1} > \hat\lambda_{i+r}.
\end{align*}
For the sake of simplicity, we have used the same $r$ as in the induction base, even though they may be different. Again, let us choose a closed curve $\Gamma^{(i)} = \Gamma_{[\hat\lambda_i-\varepsilon,\hat\lambda_i+\varepsilon]}$ where $\varepsilon>0$ is chosen small enough such that $\hat\lambda_{i-1}$ and $\hat\lambda_{i+r}$ lie on the outside of $\Gamma^{(i)}$. Similar to the induction base, we argue that if $n$ is large enough, $\mathcal R^n$ has at most $r$ different eigenvalues $\hat\lambda_{j_1^n},\hdots,\hat\lambda_{j_r^n}^n$ inside $\Gamma^{(i)}$. By the induction hypothesis, we know that $\hat\lambda_1^n,\hdots,\hat\lambda_{i-1}^n$ converge to $\hat\lambda_1,\hdots,\hat\lambda_{i-1}$ which lie outside of $\Gamma^{(i)}$ by construction, so $\hat\lambda_{j_1^n}^n,\hdots,\hat\lambda_{j_r^n}^n$ have to be distinct from those if $n$ is large enough. Again by induction, we know that the left-hand boundary of $\tilde\Gamma^{(i-1)}$ was chosen as $\tfrac{3}{4} \hat\lambda_i + \tfrac{1}{4} \hat\lambda_{i-1}$ and that for large enough $n$, the only eigenvalues of $\mathcal R^n$ inside $\tilde\Gamma^{(i-1)}$ are $\hat\lambda_{i-q}^n,\hdots,\hat\lambda_{i-1}^n$ (if $q$ was the multiplicity of $\hat\lambda_{i-1}$).  If we now choose the new curve as
\begin{align*}
    \tilde\Gamma^{(i)} := \Gamma_{[\nicefrac{3\hat\lambda_{i+1}}{4}+\nicefrac{\hat\lambda_i}{4}, \nicefrac{(\hat\lambda_i + \hat\lambda_{i-1})}{2}]}
\end{align*}
we have achieved that the interiors of $\tilde\Gamma^{(i-1)}$ and $\tilde\Gamma^{(i)}$ overlap, cf. Figure~\ref{fig:curveChoice}. 
Because of this overlap, we can conclude that the united interiors of $\tilde\Gamma^{(i-1)}$ and $\tilde\Gamma^{(i)}$ contain exactly the $\mathcal R^n$-eigenvalues $\hat\lambda_{i-q}^n,\hdots,\hat\lambda_{i-1}^n$ and $\hat\lambda_{j_1^n}^n,\hdots,\hat\lambda_{j_r^n}^n$. In other words, $\hat\lambda_{j_1^n},\hdots,\hat\lambda_{j_r^n}^n$ can only be the largest $r$ eigenvalues which are smaller than $\hat\lambda_{i-1}^n$, i.e., $\hat\lambda_i^n,\hat\lambda_{i+1}^n,\hdots,\hat\lambda_{i+r-1}$. Altogether, we have shown that $\hat\lambda_{i}^n,\hdots,\hat\lambda_{i+r}^n \to \hat\lambda_i$ as $n \to \infty$. \\
Completely analogue to the induction base, it follows that
\begin{align*}
    \sum_{j=i}^{i+r} \langle \psi,\hat\psi_j^n \rangle_X \hat\psi_j^n ~\to~ \sum_{j=i}^{i+r} \langle \psi,\hat\psi_j \rangle_X \hat\psi_j \quad \text{for all } \psi \in X \text{ as } n \to \infty.
\end{align*}
Together with the induction hypothesis, \eqref{eq:perturbationPsi} follows for $\ell = i+r$. This concludes the induction part of the proof. We still have to show \eqref{eq:perturbationLambda}. \\
We first observe that $y^k \in C(\mathscr D;X) \subset L^1(\mathscr D;X)$ since $\mathscr D$ is compact in $\mathbb R^{\mathfrak p}$. This together with \eqref{eq:perturbationIntegration} implies:
\begin{equation}
	\label{Eq2.3.17}
	\sum_{k=1}^{K}\sum_{j=1}^n \alpha_j^k\,{\| y^k(\bmu_j) \|}_X^2 \to \sum_{k=1}^{K}\int_\mathscr D{\|y^k(\bmu)\|}_X^2\,\mathrm d\bmu\quad \text{as }n\to \infty.
\end{equation}
If we utilize \eqref{SIAM:Eq-I.1.1.22} for the left and \eqref{Eq2.2.19} for the right, we can formulate \eqref{Eq2.3.17} equivalently as 
\begin{equation}
    \label{Eq2.3.18}
    \sum_{i=1}^{d^n} \hat\lambda_i^n \to \sum_{i=1}^d \hat\lambda_i \quad \text{as } n \to \infty.
\end{equation}
Now choose and fix $\ell$ such that $\hat\lambda_\ell\neq\hat\lambda_{\ell+1}$ for which we have already shown that $\hat\lambda_i^n \to \hat\lambda_i$ for all $1 \le \ell$. This together with \eqref{Eq2.3.18} finally proves \eqref{eq:perturbationLambda}.\hfill$\Box$
	
\bigskip\noindent{\bf\em Proof of Proposition~{\em\ref{Pro2.3.1}}.}
For an arbitrary $\psi \in X$ with $\|\psi\|_X=1$ we define $F:\mathscr D\to X$ by
\begin{align*}
    F(\bmu)=\sum_{k=1}^{K}{\langle y^k(\bmu),\psi \rangle}_X\,y^k(\bmu)\quad\text{for }\bmu\in\mathscr D.
\end{align*}
It follows that
\begin{equation}
    \label{Eq2.3.5a}
    \mathcal R\psi=\int_\mathscr DF(\bmu)\,\mathrm d\bmu=\int_{\bmu_\mathsf a}^{\bmu_\mathsf b}F(\bmu)\,\mathrm d\bmu=\sum_{j=1}^{n-1}\int_{\bmu_j}^{\bmu_{j+1}}F(\bmu)\,\mathrm d\bmu.
\end{equation}
The trapezoidal approximation of $\mathcal R\psi$ reads
\begin{equation}
    \label{Eq2.3.5b}
    \mathcal R^n\psi=\sum_{j=1}^n\alpha_j^nF(\bmu_j).
\end{equation}
Next, we infer from $\|\psi\|_X=1$ that
\begin{equation}
    \label{Lemma2.3.6}
    {\|F(\bmu)\|}_X^2 \le \bigg(\sum_{k=1}^{K} {\| y^k(\bmu)\|}_X^2\bigg)^2.
\end{equation}
Now, we show that $F$ belongs to $H^1(\mathscr D;X)$ and its norm is bounded independently of $\psi$. By assumption, $y^k\in C(\mathscr D;X)$ holds for $1\le k\le {K}$. Furthermore, $[0,T]$ is bounded. Thus, $y^k\in L^4(\mathscr D;X)$ which yields
\begin{align*}
    {\|F\|}_{L^2(\mathscr D;X)}^2\le\sum_{k=1}^{K}{\| y^k\|}_{L^4(\mathscr D;X)}^4 =C_1.
\end{align*}
Moreover, it holds
\begin{align*}
    F_\bmu(\bmu)= \sum_{k=1}^{K}{\langle y^k_{\bmu}(\bmu),\psi\rangle}_X\,y^k(\bmu)+{\langle y^k(\bmu),\psi\rangle}_X\,y^k_\bmu(\bmu)\quad\text{f.a.a. }\bmu\in\mathscr D,
\end{align*}
where \index{f.a.a., for almost all}{\em f.a.a.} stands for {\em for almost all}. Since the set $\mathscr D$ is one-dimensional, we have $H^1(\mathscr D;X)\hookrightarrow C(\mathscr D;X)$; see \cite[p.~269]{Eva08}, for instance. so that we have $y^k\in C(\mathscr D;X)$. Thus, we derive
\begin{align*}
    {\|F_\bmu(\bmu)\|}_X\le 2\sum_{k=1}^{K}{\| y^k(\bmu)\|}_X{\| y^k_\bmu(\bmu)\|}_X
\end{align*}
which implies
\begin{align*}
    {\|F_\bmu\|}_{L^2(\mathscr D;X)}^2&\le4\int_\mathscr D \bigg(\sum_{k=1}^{K}{\| y^k(\bmu)\|}_X{\| y^k_\bmu(\bmu)\|}_X\bigg)^2\,\mathrm d\bmu\le 4 {K} \sum_{k=1}^{K}\sum_{\nu=1}^{K}\int_\mathscr D {\|y^k(\bmu)\|}_X^2{\|y^\nu_\bmu(\bmu)\|}_X^2\,\mathrm d\bmu\\
    &\le 4 {K} \sum_{k=1}^{K}{\|y^k\|}_{C(\mathscr D;X)}^2\sum_{\nu=1}^{K} {\|y^\nu_\bmu\|}_{L^2(\mathscr D;X)}^2=C_2,
\end{align*}
where we have utilized that $y^k \in C(\mathscr D;X)$ for $1 \le k \le {K}$. Consequently,
\begin{equation}
    \label{Eq2.3.10}
    {\|F\|}_{H^1(\mathscr D;X)}=\bigg({\|F\|}_{L^2(\mathscr D;X)}^2+{\|F_\bmu\|}_{L^2(\mathscr D;X)}^2\bigg)^{1/2}\le C_3
\end{equation}
with the constant $C_3=(C_1+C_2)^{1/2}$, which is independent of $\psi$. For the next equality, we utilize the fundamental theorem of calculus for the function $F$:
\begin{align*}
    \int_{\bmu_j}^\bmu F_\bnu(\bnu)\,\mathrm d\bnu=F(\bmu)-F(\bmu_j)\quad\text{and}\quad\int_{\bmu_{j+1}}^{\bmu}F_\bnu(\bnu)\,\mathrm d\bnu=F(\bmu)-F(\bmu_{j+1})
\end{align*}
for $\bmu\in[\bmu_j,\bmu_{j+1}]$ for $j=1,\ldots,n-1$. Thus, we derive for the integral on the $j$-th subinterval:
\begin{equation}
    \label{Eq2.3.11}
    \begin{aligned}
        \int_{\bmu_j}^{\bmu_{j+1}}F(\bmu)\,\mathrm d\bmu&=\frac{1}{2}\int_{\bmu_j}^{\bmu_{j+1}}\Big(F(\bmu_j)+\int_{\bmu_j}^\bmu F_\bnu(\bnu)\,\mathrm d\bnu\Big)+\Big(F(\bmu_{j+1})+\int_{\bmu_{j+1}}^\bmu F_\bnu(\bnu)\,\mathrm d\bnu\Big)\mathrm d\bmu\\
        &=\frac{\delta\bmu_{j+1}}{2}\,\big(F(\bmu_j)+F(\bmu_{j+1})\big)\\
        &\quad+\frac{1}{2}\,\int_{\bmu_j}^{\bmu_{j+1}}\Big(\int_{\bmu_j}^\bmu F_\bnu(\bnu)\,\mathrm d\bnu+\int_{\bmu_{j+1}}^\bmu F_\bnu(\bnu)\,\mathrm d\bnu\Big)\mathrm d\bmu.
    \end{aligned}
\end{equation}
Combining \eqref{Eq2.3.5a} and \eqref{Eq2.3.11} together with the definition for $\alpha_j^n$, we get
\begin{align*}
    \mathcal R\psi&=\sum_{j=1}^{n-1}\int_{\bmu_j}^{\bmu_{j+1}}F(\bmu)\,\mathrm d\bmu=\sum_{j=1}^n\alpha_j^n F(\bmu_j)+\frac{1}{2}\sum_{j=1}^{n-1}\int_{\bmu_j}^{\bmu_{j+1}}\Big(\int_{\bmu_j}^\bmu F_\bnu(\bnu)\,\mathrm d\bnu+\int_{\bmu_{j+1}}^\bmu F_\bnu(\bnu)\,\mathrm d\bnu\Big)\mathrm d\bmu
\end{align*}
This concludes the final representation of $\mathcal R \psi$. Let us mention that the discrete variant can also be represented by the function $F$ since $\mathcal R^n \psi = \sum_{j=1}^n \alpha_j F(\bmu_j)$. Inserting these terms into \eqref{Eq2.3.5b}, we obtain
\begin{align*}
    {\big\|\mathcal R^n\psi-\mathcal R\psi\big\|}_X&=\bigg\|\frac{1}{2}\sum_{j=1}^{n-1}\int_{\bmu_j}^{\bmu_{j+1}}\Big(\int_{\bmu_j}^\bmu F_\bnu(\bnu)\,\mathrm d\bnu+\int_{\bmu_{j+1}}^\bmu F_\bnu(\bnu)\,\mathrm d\bnu\Big)\mathrm d\bmu\bigg\|_X\\
    &\le\frac{1}{2}\,\sum_{j=1}^{n-1} \bigg\|\int_{\bmu_j}^{\bmu_{j+1}}\int_{\bmu_j}^\bmu F_\bnu(\bnu)\,\mathrm d\bnu\mathrm d\bmu\bigg\|_X+\frac{1}{2}\,\sum_{j=1}^{n-1}\bigg\|\int_{\bmu_j}^{\bmu_{j+1}}\int_{\bmu_{j+1}}^\bmu F_\bnu(\bnu)\,\mathrm d\bnu\mathrm d\bmu\bigg\|_X.
\end{align*}
From the Cauchy-Schwarz inequality (cf. Lemma~\ref{App:CSIneq}) we deduce that
\begin{equation}
    \label{Eq2.3.13} 
    \begin{aligned}
        \sum_{j=1}^{n-1}\bigg\|\int_{\bmu_j}^{\bmu_{j+1}}\int_{\bmu_j}^\bmu F_\bnu(\bnu)\,\mathrm d\bnu\mathrm d\bmu \bigg\|_X&\le \sum_{j=1}^{n-1}\int_{\bmu_j}^{\bmu_{j+1}}\bigg\|\int_{\bmu_j}^\bmu F_\bnu(\bnu)\,\mathrm d\bnu\bigg\|_X\mathrm d\bmu\\
        &\le\sqrt{\Delta\bmu}\sum_{j=1}^{n-1}\bigg(\int_{\bmu_j}^{\bmu_{j+1}}\Big\|\int_{\bmu_j}^\bmu F_\bnu(\bnu)\,\mathrm d\bnu\Big\|_X^2\mathrm d\bmu\bigg)^{1/2}\\
        &\le\sqrt{\Delta\bmu}\sum_{j=1}^{n-1}\bigg(\int_{\bmu_j}^{\bmu_{j+1}}\Big(\int_{\bmu_j}^\bmu {\|F_\bnu(\bnu)\|}_X\,\mathrm d\bnu\Big)^2\mathrm d\bmu\bigg)^{1/2}\\
        &\le\Delta\bmu\sum_{j=1}^{n-1}\bigg(\int_{\bmu_j}^{\bmu_{j+1}}\int_{\bmu_j}^\bmu {\|F_\bnu(\bnu)\|}_X^2\,\mathrm d\bnu\mathrm d\bmu\bigg)^{1/2}\\
        &\le\Delta\bmu\sum_{j=1}^{n-1}\bigg( (\mu_{j+1} - \mu_j) {\|F\|}_{H^1(\mu_j,\mu_{j+1};X)}^2  \bigg)^{1/2}\\
        &\le\Delta\bmu \bigg( \sum_{j=1}^{n-1} (\mu_{j+1} - \mu_j) \bigg)^{1/2} \bigg( \sum_{j=1}^{n-1} {\|F\|}_{H^1(\mu_j,\mu_{j+1};X)}^2 \bigg)^{1/2} \\
        &= \Delta\bmu \sqrt{\bmu_\mathsf b-\bmu_\mathsf a}\,{\|F\|}_{H^1(\mathscr D;X)}.
    \end{aligned}
\end{equation}
Analogously, we derive
\begin{equation}
    \label{Eq2.3.14} 
    \sum_{j=1}^{n-1}\bigg\|\int_{\bmu_j}^{\bmu_{j+1}}\int_{\bmu_{j+1}}^\bmu F_\bnu(\bnu)\,\mathrm d\bnu\mathrm d\bmu\bigg\|_X\le \Delta\bmu \sqrt{\bmu_\mathsf b-\bmu_\mathsf a}\,{\|F\|}_{H^1(\mathscr D;X)}.
\end{equation}
From \eqref{Eq2.3.10}, \eqref{Eq2.3.13} and \eqref{Eq2.3.14} it follows that
\begin{align*}
    {\big\|\mathcal R^n\psi-\mathcal R\psi\big\|}_X\le C_4 \Delta\bmu,
\end{align*}
where $C_4=C_3(\bmu_\mathsf b-\bmu_\mathsf a)^{1/2}$ is independent of $n$ and $\psi$. Since $\Delta \bmu \to 0$ as $n \to \infty$, we conclude that
\begin{align*}
    {\|\mathcal R^n-\mathcal R\|}_{\mathscr L(X)}=\sup_{\|\psi\|_X=1}{\|\mathcal R^n\psi-\mathcal R\psi\|}_X \stackrel{n\to\infty}{\longrightarrow} 0,
\end{align*}
which gives the claim.\hfill$\Box$

\subsection{Proofs of Section~\ref{Section:PODHilbert}}
\label{SIAM:Section-2.6.4}

\noindent{\bf\em Proof of Proposition~{\em\ref{Lemma:ProjOper}}.} Since $H^\ell$ is a finite-dimensional subspace, it is also closed in the Hilbert space $H$, so Theorem \ref{thm:projectionTheorem} yields $H = H^\ell \oplus (H^\ell)^\perp$. Therefore, every element $\varphi \in H$ decomposes into $\varphi = \varphi^\ell + \varphi^\perp$ with $\varphi^\ell \in H^\ell$ and $\varphi^\perp \in (H^\ell)^\perp$. Obviously, $\{\psi_1^H,\hdots,\psi_\ell^H\}$ is an $H$-orthonormal basis of $H^\ell$, so that we can write $\varphi^\ell = \sum_{i=1}^\ell \langle \varphi,\psi_i^H \rangle_H \psi_i^H$. Lastly, Theorem \ref{thm:projectionTheorem} also states that $\varphi^\ell$ solves the minimization problem in \eqref{SIAM:Eq3.2.12a}, so it holds $\mathcal P_H^\ell \varphi = \varphi^\ell$. Analogously, one can prove that $\mathcal P_V^\ell$, $\mathcal Q_H^\ell$ and $\mathcal Q^\ell_V$ are well-defined. To prove \eqref{SIAM:Eq3.2.14a} we express $\varphi^\ell$ through the basis $\{\psi_1^H,...,\psi_\ell^H\}$ as $\varphi^\ell = \sum_{i=1}^\ell\mathrm v_i\psi_i^H$. Taking the inner product of $\varphi$ with a basis element $\psi_i^H$ yields
\begin{align*}
    {\langle\varphi,\psi_i^H\rangle}_V={\langle\varphi^\ell,\psi_i^H\rangle}_V=\sum_{j=1}^\ell{\langle \psi_j^H,\psi_i^H\rangle}_V\,\mathrm u_j 
\end{align*}
which proves \eqref{SIAM:Eq3.2.14a}. The representation \eqref{SIAM:Eq3.2.14b} is shown in the same way. Finally, Remark~\ref{SIAM:Remark_PODBasis_LinIndep} yields the unique solvability of \eqref{SIAM:Eq3.2.14a} and \eqref{SIAM:Eq3.2.14b}. The linearity and boundedness follows from the fact that the operators are orthogonal projections.\hfill$\Box$


\chapter{Reduced-Order Modeling for Evolution Problems}
\label{SIAM-Book:Section3}

In this chapter we consider linear and non-linear evolution problems which can be considered as dynamical systems in appropriate function spaces. For the numerical solution, a Galerkin approximation is introduced which leads to finite-dimensional linear or non-linear equation systems. However, in many applications these systems turn out to be of very large scale. For this reason reduced-order modeling is utilized to derive low-dimensional Galerkin schemes for the equation systems. In this book we apply the POD method in order to realize the model reduction. Of course, other methods, e.g., the reduced basis method can also be utilized. To ensure the validity of the so-called POD Galerkin approximation, a-priori and a-posteriori error estimates are derived. The chapter is organized as follows ...

\begin{itemize}
	\item In Section \ref{SIAM-Book:Section3.1}, we return to the parabolic example presented in Section \ref{SIAM-Book:Section1.2} involving already control input functions. The PDE is formulated as an abstract evolution equation using (bi-)linear forms. 
	\item This specific form of evolution equations is studied in Section \ref{SIAM-Book:Section3.2}. We present (standard) sufficient conditions ensuring existence of a unique solution for any fixed control input functions. This solution can be written as a sum of a particular solution, which depends on the given right-hand side and the initial condition, and of a control-dependent solution. The latter can be expressed as the application of a linear and bounded control-to-state operator to the chosen fixed control input function.
	\item In Section \ref{SIAM-Book:Section3.3}, we analyze the continuous POD method for the evolution problem, giving us POD vectors which are built by the solution trajectory. Using these POD vectors, we introduce a Galerkin scheme by which we project the evolution equation onto the POD subspace. In an a-priori analysis, we present error bounds between the so-called truth solution and the reduced POD solution.
	\item In Section \ref{SIAM-Book:Section3.4} the results from Sections \ref{SIAM-Book:Section3.2} and \ref{SIAM-Book:Section3.3} are transferred to a {\em semidiscrete approximation}, where the evolution problem is discretized by a Galerkin scheme. We present an {\em a-priori} and an {\em a-posteriori error analysis} to estimate the distance between the exact (unknown) soltution to the evolution problem to the associated solution to the semidiscrete approximation for any fixed control input function. In particular, the a-posteriori error estimates allow for a computation of the distance provided an approximate solution is available.
	\item In Section \ref{SIAM-Book:Section3.5}, the temporal discretization of the semidiscrete approximation is carried out. This leads to a fully discretized scheme for the evolution problem that can be realized on a computer. All results from the previous sections carry over in an analogous way, which is carried out in the section.
	\item Section \ref{SIAM-Book:Section3.6} introduces us to the more complex case when a non-linearity enters the evolution equation. Under certain assumptions, we can present a-priori and a-posteriori estimators which are similar to the linear situation from the first sections in the continuous, semidiscrete and fully discrete cases.
	\item In Section \ref{SIAM-Book:Section3.6}, we explain how POD can be applied to elliptic problems. An a-priori error analysis is also presented.
	\item To improve the readability of the chapter, all proofs have been moved to Section \ref{SIAM-Book:Section3.8}.
\end{itemize}

\section{Motivating linear evolution problem}
\label{SIAM-Book:Section3.1}
\setcounter{equation}{0}
\setcounter{theorem}{0}
\setcounter{algorithm}{0}
\setcounter{figure}{0}
\setcounter{run}{0}

Let us return to our general numerical example from Section \ref{SIAM-Book:Section1.2} in which we have encountered the \index{Equation!partial differential!parabolic}{\em parabolic partial differential equation} \eqref{SIAM:EqMotPDE1} which we will quickly repeat here:
\begin{equation}
    \label{Eq:MotEx}
    \begin{aligned}
    y_t(t,\bx)-\kappa \Delta y(t,\bx)+\bv(\bx) \cdot \nabla y(t,\bx)&=f(t,\bx)+\sum\limits_{i=1}^\mb u^{\mathsf d}_i(t)\chi_i(\bx),&&(t,\bx)\in Q,\\
    \kappa \frac{\partial y}{\partial\bn}(t,\bs)+q(\bs)y(t,\bs)&=g(t,\bs)+\sum\limits_{j=1}^\mb u^{\mathsf b}_j(t)\xi_j(\bs),&&(t,\bs)\in\Sigma,\\
    y(0,\bx)&=y_\circ(\bx),&&\bx\in\Omega.
    \end{aligned}
\end{equation}
We will now define what constitutes as a solution of \eqref{Eq:MotEx} by using appropriate function spaces. We assume that the initial heat distribution $y_\circ$ as well as the distributed and boundary heat sources $f,g$ are given bounded functions. The input or control variable $u=(u^{\mathsf d},u^{\mathsf b})$ belongs to the Hilbert space \index{Space!Hilbert!control, $\U$}$\U=L^2(0,T;\mathbb R^{m_\mathsf d} \times \mathbb R^{m_\mathsf b})$ endowed with the inner product
\begin{align*}
    {\langle u,\tilde u\rangle}_\U=\int_0^T\sum_{i=1}^\mb u_i^{\mathsf d}(t)\tilde u_i^{\mathsf d}(t)+\sum_{i=1}^\mb u_i^{\mathsf b}(t)\tilde u_i^{\mathsf b}(t)\,\mathrm dt\quad\text{for }u=(u^{\mathsf d},u^{\mathsf b}),~\tilde u=(\tilde u^{\mathsf d},\tilde u^{\mathsf b})\in\U
\end{align*}
and the induced norm $\|\cdot\|_\U=\langle\cdot\,,\cdot\rangle_\U^{1/2}$. Furthermore, the shape functions $\chi_1,\ldots,\chi_\md$, $\xi_1,\ldots,\xi_\mb$ are bounded. The advection term $\bv: \Omega \to \mathbb R^\mathfrak n$ is supposed to be bounded in each component $v_i$, $i=1,...,\mathfrak n$. 

We are interested in a weak solution $y$ to \eqref{Eq:MotEx}. Let $H=L^2(\Omega)$ and $V=H^1(\Omega)$. Then $H$ and $V$ are two separable Hilbert spaces with $V\subset H$; cf. Remark \ref{rem:gelfandTiple}. Moreover, $V$ is dense in $H$ with compact embedding. Identifying the Hilbert space $H$ with its dual space $H'$, we get the Gelfand triple
\begin{align*}
    V\hookrightarrow H\simeq H'\hookrightarrow V'.
\end{align*}
We assume
\begin{align*}
    y\in H^{1,1}(Q)=L^2(0,T;V)\cap H^1(0,T;H),
\end{align*}
where $H^{1,1}(Q)$ is endowed with the inner product
\begin{align*}
    {\langle\varphi,\phi\rangle}_{H^{1,1}(Q)}&={\langle\varphi,\phi\rangle}_{L^2(0,T;V)}+{\langle\varphi,\phi\rangle}_{H^1(0,T;H)}\\
    &=\int_0^T{\langle \varphi(t),\phi(t)\rangle}_V\,\mathrm dt+\int_0^T{\langle \varphi(t),\phi(t)\rangle}_H+{\langle\varphi_t(t),\phi_t(t)\rangle}_H\,\mathrm dt\quad\text{for }\varphi,\phi\in H^{1,1}(Q),
\end{align*}
and the associated induced norm
\begin{align*}
    {\|\varphi\|}_{H^{1,1}(Q)}=\left({\|\varphi\|}_{L^2(0,T;V)}+{\|\varphi\|}_{H^1(0,T;H)}\right)^{1/2}\quad\text{for }\varphi\in H^{1,1}(Q).
\end{align*}
Next, we multiply the first line in \eqref{Eq:MotEx} by a test function $\varphi\in V$ and utilize the integration by parts formula. From the boundary condition in \eqref{Eq:MotEx} it follows that
\begin{align*}
    &\int_\Omega f(t,\bx)\varphi(\bx)\,\mathrm d\bx + \sum_{i=1}^\md u_i^{\mathsf d}(t) \int_\Omega \chi_i(\bx)\,\mathrm d\bx \\= 
    &\int_\Omega y_t(t,\bx)\varphi(\bx)+\kappa \nabla y(t,\bx)\cdot\nabla\varphi(\bx)+\big(\bv(\bx) \cdot \nabla y(t,\bx)\big) \varphi(\bx)\,\mathrm d\bx-\kappa \int_\Gamma \frac{\partial y}{\partial \bn}(t,\bs)\varphi(\bs)\,\mathrm d\bs\\
    &=\int_\Omega y_t(t,\bx)\varphi(\bx)+\kappa \nabla y(t,\bx)\cdot\nabla\varphi(\bx)+\big(\bv(\bx) \cdot \nabla y(t,\bx)\big) \varphi(\bx)\,\mathrm d\bx\\
    &+\int_\Gamma q(\bs)y(t,\bs)\varphi(\bs) ~\mathrm d\bs - \int_\Gamma g(t,\bs) \varphi(\bs)\,\mathrm d\bs - \sum_{j=1}^\mb u^{\mathsf b}_j(t) \int_\Gamma \xi_j(\bs)\varphi(\bs)\,\mathrm d\bs
\end{align*}
for almost all $t\in(0,T]$. We define the bilinear form $a:V\times V\to\mathbb R$ by
\begin{align*}
    a(\phi,\varphi)=\kappa \int_\Omega \nabla \phi(\bx)\cdot\nabla\varphi(\bx)+\big(v(\bx)\cdot \nabla \phi(\bx) \big) \varphi(\bx)\,\mathrm d\bx+\int_\Gamma q(\bs)\phi(\bs)\varphi(\bs)\,\mathrm d\bs\quad\text{for }\phi,\varphi\in V.
\end{align*}
We prove that there exists two constant $\gamma_1>0$ and $\gamma_2\ge0$ such that
\begin{equation}
    \label{WeakCoerciveExample}
    a(\varphi,\varphi)\ge\gamma_1\,{\|\varphi\|}_V^2-\gamma_2\,{\|\varphi\|}_H^2\quad\text{for all }\varphi\in V.
\end{equation}
For that purpose we take an arbitrary $\varphi \in V$, set $v_\infty=\|\bv\|_{L^\infty(Q;\mathbb R^\mathfrak n)}$ and obtain
\begin{align*}
    a(\varphi,\varphi) &\ge \kappa \int_\Omega {|\nabla \varphi(\bx)|}_2^2+\big(\bv(\bx)\cdot\nabla\varphi(\bx)\big)\varphi(\bx)\,\mathrm d\bx\ge \kappa\big({\|\varphi\|}_V^2-{\|\varphi\|}_H^2\big)-v_\infty\,{\|\nabla\varphi\|}_{L^2(\Omega;\mathbb R^\mathfrak n)}{\|\varphi\|}_H,
\end{align*}
where the Bochner norms are recalled in Definition \ref{Definition_L2_BanachValued}. We continue to apply Young's inequality (cf. Lemma \ref{lem:youngsInequality}) to arrive at
\begin{align*}
    a(\varphi,\varphi) &\ge\kappa\big({\|\varphi\|}_V^2-{\|\varphi\|}_H^2\big)-v_\infty\,\bigg( \frac{\kappa}{2v_\infty}\,{\|\nabla \varphi\|}_{L^2(\Omega;\mathbb R^\mathfrak n)}^2+\frac{v_\infty}{2\kappa}\,{\|\varphi\|}_H^2\bigg)\ge\gamma_1\,{\|\varphi\|}_V^2-\gamma_2\,{\|\varphi\|}_H^2
\end{align*}
with $\gamma_1=\kappa/2>0$ and $\gamma_2=(\kappa^2+v_\infty^2)/(2\kappa)$. Therefore, the bilinear form $a(\cdot\,,\cdot)$ satisfies \eqref{WeakCoerciveExample}. Furthermore, we obtain with $q_\infty=\|q\|_{L^\infty(\Gamma)}$
\begin{align*}
    \big|a(\phi,\varphi)\big|&\le\kappa\int_\Omega\big|\nabla\phi(\bx)\cdot\nabla\varphi(\bx)\big|+\big|\big(\bv(\bx)\cdot\nabla\phi(\bx)\big)\varphi(\bx)\big|\,\mathrm d\bx+q_\infty\int_\Gamma \big|\phi(\bs)\varphi(\bs)\big|\,\mathrm d\bs\\
    &\le\kappa\,{\|\phi\|}_V{\|\varphi\|}_V+v_\infty\,{\|\phi\|}_V{\|\varphi\|}_V+q_\infty\,{\|\phi\|}_{L^2(\Gamma)}{\|\varphi\|}_{L^2(\Gamma)}
\end{align*}
for $\phi,\varphi\in V$. By the \index{Theorem!trace, $\gamma_\Gamma$}{\em trace theorem} (cf. \cite{Eva08}) we have
\begin{equation}
    \label{TraceTh}
    {\|\varphi\|}_{L^2(\Gamma)}\le \gamma_\Gamma\,{\|\varphi\|}_V\quad\text{for all }\varphi\in V.
\end{equation}
Thus, setting $\gamma=\kappa + \|b\|_\infty + \gamma_\Gamma^2\|q\|_\infty>0$, we have
\begin{align*}
    \big|a(\phi,\varphi)\big|\le\gamma\,{\|\phi\|}_V{\|\varphi\|}_V\quad\text{for all }\phi,\varphi\in V
\end{align*}
i.e., the bilinear form $a(\cdot\,,\cdot)$ is also {\em bounded}. 

For $t \in (0,T)$ the right-hand side $\mathcal F(t): V \to \mathbb R$ be defined as
\begin{align*}
    {\langle\mathcal F(t),\varphi\rangle}_{V',V}=\int_\Omega f(t,\bx) \varphi(\bx)\,\mathrm d\bx+\int_\Gamma g(t,\bs) \varphi(\bs)\,\mathrm d\bs \quad \text{for all } \varphi \in V.
\end{align*}
For almost all $t \in (0,T)$ and all $\varphi \in V$ we derive that
\begin{align*}
    \big|{\langle\mathcal F(t),\varphi\rangle}_{V',V}\big|\le{\|f(t,\cdot)\|}_H{\|\varphi\|}_H+\gamma_\Gamma\,{\|g(t,\cdot)\|}_{L^2(\Gamma)}{\|\varphi\|}_V.
\end{align*}
Since $f$ and $g$ are assumed to be essentially bounded and $\Omega$, $\Gamma$ have finite measure, the functions $f$ and $g$ are square integrable. The above inequality implies that it holds $\mathcal F \in L^2(0,T;V')$ with
\begin{align*}
    {\|\mathcal F\|}_{L^2(0,T;V')}\le{\|f\|}_{L^2(0,T;H)}+\gamma_\Gamma\,{\|g\|}_{L^2(0,T;L^2(\Gamma))}.
\end{align*}

Lastly, we also introduce for $u = (u^{\mathsf d}, u^{\mathsf b}) \in \U$ and $t \in (0,T)$ the linear \index{Linear operator!control, $\mathcal B$}{\em control operator}
\begin{align*}
    {\langle(\mathcal B u)(t),\varphi\rangle}_{V',V}= \sum_{i=1}^\mb u_i^\mathsf d(t)\int_\Omega \chi_i(\bx) \varphi(\bx)\,\mathrm d\bx + \sum_{j=1}^\mb u_i^{\mathsf b}(t) \int_\Gamma \xi(\bs) \varphi(\bs)\,\mathrm d\bs
\end{align*}
for all $\varphi \in V$. Clearly, the linear operator $\mathcal B:\mathscr U\to L^2(0,T;V')$ is linear. We prove that $\mathcal B$ is even continuous. For $t \in (0,T)$ and $\varphi \in V$ we get by applying the Cauchy-Schwarz inequality and using \eqref{Poincare} the estimate
\begin{align*}
    \big|{\langle (\mathcal B u)(t),\varphi\rangle}_{V',V} \big|&\le\sum_{i=1}^\mb | u^\mathsf d_i(t)| \int_\Omega\big|\chi_i(\bx)\varphi(\bx)\big|\,\mathrm d\bx+\sum_{j=1}^\mb| u^\mathsf b_j(t)|\int_\Gamma\big|\xi_j(\bs)\varphi(\bs)\big|\,\mathrm d\bs\\
    &\le\sum_{i=1}^{m_\mathsf d} | u^\mathsf d_i(t)|{\|\chi_i\|}_H{\|\varphi\|}_H+\sum_{j=1}^{m_\mathsf b}| u^\mathsf b_j(t)|{\|\xi_j\|}_{L^2(\Gamma)}{\|\varphi\|}_{L^2(\Gamma)}\\
    &\le{\|\varphi\|}_H\sum_{i=1}^\md | u^\mathsf d_i(t)|{\|\chi_i\|}_H+{\|\varphi\|}_{L^2(\Gamma)}\sum_{j=1}^\mb| u^\mathsf b_j(t)|{\|\xi_j\|}_{L^2(\Gamma)}\\
    &\le c_1c_V\,{\|\varphi\|}_V{\|u^\mathsf d\|}_{L^2(0,T;\mathbb R^{m_\mathsf d})}+c_2c_\Gamma\,{\|\varphi\|}_V{\|u^\mathsf b\|}_{L^2(0,T;\mathbb R^{m_\mathsf b})}\\
    &\le c\,{\|\varphi\|}_V\sqrt{2}\big({\|u^\mathsf d\|}_{L^2(0,T;\mathbb R^{m_\mathsf d})}^2+{\|u^\mathsf b\|}_{L^2(0,T;\mathbb R^{m_\mathsf b})}^2\big)^{1/2}=\gamma_\mathcal B\,{\|u\|}_\U{\|\varphi\|}_V
\end{align*}
with $\gamma_\mathcal B=c\sqrt{2}$, $c=\max(c_1,c_V,c_2c_\Gamma)$ and
\begin{align*}
    c_1=\bigg( \sum_{i=1}^\md{\|\chi_i\|}_H^2 \bigg)^{1/2},\quad c_2=\bigg(\sum_{j=1}^\mb{\|\xi_j\|}_{L^2(\Gamma)}^2\bigg)^{1/2} 
\end{align*}
This proves in particular that the linear operator $\mathcal B$ is bounded with
\begin{align*}
    {\|\mathcal B u\|}_{L^2(0,T;V')}\le\gamma_\mathcal B\,{\|u\|}_\U\quad\text{for all }u\in \U.
\end{align*}
Finally, the weak form of the evolution problem \eqref{Eq:MotEx} can be expressed in the form
\begin{subequations}
    \label{ExWeakForm}
    \begin{equation}
        \label{ExWeakForm-1}
        \frac{\mathrm d}{\mathrm dt}\,{\langle y(t),\varphi\rangle}_H+a(y(t),\varphi)={\langle (\mathcal F+\mathcal Bu)(t),\varphi\rangle}_{V',V}\text{ for }\varphi\in V\text{ a.e. in }(0,T]
    \end{equation}
    together with the initial condition
    \begin{equation}
        \label{ExWeakForm-2}
        {\langle y(0),\varphi\rangle}_H={\langle y_\circ,\varphi\rangle}_H\quad\text{for all }\varphi\in H.
    \end{equation}
\end{subequations}
A solution $y$ to \eqref{ExWeakForm} is called a {\em weak} or {\em variational solution} to \eqref{Eq:MotEx}.

Recall the linear space
\begin{align*}
    W(0,T)=\big\{\varphi\in L^2(0,T;V)\,\big|\,\varphi_t\in L^2(0,T;V')\big\}=L^2(0,T;V)\cap H^1(0,T;V')
\end{align*}
which is a Hilbert space endowed with the inner product
\begin{align*}
    {\langle\varphi,\phi\rangle}_{W(0,T)}=\int_0^T{\langle\varphi_t(t),\phi_t(t)\rangle}_{V'}+{\langle\varphi(t),\phi(t)\rangle}_V\,\mathrm dt\quad\text{for all }\varphi,\phi\in W(0,T)
\end{align*}
and the induced norm $\|\cdot\|_{W(0,T)}=\langle\cdot\,,\cdot\rangle_{W(0,T)}$. For more details we refer to Section~\ref{App:A.4} in the appendix. In particular, we have $H^{1,1}(Q)\hookrightarrow W(0,T)$. Due to Lemma~\ref{Lemma:HI-20}-4) we have for $y\in W(0,T)$
\begin{align*}
    \frac{\mathrm d}{\mathrm dt}\,{\langle y(t),\varphi\rangle}_H={\langle y_t(t),\varphi\rangle}_{V',V}\quad\text{for all }\varphi\in V.
\end{align*}
Thus, a solution to \eqref{ExWeakForm} do not need to be in $H^{1,1}(Q)$.

\section{The abstract linear evolution problem}
\label{SIAM-Book:Section3.2}
\setcounter{equation}{0}
\setcounter{theorem}{0}
\setcounter{figure}{0}
\setcounter{run}{0}

In this section, we introduce the abstract linear evolution problem which is mainly considered in this book. We discuss the unique solvability and finally present an a-priori error estimate for its solutions. 

\subsection{Problem formulation}
\label{SIAM-Book:Section3.2.1}

We make use of the following hypotheses.

\begin{assumption}
    \label{A1}
    Suppose that $T>0$ holds.
    \begin{enumerate}
        \item [\rm 1)] $V$ and $H$ are real, separable Hilbert spaces and $V$ is dense in $H$ with compact embedding. By $\langle\cdot\,,\cdot\rangle_H$ and $\langle\cdot\,,\cdot\rangle_V$ we denote the inner products in $H$ and $V$, respectively. We identify $H$ with its dual (Hilbert) space $H'$ by the Riesz isomorphism so that we have the \index{Gelfand triple}{\em Gelfand triple}
        \begin{align*}
            V\hookrightarrow H\simeq H'\hookrightarrow V',
        \end{align*}
        where each embedding is continous and dense. The last embedding is understood as follows: For every element $h' \in H'$ and $v \in V$, we also have $v \in H$ by the embedding $V \hookrightarrow H$, so we can define $\langle h',v \rangle_{V',V} = \langle h',v \rangle_{H',H}$. 
        \item [\rm 2)] The \index{Space!Hilbert!control, $\U$}{\em control space} $\U$ is a real Hilbert space with inner product $\langle\cdot\,,\cdot\rangle_\U$ and induced norm $\|\cdot\|_\U=(\langle\cdot\,,\cdot\rangle_\U)^{1/2}$.
        \item [\rm 3)] For almost all $t\in[0,T]$ we define a \index{Bilinear form!time-dependent, $a(t;\cdot\,,\cdot)$}time-dependent bilinear form $a(t;\cdot\,,\cdot):V\times V\to\mathbb R$ satisfying
        \begin{subequations}
            \label{SIAM:Eq3.1.1}
            \begin{align}
                \label{SIAM:Eq3.1.1-1}
                \big|a(t;\varphi,\psi)\big|&\le\gamma\,{\|\varphi\|}_V{\|\psi\|}_V&&\forall\varphi,\psi\in V \text{ a.e. in } [0,T],\\
                \label{SIAM:Eq3.1.1-2}
                a(t;\varphi,\varphi)&\ge\gamma_1\,{\|\varphi\|}_V^2-\gamma_2\,\|\varphi\|_H^2&&\forall\varphi\in V \text{ a.e. in } [0,T],
            \end{align}
        \end{subequations}
        for constants $\gamma,\gamma_2\ge0$ and $\gamma_1>0$ which do not depend on $t$. 
        \item [\rm 4)] Assume that $y_\circ\in H$, $\mathcal F\in L^2(0,T;V')$ and that $\mathcal B:\U\to L^2(0,T;V')$ is a continuous, linear \index{Linear operator!control, $\mathcal B$}{\em control operator}.
    \end{enumerate}
\end{assumption}

\begin{remark}
\label{Remark:HI-20}
\rm
\begin{enumerate}
    \item [\rm 1)] We associate with $a(t;\cdot\,,\cdot)$ a time-dependent linear and bounded operator $\mathcal A(t):V\to V'$ given by
    \begin{equation}
        \label{OperatorA(t)}
        {\langle\mathcal A(t)\varphi,\psi\rangle}_{V',V}=a(t;\varphi,\psi)\quad\text{for }\varphi,\psi\in V\text{ and }t\in[0,T].
    \end{equation}
    If the bilinear form \index{Bilinear form!time-independent, $a(\cdot\,,\cdot)$}$a$ is time-independent, i.e., $a:V\times V\to\mathbb R$, then the linear, bounded operator $\mathcal A$ introduced in \eqref{OperatorA(t)} is time-independent as well. In this case we define the \index{Domain!of a linear operator, $D(\mathcal A)$}\index{Linear operator!domain, $D(\mathcal A)$}{\em domain of $\mathcal A$} as the linear subspace
    \begin{equation}
        \label{DomainA}
        D(\mathcal A)=\big\{\varphi\in V\,\big|\,\mathcal A\varphi\in H\big\}.
    \end{equation}
    It follows that
    \begin{align*}
        D(\mathcal A) \hookrightarrow V \hookrightarrow H = H' \hookrightarrow V'
    \end{align*}
    holds, each embedding being continuous and dense, when $\mathcal A$ is bounded and $D(\mathcal A)$ is endowed with the \index{Linear operator!graph norm}{\em graph norm of $\mathcal A$}, given by 
    \begin{align*}
        {\|\varphi \|}_{\mathcal A}=\big({\|\varphi\|}_V^2+{\|\mathcal A\varphi\|}_H^2 \big)^{1/2}\quad\text{for all }\varphi\in D(\mathcal A).
    \end{align*}
    \item [\rm 2)] From Assumption~\ref{A1}-1) it follows that there exists an embedding constant $c_V>0$ satisfying \eqref{Poincare}.\hfill$\blacksquare$
    \end{enumerate}
\end{remark}

Let $\Y=W(0,T)$ denote the \index{Space!state, $\Y$}{\em state space}; cf. \eqref{EqW(0,T)}. For fixed control $u\in\U$ the state $y\in\Y$ is governed by the linear evolution problem
\begin{subequations}
    \label{SIAM:Eq3.1.6}
    \begin{align}
        \label{SIAM:Eq3.1.6a}
        \frac{\mathrm d}{\mathrm dt} {\langle y(t),\varphi \rangle}_H+a(t;y(t),\varphi)&={\langle (\mathcal F+\mathcal Bu)(t),\varphi\rangle}_{V',V}&&\hspace{-2mm}\forall\varphi\in V\text{ a.e. in }(0,T],\\
        \label{SIAM:Eq3.1.6b}
        {\langle y(0),\varphi\rangle}_H&={\langle y_\circ,\varphi\rangle}_H&&\hspace{-2mm}\forall\varphi\in H.
    \end{align}
\end{subequations}

\begin{remark}
    \label{SIAM:Remark3.1.1}
    \rm Utilizing the operator $\mathcal A(t)$ introduced in Remark~\ref{Remark:HI-20}-2) we can write \eqref{SIAM:Eq3.1.6} also as the following linear time-variant dynamical system
    \begin{align*}
        \dot y(t)+\mathcal A(t)y(t)&=(\mathcal F+\mathcal Bu)(t)&&\text{ in }V'\text{ a.e. in }(0,T],\\
        y(0)&=y_\circ && \text{ in } H
    \end{align*}
    with the notation $\dot y=y_t$.\hfill$\blacksquare$
\end{remark}

\subsection{Unique solvability and a-priori bound}
\label{SIAM-Book:Section3.2.2}

Existence of a unique solution and an a-priori estimate for the solution is stated in the following theorem which is proved in Section~\ref{SIAM-Book:Section3.8.1}.

\begin{theorem}
    \label{SIAM:Theorem3.1.1}
    Let Assumption~{\rm\ref{A1}} hold. Then for every $u\in\U$ there is a unique weak solution $y\in\Y$ satisfying \eqref{SIAM:Eq3.1.6} and
    \begin{equation}
    \label{SIAM:Eq3.1.7}
    {\|y\|}_\Y\le C\left({\|y_\circ\|}_H+{\|\mathcal F\|}_{L^2(0,T;V')}+{\|u\|}_\U\right)
    \end{equation}
    for a constant $C>0$ which is independent of $u$, $y_\circ$ and $f$. If $\mathcal F+\mathcal Bu\in L^2(0,T;H)$, $a(t;\cdot \,,\cdot)=a(\cdot\,,\cdot)$ (i.e., $a$ is independent of $t$) and $y_\circ\in V$ hold, we even have $y\in C([0,T];V) \cap H^1(0,T;H)\cap L^2(0,T;D(\mathcal A))$.
\end{theorem}

The following corollary follows directly from Theorem~\ref{SIAM:Theorem3.1.1}.

\begin{corollary}
    \label{Corollary:HI-20}
    Let Assumption~{\rm\ref{A1}} hold.
    \begin{enumerate}
        \item [\rm 1)] There exists a unique solution $\hat y\in\Y$ to
        \begin{equation}
            \label{yhatEq}
            \begin{aligned}
                \frac{\mathrm d}{\mathrm dt} {\langle\hat y(t),\varphi \rangle}_H+a(t;\hat y(t),\varphi)&={\langle \mathcal F(t),\varphi\rangle}_{V',V}&&\text{for all }\varphi\in V\text{ a.e. in }(0,T],\\
                {\langle\hat y(0),\varphi\rangle}_H&={\langle y_\circ,\varphi\rangle}_H&&\text{for all }\varphi\in H
            \end{aligned}
        \end{equation}
        satisfying
        \begin{align*}
            {\|\hat y\|}_\Y\le C\left({\|y_\circ\|}_H+{\|\mathcal F\|}_{L^2(0,T;V')}\right)
        \end{align*}
        for a constant $C>0$.
        \item [\rm 2)] Let us define the linear \index{Linear operator!solution, $\mathcal S$}{\em solution operator} $\mathcal S:\U\to\Y$ as follows: for every $u\in\U$ the function $y=\mathcal Su$ is the solution to
        \begin{equation}
            \label{eq:solutionOperator_a}
            \begin{aligned}
                \frac{\mathrm d}{\mathrm dt} {\langle y(t),\varphi \rangle}_H+a(t;y(t),\varphi)&=       {\langle (\mathcal Bu)(t),\varphi\rangle}_{V',V}&&\text{for all }\varphi\in V\text{ a.e. in }(0,T],\\
                {\langle y(0),\varphi\rangle}_H&=0&&\text{for all }\varphi\in H.
            \end{aligned}
        \end{equation}
        Then $\mathcal S$ is well-defined and bounded, i.e.,
        \begin{equation}
            \label{HI-100}
            {\|\mathcal Su\|}_\Y\le C\,{\|u\|}_\U
        \end{equation}
        for a constant $C>0$ which is independent of $u$. 
    \end{enumerate}
\end{corollary}

\begin{remark}
    \label{Remark:HI-22}
    \rm We can express the solution to \eqref{SIAM:Eq3.1.6} as $y=\hat y+\mathcal Su$. The operator $\mathcal S$ is often called {\em control-to-state operator} or {\em solution operator}. Notice that the (particular) solution $\hat y$ does not depend on the control input  $u$.\hfill$\blacksquare$
\end{remark}

It was already shown in Section \ref{SIAM-Book:Section3.1} that the weak formulation of the recurring parabolic example satisfies Assumption \ref{A1}. This yields the following corollary. 

\begin{corollary}
    \label{cor:numericalExample_solvability}
    The recurring parabolic partial differential equation \eqref{SIAM:EqMotPDE1} in its weak formulation \eqref{ExWeakForm} admits a unique solution $y \in W(0,T)$ which can be written as $y = \hat y + \mathcal S u$, where $\hat y \in W(0,T)$ is the solution to \eqref{yhatEq} and $\mathcal S: \U \to W(0,T)$ is introduced in Corollary~{\em\ref{Corollary:HI-20}}.
\end{corollary}

\begin{example}
    \label{Example:RegResult}
    \rm Let $\Omega\subset\mathbb R^\mathfrak n$, $\mathfrak n\in\{1,2,3\}$, be a bounded domain with Lipschitz-continuous boundary $\Gamma=\partial\Omega$. We set $H=L^2(\Omega)$ and $V=H^1_0(\Omega)$. The space $V$ is endowed with inner product
    \begin{align*}
        {\langle\varphi,\phi\rangle}_V=\int_\Omega\nabla\varphi\cdot\nabla\phi\,\mathrm d\bx\quad\text{for }\varphi,\phi\in V
    \end{align*}
    and the induced norm $\|\cdot\|_V=\langle\cdot\,,\cdot\rangle_V^{1/2}$. If the time-independent bilinear form $a(\cdot\,,\cdot)$ is given by
    \begin{align*}
        a(\varphi,\phi)={\langle\varphi,\phi\rangle}_V\quad\text{for }\varphi,\phi\in V,
    \end{align*}
    $\mathcal F+\mathcal Bu\in L^2(0,T;H$ holds for the given $u\in\U$ and $y_\circ\in V$ is valid. it is shown in \cite[pp.~532-533]{DL00} that the solution $y=\hat y+\mathcal Su$ belongs to the space $W(0,T)\cap L^2(0,T;H^2(\Omega))$. Since $H^2(\Omega)$ is compactly embedded in $V\hookrightarrow H$, we infer from Aubin's lemma \cite[p.~271]{Tem79} that $W(0,T)\cap L^2(0,T;H^2(\Omega))$ is compactly embedded in $L^2(0,T;V)$.\hfill$\blacklozenge$
\end{example}

\section{The continuous POD method}
\label{SIAM-Book:Section3.3}
\setcounter{equation}{0}
\setcounter{theorem}{0}
\setcounter{figure}{0}
\setcounter{run}{0}

In this section we apply a POD Galerkin approximation to the linear evolution problem \eqref{SIAM:Eq3.1.6}, but no discretization of the temporal variable. Therefore, we utilize the \index{POD method!continuous variant}{\em continuous variant of the POD method} introduced in Section~\ref{Section:ContPODHilbert}, where we distinguish two choices for $X$: $X=H$ and $X=V$. The parameter domain $\mathscr D$ is chosen to be the time interval $[0,T]$. Using the rate of convergence results presented in Section~\ref{Section:ContPODHilbert}, an a-priori error analysis is carried out for the difference between the (exact) solution to \eqref{SIAM:Eq3.1.6} and the one computed by the POD Galerkin scheme.

\subsection{POD Galerkin scheme}
\label{SIAM-Book:Section3.3.1}

We will utilize the following hypotheses.

\begin{assumption}
    \label{A2}
    \begin{enumerate}
        \item [\rm 1)] Assumption~{\rm\ref{A1}} holds.
        \item [\rm 2)] The Hilbert space $X$ denotes either $H$ or $V$. For the choice $X=V$ we additionally assume $y_\circ\in V$.
    \end{enumerate}
\end{assumption}

For given $\mathcal F\in L^2(0,T;V')$, $u\in\U$ and $y_\circ\in X$ let $y\in\Y$ be the unique solution to \eqref{SIAM:Eq3.1.6}. In the context of Section~\ref{Section:ContPODHilbert} we choose $K=1$, $\omega_1^K=1$ und $y^1=y\in L^2(0,T;X)$. Suppose that we have computed non-negative eigenvalues $\{\lambda_i\}_{i\in\mathbb I}$ and associated eigenfunctions $\{\psi_i\}_{i\in\mathbb I}$ satisfying
\begin{equation}
    \label{Eq:HI-20}
    \mathcal R\psi_i=\lambda_i\psi_i\text{ for }i\in\mathbb I,\quad\left\{
    \begin{aligned}
        &\lambda_1\ge\ldots\ge\lambda_d>\lambda_{d+1}=\ldots =0&&\text{if }d<\infty,\\
        &\lambda_1\ge\lambda_2\ge\ldots\text{ and }\lim_{i\to\infty}\lambda_i=0&&\text{if }d=\infty.
    \end{aligned}
    \right.
\end{equation}
If we want to indicate that the POD basis is computed by choosing $X=H$ or $X=V$ we write $\{\psi_i^H\}_{i\in\mathbb I}$ or $\{\psi_i^V\}_{i\in\mathbb I}$, respectively. Let $\ell\in\mathbb N$ be arbitrarily chosen with $\ell\le d$ provided $d$ is finite. The integral operator $\mathcal R$ has been introduced in \eqref{SIAM:Eq-I.1.2.5}. Then we define the POD subspace
\begin{align*}
    X^\ell=\mathrm{span}\,\big\{\psi_1,\ldots,\psi_\ell\big\}\subset X.
\end{align*}
From part 1) of Lemma~\ref{SIAM:Lemma3.2.1} we infer that $X^\ell\subset V$ holds true. If we want to consider either the specific case $X=H$ or $X=V$ we make use of the notation (cf. Section~\ref{Section:ContPODHilbert})
\begin{align*}
    H^\ell=\mathrm{span}\,\left\{\psi_1^H,\ldots,\psi_\ell^H\right\}\quad\text{and}\quad V^\ell=\mathrm{span}\,\left\{\psi_1^V,\ldots,\psi_\ell^V\right\},
\end{align*}

The POD Galerkin scheme for \eqref{SIAM:Eq3.1.6} reads as follows: find $y^\ell: [0,T] \to X^\ell$ satisfying
\begin{subequations}
    \label{SIAM:Eq3.1.6POD}
    \begin{align}
    \label{SIAM:Eq3.1.6PODa}
    \frac{\mathrm d}{\mathrm dt} {\langle y^\ell(t),\psi \rangle}_H+a(t;y^\ell(t),\psi)&={\langle (\mathcal F+\mathcal Bu)(t),\psi\rangle}_{V',V}\quad\text{for all }\psi\in X^\ell\text{ a.e. in }(0,T],\\
    \label{SIAM:Eq3.1.6PODb}
    y^\ell(0)&=\mathcal P^\ell y_\circ,
    \end{align}
\end{subequations}
where $\mathcal P^\ell:Z\to X^\ell$ stands for one of the following four orthogonal projections (cf. \eqref{SIAM:Eq3.2.12}):
\begin{align*}
    \mathcal P^\ell=\left\{
    \begin{aligned}
        &\mathcal P^\ell_H&&\text{for }X^\ell=H^\ell\text{ and } Z=H,\\
        &\mathcal P^\ell_V&&\text{for }X^\ell=H^\ell\text{ and } Z=V,\\
        &\mathcal Q^\ell_H&&\text{for }X^\ell=V^\ell\text{ and } Z=H,\\
        &\mathcal Q^\ell_V&&\text{for }X^\ell=V^\ell\text{ and } Z=V.
    \end{aligned}
    \right.
\end{align*}
We can see that in comparison to the original reduced formulation \eqref{SIAM:Eq3.1.6}, the space $V$ has been replaced by the finite-dimensional linear space $X^\ell\subset V$. In the upcoming analysis, we will therefore be required to enhance the original Gelfand triple, which is done in the following remark.

\begin{remark}
    \label{Remark:Vell_Identification}
    \rm Using the property of the \index{Gelfand triple}{\em Gelfand triple} we get
    \begin{align*}
        V^\ell \hookrightarrow V \hookrightarrow H \simeq H' \hookrightarrow V' \hookrightarrow (V^\ell)'.
    \end{align*}
    Here, for an arbitrary $\varphi \in V^\ell$ the corresponding $\varphi \in V'$ is seen as the mapping 
    \begin{align*}
        \psi\mapsto{\langle \varphi,\psi \rangle}_H={\langle \varphi,\psi \rangle}_{V',V} \quad \text{for all } \psi \in V.
    \end{align*}
    The embedding $V' \hookrightarrow (V^\ell)'$ is done by restricting a function $\varphi \in V'$ onto the domain $V^{\ell}$. By doing so, it clearly holds
    \begin{align*}
        {\| \varphi \|}_{(V^\ell)'} \leq {\| \varphi \|}_{V'}.
    \end{align*}
    For an arbitrary $\varphi \in V^\ell$, the reverse additionally holds true. To see that, le $\psi \in V$ with ${\| \psi \|}_V = 1$ be arbitrary. We can uniquely write $\psi$ as $\psi = \psi^\ell + \psi^\perp$ with $\psi^\ell \in V^\ell$ and $\psi^\perp \in (V^\ell)^{\perp_H}$, i.e. $\psi^\perp$ is in the orthogonal complement of $V^\ell$ with respect to the inner product of $H$. Then we get
    \begin{align*}
        {\langle \varphi , \psi \rangle_{V',V}} = {\langle \varphi , \psi^\ell + \psi^\perp \rangle_H} = {\langle \varphi , \psi^\ell \rangle_H} = {\langle \varphi , \psi^\ell \rangle_{V',V}} \leq {\| \varphi \|}_{(V^\ell)'}.
    \end{align*}
    Hence, 
    \begin{align}
        \label{eq:VellIsometry}
        {\| \varphi \|}_{(V^\ell)'} = {\| \varphi \|}_{V'} \qquad \text{for all } \varphi \in V^\ell
    \end{align}
    holds true.\hfill$\blacksquare$
\end{remark}

It follows from $y^\ell(t)\in X^\ell$ almost everywhere in $[0,T]$ that there exists a vector $\mathrm y^\ell(t)=(\mathrm y_i^\ell(t))\in\mathbb R^\ell$ ($t\in[0,T]$), which satisfies
\begin{equation}
    \label{PODGalAns}
    y^\ell(t)=\sum_{i=1}^\ell\mathrm y_i^\ell(t)\psi_i\quad\text{f.a.a. }t\in [0,T].
\end{equation}
Let us define the \index{Matrix!mass!POD, $\bM^\ell$}{\em mass} and (time-dependent) \index{Matrix!stiffness!POD, $\bA^\ell$}{\em stiffness matrices}
\begin{equation}
    \label{M_AMatrices}
    \bM^\ell=\big(\big({\langle \psi_j,\psi_i \rangle}_H\big)\big)\in\mathbb R^{\ell\times\ell}\quad\text{and}\quad\bA^\ell(t)=\big(\big(a(t;\psi_j,\psi_i)\big)\big)\in\mathbb R^{\ell\times\ell},
\end{equation}
respectively. Since the bilinear form $a(t;\cdot\,,\cdot)$ might not be symmetric, $\bA^\ell(t)$ can not be assumed to be symmetric. The next result is shown in Section~\ref{SIAM-Book:Section3.8.2}.

\begin{lemma}
    \label{MApdef}
    The mass matrix $\bM^\ell$ is symmetric and positive definite. For the stiffness matrix $\bA^\ell(t)$ we have
    \begin{equation}
        \label{EqvAv}
        \mathrm v^\top\bA^\ell(t)\mathrm v>0\quad\text{for all }\mathrm v\in\mathbb R^\ell \setminus \{0 \}
    \end{equation}
    provided \eqref{SIAM:Eq3.1.1-2} is satisfied with $\gamma_2=0$. In particular, $\bA^\ell(t)$ is then regular in $[0,T]$.
\end{lemma}

\begin{remark}
    \rm Since the mass matrix $\bM^\ell$ is positive definite, we can introduce the weighted inner product
    \begin{align*}
        {\langle \mathrm v^\ell,\mathrm w^\ell\rangle}_{\bM^\ell}=\big(\mathrm v^\ell\big)^\top\bM^\ell\mathrm w^\ell\quad\text{for }\mathrm v^\ell=(\mathrm v_i^\ell)^\top,\mathrm w^\ell=(\mathrm w_i^\ell)^\top\in\mathbb R^\ell.
    \end{align*}
    It follows that
    \begin{align*}
        {\langle v^\ell,w^\ell\rangle}_H={\langle \mathrm v^\ell,\mathrm w^\ell\rangle}_{\bM^\ell}\quad\text{for }v^\ell=\sum_{i=1}^\ell\mathrm v_i^\ell\psi_i,\,w^\ell=\sum_{i=1}^\ell\mathrm w_i^\ell\psi_i
    \end{align*}
    which implies also $\|v^\ell\|_H=\|\mathrm v\|_{\bM^\ell}$. Especially, the POD subspace $H^\ell$ is isometrically isomorph to the space $\mathbb R^\ell$ of the POD coefficients endowed with the weighted inner product $\langle v^\ell, w^\ell \rangle_{\bM^\ell} = (v^\ell)^\top\bM^\ell w^\ell$. Likewise, $V^\ell$ is isometrically isomorph to the space $\mathbb R^\ell$ endowed with the weighted inner product $\langle v^\ell, w^\ell \rangle_{\bS^\ell} = (v^\ell)^\top \bS^\ell w^\ell$ with the matrix $\bS^\ell=((S_{ij}^\ell))\in\mathbb R^{\ell\times\ell}$ given by $\bS_{ij}^\ell = \langle \psi_j, \psi_i \rangle_V$ for $i,j=1,...,\ell$; cf. Example \ref{SIAM:Example-I.1.1.1}-b).\hfill$\blacksquare$
\end{remark}

Inserting \eqref{PODGalAns} into \eqref{SIAM:Eq3.1.6PODa} and choosing $\psi=\psi_i$ for $1\le i\le \ell$ we derive the following system of ordinary differential equations for $\mathrm y^\ell$:
\begin{subequations}
    \label{SIAM:PODStateDis}
    \begin{equation}
        \label{SIAM:PODStateDis-a}
        \bM^\ell\dot{\mathrm y}^\ell(t)+\bA^\ell(t)\mathrm y^\ell(t)=\mathrm g^\ell(t;u)\quad\text{f.a.a. }t\in(0,T],
    \end{equation}
    where we have set
    \begin{align*}
        \mathrm g^\ell(t;u)=\left({\langle (\mathcal F+\mathcal Bu)(t),\psi_i\rangle}_{V',V}\right)\in\mathbb R^\ell\quad\text{for }u\in\U\text{ and }t\in[0,T].
    \end{align*}
    Since $\mathcal P^\ell y_\circ\in X^\ell$ holds, there is exactly one coefficient $\mathrm y_\circ^\ell=(\mathrm y_{\circ i}^\ell)\in \mathbb R^\ell$ such that
    \begin{align*}
        \mathcal P^\ell y_\circ=\sum_{i=1}^\ell\mathrm y_{\circ i}\psi_i.
    \end{align*}
    Then \eqref{SIAM:Eq3.1.6PODb} reads
    \begin{equation}
        \label{SIAM:PODStateDis-b}
        \mathrm y^\ell(0)=\mathrm y_\circ^\ell.
    \end{equation}
\end{subequations}
If $\mathrm y^\ell$ solves \eqref{SIAM:PODStateDis}, then $y^\ell$ given by \eqref{PODGalAns} is a solution to \eqref{SIAM:Eq3.1.6POD}. The other way round, the coefficient vector $\mathrm y^\ell$ solves \eqref{SIAM:PODStateDis} provided $y^\ell$ is a solution to \eqref{SIAM:Eq3.1.6POD}.

Note that \eqref{SIAM:PODStateDis} is an $\ell$-dimensional initial value problem. Since the mapping $t\mapsto \mathrm g^\ell(t;u)$ might not be continuous in $[0,T]$ for given $u\in\U$, the existence of a unique solution $\mathrm y^\ell$ to \eqref{SIAM:PODStateDis} does not follow from the theory of Picard-Lindel\"of. However, existence of a unique solution can be ensured by arguing similarly as in the proof of Theorem~\ref{SIAM:Theorem3.1.1}; cf. Section~\ref{SIAM-Book:Section3.8.2}.

\begin{theorem}
    \label{SIAM:Theorem3.1.1POD}
    Let Assumption~{\rm\ref{A2}} hold and $\{\psi_i\}_{i=1}^\ell$ solve \eqref{Eq:HI-20} for all $\ell\in\{1,\ldots,d\}$. Then for every $u\in\U$, there is a unique weak solution $y^\ell\in H^1(0,T;V)\hookrightarrow\Y=W(0,T)$ satisfying \eqref{SIAM:Eq3.1.6POD} and
    \begin{equation}
        \label{SIAM:Eq3.1.7POD}
        {\|y^\ell\|}_\mathscr Y\le C\left({\|\mathcal P^\ell y_\circ\|}_H+{\|\mathcal F\|}_{L^2(0,T;V')}+{\|u\|}_\U\right)
    \end{equation}
    for a constant $C>0$ which is independent of $\ell$, $u$, $y_\circ$ and $\mathcal F$.
\end{theorem}

\begin{remark}
    \label{Remark:HI-23}
    \rm If $\mathcal P^\ell=\mathcal P^\ell_H$ or $\mathcal P^\ell=\mathcal Q^\ell_H$ (i.e., $\mathcal P^\ell$ is an $H$-orthogonal projection), then we have
    \begin{align*}
        {\|\mathcal P^\ell y_\circ\|}_H\le{\|y_\circ\|}_H.
    \end{align*}
    If $\mathcal P^\ell=\mathcal P^\ell_V$ or $\mathcal P^\ell=\mathcal Q^\ell_V$ is chosen (i.e., $\mathcal P^\ell$ is an $H$-orthogonal projection), we conclude that $\mathcal P^\ell y_\circ\in V$ for any $y_\circ\in V$. Using the embedding inequality \eqref{Poincare} and $\|\mathcal P^\ell\|_{\mathscr L(V)}=1$ we derive
    \begin{align*}
        {\|\mathcal P^\ell y_\circ\|}_H\le c_V\,{\|\mathcal P^\ell y_\circ\|}_V\le c_V\,{\|y_\circ\|}_V
    \end{align*}
    for any initial condition satisfying $y_\circ\in V$. Consequently, if $\mathcal P^\ell$ is either an $H$- or a $V$-orthogonal projection, we infer from \eqref{SIAM:Eq3.1.7POD} that
    \begin{align*}
        {\|y^\ell\|}_\mathscr Y\le\left\{
        \begin{aligned}
            &C\left({\|y_\circ\|}_H+{\|\mathcal F\|}_{L^2(0,T;V')}+{\|u\|}_\U\right)&&\text{if }\mathcal P^\ell=\mathcal P^\ell_H\text{ or }\mathcal P^\ell=\mathcal Q^\ell_H,\\
            &\tilde C\left({\|y_\circ\|}_V+{\|\mathcal F\|}_{L^2(0,T;V')}+{\|u\|}_\U\right)&&\text{if }\mathcal P^\ell=\mathcal P^\ell_V\text{ or }\mathcal P^\ell=\mathcal Q^\ell_V\\
        \end{aligned}
        \right.
    \end{align*}
    for $\tilde C=\max\{1,c_V\}\,C$. In particular, $y^\ell$ is bounded in the $\mathscr Y$-norm by a constant that does not depend on $\ell$.\hfill$\blacksquare$
\end{remark}

Although the weak solution $y^\ell$ to \eqref{SIAM:Eq3.1.6POD} belongs to $H^1(0,T;V^\ell)$, we do not get an $H^1(0,T;V)$-estimate like \eqref{SIAM:Eq3.1.7POD} with an $\ell$-independent constant $C$. The reason is given in the next lemma which is proved in Section~\ref{SIAM-Book:Section3.8.2}.

\begin{lemma}
    \label{lem:W0T_finDim}
    Let Assumption~{\em\ref{A1}-1)} hold, but both $V$ and $H$ be of identical finite dimension $\mathfrak n<\infty$. Then the space $\mathscr Y$ is isomorphic to the space $H^1(0,T;V)$. In particular, there exists a constant $C_\mathfrak n>0$ such that the following {\em inverse inequality} holds:
    \begin{align}
        \label{eq:W0T_finDim}
        {\|\varphi\|}_{H^1(0,T;V)}\le C_\mathfrak n\,{\|\varphi\|}_\mathscr Y\quad\text{for all }\varphi\in H^1(0,T;V).
    \end{align}
\end{lemma}

\begin{remark}
    \label{RemH1WEq}
    \rm It follows from Lemma~\ref{lem:W0T_finDim}, \eqref{SIAM:Eq3.1.7POD} and Remark~\ref{Remark:HI-23} that there exists a constant $C_\ell>0$ depending on $\ell$ so that
    \begin{align*}
        {\|y^\ell\|}_{H^1(0,T;V)}\le C_\ell\left({\|y_\circ\|}_H+{\|\mathcal F\|}_{L^2(0,T;V')}+{\|u\|}_\U\right)
    \end{align*}
    holds true.\hfill$\blacksquare$
\end{remark}

The following corollary is a consequence of Theorem~\ref{SIAM:Theorem3.1.1POD}; cf. Corollary~\ref{Corollary:HI-20}.

\begin{corollary}
    \label{Corollary:HI-21}
    Let Assumption~{\rm\ref{A2}} hold.
    \begin{enumerate}
        \item [\rm 1)] There is a unique solution $\hat y^\ell\in H^1(0,T;V)\hookrightarrow\Y$ to
        \begin{align}
            \label{eq:reducedStateEq_yHat}
            \begin{aligned}
                \frac{\mathrm d}{\mathrm dt} {\langle\hat y^\ell(t),\psi \rangle}_H+a(t;\hat y^\ell(t),\psi)&={\langle\mathcal F(t),\psi\rangle}_{V',V}\quad\text{for all }\psi\in V^\ell\text{ a.e. in }(0,T],\\
                \hat y^\ell(0)&=\mathcal P^\ell y_\circ.
            \end{aligned}
        \end{align}
        Moreover,
        \begin{align*}
            {\|\hat y^\ell\|}_{\Y}\le C\left({\|\mathcal P^\ell y_\circ\|}_H+{\|\mathcal F\|}_{L^2(0,T;V')}\right)
        \end{align*}
        for a constant $C>0$ which is independent of $\ell$, $y_\circ$ and $\mathcal F$.
        \item [\rm 2)] Let us define the linear \index{Linear operator!POD solution, $\mathcal S^\ell$}{\em POD solution operator} $\mathcal S^\ell:\U\to \Y$ as follows: for every $u\in\U$ the function $\tilde y^\ell=\mathcal S^\ell u$ is the solution to
        \begin{align}
            \label{eq:reducedStateEq_Sell}
            \begin{aligned}
                \frac{\mathrm d}{\mathrm dt} {\langle\tilde y^\ell(t),\psi \rangle}_H+a(t;\tilde y^\ell(t),\psi)&={\langle (\mathcal Bu)(t),\psi\rangle}_{V',V}\quad\text{for all }\psi\in V^\ell\text{ a.e. in }(0,T],\\
                \tilde y^\ell(0)&=0.
            \end{aligned}
        \end{align}
        Then $\mathcal S^\ell$ is well-defined and continuous, i.e.,
        \begin{align*}
            {\|\mathcal S^\ell u\|}_\Y\le C\,{\|u\|}_\U
        \end{align*}
        holds for all $u \in \U$ for a constant $C>0$ which is independent of $\ell$ and $u$. Additionally, $\mathcal S^\ell u \in H^1(0,T;V)$ for all $u \in \U$.
    \end{enumerate}
\end{corollary}

\begin{remark}
    \label{Remark:HI-23-2}
    \rm Similar to Remark~\ref{Remark:HI-22} we can express the solution to \eqref{SIAM:Eq3.1.6POD} as $y^\ell=\hat y^\ell+\mathcal S^\ell u$. Again, the linear operator $\mathcal S^\ell$ is called a \index{Linear operator!control-to-state}{\em linear control-to-state operator}.\hfill$\blacksquare$
\end{remark}

\subsection{POD a-priori error analysis}
\label{SIAM-Book:Section3.3.2}

Let Assumption~\ref{A2} hold. Our next goal is to derive a-priori error estimates for the term
\begin{align*}
    {\|y-y^\ell\|}_\mathscr Y^2 = \int_0^T{\|y(t)-y^\ell(t)\|}^2_V+{\|y_t(t)-y_t^\ell(t)\|}^2_{V'}\,\mathrm dt,
\end{align*}
where $y$ and $y^\ell$ are the solutions to \eqref{SIAM:Eq3.1.6} and \eqref{SIAM:Eq3.1.6POD}, respectively, and $\mathscr Y=W(0,T)$ holds. We will later specify the snapshots which are utilized to compute the POD basis $\{\psi_i\}_{i=1}^\ell$ of rank $\ell$. Recall that we have introduced the projection $\mathcal P^\ell:Z\to X^\ell$ in \eqref{SIAM:Eq3.1.6PODb}. We shall make use of the following decomposition (cf., e.g., \cite[p.~8]{Tho97} and \cite{KV01,KV02a,KV02b})
\begin{equation}
    \label{EqDec}
    y(t)-y^\ell(t)=y(t)-\mathcal P^\ell y(t)+\mathcal P^\ell y(t)-y^\ell(t)=\varrho^\ell(t)+\vartheta^\ell(t)\quad\text{in }[0,T]\text{ a.e.}
\end{equation}
with $\varrho^\ell(t)=y(t)-\mathcal P^\ell y(t)\in (X^\ell)^\bot$ and $\vartheta^\ell(t)=\mathcal P^\ell y(t)-y^\ell(t)\in X^\ell$. It follows that
\begin{equation}
    \label{APriori-Est-1AA}
    {\|y-y^\ell\|}^2_\mathscr Y\le2 \big( {\|\varrho^\ell\|}^2_\mathscr Y+{\|\vartheta^\ell\|}^2_\mathscr Y \big).
\end{equation}

First we estimate the $\vartheta^\ell$-term. For a proof we refer the reader to Section~\ref{SIAM-Book:Section3.8.2}.

\begin{lemma}
    \label{Lemma:HI-22}
    Let Assumption~{\rm \ref{A2}} be fulfilled. Additionally, suppose that $y \in \mathscr Y$ holds. Then
    \begin{equation}
        \label{APriori-Est-5}
        {\|\vartheta^\ell\|}_\mathscr Y^2\le C\,{\|\varrho^\ell\|}_\mathscr Y^2
    \end{equation}
    for a constant $C>0$ which are independent of $\ell$, $y_\circ$, $f$ and $u$.
\end{lemma}

\begin{remark}
    \label{RemarkThetaEstimates}
    \rm
    \begin{enumerate}
        \item[1)] Combining \eqref{APriori-Est-1AA} and \eqref{APriori-Est-5} we infer that
        \begin{equation}
            \label{APriori-Est-6}
            {\|y-y^\ell\|}_\mathscr Y^2\le C_1 \,{\|\varrho^\ell\|}_\mathscr Y^2.
        \end{equation}
        for $C_1=2(1+C)>0$.
        \item [2)] Recall that $L^2(0,T;H)\hookrightarrow L^2(0,T;V')$ holds; cf. \cite[p.~471]{DL00}. Thus, there is an embedding constant $\hat c_\mathsf e>0$ satisfying
        \begin{equation}
            \label{Embi1}
            {\|\varphi\|}_{L^2(0,T;V')}\le \hat c_\mathsf e\,{\|\varphi\|}_{L^2(0,T;H)}\quad\text{for all }\varphi\in L^2(0,T;H).
        \end{equation}
        If the solution $y$ to \eqref{SIAM:Eq3.1.6} belongs to $H^1(0,T;H)\cap L^2(0,T;V)\hookrightarrow W(0,T)$, we derive from \eqref{APriori-Est-6} and \eqref{Embi1} that
        \begin{equation}
            \label{Estimate100}
            \begin{aligned}
                {\|y-y^\ell\|}_\mathscr Y^2&\le C_1 \,{\|\varrho^\ell\|}_\mathscr Y^2=C_1 \big({\|\varrho^\ell_t\|}_{L(0,T;V')}^2+{\|\varrho^\ell\|}_{L^2(0,T;V)}^2\big)\\
                &\le C_1 \big(\hat c_\mathsf e^2\,{\|\varrho^\ell_t\|}_{L(0,T;H)}^2+{\|\varrho^\ell\|}_{L(0,T;V)}^2\big)\\
                &=C_2 \big({\|\varrho^\ell_t\|}_{L(0,T;H)}^2+{\|\varrho^\ell\|}_{L(0,T;V)}^2\big)
            \end{aligned}
        \end{equation}
        holds true for $C_2=C_1\max\{1,\hat c_\mathsf e^2\}>0$.
        \item [3)] From the two embeddings $H^1(0,T;V)\hookrightarrow H^1(0,T;V')$ and $H^1(0,T;V)\hookrightarrow L^2(0,T;V)$ (cf. \cite[p.~471]{DL00}) and from $W(0,T)=H^1(0,T;V')\cap L^2(0,T;V)$ we infer that $H^1(0,T;V)\hookrightarrow W(0,T)$. Then, there is a constant $\tilde c_\mathsf e>0$ such that
        \begin{equation}
            \label{Embi2}
            {\|\varphi\|}_\mathscr Y\le \tilde c_\mathsf e\,{\|\varphi\|}_{H^1(0,T;V)}\quad\text{for all }\varphi\in H^1(0,T;V).
        \end{equation}
        Suppose that the solution $y$ to \eqref{SIAM:Eq3.1.6} belongs to $H^1(0,T;V)$. Then, it follows from \eqref{APriori-Est-6} and \eqref{Embi2} that
        \begin{equation}
            \label{Estimate101}
            {\|y-y^\ell\|}_\mathscr Y^2\le C_1 \,{\|\varrho^\ell\|}_\mathscr Y^2\le C_3\,{\|\varrho^\ell\|}_{H^1(0,T;V)}^2
        \end{equation}
        is satisfied for $C_3=\tilde c_\mathsf e^2C_1>0$.\hfill$\blacksquare$
    \end{enumerate}
\end{remark}

If $\mathcal P^\ell$ is an $H$-orthonormal projection (i.e., $\mathcal P^\ell=\mathcal P_H^\ell$ or $\mathcal P^\ell=\mathcal Q_H^\ell$), it turns out that we can omit the time derivatives of $\varrho^\ell$ in our estimate; see Section \ref{SIAM-Book:Section3.8.2} for a proof.
\begin{lemma}
    \label{APriori-Est-inV}
    Let Assumption~{\rm \ref{A2}} be fulfilled. Additionally, suppose that $y \in \mathscr Y\cap H^1(0,T;H)$ holds and that $\mathcal P^\ell$ denotes either $\mathcal P_H^\ell$ or $\mathcal Q_H^\ell$. Set $\vartheta^\ell(t)=\mathcal P^\ell y(t)-y^\ell(t)\in X^\ell$ for all $t \in [0,T]$. Then
    \begin{equation}
    \label{APriori-Est-inV-1}
        {\|\vartheta^\ell\|}_\mathscr Y^2\le \tilde C\,{\|\varrho^\ell\|}_{L^2(0,T;V)}^2
    \end{equation}
    for a constant $\tilde C>0$ which are independent of $\ell$, $y_\circ$, $\mathcal F$ and $u$. In particular, no time derivatives of $\varrho^\ell$ occur on the right-hand side of \eqref{APriori-Est-inV-1}.
\end{lemma}

Next, we discuss the choices $X=H$ and $V$ as well as different choices for the orthogonal projection $\mathcal P^\ell$. For that purpose we proceed as in Section~\ref{Section:ContPODHilbert} and introduce different notations for the POD bases in $H$ and $V$: Let the eigenvalue-eigenfunction pairs $\{(\lambda_i^H,\psi_i^H)\}_{i\in\mathbb I}$ and $\{(\lambda_i^V,\psi_i^V)\}_{i\in\mathbb I}$ satisfy
\begin{subequations}
    \label{EigProblemHV}
    \begin{align}
        \label{EigProblemHV-H}
        \mathcal R^H\psi_i^H&=\lambda_i^H\psi_i^H,&&\left\{
        \begin{aligned}
            &\lambda_1^H\ge\ldots\ge\lambda_d^H>\lambda_{d+1}^H=\ldots=0&&\text{if }d=d^H<\infty,\\
            &\lambda_1^H\ge\lambda_2^H\ge\ldots\text{ and }\lim_{i\to\infty}\lambda_i^H=0&&\text{if }d^H=\infty,
        \end{aligned}
        \right.\\
        \label{EigProblemHV-V}
        \mathcal R^V\psi_i^V&=\lambda_i^V\psi_i^V,&&\left\{
        \begin{aligned}
            &\lambda_1^H\ge\ldots\ge\lambda_d^V>\lambda_{d+1}^V=\ldots=0&&\text{if }d=d^V<\infty,\\
            &\lambda_1^H\ge\lambda_2^V\ge\ldots\text{ and }\lim_{i\to\infty}\lambda_i^V=0&&\text{if }d^V=\infty,
        \end{aligned}
        \right.
    \end{align}
\end{subequations}
where the integral operators $\mathcal R^H$ and $\mathcal R^V$ have been introduced in \eqref{OperatorsR}. Next we analyze different choices for the snapshots in the POD method. 

\subsubsection{Snapshot space with time derivatives}

Let us particularly assume that $y\in H^1(0,T;X)$ holds and that we have set $K=2$, $\omega_1^K=\omega_2^K=1$, $y^1=y$ and $y^2=y_t$, so both the solution trajectory and its derivative are considered to build the POD basis. Then we distinguish the following choices for the orthogonal projection $\mathcal P^\ell$ (cf. \eqref{SIAM:Eq3.2.12})
\begin{enumerate}
    \item [1)] {\em Case $X=H$ and $\mathcal P^\ell=\mathcal P^\ell_H$:} In this case we have $y\in H^1(0,T;H)$ so that we can apply Remark~\ref{RemarkThetaEstimates}-2). We infer for all $\ell\ge1$ from \eqref{Estimate100}, \eqref{RatePH-V} and \eqref{ErrRateH}
    \begin{align*}
        {\|y-y^\ell\|}_\mathscr Y^2&\le C_2\int_0^T \big\|y^1(t)-\mathcal P^\ell_Hy^1(t)\big\|_V^2+\big\|y^2(t)-\mathcal P^\ell_H y^2(t)\big\|_H^2\,\mathrm dt\\
        &\le C_2\sum_{k=1}^K\omega_k^K\int_0^T \big\|y^k(t)-\mathcal P^\ell_Hy^k(t)\big\|_V^2+\big\|y^k(t)-\mathcal P^\ell_H y^k(t)\big\|_H^2\,\mathrm dt\\
        &= C_2\sum_{i>\ell}\lambda_i^H\big({\|\psi_i^H\|}_V^2+1\big).
    \end{align*}
    \item [2)] {\em Case $X=H$ and $\mathcal P^\ell=\mathcal P^\ell_V$:} Again we have $y\in H^1(0,T;H)$ and can utilze Remark~\ref{RemarkThetaEstimates}-2). From \eqref{Estimate100}, \eqref{Poincare} and \eqref{RatePV-V} we find for all $\ell\ge1$
    \begin{align*}
        {\|y-y^\ell\|}_\mathscr Y^2&\le C_2\int_0^T\big\|y^1(t)-\mathcal P^\ell_V y^1(t)\big\|_V^2+\big\|y^2(t)-\mathcal P^\ell_V y^2(t)\big\|_H^2\,\mathrm dt\\
        &\le C_2\int_0^T\big\|y^1(t)-\mathcal P^\ell_V y^1(t)\big\|_V^2+c_V^2\big\|y^2(t)-\mathcal P^\ell_V y^2(t)\big\|_V^2\,\mathrm dt\\
        &\le C_4\sum_{k=1}^K\omega_k^K\int_0^T\big\|y^k(t)-\mathcal P^\ell_V y^k(t)\big\|_V^2\,\mathrm dt= C_4\sum_{i>\ell}\lambda_i^H\big\|\psi_i^H-\mathcal P^\ell_V\psi_i^H\big\|_V^2.
    \end{align*}
    with $C_4=\max\{1,c_V^2\}C_2$.
    \item [3)] {\em Case $X=V$ and $\mathcal P^\ell=\mathcal Q^\ell_H$:} Now we have $y\in H^1(0,T;V)$ so that Remark~\ref{RemarkThetaEstimates}-3) can be used. We combine \eqref{Estimate101} and \eqref{RateQH-V}. It follows for all $\ell\ge1$ that
    \begin{align*}
        {\|y-y^\ell\|}_\mathscr Y^2&\le C_3\sum_{k=1}^K\omega_k^K\int_0^T\big\|y^k(t)-\mathcal Q^\ell_H y^k(t)\big\|_V^2\,\mathrm dt\\
        &=C_3\sum_{i>\ell}\lambda_i^V{\|\psi_i^V-\mathcal Q_H^\ell\psi_i^V\|}_V^2.
    \end{align*}
    \item [4)] {\em Case $X=V$ and $\mathcal P^\ell=\mathcal Q^\ell_V$:} Due to \eqref{Estimate101} and \eqref{RateQV-V} and we derive that
    \begin{align*}
        {\|y-y^\ell\|}_\mathscr Y^2\le C_3\sum_{k=1}^K\omega_k^K\int_0^T\big\|y^k(t)-\mathcal Q^\ell_V y^k(t)\big\|_V^2\,\mathrm dt=C_3\sum_{i>\ell}\lambda_i^V.
    \end{align*}
\end{enumerate}

Notice that we have in all four cases an estimate of the form
\begin{align*}
    {\|y-y^\ell\|}_\mathscr Y^2\le C_5\sum_{k=1}^K\omega_k^K{\|y^k-\mathcal P^\ell y^k\|}_{L^2(0,T;V)}^2
\end{align*}
for the constant $C_5=\max\{2C_2,C_3,C_4\}>0$. Using \cite[Theorems~5.2 and 5.3]{Sin14} we get convergence of $y-y^\ell$ in $\mathscr Y$, where we have to suppose in the case $X=H$ and $\mathcal P^\ell=\mathcal Q^\ell_H$ that $\|\mathcal Q^\ell_H\|_{\mathscr L(V)}$ is bounded independently of $\ell$. Summarizing we have derived the following theorem.

\begin{theorem}
    \label{Th:A-PrioriError}
    Let Assumption~{\rm\ref{A2}} hold. Assume that for $u\in\U$ and $y_\circ\in X$ the unique solution $y$ to \eqref{SIAM:Eq3.1.6} satisfies $y\in H^1(0,T;X)$. In \eqref{SIAM:PellH} and \eqref{SIAM:PellV} we choose $K=2$, $\omega_1^K=\omega_2^K=1$, $y^1=y$ and $y^2=y_t$. Let the eigenvalue-eigenfunction pairs $\{(\lambda_i^H,\psi_i^H)\}_{i\in\mathbb I}$ and $\{(\lambda_i^V,\psi_i^V)\}_{i\in\mathbb I}$ satisfy \eqref{EigProblemHV}. Furthermore, $y^\ell$ denotes the solution to \eqref{SIAM:Eq3.1.6POD}. Then we have for all $\ell\ge1$ the following \index{Error estimate!a-priori!state variable}a-priori error estimate
    \begin{align*}
        {\|y-y^\ell\|}_\mathscr Y^2\le C\cdot\left\{
        \begin{aligned}
            &\sum_{i>\ell}\lambda_i^H\big({\|\psi_i^H\|}_V^2+1\big)&&\text{for }X=H,~\mathcal P^\ell=\mathcal P^\ell_H,\\
            &\sum_{i>\ell}\lambda_i^H\big\|\psi_i^H-\mathcal P^\ell_V\psi_i^H\big\|_V^2&&\text{for }X=H,~\mathcal P^\ell=\mathcal P^\ell_V,\\
            &\sum_{i>\ell}\lambda_i^V{\|\psi_i^V-\mathcal Q_H^\ell\psi_i^V\|}_V^2&&\text{for }X=V,~\mathcal P^\ell=\mathcal Q^\ell_H,\\
            &\sum_{i>\ell}\lambda_i^V&&\text{for }X=V,~\mathcal P^\ell=\mathcal Q^\ell_V
        \end{aligned}
        \right.
    \end{align*}
    for a constant $C>0$ that is independent of $\ell$, $y_\circ$, $f$ and $u$. Moreover, we have $\lim_{\ell\to\infty}\|y-y^\ell\|_{L^2(0,T;V)}=0$, where we have to suppose in the case $X=H$ and $\mathcal P^\ell=\mathcal Q^\ell_H$ that $\|\mathcal Q^\ell_H\|_{\mathscr L(V)}$ is bounded independently of $\ell$.
\end{theorem}

\begin{remark}
    \label{Remark:EstMaxNorm-SnapOhneDQ}
    \rm
    \begin{enumerate}
        \item [1)] For the choices $X=V$ and $\mathcal P^\ell=\mathcal Q^\ell_V$ the error decays as the POD aproximation error; cf. \eqref{ErrRate-V}. Moreover, the evalution of
        \begin{align*}
            \mathcal P^\ell_V y_\circ=\sum_{i=1}^\ell{\langle y_\circ,\psi_i^V\rangle}_V\,\psi_i^V
        \end{align*}
        is computationally cheap. The only issue is the computation of the inner products $\langle y_\circ,\psi_i^V\rangle_V$ for $i=1,\ldots,\ell$. However, we require that the snapshots $y$ and $y_t$ are in $L^2(0,T,V)$, so that the solution to \eqref{SIAM:Eq3.1.6} has to belong to $H^1(0,T;V)$ which can not be expected at all; cf. Theorem~\ref{SIAM:Theorem3.1.1}. Less regularity of the solution -- in fact, $y \in H^1(0,T;H)$ -- is needed for the choices $X=H$ and $\mathcal P^\ell=\mathcal P^\ell_H$. Here, the effort to compute
        \begin{align*}
            \mathcal P^\ell_H y_\circ=\sum_{i=1}^\ell{\langle y_\circ,\psi_i^H\rangle}_H\,\psi_i^H
        \end{align*}
        is just like the computation of the projection $\mathcal P^\ell_V y_\circ$. However, the decay rate for the a-priori error is damped by the $V$-norm of the basis elements $\psi_i^H$. A better decay is acchievd for the choices $X=H$ and $\mathcal P^\ell=\mathcal P^\ell_{VH^\ell}$, but now the evaluation of $\mathcal P^\ell_{VH^\ell} y_\circ$ requires an $\ell$-dimensional system solve and, in addition, the computation of the inner products $\langle y_\circ,\psi_i^V\rangle_V$ for $i=1,\ldots,\ell$.
        \item [2)] Clearly, we can also bound the difference $y-y^\ell$ in the $L^2(0,T;H)$ norm by utilizing \eqref{Poincare} in 
        \begin{align*}
            {\|y-y^\ell\|}_{L^2(0,T;H)}^2\le c_V^2\,{\|y-y^\ell\|}_{L^2(0,T;H)}^2\le c_V^2\,{\|y-y^\ell\|}_\mathscr Y^2
        \end{align*}
        and applying Theorem~\ref{Th:A-PrioriError}.\hfill$\blacksquare$
    \end{enumerate}
\end{remark}

\subsubsection{Snapshot space with arbitrary data trajectory}

Next we discuss the situation when the POD basis is computed from snapshots which are different from the solution $y$ to \eqref{SIAM:Eq3.1.6}. We choose $\mathscr D=[0,T]$ and arbitrary trajectories $\{y^k\}_{k=1}^K\subset L^2(0,T;X)$ in \eqref{Eq2.2.6}. Furthermore, the solution $y$ to \eqref{SIAM:Eq3.1.6} belongs to $H^1(0,T;X)\hookrightarrow\mathscr Y$ holds so that we can apply Remark~\ref{RemarkThetaEstimates}-2). Due to \eqref{Estimate100}, we still have
\begin{align*}
    {\|y-y^\ell\|}_\mathscr Y^2\le C_2\int_0^T \big\|y(t)-\mathcal P^\ell y(t)\big\|_V^2+\,\big\|y_t(t)-\mathcal P^\ell y_t(t)\big\|_H^2\,\mathrm dt
\end{align*}
Since the solution $y$ is not utilized to get the POD basis of rank $\ell$, we can not apply Theorem~\ref{Prop:VTopology} in order to estimate the right-hand side (as we have done to prove Theorem~\ref{Prop:VTopology}). However, Proposition~\ref{Proposition:ProjConvV} can be used. Let $X=H$ and $\mathcal P^\ell=\mathcal P^\ell_H$. If $\lambda_i^H>0$ for all $i\in\mathbb N$ we get
\begin{align*}
    \lim_{\ell\to\infty}\left({\|y(t)-\mathcal P^\ell_H y(t)\|}_V^2+{\|y_t(t)-\mathcal P^\ell_H y_t(t)\|}_H^2\right)=0\quad\text{in }[0,T]\text{ a.e.}
\end{align*}
provided $y(t)\in\mathrm{ran}\,(\mathcal Y_H)$, where
\begin{align*}
    \mathcal Y_H \phi=\sum_{k=1}^K\omega_k^K\int_0^T\phi^k(t)y^k(t)\,\mathrm d0\quad\text{for }\phi=(\phi^1,\ldots,\phi^K)\in L^2(0,T;\mathbb R^k).
\end{align*}
Next we consider again $X=H$, but $\mathcal P^\ell=\mathcal P^\ell_V$. Again, we have to assume $\lambda_i^H>0$ for all $i\in\mathbb N$ in order to apply Proposition~\ref{Proposition:ProjConvV}. Using \eqref{Poincare} we find
\begin{align*}
    &\lim_{\ell\to\infty}\left({\|y(t)-\mathcal P^\ell_V y(t)\|}_V^2+{\|y_t(t)-\mathcal P^\ell_V y_t(t)\|}_H^2\right)\\
    &\quad\le\lim_{\ell\to\infty}\left({\|\mathcal P^\ell_V y(t)-y(t)\|}_V^2+c_V^2\,{\|\mathcal P^\ell_V y_t(t)-y_t(t)\|}_V^2\right)=0\quad\text{in }[0,T]\text{ a.e.}
\end{align*}
for $y(t),\,y_t(t)\in V$ provided $\lambda_i^H>0$ for all $i\in\mathbb N$. Finally, applying \eqref{Poincare} we find
\begin{align*}
    &\lim_{\ell\to\infty}\left({\|\mathcal Q_V^\ell y(t)-y(t)\|}_V^2+{\|\mathcal Q_V^\ell y_t(t)-y_t(t)\|}_H^2\right)\\
    &\quad\le\lim_{\ell\to\infty}\left({\|\mathcal Q_V^\ell y(t)-y(t)\|}_V^2+c_V^2\,{\|\mathcal Q_V^\ell y_t(t)-y_t(t)\|}_V^2\right)=0\quad\text{in }[0,T]\text{ a.e.}
\end{align*}
for $y(t),\,y_t(t)\in V$. Consequently, we have the following result:

\begin{theorem}
    \label{Corollary:APrioriConv}
    Suppose that Assumption~{\rm\ref{A2}} holds. For $u\in\U$ and $y_\circ\in X$ let the unique solution $y$ to \eqref{SIAM:Eq3.1.6} belong to $H^1(0,T;X)$. The POD basis of rank $\ell$ be computed by solving \eqref{SIAM:PellH} or \eqref{SIAM:PellV} for an arbitrarily chosen snapshot ensemble $\{y^k\}_{k=1}^K\subset L^2(0,T;V)$, where we assume, in addition, that
    \begin{enumerate}
        \item [\rm 1)] $\lambda_i^H>0$ for all $i\in\mathbb N$, $y\in L^2(0,T;\mathrm{ran}\,(\mathcal Y_H))$ for $X=H$ and $\mathcal P^\ell=\mathcal P^\ell_H$, 
        \item [\rm 2)] $\lambda_i^H>0$ for all $i\in\mathbb N$, $y\in H^1(0,T;V)$ for $X=H$ and $\mathcal P^\ell=\mathcal P^\ell_V$. 
    \end{enumerate}
    Let $y^\ell$ denote the solution to \eqref{SIAM:Eq3.1.6POD}. Then we have $\lim_{\ell\to\infty}{\|y-y^\ell\|}_\mathscr Y=0$.
\end{theorem}

\subsubsection{Snapshot space without time derivatives}

We infer from Theorem~\ref{Th:A-PrioriError} the rate of convergence for the difference $y-y^\ell$ in the $L^2(0,T;V)$ norm provided $y$ and $y_t$ are utilized as snapshots to compute the POD bases $\{\psi_i^H\}_{i\in\mathbb I}$ and $\{\psi_i^V\}_{i\in\mathbb I}$. For a specific setting we can avoid to include the time derivatives into the snapshot ensemble.

Notice that
\begin{align*}
    {\|\varphi\|}^2_{L^2(0,T;V)}\le{\|\varphi\|}^2_\mathscr Y\quad\text{for all }\varphi\in\mathscr Y
\end{align*}
Moreover, recall that $\varrho^\ell=y-\mathcal P^\ell y$ and $\vartheta^\ell=\mathcal P^\ell y-y^\ell$. Due to Lemma~\ref{APriori-Est-inV} we have 
\begin{equation}
    \label{APriori-Est-inV-2}
    {\|y-y^\ell\|}^2_{L^2(0,T;V)}\le2 \big( {\|\varrho^\ell\|}^2_{L^2(0,T;V)}+{\|\vartheta^\ell\|}^2_{L^2(0,T;V)} \big)\le\tilde C_1\,{\|\varrho^\ell\|}^2_{L^2(0,T;V)}
\end{equation}
with the constant $\tilde C_1=2(1+\tilde C)$. Now the next theorem follows by using Theorem~\ref{Prop:VTopology}.

\begin{theorem}
    \label{Th:A-PrioriError-2}
    Let Assumption~{\rm\ref{A2}} hold. Assume that for $u\in\U$ the unique solution $y$ to \eqref{SIAM:Eq3.1.6} satisfies $y\in H^1(0,T;H)$ and that $\mathcal P^\ell$ is given either by $\mathcal P_H^\ell$ or $\mathcal Q_H^\ell$. In \eqref{SIAM:PellH} and \eqref{SIAM:PellV} we choose $K=1$, $\omega_1^K=1$ and $y^1=y$. Let the eigenvalue-eigenfunction pairs $\{(\lambda_i^H,\psi_i^H)\}_{i\in\mathbb I}$ and $\{(\lambda_i^V,\psi_i^V)\}_{i\in\mathbb I}$ satisfy \eqref{EigProblemHV}. Furthermore, $y^\ell$ denotes the solution to \eqref{SIAM:Eq3.1.6POD}. Then there exists a constant $C>0$ satisfying the following \index{Error estimate!a-priori!state variable}a-priori error estimate
    \begin{align*}
        {\|y-y^\ell\|}_{L^2(0,T;V)}^2\le C\cdot\left\{
        \begin{aligned}
            &\sum_{i>\ell}\lambda_i^H\,\big\|\psi_i^H\big\|_V^2,&&X=H,~\mathcal P^\ell=\mathcal P^\ell_H,\\
            &\sum_{i>\ell}\lambda_i^V\,\big\|\psi_i^V-\mathcal Q^\ell_H\psi_i^V\big\|_V^2,&&X=V,~\mathcal P^\ell=\mathcal Q^\ell_H.
        \end{aligned}
        \right.
    \end{align*}
\end{theorem}

\subsubsection{Operator convergence}

Recall that we have introduced the linear, bounded operators $\mathcal S$ and $\mathcal S^\ell$ in Corollaries~\ref{Corollary:HI-20}-2 and \ref{Corollary:HI-21}-2, respectively. Moreover, the linear projection operator $\mathcal P^\ell$ has been utilized in \eqref{SIAM:Eq3.1.6PODb}. We give sufficient conditions that $\mathcal S^\ell$ converges to $\mathcal S$ in $\mathscr L(\U,\Y)$ as $\ell\to\infty$. The proof of the next result is given in Section~\ref{SIAM-Book:Section3.8.2}.

\begin{theorem}
    \label{Theorem:OperatorConv}
    Let Assumption~{\rm\ref{A2}} hold. We suppose that $\U$ is a finite-dimensional Hilbert space. Moreover, let one of the following hypotheses be satisfied:
    \begin{enumerate}
        \item [\rm 1)] $X=H$, $\mathcal P^\ell=\mathcal P^\ell_H$, $\lambda_i^H>0$ for all $i\in\mathbb I$ and $y(t)$, $y_t(t)\in\mathrm{ran}\,(\mathcal Y)\subset V$ for almost all $t\in[0,T]$;
        \item [\rm 2)] $X=H$, $\mathcal P^\ell=\mathcal P^\ell_{VH^\ell}$, $\lambda_i^H>0$ holds for all $i\in\mathbb I$;
        \item [\rm 3)] Case $X=V$, $\mathcal P^\ell=\mathcal P^\ell_V$: $\lambda_i^V>0$ for all $i\in\mathbb I$.
    \end{enumerate}
    Then $\|\mathcal S^\ell-\mathcal S\|_{\mathscr L(\U,\Y)}\to0$ as $\ell\to\infty$.
\end{theorem}

\begin{remark}
    \label{Remark:SellConv}
    \rm To prove the statement of Theorem~\ref{Theorem:OperatorConv} in the case that $\U$ is an infinite-dimensional Hilbert space we additionally need that the operator $\mathcal S$ is compact and that the operator family $\{\mathcal S^\ell\}_{\ell \in \mathbb{N}}$ is uniformly compact in the sense that
    \begin{align*}
        {\| \mathcal S^\ell u^\ell - \mathcal S^\ell u_\circ \|}_{\Y} \to 0 \quad \text{as } \ell \to \infty
    \end{align*}
    holds for all sequences $\{u^\ell\}_{\ell \in \mathbb{N}} \subset \U$ with $u^\ell \rightharpoonup u_\circ$ as $\ell \to \infty$.\hfill$\blacksquare$
\end{remark}

\section{The semidiscrete approximation}
\label{SIAM-Book:Section3.4}
\setcounter{equation}{0}
\setcounter{theorem}{0}
\setcounter{figure}{0}
\setcounter{run}{0}

To compute the POD basis $\{\psi_i\}_{i=1}^\ell$ as described in Section~\ref{SIAM-Book:Section3.3}, we need the snapshots $y(t)$ for almost all $t \in [0,T]$. This is realized numerically by computing approximations for $y(t)$ using a spatial and temporal discretization method. First, we consider a Galerkin discretization method in which the snapshots are given by a finite linear combination of basis functions but the time variable is still kept continuous. In a second step we turn to the temporal discretization. We refer the reader to \cite{Tho97}, where Galerkin discretizations for evolution problems are considered.

\subsection{Galerkin discretization}
\label{Spatial discretization}
\label{Section:3.4.1}

For a given mesh size parameter $h>0$ and associated number of nodes $m=m(h)\in\mathbb N$ let $\{\varphi_i^h\}_{i=1}^m\subset V$ be linearly independent functions. Then we introduce the $m$-dimensional subspace
\begin{equation}
    \label{FESpace}
    V^h=\mathrm{span} \, \big\{\varphi_1,\ldots,\varphi_m\big\} \subset V
\end{equation}
endowed with the topology in $V$. We make use of the following hypotheses.

\begin{assumption}
    \label{A3}
    There exists a subspace $\mathscr W\hookrightarrow V$ which is dense in $V$, such that
    \begin{equation}
        \label{Eq:ApproxProp-1}
        \inf_{v^h\in V^h}\big\{{\|v-v^h\|}_H+h\,{\|v-v^h\|}_V\big\}\le C_{\mathscr W}\,h^2\,{\|v\|}_{\mathscr W}\quad\text{for any }v\in\mathscr W,
    \end{equation}
    where the constant $C_{\mathscr W}$ is independent of $h$ and $v$.
\end{assumption}

\begin{example}
    \label{FE-Tri-Plot}
    \rm Suppose that $\Omega\subset\mathbb R^2$ is a bounded, convex polygonal domain with boundary $\Gamma=\partial\Omega$. Let $\{\mathcal T^h\}_{h>0}$ denote a family of partitions of $\Omega$ into disjunct triangles such that no vertex of any triangle lies on the interior of a side of another triangle and such that the union of the closures of all triangles is equal to $\overline\Omega$; cf. Figure~\ref{Fig:FE-Tri-Plot}.
    \begin{figure}
	    \centering
	    \includegraphics[height=50mm]{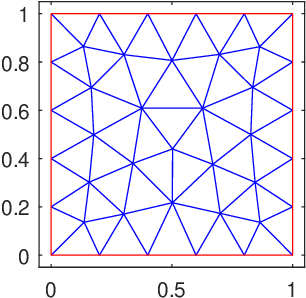}
	    \caption{Example~\ref{FE-Tri-Plot}. Finite element triangulation of the unit square $\Omega = (0,1)^2$ using $m = 41$ vertices and a maximal triangle side length of $h=0.2$.}
	    \label{Fig:FE-Tri-Plot}
    \end{figure}
    Let $h>0$ denote the maximal length of the sides of the triangulation $\mathcal T^h$. Thus, if $h$ decreases, $\mathcal T^h$ is made finer. Let the angles of the triangulation be bounded from below by a positive constant, independently of $h$. If we choose piecewise linear \index{Method!finite element, FE}{\em finite elements (FE) functions} $\{\varphi_i\}_{i=1}^m$ as a basis for $V^h$, we get \eqref{Eq:ApproxProp-1} for $\mathscr W=H^2(\Omega)\cap V$; cf. \cite[pp.~105-110]{BS08}.\hfill$\blacklozenge$
\end{example}

\begin{remark}
    \rm Estimate \eqref{Eq:ApproxProp-1} is often ensured by assuming the existence of a family of projection or interpolation operators $\{\mathcal P^h\}_{h>0}$ with $\mathcal P^h:\mathscr W\to V^h$ satisfying
    \begin{equation}
        \label{Eq:ApproxProp-2}
        {\|v-\mathcal P^hv\|}_H+h\,{\|v-\mathcal P^hv\|}_V\le \hat C_{\mathscr W}\,h^2\,{\|v\|}_\mathscr W\quad\text{for any }v\in\mathscr W,
    \end{equation}
    where the constant $\hat C_{\mathscr W}$ is independent of $h$ and $v$.\hfill$\blacksquare$
\end{remark}

We apply a standard Galerkin scheme for \eqref{SIAM:Eq3.1.6}. Thus, we look for a function $y^h$ satisfying $y^h(t)\in V^h$ in $[0,T]$ almost everywhere and
\begin{subequations}
    \label{EvProGal}
    \begin{align}
        \label{EvProGal-1}
        \frac{\mathrm d}{\mathrm dt}{\langle y^h(t),\varphi\rangle}_H+a(t;y^h(t),\varphi)&={\langle (\mathcal F+\mathcal Bu)(t),\varphi\rangle}_{V',V}\quad\text{for all }\varphi \in V^h\text{ a.e. in }(0,T],\\
        \label{EvProGal-2}
        y^h(0)&=\mathcal P^h y_\circ,
    \end{align}
\end{subequations}
where $\mathcal P^h:X\to V^h$ is an appropriate linear projection operator for the respective choice $X = V$ or $X = H$.

\begin{example}
    \label{ExampleFEProjection}
    \rm Let us present two possible choices for the projection operator $\mathcal P^h$.
    \begin{enumerate}
        \item [1)] For any $w\in H$ the element $v^h=\mathcal P^h_H w\in V^h$ has to satisfy
        \begin{align*}
            {\langle v^h,\varphi_i^h\rangle}_H={\langle w,\varphi_i^h\rangle}_H\quad\text{for }i=1,\ldots,m.
        \end{align*} 
        Utilizing the representation $v^h=\sum_{j=1}^m\mathrm v_j^h\varphi_j^h\in V^h$ the coefficient vector $\mathrm v^h=(\mathrm v_i^h)\in \mathbb R^m$ is uniquely determined as the solution to the linear system
        \begin{equation}
            \label{ProjPh-H}
            \bM^h\mathrm v^h=\mathrm w^h
        \end{equation}
        with the \index{Matrix!mass!FE, $\bM^h$}{\em mass matrix} $\bM^h=(({\langle \varphi_j^h,\varphi_i^h \rangle}_H))\in\mathbb R^{m\times m}$ and the right-hand side $\mathrm w^h=(\langle w,\varphi_i^h\rangle_H)\in\mathbb R^m$. From \eqref{ProjPh-H} we infer that
        \begin{equation}
            \label{ProjPh-Prop-H}
            {\langle \mathcal P^h_H y_\circ,\varphi_i^h\rangle}_H=\sum_{j=1}^m\mathrm v_j^h\,{\langle \varphi_j^h,\varphi_i^h\rangle}_H={\langle y_\circ,\varphi_i^h\rangle}_H
        \end{equation}
        holds for $i=1,\ldots,m$. Furthermore, we have
        \begin{align*}
            {\|\mathcal P^h_Hw\|}_H^2={\langle\mathcal P^h_Hw,\mathcal P^h_Hw\rangle}_H={\langle w,\mathcal P^h_Hw\rangle}_H\le{\|w\|}_H{\|\mathcal P^h_Hw\|}_H,
        \end{align*}
        i.e., $\|\mathcal P^h_Hw\|_H\le\|w\|_H$ for every $w\in H$. Moreover, $\|\mathcal P^h_Hw^h\|_H=\|w^h\|_H$ for every $w^h\in V^h\subset H$. This implies
        \begin{equation}
            \label{Chopin-20}
            {\|\mathcal P^h_H\|}_{\mathscr L(H)}=1.
        \end{equation}
        From $\|\mathcal P^h_Hw\|_H\le\|w\|_H$ for every $w\in H$ and \eqref{Poincare} we also infer that
        \begin{align*}
            {\|\mathcal P^h_Hw\|}_H\le c_V\,{\|w\|}_V\quad\text{for every }w\in V
        \end{align*}
        which implies that $\|\mathcal P^h_H\|_{\mathscr L(V,H)}\le c_V$.
        \item [2)] We also introduce a linear projection $\mathcal P^h_V\in L(V,V^h)$ as follows: For any $w\in V$ the element $v^h=\mathcal P^h_V w\in V^h$ is given as
        \begin{align*}
            {\langle v^h,\varphi_i^h\rangle}_V={\langle w,\varphi_i^h\rangle}_V\quad\text{for }i=1,\ldots,m.
        \end{align*} 
        From $v^h=\sum_{j=1}^m\mathrm v_j^h\varphi_j^h\in V^h$ the coefficient vector $\mathrm v^h\in \mathbb R^m$ is the unique solution to the linear system
        \begin{equation}
            \label{ProjPh-V}
            \bS^h\mathrm v^h=\mathrm w^h
        \end{equation}
        with the symmetric, positive definite \index{Matrix!stiffness!FE, $\bS^h$}{\em stiffness matrix} $\bS^h=(({\langle \varphi_j^h,\varphi_i^h \rangle}_V))\in\mathbb R^{m\times m}$ and the right-hand side $\mathrm w^h=(\langle w,\varphi_i^h\rangle_V)\in\mathbb R^m$. Utilizing \eqref{ProjPh-V} we find that
        \begin{equation}
            \label{ProjPh-Prop-V}
            {\langle \mathcal P^h_V y_\circ,\varphi_i^h\rangle}_V=\sum_{j=1}^m\mathrm v_j^h\,{\langle \varphi_j^h,\varphi_i^h\rangle}_V={\langle y_\circ,\varphi_i^h\rangle}_V
        \end{equation}
        holds for $i=1,\ldots,m$. Analogously to part 1) we derive that
        \begin{equation}
            \label{Chopin-20-aaa}
            {\|\mathcal P^h_V\|}_{\mathscr L(V)}=1.
        \end{equation}
        Applying \eqref{Chopin-20-aaa} and the embedding inequality \eqref{Poincare} we also conclude that
        \begin{equation}
            \label{Chopin-21}
            \begin{aligned}
                {\|\mathcal P^h_V\|}_{\mathscr L(V,H)}&=\sup_{{\|w\|}_V=1}{\|\mathcal P^h_Vw\|}_H\le c_V\sup_{{\|w\|}_V=1}{\|\mathcal P^h_Vw\|}_V\le c_V\,{\|\mathcal P^h_V\|}_{\mathscr L(V)}=c_V
            \end{aligned}
        \end{equation}
        holds true.\hfill$\blacklozenge$
    \end{enumerate}
\end{example}

Since $y^h(t) \in V^h$ is assumed, we have the representation
\begin{equation}
    \label{FE-Galerkin}
    y^h(t)=\sum_{i=1}^m \mathrm y^h_i(t)\varphi_i^h\in V^h
\end{equation}
with a modal coefficient vector
\begin{align*}
    \mathrm y^h(t)=\big(\mathrm y_i^h(t)\big)_{1\le i\le m}\in\mathbb R^m \quad \text{for }t \in [0,T].
\end{align*}
From \eqref{EvProGal} and \eqref{FE-Galerkin} we derive the linear system of ordinary differential equations
\begin{equation}
    \label{FineModel}
    \bM^h \dot{\mathrm y}^h(t)+\bA^h(t) \mathrm y^h(t) =\mathrm g^h(t;u)\text{ f.a.a. } t \in (0,T], \quad \bM^h \mathrm y^h(0) = \mathrm y_\circ^h
\end{equation}
with the mass matrix $\bM^h$ introduced in Example~\ref{ExampleFEProjection}, the \index{Matrix!mass!FE, $\bA^h$}{\em stiffness matrix}
\begin{equation}
    \label{FEStiffnessMatrix}
    \bA^h(t)=\big(\big(a(t;\varphi_j^h,\varphi_i^h)\big)\big)\in\mathbb R^{m\times m}\quad\text{f.a.a. }t\in[0,T]
\end{equation}
and the vectors
\begin{align*}
    \mathrm y_\circ^h=\big({\langle y_\circ,\varphi_i^h\rangle}_H\big)_{1\le i\le m},\,\mathrm g^h(t;u)=\big({\langle (\mathcal F+\mathcal Bu)(t),\varphi_i^h\rangle}_{V',V}\big)_{1\le i\le m}\in\mathbb R^m
\end{align*}
for almost all $t \in [0,T]$ and all $u\in\U$. The next result follows by similar arguments used in the proof of Theorem~\ref{SIAM:Theorem3.1.1POD}.

\begin{theorem}
    \label{Theorem:FESystem}
    Let Assumption {\rm\ref{A1}} hold. Suppose that $\{\varphi_i\}_{i=1}^m$ are linearly independent functions in $V$ and that $V^h$ is given by \eqref{FESpace}. Then for every $u\in\U$ there is a unique weak solution $y^h\in H^1(0,T;V)\hookrightarrow\Y$ satisfying \eqref{EvProGal} and the a-priori bound
    \begin{equation}
        \label{APrioriFE}
        {\|y^h\|}_{\Y}\le C\left({\|\mathcal P^hy_\circ\|}_H+{\|\mathcal F\|}_{L^2(0,T;V')}+{\|u\|}_\U\right)
    \end{equation}
    for a constant $C>0$ which is independent of $u$, $y_\circ$ and $f$.
\end{theorem}

\begin{remark}
    \rm In Example~\ref{ExampleFEProjection} we have introduced the projection $\mathcal P^h_H$. Then it follows from
    \begin{align*}
        {\|\mathcal P^h_Hy_\circ\|}_H^2={\langle\mathcal P^h_Hy_\circ,\mathcal P^h_Hy_\circ\rangle}_H={\langle y_\circ,\mathcal P^h_Hy_\circ\rangle}_H\le{\|y_\circ\|}_H{\|\mathcal P^h_Hy_\circ\|}_H
    \end{align*}
    that $\|\mathcal P^h_Hy_\circ\|_H\le\|y_\circ\|_H$ holds so that \eqref{APrioriFE} can be replaced by
    \begin{align*}
        {\|y^h\|}_{\Y}\le C\left({\|y_\circ\|}_H+{\|\mathcal F\|}_{L^2(0,T;V')}+{\|u\|}_\U\right)
    \end{align*}
    where the right-hand side is independent of $h$. In the case of $\mathcal P^h=\mathcal P^h_V$ and $y_\circ\in V$ we find that $\|\mathcal P^h_Vy_\circ\|_H\le c_V\,\|y_\circ\|_V$, where the constant $c_V$ is given by \eqref{Poincare}. Thus, we infer from \eqref{APrioriFE} that
    \begin{align*}
        {\|y^h\|}_{\Y}\le\tilde C\left({\|y_\circ\|}_H+{\|\mathcal F\|}_{L^2(0,T;V')}+{\|u\|}_\U\right)
    \end{align*}
    holds with $\tilde C=C\max(1,c_V)$.\hfill$\blacksquare$
\end{remark}

\subsection{A-priori error analysis}
\label{SIAM-Book:Section3.4.2}

Now we estimate the error between the true solution $y$ of \eqref{SIAM:Eq3.1.6} and the semidiscrete solution $y^h$ of \eqref{EvProGal}. For the proof we refer to Section~\ref{SIAM-Book:Section3.8.3}.

\begin{theorem}
    \label{Theorem:FE-AprioriError}
    Let Assumptions~\emph{\ref{A1}} and \emph{\ref{A3}} hold. We suppose that $y$ and $y^h$ are the solutions to \eqref{SIAM:Eq3.1.6} and \eqref{EvProGal}, respectively. Moreover, assume that $y\in H^1(0,T;\mathscr W)$. Then there exists a constant $C>0$ such that the following \index{Error estimate!a-priori!state variable}a-priori error estimates hold:
    \begin{equation}
        \label{Eq:Apriori-FE-8-1}
        \begin{aligned}
            \int_0^T{\|y^h(t)-y(t)\|}_H^2\,\mathrm dt&\le Ch^2\left(h^2\,{\|y\|}_{H^1(0,T;\mathscr W)}^2+{\|y\|}_{L^2(0,T;\mathscr W)}^2\right),\\
            \int_0^T{\|y^h(t)-y(t)\|}_V^2\,\mathrm dt&\le Ch^2\left(h^2\,{\|y_t\|}_{L^2(0,T;\mathscr W)}^2+{\|y\|}_{L^2(0,T;\mathscr W)}^2\right).
        \end{aligned}
    \end{equation}
\end{theorem}

\subsection{POD basis computation}
\label{SIAM-Book:Section3.4.3}

This subsection focusses on the POD basis computation utilizing the semidiscrete solution $y^h \in H^1(0,T;V^h)$ and its derivative as data vectors. Suppose that Assumption~\ref{A2} is satisfied. By
\begin{equation}
    \label{FEAnsatz}
    y^h(t)=\sum_{i=1}^m\mathrm y_i^h(t)\varphi_i\quad\text{for all }t\in[0,T],
\end{equation}
we denote the solution to \eqref{EvProGal}. In the context of Section~\ref{Section:ContPODHilbert} we choose $K=2$, $\omega_1^K=\omega_2^K=1$, $y^1=y^h$ and $y^2=y_t^h$. From Theorem~\ref{Theorem:FESystem} we conclude that $y^h\in H^1(0,T;V^h)$. Thus, $y^1,\,y^2\in L^2(0,T;X)$ for $X=H$ and
\begin{equation}
    \label{FEPODEigPro}
    \mathcal R\psi_i=\lambda_i\psi_i\text{ for }i\in\mathbb I,\quad\lambda_1\ge\lambda_2\ge\ldots\ge\lambda_d>\lambda_{d+1}=\ldots =0
\end{equation}
where the finite-rank operator $\mathcal R:X\to V^h\subset X$ is given as
\begin{align*}
    \mathcal R\psi=\int_0^T\Big({\langle y^h(t),\psi\rangle}_X\,y^h(t)+{\langle y^h_t(t),\psi\rangle}_X\,y^h_t(t)\Big)\,\mathrm dt\quad\text{for }\psi\in X
\end{align*}
and $d=\dim \big\{\mathcal R\psi\,\big|\,\psi\in X\big\}\le m$. 

Let us mention that $\mathcal R$ and therefore $\{(\lambda_i,\psi_i)\}_{i\in\mathbb I}$ as well as $d$ depend on $h$ To simplify the notation we so not indicate this dependence by a superscript $h$.

Since $d$ is finite, $\mathcal R$ is a finite rank operator. It follows from Lemma~\ref{SIAM:Lemma3.2.1}-2) that $d$ is the same for the choices $X=H$ and $X=V$. Then $\{\psi_i\}_{i=1}^\ell\subset V^h$ is a POD basis of rank $\ell$ solving
\begin{equation}
    \label{POD:FEData}
    \left\{
    \begin{aligned}
        & \min\int_0^T \Big\|y^h(t) - \sum_{i=1}^\ell {\langle y^h(t),\psi_i\rangle}_X\,\psi_i\Big\|_X^2+\Big\|y_t^h(t) - \sum_{i=1}^\ell {\langle y_t^h(t),\psi_i\rangle}_X\,\psi_i\Big\|_X^2\,\mathrm dt\\
        &\hspace{1mm}\text{s.t. }\{\psi_i\}_{i=1}^\ell\subset X\text{ and }{\langle\psi_i,\psi_j\rangle}_X=\delta_{ij} \text{ for } 1 \le i,j \le \ell.
    \end{aligned}
    \right.
\end{equation}
For $\{(\lambda_i,\psi_i)\}_{i\in\mathbb I}$ satisfying \eqref{FEPODEigPro} we obtain for every $\ell\in\{1,\ldots,d\}$
\begin{equation}
    \label{FEApprPOD}
    \int_0^T \Big\|y^h(t) - \sum_{i=1}^\ell {\langle y^h(t),\psi_i\rangle}_X\,\psi_i\Big\|_X^2+\Big\|y_t^h(t)-\sum_{i=1}^\ell{\langle y_t^h(t),\psi_i\rangle}_X\,\psi_i\Big\|_X^2\,\mathrm dt=\sum_{i=\ell+1}^d\lambda_i.
\end{equation}
From $\psi_i=\mathcal R\psi_i/\lambda_i$ for $i=1,\ldots,\ell\le d^h$, we infer that there exists a unique coefficient matrix $\bPsi=[\uppsi_1\,|\ldots\,\uppsi_\ell]\in\mathbb R^{m\times\ell}$ with the elements $\Psi_{ij}=(\uppsi_j)_i$ and columns $\bPsi_{\cdot\,,j}=\uppsi_j\in\mathbb R^m$ such that
\begin{equation}
    \label{FEAnsatzPOD}
    \psi_i=\sum_{j=1}^m\Psi_{ji}\varphi_j\quad\text{for }i=1,\ldots,\ell.
\end{equation}
Let us introduce the positive definite, symmetric \index{Matrix!weighting}{\em weighting matrix}
\begin{align*}
    \bW=((\langle\varphi_j,\varphi_i\rangle_X))\in\mathbb R^{m\times m}
\end{align*}
and the associated weighted inner product
\begin{align*}
    {\langle \mathrm v,\mathrm w\rangle}_\bW=\mathrm v^\top\bW\mathrm w\quad\text{for }\mathrm  v,\,\mathrm w\in \mathbb{R}^m.
\end{align*}
with the induced norm $|\cdot|_\bW={\langle\cdot\,,\cdot\rangle}_\bW^{1/2}$. Then we conclude from \eqref{FEAnsatz} and \eqref{FEAnsatzPOD} that
\begin{align*}
    {\langle y^h(t),\psi_j\rangle}_X\,\psi_j=\sum_{i=1}^m{\langle\mathrm y^h(t),\uppsi_i\rangle}_\bW\,\Psi_{ij}\varphi_i\quad\text{for }j=1,\ldots,\ell.
\end{align*}
Hence, we find
\begin{align*}
    \Big\|y^h(t) - \sum_{i=1}^\ell {\langle y^h(t),\psi_i\rangle}_X\,\psi_i\Big\|_X^2=\Big|\mathrm y^h(t) - \sum_{i=1}^\ell {\langle \mathrm y^h(t),\uppsi_i\rangle}_\bW\,\uppsi_i\Big|_\bW^2
\end{align*}
and
\begin{align*}
    \Big\|y_t^h(t) - \sum_{i=1}^\ell {\langle y^h_t(t),\psi_i\rangle}_X\,\psi_i\Big\|_X^2=\Big|\dot{\mathrm y}^h(t) - \sum_{i=1}^\ell {\langle \dot{\mathrm y}^h(t),\uppsi_i\rangle}_{\bW}\,\uppsi_i\Big|_\bW^2
\end{align*}
Summarizing, \eqref{POD:FEData} is equivalent to the minimization problem
\begin{equation}
    \label{POD:FEDataDiscrete}
    \left\{
    \begin{aligned}
        & \min\int_0^T \Big|\mathrm y^h(t) - \sum_{i=1}^\ell {\langle \mathrm y^h(t),\uppsi_i\rangle}_\bW\,\uppsi_i\Big|_\bW^2+\Big|\dot{\mathrm y}^h(t) - \sum_{i=1}^\ell {\langle \dot{\mathrm y}^h(t),\uppsi_i\rangle}_\bW\,\uppsi_i\Big|_\bW^2\,\mathrm dt\\
        &\hspace{1mm}\text{s.t. }\{\uppsi_i\}_{i=1}^\ell\subset \mathbb R^m\text{ and }{\langle\uppsi_i,\uppsi_j\rangle}_\bW=\delta_{ij} \text{ for } 1 \le i,j \le \ell.
    \end{aligned}
    \right.
\end{equation}
If $\{\uppsi_i\}_{i=1}^\ell$ solves \eqref{POD:FEDataDiscrete} then $\{\psi_i^h\}_{i=1}^\ell$ is a solution to \eqref{POD:FEData} with
\begin{align*}
    \psi_i=\sum_{j=1}^m\big(\uppsi_i\big)_j\varphi_i\quad\text{for }i=1,\ldots,\ell.
\end{align*}
On the other hand, if $\{\psi_i\}_{i=1}^\ell$ is a solution to \eqref{POD:FEData} then we infer from $\psi_i\in V^h$ the existence of a uniquely determined matrix $\bPsi=[\uppsi_1\,|\ldots|\,\uppsi_\ell]\in\mathbb R^{m\times\ell}$ satisfying \eqref{FEAnsatzPOD}. It follows that the columns $\{\uppsi_i\}_{i=1}^\ell$ of $\bPsi$ solve \eqref{POD:FEDataDiscrete}.

\subsection{POD Galerkin scheme}
\label{SIAM-Book:Section3.4.4}

After having computed a POD basis, this subsection illurstrates the Galerkin discretization of the semidiscrete system onto the subspace spanned by this POD basis. Suppose therefore that we have computed a POD basis $\{\psi_i^h\}_{i=1}^\ell$ of rank $\ell$ with $\ell\le d^h\le m$. Let Assumption~\ref{A2} be satisfied. We set
\begin{align*}
    X^{h\ell}=\mathrm{span}\,\big\{\psi_1,\ldots,\psi_\ell\big\}\subset X\subset V^h
\end{align*}
and
\begin{align*}
    H^{h\ell}=\mathrm{span}\,\big\{\psi_1^H,\ldots,\psi_\ell^H\big\},\quad V^{h\ell}=\mathrm{span}\,\big\{\psi_1^V,\ldots,\psi_\ell^V\big\}
\end{align*}
for the choices $X=H$ and $X=V$, respectively. For the reader's convenience we omit the indication of the dependence of the eigenvectors $\{\psi_i^H\}_{i=1}^\ell$ and $\{\psi_i^V\}_{i=1}^\ell$ on the spatial discretization by an additional superscript $h$. Then we replace \eqref{EvProGal} by the following POD Galerkin scheme: for given $u\in\U$ find $y^{h\ell}:[0,T] \to X^{h\ell}$ satisfying
\begin{subequations}
\label{EvProGal-POD}
\begin{align}
\label{EvProGal-POD-1}
\frac{\mathrm d}{\mathrm dt} \, {\langle y^{h\ell}(t),\psi\rangle}_H+a(t;y^{h\ell}(t),\psi)&={\langle (\mathcal F+\mathcal Bu)(t),\psi\rangle}_{V',V}\\
\nonumber
&\hspace{20mm}\forall\psi \in X^{h\ell}\text{ a.e. in }(0,T],\\
\label{EvProGal-POD-2}
y^{h\ell}(0)&=\mathcal P^{h\ell} y_\circ,
\end{align}
\end{subequations}
where $\mathcal P^{h\ell}:Z\to X^{h\ell}$ stands for one of the following four orthogonal projections (cf. \eqref{SIAM:Eq3.2.12}):
\begin{align*}
    \mathcal P^{h\ell}=\left\{
    \begin{aligned}
        &\mathcal P^{h\ell}_H&&\text{for }X=H^{h\ell}\text{ and } Z=H,\\
        &\mathcal P^{h\ell}_V&&\text{for }X=H^{h\ell}\text{ and } Z=V,\\
        &\mathcal Q^{h\ell}_H&&\text{for }X=V^{h\ell}\text{ and } Z=H,\\
        &\mathcal Q^{h\ell}_V&&\text{for }X=V^{h\ell}\text{ and } Z=V.
    \end{aligned}
    \right.
\end{align*}

From $y^{h\ell}(t)\in X^{h\ell}$ it follows that there exists a coefficient vector
\begin{align*}
    \mathrm y^{h\ell}(t)=\big(\mathrm y_i^{h\ell}(t)\big)_{1\le i\le\ell}\in\mathbb R^\ell\quad\text{a.e. in }[0,T]
\end{align*}
satisfying
\begin{equation}
    \label{FE-POD-Galerkin}
    y^{h\ell}(t)=\sum_{i=1}^\ell\mathrm y_i^{h\ell}(t)\psi_i\quad\text{a.e. in }[0,T].
\end{equation}
Inserting \eqref{FE-POD-Galerkin} into \eqref{EvProGal-POD-1} and using $\psi=\psi_i$, $1\le i\le \ell$, we derive the following linear system of ordinary differential equations
\begin{subequations}
    \label{FE-POD-ODE}
    \begin{equation}
        \bM^{h\ell}\dot{\mathrm y}^{h\ell}(t)+\bA^{h\ell}(t)\mathrm y^{h\ell}(t)=\mathrm g^{h\ell}(t;u)\quad\text{a.e. in }(0,T]
    \end{equation}
    with the \index{Matrix!POD!mass, $\bM^{h\ell}$}\index{Matrix!POD!stiffness, $\bA^{h\ell}$}matrices
    \begin{align*}
        \bM^{h\ell}=\big(\big({\langle\psi_j,\psi_i\rangle}_H\big)\big)\in\mathbb R^{\ell\times\ell},\quad\bA^{h\ell}(t)=\big(\big(a(t;\psi_j,\psi_i)\big)\big)\in\mathbb R^{\ell\times\ell}
    \end{align*}
    and the vector
    \begin{align*}
        \mathrm g^{h\ell}(t;u)=\big(
        {\langle (\mathcal F+\mathcal Bu)(t),\psi_i\rangle}_{V',V}\big)_{1\le i\le\ell}\in\mathbb R^\ell\quad\text{for }u\in\U\text{ and a.e. in }[0,T].
    \end{align*}
    From \eqref{EvProGal-POD-2} we infer that
    \begin{equation}
        \mathrm y^{h\ell}(0)=\mathrm y_\circ^{h\ell}=\big(\mathrm y_{\circ i}^{h\ell}\big)_{1\le i\le\ell}\in\mathbb R^\ell\quad\text{with}\quad\mathcal P^{h\ell}y_\circ=\sum_{i=1}^\ell\mathrm y_{\circ i}^{h\ell}\psi_i.
    \end{equation}
\end{subequations}

Arguing as in the proof of Theorem~\ref{SIAM:Theorem3.1.1POD} we obtain the next result.

\begin{theorem}
    \label{DWW-1}
    Suppose that Assumption~{\rm\ref{A2}} holds. Moreover, assume that $\ell\le d_h$ holds. Then for any $u\in\U$ there exists a unique weak solution $y^{h\ell}\in H^1(0,T;V)$ satisfying \eqref{EvProGal-POD}. In particular, we have the a-priori bound
    \begin{align*}
        {\|y^{h\ell}\|}_{\Y}\le C\left({\|\mathcal P^{h\ell}y_\circ\|}_H+{\|\mathcal F\|}_{L^2(0,T;V')}+{\|u\|}_\U\right)
    \end{align*}
    for a constant $C>0$ which is independent of $\ell$, $y_\circ$, $\mathcal F$ and $u$.
\end{theorem}

The next result is a direct consequence of Theorem~\ref{DWW-1}; Corollaries~\ref{Corollary:HI-20} and \ref{Corollary:HI-21}.

\begin{corollary}
    \label{Corollary:HI-30}
    Suppose that Assumption~{\rm\ref{A2}} holds. Then we have for every $\ell\in\{1,\ldots,d\}$:
    \begin{enumerate}
        \item [\rm 1)] There exists a unique solution $\hat y^{h\ell}$ satisfying $\hat y^{h\ell}(t)\in X^\ell$ in $[0,T]$ almost everywhere and
        \begin{align*}
        \frac{\mathrm d}{\mathrm dt}\,{\langle\hat y^{h\ell}(t),\psi\rangle}_H+a(t;\hat y^{h\ell}(t),\psi)&={\langle\mathcal F(t),\psi\rangle}_{V',V}\quad\text{for all }\psi\in X^\ell\text{ a.e. in }(0,T],\\
        \label{MD-5b}
        \hat y^{h\ell}(0)&=\mathcal P^{h\ell} y_\circ.
    \end{align*}
    Moreover, there exists a constant $C\ge 0$ which is independent of $\ell$, $y_\circ$ and $\mathcal F$ with the a-priori bound
    \begin{align*}
        {\|\hat y^{h\ell}\|}_{\Y}\le C\big({\|\mathcal P^{h\ell}y_\circ\|}_H+{\|\mathcal F\|}_{L^2(0,T;V)}\big).
    \end{align*}
    \item [\rm 2)] Let us define the linear operator $\mathcal S^{h\ell}:\U\to \Y$ as follows: $y^{h\ell}=\mathcal S^{h\ell}u$ satisfies
    \begin{align*}
        \frac{\mathrm d}{\mathrm dt}\,{\langle y^{h\ell}(t),\psi\rangle}_H+a(t;y^{h\ell}(t),\psi)&={\langle (\mathcal Bu)(t),\psi\rangle}_{V',V}\quad\text{for all }\psi^h\in X^{h\ell}\text{ a.e. in }(0,T],\\
        y^{h\ell}(0)&=0.
    \end{align*}
    Then $\mathcal S^{h\ell}$ is well-defined and bounded, i.e.,
    \begin{align*}
        {\|\mathcal S^{h\ell}u\|}_{\Y}\le C\,{\|u\|}_\U
    \end{align*}
    for a constant $C\ge0$ which is independent of $\ell$ and $u$. Furthermore, $\mathcal S^{h\ell}u \in H^1(0,T;V)$ holds for all $u \in \U$.
    \end{enumerate}
\end{corollary}

\begin{remark}
    \label{Remark:Shl}
    \rm Let us mention that the solution to \eqref{EvProGal-POD} can be expressed as $y^{h\ell}=\hat y^{h\ell}+\mathcal S^{h\ell}u$.\hfill$\blacksquare$
\end{remark}

\subsection{POD a-priori error analysis}
\label{SIAM-Book:Section3.4.5}

As in Section~\ref{SIAM-Book:Section3.3.2} we are interested to estimate the error between the solutions $y$ and $y^{h\ell}$  to \eqref{SIAM:Eq3.1.6} and \eqref{EvProGal-POD}, respectively. Here we make use of the triangle inequality and Theorem~\ref{Theorem:FE-AprioriError}, where we already present an a-priori error estimate for the semidiscrete solution $y^h$ to \eqref{EvProGal}. Hence, if Assumptions~\ref{A1} and \ref{A3} hold we find that
\begin{align*}
    {\|y-y^{h\ell}\|}_{L^2(0,T;V)}^2 &\le 2\,{\|y-y^h\|}_{L^2(0,T;V)}^2 + 2\,{\|y^h-y^{h\ell}\|}_{L^2(0,T;V)}^2= \mathcal O\big(h^2\big)+2\,{\|y^h-y^{h\ell}\|}_{L^2(0,T;V)}^2.
\end{align*}
Thus, our goal is to derive an a-priori error estimate for the term
\begin{align*}
    {\|y^h-y^{h\ell}\|}_{L^2(0,T;V)}^2 = \int_0^T{\|y^h(t)-y^{h\ell}(t)\|}_V^2\,\mathrm dt.
\end{align*}
We can follow the analysis carried out in Section~\ref{SIAM-Book:Section3.3.2}, where we directly apply the error formulas of Section~\ref{Section:ContPODHilbert}. Let us we make use of the following decomposition (compare \eqref{EqDec})
\begin{align*}
    y^h(t)-y^{h\ell}(t)=y^h(t)-\mathcal P^{h\ell} y^h(t)+\mathcal P^{h\ell} y^h(t)-y^{h\ell}(t)=\varrho^{h\ell}(t)+\vartheta^{h\ell}(t)\text{ a.e. in }[0,T]
\end{align*}
with $\varrho^{h\ell}(t)=y^h(t)-\mathcal P^{h\ell} y^h(t)\in (X^{h\ell})^\bot$ and $\vartheta^{h\ell}(t)=\mathcal P^{h\ell} y^h(t)-y^{h\ell}(t)\in X^{h\ell}$. To achieve that
\begin{align*}
    \vartheta^{h\ell}(0)=\mathcal P^{h\ell} y^h(0)-y^{h\ell}(0)=\mathcal P^{h\ell}\big(\mathcal P^hy_\circ\big)-\mathcal P^{h\ell} y_\circ=0
\end{align*}
we choose the operators $\mathcal P^h$ appropriately; cf. Example~\ref{ExampleFEProjection}.

\begin{remark}
    \label{Rem:ProjChoice}
    \rm
    \begin{enumerate}
        \item [1)] {\em $X=H$, $\mathcal P^{h\ell}=\mathcal P^{h\ell}_H$ and $\mathcal P^h=\mathcal P^h_H$:} We infer from \eqref{ProjPh-Prop-H} and $\psi_i^H\in V^h$ for $1\le i\le \ell$ that
        \begin{align*}
            \mathcal P^{h\ell}_H\big(\mathcal P^h_Hy_\circ\big)=\sum_{i=1}^\ell{\langle\mathcal P^h_Hy_\circ,\psi_i^H\rangle}_H\,\psi_i^H=\sum_{i=1}^\ell{\langle y_\circ,\psi_i^H\rangle}_H\,\psi_i^H=\mathcal P^{h\ell}_Hy_\circ
        \end{align*}
        which implies $\vartheta^{h\ell}(0)=0$.
        \item [2)] {\em $X=H$, $\mathcal P^{h\ell}=\mathcal P^{h\ell}_V$ and $\mathcal P^h=\mathcal P^h_V$:} We suppose that $y_\circ\in V$ holds. Then
        \begin{align*}
            \mathcal P^{h\ell}_V\big(\mathcal P^h_Vy_\circ\big)=\sum_{i=1}^\ell\mathrm v_i\psi_i^H,\quad\mathcal P^{h\ell}_Vy_\circ=\sum_{i=1}^\ell\mathrm w_i\psi_i^H
        \end{align*}
        for coefficient vectors $\mathrm v=(\mathrm v_i),\,\mathrm w=(\mathrm w_i)\in\mathbb R^\ell$. Utilizing $\psi_i^H\in V^h$ for $1\le i\le \ell$, \eqref{SIAM:Eq3.2.14b}, \eqref{ProjPh-Prop-V} and again \eqref{SIAM:Eq3.2.14b} we find
        \begin{align*}
            \big(\bS^{h\ell}\mathrm v\big)_i&=\sum_{j=1}^\ell{\langle \psi_j^H,\psi_i^H\rangle}_V\mathrm v_j={\langle \mathcal P^{h\ell}_V(\mathcal P^h_Vy_\circ),\psi_i^H\rangle}_V={\langle \mathcal P^h_Vy_\circ,\psi_i^H\rangle}_V={\langle y_\circ,\psi_i^H\rangle}_V\\
            &={\langle \mathcal P^{h\ell}_Vy_\circ,\psi_i^H\rangle}_V=\sum_{j=1}^\ell{\langle \psi_j^H,\psi_i^H\rangle}_V\mathrm w_j=\big(\bS^{h\ell}\mathrm w\big)_i\quad\text{for }i=1,\ldots,\ell,
        \end{align*}
        where $\bS^{h\ell}=((\langle \psi_j^H,\psi_i^H\rangle_V))\in\mathbb R^{\ell\times\ell}$ is invertible. This implies that $\mathrm v=\mathrm w$ holds. Consequently, $\vartheta^{h\ell}(0)=0$.
        \item [3)] {\em $X=V$, $\mathcal P^{h\ell}=\mathcal Q^{h\ell}_H$ and $\mathcal P^h=\mathcal P^h_H$:} We argue similar to case 2). There exist coefficient vectors $\mathrm v=(\mathrm v_i),\,\mathrm w=(\mathrm w_i)\in\mathbb R^\ell$ such that
        \begin{align*}
            \mathcal Q^{h\ell}_H\big(\mathcal P^h_Hy_\circ\big)=\sum_{i=1}^\ell\mathrm v_i\psi_i^V,\quad\mathcal Q^{h\ell}_Hy_\circ=\sum_{i=1}^\ell\mathrm w_i\psi_i^V
        \end{align*}
        holds. From $\psi_i^V\in V^h$ for $1\le i\le \ell$, \eqref{SIAM:Eq3.2.14a}, \eqref{ProjPh-Prop-H} and again \eqref{SIAM:Eq3.2.14a} we derive
        \begin{align*}
            \big(\bM^{h\ell}\mathrm v\big)_i={\langle \mathcal P^h_Hy_\circ,\psi_i^V\rangle}_H={\langle y_\circ,\psi_i^V\rangle}_H=\big(\bM^{h\ell}\mathrm w\big)_i\quad\text{for }i=1,\ldots,\ell,
        \end{align*}
        where $\bM^{h\ell}=((\langle \psi_j^V,\psi_i^V\rangle_H))\in\mathbb R^{\ell\times\ell}$ is invertible. Again, $\mathrm v=\mathrm w$ is satisfied which gives $\vartheta^{h\ell}(0)=0$.
        \item [4)] {\em $X=V$, $\mathcal P^{h\ell}=\mathcal Q^{h\ell}_V$ and $\mathcal P^h=\mathcal P^h_V$:} It follows from Assumption~\ref{A2}-3) that $y_\circ\in V$ holds for $X=V$. Due to \eqref{ProjPh-Prop-V} we have
        \begin{align*}
            \mathcal Q^{h\ell}_V\big(\mathcal P^h_Vy_\circ\big)=\sum_{i=1}^\ell{\langle\mathcal P^h_Vy_\circ,\psi_i^V\rangle}_V\,\psi_i^V=\sum_{i=1}^\ell{\langle y_\circ,\psi_i^V\rangle}_V\,\psi_i^V=\mathcal Q^{h\ell}_Vy_\circ.
        \end{align*}
        Again, we obtain $\vartheta^{h\ell}(0)=0$.\hfill$\blacksquare$
    \end{enumerate}
\end{remark}

Now we can argue as in the proof of Theorem~\ref{Th:A-PrioriError}.

\begin{theorem}
    \label{Th:A-PrioriError-20}
    Let Assumption~{\rm\ref{A2}} hold. We utilize the unique solution $y^h$ to \eqref{EvProGal} to compute the POD basis by solving \eqref{POD:FEDataDiscrete}. Let the eigenvalue-eigenfunction pairs $\{(\lambda_i^H,\psi_i^H)\}_{i\in\mathbb I}$ and $\{(\lambda_i^V,\psi_i^V)\}_{i\in\mathbb I}$ satisfy \eqref{FEPODEigPro} either for $X=H$ or for $X=V$, respectively. Furthermore, $y^{h\ell}$ stands for the solution to \eqref{EvProGal-POD}. Then we have for all $\ell\in\{1,\ldots,d\}$ the following \index{Error estimate!a-priori!state variable}a-priori error estimate
    \begin{align*}
        {\|y^h-y^{h\ell}\|}_\mathscr Y^2\le C\cdot\left\{
        \begin{aligned}
            &\sum_{i=\ell+1}^d\lambda_i^H\big({\|\psi_i^H\|}_V^2+1\big),&&X=H,~\mathcal P^{h\ell}=\mathcal P^\ell_H,~\mathcal P^h=\mathcal P^h_H,\\
            &\sum_{i=\ell+1}^d\lambda_i^H\,{\|\psi_i^H-\mathcal P^\ell_V\psi_i^H\|}_V^2,&&X=H,~\mathcal     P^{h\ell}=\mathcal P^{h\ell}_V,~\mathcal P^h=\mathcal P^h_V,~y_\circ\in V,\\
            &\sum_{i=\ell+1}^d\lambda_i^V\,{\|\psi_i^V-\mathcal Q_H^{h\ell}\psi_i^V\|}_V^2,&&X=V,~\mathcal P^{h\ell}=\mathcal Q^{h\ell}_H,~\mathcal P^h=\mathcal P^h_H,\\
            &\sum_{i=\ell+1}^d\lambda_i^V,&&X=V,~\mathcal P^{h\ell}=\mathcal Q^\ell_V,~\mathcal P^h=\mathcal P^h_V
        \end{aligned}
        \right.
    \end{align*}
    for a constant $C>0$ that is independent of $\ell$, $y_\circ$, $\mathcal F$ and $u$.
\end{theorem}

\subsubsection{Extensions}

Next we consider the case where the time derivative $y_t^h\in L^2(0,T;V)$ is not included in the snapshot set for the computation of the POD basis; cf. Section~\ref{SIAM-Book:Section3.3.2}. Then we determine the POD basis $\{\psi_i\}_{i=1}^\ell$ from the minimization problem
\begin{equation}
    \label{POD:FEDataDQ}
    \left\{
    \begin{aligned}
        & \min\int_0^T \Big\|y^h(t) - \sum_{i=1}^\ell {\langle y^h(t),\psi_i\rangle}_X\,\psi_i^h\Big\|_X^2\,\mathrm dt\\
        &\hspace{1mm}\text{s.t. }\{\psi_i\}_{i=1}^\ell\subset X\text{ and }{\langle\psi_i,\psi_j\rangle}_X=\delta_{ij} \text{ for } 1 \le i,j \le \ell.
    \end{aligned}
    \right.
\end{equation}
The solution to \eqref{POD:FEDataDQ} can be computed by solving
\begin{align*}
    \mathcal R\psi_i=\lambda_i\psi_i\text{ for }i\in\mathbb I,\quad\lambda_1\ge\ldots\ge\lambda_d>\lambda_{d+1}=\ldots=0
\end{align*}
with
\begin{align*}
    \mathcal R\psi=\int_0^T{\langle y^h(t),\psi\rangle}_X\,y^h(t)\,\mathrm dt\quad\text{for }\psi\in X
\end{align*}
and $d=\dim \big\{\mathcal R\psi\,\big|\,\psi\in X\big\}\le m$. Then it follows by the same arguments as in the proof of Theorem~\ref{Th:A-PrioriError-2} that there exists a constant $C>0$ that is independent of $\ell$, $y_\circ$, $f$ and $u$ satisfying the \index{Error estimate!a-priori!state variable}a-priori error estimates
\begin{equation}
    \label{Chopin-2}
    \int_0^T{\|y^h(t)-y^{h\ell}(t)\|}_V^2\,\mathrm dt\le C\cdot\left\{
    \begin{aligned}
        &\sum_{i=\ell+1}^d\lambda_i^H\,\big\|\psi_i^h\big\|_V^2,&&X=H,~\mathcal P^\ell=\mathcal P^\ell_H,~\mathcal P^h=\mathcal P^h_H,\\
        &\sum_{i=\ell+1}^d\lambda_i^V\,{\|\psi_i^V-\mathcal Q^\ell_H\psi_i^V\|}_V^2,&&X=V,~\mathcal P^\ell=\mathcal Q^\ell_H,~\mathcal P^h=\mathcal P^h_V,\\
        &&&y_\circ\in V.
    \end{aligned}
    \right.
\end{equation}
To prove \eqref{Chopin-2} we utilize that $\mathcal P^\ell_H$ and $\mathcal Q^\ell_H$ are $H$-orthonormal projections. This allows us to derive the property \eqref{Eq:SnapOhneDQ-2} for $y^h$. Let us mention that we can also derive bounds for the difference $y^h-y^{h\ell}$ in the $\Y$-norm proceeding as in the end of Section~\ref{SIAM-Book:Section3.3.2}.

\subsection{POD a-posteriori error analysis}
\label{SIAM-Book:Section3.4.6}

In general, the a-priori error analysis presented in Theorem~\ref{Th:A-PrioriError-20} can not be utilized directly in numerical realizations for two reasons:
\begin{itemize}
    \item The rate of convergence results are only available for a given $u\in\U$ if the corresponding finite element trajectories $y^h(t)$ and $y^h_t(t)$ are available. In this case, it is not meaningful to compute a POD solution $y^{h\ell}$ in addition.
    \item The constants in Theorem~\ref{Th:A-PrioriError-20} are often difficult to estimate accurately, so that the bounds may become inefficient.
\end{itemize}
For that reason, an a-posteriori error analysis is developed which allows us to estimate the error between the finite element and the POD solution for a given $u\in\U$. An essential point is that the computation of the POD basis does not require the solution $y^h(t)$. The proof of the next theorem is given in Section~\ref{SIAM-Book:Section3.8.3}.

\begin{theorem}
    \label{TheoremApostSemi}
    Suppose that Assumption~{\rm\ref{A2}} holds and that $u\in\U$. Let $\{\psi_i\}_{i=1}^\ell$ be any POD basis of rank $\ell$. Let $y^h$ and $y^{h\ell}$ be the solutions to \eqref{EvProGal} and \eqref{EvProGal-POD}, respectively. We define the residual
    \begin{align*}
        r^{h\ell}(t;u)={\langle y^{h\ell}_t(t;u)-(\mathcal F+\mathcal Bu)(t),\cdot\rangle}_{V',V}+a(t;y^{h\ell}(t;u),\cdot)\in (V^h)'
    \end{align*}
    for almost all $t\in[0,T]$. Then the following \index{Error estimate!a-posteriori!state variable}a-posteriori error estimates hold
    \begin{align*}
        {\|y^h(t)-y^{h\ell}(t)\|}_H^2\le e^{ct}\bigg({\|\mathcal P^hy_\circ-\mathcal P^\ell y_\circ\|}_H^2+\frac{1}{\gamma_1}\,{\|r^{h\ell}(\cdot\,;u)\|}_{L^2(0,t;(V^h)')}^2\bigg)
    \end{align*}
    almost everywhere in $[0,T]$ with $c=\max(2\gamma_2-\gamma_1/c_V^2,0)$ and
    \begin{align*}
        {\|y^h-y^{h\ell}\|}_{L^2(0,T;V)}^2\le \frac{1+2\gamma_2Te^{ct}}{\gamma_1}\bigg({\|\mathcal P^h y_\circ-\mathcal P^\ell y_\circ\|}_H^2+\frac{1}{\gamma_1}\,{\|r^{h\ell}(\cdot\,;u)\|}_{L^2(0,T;(V^h)')}^2\bigg).
    \end{align*}
\end{theorem}

\section{The fully discrete approximation}
\label{SIAM-Book:Section3.5}
\setcounter{equation}{0}
\setcounter{theorem}{0}
\setcounter{figure}{0}
\setcounter{run}{0}

In real computations both the solutions $y^h$ to \eqref{EvProGal} and $y^{h\ell}$ to \eqref{EvProGal-POD} can not be computed numerically for any $t\in[0,T]$. This is why we introduce a discretization of the temporal variable to get a fully discrete scheme for the evolution problem \eqref{ExWeakForm}.

\subsection{Temporal and Galerkin discretization}
\label{SIAM-Book:Section3.5.1}

Let $0=t_1 < t_2 < \ldots < t_n=T$ be a given grid in $[0,T]$ with step sizes $\delta t_j=t_j-t_{j-1}$ for $2 \le j \le n$. We set
\begin{align*}
    \delta t=\min_{2\le j\le n}\delta t_j\quad\text{and}\quad\Delta t=\max_{2\le j\le n}\delta t_j.
\end{align*}
We make use of the next hypotheses.

\begin{assumption}
    \label{A100}
    \begin{enumerate}
        \item [\em 1)] Assumptions~{\em\ref{A2}} and {\rm\ref{A3}} hold.
        \item [\em 3)] There exists a constant $\zeta>0$ which is independent of $n$, such that
        \begin{equation}
            \label{Chopin-3}
            \frac{\Delta t}{\delta t}\le\zeta,
        \end{equation}
        i.e., we have $\Delta t=\mathcal O(\delta t)$.
    \end{enumerate}
\end{assumption}

Let us introduce the finite difference quotient
\begin{align*}
    \overline\partial y^h_j=\frac{y_j^h-y_{j-1}^h}{\delta t_j}\quad\text{for }j=2,\ldots,n
\end{align*}
for the temporal derivative, where $\{y_j^h\}_{j=1}^n\subset V^h$ are approximations for the solution $y^h$ to \eqref{EvProGal} at the time instances  $t=t_j$ for $j=1,\ldots,n$. Moreover, we define for any $u\in\U$ the dual elements $\{g_j(u)\}_{j=1}^n\subset V'$ by
\begin{align*}
    g_j(u)&=\frac{2}{\delta t_1}\int_0^{\delta t_1/2}(\mathcal F+\mathcal Bu)(s)\,\mathrm ds,\\
    g_j(u)&=\frac{2}{\delta t_{j-1}+\delta t_j}\int_{t_j-\delta t_{j-1}/2}^{t_j+\delta t_j/2}(\mathcal F+\mathcal Bu)(s)\,\mathrm ds\quad\text{for }j=2,\ldots,n-1,\\
    g_n(u)&=\frac{2}{\delta t_n}\int_{T-\delta t_n/2}(\mathcal F+\mathcal Bu)(s)\,\mathrm ds.
\end{align*}
To solve \eqref{EvProGal} by time discretization, we apply the so-called \emph{$\theta$-scheme}\index{Method!$\theta$-scheme} with $\theta\in[0,1]$, for which the sequence $\{y^h_j\}_{j=1}^n \subset V^h$ has to satisfy
\begin{subequations}
    \label{EvProGal-disc}
    \begin{equation}
        \label{EvProGal-disc-1}
        {\langle \overline \partial y^h_j,\varphi^h \rangle}_H+\theta a(t_j;y^h_j,\varphi^h)+(1-\theta)a(t_{j-1};y^h_{j-1},\varphi^h)={\langle \theta g_j+(1-\theta)g_{j-1},\varphi^h\rangle}_{V',V}
    \end{equation}
    for all $\varphi^h \in V^h$, for $2 \le j \le n$ and
    \begin{equation}
        \label{EvProGal-disc-2}
        y^h_1=\mathcal P^h y_\circ,
    \end{equation}
\end{subequations}
where $\mathcal P^h:H\to V^h$ is an appropriate linear projection operator (see Example~\ref{ExampleFEProjection}). If we choose $\theta=1$, the $\theta$-method yields the {\em implicit Euler method}\index{Method!implicit Euler}, whereas the {\em explicit Euler method}\index{Method!explicit Euler} is obtained if we set $\theta=0$. Taking $\theta=1/2$ leads to the {\em Crank-Nicolson method}\index{Method!Crank-Nicolson}.

From $y_j^h\in V^h$, $j=1,\ldots,n$, we infer that
\begin{equation}
    \label{ThetaM:GalAnsatz}
    y^h_j=\sum_{l=1}^m\mathrm y_{jl}^h\varphi_l^h\in V^h\quad\text{for }1\le j\le n
\end{equation}
with coefficient vectors $\mathrm y^h_j=(\mathrm y_{jl}^h)_{1\le l\le m}\in\mathbb R^m$. Inserting \eqref{ThetaM:GalAnsatz} into \eqref{EvProGal-disc} and choosing $\varphi^h=\varphi_i^h$, $i=1,\ldots,m$, we find that the sequence $\{\mathrm y_j^h\}_{j=1}^n\subset\mathbb R^m$ of coefficient vectors satisfies the following initial value problem
\begin{subequations}
    \label{FineModel-Disc}
    \begin{align}
        \label{FineModel-Disc-1}
        \big(\bM^h+\theta\delta t_j\bA^h_j\big)\mathrm y_j^h&=\big(\bM^h+(\theta-1)\delta t_j\bA^h_{j-1}\big)\mathrm y^h_{j-1}+\delta t_j\big(\theta\mathrm g_j^h(u)+(1-\theta)\mathrm g_{j-1}^h(u)\big)\\
        \label{FineModel-Disc-2}
        \bM^h\mathrm y^h_1&=\mathrm y_\circ^h
    \end{align}
\end{subequations}
for $2 \le j \le n$, where $\bM^h$ has been introduced in Example~\ref{ExampleFEProjection}, $\bA^h_j=\bA^h(t_j)$ (cf. \eqref{FEStiffnessMatrix}) and $\mathrm g_j^h(u)=({\langle g_j(u),\varphi_i^h\rangle}_{V',V})_{1\le i\le m}$ for $1\le i\le m$. To ease the presentation we consider only the case $\theta=1$, i.e., the implicit Euler method. We refer the reader for the Crank-Nicolson method to \cite{KV02b}. The next result is proved in Section~\ref{SIAM-Book:Section3.8.4}.

\begin{theorem}
    \label{FullDiscFEModel}
    Suppose that Assumption~{\rm\ref{A100}} and $\theta=1$ hold. Then there exists a unique sequence $\{y^h_j\}_{j=1}^n \subset V^h$ to \eqref{EvProGal-disc}. If $\delta t$ is sufficiently small, we have the a-priori bounds
    \begin{equation}
        \label{TempDiscFESOL-APriori}
        {\|y_j^h\|}_H^2\le e^{-C_1(j-1)\delta t}\bigg({\|\mathcal P^hy_\circ\|}_H^2+\frac{1}{\gamma_1}\sum_{l=2}^j\delta t_l\,{\|g_l(u)\|}_{V'}^2\bigg)
    \end{equation}
    for $j=2,\ldots,n$ with the constant $C_1=2\gamma_1/c_V^2$ and
    \begin{equation}
        \label{TempDiscFESOL-APriori-2}
        {\|y_n^h\|}_H^2+\sum_{j=2}^n{\|y_j^h-y_{j-1}^h\|}_H^2+\gamma_1\sum_{j=2}^n\delta t_j\,{\|y_j^h\|}_V^2\le{\|\mathcal P^hy_\circ\|}_H^2+\frac{1}{\gamma_1}\sum_{j=2}^n\delta t_j\,{\|g_j(u)\|}_{V'}^2.
    \end{equation}
\end{theorem}

\subsection{A-priori error analysis}
\label{SIAM-Book:Section3.5.2}

In the following theorem we estimate the error between the exact solution and its approximation computed by the implicit Euler method together with a Galerkin discretization. The proof is given in Section~\ref{SIAM-Book:Section3.8.4}.

\begin{theorem}
    \label{Theorem:FEAprioriError}
    Suppose that Assumption~{\rm\ref{A100}} is satisfied and that $\theta=1$ is chosen. Assume that $y$ and $y^h$ are the solutions to \eqref{SIAM:Eq3.1.6} and \eqref{EvProGal-disc}, respectively, with $y\in H^2(0,T;\mathscr W)$. Let $\mathcal P^h$ satisfy \eqref{Eq:ApproxProp-2} and $\|\mathcal P^h\|_{\mathscr L(X,H)}\le C$ for a constant $C$ independent of $h$. Then there exists a constant $C>0$ such that the \index{Error estimate!a-priori!state variable}a-priori error estimate
    \begin{align*}
        \sum_{j=1}^n\delta t_j\,{\|y(t_j)-y^h_j\|}_V^2\le C\,\big(\Delta t^2+h^2+\Delta t^2h^2\big)
    \end{align*}
    holds.
\end{theorem}

\begin{remark}
    \rm From Example~\ref{ExampleFEProjection} it follows that
    \begin{align*}
        {\|\mathcal P^h\|}_{\mathscr L(X,H)}\le c_V\quad\left\{
        \begin{aligned}
            &\text{for }X=H\text{ and }\mathcal P^h=\mathcal P^h_H,\\
            &\text{for }X=V\text{ and }\mathcal P^h=\mathcal P^h_V.
        \end{aligned}
        \right.
    \end{align*}
    This, $\|\mathcal P^h\|_{\mathscr L(X,H)}$ be bounded independent of $h$.\hfill$\blacksquare$
\end{remark}

\subsection{POD basis computation}
\label{SIAM-Book:Section3.5.3}

Let Assumption~{\rm\ref{A100}} hold. By
\begin{equation}
    \label{FullyDiscGal}
    y^h_j=\sum_{i=1}^m\mathrm y^h_{ji}\varphi_i^h\quad\text{for }1\le j\le n
\end{equation}
we denote the Galerkin solution to \eqref{EvProGal-disc} for $\theta=1$. We introduce the coefficient vectors $\mathrm y^h_j=(\mathrm y^h_{ji})_{1\le i\le m}\in\mathbb R^m$ for $j=1,\ldots,n$. Let us choose $K=2$, $\omega_1^K=\omega_2^K=1$ and
\begin{align*}
    \alpha_1^n=\frac{\delta t_1}{2},\quad\alpha_j^n=\frac{\delta t_j+\delta t_{j+1}}{2}\text{ for }j=2,\ldots,n-1,\quad\alpha_n^n=\frac{\delta t_n}{2}.
\end{align*}
It follows from Assumption~\ref{A100}-3) and $\Delta t\ge \alpha_j^n$ that
\begin{equation}
    \label{Chopin-4}
    \delta t_j\ge\delta t=\frac{\delta t\Delta t}{\Delta t}\ge\frac{\Delta t}{\zeta}\ge\frac{\alpha_j^n}{\zeta}\quad\text{for }j=1,\ldots,n.
\end{equation}
We further choose the \index{POD space!discrete variant!snapshots}{\em snapshots}
\begin{align*}
    y_j^1=y_j^h\in V^h\subset X \text{ for }1\le j\le n\quad\text{and}\quad y_1^2=0,~y_j^2=\overline\partial y_j^h\in V^h\subset X \text{ for } j=2,\ldots,n
\end{align*}
with the difference quotient $\overline\partial y_j^h=(y_j^h-y_{j-1}^h)/\delta t_j$. Then
\begin{align*}
    &\sum_{j=1}^n\alpha_j^n\Big\|y_j^h-\sum_{i=1}^\ell{\langle y_j^h,\psi_i^n\rangle}_X\,\psi_i^n\Big\|_X^2+\sum_{j=2}^n\alpha_j^n\Big\|\overline\partial y_j^h-\sum_{i=1}^\ell{\langle \overline\partial y_j^h,\psi_i^n\rangle}_X\,\psi_i^n\Big\|_X^2\\
    &=\sum_{k=1}^2\sum_{j=1}^n\alpha_j^n\,\Big\|y_j^k-\sum_{i=1}^\ell{\langle y_j^k,\psi_i^n\rangle}_X\,\psi_i^n\Big\|_X^2.
\end{align*}
We define the \index{POD space!discrete variant!snapshot space}{\em snapshot space}
\begin{align*}
    \mathscr V^n=\Span\big\{y_j^k\in V\,\big|\,1\le j\le n\text{ and }1\le k\le 2\big\}
\end{align*}
and set $d^n=\dim\mathscr V^n\le2n$. Now we compute a POD basis of rank $\ell\in\{1,\ldots,d^n\}$ by solving
\begin{equation}
    \label{PODProb}
    \left\{
    \begin{aligned}
        &\min\sum_{k=1}^2\sum_{j=1}^n\alpha_j^n\,\Big\|y_j^k-\sum_{i=1}^\ell{\langle y_j^k,\psi_i^n\rangle}_X\,\psi_i^n\Big\|_X^2\\
        &\hspace{1mm}\text{s.t. }\{\psi_i^n\}_{i=1}^\ell\subset V^h\text{ and }{\langle\psi_i^n,\psi_j^n\rangle}_X=\delta_{ij},~1\le i,j\le \ell
    \end{aligned}
    \right.
\end{equation}
for $X=H$ and $X=V$. Utilizing the symmetric and positive weighting matrix $\bW=((\langle\varphi_j^h,\varphi_i^h\rangle_X))\in\mathbb R^{m\times m}$ and the associated weighted inner product we can express the minimization problem as
\begin{align}
    \label{PODProbDisc}
    \left\{
    \begin{aligned}
        &\min\sum_{j=1}^n\alpha_j^n\Big|\mathrm y_j^h-\sum_{i=1}^\ell{\langle \mathrm y_j^h,\uppsi_i^n\rangle}_\bW\,\uppsi_i^n\Big|_\bW^2+\sum_{j=2}^n\alpha_j^n\Big|\overline\partial\mathrm y_j^h-\sum_{i=1}^\ell{\langle \overline\partial\mathrm y_j^h,\uppsi_i^n\rangle}_\bW\,\uppsi_i^n\Big|_\bW^2\\
        &\hspace{1mm}\text{s.t. }\{\uppsi_i^n\}_{i=1}^\ell\subset\mathbb R^m\text{ and }{\langle\uppsi_i^n,\uppsi_j^n\rangle}_\bW=\delta_{ij},\quad1\le i,j\le \ell,
    \end{aligned}
    \right.
\end{align}
where now $\overline\partial\mathrm y_j^h=(\mathrm y_j^h-\mathrm y_{j-1}^h)/\delta t_j\in\mathbb R^m$. 

Recall a result from Section~\ref{Section:DiscPODHilbert}: Let the eigenvalue-eigenfunction pairs $\{(\lambda_i^{nH},\psi_i^{nH})\}_{i\in\mathbb I}$ and $\{(\lambda_i^{nV},\psi_i^{nV})\}_{i\in\mathbb I}$ satisfy
\begin{align*}
    \mathcal R_H^n\psi_i^{nH}=\lambda_i^{nH}\psi_i^{nH}\text{ for }i\in\mathbb I,&\quad\lambda_1^{nH}\ge\ldots\ge\lambda_{d^n}^H>\lambda_{d^n+1}^{nH}=\ldots=0,\\
    \mathcal R_V^n\psi_i^{nV}=\lambda_i^{nV}\psi_i^{nV}\text{ for }i\in\mathbb I,&\quad\lambda_1^{nV}\ge\ldots\ge\lambda_{d^n}^{nV}>\lambda_{d^n+1}^{nV}=\ldots=0,
\end{align*}
where
\begin{align*}
    \mathcal R_H^n\psi&=\sum_{j=1}^n\alpha_j^n\,{\langle y^h_j,\psi\rangle}_Hy_j^h+\sum_{j=2}^n\alpha_j^n\,{\langle\overline\partial y^h_j,\psi\rangle}_H\overline\partial y_j^h&&\text{for }\psi\in H,\\
    \mathcal R_V^n\psi&=\sum_{j=1}^n\alpha_j^n\,{\langle y^h_j,\psi\rangle}_Vy_j^h+\sum_{j=2}^n\alpha_j^n\,{\langle\overline\partial y^h_j,\psi\rangle}_V\overline\partial y_j^h&&\text{for }\psi\in V.
\end{align*}
Then $\{\psi_i^{nH})\}_{i=1}^\ell$ solves \eqref{PODProb} for $X=H$ and $\{\psi_i^{nV})\}_{i=1}^\ell$ is the solution to \eqref{PODProb} for $X=V$. Moreover, we find the POD approximation error as follows
\begin{equation}
    \label{POD-Error-Formula-100}
    \begin{aligned}
        \sum_{k=1}^2\sum_{j=1}^n\alpha_j^n\,\Big\|y_j^k-\sum_{i=1}^\ell{\langle y_j^k,\psi_i^{nH}\rangle}_H\,\psi_i^{nH}\Big\|_H^2&=\sum_{i=\ell+1}^{d^n}\lambda_i^{nH},\\
        \sum_{k=1}^2\sum_{j=1}^n\alpha_j^n\,\Big\|y_j^k-\sum_{i=1}^\ell{\langle y_j^k,\psi_i^{nV}\rangle}_V\,\psi_i^{nV}\Big\|_V^2&=\sum_{i=\ell+1}^{d^n}\lambda_i^{nV}.
    \end{aligned}
\end{equation}

\begin{remark}
    \rm Let us motivate why we include the finite difference quotients $\{\overline\partial y_j^h\}_{j=2}^n$ into the set of snapshots; cf. \cite[Remark~1]{KV02a}. While the finite difference quotients are contained in the span of $\{y_j^h\}_{j=1}^n$, the POD bases differ depending on whether $\{\overline\partial y_j^h\}_{j=2}^n$ are included or not. The linear dependence does not constitute a difficulty for the singular value decomposition which is required to compute the POD basis. In fact, the snapshots themselves can be linearly dependent. Secondly, in anticipation of the rate of convergence results that will be obtained further below, we note that the time derivative of $y$ in \eqref{EvProGal-1} must be approximated by the Galerkin-POD based scheme. In case the terms $\{\overline\partial y_j^h\}_{j=2}^n$ are included in the snapshot ensemble, we will be able to utilize the estimate
    \begin{equation}
        \label{DQincl-1}
        \sum_{j=1}^n\alpha_j^n\Big\|y_j^h-\sum_{i=1}^\ell{\langle y_j^h,\psi_i^h\rangle}_X\,\psi_i^h\Big\|_X^2+\sum_{j=2}^n\alpha_j^n\Big\|\overline\partial y_j^h-\sum_{i=1}^\ell{\langle \overline\partial y_j^h,\psi_i^n\rangle}_X\,\psi_i^n\Big\|_X^2=\sum_{i=\ell+1}^{d^n}\lambda_i^n.
    \end{equation}
    Otherwise, if only the snapshots $\{y_j^h\}_{j=1}^n$ are used, we only have
    \begin{align*}
        \sum_{j=1}^n\frac{\alpha_j^n}{\delta t_j^2}\,\Big\|y_j^h-\sum_{i=1}^\ell{\langle y_j^h,\psi_i^h\rangle}_X\,\psi_i^h\Big\|_X^2=\sum_{i=\ell+1}^{d^n}\lambda_i^n.
    \end{align*}
    Consequently,
    \begin{align*}
        &\sum_{j=2}^n\alpha_j^n\Big\|\overline\partial y_j^h-\sum_{i=1}^\ell{\langle \overline\partial y_j^h,\psi_i^n\rangle}_X\,\psi_i^n\Big\|_X^2\\
        &\quad=\sum_{j=1}^n\frac{\alpha_j^n}{\delta t_j^2}\,\bigg(\Big\|y_j^h-\sum_{i=1}^\ell{\langle y_j^h,\psi_i^h\rangle}_X\,\psi_i^h\Big\|_X^2-\Big\|y_{j-1}^h-\sum_{i=1}^\ell{\langle y_{j-1}^h,\psi_i^h\rangle}_X\,\psi_i^h\Big\|_X^2\bigg)\le\frac{2}{\delta t^2}\sum_{i=\ell+1}^{d^n}\lambda_i^n
    \end{align*}
    which in contrast to \eqref{DQincl-1} contains the factor $\delta t^{-2}$ on the right-hand side.\hfill$\blacksquare$
\end{remark}

To compute the POD basis numerically, we follow Section~\ref{SIAM:Section-2.1.1.4} and set
\begin{align*}
    \bY^1&=\big[\mathrm y_1^h\,|\ldots|\,\mathrm y_n^h\big]\in\mathbb R^{m\times n},&\bY^2&=\big[0\,|\,\overline\partial\mathrm y_2^h\,|\ldots|\,\overline\partial\mathrm y_n^h\big]\in\mathbb R^{m\times n},\\
    \bD&=\mathrm{diag}\,(\alpha_1^n,\ldots,\alpha_n^n)\in\mathbb R^{n\times n},&\tilde{\bD}&=\mathrm{diag}\,(\bD,\bD)\in\mathbb R^{(2n)\times(2n)},\\
    \bY&=\big[\bY^1\,|\,\bY^2\big]\in\mathbb R^{m\times(2n)},&\tilde{\bY}&=\bW^{1/2} \bY\tilde{\bD}^{1/2}.
\end{align*}
Then a solution to \eqref{PODProbDisc} is given by $\uppsi_i^n=\bW^{-1/2}\tilde \uppsi_i^n\in\mathbb R^m$, $i=1,\ldots,\ell$, where the $\tilde\uppsi_i^n$'s is obtained computing the singular value decomposition of $\tilde{\bY}$; see Algorithm~\ref{SIAM:Algorithm-I.1.1.4}.

If $\{\uppsi_i^n\}_{i=1}^\ell$ solves \eqref{PODProbDisc} then $\{\psi_i^n\}_{i=1}^\ell$ is a solution to \eqref{PODProb} with
\begin{equation}
    \label{FEAnsatzPOD-2}
    \psi_i^n=\sum_{j=1}^m\big(\uppsi_i^n\big)_j\varphi_i^h\quad\text{for }i=1,\ldots,\ell.
\end{equation}
On the other hand, if $\{\psi_i^n\}_{i=1}^\ell$ is a solution to \eqref{PODProb} then we infer from $\psi_i^n\in V^h$ the existence of a uniquely determined matrix $\bPsi=[\uppsi_1^n\,|\ldots|\,\uppsi_\ell^n]\in\mathbb R^{m\times\ell}$ satisfying \eqref{FEAnsatzPOD-2}. It follows that the columns $\{\uppsi_i^n\}_{i=1}^\ell$ of $\bPsi$ solve \eqref{PODProbDisc}.

\subsection{POD Galerkin scheme}
\label{SIAM-Book:Section3.5.4}

Suppose that we have computed a POD basis $\{\psi_i^n\}_{i=1}^\ell$ of rank $\ell\in\{1,\ldots,d^n\}$. We set
\begin{align*}
    X^{n\ell}=\mathrm{span}\,\big\{\psi_1^n,\ldots,\psi_\ell^n\big\}\subset X\subset V
\end{align*}
and
\begin{align*}
    H^{n\ell}=\mathrm{span}\,\big\{\psi_1^{nH},\ldots,\psi_\ell^{nH}\big\},\quad V^{n\ell}=\mathrm{span}\,\big\{\psi_1^{nV},\ldots,\psi_\ell^{nV}\big\}
\end{align*}
for the choices $X=H$ and $X=V$, respectively. To simplify the presentation we utilize the same time grid $\{t_j\}_{j=1}^n$ in \eqref{EvProGal-disc}, \eqref{PODProbDisc} and the POD Galerkin scheme introduced below next. We replace \eqref{FineModel-Disc} (with $\theta=1$) by the following POD Galerkin scheme: Find $\{y^{h\ell}_j\}_{j=1}^n\subset X^{n\ell}$ satisfying
\begin{subequations}
    \label{EvProGalPOD-disc}
    \begin{equation}
        \label{EvProGalPOD-disc-1}
        \begin{aligned}
            &{\langle\overline\partial y^{h\ell}_j,\psi\rangle}_H+ a(t_j;y^{h\ell}_j,\psi)={\langle g^h_j(u),\psi\rangle}_{V',V}\text{ for all }\psi \in X^{n\ell}
        \end{aligned}
    \end{equation}
    for $2\le j\le n$ and
    \begin{equation}
        \label{EvProGalPOD-disc-2}
        y^{h\ell}_1=\mathcal P^{n\ell} y_\circ.
    \end{equation}
\end{subequations}
Recall the four different choices for $\mathcal P^{n\ell}$ introduced in \eqref{Eq:ProjOperatorsDisc}.

From $\{y_j^{h\ell}\}_{j=1}^n\subset X^{n\ell}$ we infer that there are uniquely determined coefficient vectors $\mathrm y_j^{h\ell}=(\mathrm y_{ji}^{h\ell})_{1\le i\le m}\in\mathbb R^\ell$, $1\le j\le n$, satisfying
\begin{equation}
    \label{POD-Gal-Ans}
    y_j^{h\ell}=\sum_{i=1}^\ell\mathrm y_{ji}^{h\ell}\psi_i^n\quad\text{for }j=1,\ldots,n.
\end{equation}
Inserting \eqref{POD-Gal-Ans} into \eqref{EvProGalPOD-disc}, choosing $\psi=\psi_i^n$, $1\le i\le \ell$, and using the notation introduced in Section~\ref{SIAM-Book:Section3.4.4} we obtain
\begin{subequations}
    \label{PODModel-Disc}
    \begin{align}
        \label{PODModel-Disc-1}
        \big(\bM^{h\ell}+\delta t_j\bA^{h\ell}_j\big)\mathrm y_j^{h\ell}&=\bM^{h\ell}\mathrm y^{h\ell}_{j-1}+\delta t_j\mathrm g_j^{h\ell}(u)\text{ for } 2 \le j \le n,\\
        \label{PODModel-Disc-2}
        \bM^{h\ell}\mathrm y^h_1&=\mathrm y_\circ^{h\ell},
    \end{align}
\end{subequations}
where $\bA^{h\ell}_j=\bA^{h\ell}(t_j)$ and $\mathrm g^{h\ell}_j(u)=\mathrm g_j^h(u)\in\mathbb R^\ell$ for $1\le j\le n$. 

In the following theorem, an existence result and a-priori estimates for the solution $\{y^{h\ell}_j\}_{j=1}^n$ are established. Since $\bM^{h\ell}+\delta t_j\bA^{h\ell}_j$ is regular, the proof follows by the same arguments as in the proof of Theorem~\ref{FullDiscFEModel}.

\begin{theorem}
    \label{PropExUnPODDisc11}
    Suppose that Assumption~{\rm\ref{A100}} and $\theta=1$ hold. The POD basis $\{\psi_i^n\}_{i=1}^\ell$ of rank $\ell\in\{1,\ldots,d^n\}$ is computed as described in Section~{\em\ref{SIAM-Book:Section3.5.3}}. Then there exists a unique sequence $\{y^{h\ell}_j\}_{j=1}^n \subset V^h$ to \eqref{EvProGalPOD-disc} provided $\Delta t$ is sufficiently small. If $\delta t$ is sufficiently small, the following estimates are satisfied:
    \begin{equation}
        \label{TempDiscPODSOL-APriori}
        {\|y_j^{h\ell}\|}_H^2\le e^{-C_1(j-1)\delta t}\bigg({\|\mathcal P^{n\ell}y_\circ\|}_H^2+\frac{1}{\gamma_1}\sum_{l=2}^j\delta t_l\,{\|g_l^h(u)\|}_{V'}^2\bigg)
    \end{equation}
    for $j=2,\ldots,n$ with the constant $C_1=2\gamma_1/c_V^2$ and
    \begin{equation}
        \label{TempDiscPODSOL-APriori-2}
        {\|y_n^{h\ell}\|}_H^2+\sum_{j=2}^n{\|y_j^{h\ell}-y_{j-1}^{h\ell}\|}_H^2+\gamma_1\sum_{j=1}^n\delta t_j\,{\|y_j^{h\ell}\|}_V^2\le{\|\mathcal P^{n\ell}y_\circ\|}_H^2+\frac{1}{\gamma_1}\sum_{j=2}^n\delta t_j\,{\|g_j^\ell(u)\|}_{V'}^2.
    \end{equation}
\end{theorem}

\begin{remark}
    \rm We can avoid the assumption that $\Delta t$ has to be sufficiently small. If $\gamma_2=0$ holds in \eqref{SIAM:Eq3.1.1-2} (cf. Assumption~\ref{A1}-2-3)), this follows arguing as in the proof of Theorem~\ref{PropExUnPODDisc11}. If $\gamma_2>0$ holds true, we can transform our evolution problem \eqref{SIAM:Eq3.1.6} as described in Remark~\ref{Remark:OperatorS}.\hfill$\blacksquare$
\end{remark}

\subsection{POD a-priori error analysis}
\label{SIAM-Book:Section3.5.5}

Utilizing Theorem~\ref{Theorem:FEAprioriError} we find that
\begin{align*}
    \sum_{j=1}^n\alpha_j^n\,{\|y(t_j)-y^{h\ell}_j\|}_V^2&\le2\sum_{j=1}^n\alpha_j^n\,{\|y(t_j)-y^h_j\|}_V^2+2\sum_{j=1}^n\alpha_j^n\,{\|y^h_j-y^{h\ell}_j\|}_V^2\\
    &\le C\big(h^2+\Delta t^2+\Delta t^2h^2\big)+2\sum_{j=1}^n\alpha_j^n\,{\|y^h_j-y^{h\ell}_j\|}_V^2.
\end{align*}
Thus, our goal is to derive an a-priori error estimate for the term
\begin{align*}
    \sum_{j=1}^n\alpha_j^n\,{\|y^h_j-y^{h\ell}_j\|}_V^2.
\end{align*}

For the a-priori error analysis we proceed as in Section~\ref{SIAM-Book:Section3.3.2} and utilize the rate of convergence results presented in Theorem~\ref{Prop:VTopologyDisc}. Thus, we make use of the following decomposition (compare \eqref{EqDec})
\begin{align*}
    y^h_j-y^{h\ell}_j&=y^h_j-\mathcal P^{n\ell}y^h_j+\mathcal P^{n\ell}y^h_j-y^{h\ell}_j=\varrho^{h\ell}_j+\vartheta^{h\ell}_j\quad\text{for }1\le j\le n
\end{align*}
with $\varrho^{h\ell}_j=y^h_j-\mathcal P^{n\ell}y^h_j\in (X^{n\ell})^\bot$ and $\vartheta^{h\ell}_j=\mathcal P^{n\ell}y^h_j-y^{h\ell}_j\in X^{n\ell}$. Applying the notation introduced in Section~\ref{SIAM-Book:Section3.5.4} we find
\begin{equation}
    \label{Chopin-5}
    \sum_{j=1}^n\alpha_j^n\left({\|\varrho_j^{h\ell}\|}_V^2+{\|\overline\partial\varrho_j^{h\ell}\|}_V^2\right)=\sum_{k=1}^2\sum_{j=1}^n\alpha_j^n\,{\|y_j^k-\mathcal P^{n\ell}y_j^k\|}_V^2.
\end{equation}
To achieve that
\begin{equation}
    \label{Eq:FDPOD-0}
    \vartheta^{h\ell}_1=\mathcal P^{n\ell}y^h_1-y^{h\ell}_1=\mathcal P^{n\ell}\big(\mathcal P^hy_\circ\big)-\mathcal P^{n\ell}y_\circ=0
\end{equation}
we have to choose appropriate $\mathcal P^h$ and $\mathcal P^{n\ell}$ in \eqref{EvProGal-2} and \eqref{EvProGal-POD-2}, respectively. Recall the choices made in Remark~\ref{Rem:ProjChoice} in order to ensure \eqref{Eq:FDPOD-0}:
\begin{enumerate}
    \item [1)] $X=H$, $\mathcal P^{n\ell}=\mathcal P^{n\ell}_H$ and $\mathcal P^h=\mathcal P^h_H$,
    \item [2)] $X=H$, $\mathcal P^{n\ell}=\mathcal P^{n\ell}_V$ and $\mathcal P^h=\mathcal P^h_V$,
    \item [3)] $X=V$, $\mathcal P^{n\ell}=\mathcal Q^{n\ell}_H$ and $\mathcal P^h=\mathcal P^h_H$,
    \item [4)] $X=V$, $\mathcal P^{n\ell}=\mathcal Q^{n\ell}_V$ and $\mathcal P^h=\mathcal P^h_V$.
\end{enumerate}
The next lemma is shown in Section~\ref{SIAM-Book:Section3.8.4}.

\begin{lemma}
    \label{Lemma:FullDiscTheta}
    Let Assumption~{\rm\ref{A100}} and $\theta=1$ hold. Suppose that $\{y^h_j\}_{j=1}^n \subset V^h$ solves \eqref{EvProGal-disc}. The POD basis $\{\psi_i^n\}_{i=1}^\ell$ of rank $\ell\in\{1,\ldots,d^n\}$ is computed as described in Section~{\em\ref{SIAM-Book:Section3.5.3}} and $\{y^{h\ell}_j\}_{j=1}^n \subset V^h$ is the solution to \eqref{EvProGalPOD-disc}. Then $\vartheta^{h\ell}_j=\mathcal P^{h\ell}y^h_j-y^{h\ell}_j$, $1\le j\le n$ satisfies
    \begin{equation}
        \label{Chopin-6}
        \sum_{j=1}^n\alpha_j^n\,{\|\vartheta_j^{h\ell}\|}_V^2\le C\sum_{k=1}^2\sum_{j=1}^n\alpha_j^n\,{\|y_j^k-\mathcal P^{n\ell}y_j^k\|}_V^2
    \end{equation}
    with a constant depending on $\gamma$, $\gamma_1$, $c_V$, $T$ and $\zeta$.
\end{lemma}

Now, the following a-priori error follows directly from Proposition~\ref{Prop:VTopologyDisc}, \eqref{Chopin-5}, Lemma~\ref{Lemma:FullDiscTheta} and from
\begin{align*}
    \sum_{j=1}^n\alpha_j^n\,{\|y^h_j-y^{h\ell}_j\|}_V^2&\le2\sum_{j=1}^n\alpha_j^n\,{\|\vartheta^{h\ell}_j\|}_V^2+2\sum_{j=1}^n\alpha_j^n\,{\|\varrho^{h\ell}_j\|}_V^2\le C_\circ\sum_{k=1}^2\sum_{j=1}^n\alpha_j^n\,{\|y_j^k-\mathcal P^{n\ell}y_j^k\|}_V^2
\end{align*}
with $C_\circ=2(C+1)$.

\begin{theorem}
    \label{Th:A-PrioriError-200}
    Let Assumption~{\rm\ref{A100}} and $\theta=1$ hold. Suppose that $\{y^h_j\}_{j=1}^n \subset V^h$ solves \eqref{EvProGal-disc}. The POD basis $\{\psi_i^n\}_{i=1}^\ell$ of rank $\ell\in\{1,\ldots,d^n\}$ is computed as described in Section~{\em\ref{SIAM-Book:Section3.5.3}} and $\{y^{h\ell}_j\}_{j=1}^n \subset V^h$ is the solution to \eqref{EvProGalPOD-disc}. Then we have the following \index{Error estimate!a-priori!state variable}a-priori error estimates:
    \begin{enumerate}
        \item[\rm 1)] For $X=H$, $\mathcal P^{n\ell}=\mathcal P^{n\ell}_H$ and $\mathcal P^h=\mathcal P^h_H$ we have
        \begin{align*}
            \sum_{j=1}^n\alpha_j^n\,{\|y^h_j-y^{h\ell}_j\|}_V^2\le C\sum\limits_{i=\ell+1}^{d^n}\lambda_i^{nH}\,{\|\psi_i^{nH}\|}_V^2.
        \end{align*}
        \item[\rm 2)] For $X=H$, $\mathcal P^{n\ell}=\mathcal P^{n\ell}_V$ and $\mathcal P^h=\mathcal P^h_V$ it follows that 
        \begin{align*}
            \sum_{j=1}^n\alpha_j^n\,{\|y^h_j-y^{h\ell}_j\|}_V^2\le C\sum\limits_{i=\ell+1}^{d^n}\lambda_i^{nH}\,\big\|\psi_i^{nH}-\mathcal P^{n\ell}_V\psi_i^{nH}\big\|_V^2.
        \end{align*}
        \item[\rm 3)] For $X=V$, $\mathcal P^{n\ell}=\mathcal Q^{n\ell}_H$ and $\mathcal P^h=\mathcal P^h_H$ we derive
        \begin{align*}
            \sum_{j=1}^n\alpha_j^n\,{\|y^h_j-y^{h\ell}_j\|}_V^2\le C\sum\limits_{i=\ell+1}^{d^n}\lambda_i^{nV}\,\big\|\psi_i^{nV}-\mathcal Q^{n\ell}\psi_i^{nV}\big\|_V^2.
        \end{align*}
        \item[\rm 4)] For $X=V$, $\mathcal P^{n\ell}=\mathcal Q^{n\ell}_V$ and $\mathcal P^h=\mathcal P^h_V$ and  we get
        \begin{align*}
            \sum_{j=1}^n\alpha_j^n\,{\|y^h_j-y^{h\ell}_j\|}_V^2\le C\sum\limits_{i=\ell+1}^{d^n}\lambda_i^{nV}.
        \end{align*}
    \end{enumerate}
\end{theorem}

\begin{example}
    \label{Example:Ch3_aprioriestimates_example}
    \rm We want to verify the a-priori estimates in Theorem~\ref{Th:A-PrioriError-200} for the guiding model solved with a fix control $u$, which is generated randomly. To make more relevant the use of different projections, we consider $\kappa=0.02$, a time-dependent velocity field $\bv(t,\bx)= ((0.1+0.1t)(x_2-0.5),(0.1+0.1t)(0.5-x_1))$ and an initial guess $y_\circ(\bx)= x_1+x_2$. The other parameters are set as in Section~\ref{SIAM-Book:Section1.3}.
    \begin{figure}
        \begin{center}
            \includegraphics[height=50mm]{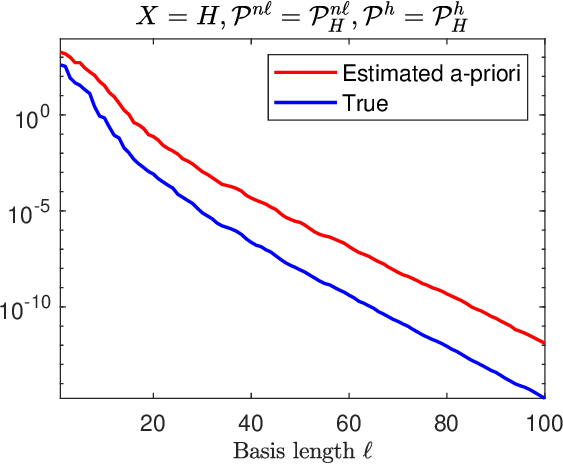}\hspace{10mm}
            \includegraphics[height=50mm]{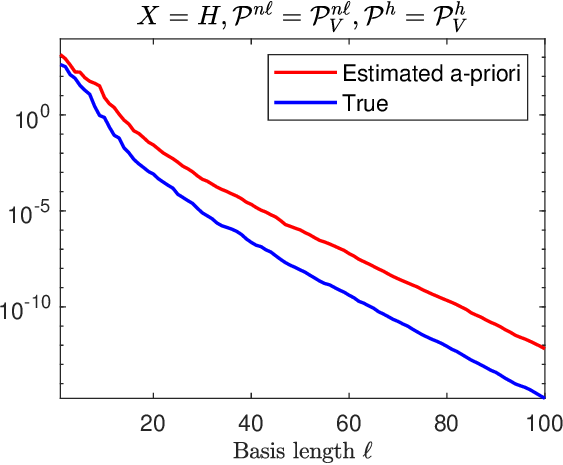} \\[4mm]
            \includegraphics[height=50mm]{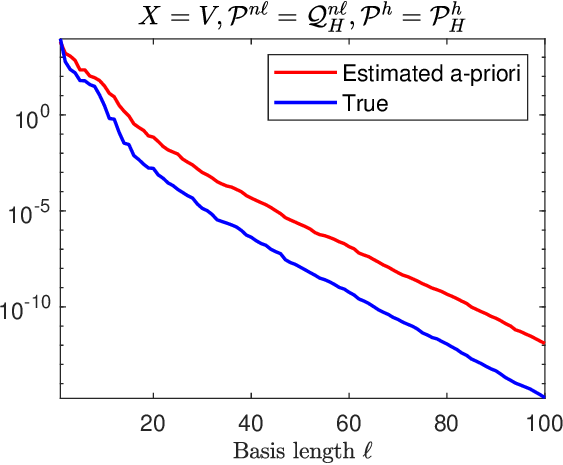}\hspace{10mm}
            \includegraphics[height=50mm]{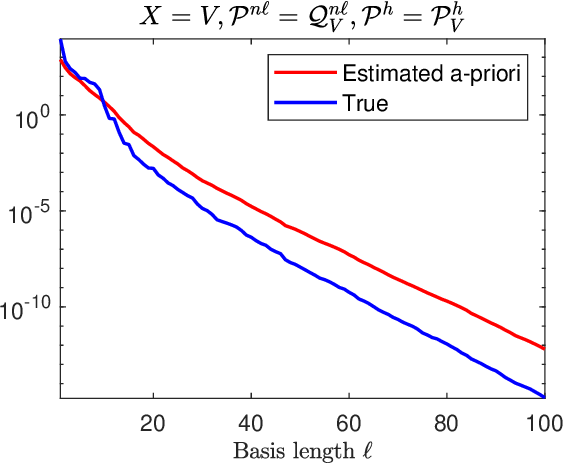}
        \end{center}
        \caption{Example~\ref{Example:Ch3_aprioriestimates_example}. Plot of the a-priori estimate and the approximation error of the reduced order model for different choices of spaces and projections (y-axis: logarithmic scale).}
        \label{fig:Ch3_aprioriestimates_example}
    \end{figure}
    As one can see from Figure~\ref{fig:Ch3_aprioriestimates_example}, the a-priori error estimate predicts the behaviour of the error as soon as $\ell$ increases. Note that the constant $C$ is not always computable and it depends on the choice of the spaces $H$ and $V$. Therefore, we report only the quantities on the right-hand side of the estimates in Theorem~\ref{Th:A-PrioriError-200}, without the multiplication of the constant $C$. This is the reason of the underestimation for small $\ell$ in case 4). Moreover, this justify the use of an a-priori estimate for predicting the qualitative behaviour of the reduced order model rather than the quantitative one and, thus, to choose the dimension $\ell$ of the reduced space.\hfill$\blacklozenge$
\end{example}

\subsubsection{Extensions}

To derive the rate of convergence results presented in Theorem~\ref{Th:A-PrioriError-200} we have utilized the snapshots $\{y^h_j\}_{j=1}^n$ and $\{\overline\partial y^h_j\}_{j=2}^n$ for the POD basis computation; cf. Section~\ref{SIAM-Book:Section3.5.3}. Next we show that we can omit the difference quotients $\{\overline\partial y^h_j\}_{j=2}^n$, but we still can derive a rate of convergence result. For that purpose we consider the $H$-orthonormal projection operators $\mathcal P^{n\ell}=\mathcal P^{n\ell}_H$ and $\mathcal P^{n\ell}=\mathcal Q^{n\ell}_H$; cf. Theorem~\ref{Th:A-PrioriError-2} and \eqref{Chopin-2}. Moreover we choose $\mathcal P^h=\mathcal P^h_H$ in \eqref{EvProGal-disc-2}. Instead of \eqref{PODProb} we compute the POD basis of rank $\ell$ by solving
\begin{align*}
    \min\sum_{j=1}^n\alpha_j^n\,\Big\|y_j^h-\sum_{i=1}^\ell{\langle y_j^h,\psi_i^n\rangle}_X\,\psi_i^n\Big\|_X^2\quad\text{s.t.}\quad\{\psi_i^n\}_{i=1}^\ell\subset V^h\text{ and }{\langle\psi_i^n,\psi_j^n\rangle}_X=\delta_{ij},~1\le i,j\le \ell.
\end{align*}

In Section~\ref{SIAM-Book:Section3.8.4} we prove the following result.

\begin{theorem}
    \label{Theorem:NoDQSnaps}
    Let Assumption~{\rm\ref{A100}} and $\theta=1$ hold. Suppose that $\{y^h_j\}_{j=1}^n \subset V^h$ solves \eqref{EvProGal-disc} with $\mathcal P^h=\mathcal P^h_H$. The POD basis $\{\psi_i^n\}_{i=1}^\ell$ of rank $\ell\in\{1,\ldots,d^n\}$ is computed as described in Section~{\em\ref{SIAM-Book:Section3.5.3}} and $\{y^{h\ell}_j\}_{j=1}^n \subset V^h$ is the solution to \eqref{EvProGalPOD-disc}. Then the \index{Error estimate!a-priori!state variable}a-priori error estimate
    \begin{align*}
        \int_0^T{\|y^h(t)-y^{h\ell}(t)\|}_V^2\,\mathrm \le~C\cdot\left\{
        \begin{aligned}
            &\sum_{i=\ell+1}^{d^n}\lambda_i^{nH}\,\big\|\psi_i^{nH}\big\|_V^2,&&X=H,~\mathcal P^{n\ell}=\mathcal P^{n\ell}_H,\\
            &\sum_{i=\ell+1}^{d^n}\lambda_i^{nV}\,\big\|\psi_i^{nV}-\mathcal Q^{n\ell}_H\psi_i^{nV}\big\|_V^2,&&X=V,~\mathcal P^{n\ell}=\mathcal Q^{n\ell}_H
        \end{aligned}
        \right.
    \end{align*}
    holds with a constant $C>0$.
\end{theorem}

\subsection{POD a-posteriori error analysis}
\label{SIAM-Book:Section3.5.6}

Suppose that $\{y^h_j\}_{j=1}^n\subset V^h$ and $\{y^{h\ell}_j\}_{j=1}^n\subset X^{h\ell}$ are the solutions to \eqref{EvProGal-disc} and \eqref{EvProGalPOD-disc}, respectively. In Theorems~\ref{Th:A-PrioriError-200} and \ref{Theorem:NoDQSnaps} we have presented a-priori error estimates for the difference $y_h^j-y_j^{h\ell}$, $j=1,\ldots,n$. However, the convergence rates are only valid if the POD basis is computed utilizing the solution $\{y^h_j\}_{j=1}^n$ and -- in case of Theorem~\ref{Th:A-PrioriError-200} -- also the difference quotients $\{\overline\partial y^h_j\}_{j=2}^n$ as snapshots. If the POD basis is computed from different snapshots, the a-priori error estimates does not hold. In this case we need an a-posteriori error estimate which indicates whether the computed POD solution $\{y^{h\ell}_j\}_{j=1}^n$ is sufficiently accurate or not.

In the next theorem we present an a-posteriori error analysis for the fully discrete Galerkin scheme \eqref{EvProGalPOD-disc}. The proof is based on the arguments used in Section~\ref{SIAM-Book:Section3.4.6}; cf. Section~\ref{SIAM-Book:Section3.8.4}.

\begin{theorem}
    \label{Prop:ApostiError}
    Let us assume that  Assumption~{\rm\ref{A100}} and $\theta=1$ hold. Suppose that $\{y^h_j\}_{j=1}^n \subset V^h$ solves \eqref{EvProGal-disc}. The POD basis $\{\psi_i^n\}_{i=1}^\ell$ of rank $\ell\in\{1,\ldots,d^n\}$ is computed as described in Section~{\em\ref{SIAM-Book:Section3.5.3}} and $\{y^{h\ell}_j\}_{j=1}^n \subset V^h$ is the POD Galerkin solution to \eqref{EvProGalPOD-disc}. If $\delta t$ is sufficiently small, the \index{Error estimate!a-posteriori!state variable}a-posteriori error estimates
    \begin{align*}
        {\|y_j^h-y_j^{h\ell}\|}_H^2&\le \frac{e^{-2\gamma_1(j-1)\delta t/c_V^2}}{\gamma_1}\,\Delta t\sum_{l=2}^{j}\,{\|r_l^{h\ell}\|}_{(V^h)'}^2\quad\text{for }2\le j\le n,\\
        \sum_{j=1}^n\alpha_j^n\,{\|y_j^h-y_j^{h\ell}\|}_V^2&\le \frac{\zeta}{\gamma_1^2}\sum_{j=2}^n\delta t_j\,{\|r_j^{h\ell}\|}_{(V^h)'}^2
    \end{align*}
    hold, where the residual vectors are given for $j=2,\ldots,n$ by 
    \begin{equation}
        \label{Shostakovich-1}
        r_j^{h\ell}={\langle\overline\partial y_j^{h\ell},\cdot\rangle}_H+a(t_j;y^{h\ell}_j,\cdot)-{\langle (\mathcal F+\mathcal Bu)(t_j),\cdot\rangle}_{V',V}\in (V^h)'.
    \end{equation}
\end{theorem}

\begin{example}
    \label{Example:Ch3_aposterioriestimates_example}
    \rm Similarly to Example~\ref{Example:Ch3_aprioriestimates_example}, we want to verify the a-posteriori error estimate in Theorem~\ref{Prop:ApostiError} for the guiding example. We set all the parameters of the problem as in Example~\ref{Example:Ch3_aprioriestimates_example}. The only difference is that we generate the snapshots with a random control $u_1$ and compute the solution of the reduced order model as well as the a-posteriori error estimate for a second random control $u_2$.
    \begin{figure}
        \begin{center}
            \includegraphics[height=50mm]{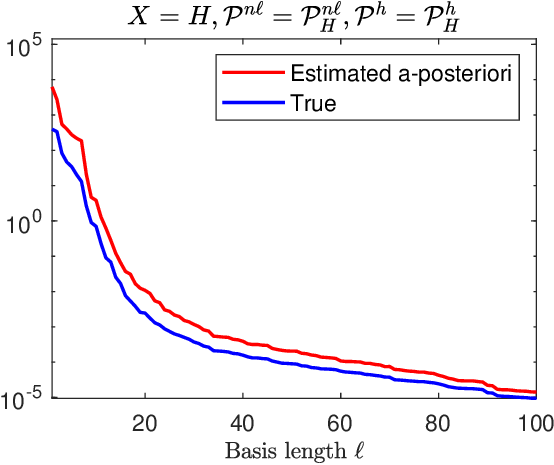}\hspace{10mm}
            \includegraphics[height=50mm]{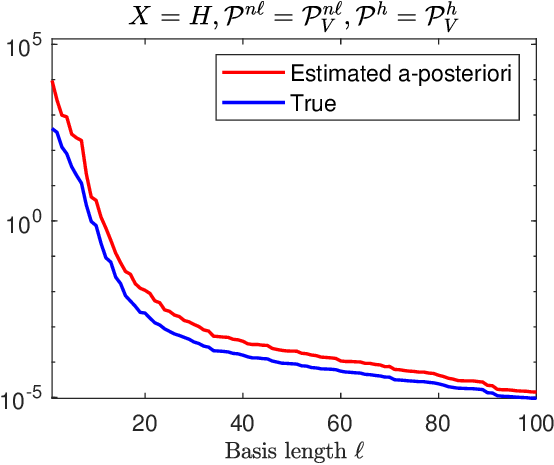}\\[4mm]
            \includegraphics[height=50mm]{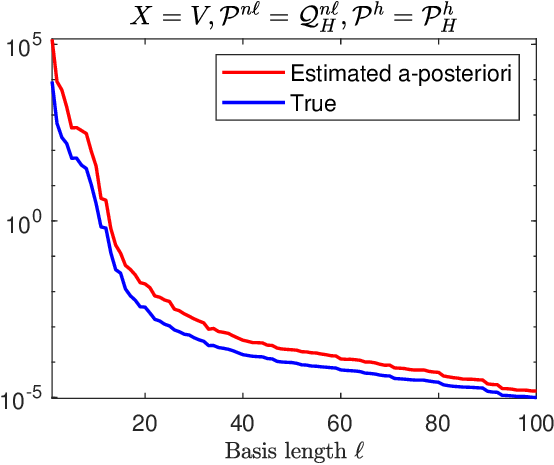}\hspace{10mm}
            \includegraphics[height=50mm]{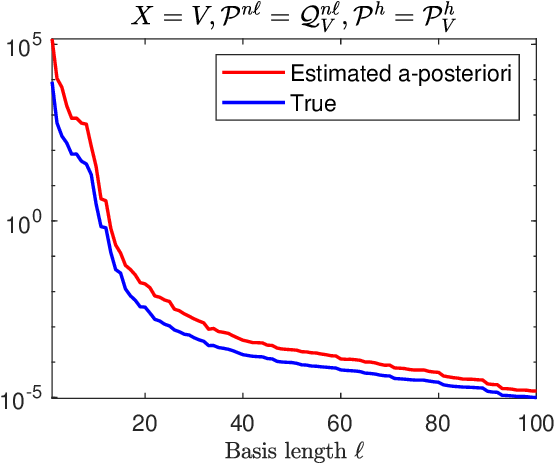}
        \end{center}
        \caption{Example~\ref{Example:Ch3_aposterioriestimates_example}. Plot of the a-posteriori estimate and the approximation error of the reduced order model for different choices of spaces and projections (y-axis: logarithmic scale).}
        \label{fig:Ch3_aposterioriestimates_example}
    \end{figure}
    According to Figure~\ref{fig:Ch3_aposterioriestimates_example} and as expected, there is a bad approximation of the full order solution of the problem for small $\ell$. This is mainly caused by the advection $\bv$ that generates several different dynamics and the fact that the two random controls $u_1$ and $u_2$ may be far from each other. For $\ell= 100$ in fact the true approximation error is still around $10^{-5}$. Another important issue is the factor of overestimation for the a-posteriori estimate, which is influenced by the quantities $\gamma_1$ ($=10^{-4}$ for this example), $\delta t_j$ ($=0.025$), $\zeta$ ($=1$) and indirectly by the velocity field $\bv$ through the residual $r_j^{h\ell}$; cf. Theorem~\ref{Prop:ApostiError}. In Figure~\ref{fig:Ch3_aposterioriestimates_example} one can note that the factor of overestimation is between $10$ and $30$ for each choice of projections $\mathcal{P}^h$ and $\mathcal{P}^{n\ell}$ and space $X$, thus still acceptable. To have a good a-posteriori error estimate is in fact crucial for a POD-based model order reduction, since it indicates whether the true approximation is expected to be good or not. For further details on this topic we refer to Chapter~\ref{SIAM-Book:Section4} and Chapter~\ref{Advanced Topics in POD Suboptimal control}. Note at last that there is no possibility of underestimation, as it is happening instead in Example~\ref{Example:Ch3_aprioriestimates_example}. This is clearly because in this example all the quantities for the a-posteriori error estimation are computable. This is another crucial point for judging the good or bad quality of an a-posteriori estimate.\hfill$\blacklozenge$
\end{example}

\section{Extension to non-linear evolution problems}
\label{SIAM-Book:Section3.6}
\setcounter{equation}{0}
\setcounter{theorem}{0}
\setcounter{figure}{0}
\setcounter{run}{0}

The application of the POD method can be extended to non-linear evolution problems. Here we focus on a non-linear problem which arises in fluid dynamics and closely follow Section III-3 in \cite{Tem88}. 

\subsection{The abstract non-linear evolution problem}
\label{SIAM-Book:Section3.6.1}

Let $V$ and $H$ be real separable Hilbert spaces and suppose that $V$ is dense in $H$ with compact embedding. As in Remark~\ref{Remark:HI-20}-1, we consider a \index{Bilinear form!time-independent, $a(\cdot\,,\cdot)$}time-independent bilinear form $a:V\times V \to \mathbb{R}$ and associate with it the linear operator $\mathcal A:V \to V'$ by
\begin{align*}
    {\langle \mathcal A\varphi,\phi \rangle}_{V',V}=a(\varphi,\phi)\quad\text{for } \varphi,\phi \in V.
\end{align*}
To be able to consider $\mathcal A$ as an operator mapping from $H$ to $H$, its domain is defined as in \eqref{DomainA}. In contrast to the linear case, we additionally introduce a linear and continuous operator $\mathcal C:V \to V'$ which maps $D(\mathcal A)$ into $H$. Moreover, let $\mathcal N: V \to V'$ be a continuous non-linear function, $y_\circ \in V$ an initial condition, $\mathcal F\in L^2(0,T;H)$ an inhomogeneity and $u$ a control element of the Hilbert space $\U$ which is mapped by a linear and bounded \index{Linear operator!contrrol, $\mathcal B$}{\em control operator} $\mathcal B: \U \to L^2(0,T;H)$ to a function $\mathcal B u \in L^2(0,T;H)$. Then we consider the non-linear evolution problem 
\begin{equation}
    \label{NonlEvPro}
    \begin{aligned}
    \frac{\mathrm d}{\mathrm dt} \, {\langle y(t),\varphi \rangle}_H + a(y(t),\varphi)+{\langle \mathcal C y(t) + \mathcal N(y(t)),\varphi\rangle}_{V',V}&={\langle(\mathcal F+\mathcal Bu)(t),\varphi\rangle}_{V',V}\\
    &\hspace{22.5mm}\text{for all }\varphi\in V\text{ a.e. in }(0,T],\\
{\langle y(0),\varphi \rangle}_H&={\langle y_\circ,\varphi\rangle}_H\quad\text{for all }\varphi \in V.
\end{aligned}
\end{equation}
Results about the existence of solutions to \eqref{NonlEvPro} depend on the concrete qualities of $a$, $\mathcal C$ and $\mathcal N$. The following hypotheses are taken from \cite{Tem88}.

\begin{assumption}
	\label{Asspt-1}
	\begin{enumerate}
        \item [\em 1)] Let $y_\circ\in V$, $\mathcal F\in L^2(0,T;H)$ and $u\in\U$ with $\mathcal Bu\in L^2(0,T;H)$.
	    \item [\em 2)] We assume that the bilinear form $a$ is bounded, i.e. there is a constant $\gamma \ge 0$ such that
        \begin{equation}
	       \label{Eq:Jan-4}
		   \big|a(\varphi,\psi)\big|\le\gamma\,{\|\varphi\|}_V{\|\psi\|}_V\quad\text{for all }\varphi,\psi\in V.
	    \end{equation}
	    Additionally, we also assume that $\mathcal A+\mathcal C$ is coercive
        \begin{equation}
	        \label{OperA&C}
		    a(\varphi,\varphi)+{\langle \mathcal C\varphi,\varphi \rangle}_{V',V} \ge\gamma_1 \, {\| \varphi \|}_V^2\quad\text{for all } \varphi \in V
	    \end{equation}
	    with a constant $\gamma_1>0$.
	    \item [\em 3)] Let further $\mathcal C$ satisfy the properties 
        \begin{equation}
	        \label{OperC}
	        \begin{aligned}
		        {\| \mathcal C \varphi \|}_H & \le c_\mathcal C \, {\| \varphi \|}_V^{1-\delta_1}{\| \mathcal A\varphi \|}_H^{\delta_1} && \text{for all } \varphi \in D(\mathcal A),\\
		        \big| {\langle \mathcal C \varphi ,\varphi \rangle}_{V',V} \big| & \le c_\mathcal C\,{\| \varphi \|}_V^{1+\delta_2} {\| \varphi \|}_H^{1-\delta_2} &&\text{for all } \varphi \in V,
	        \end{aligned}
	    \end{equation}
	    for a constant $c_\mathcal C>0$ and for $\delta_1,\delta_2 \in [0,1)$.
	    \item [\em 4)] For the non-linearity, suppose there is another bilinear and continuous function $\mathcal M: V \times V \to V'$ satisfying $\mathcal M(\varphi,\varphi) = \mathcal N(\varphi)$ for all $\varphi \in V$. Let it map from $V \times D(\mathcal A)$ respectively $D(\mathcal A)\times V$ into $H$ such that there exist constants $c_{\mathcal M}>0$ and $\delta_3,\delta_4,\delta_5 \in [0,1)$ with
        \begin{equation}
		    \label{OperN}
		    \begin{aligned}
		        {\langle \mathcal M(\varphi,\phi),\phi \rangle}_{V',V}&=0,\\
		        \big| {\langle \mathcal M(\varphi,\phi),\psi \rangle}_{V',V} \big|&\le c_\mathcal M\,{\| \varphi \|}_H^{\delta_3}{\| \varphi \|}_V^{1-\delta_3}{\| \phi \|}_V{\| \psi \|}_V^{\delta_3}{\| \psi \|}_H^{1-\delta_3},\\
		        {\| \mathcal M(\varphi,\chi) \|}_H+{\|\mathcal M(\chi,\varphi) \|}_H&\le c_\mathcal M\,{\| \varphi \|}_V{\| \chi \|}_V^{1-\delta_4}{\| \mathcal A \chi\|}_H^{\delta_4},\\
		        {\| \mathcal M(\varphi,\chi) \|}_H&\le c_\mathcal M\,{\|\varphi\|}_H^{\delta_5}{\|\varphi\|}_V^{1-\delta_5}{\|\chi\|}_V^{1-\delta_5}{\|\mathcal A\chi\|}_H^{\delta_5}
            \end{aligned}
	    \end{equation}
        for all $\varphi,\phi,\psi \in V$ and for all $\chi \in D(\mathcal A)$.
	\end{enumerate}
\end{assumption}

The following result extends Theorem~2.1 in \cite[p.~111]{Tem88} in the sense that the function $f$ was not originally assumed to be time-dependent. However, the proof is analogous in every step for our case.
\begin{theorem}
    \label{theo2-1}
    Assume that Assumption~{\em\ref{Asspt-1}} holds. Then there exists a unique solution of \eqref{NonlEvPro} satisfying
    \begin{equation}
        \label{eq2-9}
        y \in C([0,T];V) \cap L^2(0,T;D(\mathcal A)).
    \end{equation}
\end{theorem}

Let us present an example for the non-linear evolution system \eqref{NonlEvPro} under Assumption \ref{Asspt-1}.

\begin{example}
    \rm Let $\Omega$ denote a bounded domain in $\mathbb R^2$ with boundary $\Gamma$ and let $T>0$. The two-dimensional \index{Equation!Navier-Stokes}{\em Navier-Stokes equation} is given by
    \begin{subequations}
        \label{sube1}
        \begin{align}
            \label{sube1a}
            \rho~\big(y_t+(y\cdot\nabla)y\big)-\nu\Delta y+\nabla p&=f&&\text{in } Q=(0,T) \times \Omega,\\
            \label{sube1b}
            \mathrm{div}~y &=0&&\text{in } Q,
        \end{align}
        where $\rho>0$ is the density of the fluid, $\nu>0$ is the kinematic viscosity, $f$ represents volume forces and
        \begin{align*}
            (y\cdot\nabla)y=\Big(y_1\frac{\partial y_1}{\partial x_1}+y_2\frac{\partial y_1}{\partial x_2},y_1\frac{\partial y_2}{\partial x_1}+y_2\frac{\partial y_2}{\partial x_2}\Big)^\top.
        \end{align*}
        The unknowns are the velocity field $y=(y_1,y_2)$ and the pressure $p$. Together with \eqref{sube1}, one can consider non-slip boundary conditions
        \begin{equation}
            y=y_d \quad \text{on } \Sigma=(0,T) \times \Gamma
        \end{equation}
        and the initial condition
        \begin{equation}
            y(0,\cdot)=y_\circ\quad\text{in }\Omega.
        \end{equation}
    \end{subequations}
    In \cite[pp.~104-107, 116-117]{Tem88} it was shown that \eqref{sube1} can be written in the form \eqref{NonlEvPro}.\hfill$\blacklozenge$
\end{example}

\subsection{The continuous POD method}
\label{SIAM-Book:Section3.6.2}

In this subsection we introduce the continuous version of the POD method for \eqref{NonlEvPro}.

\subsubsection{POD Galerkin scheme}

Let $y\in L^2(0,T;V)$ be a solution to \eqref{NonlEvPro}. We choose the snapshot set
\begin{align*}
    \mathscr V=\bigg\{\int_0^T\varphi(t)y(t)\,\mathrm dt\,\Big|\,\varphi\in L^2(0,T)\bigg\}\subset V
\end{align*}
with dimension
\begin{align*}
    d=\left\{
    \begin{aligned}
        &\dim\mathscr V&&\text{if }\dim\mathscr V<\infty,\\
        &\infty&&\text{otherwise.}
    \end{aligned}
    \right.
\end{align*}
Note that we do not take the time derivatives as snapshots. This is often done in numerical examples. However, the inclusion of the time derivatives is straightforward; cf. Section~\ref{Section:ContPODHilbert}. Again, $X$ stands either for the separable Hilbert spaces $H$ or $V$. Then the POD basis of rank $\ell$ is given by
\begin{equation}
    \label{Eq:Jan-1}
    \left\{
    \begin{aligned}
        &\min\int_0^T\Big\|y(t)-\sum_{i=1}^\ell{\langle y(t),\psi_i\rangle}_X\,\psi_i\Big\|_X^2\,\mathrm dt\\
        &\hspace{0.4mm}\text{s.t. }\{\psi_i\}_{i=1}^\ell\subset X\text{ and }{\langle\psi_i,\psi_j\rangle}_X=\delta_{ij}\text{ for }1\le i,j\le\ell.
    \end{aligned}
    \right.
\end{equation}
A solution to \eqref{Eq:Jan-1} is given by Theorem~\ref{Theorem2.2.1}. Proceeding as in Section~\ref{Section:ContPODHilbert} we obtain the bases $\{\psi_i^H\}_{i\in\mathbb I}$ for $X=H$ and $\{\psi_i^V\}_{i\in\mathbb I}$ for $X=V$ with the associated eigenvalues $\{\lambda_i^H\}_{i\in\mathbb I}$ and $\{\lambda_i^V\}_{i\in\mathbb I}$, respectively. We set
\begin{align*}
    H^\ell=\mathrm{span}\{\psi_1^H,\ldots,\psi_\ell^H\}\quad\text{and}\quad V^\ell=\mathrm{span}\{\psi_1^V,\ldots,\psi_\ell^V\}.
\end{align*}
If we do not want to distinguish between the choices $X=H$ and $X=V$, we just write $X^\ell=\mathrm{span}\{\psi_1,\ldots,\psi_\ell\}$. Recall that $H^\ell\subset V$ follows from Lemma~\ref{SIAM:Lemma3.2.1}-1).

The POD Galerkin method for \eqref{NonlEvPro} reads as follows: find $y^\ell: [0,T] \to X^\ell$ satisfying
\begin{equation}
    \label{Eq:Jan-2}
    \begin{aligned}
        \frac{\mathrm d}{\mathrm dt} \, {\langle y^\ell(t),\psi \rangle}_H + a(y^\ell(t),\psi)+{\langle\mathcal Cy^\ell(t)+\mathcal N(y^\ell(t)),\psi\rangle}_{V',V}&={\langle(\mathcal F+\mathcal Bu)(t),\psi\rangle}_{V',V}\\
        &\quad\text{for all }\psi\in X^\ell\text{ a.e. in }(0,T],\\
        \hspace{17.5mm}y^\ell(0)&=\mathcal P^\ell y_\circ,
    \end{aligned}
\end{equation}
where $\mathcal P^\ell:X\to X^\ell$ denotes either the $H$-orthogonal projection $\mathcal P_H^\ell$ (in case of $X=H$) or the $H$-orthogonal projection $\mathcal Q_H^\ell$ (in case of $X=V$).

\begin{assumption}
    \label{Asspt-1-POD}
    \begin{enumerate}
        \item [\rm 1)] Assumption~{\rm\ref{Asspt-1}} holds.
        \item [\rm 2)] The Hilbert space $X$ denotes either $H$ or $V$.
    \end{enumerate}
\end{assumption}

Recall also that we have set $\Y=W(0,T)$. The next result is proved in Section~\ref{SIAM-Book:Section3.8.5}. 

\begin{proposition}
    \label{Prop:NavierStokes-1}
    Suppose that Assumption~{\rm\ref{Asspt-1-POD}} holds. Then the solution to \eqref{Eq:Jan-2} is unique and satisfies
    \begin{equation}
        \label{Shostakovich-3}
        {\|y^\ell\|}_\Y\le C\Big({\|\mathcal P^\ell y_\circ\|}_H+{\|f\|}_{L^2(0,T;H)}+{\|u\|}_U\Big)
    \end{equation}
    for a constant $C>0$ which is independent of $\ell$.
\end{proposition}

From $y^\ell(t)\in V^\ell$ for almost all $t\in[0,T]$, it follows that
\begin{align*}
    y^\ell(t)=\sum_{i=1}^\ell\mathrm y_i^\ell(t)\psi_i\quad\text{f.a.a. }t\in[0,T],
\end{align*}
with a coefficient vector $\mathrm y^\ell(t)=(\mathrm y_i^\ell(t))_{1\le i\le \ell}\in\mathbb R^\ell$. In Section~\ref{SIAM-Book:Section3.3.1} we have introduced the mass matrix $\bM^\ell\in\mathbb R^{\ell\times\ell}$, the vector $\mathrm g^\ell(t;u)\in\mathbb R^\ell$ for the right-hand side and the vector $\mathrm y_\circ^\ell\in\mathbb R^\ell$ for the initial condition. Moreover, we define the matrix
\begin{align*}
    \bA^\ell=\big(\big(a(\psi_j,\psi_i)+{\langle\mathcal C\psi_j,\psi_i\rangle}_{V',V}\big)\big)\in\mathbb R^{\ell\times\ell}.
\end{align*}
For the non-linearity we introduce the vector
\begin{align*}
    \mathrm n^\ell(\mathrm v)=\bigg(\Big\langle\mathcal N\Big(\sum_{j=1}^\ell\mathrm v_j\psi_j\Big),\psi_i\Big\rangle_{V',V}\bigg)_{1\le i\le\ell}
\end{align*}
for the vector $\mathrm v=(\mathrm v_i)_{1\le i\le\ell}$. Then \eqref{Eq:Jan-2} is equivalent to the following non-linear, $\ell$-dimensional system of differential equations
\begin{equation}
    \label{Eq:Jan-3}
    \begin{aligned}
        \bM^\ell\dot{\mathrm y}^\ell(t)+\bA^\ell\mathrm y^\ell(t)+\mathrm n^\ell\big(\mathrm y^\ell(t)\big)&=\mathrm g^\ell(t;u)\quad\text{for }t\in(0,T],\\
        \mathrm y^\ell(0)&=\mathrm y_\circ^\ell.
    \end{aligned}
\end{equation}
System \eqref{Eq:Jan-3} is a low-dimensional model for \eqref{NonlEvPro}.

\subsubsection{POD a-priori error analysis}

In the next theorem we present an a-priori error estimate for the term
\begin{align*}
    \int_0^T{\|y(t)-y^\ell(t)\|}_V^2\,\mathrm dt,
\end{align*}
where $y$ and $y^\ell$ solve \eqref{NonlEvPro} and \eqref{Eq:Jan-2}, respectively. The proof is based on the results of Proposition~\ref{Prop:Rate-V}, Section~\ref{SIAM-Book:Section3.3.2}, \eqref{Eq:Jan-4}, \eqref{OperC}, \eqref{OperA&C} and \eqref{OperN}. It is given in Section~\ref{SIAM-Book:Section3.8.5}.

\begin{theorem}
    \label{Th:NonlEvEqApri}
    Suppose that Assumption~{\rm\ref{Asspt-1-POD}} holds. Let $y$ and $y^\ell$ be the solutions to \eqref{NonlEvPro} and \eqref{Eq:Jan-2}, respectively, where we suppose that $y\in H^1(0,T;V)$. Then the following \index{Error estimate!a-priori!state variable}a-priori error estimate is valid:
    \begin{align*}
        \int_0^T{\|y(t)-y^\ell(t)\|}_V^2\,\mathrm dt\le C\cdot\left\{
        \begin{aligned}
            &\sum_{i>\ell}\lambda_i^H\,{\|\psi_i^H\|}_V^2,&&X=H,~\mathcal P^\ell=\mathcal P^\ell_H,\\
            &\sum_{i>\ell}\lambda_i^V\,{\|\psi_i^V-\mathcal Q^\ell_H\psi_i^V\|}_H^2,&&X=V,~\mathcal P^\ell=\mathcal Q^\ell_H
        \end{aligned}
        \right.
    \end{align*}
    for a constant $C>0$ which depends on $y$, but not on $\ell$.
\end{theorem}

\subsection{The semidiscrete approximation}
\label{SIAM-Book:Section3.6.3}

In this section we extend the approach presented in Section~\ref{SIAM-Book:Section3.4} to the non-linear evolution problem \eqref{NonlEvPro}. Our main focus lies on the derivation of the corresponding discretized systems and the application of the POD method. Therefore, we do not proof existence, uniqueness or regularity results. However, we show how the error analysis of Section~\ref{SIAM-Book:Section3.4} can be transferred to the non-linear case.

\subsubsection{Galerkin discretization}
\label{sec:nonlinearGalerkinDiscretiation}

In \eqref{FESpace} we have introduced the $m$-dimensional subspace $V^h\subset V$. We apply a standard Galerkin scheme for \eqref{NonlEvPro}. Thus, we look for a function $y^h$ satisfying $y^h(t)\in V^h$ in $[0,T]$ almost everywhere and
\begin{equation}
    \label{NonlEvProGal}
    \begin{aligned}
        \frac{\mathrm d}{\mathrm dt}\,{\langle y^h(t),\varphi^h \rangle}_H+a(y^h(t),\varphi^h)+{\langle \mathcal Cy^h(t)+\mathcal N(y^h(t)),\varphi^h\rangle}_{V',V}&={\langle(\mathcal F+\mathcal Bu)(t),\varphi^h\rangle}_{V',V}\\
        &\quad\forall\varphi^h\in V^h\text{ a.e. in }(0,T],\\
        y^h(0)&=\mathcal P^h_Hy_\circ,
    \end{aligned}
\end{equation}
where the linear projection operator $\mathcal P^h:X\to V^h$ stands for one of the two projections introduced in Example~\ref{ExampleFEProjection}. Next we derive a system of ordinary differential equations for the coefficients $\mathrm y^h(t)=(\mathrm y_i^h(t))_{1\le i\le m}\in\mathbb R^m$ of the Galerkin ansatz
\begin{align*}
    y^h(t)=\sum_{i=1}^m\mathrm y_i^h(t) \varphi_i^h \in V^h\quad\text{for }t\in[0,T].
\end{align*}
In Section~\ref{Section:3.4.1} we have introduced the mass matrix $\bM^h\in\mathbb R^{m\times m}$, the vector $\mathrm g^h(t;u)\in\mathbb R^m$ for the right-hand side and the vector $\mathrm y_\circ^h\in\mathbb R^m$ for the initial condition. Moreover, we define the matrix
\begin{align*}
    \bA^h=\big(\big(a(\varphi_j,\varphi_i)+{\langle\mathcal C\varphi_j,\varphi_i\rangle}_{V',V}\big)\big)\in\mathbb R^{m\times m}.
\end{align*}
For the non-linearity we introduce the vector
\begin{align*}
    \mathrm n^h(\mathrm v)=\bigg(\Big\langle\mathcal N\Big(\sum_{j=1}^m\mathrm v_j\varphi_j\Big),\varphi_i\Big\rangle_{V',V}\bigg)_{1\le i\le m}
\end{align*}
for the vector $\mathrm v=(\mathrm v_i)_{1\le i\le m}$. Then \eqref{Eq:Jan-2} is equivalent to the following non-linear, $m$-dimensional system of differential equations
\begin{equation}
    \label{NonlFineModel}
    \begin{aligned}
        \bM^h\dot{\mathrm y}^h(t)+\bA^h\mathrm y^h(t)+\mathrm n^h\big(\mathrm y^h(t)\big)&=\mathrm g^h(t;u)\quad\text{for }t\in(0,T],\\
        \mathrm y^h(0)&=\mathrm y_\circ^h.
    \end{aligned}
\end{equation}
Here we do not consider the existence of a unique solution to \eqref{NonlEvProGal}; cf. Section~\ref{Section:3.4.1} for the semi-discretized linear evolution problem. Moreover, we do not study the a-priori error between a solution to \eqref{NonlEvPro} and \eqref{NonlEvProGal}; cf. Section~\ref{SIAM-Book:Section3.4.2} for the semi-discretized linear evolution problem. 

\subsubsection{POD Galerkin scheme}

Suppose that there exists a unique soultion $y^h$ to \eqref{NonlEvProGal}. In the context of Section~\ref{Section:ContPODHilbert} we choose $K=1$, $\omega_1^K=1$ and $y^1=y^h$. Since $y^h\in L^2(0,T;V)$ we can apply the results derived by \cite{Sin14}; as in Section~\ref{Section:ContPODHilbert}. Suppose that $\{(\lambda_i,\psi_i)\}_{i\in\mathbb I}$ denote the eigenvalue-eigenvector pairs solving \eqref{FEPODEigPro}, where the operator $\mathcal R:X\to V^h\subset X$ is given as
\begin{align*}
    \mathcal R\psi=\int_0^T{\langle y^h(t),\psi\rangle}_X\,y^h(t)\,\mathrm dt\quad\text{for }\psi\in X
\end{align*}
and $d=\dim\,\{\mathcal R\psi\,\big|\,\psi\in X\}\le m$. Thus, $\mathcal R$ is a finite rank operator. It follows from \cite[Lemma~3.1-2)]{Sin14} that $d$ is the same for the choices $X=H$ and $X=V$. Then $\{\psi_i\}_{i=1}^\ell\subset V^h$ is a POD basis of rank $\ell$ solving
\begin{equation}
    \label{Eq:Apr-2}
    \left\{
    \begin{aligned}
        & \min\int_0^T \Big\|y^h(t) - \sum_{i=1}^\ell {\langle y^h(t),\psi_i\rangle}_X\,\psi_i\Big\|_X^2\,\mathrm dt\\
        &\hspace{1mm}\text{s.t. }\{\psi_i\}_{i=1}^\ell\subset X\text{ and }{\langle\psi_i,\psi_j\rangle}_X=\delta_{ij} \text{ for } 1 \le i,j \le \ell.
        \end{aligned}
    \right.
\end{equation}
For $\{(\lambda_i^h,\psi_i^h)\}_{i\in\mathbb I}$ satisfying \eqref{FEPODEigPro} we obtain
\begin{equation*}
    \int_0^T\Big\|y^h(t)-\sum_{i=1}^\ell {\langle y^h(t),\psi_i\rangle}_X\,\psi_i\Big\|_X^2\,\mathrm dt=\sum_{i=\ell+1}^d\lambda_i.
\end{equation*}
Following Section~\ref{Section:3.4.1} we infer that there is exactly one coefficient matrix $\bPsi=[\uppsi_1\,|\ldots\,\uppsi_\ell]\in\mathbb R^{m\times\ell}$ satisfying \eqref{FEAnsatzPOD}. We define again the positive definite and symmetric matrix $\bW=((\langle\varphi_j^h,\varphi_i^h\rangle_X))\in\mathbb R^{m\times m}$. Then \eqref{Eq:Apr-2} is equivalent with the minimization problem
\begin{equation}
    \label{Eq:Apr-4}
    \left\{
    \begin{aligned}
        & \min\int_0^T \Big|\mathrm y^h(t) - \sum_{i=1}^\ell {\langle \mathrm y^h(t),\uppsi_i\rangle}_\bW\,\uppsi_i\Big|_\bW^2\,\mathrm dt\\
        &\hspace{1mm}\text{s.t. }\{\uppsi_i\}_{i=1}^\ell\subset \mathbb R^m\text{ and }{\langle\uppsi_i,\uppsi_j\rangle}_\bW=\delta_{ij} \text{ for } 1 \le i,j \le \ell.
    \end{aligned}
    \right.
\end{equation}

Suppose that we have computed a POD basis $\{\psi_i\}_{i=1}^\ell$ of rank $\ell$ with $\ell\ll d\le m$. Again, we set
\begin{align*}
    X^\ell=\mathrm{span}\,\big\{\psi_1,\ldots,\psi_\ell\big\}
\end{align*}
and
\begin{align*}
    H^\ell=\mathrm{span}\,\big\{\psi_1^H,\ldots,\psi_\ell^H\big\},\quad V^\ell=\mathrm{span}\,\big\{\psi_1^V,\ldots,\psi_\ell^V\big\}
\end{align*}
for the choices $X=H$ and $X=V$, respectively. As before, we omit the dependence of $\{\psi_i^H\}_{i=1}^\ell$ and $\{\psi_i^V\}_{i=1}^\ell$ on $h$ for the reader's convenience. Then we replace \eqref{NonlEvProGal} by the following POD Galerkin scheme: find $y^{h\ell}: [0,T] \to X^\ell$ satisfying
\begin{equation}
    \label{NonlEvProGalPOD}
    \begin{aligned}
        \frac{\mathrm d}{\mathrm dt}\,{\langle y^{h\ell}(t),\psi\rangle}_H+a(y^{h\ell}(t),\psi)+{\langle \mathcal Cy^{h\ell}(t)+\mathcal N(y^{h\ell}(t)),\psi\rangle}_{V',V}&={\langle(\mathcal F+\mathcal Bu)(t),\psi\rangle}_{V',V}\\
        &\quad\forall\psi\in X^{h\ell}\text{ a.e. in }(0,T],\\
        y^{h\ell}(0)&=\mathcal P^{h\ell}y_\circ,
    \end{aligned}
\end{equation}
where we choose either $\mathcal P^{h\ell}=\mathcal P_H^{h\ell}$ for $X=H$ or $\mathcal P^{h\ell}=\mathcal P_{HV^{h\ell}}^{h\ell}$ for $X=V$.

We proceed as in Section~\ref{SIAM-Book:Section3.4.4} and derive the following non-linear system of ordinary differential equations: 
\begin{equation}
    \label{NonlSystPOD}
    \begin{aligned}
        \bM^{h\ell}\dot{\mathrm y}^{h\ell}(t)+\bA^{h\ell}\mathrm y^{h\ell}(t)+\mathrm n^{h\ell}\big(\mathrm y^{h\ell}(t)\big)&=\mathrm g^{h\ell}(t;u)\quad\text{for }t\in(0,T],\\
        \mathrm y^{h\ell}(0)&=\mathrm y_\circ^{h\ell},
    \end{aligned}
\end{equation}
where the mass matrix $\bM^{h\ell}\in\mathbb R^{\ell\times\ell}$ as well as the vectors $\mathrm y^{h\ell}$, $\mathrm g^{h\ell}(t;u)$, $\mathrm y_\circ^{h\ell}\in\mathbb R^\ell$ have been introduced in Section~\ref{SIAM-Book:Section3.4.4}. Moreover, we have set
\begin{align*}
    \bA^{h\ell}&=\big(\big(a(\psi_j,\psi_i)+{\langle\mathcal C\psi_j,\psi_i\rangle}_{V',V}\big)\big)\in\mathbb R^{\ell\times\ell},\\
    \mathrm n^{h\ell}(\mathrm v)&=\bigg(\Big\langle\mathcal N\Big(\sum_{j=1}^\ell\mathrm v_j\psi_j\Big),\varphi_i\Big\rangle_{V',V}\bigg)_{1\le i\le\ell}\quad\text{for the vector }\mathrm v=\big(\mathrm v_i\big)_{1\le i\le\ell}.
\end{align*}

\subsubsection{POD a-priori error analysis}

Let us suppose the following hypothesis.

\begin{assumption}
    \label{A7}
    \begin{enumerate}
        \item [\em 1)] Assumption {\rm\ref{Asspt-1-POD}} holds.
        \item [\em 2)] There exists a unique solution $y^h$ of \eqref{NonlEvProGal} satisfying
        \begin{align*}
            {\|y^h\|}_{L^\infty(0,T;V)}\le C\Big({\|\mathcal P^h y_\circ\|}_H+{\|f\|}_{L^2(0,T;H)}+{\|u\|}_\U\Big).
        \end{align*}
        \item [\em 3)] The POD basis $\{\psi_i\}_{i=1}^\ell$ of rank $\ell\in\{1,\ldots,d\}$ is computed by solving \eqref{Eq:Apr-2}.
        \item [\em 4)] There exists a unique solution $y^{h\ell}\in H^1(0,T;V)$ to the POD Galerkin scheme \eqref{NonlEvProGalPOD}.
    \end{enumerate}
\end{assumption}

Similarly as in Section~\ref{SIAM-Book:Section3.4.5} we consider the difference
\begin{align*}
    \int_0^T{\|y^h(t)-y^{h\ell}(t)\|}_V^2\,\mathrm dt
\end{align*}
and present an a-priori error estimate utilizing the error formulas of Section~\ref{Section:DiscPODHilbert} and the assuption $y^\ell \in H^1(0,T;V)$. This is stated in the next theorem which is proved in Section~\ref{SIAM-Book:Section3.8.5}.

\begin{theorem}
    \label{Th:A-PrioriError-300}
    Let Assumption~{\rm\ref{A7}} be satisfied. Then
    \begin{enumerate}
        \item[\rm 1)] For $X=H$, $\mathcal P^\ell=\mathcal P^\ell_H$ and $\mathcal P^h=\mathcal P^h_H$ we have
        \begin{align*}
            {\|y^h-y^{h\ell}\|}_{L^2(0,T;V)}^2\le C\sum\limits_{i=\ell+1}^d\lambda_i^H\,{\|\psi_i^H\big\|}_V^2.
        \end{align*}
        \item[\rm 2)] For $X=V$ and $\mathcal P^\ell=\mathcal Q^\ell_V$ and $\mathcal P^h=\mathcal P^h_V$ we derive
        \begin{align*}
            {\|y^h-y^{h\ell}\|}_{L^2(0,T;V)}^2\le C\sum\limits_{i=\ell+1}^d\lambda_i^V\,\big\|\psi_i^V-\mathcal Q^\ell_V\psi_i^V\big\|_V^2.
        \end{align*}
    \end{enumerate}
\end{theorem}

\subsubsection{POD a-posteriori error analysis}

A-posteriori error analysis\index{Error estimate!a-posteriori!state variable} can be utilized to estimate the difference $y^h-y^{h\ell}$ without the knowledge of $y^h$ and without any assumptions on the choice of the POD basis. For that reason we present also a result for our non-linear evolution problem. For almost all $t\in [0,T]$ let us define the time dependent residual $r^{h\ell}(t)\in (V^h)'$ by
\begin{equation}
    \label{Eq:timeResidual}
    r^{h\ell}(t)=a(y^{h\ell}(t),\cdot)+{\langle y_t^{h\ell}(t)+\mathcal Cy^{h\ell}(t)+\mathcal N(y^{h\ell}(t))-(\mathcal F+\mathcal Bu)(t),\cdot\rangle}_{V',V}.
\end{equation}
Then the next result holds. For a proof we refer to Section~\ref{SIAM-Book:Section3.8.5}.

\begin{proposition}
    \label{Prop:NSAposti}
    Suppose that Assumption~{\rm\ref{A7}} holds. By $y^h$ and $y^{h\ell}$ we denote the solutions to \eqref{NonlEvProGal} and \eqref{NonlEvProGalPOD}, respectively. We assume that $\|y^{h\ell}\|_{L^2(0,T;V)}$ is bounded independent of $h$ and $\ell$. Then the following \index{Error estimate!a-posteriori!state variable}a-posteriori error estimates hold for almost all $t\in[0,T]$
    \begin{align*}
        {\|y^h(t)-y^{h\ell}(t)\|}_H^2&\le \mathsf C_1^{h\ell}(t)\bigg(\mathsf R_\circ^{h\ell}+\int_0^t\mathsf R_1^{h\ell}(s)\,\mathrm ds\bigg),\\
        \int_0^t{\|y^h(s)-y^{h\ell}(s)\|}_V^2\,\mathrm ds&\le\tilde{\mathsf R}_\circ^{h\ell}+\int_0^t\tilde{\mathsf R}_1^{h\ell}(s)+\mathsf R_2^{h\ell}(s)\bigg(\mathsf R_\circ^{h\ell}+\int_0^s\mathsf R_1^{h\ell}(\tau)\,\mathrm d\tau\bigg)\,\mathrm ds,
    \end{align*}
    where
    \begin{subequations}
        \begin{align}
            \label{Alb-4a}
            \mathsf C_1^{h\ell}(t)&=\exp\bigg(\frac{2 c_\mathcal N^2}{\gamma_1}\int_0^t{\|y^{h\ell}(s)\|}_V^2\,\mathrm ds\bigg),\quad\mathsf R_\circ^{h\ell}={\|(\mathcal P^h-\mathcal P^\ell)y_\circ\|}_H^2,\\
            \label{Alb-4b}
            \mathsf R_1^{h\ell}(s)&=\frac{2}{\gamma_1}\,{\|r^{h\ell}(s)\|}_{(V^h)'}^2,\\
            \label{Alb-5}
            \tilde{\mathsf R}_\circ^{h\ell}&=\frac{\mathsf R_\circ^{h\ell}}{\gamma_1},\quad\tilde{\mathsf R}_1^{h\ell}(s)=\frac{1}{\gamma_1}\mathsf R_1^{h\ell}(s),\quad\mathsf R_2^{h\ell}(s)=\frac{2 c_\mathcal N^2}{\gamma_1^2}\,\mathsf C_1^{h\ell}(\tau)\,{\|y^{h\ell}(s)\|}_V^2
        \end{align}
    \end{subequations}
    for almost all $s\in[0,T]$.
\end{proposition}

\subsection{The fully discrete approximation}
\label{SIAM-Book:Section3.6.4}

To solve \eqref{NonlEvProGal} and \eqref{NonlEvProGalPOD} numerically, a temporal discretization is required for the interval $[0,T]$. This is done in this subsection.

\subsubsection{Temporal and Galerkin discretization}

We apply the same time discretization as in Section~\ref{SIAM-Book:Section3.5.1}: Let $0=t_1 < t_2 < \ldots < t_n=T$ be a given grid in $[0,T]$ with step sizes $\delta t_j=t_j-t_{j-1}$ for $2 \le j \le n$. We set
\begin{align*}
    \delta t=\min_{2\le j\le n}\delta t_j\quad\text{and}\quad\Delta t=\max_{2\le j\le n}\delta t_j.
\end{align*}
Recall the definition of the functions $\{g_j(u)\}_{j=1}^n\subset V'$, where
\begin{align*}
    g_j(u)&=\frac{2}{\delta t_1}\int_0^{\delta t_1/2}(\mathcal F+\mathcal Bu)(s)\,\mathrm ds,\\
    g_j(u)&=\frac{2}{\delta t_{j-1}+\delta t_j}\int_{t_j-\delta t_{j-1}/2}^{t_j+\delta t_j/2}(\mathcal F+\mathcal Bu)(s)\,\mathrm ds\quad\text{for }j=2,\ldots,n-1,\\
    g_n(u)&=\frac{2}{\delta t_n}\int_{T-\delta t_n/2}(\mathcal F+\mathcal Bu)(s)\,\mathrm ds
\end{align*}
for every $u\in\mathcal U$. To solve \eqref{NonlEvProGal} we apply the implicit Euler method. Then the sequence $\{y^h_j\}_{j=1}^n\subset V^h$ satisfies
\begin{equation}
    \label{NonlEvProGalFullyDisc}
    \begin{aligned}
        {\langle \overline\partial_ty^h_j,\varphi^h \rangle}_H+a(y^h_j,\varphi^h)+{\langle \mathcal Cy^h_j+\mathcal N(y^h_j),\varphi^h\rangle}_{V',V}&={\langle g_j(u),\varphi^h\rangle}_{V',V}\\
        &\quad\forall\varphi^h\in V^h\text{ a.e. in }(0,T],\\
        y^h(0)&=\mathcal P^h_Hy_\circ.
    \end{aligned}
\end{equation}
We interprete $y_j^h$ as approximatios for the solution $y^h$ to \eqref{NonlEvProGal} at time $t=t_j$ for $1\le j\le n$. From $y_j^h\in V^h$, $j=1,\ldots,n$, we infer that
\begin{equation}
\label{Alb-6}
    y^h_j=\sum_{i=1}^m\mathrm y_{ji}^h\varphi_i^h\in V^h\quad\text{for }1\le j\le n
\end{equation}
with coefficient vectors $\mathrm y^h_j=(\mathrm y_{ji}^h)_{1\le i\le m}\in\mathbb R^m$. Inserting \eqref{Alb-6} into \eqref{NonlEvProGalFullyDisc} and choosing $\varphi^h=\varphi_i^h$, $i=1,\ldots,m$, we find that the sequence $\{\mathrm y_j^h\}_{j=1}^n\subset\mathbb R^m$ of coefficient vectors satisfies the following non-linear initial value problem
\begin{subequations}
    \label{Nonl:FineModel-Disc}
    \begin{align}
        \label{Nonl:FineModel-Disc-1}
        \big(\bM^h+\delta t_j\bA^h\big)\mathrm y_j^h+\mathrm n^h(\mathrm y^h_j)&=\bM^h\mathrm y^h_{j-1}+\delta t_j\mathrm g_j^h(u)\text{ for } 2 \le j \le n,\\
        \label{Nonl:FineModel-Disc-2}
        \bM^h\mathrm y^h_1&=\mathrm y_\circ^h,
    \end{align}
\end{subequations}
where $\mathrm g_j^h(u)=({\langle g_j(u),\varphi_i^h\rangle}_{V',V})_{1\le i\le m}\in\mathbb R^m$ for $1\le j\le n$. The following result is proved in Section~\ref{SIAM-Book:Section3.8.5}.

\begin{proposition}
    \label{PropNonl:FullDiscFEModel}
    Let Assumption~{\em\ref{Asspt-1}} hold. If the time step $\delta t$ is sufficiently small, there exists a unique sequence $\{y^h_j\}_{j=1}^n \subset V^h$ to \eqref{NonlEvProGalFullyDisc}. Furthermore, we have the a-priori estimates
    \begin{equation}
        \label{Nonl:TempDiscFESOL-APriori}
        {\|y_j^h\|}_H^2\le e^{-C_1(j-1)\delta t}\bigg({\|\mathcal P^hy_\circ\|}_H^2+\frac{1}{\gamma_1}\sum_{l=2}^j\delta t_l\,{\|g_l(u)\|}_{V'}^2\bigg)
    \end{equation}
    for $j=2,\ldots,n$ with the constant $C_1=2\gamma_1/c_V^2$ and
    \begin{equation}
        \label{Nonl:TempDiscFESOL-APriori-2}
        {\|y_n^h\|}_H^2+\sum_{i=1}^n{\|y_j^h-y_{j-1}^h\|}_H^2+\gamma_1\sum_{j=1}^n\delta t_j\,{\|y_j^h\|}_V^2\le {\|\mathcal P^hy_\circ\|}_H^2+ \frac{1}{\gamma_1}\sum_{j=2}^n\delta t_j\,{\|g_j(u)\|}_{V'}^2.
    \end{equation}
\end{proposition}

\subsubsection{POD basis computation}

Let $y_\circ\in V$, $f\in L^2(0,T;V')$ and $u\in\U$. By
\begin{equation}
    \label{Nonl:FullyDiscGal}
    y^h_j=\sum_{i=1}^m\mathrm y^h_{ji}\varphi_i^h\quad\text{for }1\le j\le n
\end{equation}
we denote the Galerkin solution to \eqref{NonlEvProGalFullyDisc}. We introduce the coefficient vectors $\mathrm y^h_j=(\mathrm y^h_{ji})_{1\le i\le m}\in\mathbb R^m$ for $j=1,\ldots,n$. In Section~\ref{SIAM-Book:Section3.5.3} we have introduced the positive weighting parameters $\{\alpha^n_j\}_{j=1}^n$. Let us choose $\wp=1$, $\omega_1^\wp=1$ and the snapshots
\begin{align*}
    y_j^1=y_j^h\in V^h\quad\text{for }j=1,\ldots,n.
\end{align*}
The snapshot subspace is given as $\mathscr V^n=\Span\,\{y_j^h\,|\,1\le j\le n\}\subset V^h$ with dimension $d^n=\dim\mathscr V^n\le\min(n,m)$. For either $X=H$ or $X=V$ we compute a POD basis of rank $\ell\in\{1,\ldots,d^n\}$ by solving
\begin{equation}
    \label{Nonl:PODProb}
    \left\{
    \begin{aligned}
        &\min\sum_{j=1}^n\alpha_j^n\,\Big\|y_j^h-\sum_{i=1}^\ell{\langle y_j^h,\psi_i^n\rangle}_X\,\psi_i^n\Big\|_X^2\\
        &\hspace{1mm}\text{s.t. }\{\psi_i^n\}_{i=1}^\ell\subset V^h\subset X\text{ and }{\langle\psi_i^n,\psi_j^n\rangle}_X=\delta_{ij},~1\le i,j\le \ell.
    \end{aligned}
    \right.
\end{equation}
Recall that $\{\psi_i^n\}_{i=1}^\ell$ satisfies
\begin{align*}
    \mathcal R^n\psi_i^n=\lambda_i^n\psi_i^n\text{ for }1\le i\le \ell,\quad\lambda_1^n\ge\ldots\lambda_\ell^n>0,
\end{align*}
where
\begin{align*}
    \mathcal R^n\psi=\sum_{j=1}^n\alpha_i^n\,{\langle y_j^h,\psi\rangle}_X\,y^h_j\quad\text{for }\psi\in X
\end{align*}
with either $X=H$ or $X=V$. Utilizing the symmetric and positive definite weighting matrix $\bW=((\langle\varphi_j^h,\varphi_i^h\rangle_X))\in\mathbb R^{m\times m}$ and the associated weighted inner product we can express the minimization problem as
\begin{equation}
    \label{Nonl:PODProbDisc}
    \left\{
    \begin{aligned}
        &\min\sum_{j=1}^n\alpha_j^n\Big|\mathrm y_j^h-\sum_{i=1}^\ell{\langle \mathrm y_j^h,\uppsi_i\rangle}_\bW\,\uppsi_i\Big|_\bW^2\\
        &\hspace{1mm}\text{s.t. }\{\uppsi_i\}_{i=1}^\ell\subset\mathbb R^m\text{ and }{\langle\uppsi_i,\uppsi_j\rangle}_\bW=\delta_{ij},\quad1\le i,j\le \ell.
    \end{aligned}
    \right.
\end{equation}
We follow Section~\ref{SIAM:Section-2.1.1.4} and set
\begin{align*}
    \bY&=\big[\mathrm y_1^h\,|\ldots|\,\mathrm y_n^h\big]\in\mathbb R^{m\times n},&\bD&=\mathrm{diag}\,(\alpha_1^n,\ldots,\alpha_n^n)\in\mathbb R^{n\times n},\\
    \tilde{\bY}&=\bW^{1/2}\bY\tilde{\bD}^{1/2}.
\end{align*}
Then a solution to \eqref{Nonl:PODProbDisc} is given by $\uppsi_i=\bW^{-1/2}\tilde \uppsi_i\in\mathbb R^m$, $i=1,\ldots,\ell$, where the $\tilde\uppsi_i$'s solve the symmetric $m\times m$ eigenvalue problem
\begin{align*}
    \tilde{\bY}\tilde{\bY}^\top\tilde\uppsi_i=\lambda_i^n\tilde\uppsi_i\quad\text{with }\lambda_1^n\ge\ldots\ge\lambda_{d^n}^n>\lambda_{d^n+1}^n=\ldots=\lambda_m\ge0.
\end{align*}
If $\{\uppsi_i\}_{i=1}^\ell$ solves \eqref{Nonl:PODProbDisc} then $\{\psi_i^n\}_{i=1}^\ell$ is a solution to \eqref{Nonl:PODProb} with the representation \eqref{FEAnsatzPOD-2}. On the other hand, if $\{\psi_i^h\}_{i=1}^\ell$ is a solution to \eqref{Nonl:PODProb} then we infer from $\psi_i^n\in V^h$ the existence of a uniquely determined matrix $\bPsi=[\uppsi_1\,|\ldots|\,\uppsi_\ell]\in\mathbb R^{m\times\ell}$ satisfying \eqref{FEAnsatzPOD-2}. It follows that the columns $\{\uppsi_i\}_{i=1}^\ell$ of $\bPsi$ solve \eqref{Nonl:PODProbDisc}. Finally, we can quantify the POD approximation error as follows
\begin{equation}
    \label{POD-Error-Formula-100-2}
    \sum_{j=1}^n\alpha_j^n\,\Big\|y_j^h-\sum_{i=1}^\ell{\langle y_j^h,\psi_i^n\rangle}_X\,\psi_i^n\Big\|_X^2=\sum_{i=\ell+1}^{d^n}\lambda_i^n.
\end{equation}

\subsubsection{POD Galerkin scheme}

Suppose that we have computed a POD basis $\{\psi_i^n\}_{i=1}^\ell$ of rank $\ell$. Again, we set
\begin{align*}
    X^{n\ell}=\Span\,\big\{\psi_1^n,\ldots,\psi_\ell^n\big\}\subset V^h \subset X
\end{align*}
and
\begin{align*}
    H^{n\ell}=\mathrm{span}\,\big\{\psi_1^{nH},\ldots,\psi_\ell^{nH}\big\},\quad V^{n\ell}=\mathrm{span}\,\big\{\psi_1^{nV},\ldots,\psi_\ell^{nV}\big\}
\end{align*}
for the choices $X=H$ and $X=V$, respectively. To simplify the presentation we utilize the same time grid $\{t_j\}_{j=1}^n$ for \eqref{NonlEvProGalFullyDisc}, \eqref{Nonl:PODProb} and the POD Galerkin scheme introduced below in \eqref{Nonl:EvProGalPOD-disc}. Next we replace \eqref{NonlEvProGalFullyDisc} by the following POD Galerkin scheme: find $\{y^{h\ell}_j\}_{j=1}^n\subset X^{h\ell}$ satisfying
\begin{subequations}
    \label{Nonl:EvProGalPOD-disc}
    \begin{equation}
        \label{Nonl:EvProGalPOD-disc-1}
        {\langle\overline\partial y^{h\ell}_j,\psi\rangle}_H+a(y^{h\ell}_j,\psi)+{\langle\mathcal Cy_j^{h\ell}+\mathcal N(y^{h\ell}_j),\psi\rangle}_{V',V}={\langle g_j(u),\psi\rangle}_{V',V}
    \end{equation}
    for all $\psi \in X^{n\ell}$, $2\le j\le n$ and
    \begin{equation}
        \label{Nonl:EvProGalPOD-disc-2}
        y^{h\ell}_1=\mathcal P^{n\ell} y_\circ
    \end{equation}
\end{subequations}
with $\mathcal P^{n\ell}=\mathcal P^{n\ell}_H$ if $X=H$ or with $\mathcal P^{n\ell}=\mathcal Q^{n\ell}_H$ if $X=V$. In \eqref{Nonl:EvProGalPOD-disc-1} we set $\overline\partial y^{h\ell}_j=(y^{h\ell}_j-y^{h\ell}_{j-1})/\delta t_j$. From $y_j^{h\ell}\in X^{h\ell}$ we infer that there are uniquely determined coefficient vectors $\mathrm y_j^{h\ell}=(\mathrm y_{ji}^{h\ell})_{1\le i\le\ell}\in\mathbb R^\ell$ satisfying
\begin{equation}
    \label{Nonl:POD-Gal-Ans}
    y_j^{h\ell}=\sum_{i=1}^\ell\mathrm y_{ji}^{h\ell}\psi_i^n\quad\text{for }j=1,\ldots,n.
\end{equation}
Inserting \eqref{Nonl:POD-Gal-Ans} into \eqref{Nonl:EvProGalPOD-disc}, choosing $\psi=\psi_i^n$, $1\le i\le \ell$, and using the notation introduced in Section~\ref{SIAM-Book:Section3.4.4} we obtain
\begin{subequations}
    \label{Nonl:PODModel-Disc}
    \begin{align}
        \label{Nonl:PODModel-Disc-1}
        \big(\bM^{h\ell}+\delta t_j\bA^{h\ell}\big)\mathrm y_j^{h\ell}+\mathrm n^{h\ell}(\mathrm y^{h\ell}_j)&=\bM^{h\ell}\mathrm y^{h\ell}_{j-1}+\delta t_j\mathrm g_j^{h\ell}(u)\text{ for } 2 \le j \le n,\\
        \label{Nonl:PODModel-Disc-2}
        \bM^{h\ell}\mathrm y^h_1&=\mathrm y_\circ^{h\ell}
    \end{align}
\end{subequations}
and
\begin{align*}
    \mathrm n^{h\ell}(\mathrm v)=\bigg(\Big\langle\mathcal N\Big(\sum_{j=1}^\ell\mathrm v_j\psi_j^n\Big),\psi_i^n\Big\rangle\bigg)_{1\le i\le \ell}
\end{align*}
for the vector $\mathrm v=(\mathrm v_i)_{1\le i\le \ell}$. In the following proposition, existence and a-priori estimates for the solution $\{y^{h\ell}_j\}_{j=1}^n$ are established. The proof follows by the same arguments as in the proof of Proposition~\ref{PropNonl:FullDiscFEModel}.

\begin{proposition}
    \label{PropExUnPODDisc}
    Let Assumption~{\em\ref{Asspt-1-POD}} hold. The POD basis $\{\psi_i^n\}_{i=1}^\ell$ of rank $\ell\in\{1,\ldots,d\}$ is computed by solving \eqref{Nonl:PODProb}. If the time step $\delta t$ is sufficiently small, there exists a unique sequence $\{y^{h\ell}_j\}_{j=1}^n \subset X^{n\ell}$ to \eqref{Nonl:EvProGalPOD-disc}. Moreover, the following estimates are satisfied:
    \begin{align*}
        {\|y_j^{h\ell}\|}_H^2\le e^{-C_1(j-1)\delta t}\bigg({\|\mathcal P^{n\ell}y_\circ\|}_H^2+\frac{1}{\gamma_1}\sum_{l=2}^j\delta t_l\,{\|g_l(u)\|}_{V'}^2\bigg)
    \end{align*}
    for $j=2,\ldots,n$ with the constant $C_1=2\gamma_1/c_V^2$ and
    \begin{align*}
        {\|y_n^{h\ell}\|}_H^2+\sum_{j=2}^n{\|y_j^{h\ell}-y_{j-1}^{h\ell}\|}_H^2+\gamma_1\sum_{j=1}^n\delta t_j\,{\|y_j^{h\ell}\|}_V^2\le {\|\mathcal P^{n\ell}y_\circ\|}_H^2+ \frac{1}{\gamma_1}\sum_{j=2}^n\delta t_j\,{\|g_j(u)\|}_{V'}^2.
\end{align*}
\end{proposition}

\subsubsection{POD a-priori error analysis}

Let $\{y^h_j\}_{j=1}^n$ and $\{y^{h\ell}_j\}_{j=1}^n$ be solutions to \eqref{NonlEvProGalFullyDisc} and \eqref{Nonl:EvProGalPOD-disc}, respectively. Our goal is to derive an a-priori error estimate for the term
\begin{align*}
    \sum_{j=1}^n\alpha_j^n\,{\|y_j^h-y_j^{h\ell}\|}_V^2.
\end{align*}
Here we utilize the results presented in Section~\ref{Section:DiscPODHilbert} and proceed similarly as in Section~\ref{SIAM-Book:Section3.5.5}. 

In the following assumption we suppose that all hypotheses of Propositions~\ref{PropNonl:FullDiscFEModel} and \ref{PropExUnPODDisc} are satisfied so that there exist unique solutions $\{y^h_j\}_{j=1}^n$ and $\{y^{h\ell}_j\}_{j=1}^n$ to \eqref{NonlEvProGalFullyDisc} and \eqref{Nonl:EvProGalPOD-disc}, respectively.

\begin{assumption}
    \label{AssNL5}
    \begin{enumerate}
        \item [\em 1)] Suppose that Assumption~{\em\ref{Asspt-1-POD}} hold. 
        \item [\em 2)] There is a constant $\zeta>0$ which is independent of $n$, such that $\Delta t/\delta t\le\zeta$.
        \item [\em 3)] Let the POD basis $\{\psi_i^n\}_{i=1}^\ell$ of rank $\ell$ be computed from \eqref{Nonl:PODProb}.
        \item [\em 4)] Assume that $\delta t$ is sufficiently small so that there exist unique solutions $\{y^h_j\}_{j=1}^n$ and $\{y^{h\ell}_j\}_{j=1}^n$ to \eqref{NonlEvProGalFullyDisc} and \eqref{Nonl:EvProGalPOD-disc}. Moreover, the sequence $\{y^h_j\}_{j=1}^n$ satisfies
        \begin{equation}
            \label{Bach-7}
            \max_{1\le j\le n}{\|y_j^h\|}_V\le C
        \end{equation}
        with a constant $C>0$ which does not depend on $n$ and $h$.
    \end{enumerate}
\end{assumption}

The proof of the next theorem is given in Section~\ref{SIAM-Book:Section3.8.5}.

\begin{theorem}
    \label{Theorem:AprioriEstNonl}
    Suppose that Assumption~{\em\ref{AssNL5}} hold. Then
    \begin{equation}
        \label{NLEq-ApostiE}
        \sum_{j=1}^n\alpha_j^n\,{\|y_j^h-y_j^{h\ell}\|}_V^2\le C\cdot\left\{
        \begin{aligned}
            &\sum_{i=\ell+1}^{d^n}\lambda_i^{nH}\,{\|\psi_i^{nH}\|}_V^2&&X=H,~\mathcal P^{n\ell}=\mathcal P^{n\ell}_H,\\
            &\sum_{i=\ell+1}^{d^n}\lambda_i^{nV}\,{\|\psi_i^{nV}-\mathcal Q^{h\ell}_H\psi_i^{nV}\|}_V^2&&X=V,~\mathcal P^{n\ell}=\mathcal Q^{n\ell}_H
        \end{aligned}
        \right.
    \end{equation}
    with a constant $C>0$ that is independent of $n$, $h$ and $\ell$.
\end{theorem}

\subsubsection{POD a-posteriori error analysis}

Let Assumption~\ref{Asspt-1-POD} hold. We suppose that $\{y_j^h\}_{j=1}^n$ and $\{y_j^{h\ell}\}_{j=1}^n$ are solutions to \eqref{NonlEvProGalFullyDisc} and \eqref{Nonl:EvProGalPOD-disc}, respectively. To estimate the error between the (unknown) $y_j^h$ and the POD solution $y_j^h$ for $j=1,\ldots,n$ we derive an a-posteriori error estimate in the next theorem which is proved in Section~\ref{SIAM-Book:Section3.8.5}. This extends our results for the linear case presented in Section~\ref{SIAM-Book:Section3.5.6}.

\begin{theorem}
    \label{NonlEq:AposterioriProp}
    Suppose that Assumption~{\em\ref{AssNL5}} hold. Then the \index{Error estimate!a-priori!state variable}a-priori error estimate
    \begin{align*}
        {\|y_j^h-y^{h\ell}\|}_H^2\le\big(1+\hat c\Delta t\big)^{j-1}{\|(\mathcal P^h-\mathcal P^{n\ell})y_\circ\|}_H^2 + \frac{1}{\gamma_1}\sum_{l=2}^j\big(1+\hat c\Delta t\big)^{j+1-l}\delta t_l\,{\|r_l^{h\ell}\|}_{(V^h)'}^2
    \end{align*}
    for the constant $\hat c=2c_\mathcal N^2C^2/\gamma_1$ and the residual
    \begin{align*}
        r_j^{h\ell}=-\left({\langle\overline\partial y_j^{h\ell}+\mathcal Cy_j^{h\ell}+\mathcal N(y_j^{h\ell})-g_j(u),\cdot\rangle}_{V',V}+a(y_j^{h\ell},\cdot)\right)\in (V^h)'.
    \end{align*}
\end{theorem}

\subsection{The empirical interpolation method}
\label{SIAM-Book:Section3.6.5}

The ROM introduced in \eqref{NonlSystPOD} is a nonlinear system. Hence, the problem with the POD Galerkin approach is the complexity of the evaluation of the nonlinearity $\mathrm n^{h\ell}$ in \eqref{NonlSystPOD}. Setting $\bPsi=((\Psi_{ij}))\in \mathbb R^{m \times \ell}$ with $\bPsi=[\uppsi_1\,|\ldots\,\uppsi_\ell]\in \mathbb R^{m \times \ell}$ we derive for any $\mathrm y\in\mathbb R^\ell$ and $i\in\{1,\ldots,\ell\}$ that
\begin{align*}
    \big(\mathrm n^{h\ell}(\mathrm y)\big)_i&=\bigg\langle\mathcal N\Big(\sum_{j=1}^\ell\mathrm y_j\psi_j^h\Big),\psi_i \bigg\rangle_{V',V}\\
    &=\bigg\langle\mathcal N\bigg(\sum_{l=1}^m\Big(\sum_{j=1}^\ell\Psi_{lj}\mathrm y_j\Big)\varphi_l\bigg),\sum_{k=1}^m\Psi_{ki}\varphi_k\bigg\rangle_{V',V}\\
    &=\sum_{k=1}^m\Psi_{ki}\,\bigg\langle\mathcal N\bigg(\sum_{l=1}^m(\bPsi\mathrm y)_l\varphi_l\bigg),\varphi_k\bigg\rangle_{V',V}\\
    &=\sum_{k=1}^m(\bPsi^\top)_{ik}\,\big(\mathrm n^h(\bPsi\mathrm y)\big)_k=\big(\bPsi^\top\mathrm n^h(\bPsi\mathrm y)\big)_i.
\end{align*}
This implies
\begin{align*}
    \mathrm n^{h\ell}(\mathrm y)=\bPsi^\top \mathrm n^h(\bPsi\mathrm y)\quad\text{for any }\mathrm y\in\mathbb R^\ell.
\end{align*}
This can be interpreted in the way that the vector $\mathrm y\in\mathbb R^\ell$ is first expanded to a vector $\bPsi\mathrm y$ of dimension $m$, then the nonlinearity $\mathrm n^h(\bPsi\mathrm y)$ is evaluated and at last the result is reduced back to the low dimension $\ell$ of the reduced-order model. This is computationally expensive. Further this means that our reduced-order model is not independent of the full dimension $m$. Note that when applying a Newton method to the system \eqref{NonlSystPOD} the Jacobian of the nonlinearity is also needed. For any $\mathrm y\in\mathbb R^\ell$ and $i,\nu,\in\{1,\ldots,\ell\}$ we have
\begin{align*}
    \frac{\partial}{\partial\mathrm y_\nu}\,\big(\mathrm n^{h\ell}(\mathrm y)\big)_i&=\frac{\partial}{\partial\mathrm y_\nu}\,\Bigg(\sum_{k=1}^m\Psi_{ki}\,\bigg\langle\mathcal N\bigg(\sum_{l=1}^m\Big(\sum_{j=1}^\ell\Psi_{lj}\mathrm y_j\Big)\varphi_l\bigg),\varphi_k\bigg\rangle_{V',V}\Bigg)\\
    &=\sum_{k=1}^m\Psi_{ki}\,\bigg\langle\mathcal N'\bigg(\sum_{l=1}^m(\bPsi\mathrm y)_l\varphi_l\bigg)\sum_{\mu=1}^m\Psi_{l\mu}\varphi_\mu,\varphi_k\bigg\rangle_{V',V}\\
    &=\sum_{k=1}^m\Psi_{ki}\sum_{\mu=1}^m\Psi_{l\mu}\,\bigg\langle\mathcal N'\bigg(\sum_{l=1}^m(\bPsi\mathrm y)_l\varphi_l\bigg)\varphi_\mu,\varphi_k\bigg\rangle_{V',V}.
\end{align*}
Consequently,
\begin{align*}
    \nabla\mathrm n^{h\ell}(\mathrm y)=\bPsi^\top\bD^h(\bPsi\mathrm y)\bPsi\quad\text{for any }\mathrm y\in\mathbb R^\ell
\end{align*}
with
\begin{align*}
    \mathrm D^h(\bPsi\mathrm y)=\bigg(\Big\langle\mathcal N'\Big(\sum_{k=1}^m(\bPsi\mathrm y)_k\varphi_k\Big)\varphi_j,\varphi_i\Big\rangle_{V',V}\bigg)_{1 \le i,j \le m}.
\end{align*}
Again the same problem can be observed. Note that here the computation expenses are larger since the Jacobians are of dimension $m\times m$. Hence not only a vector is transformed but a matrix of full dimension. To avoid this computational expensive evaluation, the \index{Method!empirical interpolation, EIM}{\em empirical interpolation method (EIM)} was introduced \cite{BMNP04}. This method is often used in combination with the reduced basis approach \cite{Gre12}. The second approach we will investigate here is the \index{Method!discrete empirical interpolation, DEIM}{\em discrete empirical interpolation method} (DEIM) as introduced in \cite{CS09,CS10,CS12}. While the EIM implementation is based on a greedy algorithm, the DEIM strategy is based on a POD approach combined with a greedy algorithm. We will now discuss both methods. We define
\begin{align*}
    \mathrm b^h=\mathrm n^h(\bPsi\mathrm y)\in \mathbb R^m \quad \text{for }\mathrm y\in\mathbb R^\ell.
\end{align*}
Now $\mathrm b^h$ is approximated by a Galerkin ansatz utilizing $\mathsf p$ linearly independent functions $\phi_1,\ldots,\phi_\mathsf p\in\mathbb R^m$, i.e.
\begin{equation}
    \label{eq:EIM_approx}
    \mathrm b^h \approx\sum_{k=1}^\mathsf p\mathrm c_k \phi_k=\Phi\mathrm c
\end{equation}
with $\mathrm c=[\mathrm c_1\,|\ldots|\,\mathrm c_\mathsf p]^\top \in \mathbb R^\mathsf p$  and $\Phi=
[\phi_1|\ldots|\phi_\mathsf p]\in\mathbb R^{m\times\mathsf p}$. Hence we can write the approximation of $\mathrm n^{h\ell}$ as
\begin{align*}
    \mathrm n^{h\ell}(\mathrm y)=\Psi^\top\mathrm n^h(\Psi\mathrm y)=\Psi^\top \mathrm b^h\approx \Psi^\top\Phi\mathrm c.
\end{align*}
The question arising is how to compute the matrix $\bPhi$ and the vector $\mathrm c$. Let $\vec\imath\in\mathbb R^\mathsf p$ be an index vector and $\bB\in\mathbb R^{m\times\mathsf p}$ a given matrix. Then by $\bB_{\vec \imath}$ we denote the submatrix consisting of the rows of $\bB$ corresponding to the indices in $\vec\imath$. Obviously, if we choose $\mathsf p$ indices then the overdetermined system $\mathrm b^h=\bPhi\mathrm c$ can be solved by choosing $\mathsf p$ rows of $\mathrm b^h$ and $\bPhi$. Here it is assumed that the submatrix $\Phi_{\vec \imath} \in \mathbb R^{\mathsf p\times\mathsf p}$ is invertible.

Assume we have computed $\bPhi$ and $\vec\imath$ by an algorithm. Then we proceed as follows. For simplicity we introduce here the matrix $\bP=[e_{\vec\imath_1}\,|\ldots|\,e_{\vec \imath_\mathsf p}]\in \mathbb R^{m \times\mathsf p}$, where $e_{\vec \imath_i} = (0,\ldots,0,1,0,\ldots,0)^\top\in\mathbb R^m$ is a vector with all zeros and at the $\vec \imath_i$-th row a one. Note that $\bPhi_{\vec\imath}=\bP^\top \bPhi$ holds. To evaluate the approximate nonlinearity we need $\mathrm c$. Since we know $\bPhi$ and the index vector $\vec\imath$, we can compute
\begin{align*}
    \mathrm c=(\bP^\top\bPhi)^{-1}\bP^\top\mathrm b^h= \left(\bP^\top\bPhi\right)^{-1}\bP^\top\mathrm n^h(\bPsi\mathrm y)\quad\text{for }\mathrm y\in\mathbb R^\ell.
\end{align*}
Suppose that the matrix $\mathrm P$ can be moved into the nonlinearity. Then we obtain
\begin{align*}
    \bP^\top\mathrm n^h(\bPsi\mathrm y)=\big(\mathrm n^h(\bPsi\mathrm y)\big)_{\vec \imath}=\tilde{\mathrm n}^h(\bP^\top\bPsi\mathrm y)\quad\text{for }\mathrm y\in\mathbb R^\ell,
\end{align*}
where $\tilde{\mathrm n}^h:\mathbb R^{|\vec\imath|}\to\mathbb R^{|\vec\imath|}$ only depends on the vector $(\bPsi\mathrm y)_{\vec\imath}$. An extension for general nonlinearities is shown in \cite{CS10}. Let us now have a look at the computational expenses. The matrices  $\bP^\top\bPsi\in\mathbb R^{\mathsf p\times\ell}$, $\bP^\top\bPhi\in\mathbb R^{\mathsf p\times\mathsf p}$ and $\bPsi^\top\bPhi \in \mathbb R^{\ell\times\mathsf p}$ can be precomputed. All these are independent of the full dimension $m$. Additionally, during the iterations the nonlinearity only has to be evaluated at the interpolation points, i.e. only at $\mathsf p$ points. This allows the reduced-order model to be completely independent of the full dimension. Note that the used method is an interpolation and therefore is exact at the interpolation points. For the Jacobian the approach is similar.

Let us now turn to the EIM and DEIM algorithms. When \eqref{NonlFineModel} is solved, the nonlinearity $\mathrm n^h(\mathrm y^h(t))$ is evaluated for each time step $t_j$, $j=1,\ldots,n$. If these evaluations are stored, the procedure to determine $\bPhi$ and the index vector $\vec\imath$ does not involve any further evaluations of the nonlinearity. We denote by $\bF\in\mathbb R^{m \times n}$ the matrix with columns $\mathrm n^h(\mathrm y^h(t_j)) \in \mathbb R^m$ for $j=1,\ldots,n$. Next let us have a look at the two algorithms of interest and let us present some numerical results. In the algorithms $|\cdot|_\infty$ stands for the maximum norm in $\mathbb R^m$ and the operation `$\argmax$' returns the index, where
the maximum entry occurs. In Algorithm~\ref{algo:EIM} we state the empirical interpolation method (EIM) using a greedy algorithm.

\medskip
\hrule
\vspace{-3.5mm}
\begin{algorithm}[(Empirical interpolation method)]
    \label{algo:EIM}
    \vspace{-3mm}
    \hrule
    \vspace{0.5mm}
    \begin{algorithmic}[1]
        \REQUIRE $\mathsf p$ and matrix $\mathrm F=[\mathrm n^h(\mathrm y^h(t_1))\,|\ldots |\,\mathrm n^h(\mathrm y^h(t_n))]\in \mathbb R^{m\times n}$;
        \STATE $k\leftarrow\argmax_{j=1,\ldots,n}|f(t_j,y(t_j))|_\infty$;
        \STATE $\xi\leftarrow\mathrm n^h(\mathrm y^h(t_k))$;
        \STATE idx $\leftarrow\argmax_{j=1,\ldots,m}|\xi_j|$;
        \STATE $\phi_1\leftarrow \xi/\xi_{\{\text{idx}\}}$;
        \STATE $\bPhi=[\phi^1]$ and $\vec\imath=$ idx;
        \FOR{$i=2$ to $\ell^{EI}$}
            \STATE Solve $\bPhi_{\{\text{idx}\}}\mathrm c_j=\mathrm n^h(\mathrm y^h(t_j))_{\{\text{idx}\}}$ for $j = 1,\ldots,n$;
            \STATE $k\leftarrow\argmax_{j=1,\ldots,n}|\mathrm n^h(\mathrm y(t_j))-\bPhi\mathrm c_j|_\infty$;
            \STATE $\xi\leftarrow\mathrm n^h(\mathrm y^h(t_k))$;
            \STATE idx$\leftarrow\argmax_{j=1,\ldots,m}|(\xi - \bPhi\mathrm c_k)_{\{j\}}|$;
            \STATE $\phi^i\leftarrow(\xi-\bPhi\mathrm c_k)/(\xi-\bPhi\mathrm c_k)_{\{\text{idx}\}}$;
            \STATE $\bPhi\leftarrow [\bPhi|\phi^i]$ and $\vec\imath\leftarrow [\vec\imath\,|\,$idx$]$;
        \ENDFOR 
        \RETURN $\bPhi$ and $\vec \imath$
    \end{algorithmic}
    \hrule
\end{algorithm}

Note that the basis $\phi^i$, $i=1,\ldots,\mathsf p$, is chosen from the provided snapshots of $\mathrm n^h(\mathrm y^h(t))$ by scaling and shifting. The obtained basis is not orthonormal. The advantage of this method is that the submatrix $\bPhi_{\vec\imath}$ is an upper triangular matrix. Hence solving for $\mathrm c(t)$ is computationally cheap. The drawback of this method is that the computation of the basis is more expensive than the discrete empirical interpolation method (DEIM) presented in Algorithm~\ref{algo:DEIM}.

\medskip
\hrule
\vspace{-7mm}
\begin{algorithm}[(Discrete empirical interpolation method)]
    \label{algo:DEIM}
    \vspace{-3mm}
    \hrule
    \vspace{0.5mm}
    \begin{algorithmic}[1]
        \REQUIRE $\mathsf p$ and matrix $\bF=[\mathrm n^h(\mathrm y^h(t_1))\,|\ldots |\,\mathrm n^h(\mathrm y^h(t_n))]\in\mathbb R^{m \times n}$;
        \STATE Compute POD basis $\bPhi=[\phi_1,\ldots,\phi_\mathsf p]$ for $\bF$;
        \STATE idx$ \leftarrow \argmax_{j=1,\ldots,m}|(\phi_1)_{\{j\}}|$;
        \STATE $\bU=[\phi_1]$ and $\vec\imath = $ idx;
        \FOR{$i=2$ to $\mathsf p$}
            \STATE $u \leftarrow \phi_i$;
            \STATE Solve $\bU_{\vec \imath}\,\mathrm c =\mathrm u_{\vec \imath}$\,;
            \STATE $\mathrm r \leftarrow\mathrm u-\bU\mathrm c$;
            \STATE idx$ \leftarrow\argmax_{j=1,\ldots,m}|(\mathrm r)_{\{j\}}|$;
            \STATE $\bU\leftarrow[\bU|\mathrm u]$ and $\vec\imath\leftarrow[\vec \imath\,|\,$idx$]$;
        \ENDFOR 
        \RETURN $\bPhi$ and $\vec \imath$
    \end{algorithmic}
    \hrule
\end{algorithm}

The DEIM algorithm on the other hand generates the basis using the POD approach. Here the previously introduced POD approach is applied to the snapshots of the nonlinearity $\mathrm b^n(t)=\mathrm n^h(\mathrm y^h(t))$ to compute $\bPhi$. The matrix $\bPhi_{\vec\imath}$ obtained by the DEIM method has no special structure. Hence evaluating the nonlinearity using DEIM is more expensive compared to EIM. The computational cost can be reduced by precomputing a LU decomposition of $\bPhi_{\vec\imath}$. Then the evaluation of the nonlinearity using DEIM involves two solves compared
to one solve for the EIM. Further when comparing the two algorithms it can be seen that the computation for the EIM basis is more expensive compared to the DEIM basis. This can be seen when comparing line 7 in Algorithm~\ref{algo:EIM} and line 6 in Algorithm~\ref{algo:DEIM}. In each iteration of Algorithm~\ref{algo:EIM} one has to solve $n$ linear systems compared to one linear system in Algorithm~\ref{algo:DEIM}. The selection for the interpolation points in both algorithms is similar and is based on a greedy algorithm. The idea is to successively select spatial points to limit the growth of an error bound. The indices are constructed inductively from the input data. For more details we refer the reader to \cite{BMNP04,CS09}. 

\section{The POD method for elliptic problems}
\label{SIAM-Book:Section3.7}
\setcounter{equation}{0}
\setcounter{theorem}{0}
\setcounter{figure}{0}
\setcounter{run}{0}

In Section~\ref{SIAM:Section-2.1.2} we have introduced the continuous POD method an a closed and bounded parameter set $\mathscr D\subset\mathbb R^\mathfrak p$ with $\mathfrak p\in\mathbb N$. For the choice $\mathscr D=[0,T]$ we utilized the POD method to derive Galerkin POD schemes for evolution equations in Section~\ref{SIAM-Book:Section3.3}. In this section we choose a finite-dimensional set $\mathscr D$ containing parameters which enter parametrized elliptic variational problems. Then the computed POD bases are used to derive reduced-order parametrized systems.

Since we focus on evolution problems in this book, we only discuss briefly the application of the POD method to certain parametrized elliptic variational problems. For more details we refer to \cite{KV07,KV09,KV12,VH08}, for instance. Let us mention that especially in the context of elliptic variational problems the {\em reduced basis method} offers an efficient alternative for the derivation of low-order models; see, e.g., in the textbooks \cite{HRS16,PR05,QMN16}.

\subsection{The parametrized elliptic variational problem}
\label{SIAM-Book:Section3.7.1}

Let $\mathscr D \subset \mathbb R^\mathfrak p$ be a bounded and closed set. Suppose that for $\bmu \in D$ the \index{Bilinear form!parameter-dependent, $a(\bmu;\cdot\,,\cdot)$}{\em parameter-dependent bilinear form} $a(\bmu;\cdot\,,\cdot):V\times V \to\mathbb R$ satisfies
\begin{subequations}
    \label{BilinearParam}
    \begin{align}
        \label{BilinearParam-a}
        &\big|a(\bmu;\varphi,\phi)\big|\le\gamma\,{\|\varphi\|}_V{\|\phi\|}_V&&\text{for all }\varphi,\phi\in V\text{ and for }\bmu\in\mathscr D,\\
        \label{BilinearParam-b}
        & a(\bmu;\varphi,\varphi)\ge\gamma_1\,{\|\varphi\|}_V^2&&\text{for all }\varphi \in V\text{ and for }\bmu\in\mathscr D
    \end{align}
\end{subequations}
for positive constants $\gamma$, $\gamma_1$. Further, let $\mathcal F(\bmu) \in V'$ for $\bmu\in\mathscr D$ be a parameter-dependent right-hand side. For given $\mu \in \mathscr D$, we consider the variational problem: find $y=y(\bmu) \in V$ solving
\begin{equation}
    \label{Eq:ParVarPro}
    a(\bmu;y(\bmu),\varphi)= {\langle \mathcal F(\bmu),\varphi\rangle}_{V',V}\quad\text{for all }\varphi \in V.
\end{equation}

\begin{example}
    \label{Re:PaVarPro}
    \rm For $\bmu_\mathsf a,\bmu_\mathsf b \in \mathbb R$ with $\bmu_\mathsf a\le\bmu_\mathsf b$ we define the parameter set $\mathscr D=[\mu_\mathsf a,\mu_\mathsf b]$. Then we understand the parameter dependent bilinear form $a(\bmu;\cdot\,,\cdot):V \times V \to \mathbb R$ as
    \begin{align*}
        a(\bmu;\varphi,\phi)={\langle \varphi,\phi \rangle}_V+\bmu\,{\langle \varphi,\phi \rangle}_H\quad\text{for }\varphi,\phi \in V \text{ and }\bmu\in\mathscr D.
    \end{align*}
    For any $\bmu\in\mathscr D$ we infer from the embedding inequality \eqref{Poincare} that
    \begin{align*}
        |a(\bmu;\varphi,\phi)|\le\left(1+c_V^2 \, \max \{|\bmu_\mathsf a|,|\bmu_\mathsf b|\} \right){\| \varphi \|}_V{\| \phi \|}_V\quad\text{for all }\varphi,\phi \in V,
    \end{align*}
    i.e., the bilinear form $a(\bmu;\cdot \, , \cdot)$ satisfies \eqref{BilinearParam-a} with $\gamma=1+\max\{|\bmu_\mathsf a|,|\bmu_\mathsf b|\} c_V^2$. Further
    \begin{align*}
        a(\mu;\varphi,\varphi)={\|\varphi\|}_V^2+\mu\,{\|\varphi\|}_H^2\ge{\|\varphi\|}_V^2+\mu_a\,{\| \varphi \|}_H^2\quad\text{for all }\varphi \in V\text{ and }\mu\in\mathscr D.
    \end{align*}
    If $\bmu_\mathsf a \ge 0$ holds, then \eqref{BilinearParam-b} is satisfied with $\gamma_1=1$. In the case $\bmu_\mathsf a <0$ we infer from \eqref{Poincare} that
    \begin{align*}
        a(\bmu;\varphi,\varphi) \ge {\| \varphi \|}_V^2+\bmu_\mathsf a \, {\| \varphi \|}_H^2\ge \big(1+\bmu_\mathsf a c_V^2\big) {\| \varphi \|}_V^2\quad\text{for all } \varphi \in V \text{ and } \mu \in \mathscr D.
    \end{align*}
    Summarizing, \eqref{BilinearParam-b} holds if $\gamma_1=1+\min\{0,\bmu_\mathsf a c_V^2\}>0$ is fulfilled.\hfill$\blacklozenge$
\end{example}

The following theorem ensures that \eqref{Eq:ParVarPro} admits a unique solution. Its proof is given in Section~\ref{SIAM-Book:Section3.8.6}.

\begin{theorem}
    \label{Th:EllPro}
    Suppose that the parameter dependent bilinear form $a(\bmu;\cdot\,,\cdot)$ satisfies \eqref{BilinearParam} and $\mathcal F(\bmu) \in V'$ holds true for any $\bmu \in \mathscr D$. Then there exists a unique solution $y(\bmu) \in V$ to \eqref{Eq:ParVarPro} for every $\bmu \in \mathscr D$. Moreover, we have
    \begin{equation}
        \label{eq2-7}
        {\|y(\bmu)\|}_V\le\frac{1}{\gamma_1}\,{\| \mathcal F(\bmu) \|}_{V'}\quad\text{for every }\bmu \in \mathscr D.
    \end{equation}
    In particular, if the mapping $\bmu \mapsto \mathcal F(\bmu) \in V'$ belongs to $L^2(\mathscr D;V')$, the mapping $\bmu \mapsto y(\bmu)$ belongs to $L^2(\mathscr D;V)$.
\end{theorem}

\subsection{The continuous POD method}
\label{SIAM-Book:Section3.7.2}

Together with \eqref{Eq:ParVarPro} we will consider a discretized variational problem, where we apply POD for the discretization of $V$. For that purpose let $y(\bmu) \in V$ be the associated solution to \eqref{Eq:ParVarPro} for a chosen parameter $\bmu \in \mathscr D$. We suppose that $\mathcal F\in L^2(\mathscr D;V')$ holds, so that $y \in L^2(\mathscr D;V) \hookrightarrow L^2(\mathscr D;H)$ by Theorem~\ref{Th:EllPro}. Next we apply a \index{POD method!continuous variant}continuous variant of the POD method. Further, $X$ denotes either the space $V$ or the space $H$. We quickly summarize the results from Section \ref{SIAM:Section-2.1.2} for the current situation. Let us define the bounded linear operator $\mathcal Y:L^2(\mathscr D) \to X$
by
\begin{align*}
    \mathcal Y\phi=\int_\mathscr D\phi(\bmu)y(\bmu)\,\mathrm d\bmu\quad\text{for }\phi\in L^2(\mathscr D).
\end{align*}
Then the \index{POD method!continuous variant!snapshot space}snapshot space $\mathscr V$ is given as
\begin{align*}
    \mathscr V=\mathrm{ran}\,\mathcal Y\subset V.
\end{align*}
Let $d=\dim\mathscr V\le\infty$. The Hilbert space adjoint $\mathcal Y^\star:X \to L^2(\mathscr D)$ of $\mathcal Y$ is given by
\begin{align*}
    \big(\mathcal Y^\star \psi\big)(\bmu)={\langle\psi,y(\bmu)\rangle}_X\quad\text{for }\psi\in X\text{ and }\bmu \in \mathscr D.
\end{align*}
Furthermore, we find that the bounded, linear, symmetric and nonnegative operator $\mathcal R=\mathcal Y \mathcal Y^\star:X \to X$ has the form
\begin{equation}
    \label{eq3-4}
    \mathcal R\psi=\int_\mathscr D{\langle \psi,y(\bmu)\rangle}_Xy(\bmu)\,\mathrm d\mu\quad\text{for }\psi\in X.
\end{equation}
The operator $\mathcal K=\mathcal Y^\star \mathcal Y:L^2(\mathscr D) \to L^2(\mathscr D)$ is given by
\begin{equation}
    \label{opK}
    \big(\mathcal K\phi\big)(\bmu)=\int_\mathscr D{\langle y(\bnu),y(\bmu) \rangle}_X\phi(\bnu)\,\mathrm d\bnu\quad\text{for }\phi\in L^2(\mathscr D).
\end{equation}
Since the mapping $\bmu \mapsto y(\bmu) \in V$ is in $L^2(\mathscr D;V)$, we conclude that
\begin{align*}
    \int_\mathscr D\int_\mathscr D \big|{\langle y(\bnu),y(\bmu) \rangle}_X\big|^2\,\mathrm d\bnu\mathrm d\bmu<\infty.
\end{align*}
From \cite[pp. 197 and 277]{Yos95} it follows that $\mathcal K=\mathcal Y^\star \mathcal Y$ is compact and, therefore, $\mathcal R=\mathcal Y\mathcal Y^\star$ is
compact as well. From the Theorems~\ref{SIAM:Theorem-I.1.1.1} and \ref{SIAM:Theorem-I.1.1.2} it follows that there exists a complete orthonormal basis $\{\psi_i\}_{i \in \mathbb I}$ of $V$ and a sequence $\{\lambda_i\}_{i \in \mathbb I}$ of nonnegative real numbers so that for any $\ell\le d$
\begin{align*}
    \mathcal R\psi_i=\lambda_i\psi_i\text{ for }1\le i\le\ell,\quad\lambda_1\ge\lambda_2\ge\ldots\ge \lambda_\ell>0.
\end{align*}
Furthermore,
\begin{align*}
    \int_\mathscr D{\| y(\bmu) \|}_X^2\,\mathrm d\bmu=\sum_{i\in\mathbb I}\lambda_i.
\end{align*}

\begin{remark}[Methods of snapshots]
    \label{Re3-1}
    \rm Analogous to Remark~\ref{Remark2.2.1a}, we find that the bounded, linear, symmetric and nonnegative operator $\mathcal K$ defined in \eqref{opK} has the same eigenvalues $\{\lambda_i\}_{i \in \mathbb I}$ as the operator $\mathcal R$ and the eigenfunctions
    \begin{align*}
        \phi_i(t)=\frac{1}{\sqrt{\lambda_i}}~\big(\mathcal Y^\star\psi_i \big)(\bmu)=\frac{1}{\sqrt{\lambda_i}}~{\langle \psi_i, y(\bmu)\rangle}_V
    \end{align*}
    for $i \in \{ j \in \mathbb N:\lambda_j>0\}$ and almost all $\bmu\in\mathscr D$.\hfill$\blacksquare$
\end{remark}

In the following theorem we formulate properties of the eigenvalues and eigenfunctions of $\mathcal R$. Its proof follows from Theorem~\ref{Theorem2.2.1}.

\begin{theorem}
    \label{Th3-1}
    Let $\{\lambda_i\}_{i \in \mathbb I}$ and $\{\psi_i\}_{i \in \mathbb I}$ denote the eigenvalues and eigenfunctions, respectively, of the operator $\mathcal R$ introduced in \eqref{eq3-4}. Then for every $\ell \in \mathbb I$ the first $\ell$ eigenfunctions $\psi_1,\ldots,\psi_\ell \in X$ solve the minimization problem
    \begin{equation}
        \label{eq3-11}
        \left\{
        \begin{aligned}
            &\min\int_\mathscr D \Big\|y(\mu)-\sum_{i=1}^\ell {\langle y(\bmu),\psi_i\rangle}_X\,\psi_i\Big\|_X^2 \, \mathrm d\bmu\\
            &\text{s.t. }\{\psi_i\}_{i=1}^\ell\subset X\text{ and }{\langle\psi_j,\psi_i \rangle}_X=\delta_{ij}\text{ for }1\le i,j\le\ell.
        \end{aligned}
        \right.
    \end{equation}
    Moreover,
    \begin{equation}
        \label{eq3-12}
        \int_\mathscr D\Big\|y(\bmu)-\sum_{i=1}^\ell{\langle y(\bmu),\psi_i\rangle}_X\,\psi_i\Big\|_X^2\,\mathrm d\bmu=\sum_{i>\ell}\lambda_i\quad\text{for any }\ell\in\mathbb N.
    \end{equation}
\end{theorem}

We call a solution to \eqref{eq3-11} a {\em POD basis of rank $\ell$}. Recall that $\{\psi\}_{i=1}^\ell$ is also the solution to
\begin{align*}
    \left\{
    \begin{aligned}
        &\min-\int_\mathscr D \sum_{i=1}^\ell \big|{\langle y(\bmu),\psi_i\rangle}_X\big|^2\,\mathrm d\bmu\\
        &\hspace{0.8mm}\text{s.t. }\{\psi_i\}_{i=1}^\ell\subset X\text{ and }{\langle\psi_j,\psi_i \rangle}_X=\delta_{ij}\text{ for }1\le i,j\le\ell.
    \end{aligned}
    \right.
\end{align*}
Using \eqref{Eq2.2.18} we have that
\begin{align*}
    \sum_{i=1}^\ell\lambda_i=\sum_{i=1}^\ell \int_\mathscr D \big|{\langle y(\mu),\psi_i \rangle}_X\big|^2\,\mathrm d\bmu\ge \sum_{i=1}^\ell\int_\mathscr D\big|{\langle y(\bmu),\chi_i\rangle}_X\big|^2\,\mathrm d\bmu
\end{align*}
for every $\ell\in\mathbb N$, where $\{\chi_i\}_{i \in \mathbb N}$ is an arbitrary orthonormal basis in $X$.

\begin{remark}
    \label{Remark:EllProb}
    \rm
    \begin{enumerate}
        \item [1)] If a POD basis of rank $\ell$ is determined, we set $X^\ell=\Span\{\psi_1,\ldots,\psi_\ell\}$ and introduce the following POD Galerkin scheme for \eqref{Eq:ParVarPro} (cf. Section~\ref{SIAM-Book:Section3.3.1}): Find
        \begin{equation}
            \label{Eq:ParVarProPOD-1}
            y^\ell(\bmu)=\sum_{i=1}^\ell\mathrm y_i^\ell(\bmu)\psi_i,\quad\bmu\in\mathscr D,
        \end{equation}
        satisfying the variational problem
        \begin{equation}
            \label{Eq:ParVarProPOD-2}
            a(\bmu;y^\ell(\bmu),\psi)= {\langle \mathcal F(\bmu),\psi\rangle}_{V',V}\quad\text{for all }\psi \in X^\ell.
        \end{equation}
        Inserting \eqref{Eq:ParVarProPOD-1} into \eqref{Eq:ParVarProPOD-2} and choosing $\psi=\psi_i$, $1\le i\le\ell$ in \eqref{Eq:ParVarProPOD-2}, we derive the following linear system for the coefficient vector $\mathrm y(\bmu)=(\mathrm y_i(\bmu))_{1\le i\le \ell}\in\mathbb R^\ell$:
        \begin{align*}
            \bA(\bmu)\mathrm y(\bmu)=\mathrm \mathcal F(\bmu)
        \end{align*}
        with the stiffness matrix $\bA(\bmu)=((a(\bmu;\psi_j,\psi_i)))\in\mathbb R^{\ell\times\ell}$ and the right-hand side $\mathrm f(\bmu)=(\langle \mathcal F(\bmu),\psi_i\rangle_{V',V})\in\mathbb R^\ell$. Since $\bA(\bmu)$ is positive definite, existence of a unique solution $\mathrm y(\bmu)$ follows directly for any $\ell\le d$. This implies that also \eqref{Eq:ParVarProPOD-2} possesses a unique solution as well.
        \item [2)] Utilizing Proposition~\ref{Prop:VTopology} we can derive a-priori error estimates for the difference $y(\bmu)-y^\ell(\bmu)$; cf. Section~\ref{SIAM-Book:Section3.3.2}. Since we concentrate on evolution problems in this book, we does not study this error analysis here.
        \item [3)] Combining the continuous POD method for evolution problems (compare Section~\ref{SIAM-Book:Section3.3}) and for parametrized elliptic variational problems we can study parametrized evolution problems, where the parameter set is given as $\mathscr D_T=[0,T] \times \mathscr D$. To compute the POD bases the POD greedy method turns out to be very efficient; see \cite{HO08} and Section~\ref{SIAM:Section-PODGreedy}.\hfill$\blacksquare$
    \end{enumerate}
\end{remark}

\subsection{The discrete POD method}
\label{SIAM-Book:Section3.7.3}

In applications the weak solution to \eqref{Eq:ParVarPro} is not known for all parameters $\bmu \in \mathscr D$, but only for a given grid in $\mathscr D$. For that purpose we consider a \index{POD method!discrete variant}discrete variant of the POD method also lere. Let $\{\bmu_j\}_{j=1}^n$ be a grid in $\mathscr D$ and let $y_j=y(\bmu_j)$, $1\le j\le n$, denote the corresponding solutions to \eqref{Eq:ParVarPro} for the grid points $\bmu_j$. We define the snapshot set $\mathcal V^n=\mathrm{span}\,\{y_1,\ldots,y_n\}\subset V$ and set $d^n=\dim\mathscr V^n\le n$. Then we determine a POD basis of rank $\ell\le d^n$ for $\mathcal V^n$ by solving
\begin{equation}
    \label{eq4-8}
    \min\sum_{j=1}^n\alpha_j\,\bigg\|y_j-\sum_{i=1}^\ell{\langle y_j,\psi_i\rangle}_X\,\psi_i\bigg\|_X^2\text{ s.t. }\{\psi_i\}_{i=1}^\ell \subset X\text{ and }{\langle\psi_j,\psi_i\rangle}_X=\delta_{ij}\text{ for }1\le i,j\le \ell,
\end{equation}
where the $\alpha_j$'s are nonnegative weights. Let us introduce the operator $\mathcal R^n:X\to\mathscr V^n\subset X$ by
\begin{align*}
    \mathcal R^n\psi=\sum_{j=1}^n\alpha_j\,{\langle y_j,\psi\rangle}_X\,y_j \quad \text{for } \psi \in X.
\end{align*}
In contrast to $\mathcal R$ introduced in \eqref{eq3-4} the operator $\mathcal R^n$ depends on the grid  $\{\bmu_j\}_{j=1}^n$. The image space of $\mathcal R^n$ has finite dimension $d^n\le n$, whereas, in general, the image space of the operator $\mathcal R$ is infinite-dimensional. Since $\mathcal R^n$ is a linear, bounded, compact, nonnegative, selfadjoint operator, there exist eigenvalues $\{\lambda_i^n\}_{i\in\mathbb I}$ and orthonormal eigenfunctions $\{\psi_i^n\}_{i\in\mathbb I}$ satisfying
\begin{align*}
    \mathcal R^n\psi_i^n=\lambda_i^n\psi_i^n\text{ for }i=1,\ldots,d^n,\quad\lambda_1^n\ge\ldots\ge\ldots\ge\lambda_{d^n}^n>0.
\end{align*}
The solution to \eqref{eq4-8} is given by $\{\psi_i^n\}_{i=1}^\ell$ and we have
\begin{align*}
    \sum_{j=1}^n\alpha_j\,\bigg\|y_j-\sum_{i=1}^\ell{\langle y_j,\psi_i^n\rangle}_X\psi_i^n\bigg\|_X^2=\sum_{i=\ell+1}^{d^n}\lambda_i^n.
\end{align*}

\begin{remark}
    \label{Re-SnapshotPOD}
    \rm
    \begin{enumerate}
        \item [1)] \index{POD method!discrete variant!methods of snapshots}{\em Methods of snapshot} \cite{Sir87}: Let us introduce the $(n\times n)$ diagonal matrix $\bD=\mathrm{diag}\, (\alpha_1,\ldots,\alpha_n)$. We supply $\mathbb R^n$ with the weighted inner product
        \begin{align*}
            {\langle \mathrm u,\mathrm v\rangle}_\bD=\sum_{i=1}^n\alpha_i\mathrm u_i\mathrm v_i=\mathrm u^\top\bD\mathrm v \quad\text{for }\mathrm u=(u_i),\,\mathrm v=(v_i)\in \mathbb R^n.
        \end{align*}
        If the $\alpha_i$'s are quadrature weights corresponding to the parameter grid $\{\bmu_i \}_{i=1}^n$ then the inner product $\langle \cdot\,,\cdot\rangle_\bD$ is a discrete version of the inner product in $L^2(\mathscr D)$. We define the symmetric nonnegative matrix $\mathcal K^n\in\mathbb R^{n \times n}$ with the elements $\langle y_i,y_j\rangle_X$, $1 \le i,j\le n$, and consider the eigenvalue problem
        \begin{equation}
            \label{Snapshot-POD1}
            \mathcal K^n\phi_i^n=\lambda_i^n\phi_i^n,~1\le i\le\ell\quad\text{and}\quad {\langle\phi_i^n,\phi_j^n\rangle}_\bD=\delta_{ij},~1\le i,j\le\ell.
        \end{equation}
        From the singular value decomposition it follows that $\mathcal K^n$ has the same eigenvalues $\{\lambda_i^n\}_{i\in\mathbb I}$ as the operator $\mathcal R^n$. Furthermore, the POD basis functions are given by the formula
        \begin{equation}
            \label{Snapshot-POD2}
            \psi_i=\frac{1}{\sqrt{\lambda_i^n}}\sum_{j=1}^n\alpha_j\Phi_{ji}^n y_j\quad\text{for } i=1,\ldots,\ell\le d^n, 
        \end{equation}
        where $\Phi_{ji}^n$ denotes the $j$-th component of the eigenvector $\phi_i^n\in\mathbb R^n$.
        \item [2)] Usually, the snapshots $\{y_j\}_{j=1}^n$ has to computed by a numerical approximation scheme. Here we can introduce a Galerkin discretization as in Section~\ref{Spatial discretization}.
        \item [3)] Following Section~\ref{SIAM-Book:Section3.5.4} we can introduce a POD Galerkin scheme and discuss an a-priori error analysis utilizing Proposition~\ref{Prop:VTopologyDisc}; cf. Section~\ref{SIAM-Book:Section3.5.5}. However, we will not do it in this book, because our focus is on evolution problems; cf. Remark~\ref{Remark:EllProb}-2).\hfill$\blacksquare$
    \end{enumerate}
\end{remark}

\subsection{Extension to nonlinear problems}
\label{SIAM-Book:Section3.7.4}

Let us give a few comments how to deal with parametrized nonlinear elliptic variational problems. Suppose that for any $\bmu\in\mathscr D$ the mapping $\mathcal N(\bmu;\cdot):V \to V'$ is a nonlinear, locally Lipschitz-continuous mapping satisfying
\begin{equation}
    \label{MonN}
    {\langle \mathcal N(\bmu;\phi)-\mathcal N(\bmu;\varphi),\phi-\varphi\rangle}_{V',V} \ge 0 \quad \text{for all } \phi,\varphi \in V \text{ and for all }\bmu \in \mathscr D,
\end{equation}
i.e., $\mathcal N(\bmu;\cdot)$ is monotone for any $\bmu \in \mathscr D$. Instead of \eqref{Eq:ParVarPro} we consider
\begin{equation}
    \label{Eq:ParVarProNonl}
    a(\bmu;y(\bmu),\varphi) + {\langle \mathcal N(\bmu;y(\bmu)),\varphi \rangle}_{V',V}={\langle\mathcal f(\bmu),\varphi\rangle}_{V',V} \quad \text{for all }\varphi \in V.
\end{equation}

\begin{example}
    \rm Let us give an example for a semilinear problem satisfying \eqref{MonN}. Suppose that $\Omega \subset \mathbb R^\mathfrak n$, $\mathfrak n \in \{1,2,3\}$, is a bounded and open set with Lipschitz-continuous boundary $\Gamma=\partial\Omega$. We consider the \index{Equation!partial differential!elliptic}{\em elliptic partial differential equation}
    \begin{equation}
        \label{GinzEq}
        -\mu_1\Delta y+\mu_2y+\mu_3y^3=f\text{ in }\Omega\quad\text{and}\quad\mu_1\frac{\partial y}{\partial\bn}+\mu_4y=g \text{ on } \Gamma,
    \end{equation}
    where $f \in L^2(\Omega)$, $g\in L^2(\Gamma)$ and
    \begin{align*}
        \mathscr D=\big\{\bmu=(\mu_1,\ldots,\mu_4)\in\mathbb R^4\, \big|\,\bmu_\mathsf a \le\bmu\le\bmu_\mathsf b\text{ in } \mathbb R^4\big\}
    \end{align*}
    with $0<\mu_a\le\mu_b$. A {\em weak solution} to \eqref{GinzEq} satisfies $y \in V=H^1(\Omega)$ and
    \begin{equation}
        \label{GinzEqVar}
        \int_\Omega\mu_1\nabla y\cdot\nabla\varphi+\big(\mu_2y+\mu_3y^3\big)\varphi\,\mathrm d\bx+\int_\Gamma\mu_4y\varphi \, \mathrm d\bs=\int_\Omega f\varphi\,\mathrm d\bx +\int_\Gamma g\varphi\,\mathrm d\bs
    \end{equation}
    for all $\varphi \in V$. Next we express \eqref{GinzEqVar} in the form \eqref{Eq:ParVarProNonl}. For that purpose we utilize the parametrized bilinear form $a(\mu;\cdot\,,\cdot):V \times V \to \mathbb R$ given by
    \begin{align*}
        a(\bmu;\phi,\varphi)=\int_\Omega\mu_1\nabla\phi\cdot\nabla\varphi+\mu_2\phi\varphi\,\mathrm d\bx+\int_\Gamma\mu_4\phi\varphi\, \mathrm d\bs\quad\text{for }\phi,\varphi\in V\text{ and }\bmu\in\mathscr D.
    \end{align*}
    Notice that this bilinear form satisfies \eqref{BilinearParam}. Moreover, let the parameter independent right-hand side be given as
    \begin{align*}
        {\langle\mathcal F,\varphi \rangle}_{V',V}=\int_\Omega f\varphi\,\mathrm d\bx+\int_\Gamma g\varphi\,\mathrm d\bs\quad\text{for }\varphi\in V.
    \end{align*}
    Finally, we define the nonlinearity
    \begin{align*}
        {\langle \mathcal N(\bmu;\varphi),\phi \rangle}_{V',V}=\mu_3\int_\Omega\varphi^3\phi \, \mathrm d\bx\quad \text{for }\varphi,\phi \in V \text{ and }\bmu\in\mathscr D.
    \end{align*}
    Then a weak solution to \eqref{GinzEq} satisfies the variational formulation \eqref{Eq:ParVarProNonl}. Recall that $\varphi \in V$ implies $\varphi \in L^6(\Omega)$. Consequently, $\varphi^3 \in H=L^2(\Omega) \subset V'$. Let $\phi,\varphi \in V$ and $\chi=\phi-\varphi \in V$. From $\mu_3 \ge \mu_a >0$ we infer that
    \begin{align*}
        {\langle \mathcal N(\bmu;\phi)-\mathcal N(\bmu;\varphi),\chi \rangle}_{V',V}&=\int_\Omega \mu_3 \big(\phi^3-\varphi^3\big)\chi\, \mathrm d\bx=\mu_3\int_\Omega \bigg(\int_0^1 3(\varphi+\tau\chi)^2\chi\mathrm d\tau\big)\chi\,\mathrm d\bx\\
       &=\mu_3\int_0^1 \int_\Omega (\varphi+\tau\chi)^2\chi^2 d\tau\,\mathrm d\bx\ge 0
    \end{align*}
    holds true. Thus, \eqref{MonN} is satisfied.\hfill$\blacklozenge$
\end{example}

\begin{remark}
    \rm If a solution $y(\bmu)$ to \eqref{Eq:ParVarProNonl} is given, then a POD basis can be computed as described above for the linear elliptic variational problem \eqref{Eq:ParVarPro}.\hfill$\blacksquare$
\end{remark}

\section{Proofs of Section~\ref{SIAM-Book:Section3}}
\label{SIAM-Book:Section3.8}
\setcounter{equation}{0}
\setcounter{theorem}{0}
\setcounter{figure}{0}
\setcounter{run}{0}

\subsection{Proofs of Section~\ref{SIAM-Book:Section3.2}}
\label{SIAM-Book:Section3.8.1}

\noindent{\bf\em Proof of Theorem~{\em\ref{SIAM:Theorem3.1.1}}.}
For a proof of the existence of a unique weak solution we refer to \cite[pp.~512-520]{DL00}. The regularity result follows from \cite[pp.~532-533]{DL00} and \cite[pp.~360-364]{Eva08}. To prove the a-priori error estimate \eqref{SIAM:Eq3.1.7} we choose $\varphi=y(t)$ in \eqref{SIAM:Eq3.1.6a}. Utilizing the coercivity estimate \eqref{SIAM:Eq3.1.1-2} and Young's inequality we derive
\begin{align*}
    \frac{1}{2}\,\frac{\mathrm d}{\mathrm dt}\,{\|y(t)\|}_H^2+\gamma_1\,{\|y(t)\|}_V^2 - \gamma_2 \|y(t)\|_H^2&\le\big({\|\mathcal F(t)\|}_{V'}+{\|(\mathcal Bu)(t)\|}_{V'}\big){\|y(t)\|}_V\\
    &\le\frac{1}{2\gamma_1}\big({\|\mathcal F(t)\|}_{V'}^2+{\|(\mathcal Bu)(t)\|}_{V'}^2\big)+\frac{\gamma_1}{2}\,{\|y(t)\|}_V^2,
\end{align*}
which implies
\begin{equation}
    \label{EstAprioriState-1}
    \frac{\mathrm d}{\mathrm dt}\,{\|y(t)\|}_H^2+\gamma_1\,{\|y(t)\|}_V^2\le 2\gamma_2\|y(t)\|^2_H + \frac{2}{\gamma_1}\big({\|\mathcal F(t)\|}_{V'}^2+{\|(\mathcal Bu)(t)\|}_{V'}^2\big).
\end{equation}
We insert the embedding estimate \eqref{Poincare} and infer that
\begin{align*}
\frac{\mathrm d}{\mathrm dt}\,{\|y(t)\|}_H^2\le c_1\,{\|y(t)\|}_H^2+\frac{2}{\gamma_1}\big({\|\mathcal F(t)\|}_{V'}^2+{\|(\mathcal Bu)(t)\|}_{V'}^2\big)
\end{align*}
with $c_1=\max(2\gamma_2 - \gamma_1/c_V^2,0) \ge 0$. By applying the Gronwall Lemma (Proposition~\ref{Gronwall}), we find
\begin{align}
    \label{EstAPrioriState-1b}
    \begin{aligned}
        {\|y(t)\|}_H^2&\le e^{c_1t}\bigg({\|y_\circ\|}_H^2+\int_0^t\frac{2}{\gamma_1}\big({\|\mathcal F(s)\|}_{V'}^2+{\|(\mathcal Bu)(s)\|}_{V'}^2\big)\,\mathrm ds\bigg)\\
        &\le e^{c_1t}\bigg({\|y_\circ\|}_H^2+\frac{2}{\gamma_1}\left({\|\mathcal F\|}_{L^2(0,T;V')}^2+{\|\mathcal B\|}_{\mathscr L(U;L^2(0,T;V'))}^2{\|u\|}_\U^2\right)\bigg)
    \end{aligned}
\end{align}
which implies
\begin{equation}
\label{EstAprioriState-2}
{\|y(t)\|}_H^2\le e^{c_1t}\Big({\|y_\circ\|}_H^2+c_2\big({\|\mathcal F\|}_{L^2(0,T;V')}^2+{\|u\|}_\U^2\big)\Big)
\end{equation}
with $c_2=2\max(1,\|\mathcal B\|_{\mathscr L(U;L^2(0,T;V'))}^2)/\gamma_1>0$. Integrating \eqref{EstAprioriState-1} over $[0,T]$ and dropping the nonnegative term $\|y(T)\|_H^2$, it follows that
\begin{align*}
    {\|y\|}_{L^2(0,T;V)}^2&~\le\frac{1}{\gamma_1}\bigg({\|y_\circ\|}_H^2+c_2\,\big({\|\mathcal F\|}_{L^2(0,T;V')}^2+{\|u\|}_\U^2\big) + 2\gamma_2 \int_0^T \|y(t)\|_H^2 ~\mathrm dt \bigg)
\end{align*}
We insert the estimate \eqref{EstAPrioriState-1b}, which results in
\begin{equation}
    \label{EstAprioriState-3}
    \begin{aligned}
        {\|y\|}_{L^2(0,T;V)}^2 \le c_3\big({\|y_\circ\|}_H^2+{\|\mathcal F\|}_{L^2(0,T;V')}^2+{\|u\|}_\U^2\big),
    \end{aligned}
\end{equation}
where the constant is given by
\begin{align*}
	c_3 = \frac{(1+2\gamma_2 T e^{c_1T})\max(1,c_2)}{\gamma_1}.
\end{align*}
Let $\varphi \in V$ with $\|\varphi\|_V = 1$ be arbitrary. Using \eqref{SIAM:Eq3.1.6a} and \eqref{SIAM:Eq3.1.1-1}, we estimate
\begin{align*}
    {\langle y_t(t),\varphi\rangle}_{V',V} & = {\langle (\mathcal F+\mathcal Bu)(t),\varphi\rangle}_{V',V}-a(t;y(t),\varphi)\\
    & \le {\|(\mathcal F+\mathcal Bu)(t)\|}_{V'}{\|\varphi\|}_V+\gamma\,{\|y(t)\|}_V{\|\varphi\|}_V\\
    & \le {\|\mathcal F(t)\|}_{V'}+{\|(\mathcal B u)(t)\|}_{V'} +\gamma\,{\|y(t)\|}_{V}
\end{align*}
for almost all $t \in (0,T)$. Thus, it holds
\begin{align*}
    {\| y_t(t) \|}_{V'} \leq {\|\mathcal F(t)\|}_{V'}+{\|(\mathcal B u)(t)\|}_{V'} +\gamma\,{\|y(t)\|}_{V}
\end{align*}
for almost all $t \in (0,T)$. Consequently, using \eqref{EstAprioriState-3} yields
\begin{equation}
    \label{EstAprioriState-4}
    \begin{aligned}
        {\|y_t\|}_{L^2(0,T;V')}^2 & = \int_0^T {\| y_t(t) \|}_{V'}^2\,\mathrm dt\le \int_0^T 3 {\|\mathcal F(t)\|}_{V'}^2 + 3\,{\|(\mathcal B u)(t)\|}_{V'}^2 + 3\gamma\,{\|y(t)\|}_{V}^2 \, \mathrm dt\\
        & = 3\,{\|\mathcal F\|}_{L^2(0,T;V')}^2 + 3\,{\| \mathcal B u \|}_{L^2(0,T;V')}^2 + 3\gamma\,{\|y\|}_{L^2(0,T;V)}^2 \\
        & \leq 3\,{\|\mathcal F\|}_{L^2(0,T;V')}^2 + 3\,{\|\mathcal B\|}_{\mathscr L(\U;L^2(0,T;V'))}^2{\|u\|}_\U^2 + 3\gamma\,{\|y\|}_{L^2(0,T;V)}^2 \\
        & \leq c_4\big({\|y_\circ\|}_H^2+{\|\mathcal F\|}_{L^2(0,T;V')}^2+{\|u\|}_\U^2\big)
    \end{aligned} 
\end{equation}
with $c_4= 3 \gamma c_3 + 3\max(1,\|\mathcal B\|_{\mathscr L(U;L^2(0,T;V'))}^2)$. From \eqref{EstAprioriState-3} and \eqref{EstAprioriState-4} we get
\begin{align*}
    {\|y\|}_\Y&=\big({\|y_t\|}_{L^2(0,T;V')}^2+{\|y\|}_{L^2(0,T;V)}^2\big)^{1/2}\\
    &\le \left(c_4+c_3\right)^{1/2} \big({\|y_\circ\|}_H^2+{\|\mathcal F\|}_{L^2(0,T;V')}^2+{\|u\|}_\U^2\big)^{1/2} \\
    & \leq \left(c_4+c_3\right)^{1/2} \big({\|y_\circ\|}_H +{\|\mathcal F\|}_{L^2(0,T;V')} +{\|u\|}_\U\big).
\end{align*}
Setting $C=\left(c_4+c_3\right)^{1/2}$ we obtain the a-priori estimate \eqref{SIAM:Eq3.1.7}.\hfill$\Box$

\subsection{Proofs of Section~\ref{SIAM-Book:Section3.3}}
\label{SIAM-Book:Section3.8.2}

\noindent{\bf\em Proof of Lemma~{\em\ref{MApdef}}.} Clearly, $\bM^\ell$ is symmetric. Let $\mathrm v=(\mathrm v_i)^\top\in\mathbb R^\ell$ be chosen arbitrarily. We set
\begin{align*}
    v^\ell=\sum_{i=1}^\ell\mathrm v_i\psi_i\in X^\ell\subset V.
\end{align*}
Then we derive
\begin{align*}
    \mathrm v^\top\bM^\ell\mathrm v=\sum_{i=1}^\ell\sum_{j=1}^\ell\mathrm v_i\mathrm v_j\,{\langle\psi_j,\psi_i\rangle}_H={\|v^\ell\|}_H^2\ge 0.
\end{align*}
Further, $\mathrm v^\top\bM^\ell\mathrm v=0$ holds if and only if $\mathrm v=0$ is true because $\{\psi_1,\hdots,\psi_\ell\}$ are linearily independent. Thus, the matrix $\bM^\ell$ is positive definite. From \eqref{SIAM:Eq3.1.1-2} and the assumption that $\gamma_2=0$, we derive
\begin{align*}
    \mathrm v^\top\bA^\ell(t)\mathrm v=\sum_{i=1}^\ell\sum_{j=1}^\ell\mathrm v_i\mathrm v_ja(\psi_j,\psi_i)\ge\gamma_1\,{\|v^\ell\|}_V^2\ge0
\end{align*}
which implies \eqref{EqvAv}. From \eqref{EqvAv} we infer that $\bA^\ell(t)\mathrm v\neq0$ holds for any $\mathrm v\in\mathbb R^\ell \setminus \{0\}$ in $[0,T]$. Thus, the columns of $\bA^\ell(t)$ are linearly independent, so that $\bA^\ell(t)$ is a regular matrix in $[0,T]$.\hfill$\Box$

\medskip\noindent{\bf\em Proof of Theorem~{\em\ref{SIAM:Theorem3.1.1POD}}.}
Since $V^\ell$ is a finite-dimensional subspace of $V$, the proof of the existence of a unique solution $y\in H^1(0,T;V^\ell)\subset H^1(0,T;V)$ to \eqref{SIAM:Eq3.1.6POD} follows again from \cite[pp.~512-520]{DL00}. To prove the a-priori estimate we proceed as in the proof of Theorem~\ref{SIAM:Theorem3.1.1} and derive
\begin{equation}
    \label{EstAprioriState-3-1}
    {\|y^\ell\|}_{L^2(0,T;V)}\le c_1\big({\|\mathcal P^\ell y_\circ\|}_H+{\|\mathcal F\|}_{L^2(0,T;V')}+{\|u\|}_\U\big);
\end{equation}
where $c_1=\max(1/\sqrt\gamma_1,2/\gamma_1,2\|\mathcal B\|_{\mathscr L(\U;L^2(0,T;V'))}/\gamma_1)>0$. From \eqref{SIAM:Eq3.1.1-1}, \eqref{SIAM:Eq3.1.6PODa} and $\mathcal F,\mathcal Bu\in L^2(0,T;V')$  we infer that
\begin{equation}
    \label{yellVprime}
    y_t^\ell(t)={\langle (\mathcal F+\mathcal Bu)(t),\cdot\rangle}_{V',V}-a(y^\ell(t),\cdot)\text{ in }(V^\ell)'
\end{equation}
holds for almost all $t \in [0,T]$. Hence, \eqref{yellVprime} and \eqref{SIAM:Eq3.1.1-1} imply that
\begin{align*}
    {\langle y_t^\ell(t) , \varphi \rangle_{V',V}} & = {\langle (\mathcal F+\mathcal Bu)(t),\varphi\rangle}_{V',V}-a(y^\ell(t),\varphi) \\
    & \leq {\| \mathcal F(t) \|}_{V'} + {\| (\mathcal Bu)(t) \|}_{V'} + \gamma \, {\|y^\ell(t)\|}_V
\end{align*}
holds for all $\varphi \in V^\ell$ with ${\| \varphi \|}_V = 1$ and almost all $t \in [0,T]$. Using \eqref{eq:VellIsometry} this implies
\begin{equation}
    \label{EstAprioriState-3-1a}
    \begin{aligned}
    {\| y_t^\ell(t) \|}_{V'} = \| y_t^\ell(t) \|_{(V^\ell)'} \leq {\| \mathcal F(t) \|}_{V'} + {\| (\mathcal Bu)(t) \|}_{V'} + \gamma \, {\|y^\ell(t)\|}_V.
    \end{aligned}
\end{equation}
Hence, using \eqref{EstAprioriState-3-1} and \eqref{EstAprioriState-3-1a} yields
\begin{equation}
    \label{AppEstimate}
    \begin{aligned}
        {\|y_t^\ell\|}_{L^2(0,T;V')}^2 &= \int_0^T {\| y_t^\ell(t) \|}_{V'}^2 \,\mathrm dt \\
        &\leq \int_0^T 3\,{\| \mathcal F(t) \|}_{V'}^2 + 3\,{\| (\mathcal Bu)(t) \|}_{V'}^2 + 3\gamma \, {\|y^\ell(t)\|}_V^2 \,\mathrm dt \\
        &\leq 3\,{\|\mathcal F\|}_{L^2(0,T;V')}^2 + 3\,{\|\mathcal B\|}_{\mathscr L(\U,L^2(0,T;V'))}^2 {\|u\|}_\U^2 + 3 \gamma \, {\|y^\ell\|}_{L^2(0,T;V)}^2 \\
        &\leq c_2\big({\|\mathcal P^\ell y_\circ\|}_H^2 + {\|\mathcal F\|}_{L^2(0,T;V')}^2+{\|u\|}_\U^2\big),
    \end{aligned}
\end{equation}
where $c_2=9\gamma\,c_1+3\max(1,\|\mathcal B\|_{\mathscr L(\U,L^2(0,T;V'))}^2)$. Now, we directly get
\begin{align*}
    {\|y^\ell\|}_{W(0,T)}\le C\left({\|\mathcal P^\ell y_\circ\|}_H+{\|\mathcal F\|}_{L^2(0,T;V')}+{\|u\|}_\U\right)
\end{align*}
follows directly with $C = (3 c_1^2 + c_2)^{1/2}$. This implies \eqref{SIAM:Eq3.1.7POD} due to the isomorphism between $H^1(0,T;V)$ and $\Y$, cf. Lemma \ref{lem:W0T_finDim}. \hfill$\Box$

\medskip\noindent{\bf\em Proof of Lemma~{\em\ref{lem:W0T_finDim}}.} Let $y \in W(0,T)$ be given. Since $V$ and $H$ have the same finite dimension $\mathfrak n$, $V'$ has this dimension as well and the linear and injective embeddings $V \hookrightarrow H \hookrightarrow V'$ are thus surjective as well. For almost every $t \in (0,T)$, it holds $y_t(t) \in V'$ and, due to the surjectivity, there exists $\tilde y(t) \in V$ whose Gelfand embedding into $V'$ is given by $\langle \tilde{y}(t),\cdot \rangle_H$ and is identical to $y_t(t)$. Because of $y \in H^1(0,T;V')\hookrightarrow W(0,T)$, it holds by definition of this Sobolev space:
\begin{align*}
    -\int_0^T \langle y(t),\cdot \rangle_H \xi'(t) ~\mathrm dt &= \int_0^T y_t(t) \xi(t) ~\mathrm dt = \int_0^T \langle \tilde{y}(t),\cdot \rangle_H \xi(t) ~\mathrm dt \quad \text{in } V',~ \text{for all } \xi \in C_0^\infty(0,T).
\end{align*}
Alternatively, we can write this as
\begin{align*}
    \bigg\langle-\int_0^Ty(t)\xi'(t)\,\mathrm dt,\varphi\bigg\rangle_H=\bigg\langle\int_0^T \tilde y(t)\xi(t)\,\mathrm dt,\varphi\bigg\rangle_H \quad\text{for all } \xi \in C_0^\infty(0,T),\,\varphi \in V.
\end{align*}
Since $V \hookrightarrow H$ is dense (even surjective) and injective, this yields
\begin{equation}
    \label{eq:appFinDim_1}
    -\int_0^T y(t) \xi'(t) ~\mathrm dt = \int_0^T \tilde{y}(t) \xi(t)\,\mathrm dt\quad\text{in }V\text{ for all }\xi\in C_0^\infty(0,T).
\end{equation}
This proves that the $V$-derivative of $y$ is given by $\tilde y$. After we have shown that $\tilde y \in L^2(0,T;V)$, this means $y \in H^1(0,T;V)$ with $y_t = \tilde y$. First of all, the inverse of the surjective embedding $V \hookrightarrow H$ is also continuous as a linear mapping between finite-dimensional spaces. The same holds for the inverse of $H \hookrightarrow V'$. In other words, there is $C_\mathfrak n>0$ such that the {\em inverse inequality}
\begin{equation}
    \label{InvIneq}
    {\|\varphi\|}_V\le C_\mathfrak n\,{\| \varphi \|}_H\quad\text{for all }\varphi\in H
\end{equation}
and $\| \varphi \|_H \le C_\mathfrak n\,\| \langle \cdot, \varphi \rangle_H \|_{V'}$ for all $\varphi \in H$. Notice that \eqref{InvIneq} does not hold if $V$ and $H$ have infinite dimensions, which implies that the constant $C_\mathfrak n$ depends on the dimension $\mathfrak n$ of $H$ and $V$. In particular, $\lim_{\mathfrak n\to\infty}C_\mathfrak n=\infty$. Therefore, 
\begin{equation}
    \label{eq:appFinDim_2}
    {\| \tilde y(t)\|}_V\le C_\mathfrak n\,{\|\tilde y(t)\|}_H\le C_\mathfrak n^2\,{\|{\langle \cdot\,,\tilde y(t)\rangle}_H\|}_{V'}=C_\mathfrak n^2\,{\|y_t(t)\|}_{V'}.
\end{equation}
From $y_t \in L^2(0,T;V')$, it follows $\tilde y \in L^2(0,T;V)$. We have therefore defined a mapping 
\begin{align*}
    W(0,T) \to H^1(0,T;V), \quad y \mapsto y
\end{align*}
which is obviously linear and, because of \eqref{eq:appFinDim_2}, continuous. Because of the nature of the mapping and the canonical injection $H^1(0,T;V) \hookrightarrow W(0,T)$, it is also bijective, making it an isomorphism between Hilbert spaces.\hfill$\Box$

\medskip\noindent{\bf\em Proof of Lemma~{\em\ref{Lemma:HI-22}}.} From \eqref{SIAM:Eq3.1.6}, \eqref{SIAM:Eq3.1.6POD}, \eqref{SIAM:Eq3.1.1-1} and $\varrho^\ell=\mathcal P^\ell y-y$ we get for any $\psi\in X^\ell$
\begin{equation}
    \label{APriori-Est-1c}
    \begin{aligned}
        &{\langle \vartheta_t^\ell(t),\psi\rangle}_{V',V}+a(t;\vartheta^\ell(t),\psi)={\langle \vartheta_t^\ell(t),\psi\rangle}_{V',V}+a(t;\vartheta^\ell(t),\psi)\\
        &\quad={\langle\mathcal P^\ell y_t(t),\psi\rangle}_{V',V}+a(t;\mathcal P^\ell y(t),\psi)-{\langle y_t^\ell(t),\psi\rangle}_{V',V}-a(t;y^\ell(t),\psi)\\
        &\quad={\langle(\mathcal P^\ell y_t-y_t)(t),\psi\rangle}_{V',V}+a(t;\mathcal P^\ell y(t)-y(t),\psi)\\
        &\quad\le\big({\|\varrho^\ell_t(t)\|}_{V'}+\gamma\,{\|(\varrho^\ell(t)\|}_V\big){\|\psi\|}_V
    \end{aligned}
\end{equation}
Choosing $\psi=\vartheta^\ell(t)\in X^\ell$ for almost all $t \in [0,T]$ applying Lemma~\ref{Lemma:HI-20}-4) and \eqref{SIAM:Eq3.1.1-2} we find
\begin{align*}
    &\frac{1}{2}\,\frac{\mathrm d}{\mathrm dt}\,{\|\vartheta^\ell(t)\|}_H^2+\gamma_1\,{\|\vartheta^\ell(t)\|}_V^2-\gamma_2\|\vartheta^\ell(t)\|_H^2\le\big({\|\varrho^\ell_t(t)\|}_{V'}+\gamma\,{\|\varrho^\ell(t)\|}_V\big){\|\vartheta^\ell(t)\|}_V
\end{align*}
in $[0,T]$ almost everywhere. Applying two times Young's inequality (cf. Lemma~\ref{lem:youngsInequality}) with $\varepsilon=\gamma_1/4$ and multiplying by two we get
\begin{equation}
    \label{APriori-Est-3}
    \frac{\mathrm d}{\mathrm dt}\,{\|\vartheta^\ell(t)\|}_H^2+\gamma_1\,{\|\vartheta^\ell(t)\|}_V^2\le\,2\gamma_2\,{\|\vartheta^\ell(t)\|}_H^2+\frac{2}{\gamma_1}\big({\|\varrho^\ell_t(t)\|}_{V'}^2+\gamma^2\,{\|\varrho^\ell(t)\|}_V^2\big)
\end{equation}
In a next step, we use the embedding estimate \eqref{Poincare} for the second summand on the left-hand side and set
\begin{align*}
    c_1=\max\Big(2\gamma_2-\frac{\gamma_1}{c_V^2},~0 \Big),
\end{align*}
which gives us
\begin{align*}
    \frac{\mathrm d}{\mathrm dt}\,{\|\vartheta^\ell(t)\|}_H^2\le c_1\,{\|\vartheta^\ell(t)\|}_H^2+c_2\big({\|\varrho^\ell_t(t)\|}_{V'}^2+{\|\varrho^\ell(t)\|}_V^2\big)
\end{align*}
for $c_2=2\max\{1,\gamma^2\}/\gamma_1>0$. Now we apply the Gronwall lemma (Proposition~\ref{Gronwall}) and use $\vartheta^\ell(0)=0$ to get
\begin{equation}
    \label{APriori-Est-4}
    \big\|\vartheta^\ell(t)\big\|^2_H\le c_3\int_0^t{\|\varrho^\ell_t(s)\|}_{V'}^2+{\|\varrho^\ell(s)\|}_V^2\,\mathrm ds
    \le c_3\,{\|\varrho^\ell\|}_\mathscr Y^2\quad\text{a.e. in }[0,T]
\end{equation}
with $c_3=c_2\exp(c_1T)>0$. From \eqref{APriori-Est-4} it follows that
\begin{equation}
    \label{ThetaW(0,T)-1}
    {\|\vartheta^\ell\|}_{C([0,T];H)}^2\le c_3\,{\|\varrho^\ell\|}_\mathscr Y^2.
\end{equation}
We continue by integrating \eqref{APriori-Est-3} over $[0,T]$ dropping the non-negative term $\|\vartheta^\ell(T)\|_H^2$ and using again $\vartheta^\ell(0)=0$ in $H$. This gives us
\begin{align*}
    \gamma_1\,{\|\vartheta^\ell\|}_{L^2(0,T;V)}^2 \le2\gamma_2\int_0^T \|\vartheta^\ell(t)\|_H^2\,\mathrm dt+ c_2\,{\|\varrho^\ell\|}_\mathscr Y^2.
\end{align*} 
Inserting the estimate \eqref{ThetaW(0,T)-1} for $\|\vartheta^\ell(t)\|_H^2$ we derive
\begin{equation}
    \label{ThetaW(0,T)-2}
    {\|\vartheta^\ell\|}_{L^2(0,T;V)}^2\le c_4\,{\|\varrho^\ell\|}_\mathscr Y^2
\end{equation}
for the constant $c_4=(2\gamma_2c_3T+c_2)/\gamma_1>0$. Using \eqref{SIAM:Eq3.1.1-1} and \eqref{EqDec} we find that
\begin{align*}
    {\langle \vartheta_t^\ell(t),\psi\rangle}_{V',V}& ={\langle \mathcal P^\ell y_t(t),\psi\rangle}_{V',V}-{\langle y^\ell_t(t),\psi\rangle}_{V',V}\\
    &= {\langle\mathcal P^\ell y_t(t)-y_t(t),\psi\rangle}_{V',V}-a(t;y(t)-y^\ell(t),\psi)\\
    &\le {\|\mathcal P^\ell y_t(t)-y_t(t)\|}_{V'}{\|\psi\|}_V+\gamma\,{\|\varrho^\ell(t)+\vartheta^\ell(t)\|}_V{\|\psi\|}_V \\
    &\le c_5\left({\|\varrho_t^\ell(t)\|}_{V'}+{\|\varrho^\ell(t)\|}_{V}+{\|\vartheta^\ell(t)\|}_{V}\right)
\end{align*}
with $c_5=\max\{1,\gamma\}>0$. Hence, applying \eqref{eq:VellIsometry} yields
\begin{equation*}
    {\| \vartheta_t^\ell(t) \|}_{V'} = {\| \vartheta_t^\ell(t) \|}_{(V^\ell)'} \leq c_5 \left({\|\varrho_t^\ell(t)\|}_{V'}+{\|\varrho^\ell(t)\|}_{V}+{\|\vartheta^\ell(t)\|}_{V}\right).
\end{equation*}
Thus, \eqref{ThetaW(0,T)-2} implies
\begin{align*}
    {\| \vartheta_t^\ell \|}_{L^2(0,T;V')}^2  = \int_0^T {\| \vartheta_t^\ell(t) \|}_{V'}^2 \,\mathrm dt\leq 3 c_5^2 \big({\|\varrho_t^\ell\|}_\mathscr Y^2+{\|\vartheta^\ell\|}_{L^2(0,T,V)}^2\big)\le c_6\,{\|\varrho^\ell\|}_\mathscr Y
\end{align*}
for $c_6=3c_5^2\max\{1,c_4\}>0$. Together with \eqref{ThetaW(0,T)-2} this implies
\begin{equation}
    \label{APriori-Est-10}
    {\|\vartheta^\ell\|}_{\Y}^2={\|\vartheta^\ell\|}_{L^2(0,T;V)}^2+{\|\vartheta^\ell_t\|}_{L^2(0,T;V')}^2\le (c_4+c_6)\,{\|\varrho^\ell\|}_\mathscr Y.
\end{equation}
Thus, \eqref{APriori-Est-5} follows from \eqref{APriori-Est-10} with $C=c_4+c_6$.\hfill$\Box$

\medskip\noindent{\bf\em Proof of Lemma~{\em\ref{APriori-Est-inV}}.} We know that $\mathcal P^\ell$ is an $H$-orthonormal projection in $X=H$ or $X=V$. For arbitary $\psi, \phi \in H$, we set $\psi^H = \mathcal P^\ell\psi$, $\phi^H=\mathcal P^\ell\phi\in X^\ell$ and $\psi^\bot=\psi-\psi^H$, $\phi^\bot=\phi-\phi^H$. It follows  that $\langle\psi^H,\phi^\bot\rangle_H=0$ and $\langle\psi^\bot,\phi^H\rangle_H=0$ holds true. Hence,
\begin{equation}
    \label{Eq:SnapOhneDQ-1}
    \begin{aligned}
        {\langle\mathcal P^\ell\psi,\phi\rangle}_H={\langle\psi^H,\phi^H+\phi^\bot\rangle}_H&={\langle\psi^H,\phi^H\rangle}_H+{\langle\psi^H,\phi^\bot\rangle}_H\\
        &={\langle\psi^H,\phi^H\rangle}_H+{\langle\psi^\bot,\phi^H\rangle}_H={\langle\psi,\mathcal P^\ell\phi\rangle}_H.
    \end{aligned}
\end{equation}
Following \cite{Sin14} and using $y\in H^1(0,T;H)$ we observe that 
\begin{align*}
    {\langle \mathcal P^\ell y_t(t),\varphi\rangle}_H=\frac{\mathrm d}{\mathrm dt}\,{\langle \mathcal P^\ell y(t),\varphi\rangle}_H=\frac{\mathrm d}{\mathrm dt}{\langle y(t),\mathcal P^\ell\varphi\rangle}_H={\langle y_t(t),\mathcal P^\ell\varphi\rangle}_H
\end{align*}
for all $\varphi\in H$ and for almost all $t \in [0,T]$. In particular, we find that
\begin{equation}
    \label{Eq:SnapOhneDQ-2}
    {\langle \mathcal P^\ell y_t(t),\psi\rangle}_H={\langle y_t(t),\psi\rangle}_H\quad\text{for all }\psi\in X^\ell\text{ and in }[0,T]\text{ a.e.}
\end{equation}
Thus, we can stray from the proof of Lemma \ref{Lemma:HI-22} by making a different estimate than \eqref{APriori-Est-1c} in this case:
\begin{equation}
    \label{Eq:SnapOhneDQ-3a}
    \begin{aligned}
        {\langle \vartheta_t^\ell(t),\psi\rangle}_H+a(t;\vartheta^\ell(t),\psi)&={\langle y_t(t),\psi\rangle}_H+a(t;\mathcal P^\ell y(t),\psi)-{\langle y_t^\ell(t),\psi\rangle}_H-a(t;y^\ell(t),\psi)\\
        &=a(t;(\mathcal P^\ell y-y)(t),\psi)\le\gamma\,{\|(\mathcal P^\ell y-y)(t)\|}_V{\|\psi\|}_V
    \end{aligned}
\end{equation}
for all $\psi\in X^\ell$ and in $[0,T]$ almost everywhere. Note that no time derivatives of $\mathcal P^\ell y - y$ occur in this case, in contrast to \eqref{APriori-Est-1c}. The same computations as in the proof of Lemma \ref{Lemma:HI-22} and Remark~\ref{RemarkThetaEstimates}-2) lead to \eqref{APriori-Est-inV-1}.\hfill$\Box$

\medskip\noindent{\bf\em Proof of Theorem~{\em\ref{Theorem:OperatorConv}}.} Define the set $\mathscr S=\{u\in\U\,|\,\|u\|_\U=1\}$. Note that
\begin{align*}
    {\|\mathcal S^\ell-\mathcal S\|}_{\mathscr L(\U,\Y)}=\sup_{u\in\mathscr S}{\|\mathcal S^\ell u-\mathcal S u\|}_\Y.
\end{align*}
Consequently, $\mathcal S^\ell$ converges to $\mathcal S$ in $\mathscr L(\U,\Y)$ if and only if $\mathcal S^\ell u$ converges to $\mathcal S u$ uniformly on $\mathscr S$. \\ 
Suppose that $\{\mathcal S^\ell u\}_{\ell\in\mathbb N}$ does not converge uniformly for all $u \in \mathscr S$. Then there is an $\varepsilon>0$ such that for any $\ell_\circ=\ell_\circ(\varepsilon)\in\mathbb N$ there is an $\ell\ge\ell_\circ$ and a $u^\ell\in\mathscr S$ satisfying
\begin{equation*}
    {\|\mathcal S^\ell u^\ell -\mathcal S u^\ell \|}_\Y\ge 3\varepsilon,
\end{equation*}
i.e., we can choose a sequence $\{u^n\}_{n \in \mathbb{N}} \subset \mathscr S$ and a subsequence $\{\mathcal S^{\ell_n}\}_{n \in \mathbb{N}}$ of $\{\mathcal S^{\ell}\}_{\ell \in \mathbb{N}}$ such that 
\begin{equation}
    \label{Th:SellConv-1}
    {\|\mathcal S^{\ell_n} u^n -\mathcal S u^n \|}_\Y\ge 3\varepsilon
\end{equation}
holds for all $n \in \mathbb N$. Since $\U$ is finite-dimensional the set $\mathscr S$ is compact. Thus, without loss of generality we can assume that there is a limit point $u_\circ \in \mathscr S$ with $u^n \to u_\circ$ in $\U$ as $n \to \infty$. \\
The continuity of the operator $\mathcal{S}$ immediately implies that there is a $n_1 \in \mathbb{N}$ such that
\begin{equation}
    \label{Th:SellConv-2}
    {\| \mathcal S u^n - \mathcal S u_\circ \|}_{\Y} < \varepsilon \quad \text{for all } n \geq n_1.
\end{equation}
Since $\mathcal S^\ell$ is linear and bounded independently of $\ell$, there exists a constant $C>0$ such that
\begin{align*}
    {\|\mathcal S^\ell u_1-\mathcal S^\ell u_2\|}_\Y\le C\,{\|u_1-u_2\|}_\U
\end{align*}
for all $\ell \in \mathbb{N}$ and for all $u_1,u_2\in\U$. Setting $\delta=\varepsilon/C>0$ it holds
\begin{equation*}
    {\|\mathcal S^\ell u_1 -\mathcal S^\ell u_2 \|}_\Y < \varepsilon
\end{equation*}
for all $\ell \in \mathbb{N}$ and for all $u_1,u_2\in\U$ with ${\|u_1-u_2\|}_\U<\delta$.
Using again that $u^n \to u_\circ$ in $\U$ as $n \to \infty$ we get the existence of a $n_2 \in \mathbb{N}$ such that
\begin{align*}
    {\| u^n - u_\circ \|}_\U \leq \frac{\varepsilon}{C} \quad \text{for all } n \geq n_2.
\end{align*}
In total, we thus get
\begin{equation}
    \label{Th:SellConv-3}
    {\|\mathcal S^{\ell_n} u^n -\mathcal S^{\ell_n} u_\circ \|}_\Y < \varepsilon \quad \text{for all } n \geq n_2.
\end{equation}
Moreover, applying Theoremy~\ref{Corollary:APrioriConv} yields the existence of a $n_3 \in \mathbb{N}$ such that
\begin{equation}
    \label{Th:SellConv-4}
    {\| \mathcal S^{\ell_n} u_\circ - \mathcal S u_\circ \|}_{\Y} < \varepsilon \quad \text{for all } n \geq n_3.
\end{equation}
Therefore, we derive from \eqref{Th:SellConv-2}, \eqref{Th:SellConv-3} and \eqref{Th:SellConv-4} that
\begin{align*}
    {\|\mathcal S^{\ell_n} u^n -\mathcal S u^n \|}_\Y & \le {\|\mathcal S^{\ell_n} u^n -\mathcal S^{\ell_n} u_\circ \|}_\Y + {\|\mathcal S^{\ell_n} u_\circ -\mathcal  S u_\circ \|}_\Y + {\|\mathcal S u_\circ - \mathcal S u^n \|}_\Y \\
    & \le \varepsilon + \varepsilon + \varepsilon = 3 \varepsilon
\end{align*}
which contradicts \eqref{Th:SellConv-1}\hfill$\Box$

\subsection{Proofs of Section~\ref{SIAM-Book:Section3.4}}
\label{SIAM-Book:Section3.8.3}

\noindent{\bf\em Proof of Theorem~{\em\ref{TheoremApostSemi}}.} We set $e^{h\ell}(t;u)=y^h(t;u)-y^{h\ell}(t;u)\in V^h$ for the error between $y^h$ and $y^{h\ell}$. Then we derive from \eqref{EvProGal} and \eqref{EvProGal-POD} that for $t\in[0,T]$ and $\varphi\in V^h$
\begin{align*}
    {\langle e^{h\ell}_t(t;u),\varphi\rangle}_{V',V}+a(t;e^{h\ell}(t;u),\varphi)=-{\langle r^{h\ell}(t;u),\varphi\rangle}_{V',V}
\end{align*}
with the residual
\begin{align*}
    r^{h\ell}(t;u)={\langle y^{h\ell}_t(t;u)-(\mathcal F+\mathcal Bu)(t),\cdot\rangle}_{V',V}+a(t;y^{h\ell}(t;u),\cdot)\in (V^h)'.
\end{align*}
Choosing $\varphi=e^{h\ell}(t;u)\in V^h$ we derive from \eqref{SIAM:Eq3.1.1-2} that
\begin{align*}
    &\frac{1}{2}\,\frac{\mathrm d}{\mathrm dt}\,{\|e^{h\ell}(t;u)\|}_H^2+\gamma_1\,{\|e^{h\ell}(t;u)\|}_V^2 - \gamma_2 \|e^{h\ell}(t;u)\|_H^2 \\
&\quad\le{\langle e^{h\ell}_t(t;u),e^{h\ell}(t;u)\rangle}_{V',V}+a(t;e^{h\ell}(t;u),e^{h\ell}(t;u)) \\
&\quad=-{\langle r^{h\ell}(t;u),e^{h,\ell}(t;u) \rangle}_{V',V}
\end{align*}
Applying the Young inequality we have
\begin{align*}
    &\frac{1}{2}\,\frac{\mathrm d}{\mathrm dt}\,{\|e^{h\ell}(t;u)\|}_H^2+\gamma_1\,{\|e^{h\ell}(t;u)\|}_V^2 - \gamma_2 {\|e^{h\ell}(t;u)\|}_H^2 \\
    &\quad\le{\|r^{h\ell}(t;u)\|}_{(V^h)'}{\|e^{h\ell}(t;u)\|}_V\\
    &\quad\le\frac{1}{2\gamma_1}\,{\|r^{h\ell}(t;u)\|}_{(V^h)'}^2+\frac{\gamma_1}{2}\,{\|e^{h\ell}(t;u)\|}_V^2,
\end{align*}
which implies
\begin{equation}
    \label{ApostErr}
    \frac{\mathrm d}{\mathrm dt}\,{\|e^{h\ell}(t;u)\|}_H^2+\gamma_1\,{\|e^{h\ell}(t;u)\|}_V^2-2\gamma_2{\|e^{h\ell}(t;u)\|}_H^2 \le \frac{1}{\gamma_1}\,{\|r^{h\ell}(t;u)\|}_{(V^h)'}^2.
\end{equation}
From the embedding inequality \eqref{Poincare}, this yields
\begin{align}
	\frac{\mathrm d}{\mathrm dt}\,{\|e^{h\ell}(t;u)\|}_H^2 \le c \|e^{h\ell}(t;u)\|_H^2 + \frac{1}{\gamma_1}\,{\|r^{h\ell}(t;u)\|}_{(V^h)'}^2.
\end{align}
with the nonnegative constant $c = \max(2\gamma_2-\gamma_1/c_V^2,0)$. From applying the Gronwall Lemma (Proposition~\ref{Gronwall}), it follows that
\begin{align}
	\label{eq:ApostErr_diffH}
    {\|e^{h\ell}(t;u)\|}_H^2\le e^{ct}\bigg({\|e^{h\ell}(0;u)\|}_H^2+\frac{1}{\gamma_1}\int_0^t{\|r^{h\ell}(s;u)\|}_{(V^h)'}^2\,\mathrm ds\bigg)
\end{align}
almost everywhere in $[0,T]$. Integrating \eqref{ApostErr} over $[0,T]$ we find that
\begin{align*}
    &\int_0^T{\|e^{h\ell}(t;u)\|}_V^2\,\mathrm dt \\
    &\quad\le\frac{1}{\gamma_1}\bigg( \|e^{h\ell}(0;u)\|_H^2 + \frac{1}{\gamma_1} \int_0^T \|r^{h\ell}(t;u)\|_{(V^h)'}^2 ~\mathrm dt + 2 \gamma_2 \int_0^T \|e^{h\ell}(t;u)\|_H^2 ~\mathrm dt \bigg) \\
    &\quad\le\frac{1+2\gamma_2 T e^{ct}}{\gamma_1}\bigg({\|e^{h\ell}(0;u)\|}_H^2+\frac{1}{\gamma_1}\int_0^T{\|r^{h\ell}(t;u)\|}_{(V^h)'}^2\,\mathrm dt\bigg),
\end{align*}
where we have again used the estimate \eqref{eq:ApostErr_diffH} in the integral after $\gamma_2$. Inserting $e^{h\ell}=y^h(\cdot\,;u)-y^{h\ell}(\cdot\,;u)$, using the initial conditions for $y^h(0)$ and $y^{h\ell}(0)$ we get the claim.\hfill$\Box$

\subsection{Proofs of Section~\ref{SIAM-Book:Section3.5}}
\label{SIAM-Book:Section3.8.4}

\noindent{\bf\em Proof of Theorem~{\em\ref{FullDiscFEModel}}.} Clearly, $y_1^h$ is given by the projected initial condition $\mathcal P^h y_\circ\in V^h$. Furthermore, $y_j^h\in V^h$ is is given for $j=2,\hdots,n$ by \eqref{ThetaM:GalAnsatz} with the coordinate vectors $\mathrm y_{j\ell}^h$ solving \eqref{FineModel-Disc-1}. The matrix $\bM^h+\delta t_j\bA^h(t_j)$ is regular for sufficiently small $\Delta t$. Consequently, there exists a unique $y_j^h$. Next we derive the a-priori error bound for the solution to \eqref{EvProGal-disc}. We choose $\varphi^h=y_j^h\in V^h$. It follows from Lemma~\ref{app_lem_innerProd1} (with $\varphi=y_{j-1}^h$ and $\phi=y_j^h$) that
\begin{align*}
    {\langle y_{j-1}^h-y_j^h,y_j^h\rangle}_H=\frac{1}{2}\,{\|y_{j-1}^h\|}_H^2-\frac{1}{2}\,{\|y_j^h-y_{j-1}^h\|}_H^2-\frac{1}{2}\,{\|y_j^h\|}_H^2
\end{align*}
which implies
\begin{align*}
    {\langle y_j^h-y_{j-1}^h,y_j^h\rangle}_H=\frac{1}{2}\,{\|y_j^h\|}_H^2+\frac{1}{2}\,{\|y_j^h-y_{j-1}^h\|}_H^2-\frac{1}{2}\,{\|y_{j-1}^h\|}_H^2.
\end{align*}
Plugging the last identity into \eqref{EvProGal-disc-1} and using \eqref{SIAM:Eq3.1.1-2} yields
\begin{align*}
    &\frac{1}{2\delta t_j}\,\big({\|y_j^h\|}_H^2-{\|y_{j-1}^h\|}_H^2\big)+\frac{1}{2\delta t_j}\,{\|y_j^h-y_{j-1}^h\|}_H^2+\gamma_1\,{\|y_j^h\|}_V^2\\
    &\quad\le{\|g_j(u)\|}_{V'}{\|y_j^h\|}_V\le\frac{1}{2\gamma_1}\,{\|g_j(u)\|}^2_{V'}+\frac{\gamma_1}{2}\,{\|y_j^h\|}_V^2.
\end{align*}
Multiplication by $2\delta t_j$ gives
\begin{equation}
    \label{AprEst100}
    {\|y_j^h\|}_H^2-{\|y_{j-1}^h\|}_H^2+{\|y_j^h-y_{j-1}^h\|}_H^2+\gamma_1\delta t_j\,{\|y_j^h\|}_V^2\le\frac{\delta t_j}{\gamma_1}\,{\|g_j(u)\|}_{V'}^2
\end{equation}
for $j=2,\ldots,n$. From \eqref{Poincare} and $\delta t_{j}\ge\delta t>0$ it follows that
\begin{equation}
    \label{AprEst100ef}
    \big(1+c_1\delta t\big)\,{\| y_j^h\|}_H^2 \leq {\|y_j^h\|}_H^2 + \gamma_1\delta t_j\,{\|y_j^h\|}_V^2
\end{equation}
with $c_1=\gamma_1/c_V^2$. Now, combining \eqref{AprEst100ef} and \eqref{AprEst100} yields after summation over $j=2,\ldots,n$
\begin{equation}
    \label{AprEst100a}
    \begin{aligned}
        {\|y_j^h\|}_H^2&\le\frac{1}{1+c_1\delta t}\,{\|y_{j-1}^h\|}_H^2+\frac{\delta t_j}{\gamma_1(1+c_1\delta t_j)}\,{\|g_j(u)\|}_{V'}^2\\
        &\le\bigg(\frac{1}{1+c_1\delta t}\bigg)^{j-1}\bigg({\|y_1^h\|}_H^2+\frac{1}{\gamma_1}\sum_{l=2}^j\delta t_l\,{\|g_l(u)\|}_{V'}^2\bigg).
    \end{aligned}
\end{equation}
Suppose that
\begin{equation}
    \label{DeltaTCond}
    \delta t\le\frac{1}{2c_1}.
\end{equation}
Then we infer $1+c_1\delta t\le1/2$ and
\begin{align*}
    \frac{1}{1+c_1\delta t}=1-\frac{c_1\delta t}{1+c_1\delta t}\le1-2c_1\delta t.
\end{align*}
Due to \eqref{DeltaTCond} we also have $-2c_1\delta t>-1$ such that
\begin{equation}
    \label{DeltaTCond-2}
    \bigg(\frac{1}{1+c_1\delta t}\bigg)^{j-1}\le\big(1-2c_1\delta t\big)^{j-1}\le e^{-c_2(j-1)\delta t}
\end{equation}
with $c_2=2c_1$. Thus, \eqref{EvProGal-disc-2}, \eqref{AprEst100a}, \eqref{DeltaTCond} and \eqref{DeltaTCond-2} imply
\begin{equation}
    \label{AprEst100b}
    \begin{aligned}
        {\|y_j^h\|}_H^2&\le e^{-c_2(j-1)\delta t}\bigg({\|y_1^h\|}_H^2+\frac{1}{\gamma_1}\sum_{l=2}^j\delta t_l\,{\|g_l(u)\|}_{V'}^2\bigg)\\
        &\le e^{-c_2(j-1)\delta t}\bigg({\|\mathcal P^hy_\circ\|}_H^2+\frac{1}{\gamma_1}\sum_{l=2}^j\delta t_l\,{\|g_l(u)\|}_{V'}^2\bigg)
    \end{aligned}
\end{equation}
for $j=2,\ldots,n$, which gives \eqref{TempDiscFESOL-APriori} for $j=2,\ldots,n$. Summation upon $j$ in \eqref{AprEst100} yields \eqref{TempDiscFESOL-APriori-2}.\hfill$\Box$

\medskip\noindent{\bf\em Proof of Theorem~{\em\ref{Theorem:FEAprioriError}}.} Let $y$ and $\{y^h_j\}_{j=1}^n$ be the solutions to \eqref{SIAM:Eq3.1.6} and \eqref{EvProGal-disc}, respectively. It was assumed that both $y(t_j)$ and $y_t(t_j)$ belong to $\mathscr W$ for $j=1,\ldots,n$. To estimate the error we proceed analogously as in the proof of Theorem~\ref{Theorem:FE-AprioriError}. Therefore, we introduce the decomposition
\begin{align*}
    y(t_j)-y^h_j=y(t_j)-\mathcal P^hy(t_j)+\mathcal P^hy(t_j)-y_j^h=\varrho_j^h+\vartheta_j^h\quad\text{for }j=1,\ldots,n,
\end{align*}
where we have set $\varrho_j^h=y(t_j)-\mathcal P^hy(t_j)\in V$ and $\vartheta_j^h=\mathcal P^hy(t_j)-y_j^h\in V^h$. From \eqref{Eq:ApproxProp-2} it follows that there exist two constants $C_1,C_2>0$ satisfying
\begin{equation}
    \label{Chopin-30}
    {\|\varrho_j^h\|}_H^2\le C_1h^4\,{\|y(t_j)\|}^2_{\mathscr W}\quad\text{and}\quad{\|\varrho_j^h\|}^2_V\le C_2h^2\,{\|y(t_j)\|}^2_{\mathscr W}
\end{equation}
for $j=1,\ldots,n$. Moreover, \eqref{Eq:ApproxProp-2} implies that
\begin{equation}
    \label{Tosca-1}
    {\|y_t(t_j)-\mathcal P^hy_t(t_j)\|}_H\le C_1h^4\,{\|y_t(t_j)\|}^2_{\mathscr W}
\end{equation}
for $j=1,\ldots,n$. We use the notations $\overline\partial\vartheta_j^h=(\vartheta_j^h-\vartheta_{j-1}^h)/\delta t_j$ and $\overline\partial y(t_j)=(y(t_j)-y(t_{j-1}))/\delta t_j$ for $j=2,\ldots,n$. From \eqref{SIAM:Eq3.1.6}, \eqref{EvProGal-disc} and \eqref{SIAM:Eq3.1.1-1} we derive for every $\varphi^h\in V$
\begin{equation}
    \label{Chopin-31}
    \begin{aligned}
        {\langle\overline\partial\vartheta_j^h,\varphi^h\rangle}_H+a(t_j;\vartheta_j^h,\varphi^h)&={\langle\mathcal P^h\overline\partial y(t_j)-y_t(t_j),\varphi^h\rangle}_H+a(t_j;\mathcal P^hy(t_j)-y(t_j),\varphi^h)\\
        &\le{\|\mathcal P^h\overline\partial y(t_j)-y_t(t_j)\|}_H{\|\varphi^h\|}_H+\gamma\,{\|\varrho_j^h\|}_V{\|\varphi^h\|}_V.
    \end{aligned}
\end{equation}
Choosing $\varphi^h=\vartheta_j^h\in V^h$, using \eqref{SIAM:Eq3.1.1-2}, Young's inequality and Lemma~\ref{app_lem_innerProd1} we infer from \eqref{Chopin-31} that
\begin{equation}
    \label{Chopin-32}
    {\|\vartheta_j^h\|}_H^2+\gamma_1\delta t_j\,{\|\vartheta_j^h\|}_V^2\le{\|\vartheta_{j-1}^h\|}_H^2+\delta t_j \left( \frac{2c_V^2}{\gamma_1}\,{\|\mathcal P^h\overline\partial y(t_j)-y_t(t_j)\|}_H^2 + \frac{2\gamma^2}{\gamma_1} \,{\|\varrho_j^h\|}_V^2 \right)
\end{equation}
for $j=2,\ldots,n$, where the constant $c_V$ has been introduced in \eqref{Poincare}. Note that
\begin{equation}
    \label{Tosca-2}
    {\|y_t(t_j)-\overline\partial y(t_j)\|}_H\le\frac{1}{\delta t_j}\bigg\|\int_{t_{j-1}}^{t_j}(t-t_{j-1})y_{tt}(t)\,\mathrm dt\bigg\|_H\le\sqrt{\frac{\delta t_j}{3}}{\|y_{tt}\|}_{L^2(t_{j-1},t_j;H)}.
\end{equation}
From \eqref{Tosca-1} and \eqref{Tosca-2} we have
\begin{align*}
    {\|\mathcal P^h\overline\partial y(t_j)-\overline\partial y(t_j)\|}_H^2&\le3\,{\|\mathcal P^h\overline\partial y(t_j)-\mathcal P^hy_t(t_j)\|}_H^2+3\,{\|\mathcal P^hy_t(t_j)-y_t(t_j)\|}_H^2+3\,{\|y_t(t_j)-\overline\partial y(t_j)\|}_H^2\\
    &\le\big(1+{\|\mathcal P^h\|}_{\mathscr L(X,H)}^2\big)\delta t_j{\|y_{tt}\|}_{L^2(t_{j-1},t_j;H)}^2+3\,C_1h^4\,{\|y_t(t_j)\|}^2_{\mathscr W}
\end{align*}
for $j=2,\ldots,n$. Hence, there is a constant $C_3>0$ such that
\begin{equation}
    \label{Tosca-3}
    {\|\mathcal P^h\overline\partial y(t_j)-\overline\partial y(t_j)\|}_H^2\le C_3\,\big(\delta t_j\,{\|y_{tt}\|}_{L^2(t_{j-1},t_j;H)}^2+h^4\,{\|y_t(t_j)\|}^2_{\mathscr W}\big)
\end{equation}
for $j=2,\ldots,n$. Thus, combining \eqref{Chopin-30}, \eqref{Chopin-32}, \eqref{Tosca-2} and \eqref{Tosca-3} yields
\begin{equation}
    \label{Tosca-3ab}
    \begin{aligned}
        &6{\|\vartheta_j^h\|}_H^2+\gamma_1\delta t_j\,{\|\vartheta_j^h\|}_V^2\\
        &\quad\le{\|\vartheta_{j-1}^h\|}_H^2 + C_4\delta t_j\,\big(\delta t_j\,{\|y_{tt}\|}_{L^2(t_{j-1},t_j;H)}^2+h^4\,{\|y_t(t_j)\|}^2_{\mathscr W}+h^2\,{\|y(t_j)\|}^2_{\mathscr W}\big)
    \end{aligned}
\end{equation}
for $j=2,\ldots,n$ and a suitable constant $C_4>0$. Further, $\vartheta_1^h=0$ holds. From numerical integration it follows that
\begin{align*}
    \sum_{i=1}^n\delta t_i\,\big(h^2\,{\|y_t(t_i)\|}^2_{\mathscr W}+{\|y(t_j)\|}^2_{\mathscr W}\big)\le h^2\,{\|y_t\|}_{L^2(0,T;\mathscr W)}^2+{\|y\|}_{L^2(0,T;\mathscr W)}^2+C_5\Delta t^2
\end{align*}
for a constant $C_5>0$ provided $y_{tt}\in L^2(0,T;\mathscr W)$ holds. Summation of \eqref{Tosca-3ab} over $i=2,\ldots,j$ for $j\in\{2,\ldots,n\}$ implies that
\begin{align*}
    {\|\vartheta_j^h\|}_H^2 + \gamma_1 \, \sum_{i=1}^j \delta t_i {\|\vartheta_j^h\|}_V^2&\le C_4\,\bigg(\Delta t^2\,{\|y_{tt}\|}_{L^2(0,T;H)}^2+h^2\big(h^2\,{\|y\|}_{L^2(0,T;\mathscr W)}^2+{\|y_t\|}_{L^2(0,T;\mathscr W)}^2+C_5\Delta t^2\big)\bigg)\\
    &\le C_6\,\big(\Delta t^2+h^2+\Delta t^2h^2\big)
\end{align*}
for $j=2,\ldots,n$ and for a constant $C_6>0$ depending on $y,y_t$ and $y_{tt}$. Thus, we have on the one hand
\begin{align*}
    \sum_{j=1}^n\delta t_j\,{\|y(t_j)-y^h_j\|}_H^2&\le2\sum_{j=1}^n\delta t_j\,{\|\varrho_j^h\|}_H^2+2\sum_{j=1}^n\delta t_j\,{\|\vartheta_j^h\|}_H^2\le C_7\,\big(\Delta t^2+h^2+\Delta t^2h^2\big)
\end{align*}
for a constant $C_7>0$ and on the other hand with the same argument
\begin{align*}
    \sum_{j=1}^n\delta t_j\,{\|y(t_j)-y^h_j\|}_V^2\le C_{8}\,\big(\Delta t^2+h^2+\Delta t^2h^2\big).
\end{align*}
for a constant $C_8>0$.\hfill$\Box$

\medskip\noindent{\bf\em Proof of Lemma~{\em\ref{Lemma:FullDiscTheta}}.} We use the notation $\overline\partial\vartheta_j^{h\ell}=(\vartheta_j^{h\ell}-\vartheta_{j-1}^{h\ell})/\delta t_j$ for $2\le j\le n$. From \eqref{EvProGal-disc}, \eqref{EvProGalPOD-disc}, \eqref{SIAM:Eq3.1.1-1} and $\varrho_j^{h\ell}=y^h_j-\mathcal P^{h\ell}y^h_j$ for $j=1,\ldots,n$, we obtain for every $\psi\in V$
\begin{equation}
    \label{Eq:FDPOD-1}
    \begin{aligned}
        {\langle\overline\partial\vartheta_j^{h\ell},\psi\rangle}_H+a(t_j;\vartheta_j^{h\ell},\psi)
        &={\langle\overline\partial \mathcal P^{h\ell}y^h_j,\psi\rangle}_H+a(t_j;\mathcal P^{h\ell}y^h_j,\psi)-{\langle\overline\partial y^{h\ell}_j,\psi\rangle}_H-a(t_j;y^{h\ell}_j,\psi)\\
        &={\langle\mathcal P^{h\ell}\overline\partial y^h_j-\overline\partial y^h_j,\psi\rangle}_H+a(t_j;\mathcal P^{h\ell}y^h_j-y^h_j,\psi)\\
        &\le{\|\overline\partial \varrho_j^{h\ell}\|}_H{\|\psi\|}_H+\gamma\,{\|\varrho_j^{h\ell}\|}_V{\|\psi\|}_V\quad\text{for }2\le j\le n.
    \end{aligned}
\end{equation}
Let us choose $\psi=\vartheta_j^{h\ell}$ in \eqref{Eq:FDPOD-1}. Utilizing \eqref{SIAM:Eq3.1.1-2}, Lemma~\ref{app_lem_innerProd1} and the Young inequality we get
\begin{equation}
    \label{Eq:FDPOD-2}
    \begin{aligned}
        &{\|\vartheta_j^{h\ell}\|}_H^2-{\|\vartheta_{j-1}^{h\ell}\|}_H^2+{\|\vartheta_j^{h\ell}-\vartheta_{j-1}^{h\ell}\|}_H^2+2\gamma_1\delta t_j\,{\|\vartheta_j^{h\ell}\|}_V^2\\
        &\le\delta t_j\bigg(\frac{2c_V^2}{\gamma_1}\,{\|\overline\partial\varrho_j^{h\ell}\|}_H^2+\frac{\gamma_1}{2c_V^2}\,{\|\vartheta_j^{h\ell}\|}_H^2+\frac{\gamma_1}{2}\,{\|\vartheta_j^{h\ell}\|}_V^2+\frac{2\gamma^2}{\gamma_1}\,{\|\varrho_j^{h\ell}\|}_V^2\bigg).
    \end{aligned}
\end{equation}
From \eqref{Poincare} and \eqref{Eq:FDPOD-2} we find
\begin{align*}
    {\|\vartheta_j^{h\ell}\|}_H^2\le{\|\vartheta_{j-1}^{h\ell}\|}_H^2+c_1\delta t_j\big({\|\varrho_j^{h\ell}\|}_V^2+{\|\overline\partial\varrho_j^{h\ell}\|}_V^2\big)
\end{align*}
with $c_1=2\max(c_V^4,\gamma^2)/\gamma_1$. Recall that $\vartheta_1^{h\ell}=0$. By summation on $j$ we find
\begin{align*}
    {\|\vartheta_j^{h\ell}\|}_H^2\le c_1\sum_{l=2}^j\delta t_l\left({\|\varrho_l^{h\ell}\|}_V^2+{\|\overline\partial\varrho_l^{h\ell}\|}_V^2\right)\le c_1\sum_{l=1}^n\delta t_l\left({\|\varrho_l^{h\ell}\|}_V^2+{\|\overline\partial\varrho_l^{h\ell}\|}_V^2\right).
\end{align*}
Using $\delta t_l\le2\alpha_l^n$ and \eqref{Chopin-5} we have
\begin{align*}
    {\|\vartheta_j^{h\ell}\|}_H^2\le 2c_1\sum_{l=1}^n\alpha_l^n \left({\|\varrho_l^{h\ell}\|}_V^2+{\|\overline\partial\varrho_l^{h\ell}\|}_V^2\right)=2c_1\sum_{k=1}^2\sum_{l=1}^n\alpha_l^n\,{\|y_l^k-\mathcal P^{h\ell}y_l^k\|}_V^2.
\end{align*}
Recall that $\sum_{j=1}^n\alpha_j=T$. Consequently,
\begin{align*}
    \sum_{j=1}^n\alpha_j^n\,{\|\vartheta_j^{h\ell}\|}_H^2\le 2c_1T\sum_{k=1}^2\sum_{j=1}^n\alpha_j^n\,{\|y_j^k-\mathcal P^{h\ell}y_j^k\|}_V^2.
\end{align*}
Setting $c_2=2(c_V^2+2c_1T)$ and using \eqref{Poincare}, \eqref{Chopin-5} we find that
\begin{align*}
    &\sum_{j=1}^n\alpha_j^n\,{\|y^h_j-y^{h\ell}_j\|}_H^2\le 2\sum_{j=1}^n\alpha_j^n\,{\|\varrho^{h\ell}_j\|}_H^2+2\sum_{j=1}^n\alpha_j^n\,{\|\vartheta^{h\ell}_j\|}_H^2\\
    &\quad\le 2c_V^2\sum_{j=1}^n\alpha_j^n\,{\|\varrho^{h\ell}_j\|}_V^2+4c_1T\sum_{k=1}^2\sum_{j=1}^n\alpha_j^n\,{\|y_j^k-\mathcal P^{h\ell}y_j^k\|}_V^2\le c_2\sum_{k=1}^2\sum_{j=1}^n\alpha_j^n\,{\|y_j^k-\mathcal P^{h\ell}y_j^k\|}_V^2.
\end{align*}
From \eqref{Chopin-4}, \eqref{Eq:FDPOD-2} and $\vartheta_1^{h\ell}=0$ we also derive
\begin{align*}
    &\sum_{j=1}^n\alpha_j^n\,{\|\vartheta_j^{h\ell}\|}_V^2\le \zeta \sum_{j=1}^n\delta t_j\,{\|\vartheta_j^{h\ell}\|}_V^2= \zeta \sum_{j=2}^n\delta t_j\,{\|\vartheta_j^{h\ell}\|}_V^2\\
    &\quad\le \frac{c_1\zeta}{\gamma_1}\sum_{j=2}^n\delta t_j\,\left({\|\overline\partial\varrho_j^{h\ell}\|}_V^2+{\|\varrho_j^{h\ell}\|}_V^2\right) \leq C\sum_{k=1}^2\sum_{j=1}^n\alpha_j^n\,{\|y_j^k-\mathcal P^{h\ell}y_j^k\|}_V^2
\end{align*}
which gives \eqref{Chopin-6} with $C= 2 c_1 \zeta/\gamma_1$.\hfill$\Box$

\medskip\noindent{\bf\em Proof of Theorem~{\em\ref{Theorem:NoDQSnaps}}.} Let $\psi\in X^{n\ell}$. Then $\mathcal P^{n\ell}\psi=\psi$ holds and using \eqref{Eq:SnapOhneDQ-1} we find
\begin{equation}
    \label{Bach-100}
    {\langle \mathcal P^{n\ell}\overline\partial y_j^h,\psi\rangle}_H=\frac{1}{\delta t_j}\,{\langle\mathcal P^{n\ell}(y^h_j-y_{j-1}^h),\psi\rangle}_H=\frac{1}{\delta t_j}\,{\langle y_j^h-y_{j-1}^h,\mathcal P^{n\ell}\psi\rangle}_H={\langle\overline\partial y_j^h,\psi\rangle}_H 
\end{equation}
for $j=2,\ldots,n$. Thus, we can modify the estimation of $\vartheta_j^{h\ell}$. Instead of \eqref{Eq:FDPOD-1} we find that for any $\psi\in X^{n\ell}$
\begin{align*}
    {\langle\overline\partial\vartheta_j^{h\ell},\psi\rangle}_H+a(t_j;\vartheta_j^{h\ell},\psi)&={\langle\mathcal P^{n\ell}\overline\partial y^h_j-\overline\partial y^h_j,\psi\rangle}_H+a(t_j;\mathcal P^{n\ell}y^h_j-y^{h\ell}_j,\psi)\\
    &=a(t_j;\mathcal P^{n\ell}y^h_j-y^{h\ell}_j,\psi)\le\gamma\,{\|\varrho_j^{h\ell}\|}_V{\|\psi\|}_V
\end{align*}
for $j=2,\ldots,n$. Thus, we derive instead of \eqref{Eq:FDPOD-2}
\begin{align*}
    {\|\vartheta_j^{h\ell}\|}_H^2-{\|\vartheta_{j-1}^{h\ell}\|}_H^2+{\|\vartheta_j^{h\ell}-\vartheta_{j-1}^{h\ell}\|}_H^2+2\gamma_1\delta t_j\,{\|\vartheta_j^{h\ell}\|}_V^2\le\delta t_j\bigg(\gamma_1\,{\|\vartheta_j^{h\ell}\|}_V^2+\frac{\gamma^2}{\gamma_1}\,{\|\varrho_j^{h\ell}\|}_V^2\bigg)
\end{align*}
for $j=2,\ldots,n$. Consequently,
\begin{align*}
    \sum_{j=1}^n\alpha_j^n\,{\|\vartheta_j^{h\ell}\|}_H^2\le \frac{\gamma^2T}{\gamma_1} \sum_{j=1}^n\alpha_j^n\,{\|y_j^1-\mathcal P^{n\ell}y_j^1\|}_V^2
\end{align*}
and
\begin{align*}
    \sum_{j=1}^n\alpha_j^n\,{\|\vartheta_j^{h\ell}\|}_V^2\le \zeta\sum_{j=2}^n\delta t_j\,{\|\vartheta_j^{h\ell}\|}_V^2\le c_4 \sum_{j=1}^n\alpha_j^n\,{\|y_j^1-\mathcal P^{n\ell}y_j^1\|}_V^2
\end{align*}
with $c_4=2 \zeta \gamma^2/\gamma_1^2$. Now we argue as in the proof of Theorem~\ref{Th:A-PrioriError-200}.\hfill$\Box$

\medskip\noindent{\bf\em Proof of Theorem~{\em\ref{Prop:ApostiError}}.} We set $e_j^{h\ell}=y^h_j-y^{h\ell}_j\in V^h$ for $j=1,\ldots,n$ and $\overline\partial e_j^{h\ell}=(e_j^{h\ell}-e_{j-1}^{h\ell})/\delta t_j\in V^h$ for $j=2,\ldots,n$. Then we conclude that for every $\varphi^h\in V^h$ we have
\begin{align*}
    &{\langle \overline\partial e_j^{h\ell},\varphi^h\rangle}_H+a(t_j;e_j^{h\ell},\varphi^h)={\langle \overline\partial y_j^h-\overline\partial y_j^{h\ell},\varphi^h\rangle}_H+a(t_j;y^h_j-y^{h\ell}_j,\varphi^h)\\
    &\quad={\langle (\mathcal F+\mathcal Bu)(t_j),\varphi^h\rangle}_{V',V}-{\langle\overline\partial y_j^{h\ell},\varphi^h\rangle}_H-a(t_j;y^{h\ell}_j,\varphi^h)=-{\langle r_j^{h\ell},\varphi^h\rangle}_{(V^h)',V^h}
\end{align*}
for $j=2,\ldots,n$ with the residual $r_j^{h\ell}$ given by \eqref{Shostakovich-1} for $j=2,\ldots,n$. Choosing $\varphi^h=e_j^{h\ell}\in V^h$ and using $\|e_j^{h\ell}\|_{V^h}=\|e_j^{h\ell}\|_V$ for $\varphi\in V^h$ we find
\begin{align*}
    {\langle \overline\partial e_j^{h\ell},e_j^{h\ell}\rangle}_H+a(t_j;e_j^{h\ell},e_j^{h\ell})\le{\|r_j^{h\ell}\|}_{(V^h)'}{\|e_j^{h\ell}\|}_V
\end{align*}
for $j=2,\ldots,n$. From \eqref{SIAM:Eq3.1.6b}, Lemma~\ref{app_innerProd1} and Young's inequality we infer that
\begin{align*}
    &{\|e_j^{h\ell}\|}_H^2-{\|e_{j-1}^{h\ell}\|}_H^2+{\|e_j^{h\ell}-e_{j-1}^{h\ell}\|}_H^2+2\gamma_1\delta t_j\,{\|e_j^{h\ell}\|}_V^2 \\
    &\quad\le\frac{\delta t_j}{\gamma_1}\,{\|r_j^{h\ell}\|}_{(V^h)'}^2+\gamma_1\delta t_j\,{\|e_j^{h\ell}\|}_V^2\quad\text{for }2\le j\le n.
\end{align*}
This yields
\begin{equation}
    \label{IE-2}
    {\|e_j^{h\ell}\|}_H^2+\gamma_1\delta t_j\,{\|e_j^{h\ell}\|}_V^2\le{\|e_{j-1}^{h\ell}\|}_H^2+\frac{\delta t_j}{\gamma_1}\,{\|r_j^{h\ell}\|}_{(V^h)'}^2\quad\text{for }2\le j\le n.
\end{equation}
Using \eqref{Poincare} and setting $c_1=\gamma_1/c_V^2$ it follows from \eqref{IE-2} that
\begin{align*}
    {\|e_j^{h\ell}\|}_H^2\le\frac{1}{1+c_1\delta t_j}\,{\|e_{j-1}^{h\ell}\|}_H^2+\frac{\delta t_j}{\gamma_1(1+c_1\delta t_j)}\,{\|r_j^{h\ell}\|}_{(V^h)'}^2\quad\text{for }2\le j\le n.
\end{align*}
From $e_1^{h\ell}=0$, \eqref{DeltaTCond}, \eqref{DeltaTCond-2} and by summation upon $j$ we derive
\begin{align*}
    {\|e_j^{h\ell}\|}_H^2\le \frac{e^{-c_2(j-1)\delta t}}{\gamma_1}\Delta t\sum_{l=2}^{j}\,{\|r_l^{h\ell}\|}_{(V^h)'}^2\quad\text{for }2\le j\le n.
\end{align*}
with $c_2=2c_1$. Arguing as in \eqref{Chopin-6}, using again $e_1^{h\ell}=0$ and also \eqref{IE-2} the  second estimate follows.\hfill$\Box$

\subsection{Proofs of Section~\ref{SIAM-Book:Section3.6}}
\label{SIAM-Book:Section3.8.5}

\noindent{\bf\em Proof of Proposition~{\em\ref{Prop:NavierStokes-1}}.} Let us insert $\psi = y^\ell(t)$ into the Galerkin system \eqref{Eq:Jan-2}, use Assumption~\ref{Asspt-1} and apply Young's inequality to get for almost all $t \in [0,T]$
\begin{align*}
	\frac{1}{2} \frac{\mathrm d}{\mathrm dt} \, \|y^\ell(t)\|_H^2 + \gamma_1 \|y^\ell(t)\|_V^2 - \gamma_2 \|y^\ell(t)\|_H^2 &\le \|(\mathcal F+ \mathcal B u)(t) \|_{V'} \|y^\ell(t)\|_V \\
	&\le \frac{\gamma_1}{2} \|y^\ell(t)\|_V^2 + \frac{1}{2\gamma_1} \|(\mathcal F+\mathcal Bu)(t)\|_{V'}^2.
\end{align*}
Consequently, we obtain \eqref{EstAprioriState-1}, so that we can proceed as in the proof of Theorem~\ref{SIAM:Theorem3.1.1} to get a bound for $\|y^\ell\|_{L^2(0,T;V)}$; cf. \eqref{EstAprioriState-3}. To estimate $\|y_t\|_{L^2(0,T;V')}$ we argue as in the proof of Theorem~\ref{SIAM:Theorem3.1.1POD}; cf. \eqref{AppEstimate}. Hence, we obtain \eqref{Shostakovich-3}. Next we prove the uniqueness of the solution. Suppose that $\hat y^\ell$ and $\tilde y^\ell$ are two solutions to \eqref{Eq:Jan-2}. Then $y^\ell=\hat y^\ell-\tilde y^\ell$ satisfies
\begin{align*}
    \frac{\mathrm d}{\mathrm dt} \, {\langle y^\ell(t),\psi \rangle}_H + a(y^\ell(t),\varphi)+{\langle\mathcal Cy^\ell(t),\psi\rangle}_{V',V}\quad=-\langle \mathcal N(\hat y^\ell(t)) - \mathcal N(\tilde y^\ell(t)), \psi \rangle_{V',V}
\end{align*}
for all $\psi\in X^\ell$ a.e. in $(0,T]$. Choosing $\psi=y^\ell(t)\in X^\ell$ for almost all $t \in [0,T]$, utilizing \eqref{OperA&C}, \eqref{OperN} and Young's inequality, it follows that
\begin{align*}
    \frac{1}{2}\,\frac{\mathrm d}{\mathrm dt} \, {\|y^\ell(t)\|}_H^2+\gamma_1\,{\|y^\ell(t)\|}_V^2\le c_\mathcal M\,{\|\tilde y^\ell(t)\|}_V{\|y^\ell(t)\|}_H{\|y^\ell(t)\|}_V\le \frac{\gamma_1}{2}\,{\|y^\ell(t)\|}^2_V+\frac{c_\mathcal M^2}{2\gamma_1}\,{\|\tilde y^\ell(t)\|}^2_V{\|y^\ell(t)\|}^2_H
\end{align*}
for almost all $t\in[0,T]$. This especially implies
\begin{equation}
    \label{Shostakovich-2}
    \frac{\mathrm d}{\mathrm dt} \, {\|y^\ell(t)\|}_H^2\le\Big( 2\gamma_2 + \frac{c_\mathcal M^2}{\gamma_1}\,{\|\tilde y^\ell(t)\|}_V^2 \Big){\|y^\ell(t)\|}^2_H\quad\text{f.a.a. }t\in[0,T].
\end{equation}
By the Gronwall lemma (cf. Proposition~\ref{Gronwall}) and $y^\ell(0)=\hat y(0)-\tilde y(0)=y_\circ-y_\circ=0$ we obtain $\|y^\ell(t)\|_H^2 \le 0$ for almost all $t\in[0,T]$. Thus, $\hat y^\ell(t)=\tilde y^\ell(t)$ holds true for almost all $t\in[0,T]$.\hfill$\Box$

\medskip\noindent{\bf\em Proof of Theorem~{\em\ref{Th:NonlEvEqApri}}.} Recall that $X^\ell=\mathrm{span}\,\{\psi_1,\ldots,\psi_\ell\}\subset V$ for the two choices $X=H$ and $X=V$. Again we shall use the decomposition
\begin{equation}
    \label{Eq:Apr-1}
    y(t)-y^\ell(t)=y(t)-\mathcal P^\ell y(t)+\mathcal P^\ell y(t)-y^\ell(t)=\varrho^\ell(t)+\vartheta^\ell(t),
\end{equation}
where $\varrho^\ell(t)=y(t)-\mathcal P^\ell y(t)$ and $\vartheta^\ell(t)=\mathcal P^\ell y(t)-y^\ell(t)\in X^\ell$. It follows that
\begin{equation}
    \label{APriori-Est-1}
    \int_0^T{\|y(t)-y^\ell(t)\|}^2_V\,\mathrm dt\le2\int_0^T{\|\varrho^\ell(t)\|}^2_V\,\mathrm dt+2\int_0^T{\|\vartheta^\ell(t)\|}^2_V\,\mathrm dt.
\end{equation}
First we estimate the $\vartheta^\ell$-term. Since $\mathcal P^\ell$ is an $H$-orthonormal projection either for $X=H$ or for $X=V$, we can apply \eqref{Eq:SnapOhneDQ-2}. From \eqref{NonlEvPro}, \eqref{Eq:Jan-2}, \eqref{Eq:Jan-4} and $\varrho^\ell(t)=y(t)-\mathcal P^\ell y(t)$ it follows that
\begin{align*}
    &{\langle \vartheta_t^\ell(t),\psi\rangle}_H+a(\vartheta^\ell(t),\psi)+{\langle\mathcal C\vartheta^\ell(t),\psi\rangle}_{V',V}\\
    &\quad={\langle y_t(t),\psi\rangle}_H+a(t;\mathcal P^\ell y(t),\psi)+{\langle\mathcal C\mathcal P^\ell y(t),\psi\rangle}_{V',V}-{\langle y_t^\ell(t),\psi\rangle}_H-a(t;y^\ell(t),\psi)-{\langle\mathcal Cy^\ell(t),\psi\rangle}_{V',V}\\
    &\quad=a((\mathcal P^\ell y-y)(t),\psi)+{\langle\mathcal C(\mathcal P^\ell y-y)(t),\psi\rangle}_{V',V}+{\langle\mathcal N(y^\ell(t))-\mathcal N(y(t)),\psi\rangle}_{V',V}\\
    &\quad\le\gamma\,{\|\varrho^\ell(t)\|}_V{\|\psi\|}_V+{\|\mathcal C\varrho^\ell(t)\|}_{V'}{\|\psi\|}_V+\big|{\langle\mathcal N(y^\ell(t))-\mathcal N(y(t)),\psi\rangle}_{V',V}\big|.
\end{align*}
Thus, choosing $\psi=\vartheta^\ell(t)\in X^\ell$ and using \eqref{OperA&C} we have
\begin{equation}
    \label{sie80a}
    \begin{aligned}
        \frac{1}{2} \frac{\mathrm d}{\mathrm dt}\,{\|\vartheta^\ell(t)\|}^2_H+\gamma_1\,{\|\vartheta^\ell(t)\|}_V^2&\le\gamma\,{\|\varrho^\ell(t)\|}_V{\|\vartheta^\ell(t)\|}_V+{\|\mathcal C\varrho^\ell(t)\|}_{V'}{\|\vartheta^\ell(t)\|}_V\\
        &\quad+\big|{\langle\mathcal N(y^\ell(t))-\mathcal N(y(t)),\vartheta^\ell(t)\rangle}_{V',V}\big|.
    \end{aligned}
\end{equation}
Applying Young's inequality it follows that
\begin{equation}
    \label{sie80b}
    \gamma\,{\|\varrho^\ell(t)\|}_V{\|\vartheta^\ell(t)\|}_V\le\frac{\gamma_1}{10}\,{\|\vartheta^\ell(t)\|}_V^2+c_0\,{\|\varrho^\ell(t)\|}_V^2
\end{equation}
and
\begin{equation}
    \label{sie80}
    \begin{aligned}
        {\|\mathcal C\varrho^\ell(t)\|}_{V'}{\|\vartheta^\ell(t)\|}_V&\le {\|\mathcal C\|}_{\mathscr L(V,V')}\,{\| \varrho^\ell(t)\|}_V{\|\vartheta^\ell(t)\|}_V\\
        &\le\frac{\gamma_1}{10}\,{\|\vartheta^\ell(t)\|}_V^2+c_0\,{\|\varrho^\ell(t)\|}_V^2
    \end{aligned}
\end{equation}
for a constant $c_0>0$ depending on $\|\mathcal C\|_{\mathscr L(V,V')}$, $\gamma$ and $\gamma_1$. We proceed by estimating the nonlinear terms on the right-hand side of \eqref{sie80a}. Note that
\begin{equation}
    \label{eq4-36}
    \mathcal N(y^\ell(t))-\mathcal N(y(t))=\mathcal M(y(t),y^\ell(t)-y(t))+\mathcal N(y^\ell(t)-y(t))+\mathcal M(y^\ell(t)-y(t),y(t)).
\end{equation}
Applying \eqref{OperN}, \eqref{Poincare}, and Young's inequality (with $p=q=2$ and $\varepsilon=2\gamma_1/5$) we obtain satisfying
\begin{equation}
    \label{eq4-38}
    \begin{aligned}
        \big|{\langle\mathcal M(y(t),y^\ell(t)-y(t)),\vartheta^\ell(t)\rangle}_{V',V}\big|&=\big|{\langle\mathcal M(y(t),\varrho^\ell(t)),\vartheta^\ell(t)\rangle}_{V',V}\big|\\
        &\le c_\mathcal Mc_V^{\delta_3}\,{\| y \|}_{C([0,T];V)}{\|\varrho^\ell(t)\|}_V{\|\vartheta^\ell(t) \|}_V^{\delta_3}{\|\vartheta^\ell(t)\|}_H^{1-\delta_3}\\
        &\le c_1\,{\|\varrho^\ell(t)\|}_V^2+\frac{\gamma_1}{5}\,{\|\vartheta^\ell(t) \|}_V^{2\delta_3}{\|\vartheta^\ell(t)\|}_H^{2(1-\delta_3)}
    \end{aligned}
\end{equation}
with $c_1=5c_\mathcal M^2c_V^{2\delta_3}\,\| y \|_{C([0,T];V)}^2/(4\gamma_1)$. Again applying Young's inequality (with $p=1/\delta_3$, $q=1/(1-\delta_3)$ and $\varepsilon=1/(2\delta_3)$) we find
\begin{equation}
    \label{eq4-38-1}
    {\|\vartheta^\ell(t) \|}_V^{2\delta_3}{\|\vartheta^\ell(t)\|}_H^{2(1-\delta_3)}\le \frac{1}{2}\,{\|\vartheta^\ell(t) \|}_V^2+ \frac{1}{c_2}\,{\|\vartheta^\ell(t) \|}_H^2+c_1\,{\|\varrho^\ell(t)\|}_V^2
\end{equation}
with $c_2=(\delta_3/2)^{\delta_3/(1-\delta_3)}/(1-\delta_3)$. Combining \eqref{eq4-38} and \eqref{eq4-38-1} we derive that
\begin{equation}
    \label{eq4-38a}
    \begin{aligned}
        \big|{\langle\mathcal M(y(t),y^\ell(t)-y(t)),\vartheta^\ell(t)\rangle}_{V',V}\big|&\le c_1\,{\|\varrho^\ell(t)\|}_V^2+\frac{\gamma_1}{5}\,\bigg(\frac{1}{2}\,{\|\vartheta^\ell(t) \|}_V^2+ \frac{1}{c_2}\,{\|\vartheta^\ell(t) \|}_H^2\bigg)\\
        &=\frac{\gamma_1}{10}\,{\|\vartheta^\ell(t) \|}_V^2+\frac{\gamma_1}{5c_2}\,{\|\vartheta^\ell(t) \|}_H^2+ c_1\,{\|\varrho^\ell(t)\|}_V^2.
    \end{aligned}
\end{equation}
Again utilizing \eqref{OperN}, Young's inequality and \eqref{Poincare} we find that there exist constants $c_3,c_4>0$ such that
\begin{equation}
    \label{eq4-39}
    \begin{aligned}
        &\big|{\langle\mathcal M(y^\ell(t)-y(t),y(t)),\vartheta^\ell(t)\rangle}_{V',V}\big|\\
        &=\big|{\langle\mathcal M(\vartheta^\ell(t),y(t))+\mathcal M(\varrho^\ell(t),y(t)),\vartheta^\ell(t)\rangle}_{V',V}\big|\\
        &\le  c_\mathcal M\,{\|y\|}_{C([0,T];V)} \left({\|\vartheta^\ell(t)\|}_H {\|\vartheta^\ell(t)\|}_V+c_V^{\delta_3}\,{\|\varrho^\ell(t)\|}_V{\|\vartheta^\ell(t)\|}_H^{1-\delta_3}\,{\| \vartheta^\ell(t)\|}_V^{\delta_3}\right)\\
        &\le \frac{\gamma_1}{10}\,{\|\vartheta^\ell(t)\|}_V^2+c_3\,{\|\vartheta^\ell(t)\|}_H^2+c_4\,{\|\varrho^\ell(t)\|}_V^2.
    \end{aligned}
\end{equation}
By assumption, $y$ belongs to $L^\infty(0,T;V)$. Due to \eqref{Poincare} we have
\begin{align*}
    {\|\varrho^\ell(t)\|}_H\le c_V\,{\|y(t)-\mathcal P^\ell y(t)\|}_V\le c_V\,{\|y(t)\|}_V\quad\text{f.a.a. }t\in[0,T].
\end{align*}
Therefore, there exists a constant $c_5>0$ such that
\begin{equation}
    \label{eq4-41}
    \esssup_{t\in[0,T]}\Big({\|\varrho^\ell(t)\|}_H^{\delta_3}{\| \varrho^\ell(t)\|}_V^{1-\delta_3},{\|\varrho^\ell(t)\|}_V\Big)\le c_5.
\end{equation}
Using \eqref{Eq:Apr-1} and \eqref{OperN} we conclude
\begin{equation}
    \label{eq4-40}
    \begin{aligned}
        &{\langle\mathcal N(y^\ell(t)-y(t)),\vartheta^\ell(t)\rangle}_{V',V}\\
        &\qquad=-{\langle\mathcal M(\vartheta^\ell(t),\varrho^\ell(t))-\mathcal M(\varrho^\ell(t),\varrho^\ell(t)),\vartheta^\ell(t)\rangle}_{V',V}.
    \end{aligned}
\end{equation}
Applying \eqref{OperA&C}, \eqref{eq4-41}, \eqref{eq4-40} and Young's inequality we find
\begin{equation}
    \label{eq4-42}
    \begin{aligned}
        &\big|{\langle\mathcal N(y^\ell(t)-y(t)),\vartheta^\ell(t)\rangle}_{V',V}\big|\\
        &\quad\le c_\mathcal Mc_5\,\Big({\|\vartheta^\ell(t)\|}_H{\|\vartheta^\ell(t)\|}_V+{\|\varrho^\ell(t)\|}_V{\|\vartheta^\ell(t)\|}_H^{1-\delta_3}{\|\vartheta^\ell(t)\|}_V^{\delta_3}\Big)\\
        &\quad\le\frac{\gamma_1}{10}\,{\| \vartheta^\ell(t) \|}_V^2+ c_6\,{\|\vartheta^\ell(t)\|}_H^2+c_7\,{\|\varrho^\ell(t)\|}_V^2
    \end{aligned}
\end{equation}
for two constants $c_6,c_7>0$. From \eqref{sie80a}-\eqref{eq4-42}
\begin{equation}
    \label{Eq:Jan-5}
    \frac{\mathrm d}{\mathrm dt}\,{\|\vartheta^\ell(t)\|}^2_H+\frac{\gamma_1}{2}\,{\|\vartheta^\ell(t)\|}_V^2\le c_8\,\,{\|\vartheta^\ell(t)\|}_H^2+c_9\,{\|\varrho^\ell(t)\|}_V^2,
\end{equation}
where $c_8=c_1+c_3+c_6$ and $c_9=2c_0+c_2+c_4+c_7$. Let
\begin{align*}
    c_{10}=\left\{
    \begin{aligned}
        &0&&\text{if }1-\frac{\gamma_1}{2c_V^2}\le0,\\
        &1-\frac{\gamma_1}{2c_V^2}&&\text{otherwise}
    \end{aligned}
    \right.
\end{align*}
Thus, from $\vartheta^\ell(0)=0$, \eqref{Poincare}, \eqref{Eq:Jan-5} and Proposition~\ref{Gronwall} we obtain
\begin{equation}
    \label{eq4-43a}
    {\|\vartheta^\ell(t)\|}^2_H\le c_{11}\int_0^t{\|y(s)-\mathcal P^\ell y(s)\|}_V^2\,\mathrm ds\quad\text{in }[0,T]\text{ a.e.},
\end{equation}
for $c_{11}=c_{12}\exp(c_{10}T)$ and $c_{12}=\max(c_8,c_9)$. Thus,
\begin{equation}
    \label{eq4-43}
    \int_0^T{\|\vartheta^\ell(t)\|}^2_H\,\mathrm dt\le c_{13}\int_0^T{\|y(t)-\mathcal P^\ell y(t)\|}_V^2\,\mathrm dt\quad\text{in }[0,T]\text{ a.e.}
\end{equation}
with $c_{13}=c_{11}T$. Integrating \eqref{Eq:Jan-5} over $[0,T]$ and using $\vartheta^\ell(0)=0$ as well as \eqref{eq4-43a} it follows that
\begin{equation}
    \label{Eq:Jan-6}
    \int_0^T{\|\vartheta^\ell(t)\|}^2_V\,\mathrm dt\le\frac{2}{\gamma_1}\int_0^T c_8\,{\|\vartheta^\ell(t)\|}_H^2+c_9\,{\|\varrho^\ell(t)\|}_V^2\,\mathrm dt\le c_{14}\int_0^T{\|y(t)-\mathcal P^\ell y(t)\|}_V^2\,\mathrm dt
\end{equation}
with $c_{14}=2\max(c_8c_{11}T,c_9)/\gamma_1$. Thus, we find from \eqref{APriori-Est-1} that
\begin{align*}
    \int_0^T{\|y(t)-y^\ell(t)\|}_V^2\,\mathrm dt\le2\big(1+c_{14}\big)\int_0^T{\|y(t)-\mathcal P^\ell y(t)\|}_V^2\,\mathrm dt
\end{align*}
for the $H$-orthogonal projection operators $\mathcal P^\ell=\mathcal P^\ell_H$ and $\mathcal P^\ell=\mathcal Q^\ell_H$. From Proposition~\ref{Prop:Rate-V} we finally infer the a-priori error estimate.\hfill$\Box$

\medskip\noindent{\bf\em Proof of Theorem~{\em\ref{Th:A-PrioriError-300}}.} For the a-priori error analysis we follow the proof of Theorem~\ref{Th:NonlEvEqApri}. Thus, we make use of the following decomposition
\begin{align*}
    y^h(t)-y^{h\ell}(t)&=y^h(t)-\mathcal P^\ell y^h(t)+\mathcal P^\ell y^h(t)-y^{h\ell}(t)=\varrho^{h\ell}(t)+\vartheta^{h\ell}(t)\quad\text{in }[0,T]\text{ a.e.}
\end{align*}
with $\varrho^{h\ell}(t)=y^h(t)-\mathcal P^\ell y^h(t)\in (X^\ell)^\bot$ and $\vartheta^{h\ell}(t)=\mathcal P^\ell y^h(t)-y^{h\ell}(t)\in X^\ell$. It is shown in Section~\ref{SIAM-Book:Section3.3.2} that $\vartheta^{h\ell}(0)=0$ holds for the choices $\mathcal P^\ell=\mathcal P^\ell_H$, $\mathcal P^h=\mathcal P^h_H$ and $\mathcal P^\ell=\mathcal Q^\ell_H$, $\mathcal P^h=\mathcal P^h_V$. Now we argue as in the proof of Theorem~\ref{Th:NonlEvEqApri}.\hfill$\Box$

\medskip\noindent{\bf\em Proof of Proposition~{\em\ref{Prop:NSAposti}}.} Let $y_\circ\in V$, $f\in L^2(0,T;H)$ and $u\in\U$ with $\mathcal Bu\in L^2(0,T;H)$. Suppose that Assumptions~\ref{Asspt-1} and \ref{A7} hold. We set $e^{h\ell}(t)=y^h(t)-y^{h\ell}(t)\in V^h$ for the error between $y^h$ and $y^{h\ell}$. Then we derive from \eqref{NonlEvProGal} and \eqref{NonlEvProGalPOD} that for almost all $t\in[0,T]$ and for every $\varphi^h\in V^h$
\begin{equation}
    \label{Alb-1}
    \begin{aligned}
        &{\langle e_t^{h\ell}(t),\varphi^h\rangle}_{V',V}+a(e^{h\ell}(t),\varphi^h)+{\langle\mathcal Ce^{h\ell}(t),\varphi^h\rangle}_{V',V}\\
        &\quad={\langle y_t^h(t),\varphi^h\rangle}_{V',V}+a(y^h(t),\varphi^h)+{\langle\mathcal Cy^h(t),\varphi^h\rangle}_{V',V}\\
        &\qquad-{\langle y_t^{h\ell}(t),\varphi^h\rangle}_{V',V}-a(y^{h\ell}(t),\varphi^h)-{\langle\mathcal Cy^{h\ell}(t),\varphi^h\rangle}_{V',V}\\
        &\quad={\langle\mathcal N(y^{h\ell}(t))-\mathcal N(y^h(t)),\varphi^h\rangle}_{V',V}-{\langle r^{h\ell}(t),\varphi^h\rangle}_{V',V}
    \end{aligned}
\end{equation}
with the time dependent residual $r^{h\ell}(t)$ defined in \eqref{Eq:timeResidual}. Choosing $\varphi^h=e^{h\ell}(t)\in V^h$ in \eqref{Alb-1} we derive from \eqref{OperA&C} that
\begin{align*}
    {\langle e_t^{h\ell}(t),e^{h\ell}(t)\rangle}_{V',V}+a(e^{h\ell}(t),e^{h\ell}(t))+{\langle\mathcal Ce^{h\ell}(t),e^{h\ell}(t)\rangle}_{V',V}\ge\frac{1}{2}\frac{\mathrm d}{\mathrm dt}\,{\|e^{h\ell}(t)\|}_H^2+\gamma_1\,{\|e^{h\ell}(t)\|}_V^2
\end{align*}
for almost all $t\in[0,T]$. Recall that the operator $\mathcal N$ is bilinear. Moreover, $\mathcal N(y^{h\ell}(t))=\mathcal N(y^{h\ell}(t),y^{h\ell}(t))$ and $\mathcal N(y^h(t))=\mathcal N(y^h(t),y^h(t))$ hold for every $t\in[0,T]$. Thus, we have
\begin{equation}
    \label{Queen-100}
    \begin{aligned}
        \mathcal N(y^{h\ell}(t))-\mathcal N(y^h(t))&=\mathcal M(y^{h\ell}(t),y^{h\ell}(t))-\mathcal M(y^h(t),y^{h\ell}(t))\\
        &\quad-\big(\mathcal M(y^h(t),y^h(t))-\mathcal M(y^h(t),y^{h\ell}(t))\big)\\
        &=-\mathcal M(e^{h\ell}(t),y^{h\ell}(t))-\mathcal M(y^h(t),e^{h\ell}(t))
    \end{aligned}
\end{equation}
for almost all $t\in[0,T]$. By \eqref{OperN} we have
\begin{align*}
    {\langle\mathcal M(y^h(t),e^{h\ell}(t)),e^{h\ell}(t)\rangle}_{V',V}=0
\end{align*}
for almost all $t\in[0,T]$. Consequently, \eqref{Alb-1}, \eqref{OperN} and the Young inequality imply that
\begin{align*}
    \frac{1}{2}\frac{\mathrm d}{\mathrm dt}\,{\|e^{h\ell}(t)\|}_H^2+\gamma_1\,{\|e^{h\ell}(t)\|}_V^2&\le{\|r^{h\ell}(t)\|}_{(V^h)'}{\|e^{h\ell}(t)\|}_V+\big|{\langle\mathcal M(e^{h\ell}(t),y^{h\ell}(t)),e^{h\ell}(t)\rangle}_{V',V}\big|\\
    &\le\frac{1}{4\gamma_1}\,{\|r^{h\ell}(t)\|}_{(V^h)'}^2+\frac{\gamma_1}{4}\,{\|e^{h\ell}(t)\|}_V^2+c_\mathcal M\,{\|e^{h\ell}(t)\|}_H{\|y^{h\ell}(t)\|}_V{\|e^{h\ell}(t)\|}_V\\
    &\le\frac{1}{\gamma_1}\,{\|r^{h\ell}(t)\|}_{(V^h)'}^2+\frac{\gamma_1}{2}\,{\|e^{h\ell}(t)\|}_V^2+\frac{c_\mathcal M^2}{\gamma_1}\,{\|y^{h\ell}(t)\|}_V^2{\|e^{h\ell}(t)\|}_H^2
\end{align*}
which implies 
\begin{equation}
    \label{Alb-2}
    \frac{\mathrm d}{\mathrm dt}\,{\|e^{h\ell}(t)\|}_H^2+\gamma_1\,{\|e^{h\ell}(t)\|}_V^2\le\frac{2}{\gamma_1}\Big({\|r^{h\ell}(t)\|}_{(V^h)'}^2+c_\mathcal M^2\,{\|y^{h\ell}(t)\|}_V^2{\|e^{h\ell}(t)\|}_H^2\Big)
\end{equation}
for almost all $t\in[0,T]$. From $e^{h\ell}(0)=y^h(0)-y^{h\ell}(0)=(\mathcal P^h-\mathcal P^{h\ell})y_\circ$ and Proposition~\ref{Gronwall} it follows that
\begin{equation}
    \label{Alb-3}
    {\|e^{h\ell}(t)\|}_H^2\le \mathsf C_1^{h\ell}(t)\bigg(\mathsf R_\circ^{h\ell}+\int_0^t\mathsf R_1^{h\ell}(s)\,\mathrm ds\bigg)\quad\text{for }t\in[0,T]\text{ a.e.}
\end{equation}
with the constants defined in \eqref{Alb-4a} and \eqref{Alb-4b}. Integrating \eqref{Alb-2} over $[0,t]$ with $t\in(0,T]$ and using \eqref{Alb-3} we find that
\begin{align*}
    \int_0^t{\|e^{h\ell}(s)\|}_V^2\,\mathrm ds&\le \tilde{\mathsf R}_\circ^{h\ell}+\int_0^t\tilde{\mathsf R}_1^{h\ell}(s)+\frac{2 c_\mathcal M^2}{\gamma_1^2}\,{\|y^{h\ell}(s)\|}_V^2{\|e^{h\ell}(s)\|}_H^2\,\mathrm ds\\
    &\le\tilde{\mathsf R}_\circ^{h\ell}+\int_0^t\tilde{\mathsf R}_1^{h\ell}(s)+\mathsf R_2^{h\ell}(s)\bigg(\mathsf R_\circ^{h\ell}+\int_0^s\mathsf R_1^{h\ell}(\tau)\,\mathrm d\tau\bigg)\,\mathrm ds
\end{align*}
with the constant in \eqref{Alb-5}.\hfill$\Box$

\medskip\noindent{\bf\em Proof of Proposition~{\em\ref{PropNonl:FullDiscFEModel}}.} We follow the arguments of the proof of Theorem~4.2 in \cite{KV02a}. Existence of a solution $\{y_j^h\}_{j=1}^n$ can be proved by using the Schauder fixed point theorem, see \cite[p.~222]{GT77}, for instance. For that purpose we define $z^h=\mathcal T_j^hw^h$ via the mappings $\mathcal T_j^h:V^h\to V^h$, $j=1,\ldots,n$, as follows: $z^h\in V^h$ is the solution to
\begin{equation}
    \label{eq4-11c}
    \begin{aligned}
        &{\langle z^h,\varphi^h\rangle}_H+\delta t_j\Big(a(z^h,\varphi^h)+{\langle \mathcal Cz^h+\mathcal M(w^h,z^h),\varphi^h\rangle}_{V',V}\Big)\\
        &\qquad={\langle y_{j-1}^h,\varphi^h\rangle}_H+\delta t_j\,{\langle g_j(u),\varphi^h\rangle}_{V',V}\quad\text{for all }\varphi^h\in V^h.
    \end{aligned}
\end{equation}
The bilinear form
\begin{align*}
    {\langle \cdot \, ,\cdot \rangle}_H+\delta t_j\left(a(\cdot \,,\cdot)+{\langle \mathcal C(\cdot)+\mathcal M(w^h,\cdot),\cdot\rangle}_{V',V}\right)
\end{align*}
is continuous and coercive in $V^h\times V^h$ by \eqref{OperC}-\eqref{OperN}. The existence and uniqueness of a solution to \eqref{eq4-11c} can thus be shown by the Lax-Milgram theorem. The fixed points of $\mathcal T_j^h$ are the solutions of \eqref{NonlEvProGalFullyDisc}. Taking $\varphi^h=z^h$ we infer from \eqref{eq4-11c}, \eqref{OperA&C}, \eqref{OperN} and Young's inequality that
\begin{align*}
    {\|z^h\|}_H^2+\gamma_1\delta t_j\,{\|z^h\|}_V^2&\le\frac{1}{2}\,{\|y^h_{j-1}\|}_H^2+\frac{1}{2}\,{\|z^h\|}_H^2+\frac{\delta t_j}{2\gamma_1}\,{\|g_j(u)\|}_{V'}^2+\frac{\gamma_1\delta t_j}{2}\,{\|z^h\|}_V^2
\end{align*}
which implies that
\begin{align*}
    {\|z^h\|}_H^2+\gamma_1\delta t_j\,{\|z^h\|}_V^2\le{\|y^h_{j-1}\|}_H^2+\frac{\delta t_j}{\gamma_1}\,{\|g_j(u)\|}_{V'}^2.
\end{align*}
Consequently,
\begin{equation}
    \label{eq4-11e}
    {\|z^h\|}_V\le\frac{1}{\gamma_1^2\delta t_j}\,\Big(\gamma_1{\|y^h_{j-1}\|}_H^2 +\delta t_j\,{\|g_j(u)\|}_{V'}^2\Big).
\end{equation}
Let us introduce the set
\begin{align*}
    \mathscr M_j^h=\bigg\{w^h\in V^h\,\big|\,{\|w^h\|}_V^2\le\frac{1}{\gamma_1^2\delta t_j}\,\Big(\gamma_1{\|y^h_{j-1}\|}_H^2 +\delta t_j\,{\|g_j(u)\|}_{V'}^2\Big)\bigg\}.
\end{align*}
From \eqref{eq4-11e} we infer that $\mathcal T_j^h$ maps $\mathscr M_j^h$ into itself. Since $\mathscr M_j^h$ is a closed ball in $V^h$, the set $\mathscr M_j^h$ is bounded, closed and convex. Since the image of $\mathcal T_j^h$ is finite dimensional, $\mathcal T_j^h$ is compact. Thus, the existence of a fixed point $y_j^h$ follows from the Schauder fixed point theorem. Next we show that the solution to \eqref{NonlEvProGalFullyDisc} is unique. For that purpose we assume that the two sequences $\{y_j^h\}_{j=1}^n$, $\{\tilde y_j^h\}_{j=1}^n$ in $V^h$ are solutions of \eqref{NonlEvProGalFullyDisc}. Then $z_j^h= y_j^h-\tilde y_j^h\in V^h$ solves
\begin{align*}
    {\langle z_j^h,\varphi^h\rangle}_H+\delta t_j\left(a(z_j^h,\varphi^h)+{\langle \mathcal Cz_j^h,\varphi^h\rangle}_{V',V}\right)=\delta t_j\,{\langle\mathcal N(\tilde y^h_j)-\mathcal N(y_j^h),\varphi\rangle}_{V',V}
\end{align*}
for all $\varphi\in V^h$. Setting $\varphi=z_j^h$ and using \eqref{OperA&C}, \eqref{OperN} and Young's inequality we obtain
\begin{align*}
    {\|z_j^h\|}_H^2+\gamma_1\delta t_j\,{\|z_j^h\|}_V^2 & \le\delta t_j\,{\langle \mathcal N(\tilde y_j^h)-\mathcal N(y_j^h),z_j^k\rangle}_{V',V}=-\delta t_j\,{\langle\mathcal M(z_j^h,\tilde y_j^h)+\mathcal M(y_j^h,z_j^h),z_j^h\rangle}_{V',V}\\
    &=-\delta t_j~{\langle\mathcal M(z_j^h,\tilde y_j^h),z_j^h \rangle}_{V',V}\le c_\mathcal M \delta t_j\,{\|\tilde y_j^h\|}_V{\|z_j^h \|}_H{\|z_j^h \|}_V\\
    & \le {\|z_j^h\|}_H^2+ \frac{c_\mathcal M^2 \delta t_j^2}{4}~{\|\tilde y_j^h\|}_V^2{\|z_j^h\|}_V^2.
\end{align*}
It follows that
\begin{align*}
    \Big(1-\frac{c_\mathcal M^2 \delta t_j}{4\gamma_1}\,{\|\tilde y_j^h\|}_V^2\Big)\,{\|z_j^h\|}_V^2\le 0.
\end{align*}
Let $c=\max\{\| \tilde y_j^h\|_V\,|\,j=1,\ldots,n\}$. Then $z_j^h=0$ and hence $y_j^h=\tilde y_j^h$ provided that $\Delta t\le4\gamma_1/(c^2c_\mathcal M^2)$.\hfill\\
To prove the estimates \eqref{Nonl:TempDiscFESOL-APriori} and \eqref{Nonl:TempDiscFESOL-APriori-2} we take $\psi=y_j^h$ in \eqref{NonlEvProGalFullyDisc}. Due to \eqref{OperA&C}, \eqref{OperN} and Lemma~\ref{app_lem_innerProd1} we obtain
\begin{align*}
    {\|y_j^h\|}_H^2-{\|y_{j-1}^h\|}_H^2+{\|y_j^h-y_{j-1}^h\|}_H^2 +2\gamma_1\delta t_j\,{\|y_j^h\|}_V^2\le 2\delta t_j~{\| g_j(u) \|}_{V'} {\|y_j^h\|}_V.
\end{align*}
Using Young's inequality it follows that
\begin{equation}
    \label{eq4-14}
    {\|y_j^h\|}_H^2+{\|y_j^h-y_{j-1}^h\|}_H^2+\gamma_1\delta t_j\,{\| y_j^h\|}_V^2\le{\|y_{j-1}^h\|}_H^2+\frac{\delta t_j}{\gamma_1}\,{\|g_j(u) \|}_{V'}^2.
\end{equation}
From \eqref{eq4-14} and \eqref{Poincare} we infer that
\begin{align*}
    (1+c_1 \delta t_j)\,{\|y_j^h\|}_H^2\le{\|y_{j-1}^h\|}_H^2+\frac{\delta t_j}{\gamma_1}\,{\|g_j(u) \|}_{V'}^2,
\end{align*}
where $c_1=\gamma_1/c_V^2$, which yields
\begin{align*}
    {\|y_j^h\|}_H^2&\le\frac{1}{1+c_1\delta t}~{\|y_{j-1}^h\|}_H^2+\frac{\delta t_j}{\gamma_1(1+c_1\delta t_j)}\,{\|g_j(u) \|}_{V'}^2\\
    &\le\bigg(\frac{1}{1+c_1\delta t}\bigg)^{j-1}\bigg({\|y_1^h\|}_H^2+\frac{1}{\gamma_1}\sum_{l=2}^j\delta t_l\,{\|g_l(u)\|}_{V'}^2\bigg)
\end{align*}
which coincides \eqref{AprEst100a}. Hence we can proceed as in the proof of Theorem~\ref{FullDiscFEModel}. Suppose that $\delta t\le1/C_1$ holds with $C_1=2c_1$ (cf. \eqref{DeltaTCond}). Then we get directly \eqref{Nonl:TempDiscFESOL-APriori} for $j=2,\ldots,n$. Summation upon $j$ in \eqref{eq4-14} we obtain \eqref{Nonl:TempDiscFESOL-APriori-2}.\hfill$\Box$

\medskip\noindent{\bf\em Proof of Theorem~{\em\ref{Theorem:AprioriEstNonl}}.} Again we utilize the decomposition
\begin{align*}
    y_j^h-y_j^{h\ell}=y_j^h-\mathcal P^{n\ell}y_j^h+\mathcal P^{n\ell}y_j^h-y_j^{h\ell}=\varrho_j^{h\ell}+\vartheta_j^{h\ell}\quad\text{for }1\le j\le n
\end{align*}
with $\varrho_j^{h\ell}=y_j^h-\mathcal P^{n\ell}y_j^h \in (X^{n\ell})^\bot$ and $\vartheta_j^{h\ell}=\mathcal P^{n\ell}y_j^h-y_j^{h\ell}\in X^{n\ell}$. Here, we consider the $H$-orthonormal projections. Thus, we take $\mathcal P^{n\ell}=\mathcal P^{n\ell}_H$ for $X=H$ and $\mathcal P^{n\ell}=\mathcal Q^{h\ell}_H$ for $X=V$. It follows from Proposition~\ref{Prop:VTopologyDisc} that
\begin{equation}
    \label{Bach-1}
    \sum_{j=1}^n\alpha_j^n\,{\|\varrho_j^{h\ell}\|}^2_V=\sum_{j=1}^n\alpha_j^n\,{\|y_j^h-\mathcal P_H^{n\ell}y_j^h\|}^2_V=\sum_{i=\ell+1}^{d^n}\lambda_i^{nH}\,{\|\psi_i^{nH}\|}_V^2
\end{equation}
for $X=H$ and
\begin{equation}
    \label{Bach-2}
    \sum_{j=1}^n\alpha_j^n\,{\|\varrho_j^{h\ell}\|}^2_V=\sum_{j=1}^n\alpha_j^n\,{\|y_j^h-\mathcal Q_H^{n\ell}y_j^h\|}^2_V=\sum_{i=\ell+1}^{d^n}\lambda_i^{nV}\,{\|\psi_i^{nV}-\mathcal P^{n\ell}\psi_i^{nV}\|}_V^2
\end{equation}
for $X=V$. Now we follow Lemma 4.5 in \cite{KV02a}. Since $y_1^h=\mathcal P^h_Hy_\circ$ holds, we have $\vartheta_1^{h\ell}=0$ for our chosen projections $\mathcal P^{n\ell}$. Using the notaion $\overline\partial\vartheta_j^{h\ell}=(\vartheta_j^{h\ell}-\vartheta_{j-1}^{h\ell})/\delta t_j$, $j=2,\ldots,n$, we obtain for every $\psi\in X^{n\ell}$
\begin{equation}
    \label{Bach-3}
    \begin{aligned}
        {\langle\overline\partial\vartheta_j^{h\ell},\psi\rangle}_H+a(\vartheta_j^{h\ell},\psi)+{\langle \mathcal C \vartheta_j^{h\ell},\psi\rangle}_{V',V}&={\langle\overline\partial \mathcal P^{n\ell}y_j^h,\psi\rangle}_H+a(\mathcal P^{n\ell}y_j^h,\psi)+{\langle\mathcal C\mathcal P^{n\ell}y_j^h,\psi\rangle}_{V',V}\\
        &\quad-\big({\langle\overline\partial y_j^{h\ell},\psi\rangle}_H+a(y_j^{h\ell},\psi)+{\langle\mathcal Cy_j^{h\ell},\psi\rangle}_{V',V}\big).
    \end{aligned}
\end{equation}
From \eqref{Bach-100} it follows that
\begin{align*}
    {\langle\overline\partial \mathcal P^{n\ell}y_j^h,\psi^h\rangle}_H={\langle\overline\partial y_j^h,\psi\rangle}_H\quad\text{for }j=2,\ldots,n.
\end{align*}
Thus, we infer from \eqref{Bach-3}, \eqref{NonlEvProGalFullyDisc} and \eqref{Nonl:EvProGalPOD-disc} that
\begin{equation}
    \label{Bach-4}
    \begin{aligned}
        &{\langle\overline\partial\vartheta_j^{h\ell},\psi\rangle}_H+a(\vartheta_j^{h\ell},\psi)+{\langle\vartheta_j^{h\ell},\psi\rangle}_{V',V}\\
        &\quad={\langle\mathcal N(y_j^{h\ell})-\mathcal N(y_j^h),\psi\rangle}_{V',V}-a(\varrho^{h\ell}_j,\psi)-{\langle\mathcal C\varrho^{h\ell}_j,\psi\rangle}_{V',V}
    \end{aligned}
\end{equation}
for $j=2,\ldots,n$. Choosing $\psi=\vartheta^{h\ell}_j\in X^{n\ell}$ and utilizing Lemma~\ref{app_lem_innerProd1} and \eqref{OperA&C} we derive from \eqref{Bach-4} that
\begin{equation}
    \label{Bach-5}
    \begin{aligned}
        &{\|\vartheta_j^{h\ell}\|}_H^2-{\|\vartheta_{j-1}^{h\ell}\|}_H^2+{\|\vartheta_j^{h\ell}-\vartheta_{j-1}^{h\ell}\|}_H^2+2\gamma_1\delta t_j\,{\|\vartheta_j^{h\ell}\|}_V^2\\
        &\quad\le2\delta t_j\big(\gamma\,{\|\varrho_j^{h\ell}\|}_V{\|\vartheta_j^{h\ell}\|}_V+{\|\mathcal C\varrho_j^{h\ell}\|}_{V'}{\|\vartheta_j^{h\ell}\|}_V\big)+2\delta t_j\,\big|{\langle\mathcal N(y_j^{h\ell})-\mathcal N(y_j^h),\vartheta_j^{h\ell}\rangle}_{V',V}\big|
    \end{aligned}
\end{equation}
for $j=2,\ldots,n$. Applying Young's inequality it follows that
\begin{equation}
    \label{Bach-5a}
    \gamma\,{\|\varrho_j^{h\ell}\|}_V{\|\vartheta_j^{h\ell}\|}_V\le\frac{\gamma_1}{5}\,{\|\vartheta_j^{h\ell}\|}_V^2+c_1\,{\|\varrho_j^{h\ell}\|}_V^2
\end{equation}
for a non-negative constant $c_1$ depending on $\gamma$, $\gamma_1$ and
\begin{equation}
    \label{Bach-6}
    {\|\mathcal C\varrho_j^{h\ell}\|}_{V'}{\|\vartheta_j^{h\ell}\|}_V\le{\|\mathcal C\|}_{\mathscr L(V,V')}{\|\varrho_j^{h\ell}\|}_V{\|\vartheta_j^{h\ell}\|}_V\le\frac{\gamma_1}{5}\,{\|\vartheta_j^{h\ell}\|}_V^2+c_2\,{\|\varrho_j^{h\ell}\|}_V^2
\end{equation}
for a constant $c_2\ge0$ depending on $\gamma_1$ and $\|\mathcal C\|_{\mathscr L(V,V')}$. We proceed by estimating the nonlinear terms on the right-hand side of \eqref{Bach-5}. Note that
\begin{align*}
    \mathcal N(y_j^{h\ell})-\mathcal N(y_j^h)=\mathcal M(y_j^h,y_j^{h\ell}-y_j^h)+\mathcal M(y_j^{h\ell}-y_j^h)+\mathcal M(y_j^{h\ell}-y_j^h,y_j^h).
\end{align*}
By Assumption~\ref{AssNL5}-4) there exists a constant $c_3=C$ which is independent of $h$ and $n$ such that \eqref{Bach-7} holds. Using \eqref{OperN}, \eqref{Poincare}, Young's inequality and \eqref{Bach-7} we conclude that there are two non-negative constants $c_4$, $c_5$, which depend on $\gamma_1$, $c_V$, $c_3$, $c_\mathcal N$ and $\delta_3$, so that
\begin{equation}
    \label{Bach-8}
    \begin{aligned}
        &\big|{\langle\mathcal M(y_j^h,y_j^{h\ell}-y_j^h),\vartheta_j^{h\ell}\rangle}_{V',V}\big|\\
        &\quad=\big|{\langle\mathcal M(y_j^h,\varrho_j^{h\ell})+\mathcal M(y_j^h,\vartheta_j^{h\ell}),\vartheta_j^{h\ell}\rangle}_{V',V}\big|=\big|{\langle\mathcal M(y_j^h,\varrho_j^{h\ell}),\vartheta_j^{h\ell}\rangle}_{V',V}\big|\\
        &\quad\le c_\mathcal Nc_V^{\delta_3}\,{\|y_j^h\|}_V{\|\varrho_j^{h\ell}\|}_V{\|\vartheta_j^{h\ell}\|}_H^{1-\delta_3}{\|\vartheta_j^{h\ell}\|}_V^{\delta_3}\le\frac{\gamma_1}{5}\,{\|\vartheta_j^{h\ell}\|}_V^2+c_4\,{\|\vartheta_j^{h\ell}\|}_H^2+c_5\,{\|\varrho_j^{h\ell}\|}_V^2.
    \end{aligned}
\end{equation}
Again applying \eqref{OperN}, \eqref{Bach-7}, \eqref{Poincare} and Young's inequality we infer that there exist two non-negative constants $c_6$, $c_7$, which depend on $\gamma_1$, $c_V$, $c_3$, $c_\mathcal N$ and $\delta_3$, satisfying
\begin{equation}
    \label{Bach-9}
    \begin{aligned}
        \big|{\langle\mathcal M(y_j^{h\ell}-y_j^h,y_j^h),\vartheta_j^{h\ell}\rangle}_{V',V}\big|&= \big|{\langle\mathcal M(\varrho_j^{h\ell},y_j^h)+\mathcal M(\vartheta_j^{h\ell},y_j^h),\vartheta_j^{h\ell}\rangle}_{V',V}\big|\\
        &\le c_\mathcal M\,{\|y_j^h\|}_V\Big({\|\vartheta_j^{h\ell}\|}_H{\|\vartheta_j^{h\ell}\|}_V+c_V^{\delta_3}\,{\|\varrho_j^{h\ell}\|}_V{\|\vartheta_j^{h\ell}\|}_H^{1-\delta_3}{\|\vartheta_j^{h\ell}\|}_V^{\delta_3}\Big)\\
        &\le \frac{\gamma_1}{5}{\|\vartheta_j^{h\ell}\|}_V^2+c_6\,{\|\vartheta_j^{h\ell}\|}_H^2+c_7\,{\|\varrho_j^{h\ell}\|}_V^2.
    \end{aligned}
\end{equation}
By \eqref{OperN} we get
\begin{align*}
    {\langle\mathcal N(y_j^{h\ell}-y_j^h),\vartheta_j^{h\ell}\rangle}_{V',V}={\langle\mathcal M(\vartheta_j^{h\ell},\varrho_j^{h\ell})+\mathcal N(\varrho_j^{h\ell}),\vartheta_j^{h\ell}\rangle}_{V',V}.
\end{align*}
From \eqref{Bach-7} it follows that there is a non-negative constant $c_8$ such that
\begin{equation}
    \label{Bach-10}
    \max_{1\le j\le n}\Big({\|\varrho_j^{h\ell}\|}_H^{\delta_3}{\|\varrho_j^{h\ell}\|}_V^{1-\delta_3},{\|\varrho_j^{h\ell}\|}_V\Big)\le c_8.
\end{equation}
Using \eqref{OperN}, \eqref{Bach-10} and Young's inequality we get
\begin{equation}
    \label{Bach-11}
    \begin{aligned}
        \big|{\langle\mathcal N(y_j^{h\ell}-y_j^h),\vartheta_j^{h\ell}\rangle}_{V',V}\big|&\le c_\mathcal Mc_8\Big({\|\vartheta_j^{h\ell}\|}_H{\|\vartheta_j^{h\ell}\|}_V+{\|\varrho_j^{h\ell}\|}_V{\|\vartheta_j^{h\ell}\|}_H^{1-\delta_3}{\|\vartheta_j^{h\ell}\|}_V^{\delta_3}\Big)\\
        &\le \frac{\gamma_1}{5}{\|\vartheta_j^{h\ell}\|}_V^2+c_9\,{\|\vartheta_j^{h\ell}\|}_H^2+c_{10}\,{\|\varrho_j^{h\ell}\|}_V^2
    \end{aligned}
\end{equation}
for non-negative constants $c_9$ and $c_{10}$ depending on $\gamma_1$, $c_1$, $c_\mathcal M$ and $\delta_3$. Combining \eqref{Bach-5}, \eqref{Bach-6}, \eqref{Bach-8}, \eqref{Bach-9}, \eqref{Bach-11} we deduce the estimate
\begin{equation}
    \label{Bach-12a}
    {\|\vartheta_j^{h\ell}\|}_H^2-{\|\vartheta_{j-1}^{h\ell}\|}_H^2+\gamma_1\delta t_j\,{\|\vartheta_j^{h\ell}\|}_V^2\le\delta t_j\big(c_{11}\,{\|\vartheta_j^{h\ell}\|}_H^2+c_{12}\,{\|\varrho_j^{h\ell}\|}_V^2\big)
\end{equation}
with $c_{11}=2 (c_4+c_6+c_9)$ and $c_{12}=2 (c_1+c_2+c_5+c_7+c_{10})$. Consequently,
\begin{align*}
    {\|\vartheta_j^{h\ell}\|}_H^2\le{\|\vartheta_{j-1}^{h\ell}\|}_H^2+\delta t_j\big(c_{11}\,{\|\vartheta_j^{h\ell}\|}_H^2+c_{12}\,{\|\varrho_j^{h\ell}\|}_V^2\big)
\end{align*}
Suppose that
\begin{equation}
    \label{Bach-13}
    \Delta t\le\frac{1}{2c_{11}}.
\end{equation}
From \eqref{Bach-13} we have $1/2 \leq 1-c_{11}\delta t_j < 1$ and
\begin{equation}
    \label{Bach-14}
    \frac{1}{1-c_{11}\delta t_j}\le\frac{1}{1-c_{11}\Delta t}\le1+2c_{11}\Delta t
\end{equation}
for $j=2,\ldots,n$. From \eqref{Bach-13}, \eqref{Bach-14} and by assumption we infer that
\begin{equation}
    \label{Bach-15}
    \begin{aligned}
        {\|\vartheta_j^{h\ell}\|}_H^2&\le(1+2c_{11}\Delta t)\Big({\|\vartheta_{j-1}^{h\ell}\|}_H^2+c_{12}\delta t_j\,{\|\varrho_j^{h\ell}\|}_V^2\Big)\\
        &\le\big(1+2c_{11}\Delta t\big)^{j-1}\,{\|\vartheta_1^{h\ell}\|}_H^2+c_{12}\sum_{l=2}^j\delta t_j\big(1+2c_{11}\Delta t\big)^{j+1-l}\,{\|\varrho_l^{h\ell}\|}_V^2
    \end{aligned}
\end{equation}
for $j=2,\ldots,n$. Due to Assumption~\ref{A100}-3) the quotient $\Delta t/\delta t$ is uniformly bounded by $\zeta$. Thus $0\le2c_{11}\Delta t\le c_{13}\delta t$ with $c_{13}=2c_{11}\zeta$. Then by \eqref{app_bernIneq} with $x=c_{13}(j-1)\delta t$
\begin{align*}
    \big(1+2c_{11}\Delta t\big)^{j-1}\le\bigg(1+\frac{c_{13}(j-1)\delta t}{j-1}\bigg)^{j-1}\le e^{c_{13}(j-1)\delta t},
\end{align*}
where the right-hand side is bounded by $e^{c_{13} T}$ independently of $j$. Consequently we derive from \eqref{Bach-15} and $\vartheta_1^{h\ell}=0$ that
\begin{equation}
    \label{Bach-16}
    {\|\vartheta_j^{h\ell}\|}_H^2\le c_{12}e^{c_{13}(j-1)\delta t}\sum_{l=2}^j\delta t_j\,{\|\varrho_l^{h\ell}\|}_V^2\le c_{14}\sum_{l=1}^n\delta t_j\,{\|\varrho_l^{h\ell}\|}_V^2.
\end{equation}
where we have used that $\vartheta_1^{h\ell}=0$ holds and $c_{14}=c_{12}e^{c_{13}T}$. By summation on $j=2,\ldots,n$ we infer from $\vartheta_1^{h\ell}=0$, $\sum_{j=1}^n\alpha_j^n=T$, \eqref{Bach-12a} and
\begin{align*}
    \delta t_j\le2\alpha^n_j,\quad\delta t_j\ge\frac{\delta t}{\Delta t}\,\Delta t\ge\frac{\Delta t}{\zeta}\ge\frac{\alpha^n_j}{\zeta}\quad\text{for }j=2,\ldots,n
\end{align*}
that
\begin{align*}
    {\|\vartheta_n^{h\ell}\|}_H^2+\frac{\gamma_1}{\zeta}\sum_{j=1}^n\alpha^n_j\,{\|\vartheta_j^{h\ell}\|}_V^2&\le {\|\vartheta_n^{h\ell}\|}_H^2+\gamma_1\sum_{j=1}^n\delta t_j\,{\|\vartheta_j^{h\ell}\|}_V^2\\
    &\le c_{11}\sum_{j=1}^n\delta t_j\,{\|\vartheta_j^{h\ell}\|}_H^2+c_{12}\sum_{j=1}^n\delta t_j\,{\|\varrho_j^{h\ell}\|}_V^2\\
    &\le 2c_{11}\sum_{j=1}^n\alpha^n_j\,{\|\vartheta_j^{h\ell}\|}_H^2+2c_{12}\sum_{j=1}^n\alpha^n_j\,{\|\varrho_j^{h\ell}\|}_V^2\\
    &\le 2c_{11}\sum_{j=1}^n\alpha^n_j\,\bigg(c_{14}\sum_{l=1}^n\delta t_l\,{\|\varrho_l^{h\ell}\|}_V^2\bigg)+2c_{12}\sum_{j=1}^n\alpha^n_j\,{\|\varrho_j^{h\ell}\|}_V^2\\
    &\le 4c_{11}c_{14}T\,\bigg(\sum_{l=1}^n\alpha^n_l\,{\|\varrho_l^{h\ell}\|}_V^2\bigg)+2c_{12}\sum_{j=1}^n\alpha^n_j\,{\|\varrho_j^{h\ell}\|}_V^2\\
    &\le \big(4c_{11}c_{14}T+2c_{12}\big)\sum_{j=1}^n\alpha^n_j\,{\|\varrho_j^{h\ell}\|}_V^2\le c_{15}\sum_{j=1}^n\alpha^n_j\,{\|\varrho_j^{h\ell}\|}_V^2
\end{align*}
with $c_{15}=2(c_{11}c_{14}T+c_{12})$. Thus, we have
\begin{equation}
    \label{Bach-17}
    \sum_{j=1}^n\alpha^n_j\,{\|\vartheta_j^{h\ell}\|}_V^2\le c_{16}\sum_{j=1}^n\alpha^n_j\,{\|\varrho_j^{h\ell}\|}_V^2
\end{equation}
with $c_{16}=c_{15}\zeta/\gamma_1$. Combining \eqref{Bach-1}, \eqref{Bach-17} we get
\begin{align*}
    \sum_{j=1}^n\alpha_j^n\,{\|y_j^h-y_j^{h\ell}\|}_V^2\le 2\sum_{j=1}^n\alpha_j^n\,\left({\|\varrho_j^{h\ell}\|}_V^2+{\|\vartheta_j^{h\ell}\|}_V^2\right)\le C\sum_{j=1}^n\alpha^n_j\,{\|\varrho_j^{h\ell}\|}_V^2
\end{align*}
with $C=2(1+c_{16})$. Now \eqref{NLEq-ApostiE} follows from \eqref{RatePH-VDisc} and \eqref{RateQH-VDisc}.\hfill$\Box$

\medskip\noindent{\bf\em Proof of Theorem~{\em\ref{NonlEq:AposterioriProp}}.} We set $e^{h\ell}_j=y_j^h-y_j^{h\ell}\in V^h$ for $j=1,\ldots,n$ and $\overline\partial e_j^{h\ell}=(e^{h\ell}_j-e^{h\ell}_{j-1})/\delta t_j$ for $j=2,\ldots,n$. Then we conclude that for every $\varphi\in V^h$
\begin{equation}
    \label{Queen-1}
    \begin{aligned}
        &{\langle\overline\partial e^{h\ell}_j,\varphi^h\rangle}_H+a(e^{h\ell}_j,\varphi^h)+{\langle\mathcal Ce^{h\ell}_j,\varphi^h\rangle}_{V',V}\\
        &\quad={\langle g_j(u)+\mathcal N(y_j^{h\ell})-\mathcal N(y_j^h)-\overline\partial y_j^{h\ell}-\mathcal Cy_j^{h\ell}-\mathcal N(y_j^{h\ell}),\varphi^h\rangle}_{V',V}-a(y_j^{h\ell},\varphi^h)\\
        &\quad={\langle\mathcal N(y_j^{h\ell})-\mathcal N(y_j^h),\varphi^h\rangle}_{V',V}+{\langle r_j^{h\ell},\varphi^h\rangle}_{V',V}
    \end{aligned}
\end{equation}
with the residual
\begin{align*}
    r_j^{h\ell}=-\left({\langle\overline\partial y_j^{h\ell}+\mathcal Cy_j^{h\ell}+\mathcal N(y_j^{h\ell})-g_j(u),\cdot\rangle}_{V',V}+a(y_j^{h\ell},\cdot)\right)\in (V^h)'
\end{align*}
for $j=2,\ldots,n$. We argue as for the a-posteriori error analysis carried out in Section~\ref{SIAM-Book:Section3.6.3}. Choosing $\varphi^h=e^{h\ell}_j\in V^h$ in \eqref{Queen-1} we derive from \eqref{OperA&C} and Lemma~\ref{app_lem_innerProd1} that
\begin{equation}
    \label{Queen-2}
    \begin{aligned}
    &2\,{\langle\overline\partial e_j^{h\ell},e_j^{h\ell}\rangle}_H+2a(e_j^{h\ell},e_j^{h\ell})+2\,{\langle\mathcal Ce_j^{h\ell},e_j^{h\ell}\rangle}_{V',V}\\
    &\quad\ge\frac{1}{\delta t_j}\left({\|e_j^{h\ell}\|}_H^2-{\|e_{j-1}^{h\ell}\|}_H^2+{\|e_j^{h\ell}-e_{j-1}^{h\ell}\|}_H^2\right)+2\gamma_1\,{\|e_j^{h\ell}\|}_V^2
    \end{aligned}
\end{equation}
for $j=2,\ldots,n$. Since $\mathcal M$ is bilinear, we have
\begin{align*}
    \mathcal N(y^{h\ell}_j)-\mathcal N(y^h_j)=-\mathcal M(e^{h\ell}_j,y^{h\ell}_j)-\mathcal M(y^h_j,e^{h\ell}_j)
\end{align*}
for $j=2,\ldots,n$ (cf. \eqref{Queen-100}). From \eqref{OperN} we derive that
\begin{align*}
    {\langle\mathcal M(y^h_j,e^{h\ell}_j),e_j^{h\ell}\rangle}_{V',V}=0
\end{align*}
for $j=2,\ldots,n$ and proceeding as in Section~\ref{SIAM-Book:Section3.6.3} we find
\begin{align*}
    \big|{\langle\mathcal M(e^{h\ell}_j,y^{h\ell}_j),e_j^{h\ell}\rangle}_{V',V}\big|\le c_\mathcal M\,{\|e_j^{h\ell}\|}_H{\|y_j^{h\ell}\|}_V{\|e_j^{h\ell}\|}_V
\end{align*}
for $j=2,\ldots,n$. Thus, \eqref{Queen-1} and \eqref{Queen-2} imply that
\begin{align*}
    &{\|e_j^{h\ell}\|}_H^2-{\|e_{j-1}^{h\ell}\|}_H^2+{\|e_j^{h\ell}-e_{j-1}^{h\ell}\|}_H^2+2\gamma_1\delta t_j\,{\|e_j^{h\ell}\|}_V^2\\
    &\quad\le2\delta t_j\left({\langle r^{h\ell}_j,e_j^{h\ell}\rangle}_{V',V}+{\langle\mathcal N(y_j^{h\ell})-\mathcal N(y_j^h),e_j^{h\ell}\rangle}_{V',V}\right)\\
    &\quad\le2\delta t_j\,\left({\|r_j^{h\ell}\|}_{(V^h)'}{\|e_j^{h\ell}\|}_V+c_\mathcal M\,{\|e_j^{h\ell}\|}_H{\|y_j^{h\ell}\|}_V{\|e_j^{h\ell}\|}_V\right)
\end{align*}
for $j=2,\ldots,n$. Consequently, by using Young's inequality we find that
\begin{align*}
    {\|e_j^{h\ell}\|}_H^2-{\|e_{j-1}^{h\ell}\|}_H^2+\gamma_1\delta t_j\,{\|e_j^{h\ell}\|}_V^2\le\frac{\delta t_j}{\gamma_1}\,\bigg({\|r_j^{h\ell}\|}_{(V^h)'}^2+c_\mathcal M^2\,{\|y_j^{h\ell}\|}_V^2{\|e_j^{h\ell}\|}_H^2\bigg).
\end{align*}
for $j=2,\ldots,n$. By assumption there exists a constant $C>0$ which is independent of $h$ and $n$ such that
\begin{align*}
    \max_{1\le j\le n}\,{\|y_j^{h\ell}\|}_V^2\le C^2.
\end{align*}
We set $c_2=c_\mathcal M^2C^2/\gamma_1$. If $\Delta t\le1/(2c_2)$ holds, we derive
\begin{align*}
    \frac{1}{1-c_2\delta t_j}\le1+2c_2\Delta t.
\end{align*}
Consequently, we get upon summation over $j$ 
\begin{align*}
    {\|e_j^{h\ell}\|}_H^2\le\big(1+2c_2\Delta t\big)^{j-1}{\|e_1^{h\ell}\|}_H^2 + \frac{1}{\gamma_1}\sum_{l=2}^j\big(1+2c_2\Delta t\big)^{j+1-l}\delta t_l\,{\|r_l^{h\ell}\|}_{(V^h)'}^2.
\end{align*}
for $j\in\{2,\ldots,n\}$ which gives the claim.\hfill$\Box$ 

\subsection{Proof of Section~\ref{SIAM-Book:Section3.7}}
\label{SIAM-Book:Section3.8.6}

\noindent{\bf\em Proof of Theorem~{\em \ref{Th:EllPro}}.}
Since the bilinear form $a(\bmu;\cdot\,,\cdot)$ is bounded and coercive on $V \times V$ for every parameter $\mu \in \mathscr D$, the existence of a unique solution to \eqref{Eq:ParVarPro} follows directly from the Lax-Milgram lemma;
see \cite{Eva08}, for instance. Next we prove the a-priori estimate \eqref{eq2-7}. For that purpose we take $\varphi =y \in V$ in \eqref{Eq:ParVarPro}. It follows that 
\begin{align*}
    \gamma_1\,{\| y \|}_V^2\le a(\bmu;y,y)={\langle \mathcal F(\bmu),y \rangle}_{V',V}\le{\|\mathcal F(\bmu)\|}_{V'}{\|y\|}_V,
\end{align*}
which gives the claim.\hfill$\Box$


\chapter{Linear-Quadratic Optimal Control Problems}
\label{SIAM-Book:Section4}

\section{The abstract linear-quadratic control problem}
\label{SIAM-Book:Section4.1}
\setcounter{equation}{0}
\setcounter{theorem}{0}

In this section we introduce our abstract optimal control problem, which is a quadratic programming problem in a Hilbert space. Throughout this section we make use of the following assumption; cf. Assumption~\ref{A1}.
\begin{assumption}
    \label{A8}
    Let $T>0$ hold.
    \begin{enumerate}
        \item [\rm 1)] $V$ and $H$ are real, separable Hilbert spaces and $V$ is dense in $H$ with compact embedding.  We identify $H$ with its dual (Hilbert) space $H'$ so that we have the Gelfand triple
        \begin{align*}
            V\hookrightarrow H\simeq H'\hookrightarrow V',
        \end{align*}
        where each embedding is continuous and dense. Here, the embedding $H' \hookrightarrow V'$ is given by $\langle \varphi', \varphi \rangle_{V,V} = \langle \varphi', \varphi \rangle_{H',H}$ for all $\varphi' \in H'$, $\varphi \in V$.
        \item [\rm 2)] The \index{Space!Hilbert!control, $\U$} {\em control space} $\U$ denotes either the Euclidean (Hilbert) space $\mathbb R^\mathsf m$ or the (Hilbert) space $L^2(\mathscr D;\mathbb R^\mathsf m)$, where $\mathscr D\subset\mathbb R^{\mathsf m_\circ}$ is an open, bounded set with $\mathsf m_\circ\in\mathbb N$.
        \item [\rm 3)] For almost all $t\in[0,T]$, let $a(t;\cdot\,,\cdot):V\times V\to\mathbb R$ be a \index{Bilinear form!time-dependent, $a(t;\cdot\,,\cdot)$}time-dependent bilinear form satisfying
        \begin{subequations}
            \label{HI-101}
            \begin{align}
                \label{HI-101-1}
                \big|a(t;\varphi,\psi)\big|&\le\gamma\,{\|\varphi\|}_V{\|\psi\|}_V&&\forall\varphi,\psi\in V \text{ a.e. in } [0,T],\\
                \label{HI-101-2}
                a(t;\varphi,\varphi)&\ge\gamma_1\,{\|\varphi\|}_V^2&&\forall\varphi\in V \text{ a.e. in } [0,T]
            \end{align}
        \end{subequations}
        for constants $\gamma\ge0$ and $\gamma_1>0$ which do not depend on $t$.
        \item [\rm 4)] Assume that $y_\circ\in H$, $\mathcal F\in L^2(0,T;V')$ and that $\mathcal B:\U\to L^2(0,T;V')$ is a continuous, linear and injective \index{Linear operator!control, $\mathcal B$}{\em control operator}.
    \end{enumerate}
\end{assumption}

\begin{remark}
    \label{Remark:HI-30}
    \rm In $\U=\mathbb R^\mathsf m$ we understand the inequality $u\le\tilde u$ componentwise for $u=(u_j)_{1\le j\le\mathsf m}$, $\tilde u=(\tilde u_j)_{1\le j\le\mathsf m}\in\U$, i.e., we have $u_j\le\tilde u_j$ for $j=1,\ldots,\mathsf m$. If $\U=L^2(\mathscr D;\mathbb R^\mathsf m)$ holds, we interpret $u\le\tilde u$ for $u=(u_j)_{1\le j\le\mathsf m}$, $\tilde u=(\tilde u_j)_{1\le j\le\mathsf m}\in\U$ as $u(s)\le\tilde u(s)$ f.a.a. $s\in\mathscr D$.\hfill$\blacksquare$
\end{remark}

\subsection{Problem formulation}
\label{SIAM-Book:Section4.1.1}

First we formulate our state equation which coincides with the linear evolution problem introduced in \eqref{SIAM:Eq3.1.6}. For $u\in\U$ we consider the state equation
\begin{equation}
    \label{GVLuminy:Eq4.1.1}
    \begin{aligned}
        \frac{\mathrm d}{\mathrm dt} {\langle y(t),\varphi \rangle}_H+a(t;y(t),\varphi)&={\langle (\mathcal F+\mathcal Bu)(t),\varphi\rangle}_{V',V}~\forall\varphi\in V\text{ a.e. in }(0,T],\\
        {\langle y(0),\varphi\rangle}_H&={\langle y_\circ,\varphi\rangle}_H\hspace{20mm}\forall\varphi\in H.
    \end{aligned}
\end{equation}

\begin{remark}
    \label{Remark:OperatorS}
    \rm To simplify the presentation we have replaced \eqref{SIAM:Eq3.1.1-2} by the stronger assumption \eqref{HI-101-2} in Assumption~\ref{A8}. However, in some applications \eqref{HI-101-2} can not be guaranteed, but the bilinear form $a(t;\cdot\,,\cdot)$ satisfies weaker estimate \eqref{SIAM:Eq3.1.1-2}; cf. Section~\ref{SIAM-Book:Section3.1}. Then, we define the \index{Bilinear form!time-dependent, $a(t;\cdot\,,\cdot)$}time-dependent bilinear form
    \begin{align*}
        \tilde a(t;\varphi,\psi)=a(t;\varphi,\psi)+\gamma_2\,{\langle\varphi,\psi\rangle}_H\quad\text{for }\varphi,\psi\in V\text{ a.e. in }[0,T].
    \end{align*}
    Clearly, $\tilde a$ satisfies \eqref{SIAM:Eq3.1.1-2}. By \eqref{Poincare} we derive that
    \begin{align*}
        \big|\tilde a(t;\varphi,\psi)\big|\le \gamma\,{\|\varphi\|}_V{\|\psi\|}_V+\gamma_2\,{\|\varphi\|}_H{\|\psi\|}_H\le\tilde\gamma\,{\|\varphi\|}_V{\|\psi\|}_V\quad\text{for all }\varphi,\psi\in V
    \end{align*}
    with $\tilde\gamma=\gamma+\gamma_2c_V^2>0$. Thus, the bilinear form $\tilde a(t;\cdot\,,\cdot)$ is also bounded. Let 
    \begin{align*}
        \tilde{\mathcal F}(t)=e^{-\gamma_2t}\mathcal F(t),\quad(\tilde{\mathcal B}u)(t)=e^{-\gamma_2t}(\mathcal Bu)(t) \quad \text{f.a.a. } t \in [0,T].
    \end{align*}
    Then it clearly holds $\tilde{\mathcal F}\in L^2(0,T;V')$, $\tilde{\mathcal B}\in\mathscr L(\U;L^2(0,T;V'))$. Instead of \eqref{GVLuminy:Eq4.1.1} we consider the evolution problem
    \begin{equation}
        \label{GVLuminy:Eq4.1.1tilde}
        \begin{aligned}
            \frac{\mathrm d}{\mathrm dt} {\langle \tilde y(t),\varphi \rangle}_H+a(t;\tilde y(t),\varphi)&={\langle (\tilde{\mathcal F}+\tilde{\mathcal B}u)(t),\varphi\rangle}_{V',V}&&\text{for all }\varphi\in V\text{ a.e. in }(0,T],\\
            {\langle \tilde y(0),\varphi\rangle}_H&={\langle y_\circ,\varphi\rangle}_H&&\text{for all }\varphi\in H.
        \end{aligned}
    \end{equation}
    for any $u\in\U$ and $y_\circ\in H$. Due to Theorem~\ref{SIAM:Theorem3.1.1} there exists a unique solution $\tilde y\in\Y$ to \eqref{GVLuminy:Eq4.1.1tilde}. Then we define 
    \begin{align*}
        y(t)=e^{\gamma_2t}\tilde y(t) \quad \text{f.a.a. } t \in [0,T].
    \end{align*}
    Now we find for all $\varphi\in V$
    \begin{align*}
        &\frac{\mathrm d}{\mathrm dt} {\langle y(t),\varphi \rangle}_H+a(t;y(t),\varphi)=\gamma_2\,{\langle y(t),\varphi\rangle}_H+e^{\gamma_2t}\frac{\mathrm d}{\mathrm dt} {\langle\tilde y(t),\varphi \rangle}_H+a(t;y(t),\varphi)\\
        &\quad=\tilde a(t;y(t),\varphi)+e^{\gamma_2t}\Big({\langle (\tilde f+\tilde{\mathcal B}u)(t),\varphi\rangle}_{V',V}-\tilde a(t;\tilde y(t),\varphi)\big)={\langle (\mathcal F+\mathcal Bu)(t),\varphi\rangle}_{V',V}.
    \end{align*}
    As $y(0)=\tilde y(0)$ holds, $y$ solves \eqref{SIAM:Eq3.1.6}.\hfill$\blacksquare$
\end{remark}

Since Assumption~\ref{A8} holds, Theorem~\ref{SIAM:Theorem3.1.1} ensures the existence of a unique solution $y$ belonging to the \index{Space!state, $\Y$}{\em state space} $\Y=W(0,T)$ to \eqref{GVLuminy:Eq4.1.1}. In Corollary~\ref{Corollary:HI-20} we have introduced the particular solution $\hat y\in\Y$ as well as the linear, bounded solution operator $\mathcal S$. Then the solution to \eqref{GVLuminy:Eq4.1.1} can be expressed as $y=\hat y+\mathcal Su$. 

\begin{remark}
    \label{Remark:HI-33}
    \rm Since the operator $\mathcal{B}$ is injective, we can conclude that the operator $\mathcal{S}$ is injective as well: Let $u \in \U$ be given with $\mathcal S u = 0$, then it immediately follows from \eqref{eq:solutionOperator_a} that $\langle (\mathcal B u)(t), \varphi \rangle_{V' \times V}=0$ for all $\varphi \in V$. This is equivalent to $\mathcal B u = 0$ in $L^2(0,T;V')$, which yields $u=0$ due to the fact that $\mathcal B$ is injective. Therefore, $\mathcal S$ is injective. \hfill$\blacksquare$
\end{remark}

\begin{definition}
    \label{Definition:HI-100}
    We define the Hilbert space
    \begin{align*}
        \mathscr X=\Y\times\U
    \end{align*}
    endowed with the natural product topology, i.e., with the inner product
    \begin{align*}
        {\langle x,\tilde x\rangle}_\X={\langle y,\tilde y\rangle}_\Y+{\langle u,\tilde u\rangle}_\U\quad \text{for }x=(y,u),~\tilde x=(\tilde y,\tilde u)\in\X
    \end{align*}
    and the norm $\|x\|_\X=(\|y\|_\Y^2+\|u\|_\U^2)^{1/2}$ for $x=(y,u)\in\X$.
\end{definition}

In a next step we introduce the quadratic objective function. For given desired states $(\ydQ,\ydT) \in L^2(0,T;H) \times H$, a nominal control $\un\in\U$ and non-negative weighting parameters $\sigma_1$, $\sigma_2$ and $\sigma$ the \index{Cost functional!quadratic}{\em quadratic cost functional} $J:\X\to\mathbb R$ is given by
\begin{equation}
    \label{GVLuminy:Eq4.1.5}
    J(x)=\frac{\sigma_1}{2}\int_0^T{\|y(t)-\ydQ(t)\|}_H^2\,\mathrm dt+\frac{\sigma_2}{2}\,{\|y(T)-\ydT\|}_H^2+\frac{\sigma}{2}\,{\|u-\un\|}_\U^2
\end{equation}
for $x=(y,u)\in\X$. Then given the control bounds $u_\mathsf a,\,u_\mathsf b \in \U$ with $u_a \leq u_b$ we consider the quadratic programming problem
\begin{equation}
    \label{GVLuminy:Eq4.1.6a}
    \tag{\bf P}
    \min J(x)\quad\mbox{subject to (s.t.)}\quad x = (y,u) \in \X,~y = \hat y+\mathcal Su \text{ and } u_\mathsf a \leq u \leq u_\mathsf b.
\end{equation}

\begin{definition}
    \label{Definition:HI-1001}
    The \index{Set!of admissible controls, $\Uad$}{\em set of admissible controls} is given as
    \begin{align*}
        \Uad=\big\{u\in\U\,\big|\,u_\mathsf a\le u\le u_\mathsf b\text{ in }\U\big\}.
    \end{align*}
    By $\Xad\subset\X$ we denote the admissible set for the optimization problem as
    \begin{align*}
        \Xad=\big\{(y,u)\in\X\,\big|\,y=\hat y+\mathcal Su\text{ and }u\in\Uad\big\}.
    \end{align*}
\end{definition}

With this definition we can write the \index{Problem!quadratic programming}{\em quadratic programming problem} in the form
\begin{equation}
    \label{GVLuminy:Eq4.1.6}
    \tag{\bf P}
    \min J(x)\quad\mbox{subject to (s.t.)}\quad x\in\Xad.
\end{equation}

\begin{remark}
    \label{Remark:HI_101}
    \rm From $x=(y,u)\in\Xad$ we infer that $y=\hat y+\mathcal Su$ holds. Hence, $y$ is a dependent variable and \eqref{GVLuminy:Eq4.1.6} is called a \index{Problem!optimal control!linear-quadratic}{\em linear-quadratic optimal control problem}. Moreover, we call $y$ the {\em state variable} and $u$ the {\em control variable}.\hfill$\blacksquare$
\end{remark}

We summarize the framework of the optimal control problem in the following assumption.

\begin{assumption}
    \label{A9}
    \begin{enumerate}
        \item [\rm 1)] The pair $(\ydQ,\ydT)$ of desired states belongs to $L^2(0,T;H)\times H$ and $\un\in\U$ is a nominal control. Moreover, $\sigma_1$, $\sigma_2$ and $\sigma$ are non-negative weighting parameters.
        \item [\rm 2)] For $u_\mathsf a,u_\mathsf b\in\U$ it holds $u_\mathsf a \le u_\mathsf b$.
        \item [\rm 3)] If $\Uad$ is unbounded (i.e., at least one of the components of $u_\mathsf a$ or $u_\mathsf b$ is $-\infty$ or $\infty$, respectively) we have $\sigma>0$.
    \end{enumerate}
\end{assumption}
 
\begin{remark}
    \label{Remark:HI-100}
    \rm
    \begin{enumerate}
        \item [1)] The set $\Uad$ is closed and convex. Since $\mathcal S\in\mathscr L(\U,\Y)$ holds, $\Xad$ is closed and convex as well.
        \item [2)] From Assumption~\ref{A9}-2) and Theorem~\ref{SIAM:Theorem3.1.1} we conclude that the set $\Xad$ is nonempty.
        \item [3)] More general quadratic cost functionals can be treated as well. For instance, the presentation can be extended to the following \index{Cost functional!quadratic}quadratic cost functional:
        \begin{align*}
            J(x)=\frac{\sigma_1}{2}\,{\|\mathcal C_1y-\ydQ\|}_{\W_1}^2+\frac{\sigma_2}{2}\, {\|\mathcal C_2(y(T))-\ydT\|}_{\W_2}^2+\frac{\sigma}{2}\,{\|u-\un\|}_\U^2
        \end{align*}
        for $x=(y,u)\in\X$, where $\W_1$, $\W_2$ are Hilbert spaces, $\mathcal C_1:\Y\to\W_1$, $\mathcal C_2:H\to\W_2$ are bounded, \index{Linear operator!observation} linear {\em observation operators}, $(\ydQ,\ydT)\in \W_1\times\W_2$, $\un\in\U$ are given data functions, $\sigma_1$, $\sigma_2\geq0$ and $\sigma>0$. \hfill$\blacksquare$
    \end{enumerate}
\end{remark}

\begin{definition}
    \label{Definition:HI-101}
    A point $\bar x=(\bar y,\bar u)\in\X$ is called a \index{Optimal solution!global}{\em globally optimal solution} to \eqref{GVLuminy:Eq4.1.6} if $\bar x\in\Xad$ and $J(\bar x)\le J(x)$ for all $x\in\Xad$ hold. We say that $\bar x=(\bar y,\bar u)\in\X$ is a \index{Optimal solution!local}{\em locally optimal solution} to \eqref{GVLuminy:Eq4.1.6} if there exists an open set $\mathscr N(\bar x)\subset\X$ containing $\bar x$ such that $\bar x\in\Xad$ and $J(\bar x)\le J(x)$ for all $x\in\Xad\cap\mathscr N(\bar x)$ are satisfied. 
\end{definition}

The proof of the next lemma is given in Section~\ref{SIAM-Book:Section4.7.1}.

\begin{lemma}
    \label{Lemma:HI-100}
    \begin{enumerate}
        \item [\rm 1)] Let Assumption~{\rm\ref{A9}-1)} hold. Then $J$ is convex.
        \item [\rm 2)] Suppose that Assumptions~{\rm\ref{A8}} and {\rm\ref{A9}} are satisfied. Then every locally optimal solution to \eqref{GVLuminy:Eq4.1.6} is a globally optimal one.
    \end{enumerate}
\end{lemma}

Utilizing the relationship $y=\hat y+\mathcal Su$ we define a so-called \index{Cost functional!reduced}{\em reduced cost functional} $\hat J:\U\to\mathbb R$ by
\begin{align*}
    \hat J(u)=J(\hat y+\mathcal Su,u)\quad \text{for }u\in\U.
\end{align*}
Then we consider the \index{Problem!optimal control!reduced}{\em reduced optimal control problem}:
\begin{equation}
    \label{GVLuminy:Eq4.1.9}
    \tag{$\mathbf{\hat P}$}
    \min\hat J(u)\quad\text{s.t.}\quad u\in\Uad.
\end{equation}

The proof of the next lemma is straightforward. Therefore, we omit the proof here.

\begin{lemma}
    \label{Lemma:HI-101}
    Let Assumptions~{\rm\ref{A8}} and {\rm\ref{A9}} hold. Then  if $\bar u$ is the optimal solution to \eqref{GVLuminy:Eq4.1.9}, then $\bar x=(\hat y+\mathcal S\bar u,\bar u)$ is the optimal solution to \eqref{GVLuminy:Eq4.1.6}. On the other hand, if $\bar x=(\bar y,\bar u)$ is the solution to \eqref{GVLuminy:Eq4.1.6}, then $\bar u$ solves \eqref{GVLuminy:Eq4.1.9}.
\end{lemma}

\begin{example}
    \label{GVLuminy:Example4.1.1}
    \rm An example for \eqref{GVLuminy:Eq4.1.6} is the one given in Sections~\ref{SIAM-Book:Section1.2} and \ref{SIAM-Book:Section3.1}. We discuss the presented theory for this application. Let $\Omega\subset\mathbb R^\mathfrak n$, $\mathfrak n\in\{1,2,3\}$, be an open and bounded domain with Lipschitz-continuous boundary $\Gamma=\partial \Omega$. For $T>0$ we set $Q=(0,T) \times \Omega$ and $\Sigma=(0,T) \times \Gamma$. We choose $H=L^2(\Omega)$ and $V=H_0^1(\Omega)$ endowed with the usual inner products
    \begin{align*}
        {\langle \varphi,\psi \rangle}_H=\int_\Omega \varphi\psi\,\mathrm d\boldsymbol x, \quad {\langle \varphi,\psi \rangle}_V=\int_\Omega \varphi\psi+\nabla \varphi \cdot \nabla \psi\,\mathrm d\boldsymbol x 
    \end{align*}
    and their induced norms, respectively. We choose $\mathscr D=(0,T)$ and $\U=L^2(0,T;\mathbb R^\mathsf m)$. Let $\chi_i \in H$, $1 \le i \le\mathsf m$, denote given control shape functions. Then  for given control $u \in\U$, initial condition $y_\circ \in H$, diffusion coefficient $c(\bx)\in L^\infty(\Omega)$ with $c(\bx)\geq c_0>0$ for all $\bx\in\Omega$, reaction coefficient $d(\bx)\in L^\infty(\Omega)$ with $d(\bx)\geq d_0>0$ for all $\bx\in\Omega$, non-negative $q\in L^\infty(\Gamma)$, inhomogeneity $f \in L^2(0,T;H)$ and boundary inhomogeneity $g \in L^2(0,T;L^2(\Gamma))$ we consider the linear parabolic \index{Equation!partial differential!diffusion-reaction}\index{Equation!partial differential!parabolic}{\em diffusion-reaction equation} \eqref{SIAM:EqMotPDE1}
    \begin{equation*}
        \begin{aligned}
            y_t(t,\bx)-\kappa\Delta y(t,\bx))+\bv(\bx) \cdot \nabla y(t,\bx)&=f(t,\bx)+\sum_{i=1}^\md u^\md_i(t)\chi_i(\bx),&&(t,\bx)\in Q,\\
            \kappa\frac{\partial y}{\partial\bn}(t,\bs)+q(\bs)y(t,\bs)&=g(t,\bs)+\sum_{j=1}^\mb u^{\mathsf b}_j(t)\xi_{j}(\bs),&&(t,\bs)\in\Sigma,\\
            y(0,\bx)&=y_\circ(0),&&\bx\in\Omega.
        \end{aligned}
    \end{equation*}
    We introduce the time-independent, symmetric bilinear form
    \begin{align*}
        a(\varphi,\psi)=\kappa\int_\Omega \nabla \varphi \cdot \nabla \psi+\big(\bv\cdot\nabla\varphi\big)\psi \, \mathrm d\boldsymbol x + \int_\Gamma q \varphi \psi \, \mathrm d s \quad \text{ for } \varphi, \psi \in V
    \end{align*}
    the bounded time-dependent linear operator $\mathcal F(t): V\to V'$ such that
    \begin{align*}
        {\langle \mathcal F(t),\varphi\rangle}_{V',V} = \int_\Omega f(t)\varphi \, \mathrm d\bx + \int_\Gamma g(t) \varphi \mathrm d s \quad \text{for } \varphi\in V \text{ a.e. in } (0,T)
    \end{align*}
    and the bounded, linear and injective operator $\mathcal B:\U \to L^2(0,T;H) \hookrightarrow L^2(0,T;V')$ as
    \begin{align*}
        (\mathcal B u)(t,\boldsymbol x)=\sum_{i=1}^\md u^{\mathsf d}_i(t)\chi_i(\boldsymbol x)\,+\,\sum_{j=1}^\mb u^{\mathsf b}_j(t)\xi_{j}(\bx)\quad \text{for } (t,\boldsymbol x) \in Q\text{ a.e. and } u \in\U.
    \end{align*}
    Applying Cauchy-Schwarz inequality and Theorem~\ref{TraceTh}, one can derive the following estimate:
    \begin{align*}
        \big| a(\varphi,\psi)\big| & \leq {\|c\|}_{L^\infty(\Omega)}{\|\varphi\|}_V{\|\psi\|}_V+{\|d\|}_{L^\infty(\Omega)}{\|\varphi\|}_H{\|\psi\|}_H+\gamma_\Gamma^2{\|q\|}_{L^\infty(\Gamma)}{\|\varphi\|}_V{\|\psi\|}_V\leq \gamma\,{\|\varphi\|}_V{\|\psi\|}_V 
    \end{align*}
    with $\gamma=\|c\|_{L^\infty(\Omega)}+\|d\|_{L^\infty(\Omega)}+\gamma_\Gamma^2\|q\|_{L^\infty(\Gamma)}$. Moreover, we have
    \begin{align*}
        a(\varphi,\varphi) \geq \gamma_1 \|\varphi\|_V^2
    \end{align*}
    with $\gamma_1=\min(c_0,d_0)$. Hence, $a(\cdot,\cdot)$ satisfies \eqref{HI-101-1} and \eqref{HI-101-2}. From Remark~\ref{Remark:OperatorS} it follows that the weak formulation of \eqref{SIAM:EqMotPDE1} can be expressed in the form \eqref{SIAM:Eq3.1.6}. Moreover, the unique weak solution to \eqref{SIAM:EqMotPDE1} belongs to the space $L^\infty(0,T;V)$ provided $y_\circ \in V$ holds.\hfill$\blacklozenge$
\end{example}

\subsection{Existence of a unique optimal solution}
\label{SIAM-Book:Section4.1.2}

Next, we will review an existence result for quadratic optimization problems in Hilbert spaces. While we have previously denoted the cost function in this chapter with $\hat J$, the following result from \cite[pp.~50-51]{Tro10} will use a more general cost function $\mathcal J$.

\begin{theorem}
    \label{GVLuminy:Theorem4.2.1}
    Suppose that $\tU$ and $\Ha$ are given Hilbert spaces with norms $\|\cdot\|_\tU$ and $\|\cdot\|_\Ha$, respectively. Furthermore, let $\tUad\subset\tU$ be non-empty, closed, convex and $\yd\in \Ha$, $\tilde u^\mathsf n\in\tU$, $\sigma\geq 0$. Moreover, let either $\tUad$ be bounded or $\sigma>0$ hold. The mapping $\mathcal G\,:\tU\to \Ha$ is assumed to be linear and continuous. Then there exists an optimal control $\bar u$ solving
    \begin{equation}
        \label{GVLuminy:Eq4.2.1}
        \min\mathcal J(u):=\frac{1}{2}\,{\|\mathcal Gu-\yd\|}_\Ha^2+\frac{\sigma}{2}\,{\|u-\tilde u^\mathsf n\|}_\tU^2\quad\text{s.t.}\quad u\in\tUad.
    \end{equation}
    If $\sigma>0$ holds or if $\mathcal G$ is injective, then $\bar u$ is uniquely determined. 
\end{theorem}

\begin{remark}
    \label{GVLuminy:Remark4.2.1}
    \rm In the proof of Theorem \ref{GVLuminy:Theorem4.2.1} it is only used that $\mathcal J$ is continuous and convex and that, in the case of $\tUad$ being unbounded, it holds $\lim_{{\| u \|}_\tU \to \infty} \mathcal{J}(u) = \infty$. Therefore, the existence of an optimal control follows for general cost functionals $\mathcal J:\tU \to \mathbb R$ with a Hilbert space $\tU$ fulfilling these properties.\hfill$\blacksquare$
\end{remark}

Next we can use Theorem~\ref{GVLuminy:Theorem4.2.1} to obtain an existence result for the optimal control problem \eqref{GVLuminy:Eq4.1.9}, which implies the existence of an optimal solution to \eqref{GVLuminy:Eq4.1.6}. For a proof we refer the reader to Section~\ref{SIAM-Book:Section4.7.1}.

\begin{theorem}
    \label{GVLuminy:Theorem4.2.2}
    Let Assumptions~{\rm\ref{A8}} and {\rm\ref{A9}} be satisfied. Then \eqref{GVLuminy:Eq4.1.9} possesses a solution $\bar u$. If $\sigma_1 > 0$ or $\sigma > 0$ the solution is unique.
\end{theorem}

\subsection{First-order sufficient optimality conditions}
\label{SIAM-Book:Section4.1.3}

Existence of a unique solution to \eqref{GVLuminy:Eq4.2.1} has been investigated in Section~\ref{SIAM-Book:Section4.1.2}. In this section we characterize the solution to \eqref{GVLuminy:Eq4.2.1} by first-order sufficient optimality conditions. To derive first-order conditions we require the notion of G\^ateaux- and Fr\'echet-derivatives in function spaces; see Definition~\ref{Definition:GFDer}.

\begin{remark}
    \label{GVLuminy:Remark4.3.1}
    \rm Let $\Ha$ be a real Hilbert space and $F:\Ha\to\mathbb R$ be G\^{a}teaux-differentiable at $u\in\mathscr H$. Then its G\^ateaux derivative $F'(u)$ at $u$ belongs to $\Ha'=\mathscr L(\Ha,\mathbb R)$. Due to the Riesz theorem there exists a unique element $\nabla F(u)\in\mathscr H$ satisfying
    \begin{align*}
        {\langle\nabla F(u),v\rangle}_\Ha={\langle F'(u),v\rangle}_{\Ha',\Ha}\quad\text{for all }v\in\Ha.
    \end{align*}
    We call $\nabla F(u)$ the {\em (G\^ateaux) gradient} of $F$ at $u$.\hfill$\blacksquare$
\end{remark}

For a proof of the next lemma we refer to Section~\ref{SIAM-Book:Section4.7.1}.

\begin{lemma}
    \label{Lemma:GradRepr}
    Let all hypotheses of Theorem~{\rm\ref{GVLuminy:Theorem4.2.1}} hold. Then $\mathcal J$ defined in \eqref{GVLuminy:Eq4.2.1} is Gateaux-differentiable in $\tU$. Its Gateaux derivative is given as
    \begin{equation}
        \label{HI-1}
        \mathcal J'(u)={\langle\mathcal Gu-\yd,\mathcal G(\cdot)\rangle}_\mathscr H+\sigma\,{\langle u-\tilde u^\mathsf n,\cdot\rangle}_\tU\quad \text{for }u\in\tU.
    \end{equation}
\end{lemma}

\begin{remark}
\rm
The gradient $\nabla \mathcal J$ of $\mathcal J$ introduced in Remark~\ref{GVLuminy:Remark4.3.1} is given as
\begin{equation}
\label{HI-2}
\nabla\mathcal J(u)=\mathcal G^\star(\mathcal Gu-\yd)+\sigma\,(u-\tilde u^\mathsf n)\quad \text{for }u\in\tU
\end{equation}
with the adjoint operator $\mathcal G^\star\in\mathscr L(\Ha,\tU)$.\hfill$\blacksquare$
\end{remark}

In the next theorem, which is proved in Section~\ref{SIAM-Book:Section4.7.1}, we formulate the first-order necessary optimality conditions for \eqref{GVLuminy:Eq4.2.1}. An essential assumption is that the set $\tUad$ of admissible controls is non-empty and convex.

\begin{theorem}[First-order necessary optimality conditions]
    \label{GVLuminy:Theorem4.3.1}
    Let\index{Optimality conditions!first-irder!necessary} all assumptions of Theorem~{\rm\ref{GVLuminy:Theorem4.2.1}} hold. Suppose that $\bar u\in\tUad$ is a solution to \eqref{GVLuminy:Eq4.2.1}. Then the following variational inequality holds
    \begin{equation}
        \label{GVLuminy:Eq4.3.6} 
            {\langle \nabla \mathcal J(\bar u),u-\bar u\rangle}_\tU\geq0\quad\text{for all }u \in \tUad,
    \end{equation}
    where the gradient of $\mathcal J$ is given by \eqref{HI-2}. If $\bar u\in\tUad$ solves \eqref{GVLuminy:Eq4.3.6}, then $\bar u$ is a solution to \eqref{GVLuminy:Eq4.2.1}. 
\end{theorem}

\begin{remark}
    \label{Remark:HI-31}
    \rm Combining \eqref{HI-2} and \eqref{GVLuminy:Eq4.3.6} we find that the variational inequality
    \begin{equation}
        \label{GVLuminy:Eq4.3.8}
        {\langle \mathcal G\bar u-\yd,\mathcal G(u-\bar u)\rangle}_\Ha+{\langle\sigma\,(\bar{u}-\un),u-\bar u \rangle}_\U\ge 0\quad\text{for all }u\in\tUad
    \end{equation}
    is a first-order sufficient optimality condition for \eqref{GVLuminy:Eq4.2.1}.\hfill$\blacksquare$
\end{remark}

From now on we assume the situation given in the proof of Theorem~\ref{GVLuminy:Theorem4.2.2}, i.e. the operator $\mathcal{G}$ has the form like in \eqref{GVLuminy:Eq4.2.2} and the cost function is specifically given by $\hat J$ instead of the general $\mathcal J$ from \eqref{GVLuminy:Theorem4.2.2}. Our aim is to find a representation of the gradient $\nabla \mathcal{J}$ in \eqref{HI-2} which can be evaluated numerically. To do this, we need to find an equivalent, more practical representation of the adjoint operator $\mathcal G^\star$. This task will be completed in Lemma \ref{Lemma:HI-30}, but first we require some preparations in the form of the the so-called adjoint equation. A proof of the next lemma is given in Section~\ref{SIAM-Book:Section4.7.1}.

\begin{lemma}
    \label{Lemma:HI-31}
    Suppose that Assumptions~{\rm\ref{A8}} and {\rm\ref{A9}} hold. Let us define the linear operator $\mathcal A:\U\to\Y$ as follows: For given $u\in\U$ the function $p=\mathcal Au\in\Y$ is the unique solution to
    \begin{equation}
        \label{GVLuminy:Eq4.3.12}
        \begin{aligned}
            -\frac{\mathrm d}{\mathrm dt}\,{\langle p(t),\varphi\rangle}_H+a(t;\varphi,p(t))&=-\sigma_1\,{\langle (\mathcal Su)(t),\varphi\rangle}_H&&\text{for all }\varphi\in V\text{ a.e. in }[0,T),\\
            p(T)&=-\sigma_2\,(\mathcal Su)(T)&&\text{in } H.
        \end{aligned}
    \end{equation}
    \begin{enumerate}
        \item [\em 1)] The operator $\mathcal A$ is well-defined and bounded. 
        \item [\em 2)] We set $y=\mathcal Su \in W_0(0,T)$, $w=\mathcal Sv\in W_0(0,T)$, and $p=\mathcal Av\in\Y$. Then 
        \begin{align*}
            \int_0^T {\langle (\mathcal Bu)(t),p(t)\rangle}_{V',V}\,\mathrm dt=-\int_0^T\sigma_1\,{\langle w(t),y(t)\rangle}_H\,\mathrm dt-\sigma_2\,{\langle w(T),y(T)\rangle}_H.
        \end{align*}
    \end{enumerate}
\end{lemma}

\begin{remark}
    \label{Remark:HI-105}
    \rm We define $\hat p\in\Y$ as the unique solution to
    \begin{equation}
        \label{GVLuminy:Eq4.3.15}
        \begin{aligned}
            -\frac{\mathrm d}{\mathrm dt}\,{\langle \hat p(t),\varphi\rangle}_H+a(t;\varphi,\hat p(t))&=\sigma_1\,{\langle\ydQ(t)-\hat y(t),\varphi\rangle}_H&&\text{for all }\varphi\in V\text{ a.e. in }[0,T),\\
            \hat p(T)&=\sigma_2\,(\ydT-\hat y(T))&&\text{in } H.
        \end{aligned}
    \end{equation}
    Then for every $u\in\U$ the function $p=\hat p+\mathcal Au$ is the unique solution to
    \begin{subequations}
        \label{GVLuminy:Eq4.3.16}
        \begin{align}
            \label{GVLuminy:Eq4.3.16a}
            -\frac{\mathrm d}{\mathrm dt}\,{\langle p(t),\varphi\rangle}_H+a(t;\varphi,p(t))&=\sigma_1\,{\langle\ydQ(t)-y(t),\varphi\rangle}_H&&\text{for all }\varphi\in V\text{ a.e. in }[0,T),\\
            \label{GVLuminy:Eq4.3.16b}
            p(T)&=\sigma_2\,(\ydT-y(T))&&\text{in } H
        \end{align}
    \end{subequations}
    with $y=\hat y+\mathcal Su$.\hfill$\blacksquare$
\end{remark}

Moreover, we have the following result which yields a computable formula for the gradient of the reduced cost functional. For the proof we refer the reader to Section~\ref{SIAM-Book:Section4.7.1}.

\begin{lemma}
    \label{Lemma:HI-30}
    Let Assumptions~{\rm\ref{A8}} and {\rm\ref{A9}} hold. Let us define the two linear, bounded operators $\Theta:W_0(0,T)\to W_0(0,T)'$ and $\Xi:L^2(0,T;H)\times H \to W_0(0,T)'$ by
    \begin{equation}
        \label{GVLuminy:Eq4.3.10}
        \begin{aligned}
            {\langle\Theta\varphi,\phi\rangle}_{W_0(0,T)',W_0(0,T)}&=\int_0^T\sigma_1\,{\langle\varphi(t),\phi(t)\rangle}_H\,\mathrm dt+\sigma_2\,{\langle\varphi(T),\phi(T)\rangle}_H,\\
            {\langle\Xi z,\phi\rangle}_{W_0(0,T)',W_0(0,T)}&=\int_0^T\sigma_1\,{\langle z_1(t),\phi(t)\rangle}_H\,\mathrm dt+\sigma_2\,{\langle z_2,\phi(T)\rangle}_H
        \end{aligned}
    \end{equation}
    for $\varphi,\,\phi\in W_0(0,T)$ and $z=(z_1,z_2)\in L^2(0,T;H)\times H$.
    \begin{enumerate}
        \item [\em 1)] We have
        \begin{equation}
            \label{GVLuminy:Eq4.3.11}
            \mathcal G^\star(\mathcal G u-\yd)=\mathcal S'\Theta\mathcal S u-\mathcal S'\Xi(y_1^\mathsf d-\hat y,y_2^\mathsf d-\hat y(T))
        \end{equation}
        for all $u \in \U$, where $\mathcal S':\Y'\to\U'\simeq\U$ is the dual operator of $\mathcal S$ and $\mathcal G^\star \in \mathscr L(\mathscr H, \tilde \U)$ is the Hilbert adjoint of $\mathcal G$ defined in \eqref{GVLuminy:Eq4.2.2}. 
        \item [\em 2)] Then 
        \begin{align}
            \label{eq:BPrimeA}
            \mathcal B'\mathcal A=-\mathcal S'\Theta\mathcal S\text{ in }\mathscr L(\U)\quad\text{and}\quad\mathcal B'\hat p=\mathcal S'\Xi(\ydQ-\hat y,\ydT-\hat y(T)),
        \end{align}
        where the linear and bounded operator $\mathcal A$ has been defined in Lemma~{\em\ref{Lemma:HI-31}}.
    \end{enumerate}
\end{lemma}

\begin{remark}
    \label{rem:concreteGradient}
    \rm Combining the previous lemmata, we can now express the gradient of the function $\hat J$ in \eqref{GVLuminy:Eq4.2.1} in a more explicit way. In fact, applying \eqref{HI-2}, \eqref{GVLuminy:Eq4.3.11} and \eqref{eq:BPrimeA} we get
    \begin{equation}
        \label{GradientJhat}
        \begin{aligned}
            \nabla \hat J(u)&=\mathcal G^\star(\mathcal Gu-\yd)+\sigma (u-\tilde u^\mathsf n) \\
            &=\mathcal S' \Theta \mathcal S u - \mathcal S' \Xi (y_1^\mathsf d-\hat y, y_2^\mathsf d - \hat y(T)) + \sigma(u - \tilde u^\mathsf n) \\
            &=- \mathcal B' \mathcal A u - \mathcal B' \hat p + \sigma (u-\tilde u^\mathsf n)=-\mathcal B' p + \sigma(u-\tilde u^\mathsf n),
        \end{aligned}
    \end{equation}
    where we have defined the adjoint variable $p = \hat p + \mathcal A u$ that solves the adjoint equation \eqref{GVLuminy:Eq4.3.16}.\hfill$\blacksquare$
\end{remark}

The following result is just Theorem~\ref{GVLuminy:Theorem4.3.1}, where we have explicitly used the representation of the gradient $\nabla \hat J$ from Remark \ref{rem:concreteGradient} for the optimal variable $u = \bar u$.

\begin{theorem}
\label{GVLuminy:Theorem4.3.2}
Suppose that Assumptions~{\rm\ref{A8}} and {\rm\ref{A9}} hold. Then $\bar x=(\bar y,\bar u)$ is a solution to \eqref{GVLuminy:Eq4.1.6} if and only if $\bar x$ together with the adjoint variable $\bar p$ satisfy the first-order optimality system
\begin{subequations}
    \label{GVLuminy:Eq4.3.20}
    \begin{align}
        \label{GVLuminy:Eq4.3.20-1}
        &\bar y=\hat y+\mathcal S\bar u,\qquad\bar p=\hat p+\mathcal A\bar u,\qquad \bar u\in\Uad,\\
        \label{GVLuminy:Eq4.3.20-2}
        &{\langle\sigma(\bar u-\un)-\mathcal B'\bar p,u-\bar u\rangle}_\U\ge 0\quad \forall u\in\Uad.
    \end{align}
\end{subequations}
\end{theorem}

It turns out that \eqref{GVLuminy:Eq4.3.20} is only one of several equivalent optimality conditions. 
If we apply \cite[Lemma 1.12]{HPUU09} to our problem, we achieve the following result.

\begin{proposition}
    \label{Pro:HI-1}
    Let Assumptions~{\rm\ref{A8}} and {\rm\ref{A9}} hold. Due to \eqref{GradientJhat} the gradient $\nabla\hat J$ of the reduced cost functional $\hat J$ at $\bar u$ is given by $\nabla\hat J(\bar u)=\sigma(\bar u-\un)-\mathcal B'\bar p\in\U$ with $\bar p=\hat p+\mathcal A\bar u$. Then the following four statements are equivalent:
    \begin{enumerate}
        \item [\rm 1)] $\bar u\in\Uad$ and
        \begin{align*}
            {\langle\nabla\hat J(\bar u),u-\bar u\rangle}_\U\ge0\quad\text{for all }u\in\Uad
        \end{align*}
        \item [\rm 2)] $\bar u\in\Uad$ and in case of $\U=\mathbb R^\mathsf m$ 
        \begin{align*}
            (\nabla\hat J(\bar u))_i=\left\{
            \begin{aligned}
                &=0&&\text{if }u_{\mathsf ai}<\bar u_i<u_{\mathsf bi},\\
                &\ge0&&\text{if }u_{\mathsf ai}=\bar u_i<u_{\mathsf bi},\\
                &\le0&&\text{if }u_{\mathsf ai}<\bar u_i=u_{\mathsf bi},
            \end{aligned}
            \right.\quad1\le i\le\mathsf m,
        \end{align*}
        in case of $\U=L^2(\mathscr D;\mathbb R^\mathsf m)$
        \begin{align*}
            (\nabla\hat J(\bar u))_i(s)=\left\{
            \begin{aligned}
                &=0&&\text{if }u_{\mathsf ai}(s)<\bar u_i(s)<u_{\mathsf bi}(s),\\
                &\ge0&&\text{if }u_{\mathsf ai}(s)=\bar u_i(s)<u_{\mathsf bi}(s),\\
                &\le0&&\text{if }u_{\mathsf ai}(s)<\bar u_i(s)=u_{\mathsf bi}(s),
            \end{aligned}
            \right.\quad1\le i\le\mathsf m\text{ a.e. in }\mathscr D,
        \end{align*}
        respectively.
        \item [\rm 3)] The two functions
        \begin{align*}
            \bar\mu_\mathsf a=\max\big\{\nabla\hat J(\bar u),0\big\}\text{ in }\U\quad\text{and}\quad\bar\mu_\mathsf b=-\min\big\{\nabla\hat J(\bar u),0\big\}\text{ in }\U
        \end{align*}
        satisfy
        \begin{align*}
            \bar u\in\Uad,\quad\bar\mu_a,\bar\mu_b\ge0\text{ in }\U,\quad{\langle\bar\mu_\mathsf a,\bar u-u_\mathsf a\rangle}_\U={\langle\bar\mu_\mathsf b,u_\mathsf b-\bar u\rangle}_\U=0.
        \end{align*}
        We call $\bar\mu_\mathsf a$ and $\bar\mu_\mathsf b$ {\em Lagrange multipliers} associated with the inequality constraints $\bar u-u_\mathsf a\le0$ in $\U$ and $\bar u-u_\mathsf b\le0$ in $\U$, respectively.
        \item [\rm 4)] For any $\eta>0$ we have
        \begin{align*}
            \bar u=\mathcal P_\Uad(\bar u-\eta\nabla\hat J(\bar u)).
        \end{align*}
    \end{enumerate}
\end{proposition} 

\subsection{A Lagrangian-based approach for \eqref{GVLuminy:Eq4.1.6}}
\label{SIAM-Book:Section4.1.4}

Until now, we have dealt with the optimal control problem \eqref{GVLuminy:Eq4.1.6} by using the solution operator $\mathcal S$. This approach allowed us to introduce the reduced cost function $\hat J$ which led to the equivalent problem \eqref{GVLuminy:Eq4.1.9} in which the state equation $y=\hat y+\mathcal S u$ did not explicitly occur. In this section, we will present an equivalent, Lagrangian-based approach by writing the state equation as a specific equality constraint $e(y,u)=0$. To do so, let us define the product Hilbert space
\begin{align*}
    \mP=L^2(0,T;V)\times H
\end{align*}
endowed with the canonical inner product. We identify $\mP'$ with $L^2(0,T;V')\times H$ and introduce the affine-linear operator $e=(e_1,e_2):\X\to\mP'$ 
\begin{equation}
    \label{ST-1}
    \begin{aligned}
        {\langle e_1(x),p^1\rangle}_{L^2(0,T;V'),L^2(0,T;V)}&=\int_0^T{\langle y_t(t),p^1(t)\rangle}_{V',V}+a(t;y(t),p^1(t))\,\mathrm dt\\
        &\quad-\int_0^T{\langle(\mathcal F+\mathcal Bu)(t),p^1(t)\rangle}_{V',V}\,\mathrm dt,\\
        {\langle e_2(x), p^2 \rangle}_H &= \langle y(0)-y_\circ, p^2 \rangle_H
    \end{aligned}
\end{equation}
for $x=(y,u)\in\X$ and $p = (p^1,p^2) \in \mP$. Then \eqref{GVLuminy:Eq4.1.1} can be expressed equivalently by the equality constraint
\begin{align*}
    e(x)=0\quad\text{in }\mP'.
\end{align*}
In particular, we infer that
\begin{equation}
    \label{HI-4}
    e(\hat y+\mathcal Su,u)=0\quad\text{for all }u\in\Uad.
\end{equation}
Now \eqref{GVLuminy:Eq4.1.6} can be written as a quadratic programming problem with equality and inequality constraints:
\begin{equation}
    \label{QP}
    \tag{\bf{QP}} 
    \min J(x)\quad\text{s.t.}\quad x=(y,u)\in\X,~e(x)=0\text{ in }\mP'\text{ and }u\in\Uad.
\end{equation}
Let $\bar u\in\Uad$ be the unique solution to \eqref{GVLuminy:Eq4.1.9}. Then $\bar x=(\hat y+\mathcal S\bar u,\bar u)\in \Xad$ is the unique solution to \eqref{QP}. The Lagrange function for \eqref{QP} is defined as
\begin{align*}
    \mathcal L(x,p)=J(x)+{\langle e(x),p\rangle}_{\mP',\mP}=J(x)+{\langle e_1(x),p^1\rangle}_{L^2(0,T;V'),L^2(0,T;V)}+{\langle e_2(x),p^2\rangle}_H
\end{align*}
for $x=(y,u)\in\X$ and $p=(p^1,p^2)\in\mP$. It follows from \eqref{HI-4} that
\begin{equation}
    \label{HI-5}
    \hat J(u)=J(\hat y+\mathcal Su,u)=\mathcal L(\hat y+\mathcal Su,u,p)
\end{equation}
for every $u\in\Uad$ and $p\in\mP$.

To characterize the solution to \eqref{QP} we derive first-order optimality conditions. Therefore, we will show that the Lagrangian is Fr\'echet-differentiable; cf. Definition~\ref{Definition:GFDer}. For a proof we refer to Section~\ref{SIAM-Book:Section4.7.1}.

\begin{lemma}
    \label{Lemma:Jprime}
    Let Assumption~{\rm\ref{A9}} be valid. The mapping $J$ is continuously Fr\'echet-differentiable from $\X$ to $\mathbb R$ and
    \begin{equation}
        \label{Custer-1}
        J'(x)x^\delta=\sigma_1\int_0^T{\langle y(t)-\ydQ(t),y^\delta(t)\rangle}_H\,\mathrm dt+\sigma_2\,{\langle y(T)-\ydT,y^\delta(T)\rangle}_H+\sigma\,{\langle u-\un,u^\delta\rangle}_\U
    \end{equation}
    for $x=(y,u)\in\X$ and $x^\delta=(y^\delta,u^\delta)\in\X$.
\end{lemma}

\begin{remark}
    \label{Remark:HI-106}
    \rm Note that $J'(x)$ is an element in $\X'$ for any $x\in\X$. Thus, we can also write $J'(x)x^\delta=\langle J'(x),x^\delta\rangle_{\X',\X}$ for every $x^\delta\in\X$.\hfill$\blacksquare$
\end{remark}

In Section~\ref{SIAM-Book:Section4.7.1} the following result is proved.

\begin{lemma}
    \label{Lemma:EGprime}
    Let Assumption~{\rm\ref{A8}} hold. Then the operator $e:\X\to\Lambda'$ is continuously Fr\'echet-differentiable in $\X$. In particular, we have
    \begin{equation}
        \label{ST-2}
        \begin{aligned}
            {\langle e_1'(x)x^\delta,p^1\rangle}_{L^2(0,T;V'),L^2(0,T;V)}&=\int_0^T{\langle y^\delta_t(t),p^1(t)\rangle}_{V',V}+a(t;y^\delta(t),p^1(t))\,\mathrm dt\\
            &\quad-\int_0^T{\langle(\mathcal Bu^\delta)(t),p^1(t)\rangle}_{V',V}\,\mathrm dt,\\
            e_2'(x)x^\delta&=y^\delta(0)
        \end{aligned}
    \end{equation}
    for $x=(y,u)\in\X$, $x^\delta=(y^\delta,u^\delta)\in\X$ and $p^1\in L^2(0,T;V)$.
\end{lemma}

From Lemmas~\ref{Lemma:Jprime} and \ref{Lemma:EGprime} and the Fr\'echet-differentiability of the mapping $(\mP' \times \mP) \ni (f,p) \mapsto {\langle f,p \rangle}_{\mP',\mP}$ it follows that the Lagrangian is continuously Fr\'echet-differentiable in $\mathscr X\times\Lambda$. In particular, we infer from \eqref{HI-5} that
\begin{equation}
    \label{HI-8}
    \begin{aligned}
        \hat J'(u)u^\delta&=\mathcal L_y(\hat y+\mathcal Su,u,p)\mathcal Su^\delta+\mathcal L_u(\hat y+\mathcal Su,u,p)u^\delta\\
        &={\langle\mathcal L_y(\hat y+\mathcal Su,u,p),\mathcal Su^\delta\rangle}_{\Y',\Y}+{\langle\mathcal L_u(\hat y+\mathcal Su,u,p),u^\delta\rangle}_{\U',\U}
    \end{aligned}
\end{equation}
for every $u\in\U$, $p\in\mP$ and $u^\delta\in\U$, where $\mathcal L_y$ and $\mathcal L_u$ denote the partial Fr\'echet derivatives with respect to $y$ and $u$, respectively.

Let $x\in\X$ be chosen arbitrarily. It follows from Lemma~\ref{Lemma:EGprime} that
\begin{align*}
    e'(x)=\big(e_y(x)\,\big|\,e_u(x)\big):\X\to\mP'
\end{align*}
is a bounded linear operator, where the partial Fr\'echet-derivatives $e_y=(e_{1y},e_{2y}):X\to\mP'$ and $e_u=(e_{1u},e_{2u}):\Y\to\mP'$ are given as
\begin{align*}
    {\langle e_{1y}'(x)y^\delta,p^1\rangle}_{L^2(0,T;V'),L^2(0,T;V)}&=\int_0^T{\langle y^\delta_t(t),p^1(t)\rangle}_{V',V}+a(t;y^\delta(t),p^1(t))\,\mathrm dt,\\
    e_{2y}'(x)y^\delta&=y^\delta(0),\\
    {\langle e_{1u}'(x)u^\delta,p^1\rangle}_{L^2(0,T;V'),L^2(0,T;V)}&=-\int_0^T{\langle(\mathcal Bu^\delta)(t),p^1(t)\rangle}_{V',V}\,\mathrm dt,\\
    e_{2u}'(x)u^\delta&=0
\end{align*}
for $y^\delta \in \Y$, $u^\delta\in\U$ and $p^1\in L^2(0,T;V)$.

In the following lemma, which is proved in Section~\ref{SIAM-Book:Section4.7.1}, it is proved that the linear operator $e_y(x)$ possesses a bounded inverse for all $x \in \X$. This implies that the mapping $e$ satisfies a constrained qualification condition; cf. \cite{Lue69,MZ79,NW06}.

\begin{lemma}
    \label{Lemma:EprimeReg}
    Let Assumption~{\rm\ref{A8}} hold and $x=(y,u)\in\X$. Then the linear and bounded operator $e_y(x):\Y\to\mP'$ is bijective with a bounded inverse $e_y(x)^{-1}$. In particular, $e'(x)$ is surjective.
\end{lemma}

\begin{remark}
    \label{Remark:HI-35}
    \rm Differentiating \eqref{HI-4} with respect to $u\in\U$ in direction $u^\delta\in\U$ we find
    \begin{align*}
        e_y(\hat y+\mathcal Su,u)\mathcal Su^\delta+e_u(\hat y+\mathcal Su,u)u^\delta=0\quad\text{in }\mP'.
    \end{align*}
    Since the inverse operator $e_y(x)^{-1}:\mP'\to\Y$ exists and is bounded for every $x\in\X$, we infer that $y^\delta=\mathcal Su^\delta$ satisfies
    \begin{align*}
        y^\delta=-e_y(\hat y+\mathcal Su,u)^{-1}e_u(\hat y+\mathcal Su,u)u^\delta\quad\text{in }\Y.
    \end{align*}
    Moreover, if $\mathscr J_\mP: \mP \to \mP'$ and $\mathscr J_\Y: \Y \to \Y'$ denote the Riesz isomorphisms, then the operator $\mathscr J_\mP^{-1} e_y(x): \Y \to \mP$ is an isomorphism of Hilbert spaces, so it has a Hilbert adjoint $(\mathscr J_\mP^{-1} e_y(x))^\star: \mP \to \Y$ which is likewise bijective with inverse $(\mathscr J_\mP^{-1} e_y(x))^{\star,-1}: \Y \to \mP$.\\
    Lastly, since $\mP$ is a Hilbert space, it is reflexive and we can identify the \index{Space!bidual}{\em bidual space} $(\mP')'$ with $\mP$ itself. Therefore, the dual operator $e_y(x)'$ can be considered a mapping from $\mP$ to $\Y'$. Furthermore, it is bijective and its the inverse is given by
    \begin{align*}
        \big(e_y(x)'\big)^{-1} = (\mathscr J_\mP^{-1} e_y(x))^{\star,-1} \mathscr J_\Y^{-1}: ~\Y' \to \mP
    \end{align*}
    which can be shown in a straightforward way.\hfill $\blacksquare$
\end{remark}

The next proposition yields a characterization of the first-order optimality condition $\mathcal L_y(x,p)=0$ in $\Y'$. For a proof we refer to Section~\ref{SIAM-Book:Section4.7.1}.

\begin{proposition}
    \label{Prop:HI-30}
    Let Assumption~{\rm\ref{A8}} and {\rm\ref{A9}} hold, $u\in\U$, $y(u)=\hat y+\mathcal Su\in\Y$ and $x(u)=(y(u),u)\in\X$. Then 
    \begin{equation}
        \label{Eq:HI-35}
        \mathscr L_y(x(u),p(u))=0\quad\text{in }\Y'
    \end{equation}
    is equivalent with
    \begin{equation}
        \label{HI-7a}
        p(u)=\big(p^1(u),p^2(u)\big)=-\big(e_y(y(u),u)'\big)^{-1}J_y(y(u),u)\in\mP.
    \end{equation}
    Moreover, $p=(p^1,p^2)$ defined by \eqref{HI-7a} satisfies $p^1=\hat p+\mathcal Au$ and $p^2=p^1(0)$.
\end{proposition}

For now, let us quickly summarize what we have proven altogether for a variable $\bar x = (\hat y + \mathcal S \bar u, \bar u)$.
\begin{itemize}
	\item Theorem \ref{GVLuminy:Theorem4.3.1} stated that $\bar x$ solves the quadratic programming problem \eqref{QP} if and only if $\bar u \in \Uad$ and
	\begin{align*}
	{\langle \nabla \hat J(\bar u),u -\bar u\rangle}_\U \ge 0 \quad\text{for all }u \in \Uad.
	\end{align*}
	\item Associated with the primal variable $\bar x$ the adjoint variable $\bar p$ is defined like in \eqref{HI-7a} to $\bar p=-(e_y(\bar x)')^{-1} J_y(\bar x) \in \mP$. It was shown in Proposition \ref{Prop:HI-30} that this variable is given by $\bar p = (\bar p^1,\bar p^2)$ with $\bar p^1 = \hat p + \mathcal A u$ (the solution to the adjoint equation \eqref{GVLuminy:Eq4.3.16}) and $\bar p^2=\bar p^1(0)$.
	\item By Proposition \ref{Prop:HI-30}, we know that the variable $\bar p$ also satisfies $\mathcal L_y(\bar x, \bar p) = 0$ in $\Y'$. Inserting this into the derivative representation \eqref{HI-8} yields $\hat J'(\bar u)=\mathcal L_u(\bar x, \bar p)$ in $\U'$. 
\end{itemize}
We summarize the above observations in the following theorem. 

\begin{theorem}
    Let Assumption~{\rm\ref{A8}} and {\rm\ref{A9}} hold. Then $\bar x=(\hat y+\mathcal S\bar u,\bar u)$ is a solution to \eqref{QP} if and only if the first-order optimality system holds for $\bar p=(\bar p^1,\bar p^2)$ with $\bar p^1=\hat p+\mathcal A\bar u$ and $\bar p^2=\bar p^1(0)$:
    \begin{align*}
        \mathcal L_y(\bar x,\bar p)=0\text{ in }\Y',\quad\mathcal L_u(\bar x,\bar p)(u-\bar u)\ge0\text{ for all }u\in\Uad\quad\text{and}\quad \bar u\in\Uad.
    \end{align*}
\end{theorem}

\subsection{Second-order derivatives}
\label{SIAM-Book:Section4.1.5}

In this section we study second derivatives (see Definition~\ref{def:secondOrderDerivatives}) of the cost functional and of the constraints. 

\begin{lemma}
    \label{Lemma:HI-106}
    Let Assumption~{\rm\ref{A8}} be satisfied. Then $J:\X\to\mathbb R$ is twice continuously Fr\'echet-differentiable in $\X$ and
    \begin{equation}
        \label{HI-106}
        \begin{aligned}
            {\langle J''(x)x^\delta,\tilde x^\delta\rangle}_{\X',\X}&=\sigma_1\int_0^T{\langle \tilde y^\delta(t),y^\delta(t)\rangle}_H\,\mathrm dt+\sigma_2\,{\langle\tilde y^\delta(T),y^\delta(T)\rangle}_H+\sigma\,{\langle\tilde u^\delta,u^\delta\rangle}_\U
        \end{aligned}
    \end{equation}
    for $x=(y,u)\in\X$, $x^\delta=(y^\delta,u^\delta)\in\X$ and $\tilde x^\delta=(\tilde y^\delta,\tilde u^\delta)\in\X$. Furthermore, the operator $e:\X\to\mP'$ is twice continuously Fr\'echet-differentiable and its second derivative is zero.
\end{lemma}

For a proof of Lemma~\ref{Lemma:HI-106} we refer to Section~\ref{SIAM-Book:Section4.7.1}. It follows from the lemma that the Lagrange function $\mathcal L$ is twice continuously Fr\'echet-differentiable. In particular, we find that
\begin{equation}
    \label{HI-200}
    {\langle\mathcal L_{xx}(x,p)x^\delta,x^\delta\rangle}_{\X',\X}=\sigma_1\int_0^T{\|y^\delta(t)\|}^2_H\,\mathrm dt+\sigma_2\,{\|y^\delta(T)\|}_H^2+\sigma\,{\|u^\delta\|}_\U^2
\end{equation}
for $x=(y,u)\in\X$, $p\in\mP$ and $x^\delta=(y^\delta,u^\delta)\in\X$. Recall that for every $x\in\X$, the kernel of $e'(x)$ is defined as
\begin{align*}
    \mathrm{ker}\,e'(x)=\big\{x^\delta\in\X\,\big|\,e'(x)x^\delta=0\text{ in }\mP'\big\};
\end{align*}
cf. Definition~\ref{SIAM:Definition-I.1.1.1}. Elements in $\mathrm{ker}\,e'(x)$ can be characterized as follows.

\begin{lemma}
    \label{Lemma:HI-200}
    Let Assumption~{\rm\ref{A8}} hold. Suppose that $x=(y,u)\in\X$ and $x^\delta=(y^\delta,u^\delta)\in\mathrm{ker}\,e'(x)\subset\X$. Then 
    \begin{equation}
        \label{KernelEstimate}
        {\|y^\delta\|}_\Y\le c_\mathsf{ker}\,{\|u^\delta\|}_\U
    \end{equation}
    for a constant $c_\mathsf{ker}>0$.
\end{lemma}

For a proof of Lemma~\ref{Lemma:HI-200} we refer to Section~\ref{SIAM-Book:Section4.7.1}. Estimate \eqref{KernelEstimate} can be utilized to proof the \index{Optimality conditions!second-order!sufficient}{\em second-order sufficient optimality condition} for \eqref{QP}. This is stated in the next proposition which is proved in Section~\ref{SIAM-Book:Section4.7.1}.

\begin{proposition}
    \label{Proposition:HI-110}
    Let Assumptions~{\rm\ref{A8}} and {\rm\ref{A9}} hold. Suppose that $\sigma>0$ is satisfied. Then the second Fr\'echet derivative $\mathcal L_{xx}(x,p)$ with respect to $x$ is coercive on $\mathrm{ker}\,e'(x)$ for all $x\in\X$ and $p\in\mP$, i.e., there exists a constant $\kappa>0$ such that
    \begin{align*}
        {\langle\mathcal L_{xx}(x,p)x^\delta,x^\delta\rangle}_{\X',\X}\ge\kappa\,{\|x^\delta\|}_\X^2\quad \text{for all } x^\delta\in\mathrm{ker}\,e'(x)
    \end{align*}
    is valid.
\end{proposition}

\begin{remark}
    \label{Remark:HI-200}
    \rm Let Assumptions~{\rm\ref{A8}} and {\rm\ref{A9}} hold. Suppose that $\bar x=(\bar y,\bar u)$ and $\bar p=(\bar p^1,\bar p^2)$ solve the first order system
    \begin{align*}
        \mathcal L_y(\bar x,\bar p)=0\text{ in }\Y',\quad\mathcal L_u(\bar x,\bar p)(u-\bar u)\ge0\text{ for all }u\in\Uad\quad\text{and}\quad \bar u\in\Uad.
    \end{align*}
    Then we infer from Proposition~\ref{Proposition:HI-110} that $\bar x$ satisfies a second-order sufficient optimality condition for \eqref{QP}. In particular, $\bar x$ is the unique global solution to \eqref{QP}.\hfill$\blacksquare$
\end{remark}

\section{The continuous POD method for the optimal control problem}
\label{SIAM-Book:Section4.2}
\setcounter{equation}{0}
\setcounter{theorem}{0}

In this section, we introduce POD Galerkin approximations to the state and adjoint variables in the optimal control problem \eqref{GVLuminy:Eq4.1.6a}. The linear state and adjoint solution operators $\mathcal S$, $\mathcal A$ from the previous section are then replaced by lower-dimensional counterparts $\mathcal S^\ell$, $\mathcal A^\ell$. We study the convergence of these POD discretizations, where we make use of the analysis in \cite{HV08,KV01,KV02a,KV02b,Sin14,TV09}. For a general introduction we also refer the reader to the survey papers \cite{BSV14,HV05,SV10}.

\subsection{POD Galerkin schemes for the dual equation}
\label{SIAM-Book:Section4.2.1}

Most numerical solution strategies for \eqref{GVLuminy:Eq4.1.6a} rely on repeatedly solving systems like the state equation $y = \hat y + \mathcal S u$ or the adjoint equation $p = \hat p + \mathcal A u$. If the underlying spaces $H$, $V$ are highly-dimensional, this requires a lot of computational effort. To work around this problem, we can compute a POD basis $\{\psi_i\}_{i=1}^\ell$ of rank $\ell$ and then perform a Galerkin projection onto the according POD space. As a result, the high-dimensional problem \eqref{GVLuminy:Eq4.1.6a} is replaced by a lower-dimensional one. In this section, we we will perform this projection and discuss consequences for the primal and dual variables. 

Recall that in \eqref{SIAM:Eq3.2.12a} and \eqref{SIAM:Eq3.2.12b} we have defined the projections
\begin{align*}
    &\mathcal P^\ell_H:H\to H^\ell,&&v^\ell=\mathcal P^\ell_H\varphi&&\text{solves }\min_{w^\ell\in H^\ell}{\|\varphi-w^\ell\|}_H\quad \text{ for }\varphi\in H,\\
    &\mathcal Q^\ell_H:H\to V^\ell,&&v^\ell=\mathcal Q^\ell_H\varphi&&\text{solves }\min_{w^\ell\in V^\ell}{\|\varphi-w^\ell\|}_H\quad \text{ for }\varphi\in H.
\end{align*}
We require some data trajectories $\{y^k(t)\,|\,0 \le t \le T\}$ ($1 \le k \le K$) to compute a POD basis. This is done as follows:
\begin{enumerate}
    \item [1)] Choose an admissible control $u\in\Uad$ and compute the solutions $y=\hat y+\mathcal Su$ and $p=\hat p+\mathcal Au$ to \eqref{GVLuminy:Eq4.1.1} and \eqref{GVLuminy:Eq4.3.16}, respectively. We assume that $y,p\in H^1(0,T;V)$ holds which implies that $y,p\in H^1(0,T;X)$.
    \item [2)] Set $K=2$, $\omega_i^K=1$ for $i=1,2$, $y^1=y$ and $y^2=p$, i.e., the \index{POD method!snapshot space}{\em snapshot space} is given as
    \begin{equation}
        \label{HI-203}
        \mathscr V=\mathrm{span}\,\bigg\{\int_0^Tw_1(t)y(t)+w_2(t)p(t)\,\mathrm dt\,\Big|\,w_1,w_2\in L^2(0,T)\bigg\}\subset X
    \end{equation}
    with the dimension
    \begin{align*}
        d=\left\{
        \begin{aligned}
            &\dim\mathscr V&&\text{if }\dim\mathscr V<\infty,\\
            &\infty&&\text{otherwise}.
        \end{aligned}
        \right.
    \end{align*}
    \item [3)] For any $\ell\in\mathbb N$ with $\ell\le d$ solve (cf. \pageref{SIAM:PellH})
    \begin{equation}
        \tag{$\mathbf P^\ell_H$}
        \label{PellH}
        \left\{
        \begin{aligned}
            &\min\sum_{k=1}^K \omega_k^K\int_0^T \Big\| y^k(t)-\sum_{i=1}^\ell {\langle y^k(t),\psi_i\rangle}_H\,\psi_i\Big\|_H^2\,\mathrm dt\\
            &\hspace{0.5mm}\text{s.t. } \{\psi_i\}_{i=1}^\ell\subset H\text{ and }{\langle\psi_i,\psi_j\rangle}_H=\delta_{ij},~1 \le i,j \le \ell
        \end{aligned}
        \right.
    \end{equation}
    or
    \begin{equation}
        \tag{$\mathbf P^\ell_V$}
        \label{PellV}
        \left\{
        \begin{aligned}
            &\min\sum_{k=1}^K \omega_k^K\int_0^T \Big\| y^k(t)-\sum_{i=1}^\ell {\langle y^k(t),\psi_i\rangle}_V\,\psi_i\Big\|_V^2\,\mathrm dt\\
            &\hspace{0.5mm}\text{s.t. } \{\psi_i\}_{i=1}^\ell\subset V\text{ and }{\langle\psi_i,\psi_j\rangle}_V=\delta_{ij},~1 \le i,j \le \ell.
        \end{aligned}
        \right.
    \end{equation}
    to obtain POD bases $\{\psi_i^H\}_{i=1}^\ell$ or $\{\psi^V_i\}_{i=1}^\ell$, respectively, with rank $\ell$. According to the a-priori error estimates (cf. \eqref{Chopin-2}) the number $\ell$ should be chosen as follows: for given tolerance $\varepsilon_\mathsf{POD}>0$
    \begin{align*}
        \xi\sum_{i=\ell+1}^d\lambda_i^H\,\big\|\psi_i^H\big\|_V^2&<\varepsilon_\mathsf{POD}\text{ if }X=H\quad\text{and}\quad\xi\sum_{i=\ell+1}^d\lambda_i^V\,{\|\psi_i^V-\mathcal Q^\ell_H\psi_i^V\|}_V^2<\varepsilon_\mathsf{POD}\text{ if }X=V
    \end{align*}
    with a factor $\xi\in (0,1)$.
    \item [4)] Set $H^\ell=\mathrm{span}\,\big\{\psi_1^H,\ldots,\psi_\ell^H\big\}$ or $V^\ell=\mathrm{span}\,\big\{\psi_1^V,\ldots,\psi_\ell^V\big\}$.
\end{enumerate}

In other words, we choose the data trajectories for the POD basis as the primal and dual variables $y$ and $p$ belonging to a particular control $u$. As long as an optimization algorithm remains in close proximity to this control, the POD surrogate model will usually be a good approximation of the original problem. 

\begin{remark}
    \label{Remark:HI-201}
    \rm
    \begin{enumerate}
        \item [1)] Recall that the pairs $\{(\lambda_i^H,\psi_i^H)\}_{i=1}^\ell$ and $\{(\lambda_i^V,\psi_i^V)\}_{i=1}^\ell$ satisfy the eigenvalue problems \eqref{EigProblemHV-H} and \eqref{EigProblemHV-V}, respectively.
        \item [2)] From Lemma~\ref{SIAM:Lemma3.2.1}-1) it follows that $H^\ell\subset V$ holds true.
        \item [3)] If we do not distinguish between the choice $X=H$ or $X=V$ we only write $\{\psi_i\}_{i=1}^\ell$ for the computed POD basis of rank $\ell$ and define the subspace $X^\ell=\mathrm{span}\,\big\{\psi_1,\ldots,\psi_\ell\big\}\subset V$.
        \item [4)] Note that we do not include time derivatives in our snapshot set. Of course, we can also enlarge the snapshot space by setting $\wp=4$, $\omega_i^\wp=1$ for $i=1,\ldots,4$, $y^1=y$, $y^2=y_t$, $y^3=p$ and $y^4=p_t$, i.e., the snapshot space is given as
        \begin{align*}
            \mathscr V=\bigg\{\sum_{k=1}^4\int_0^T\phi_k(t)y^k(t)\,\mathrm dt\in X \Big| \phi_k \in L^2(0,T) \text{ for } k=1,\dots,4 \bigg\}.
        \end{align*}
        In that case the obtained rates of a-priori convergence results are based on Theorem~\ref{Th:A-PrioriError}.\hfill$\blacksquare$
    \end{enumerate}
\end{remark}

Recall that in \eqref{SIAM:Eq3.2.12a} and \eqref{SIAM:Eq3.2.12b} we have defined the projections
\begin{align*}
    &\mathcal P^\ell_H:H\to H^\ell,&&v^\ell=\mathcal P^\ell_H\varphi&&\text{solves }\min_{w^\ell\in H^\ell}{\|\varphi-w^\ell\|}_H\quad \text{ for }\varphi\in H,\\
    &\mathcal Q^\ell_H:H\to V^\ell,&&v^\ell=\mathcal Q^\ell_H\varphi&&\text{solves }\min_{w^\ell\in V^\ell}{\|\varphi-w^\ell\|}_H\quad \text{ for }\varphi\in H.
\end{align*}
We will utilize the results presented in Proposition~\ref{Prop:VTopology} and Section~\ref{SIAM-Book:Section3.3.2}. For the readers convenience we sum up the required hypotheses in the next

\begin{assumption}
    \label{A10}
    \begin{enumerate}
        \item [\em 1)] Let Assumptions~{\em\ref{A8}} and {\em\ref{A9}} hold. Moreover, $y=\hat y+\mathcal Su$ and $p=\hat p+\mathcal Au$ belong to the space $H^1(0,T;V)$ for any $u\in\Uad$.
        \item [\em 2)] The Hilbert space $X$ denotes either $H$ or $V$. In case of $X=V$ we suppose that $y_\circ\in V$ is valid.
        \item [\em 3)] We set $K=2$, $\omega_1^K=\omega_2^K=1$ and $y^1=y$, $y^2=p$. The snapshot space $\mathscr V$ is given by \eqref{HI-203} with dimension $d\in\mathbb N\cup\{\infty\}$. If $X=H$ is chosen, we suppose that $d\in\mathbb N$ is satisfied, i.e., $\lambda_i^H=0$ for $i>d$. For $\ell\in\mathbb N$ with $\ell\le d$ let $\{\psi_i^H\}_{i=1}^\ell$ and $\{\psi_i^V\}_{i=1}^\ell$ be solutions to \eqref{PellH} and \eqref{PellV}, respectively.
        \item [\em 4)] In case of $X=H$, let $\mathcal P^\ell$ denote $\mathcal P^\ell_H$. If $X=V$, let it be given by $\mathcal Q^\ell_H$.
    \end{enumerate}
\end{assumption}

The POD Galerkin scheme for the state equation \eqref{GVLuminy:Eq4.1.1} has already been introduced in Section \ref{SIAM-Book:Section3.3.1}. If $y = \hat y + \mathcal S u$ denotes the full-order state variable, its POD-Galerkin approximation is defined as $y^\ell = \hat y^\ell + \mathcal S^\ell u$ and solves
\begin{subequations}
    \label{eq:reducedStateEq}
    \begin{align}
        \label{eq:reducedStateEqa}
        \frac{\mathrm d}{\mathrm dt} {\langle y^\ell(t),\psi \rangle}_H+a(t;y^\ell(t),\psi)&={\langle (\mathcal F+\mathcal Bu)(t),\psi\rangle}_{V',V}\quad\text{for all }\psi\in X^\ell\text{ a.e. in }(0,T],\\
        \label{eq:reducedStateEqb}
        y^\ell(0)&=\mathcal P^\ell y_\circ,
    \end{align}
\end{subequations}
where $\hat y^\ell$ and $\mathcal S^\ell u$ solve \eqref{eq:reducedStateEq_yHat} and \eqref{eq:reducedStateEq_Sell}, respectively. The operator $\mathcal S^\ell: \U\to H^1(0,T;V)$ was shown in Corollary \ref{Corollary:HI-21} to be linear and continuous.  If Assumption~\ref{A10} hold, we can infer from $\|\mathcal P^\ell y_\circ\|_H\le\|y_\circ\|_H$ for all $y_\circ\in H$ and from Theorem~\ref{SIAM:Theorem3.1.1POD} that
\begin{equation}
    \label{Ida-1}
    {\|y^\ell\|}_\Y\le C_y\left({\|y_\circ\|}_H+{\|\mathcal F\|}_{L^2(0,T;V')}+{\|u\|}_\U\right)
\end{equation}
for a constant $C_y>0$ which is independent of $\ell$. Combining Theorems~\ref{Th:A-PrioriError} and ~\ref{Th:A-PrioriError-2}, there also exists a constant $C>0$ satisfying the \index{Error estimate!a-priori!state variable}a-priori error estimate
\begin{align*}
    {\|y-y^\ell\|}_\Y^2\le C\cdot\left\{
    \begin{aligned}
        &\sum_{i>\ell}\lambda_i^H\,\big\|\psi_i^H\big\|_V^2,&&\text{for }X=H\\
        &\sum_{i>\ell}\lambda_i^V\,\big\|\psi_i^V-\mathcal P^\ell\psi_i^V\big\|_V^2,&&\text{for }X=V.
    \end{aligned}
    \right.
\end{align*}
In this subsection, we will perform a analogous Galerkin projection for the adjoint equation \eqref{GVLuminy:Eq4.3.16} which was used to express the gradient of the function $\hat{J}$. First of all, we make the ansatz 
\begin{equation}
    \label{DualGalAn}
    p^\ell(t)=\sum_{i=1}^\ell\mathrm p_i^\ell(t)\psi_i\in X^\ell,\quad t\in[0,T],
\end{equation}
and replace $y,p$ by $y^\ell, p^\ell$ in \eqref{GVLuminy:Eq4.3.16}. Second, we consider the equation \eqref{GVLuminy:Eq4.3.16a} in $(X^\ell)'$ instead of $V'$ and project the terminal data in \eqref{GVLuminy:Eq4.3.16b} onto $X^\ell$ which leads to the following problem:
\begin{subequations}
    \label{SIAM:DualPOD}
    \begin{align}
        \label{SIAM:DualPODa}
        -\frac{\mathrm d}{\mathrm dt} {\langle p^\ell(t),\psi \rangle}_H+a(t;\psi,p^\ell(t))&=\sigma_1\,{\langle \ydQ(t)-y^\ell(t),\psi\rangle}_H\quad\text{for all }\psi\in X^\ell\text{ a.e. in }[0,T),\\
        \label{SIAM:DualPODb}
        p^\ell(T)&=\sigma_2\big(\mathcal P^\ell\ydT-y^\ell(T)\big).
    \end{align}
\end{subequations}
Note that $\tilde p^\ell(t)=p^\ell(T-t)$, $t\in[0,T]$, satisfies
\begin{subequations}
    \label{SIAM:DualPODt}
    \begin{align}
        \label{SIAM:DualPODta}
        \frac{\mathrm d}{\mathrm dt} {\langle \tilde p^\ell(t),\psi \rangle}_H+\tilde a(t;\psi,\tilde p^\ell(t))&=\sigma_1\,{\langle\ydQ(t)-y^\ell(T-t),\psi\rangle}_H\quad\text{for all }\psi\in X^\ell\text{ a.e. in }(0,T],\\
        \label{SIAM:DualPODtb}
        \tilde p^\ell(T)&=\sigma_2\big(\mathcal P^\ell\ydT-y^\ell(T)\big).
    \end{align}
\end{subequations}
We can observe that \eqref{SIAM:DualPODt} has a very similar structure to \eqref{eq:reducedStateEq}, so existence and uniqueness of solutions along with a boundedness result in the manner of \eqref{Ida-1} can be shown exactly as the proof of Theorem~\ref{SIAM:Theorem3.1.1POD}, which gives the following result.

\begin{proposition}
    \label{Prop:DualPODApriori}
    Let Assumptions~{\rm\ref{A10}} hold. Then for every $u\in\U$ there exists a unique $p^\ell\in H^1(0,T;X^\ell)\hookrightarrow\Y$ to \eqref{SIAM:DualPOD} satisfying
    \begin{align*}
    {\|p^\ell\|}_\Y\le C\big(\sigma_2\,{\|\ydT-y^\ell(T)\|}_H+\sigma_1\,{\|\ydQ-y^\ell\|}_{L^2(0,T;H)}\big)
    \end{align*}
    for a constant $C\ge0$ which is independent of $\ell$ and for $y^\ell=\hat y^\ell+\mathcal S^\ell u$.
\end{proposition}

The following corollary follows directly from Proposition~\ref{Prop:DualPODApriori}. For more details we refer to Section~\ref{SIAM-Book:Section4.7.2}.

\begin{corollary}
    \label{Corollary:HI-40}
    Let all assumptions of Proposition~{\rm\ref{Prop:DualPODApriori}} hold.
    \begin{enumerate}
        \item [\rm 1)] There exists a unique solution $\hat p^\ell \in H^1(0,T;X^\ell)\hookrightarrow\Y$ to
        \begin{align*}
            -\frac{\mathrm d}{\mathrm dt} {\langle\hat p^\ell(t),\psi \rangle}_H+a(t;\psi,\hat p^\ell(t))&=\sigma_1\,{\langle\ydQ(t)-\hat y^\ell(t),\psi\rangle}_H\quad\text{for all }\psi\in X^\ell\text{ a.e. in }[0,T),\\
            \hat p^\ell(T)&=\sigma_2\big(\mathcal P^\ell\ydT-\hat y^\ell(T)\big).
        \end{align*}
        In particular, we have
        \begin{align*}
            {\|\hat p^\ell\|}_\Y\le C\big(\sigma_2\,{\|\ydT-\hat y^\ell(T)\|}_H+\sigma_1\,{\|\ydQ-\hat y^\ell\|}_{L^2(0,T;H)}\big)
        \end{align*}
        for a constant $C\ge0$ which is independent of $\ell$.
        \item [\rm 2)] Let us introduce the linear operator $\mathcal A^\ell:\U\to \Y$ as follows: $p^\ell=\mathcal A^{\ell}u$, $u\in\U$, satisfies
        \begin{align*}
            -\frac{\mathrm d}{\mathrm dt} {\langle p^\ell(t),\psi \rangle}_H+a(t;\psi,p^\ell(t))&=-\sigma_1\,{\langle (\mathcal S^\ell u)(t),\psi\rangle}_H\quad\text{for all }\psi\in X^\ell\text{ a.e. in }[0,T),\\
            p^\ell(T)&=-\sigma_2(\mathcal S^\ell u)(T).
        \end{align*}
        Then $\mathcal A^{\ell}$ is well-defined and bounded. In particular,
        \begin{align}
            \label{eq:AellContinuity1}
            {\|\mathcal A^\ell u\|}_\Y\le C\big(\sigma_2\,{\|(\mathcal S^\ell u)(T)\|}_H+\sigma_1\,{\|\mathcal S^\ell u\|}_{L^2(0,T;H)}\big)
        \end{align}
        for a constant $C\ge0$ which is independent of $\ell$.
    \end{enumerate}
\end{corollary}

\begin{remark}
    \label{Remark:KoH-1}
    \rm
    \begin{enumerate}
        \item [1)] The solution to \eqref{SIAM:DualPOD} can be expressed as $p^\ell=\hat p^\ell+\mathcal A^\ell u$.
        \item [2)] The equation \eqref{SIAM:DualPODa} can be expressed as an $\ell$-dimensional system of ordinary differential equations. For that purpose we introduce the $\ell$-dimensional vectors
        \begin{align*}
            \mathrm p^\ell(t)=\big(\mathrm p_i^\ell(t)\big)_{1\le i\le\ell},\quad\mathrm y_1^\ell(t)=\big({\langle\ydQ(t),\psi_i\rangle}_H\big)_{1\le i\le\ell},\quad\mathrm y_2^\ell=\big(\mathrm y_{2 i}^\ell\big)_{1\le i\le\ell},
        \end{align*}
        where
        \begin{align*}
            \mathcal P^\ell \ydT=\sum_{i=1}^\ell\mathrm y_{2 i}^\ell\psi_i\in X^\ell.
        \end{align*}
        Recall that the matrices $\bM^\ell\in\mathbb R^{\ell\times\ell}$ and $\bA^\ell(t)\in\mathbb R^{\ell\times\ell}$ have been defined in \eqref{M_AMatrices}. Then inserting \eqref{DualGalAn} as well as \eqref{PODGalAns} into \eqref{SIAM:DualPODa} and choosing $\psi=\psi_i$, $i=1,\ldots,\ell$, we derive
        \begin{subequations}
            \label{SIAM:DualPODODE}
            \begin{equation}
                \label{SIAM:DualPODODEa}
                -\bM^\ell\dot{\mathrm p}^\ell(t)+\bA^\ell(t)^\top\mathrm p^\ell(t)=\sigma_1\big(\mathrm y_1^\ell(t)-\bM^\ell\mathrm y^\ell(t)\big),\quad t\in[0,T).
            \end{equation}
            Moreover, \eqref{SIAM:DualPODb} implies that
            \begin{equation}
            \label{SIAM:DualPODODEb}
            \mathrm p^\ell(T)=\sigma_2\big(\mathrm y_2^\ell-\mathrm y^\ell(T)\big)
        \end{equation}
        holds true.\hfill$\blacksquare$
    \end{subequations}
    \end{enumerate}
\end{remark}

In the next corollary we present an a-priori bound for $p^\ell$. The proof is given in Section~\ref{SIAM-Book:Section4.7.2}.

\begin{corollary}
    \label{Cor:pell}
    Let all assumptions of Proposition~{\rm\ref{Prop:DualPODApriori}} hold. Then for $p^{\ell} = \hat{p} + \mathcal{A}^{\ell}(u)$ it holds
    \begin{align*}
        {\|p^\ell\|}_\Y\le C_p\,\big({\|\ydT\|}_H+{\|\ydQ\|}_{L^2(0,T;H)}+{\|y_\circ\|}_H+{\|f\|}_{L^2(0,T;V')}+{\|u\|}_\U\big)
    \end{align*}
    for a constant $C_p\ge0$ which is independent of $\ell$.
\end{corollary}

\subsection{POD a-priori error analysis for the dual equation}
\label{SIAM-Book:Section4.2.2}

By Assumption \ref{A10}-4), $\mathcal P^\ell$ is either given by $\mathcal P^\ell_H$ or $\mathcal Q_H$. Therefore, it is always an $H$-orthonormal projection. Our next goal is to derive a-priori error estimates for the term
\begin{align*}
    \int_0^T {\|p(t)-p^\ell(t)\|}_V^2\,\mathrm dt,
\end{align*}
where $p$ and $p^\ell$ are the solutions to \eqref{GVLuminy:Eq4.3.16} and \eqref{SIAM:DualPOD}, respectively. We suppose that $p\in H^1(0,T;V)$ holds. Recall that $X^\ell=\mathrm{span}\,\{\psi_1,\ldots,\psi_\ell\}\subset V$ for the two choices $X=H$ and $X=V$. Analogously to Section~\ref{SIAM-Book:Section3.3.2} we will make use of the following decomposition (cf. \eqref{EqDec})
\begin{align*}
    p(t)-p^\ell(t)=p(t)-\mathcal P^\ell p(t)+\mathcal P^\ell p(t)-p^\ell(t)=\rho^\ell(t)+\theta^\ell(t)
\end{align*}
with $\rho^\ell(t)=p(t)-\mathcal P^\ell p(t)\in (X^\ell)^\bot$ and $\theta^\ell(t)=\mathcal P^\ell p(t)-p^\ell(t)\in X^\ell$. Note that $(X^\ell)^\perp$ denotes the orthogonal space to $X^\ell$ with respect to the $H$-inner product because $\mathcal P^\ell$ is an $H$-orthogonal projection. The above decomposition gives
\begin{equation}
    \label{Dual:APriori-Est-1}
    \int_0^T {\|p(t)-p^\ell(t)\|}_V^2\,\mathrm dt\le2\int_0^T {\|\rho^\ell(t)\|}_V^2\,\mathrm dt+2\int_0^T {\|\theta^\ell(t)\|}_V^2\,\mathrm dt.
\end{equation}

\begin{lemma}
    \label{Lemma:HI-50}
    Let Assumption~{\rm\ref{A10}} hold. Let $\theta^\ell=\mathcal P^\ell p-p^\ell$, where $p$ and $p^\ell$ are the solutions to \eqref{GVLuminy:Eq4.3.16} and \eqref{SIAM:DualPOD}, respectively. Then 
    \begin{equation}
        \label{Albi-4a}
        {\|\theta^\ell\|}_{L^2(0,T;V)}^2\le c_\theta\Big({\|\mathcal P^\ell p-p\|}_{L^2(0,T;V)}^2+{\|\mathcal P^\ell y-y\|}_\Y^2\Big)
    \end{equation}
    for a constant $c_\theta>0$ which is independent of $\ell$.
\end{lemma}

\begin{theorem}
    \label{Th:DualA-PrioriError-2}
    Assume that Assumption~{\rm\ref{A10}} hold. Let $\|\mathcal P^\ell_{HV^\ell}\|_{\mathscr L(V)}$ be bounded independently of $\ell$. Suppose that $p$ and $p^\ell$ are the solutions to \eqref{GVLuminy:Eq4.3.16} and \eqref{SIAM:DualPOD}, respectively. Then there exists a constant $C>0$ satisfying the \index{Error estimate!a-priori!dual variable}a-priori error estimate
    \begin{align*}
        \int_0^T\big\|p(t)-p^\ell(t)\big\|_V^2\,\mathrm dt\le C\cdot\left\{
        \begin{aligned}
            &\sum_{i>\ell}\lambda_i^H\,\big\|\psi_i^H\big\|_V^2,&&\text{for }X=H,\\
            &\sum_{i>\ell}\lambda_i^V\,\big\|\psi_i^V-\mathcal Q^\ell_H\psi_i^V\big\|_V^2,&&\text{for }X=V.
        \end{aligned}
        \right.
    \end{align*}
\end{theorem}

\begin{remark}
    \label{GVLuminy:Remark4.4.2}
    \rm Let us discuss variants of the a-priori error estimate for the dual variable provided different snapshots are chosen.
    \begin{enumerate}
        \item [1)] In Assumption~\ref{A10}-3) we could also choose $K=1$, $\omega_1^K=1$ and
        \begin{align*}
            \mathscr V=\mathrm{span}\,\bigg\{\int_0^T\phi(t)y(t)\,\mathrm dt\,\big|\,\phi\in L^2(0,T)\bigg\}\subset X,
        \end{align*}
        i.e., the adjoint variable $p$ is not included in the snapshot space. Then we only get the a-priori error estimates
        \begin{align*}
            \int_0^T\big\|p(t)-p^\ell(t)\big\|_V^2\,\mathrm dt\le C\bigg({\|\mathcal P^\ell p-p\|}_{L^2(0,T;V)}^2+\sum_{i>\ell}\lambda_i^H\,\big\|\psi_i^H\big\|_V^2\bigg)
        \end{align*}
        for $X=H$ and
        \begin{align*}
            \int_0^T\big\|p(t)-p^\ell(t)\big\|_V^2\,\mathrm dt\le C\bigg({\|\mathcal P^\ell p-p\|}_{L^2(0,T;V)}^2+\sum_{i>\ell}\lambda_i^V\,\big\|\psi_i^V-\mathcal P^\ell\psi_i^V\big\|_V^2\bigg)
        \end{align*}
        for $X=V$. Thus, the inclusion of adjoint information into the snapshot ensemble improves the approximation quality. This is also observed for nonlinear problems; see \cite{DV01}.
        \item [2)] Now we choose $X=V$, $K=4$, $y^1=y$, $y^2=y_t$, $y^3=p$, $y^4=p_t$, $\omega_k^K=1$ for $k=1,\ldots,4$ and
        \begin{align*}
            \mathscr V=\mathrm{span}\,\bigg\{\sum_{i=1}^4\int_0^T\phi_k(t)y^k(t)\,\mathrm dt\,\Big|\,\phi_k\in L^2(0,T)\text{ for }k=1,\ldots,4\bigg\}\subset X,
        \end{align*}
        i.e. the time derivatives are included in the snapshot space. Moreover, we take the projection $\mathcal P^\ell=\mathcal P^\ell_V$ with
        \begin{align*}
            \mathcal P^\ell_V:V\to V^\ell,\quad v^\ell=\mathcal P^\ell_V\varphi\text{ solves }\min_{w^\ell\in V^\ell}{\|\varphi-w^\ell\|}_V\quad \text{ for }\varphi\in V.
        \end{align*}
        Then the following a-priori error estimate (cf. Theorem~\ref{Th:A-PrioriError})
        \begin{align*}
            \int_0^T\big\|y(t)-y^\ell(t)\big\|_V^2+\big\|p(t)-p^\ell(t)\big\|_V^2\,\mathrm dt\le C\sum\limits_{i>\ell}\lambda_i^V
        \end{align*}
        holds for a constant $C>0$. Consequently, the rate of convergence for the state and adjoint variable is given by the decay of the eigenvalues.\hfill$\blacksquare$
    \end{enumerate}
\end{remark}

For every $u\in\Uad$ we set $p=\hat p+\mathcal Au$ and $p^\ell=\hat p+\mathcal A^\ell u$. We infer from Theorem~\ref{Th:DualA-PrioriError-2}
\begin{equation}
    \label{ConvDualVar}
    {\|p-p^\ell\|}_{L^2(0,T;V)}={\|(\mathcal A-\mathcal A^\ell)u\|}_{L^2(0,T;V)}\to0\quad\text{for }\ell\to\infty.
\end{equation}

Recall that we have introduced the linear, bounded operators $\mathcal A$ and $\mathcal A^\ell$ in \eqref{GVLuminy:Eq4.3.12} and Corollary~\ref{Corollary:HI-40}-2, respectively. We give sufficient conditions that $\mathcal A^\ell$ converges to $\mathcal A$ in $\mathscr L(\U,\Y)$ as $\ell\to\infty$. The proof of the result can be shown as in the proof of Theorem~\ref{Theorem:OperatorConv}.

\begin{theorem}
    \label{Theorem:DualOperatorConv}
    Let Assumptions~{\rm\ref{A8}}, {\rm\ref{A9}} and {\rm\ref{A10}} hold. Assume that $\|\mathcal Q^\ell_H\|_{\mathscr L(V)}$ is bounded independently of $\ell$. We suppose that either $\U$ is a finite-dimensional Hilbert space or that the operator $\mathcal A$ is even compact. Moreover, let one of the following assumptions be satisfied:
    \begin{enumerate}
        \item [\rm 1)] $X=H$, $\mathcal P^\ell=\mathcal P^\ell_H$, $\lambda_i^H>0$ for all $i\in\mathbb I$ and $y(t)$, $p(t)\in\mathrm{ran}\,(\mathcal Y)\subset V$ f.a.a. $t\in[0,T]$, where the linear and bounded operator $\mathcal Y:L^2(0,T;\mathbb R^2)\to X$ is given as
        \begin{align*}
            \mathcal Y\phi=\int_0^T\phi^1(t)y(t,\cdot)+\phi^2(t)p(t,\cdot)\,\mathrm dt\quad \text{for }\phi=(\phi^1,\phi^2)\in L^2(0,T;\mathbb R^2);
        \end{align*}
        \item [\rm 2)] $X=H$, $\mathcal P^\ell=\mathcal P^\ell_V$ and $\lambda_i^H>0$ holds for all $i\in\mathbb I$;
        \item [\rm 3)] $X=V$, $\mathcal P^\ell=\mathcal Q^\ell_V$ and $\lambda_i^V>0$ holds for all $i\in\mathbb I$.
    \end{enumerate}
    Then $\|\mathcal A^\ell-\mathcal A\|_{\mathscr L(\U,\Y)}\to0$ as $\ell\to\infty$.
\end{theorem}

\subsection{POD Galerkin scheme for the optimality system}
\label{SIAM-Book:Section4.2.3}

In this section, we will present the POD Galerkin approximation for the optimal control problem \eqref{GVLuminy:Eq4.1.9}. It is defined as follows:
\begin{equation}
    \tag{$\mathbf {\hat P}^\ell$}
    \label{GVLuminy:Eq4.4.12}
    \min \hat J^\ell(u) \quad \text{s.t.} \quad u \in\Uad,
\end{equation}
where the cost is given by $\hat J^\ell(u)=J(\hat y^\ell+\mathcal S^\ell u,u)$ for $u \in\U$. Due to its similar structure, many results for \eqref{GVLuminy:Eq4.1.9} from Section \ref{SIAM-Book:Section4.1} translate to \eqref{GVLuminy:Eq4.4.12} and the proofs are almost identical. Therefore, we will leave it to the interested reader to proof the results below as an exercise; cf. Theorem~\ref{GVLuminy:Theorem4.2.2}, Lemma~\ref{Lemma:HI-30} and Remark~\ref{rem:concreteGradient}.
\begin{theorem}
    Let Assumption~{\em\ref{A10}} be satisfied.
    \begin{enumerate}
        \item [\em 1)] Problem \eqref{GVLuminy:Eq4.4.12} possesses a solution $\bar u$. If $\sigma_1 > 0$ or $\sigma > 0$ the solution is unique.
        \item [\em 2)] Let $\Theta: W_0(0,T) \to W_0(0,T))'$ and $\Xi: L^2(0,T;H) \times H \to W_0(0,T)'$ by given by \eqref{GVLuminy:Eq4.3.10}. Then it holds
        \begin{align*}
            \mathcal B' \mathcal A^\ell=-(\mathcal S^\ell)' \Theta \mathcal S^\ell\text{ in }L(\U),\quad\mathcal B'\hat p=(\mathcal S^\ell)' \Xi\big(y_1^\mathsf d-\hat y^\ell,y_2^\mathsf d-\hat y^\ell(T)\big)\text{ in }\U.
        \end{align*}
        \item [\em 3)] The gradient of $\hat J^\ell$ is given by
        \begin{align*}
            \nabla\hat J^\ell(u)=-\mathcal B' p^\ell+\sigma(u-\tilde u^\mathsf n)\quad\text{with }p^\ell=\mathcal A^\ell u+\hat p^\ell.
        \end{align*}
    \end{enumerate}
\end{theorem}

We will assume that \eqref{GVLuminy:Eq4.4.12} admits a unique solution $\bar u^\ell$. Similar to the non-reduced case, a first-order sufficient optimality condition is given by the variational inequality
\begin{equation}
    \label{GVLuminy:Eq4.4.13}
    {\langle \sigma(\bar u^\ell-\un)-\mathcal B' \bar p^\ell,u-\bar u^\ell\rangle}_\U \ge 0 \quad \text{for all }u \in\Uad
\end{equation}
with $\bar p^\ell=\hat p^\ell+\mathcal A^\ell\bar u^\ell$. The next theorem we present an a-priori result for the control variable. Its proof is given in Section~\ref{SIAM-Book:Section4.7.2}.

\begin{theorem}
    \label{GVLuminy:Theorem4.4.2}
    Assume that Assumptions~{\rm\ref{A10}} and that $\bar u$ is the solution to \eqref{GVLuminy:Eq4.1.9}. Let $\|\mathcal Q^\ell_H\|_{\mathscr L(V)}$ be bounded independently of $\ell$. Furthermore, $\bar u^\ell$ denotes the unique solution to \eqref{GVLuminy:Eq4.4.13}. Then there exists a constant $C>0$ depending on $c_V$, $\gamma$, $\gamma_1$, $T$, $\sigma_1$, $\sigma_2$ and the norm $\|\mathcal B'\|_{\mathscr L(L^2(0,T;V),\U)}$ such that the \index{Error estimate!a-priori!control variable}a-priori error estimate
    \begin{equation}
        \label{GVLuminy:Eq4.4.15}
        \big\|\bar u-\bar u^\ell\big\|_\U^2 \le C\cdot\left\{
        \begin{aligned}
            &\sum_{i>\ell}\lambda_i^H\,\big\|\psi_i^H\big\|_V^2&&\text{for }X=H,\\
            &\sum_{i>\ell}\lambda_i^V\,\big\|\psi_i^V-\mathcal Q_H^\ell\psi_i^V\big\|_V^2&&\text{for }X=V
        \end{aligned}
        \right.
    \end{equation}
    holds.
\end{theorem}

\begin{remark}
    \label{Remark:HI-202}
    \rm From Remark~\ref{GVLuminy:Remark4.4.2} we obtain the following consequences:
    \begin{enumerate}
        \item [1)] Instead of Assumption~\ref{A10}-2) let us choose $K=1$, $\omega_1^K=1$ and
        \begin{align*}
            \mathscr V=\mathrm{span}\,\bigg\{\int_0^T\phi(t)y(t)\,\mathrm dt\,\Big|\,\phi\in L^2(0,T)\bigg\}\subset X
        \end{align*}
        i.e., the adjoint variable $p$ is not included in the snapshot space. Then we only get the a-priori error estimates
        \begin{align*}
            {\|\bar u-\bar u^\ell\|}_\U^2\le C\bigg({\|\mathcal P^\ell p-p\|}_{L^2(0,T;V)}^2+\sum_{i>\ell}\lambda_i^H\,\big\|\psi_i^H\big\|_V^2\bigg)
        \end{align*}
        for $X=H$ and
        \begin{align*}
            {\|\bar u-\bar u^\ell\|}_\U^2\le C\bigg({\|\mathcal P^\ell p-p\|}_{L^2(0,T;V)}^2+\sum_{i>\ell}\lambda_i^V\,\big\|\psi_i^V-\mathcal Q_H^\ell\psi_i^V\big\|_V^2\bigg)
        \end{align*}
        for $X=V$. Thus, the inclusion of adjoint information into the snapshot ensemble improves the approximation quality for the control.
        \item [2)] Now we choose $X=V$, $K=4$, $y^1=y$, $y^2=y_t$, $y^3=p$, $y^4=p_t$, $\omega_k^K=1$ for $k=1,\ldots,4$ and
        \begin{align*}
            \mathscr V=\mathrm{span}\,\bigg\{\sum_{k=1}^4\int_0^T\phi_k(t)y^k(t)\,\mathrm dt\,\Big|\,\phi_k\in L^2(0,T)\text{ for }k=1,\ldots,4\bigg\}\subset X,
        \end{align*}
        i.e. the time derivatives are included in the snapshot space. Then (cf. Theorem~\ref{Th:A-PrioriError})
        \begin{align*}
            {\|\bar u-\bar u^\ell\|}_\U^2\le C\sum\limits_{i>\ell}\lambda_i^V
        \end{align*}
        for a constant $C>0$ and $\mathcal P^\ell=\mathcal Q^\ell_V$. Thus, the rate of convergence for the control variable is given by the decay of the eigenvalues.\hfill$\blacksquare$
    \end{enumerate}
\end{remark}

\subsection{POD a-posteriori error analysis for the optimality system}
\label{SIAM-Book:Section4.2.4}

In this section, we fix the more general control space $\U = L^2(0,T;\mathbb R^m)$. All results transfer in a natural way to the simpler case of $\U = \mathbb R^m$. \\
In \cite{TV09}, POD a-posteriori error estimates are presented which can be applied to our optimal control problem as well. Based on a perturbation method (cf. \cite{DHPY95}) it is deduced how far the suboptimal control $\bar u^\ell$ is from the (unknown) exact optimal control $\bar u$. Thus, our goal is to estimate the norm $\|\bar u-\bar u^\ell\|_\U$ without the knowledge of the optimal solution $\bar u$. In general, $\bar u^\ell\neq\bar u$ holds, so that $\bar u^\ell$ does not satisfy the variational inequality \eqref{GVLuminy:Eq4.3.20-2}. However, we can define a function $\zeta^\ell\in\U$ such that
\begin{equation}
    \label{GVLuminy:Eq4.5.1}
    {\langle\sigma(\bar u^\ell-\un)-\mathcal B'\tilde p^\ell+\zeta^\ell,u-\bar u^\ell\rangle}_\U\ge0\quad\text{for all }u\in\Uad,
\end{equation}
with $\tilde p^\ell=\hat p+\mathcal A\bar u^\ell$. Therefore, $\bar u^\ell$ satisfies the optimality condition of the perturbed parabolic optimal control problem
\begin{align*}
    \min_{u\in \Uad}\tilde J(u)=J(\hat y+\mathcal Su,u)+{\langle \zeta^\ell,u\rangle}_\U
\end{align*}
with ``perturbation'' $\zeta^\ell$. The smaller $\zeta^\ell$ is, the closer $\bar u^\ell$ is to $\bar u$. Next we estimate $\|\bar u-\bar u^\ell\|_\U$ in terms of $\|\zeta^\ell\|_\U$. By Lemma~\ref{Lemma:HI-30}-2) we have
\begin{equation}
    \label{GVLuminy:Eq4.5.3}
    \mathcal B' \big(\bar p-\tilde p^\ell\big)=\mathcal B' \mathcal A\big(\bar u-\bar u^\ell\big)=-\mathcal S'\Theta \mathcal S\big(\bar u-\bar u^\ell\big)=\mathcal S'\Theta\big(\tilde y^\ell-\bar y\big)
\end{equation}
with $\tilde y^\ell=\hat y+\mathcal S\bar u^\ell$. Choosing $u=\bar u^\ell$ in \eqref{GVLuminy:Eq4.3.20-2}, $u=\bar u$ in \eqref{GVLuminy:Eq4.5.1} and using \eqref{GVLuminy:Eq4.5.3} we obtain
\begin{align*}
    0&\le{\langle-\sigma(\bar u-\bar u^\ell)+\mathcal B'(\bar p-\tilde p^\ell)+\zeta^\ell,\bar u-\bar u^\ell\rangle}_\U\\
    &=-\sigma\,{\|\bar u-\bar u^\ell\|}^2_\U+{\langle\mathcal S'\Theta(\tilde y^\ell-\bar y),\bar u-\bar u^\ell\rangle}_\U+{\langle\zeta^\ell,\bar u-\bar u^\ell\rangle}_\U\\
    &=-\sigma\,{\|\bar u-\bar u^\ell\|}^2_\U-{\langle\Theta(\bar y-\tilde y^\ell),\bar y-\tilde y^\ell\rangle}_{W_0(0,T)',W_0(0,T)}+{\langle\zeta^\ell,\bar u-\bar u^\ell\rangle}_\U\\
    &\le-\sigma\,{\|\bar u-\bar u^\ell\|}^2_\U+{\langle\zeta^\ell,\bar u^\ell-\bar u^\ell\rangle}_\U\le-\sigma\,{\|\bar u-\bar u^\ell\|}^2_\U+{\|\zeta^\ell\|}_\U{\|\bar u-\bar u^\ell\|}_\U.
\end{align*}
Hence, we get the a-posteriori error estimation
\begin{align*}
    {\|\bar u-\bar u^\ell\|}_\U \le \frac{1}{\sigma} \, {\|\zeta^\ell\|}_\U.
\end{align*}
This estimate is stated in the following theorem, where, in addition, an appropriate perturbation function $\zeta^\ell$ is chosen. For the proof we refer to Section~\ref{SIAM-Book:Section4.7.2}.

\begin{theorem}
    \label{GVLuminy:Theorem4.5.1}
    Suppose that Assumption~{\rm\ref{A10}} holds. Let $\hat y,\hat p\in H^1(0,T;V)$ and $u\in\U = L^2(0,T;\mathbb R^m)$ be arbitrarily given so that $\mathcal Su,\,\mathcal Au \in H^1(0,T;V)\setminus\{0\}$. To compute a POD basis $\{\psi_i\}_{i=1}^\ell$ of rank $\ell$ we choose $\wp=2$, $y^1=\hat y+\mathcal Su$ and $y^2=\hat p+\mathcal Au$. Define the function $\zeta^\ell\in\U$ by
    \begin{align*}
        \zeta_i^\ell(t)=\left\{
        \begin{aligned}
            &-\min(0,\xi_i^\ell(s))&&\text{a.e. in }\mathscr A^\ell_{\mathsf ai}=\big\{s\in\mathscr D\,|\bar u^\ell_i(s)=u_{\mathsf ai}(s)\big\},\\
            &-\max(0,\xi^\ell_i(s))&&\text{a.e. in }\mathscr A^\ell_{\mathsf bi}=\big\{s\in\mathscr D\,|\bar u^\ell_i(s)=u_{\mathsf bi}(s)\big\},\\
            &-\xi_i^\ell(t) && \text{a.e. in }\mathscr D\setminus\big(\mathscr A^\ell_{\mathsf ai}\cup\mathscr A_{\mathsf bi}^\ell\big),
        \end{aligned}
        \right.
    \end{align*}
    where $\xi^\ell=\sigma(\bar u^\ell-\un)-\mathcal B' (\hat p+\mathcal A\bar u^\ell)$ in $\U$. Then the \index{Error estimate!a-posteriori!control variable}a-posteriori error estimate
    \begin{equation}
        \label{GVLuminy:Eq4.5.7a}
        {\|\bar u-\bar u^\ell\|}_\U\le\frac{1}{\sigma}\,{\|\zeta^\ell\|}_\U
    \end{equation}
    holds. In particular, $\lim \limits_{\ell \to \infty}\big\|\zeta^\ell\big\|_\U=0$.
\end{theorem}

\begin{remark}
    \label{GVLuminy:Remark4.5.1}
    \rm
    \begin{enumerate}
        \item [1)] Theorem~\ref{GVLuminy:Theorem4.5.1} shows that $\|\zeta^\ell\|_\U$ can be expected smaller than any tolerance $\epsilon>0$ provided that $\ell$ is taken sufficiently large. Motivated by this result we set up Algorithm~4.2.1. Note that the quality of the POD Galerkin scheme is improved by only increasing the number of POD basis elements. Another approach is to update the POD basis in the optimization process; see, e.g., \cite{AH02,AFS00,KV08}. 
        \item [2)] In \cite{SV13} POD a-posteriori error estimates are tested numerically for a linear-quadratic optimal control problem. It turns out that in certain cases only an increase of the number of POD ansatz functions does not decrease the error in the reduced-order solution satisfyingly. In this case a change of the POD basis is needed; see, \cite{KV08,Vol11}, for instance.
        \item [3)] Let us refer to \cite{KV12}, where POD a-posteriori error estimates are combined with a sequential quadratic programming method in order to solve a nonlinear PDE constrained optimal control problem. Furthermore, the presented analysis for linear-quadratic problems can be extended to semilinear optimal control problems by a second-order analysis; see in \cite{KTV13}.\hfill$\blacksquare$
    \end{enumerate}
\end{remark}

\medskip
\hrule
\vspace{-3.5mm}
\begin{algorithm}[(POD method with a-posteriori estimator)]
    \label{GVLuminy:Algorithm4.5.1}
    \vspace{-3mm}
    \hrule
    \vspace{0.5mm}
    \begin{algorithmic}[1]
        \REQUIRE Initial control $u^{\ell 0}\in\U$, initial number $\ell$ for the POD ansatz functions, a maximal number $\ell_\mathsf{max}>\ell$ of POD ansatz functions, maximal number $k_\mathsf{max}$ of iterations and a stopping tolerance $\epsilon>0$.
        \STATE Set $k=0$ and determine $\hat y$, $\hat p$, $y^1=\mathcal Su^{\ell 0}$, $y^2=\mathcal Au^{\ell 0}$.
        \STATE Compute a POD basis $\{\psi_i\}_{i=1}^{\ell_\mathsf{max}}$ choosing $y^1$ and $y^2$.
        \REPEAT
            \STATE Establish the POD Galerkin discretization using $\{\psi_i\}_{i=1}^\ell$.
            \STATE Set $k=k+1$ and call Algorithm~4.2 to compute suboptimal control $\bar u^{\ell k}$.
            \STATE Determine $\zeta^\ell$ according to Theorem~\ref{GVLuminy:Theorem4.5.1} and compute $\epsilon_\mathsf{ape}=\|\zeta^\ell\|_\U/\sigma$.
            \IF {$\epsilon_\mathsf{ape}<\epsilon$ {\bf or} $\ell=\ell_\mathsf{max}$}
                \STATE Return $\ell$ and suboptimal control $\bar u^{\ell k}$ and STOP.
            \ENDIF
            \STATE Set $\ell=\ell+1$.
        \UNTIL{$\ell>\ell_\mathsf{max}$ {\bf or} $k=k_\mathsf{max}$}
    \end{algorithmic}
    \hrule
\end{algorithm}

\section{The semidiscrete approximation}
\label{SIAM-Book:Section4.3}
\setcounter{equation}{0}
\setcounter{theorem}{0}

While the last section dealt with POD approximations of the infinite-dimensional optimal control problem \eqref{GVLuminy:Eq4.1.6}, in this section we deal with the more practically relevant situation where the state equation undergoes a semidiscretization, replacing the spaces $V$ and $H$ in the Gelfand triple by a finite-dimensional space $V^h$. We then study a POD approximation of the semidiscretized system, which is in practice used to further reduce the computational complexity of the problem. 

\subsection{Galerkin discretization}
\label{SIAM-Book:Section4.3.1}

In Section~\ref{Section:3.4.1} we have introduced the $m$-dimensional subspace
\begin{align*}
    V^h=\mathrm{span}\,\big\{\varphi_1^h,\ldots,\varphi_m^h\big\}\subset V.
\end{align*}
Suppose that $\hat y:[0,T]\to V^h$ satisfies
\begin{subequations}
    \label{MD-1.1}
    \begin{align}
        \label{MD-1a}
        \frac{\mathrm d}{\mathrm dt} {\langle\hat y^h(t),\varphi^h\rangle}_H+a(t;\hat y^h(t),\varphi^h)&={\langle \mathcal F(t),\varphi^h\rangle}_{V',V}\quad\text{for all }\varphi^h\in V^h\text{ a.e. in }(0,T],\\
        \label{MD-1b}
        \hat y^h(0)&=\mathcal P^h y_\circ,
    \end{align}
\end{subequations}
where $\mathcal P^h:X\to V^h$ is an appropriate linear projection operator and $X$ stands either for the space $H$ or for the space $V$. In Example~\ref{ExampleFEProjection} we have presented two possible choices for $\mathcal P^h$. It follows from Theorem~\ref{Theorem:FESystem} that there exists a unique solution $\hat y^h$ to \eqref{MD-1.1} satisfying $\hat y^h\in H^1(0,T;V)$ and
\begin{align*}
    {\|\hat y^h\|}_\Y\le C\Big({\|\mathcal P^h y_\circ\|}_H+{\|\mathcal F\|}_{L^2(0,T;V')}\Big)
\end{align*}
for a constant $C>0$ which is independent of $y_\circ$ and $f$. Moreover, we define the linear solution operator $\mathcal S^h:\U\to \Y$ as follows: for given $u\in\U$ let $y^h=\mathcal S^hu$ be the solution to 
\begin{subequations}
    \label{MD-2}
    \begin{align}
        \label{MD-2a}
        \frac{\mathrm d}{\mathrm dt} {\langle y^h(t),\varphi^h \rangle}_H+a(t;y^h(t),\varphi^h)&={\langle (\mathcal Bu)(t),\varphi^h\rangle}_{V',V}\quad\text{for all }\varphi^h\in V^h\text{ a.e. in }(0,T],\\
        \label{MD-2b}
        y^h(0)&=0.
    \end{align}
\end{subequations}
Due to Theorem~\ref{Theorem:FESystem} the operator $\mathcal S^h$ is well-defined and bounded. Moreover, $y^h=\hat y^h+\mathcal S^hu$, $u\in\U$, is the unique solution to 
\begin{subequations}
    \label{MD-3}
    \begin{align}
        \label{MD-3a}
        \frac{\mathrm d}{\mathrm dt} {\langle y^h(t),\varphi^h \rangle}_H+a(t;y^h(t),\varphi^h)&={\langle (\mathcal F+\mathcal Bu)(t),\varphi^h\rangle}_{V',V}\quad\text{for all }\varphi^h\in V^h\text{ a.e. in }(0,T],\\
        \label{MD-3b}
        y^h(0)&=\mathcal P^h y_\circ,
    \end{align}
\end{subequations}
which has already been introduced in \eqref{EvProGal}. Now, we define the following Galerkin approximation for the reduced cost functional:
\begin{align*}
    \hat J^h(u)=J(\hat y^h+\mathcal S^hu,u) \quad \text{for }u\in\U.
\end{align*}
The Galerkin approximation for \eqref{GVLuminy:Eq4.1.9} reads
\begin{equation}
    \label{PhatFE}
    \tag{$\mathbf{\hat P}^h$}
    \min \hat J^h(u)\quad\text{s.t.}\quad u\in\Uad.
\end{equation}
Problem \eqref{PhatFE} is a linear-quadratic optimal control problem. Since $\mathcal S^h$ is a linear and bounded operator, existence and uniqueness of a solution $\bar u^h\in\Uad$ to \eqref{PhatFE} follows from Theorem~\ref{GVLuminy:Theorem4.2.1}. Moreover, a first-order sufficient optimality condition is given by
\begin{equation}
    \label{MD-1}
    {\langle\sigma(\bar u^h-\un)-\mathcal B'\bar p^h,u-\bar u^h\rangle}_\U\ge0\quad\text{for all }u\in\Uad,
\end{equation}
where $\bar p^h$ solves the adjoint equation
\begin{subequations}
    \label{MD-2-2}
    \begin{align}
        \label{MD-2-2a}
        -\frac{\mathrm d}{\mathrm dt}\,{\langle\bar p^h(t),\varphi^h\rangle}_H+a(t;\varphi^h,\bar p^h(t))&=\sigma_1\,{\langle\ydQ(t)-\bar y^h(t),\varphi^h\rangle}_H\quad\text{for all }\varphi^h\in V^h\text{ a.e. in }[0,T),\\
        \label{MD-2-2b}
        \bar p^h(T)&=\sigma_2\big(\mathcal P^h\ydT-\bar y^h(T)\big)
    \end{align}
\end{subequations}
and $\bar y^h=\hat y^h+\mathcal S^h\bar u$ holds. In \eqref{MD-1} we have utilized that the gradient $\nabla\hat J^h$ of the reduced cost functional $\hat J^h$ is given by the formula
\begin{equation}
    \label{GradientJhat-h}
    \nabla\hat J^h(u)=\sigma(u-\un)-\mathcal B'(\hat p^h+\mathcal A^hu)\quad\text{for }u\in\Uad.
\end{equation}
Note that \eqref{MD-2-2} is the Galerkin scheme for the adjoint equation \eqref{GVLuminy:Eq4.3.16}. For $\bar p^h\in H^1(0,T;V^h)$ we have the representation
\begin{equation}
    \label{Equation:ch4_GalerkinAdjoint}
    \bar p^h(t)=\sum_{j=1}^m\bar{\mathrm p}_j^h(t)\varphi_j^h\quad\text{f.a.a. }t\in[0,T]
\end{equation}
with a coefficient function $\bar{\mathrm p}^h=(\bar{\mathrm p}_1^h,\ldots,\bar{\mathrm p}_m^h):[0,T]\to\mathbb R^m$. Moreover, we introduce the two $m$-dimensional vectors
\begin{align*}
    \mathrm y_1^h(t)=\big({\langle\ydQ(t),\varphi^h_i\rangle}_H\big)_{1\le i\le m},~t\in[0,T],\quad\text{and}\quad\mathrm y_2^h=\big(y^h_{2 i}\big)_{1\le i\le m},
\end{align*}
where
\begin{align*}
    \mathcal P^h \ydT=\sum_{i=1}^h\mathrm y_{2 i}^h\varphi^h_i\in V^h.
\end{align*}
Recall that the matrices $\bM^h\in\mathbb R^{m\times m}$ and $\bA^h(t)\in\mathbb R^{m\times m}$ have been defined in Example~\ref{ExampleFEProjection} and \eqref{FEStiffnessMatrix}, respectively. Then  inserting \eqref{Equation:ch4_GalerkinAdjoint} into \eqref{MD-2-2a} and choosing $\varphi^h=\varphi^h_i$, $i=1,\ldots,m$, we derive the following end-value problem
\begin{equation}
    \label{MD-4}
    \begin{aligned}
        -\bM^h\dot{\bar{\mathrm p}}^h(t)+\bA^{h,\top}(t)\bar{\mathrm p}^h(t)&=\sigma_1\big(\mathrm y_1^h(t)-\bM^h\bar{\mathrm y}^h(t)\big)&&\text{for }t\in(0,T],\\
        \bar{\mathrm p}^h(T)&=\sigma_2\big(\mathrm y_2^h-\bM^h\bar{\mathrm y}^h(T)\big),
    \end{aligned}
\end{equation}
where $\bar{\mathrm y}^h:[0,T]\to\mathbb R^m$ is the solution to \eqref{FineModel} for the choice $u=\bar u^h$.

Let $\hat p^h\in H^1(0,T;V^h)$ be the unique solution to 
\begin{subequations}
    \label{MD-5}
    \begin{align}
        \label{MD-5a}
        -\frac{\mathrm d}{\mathrm dt}\,{\langle\hat p^h(t),\varphi^h\rangle}_H+a(t;\varphi^h,\hat p^h(t))&=\sigma_1\,{\langle\ydQ(t)-\hat y^h(t),\varphi^h\rangle}_H\quad\text{for all }\varphi\in V^h\text{ a.e. in }[0,T),\\
        \label{MD-5b}
        \hat p^h(T)&=\sigma_2\big(\mathcal P^h\ydT-\hat y^h(T)\big).
    \end{align}
\end{subequations}
Moreover, we introduce the linear and bounded operator $\mathcal A^h:\U\to\Y$ as follows: for given $u\in\Uad$ the function $p^h=\mathcal A^hu$ is the unique solution to
\begin{subequations}
    \label{MD-6}
    \begin{align}
        \label{MD-6a}
        -\frac{\mathrm d}{\mathrm dt}\,{\langle p^h(t),\varphi^h\rangle}_H+a(t;\varphi^h,p^h(t))&=-\sigma_1\,{\langle(\mathcal S^hu)(t),\varphi^h\rangle}_H\quad\text{for all }\varphi\in V^h\text{ a.e. in }[0,T),\\
        \label{MD-6b}
        p^h(T)&=-\sigma_2(\mathcal S^hu)(T).
    \end{align}
\end{subequations}
Then the solution to \eqref{MD-2-2} can be expressed as $\bar p=\hat p+\mathcal A^h\bar u$.

\begin{example} \label{Example:ch4_IntroducingOptimalControlProblem}
    \rm Now we introduce an example for an optimal control problem, which will be used for the numerical tests throughout this section. For the PDE we take our guiding model, which was described in Example~\ref{Example:Ch3_aprioriestimates_example} together with the initial state $y_\circ(\bx) = 17$. The reduced cost function $\hat{J}^h: \Uad \to \mathbb{R}$ is given by 
    \begin{equation}
        \label{Equation:ch4_ExampleOptimalControlProblem_CostFunction}
        \begin{aligned}
            \hat{J}^h(u)=& \frac{\sigma_1}{2}\int_0^T{\|(S^hu)(t) + \hat{y}^h(t) -\ydQ(t)\|}_H^2\,\mathrm dt\\
            &\quad+\frac{\sigma_2}{2}\,{\|(S^hu)(T) + \hat{y}^h(T)-\ydT\|}_H^2+\frac{\sigma}{2}{\|u-\un\|}_\U^2,
        \end{aligned}
    \end{equation}
    where we set $\sigma_1 = 1$, $\sigma_2 = 0.1$ and $\sigma = 0.001$. Moreover, we choose the desired states $\ydQ(t) := 17$ for all $t \in (0,T)$ and $\ydT := 17$, as well as the desired control $\un(t) = 0$ for all $t \in (0,T)$. In other words the goal of our optimal control problem is to maintain the constant initial temperature of $17^\circ$ in the entire room in spite of the heat exchange with the outer world and the advection in the room. \hfill$\blacklozenge$ 
\end{example}
 
\subsection{POD basis computation}
\label{SIAM-Book:Section4.3.2}

Suppose that Assumptions~\ref{A8} and \ref{A9} hold. Let $u\in\Uad$ be given, $y^h=\hat y^h+\mathcal S^hu$ and $p^h=\hat p^h+\mathcal A^hu$. Then we have $y^h,p^h\in H^1(0,T;V)$; cf. Assumption~\ref{A10}-1). In the context of Section~\ref{Section:ContPODHilbert} we choose $K=2$, $\omega_1^K=\omega_2^K=1$, $y^1=y^h$ and $y^2=p^h$ so that the \index{POD method!snapshot space}{\em snapshot space} is given by
\begin{equation}
    \label{HI-215}
    \mathscr V=\mathrm{span}\,\bigg\{\int_0^T\phi_1(t)y^h(t)+\phi_2(t)p^h(t)\,\mathrm dt\,\Big|\,\phi_1,\phi_2\in L^2(0,T)\bigg\}\subset V,
\end{equation}
with finite dimension $d=\dim\mathscr V\le m$. Therefore, it follows from \cite[Lemma~3.1-2)]{Sin14} that $d$ is the same for the choices $X=H$ and $X=V$. Suppose that $\{\psi_i\}_{i=1}^\ell\subset V^h$ is a POD basis of rank $\ell\in\{1,\ldots,d\}$ solving
\begin{equation}
    \label{POD:MD-7}
    \left\{
    \begin{aligned}
        & \min\int_0^T \Big\|y^h(t) - \sum_{i=1}^\ell {\langle y^h(t),\psi_i\rangle}_X\,\psi_i\Big\|_X^2+\Big\|p^h(t) - \sum_{i=1}^\ell {\langle p^h(t),\psi_i\rangle}_X\,\psi_i\Big\|_X^2\,\mathrm dt\\
        &\hspace{1mm}\text{s.t. }\{\psi_i\}_{i=1}^\ell\subset X\text{ and }{\langle\psi_i,\psi_j\rangle}_X=\delta_{ij} \quad \text{ for } 1 \le i,j \le \ell.
    \end{aligned}
    \right.
\end{equation}
Recall that the solution to \eqref{POD:MD-7} is given by the solution to the eigenvalue problem \eqref{FEPODEigPro}, where the operator $\mathcal R:X\to V^h\subset X$ is given as
\begin{align*}
    \mathcal R\psi=\int_0^T\Big({\langle y^h(t),\psi\rangle}_X\,y^h(t)+{\langle p^h(t),\psi\rangle}_X\,p^h(t)\Big)\,\mathrm dt\quad \text{for }\psi\in X
\end{align*}
In particular, $\lambda_i=0$ for all $i>d$; cf. Assumption~\ref{A10}-3). Then we obtain
\begin{equation}
    \label{MD-8}
    \int_0^T \Big\|y^h(t) - \sum_{i=1}^\ell {\langle y^h(t),\psi_i\rangle}_X\,\psi_i\Big\|_X^2 +\Big\|p^h(t) - \sum_{i=1}^\ell {\langle p^h(t),\psi_i\rangle}_X\,\psi_i\Big\|_X^2\,\mathrm dt =\sum_{i=\ell+1}^d\lambda_i.
\end{equation}
We set
\begin{align*}
    X^\ell=\mathrm{span}\,\big\{\psi_1,\ldots,\psi_\ell\big\}\subset V\subset X
\end{align*}
and
\begin{align*}
    H^\ell=\mathrm{span}\,\big\{\psi_1^H,\ldots,\psi_\ell^H\big\},\quad V^\ell=\mathrm{span}\,\big\{\psi_1^V,\ldots,\psi_\ell^V\big\}
\end{align*}
for the choices $X=H$ and $X=V$, respectively. For the reader's convenience we omitted the dependence of the eigenvectors $\{\psi_i^H\}_{i=1}^\ell$ and $\{\psi_i^V\}_{i=1}^\ell$ on the spatial discretization by an index $h$. Finally, we replace Assumption~\ref{A10} by the following one.

\begin{assumption}
    \label{A11}
    \begin{enumerate}
        \item [\em 1)] Let Assumptions~{\em\ref{A8}} and {\em\ref{A9}} hold.
        \item [\em 2)] The Hilbert space $X$ denotes either $H$ or $V$. In case of $X=V$ we suppose that $y_\circ\in V$ is valid.
        \item [\em 3)] We set $K=2$, $\omega_1^K=\omega_2^K=1$ and $y^1=y^h$, $y^2=p^h$. The snapshot space $\mathscr V$ is given by \eqref{HI-215} with finite dimension $d\le m$. Let $\{\psi_i^H\}_{i=1}^\ell$ and $\{\psi_i^V\}_{i=1}^\ell$ be solutions to \eqref{PellH} and \eqref{PellV}, respectively.
        \item [\em 4)] In case of $X=H$, let $\mathcal P^{h\ell}$ denote $\mathcal P^{h\ell}_H$. If $X=V$, let it be given by $\mathcal Q^{h\ell}_H$.
    \end{enumerate}
\end{assumption}

\subsection{POD Galerkin scheme for the dual equation}
\label{SIAM-Book:Section4.3.3}

Suppose that  Assumption~\ref{A11} hold. In Section~\ref{SIAM-Book:Section3.4.4} we have introduced a POD Galerkin scheme for the state equation. From Theorem~\ref{DWW-1} and $\|\mathcal P^\ell y_\circ\|_H\le\|y_\circ\|_H$ we infer that
\begin{equation}
\label{HI-220}
{\|y^{h\ell}\|}_\Y\le C\left({\|y_\circ\|}_H+{\|\mathcal F\|}_{L^2(0,T;V')}+{\|u\|}_\U\right)
\end{equation}
for a constant $C>0$ which is independent of $y_\circ$, $\mathcal F$, $u$ and $\ell$. Due to \eqref{Chopin-2} there exists a constant $C>0$ satisfying the \index{Error estimate!a-priori!state variable}a-priori error estimate
\begin{align*}
    \int_0^T{\|y^h(t)-y^{h\ell}(t)\|}_V^2\,\mathrm dt\le C\cdot\left\{
    \begin{aligned}
        &\sum_{i=\ell+1}^d\lambda_i^H\,\big\|\psi_i^H\big\|_V^2&&\text{for }X=H,\\
        &\sum_{i=\ell+1}^d\lambda_i^V\,\big\|\psi_i^V-\mathcal Q_H^{h\ell}\psi_i^V\big\|_V^2&&\text{for }X=V,
    \end{aligned}
    \right.
\end{align*}
where $y^h=\hat y^h+\mathcal S^hu$ and $y^{h\ell}=\hat y^{h\ell}+\mathcal S^{h\ell}u$ denote the solutions to \eqref{EvProGal} and \eqref{EvProGal-POD}, respectively.

Analogously, we proceed for the dual equation. Making the ansatz
\begin{equation}
    \label{DualGalAn-2}
    p^{h\ell}(t)=\sum_{i=1}^\ell\mathrm p_i^{h\ell}(t)\psi_i^h \in X^{h\ell},\quad t\in[0,T]
\end{equation}
and plugging this into \eqref{MD-2-2} leads to the following POD Galerkin scheme:
\begin{subequations}
    \label{LQR:Dual-1}
    \begin{align}
        \label{LQR:Dual-1a}
        -\frac{\mathrm d}{\mathrm dt}\,{\langle p^{h\ell}(t),\psi\rangle}_H+a(t;\psi,p^{h\ell}(t))&=\sigma_1\,{\langle\ydQ(t)-y^{h\ell}(t),\psi\rangle}_{V',V}\text{ for all }\psi\in X^\ell\text{ a.e. in }[0,T),\\
        \label{LQR:Dual-1b}
        \hat p^{h\ell}(T)&=\sigma_2\big(\mathcal P^{h\ell}\ydT-y^{h\ell}(T)\big)
    \end{align}
\end{subequations}
with $y^{h\ell}=\hat y^{h\ell}+\mathcal S^{h\ell}u$. Let us also define the $\ell$-dimensional vectors
\begin{align*}
    \mathrm p^{h\ell}(t)=\big(\mathrm p_i^{h\ell}(t)\big)_{1\le i\le \ell},\quad\mathrm y_1^{h\ell}=\big({\langle\ydQ(t),\psi_i\rangle}_H\big)_{1\le i\le\ell},\quad\mathrm y_2^{h\ell}=\big(\mathrm y_{2 i}^{h\ell}\big)_{1\le i\le \ell},
\end{align*}
where
\begin{align*}
    \mathcal P^{h\ell}\ydT=\sum_{i=1}^\ell\mathrm y_{2 i}^{h\ell}\psi_i\in X^{h\ell}.
\end{align*}
In Section~\ref{SIAM-Book:Section3.4.4} we have introduced the $(\ell\times\ell)$-matrices $\mathrm M^{h\ell}$ and $\mathrm A^{h\ell}(t)$, $t\in[0,T]$. Then \eqref{LQR:Dual-1} can be expressed as the following system of ordinary differential equations
\begin{equation}
    \label{LQR:Dual-2}
        \begin{aligned}
        -\bM^{h\ell}\dot{\mathrm p}^{h\ell}(t)+\bA^{h\ell}(t)^\top\mathrm p^{h\ell}(t)&=\sigma_1\big(\mathrm y_1^{h\ell}-\bM^{h\ell}\mathrm y^{h\ell}(t)\big),\quad t\in[0,T),\\
        \mathrm p^{h\ell}(T)&=\sigma_2\big(\mathrm y_2^{h\ell}-\mathrm y^{h\ell}(T)\big),
    \end{aligned}
\end{equation}
where the vector $\mathrm y^{h\ell}(t)$, $t\in[0,T]$, solves \eqref{FE-POD-ODE}. Note that $\tilde p^{h\ell}(t)=p^{h\ell}(T-t)$, $t\in[0,T]$, satisfies
\begin{subequations}
    \label{LQR:Dual-3}
    \begin{align}
        \label{LQR:Dual-3a}
        \frac{\mathrm d}{\mathrm dt} {\langle \tilde p^{h\ell}(t),\psi^h\rangle}_H+\tilde a(t;\psi^h,\tilde p^{h\ell}(t))&=\sigma_1\,{\langle\ydQ(t)-y^{h\ell}(T-t),\psi^h\rangle}_H\text{ for all }\psi^h\in X^\ell\text{ a.e. in }(0,T],\\
        \label{LQR:Dual-3b}
        \tilde p^{h\ell}(T)&=\sigma_2\big(\mathcal P^{h\ell}\ydT-y^{h\ell}(T)\big).
    \end{align}
\end{subequations}
Now, existence and uniqueness of a solution $\tilde p$ to \eqref{LQR:Dual-3} can be shown by similar arguments as in the proof of Theorem~\ref{DWW-1}.  Thus, \eqref{LQR:Dual-1} is also uniquely solvable and $p^{h\ell}(t)=\tilde p^{h\ell}(T-t)$, $t\in[0,T]$, is the unique solution. Summarizing, we obtain the next result.

\begin{proposition}
    \label{Prop:LQR-1}
    Let Assumption~{\rm\ref{A11}} hold. Then  for every $u\in\U$ there exists a unique $p^{h\ell}\in H^1(0,T;V)\hookrightarrow\Y$ solving \eqref{LQR:Dual-1} and satisfying
    \begin{align*}
        {\|p^{h\ell}\|}_\Y\le C\big(\sigma_2\,{\|\ydT-y^{h\ell}(T)\|}_H+\sigma_1\,{\|\ydQ-y^{h\ell}\|}_{L^2(0,T;H)}\big)
    \end{align*}
    for a constant $C\ge0$ which is independent of $\ell$ and $y^{h\ell}=\hat y^{h\ell}+\mathcal S^{h\ell} u$.
\end{proposition}

The following corollary follows by analogous arguments utilized in the proof of Corollary~\ref{Corollary:HI-40} and in Remark~\ref{Remark:KoH-1}-1).

\begin{corollary}
    \label{Corollary:LQR-1}
    Let Assumption~{\rm\ref{A11}} hold.
    \begin{enumerate}
        \item [\rm 1)] There exists a unique solution $\hat p^{h\ell} \in H^1(0,T;V)\hookrightarrow\Y$ to
        \begin{equation}
            \label{LQR:Dual-0}
            \begin{aligned}
                \hspace{-5mm}-\frac{\mathrm d}{\mathrm dt}\,{\langle\hat p^{h\ell}(t),\psi\rangle}_H+a(t;\psi,\hat p^{h\ell}(t))&=\sigma_1\,{\langle\ydQ(t)-\hat y(t),\psi\rangle}_{V',V}\text{ for all }\psi\in X^\ell\text{ a.e. in }[0,T),\\
                \hat p^{h\ell}(T)&=\sigma_2\mathcal P^\ell\ydT-\hat y^{h\ell}(T).
            \end{aligned}
        \end{equation}
        In particular, we have
        \begin{align*}
            {\|\hat p^{h\ell}\|}_\Y\le C\big(\sigma_2\,{\|\ydT-\hat y^{h\ell}(T)\|}_H+\sigma_1\,{\|\ydQ-\hat y^{h\ell}\|}_{L^2(0,T;H)}\big)
        \end{align*}
        for a constant $C\ge0$ which is independent of $\ell$.
        \item [\rm 2)] Let us introduce the linear operator $\mathcal A^{h\ell}:\U\to \Y$ as follows: $p^{h\ell}=\mathcal A^{h\ell}u$, $u\in\U$, satisfies
        \begin{align*}
            -\frac{\mathrm d}{\mathrm dt}\,{\langle p^{h\ell}(t),\psi\rangle}_H+a(t;\psi,p^{h\ell}(t))&=-\sigma_1\,{\langle (\mathcal S^{h\ell}u)(t),\psi\rangle}_{V',V}\quad\text{for all }\psi\in X^\ell\text{ a.e. in }[0,T),\\
            p^{h\ell}(T)&=-\sigma_2(\mathcal S^{h\ell}u)(T).
        \end{align*}
        Then $\mathcal A^{h\ell}$ is well-defined and bounded. In particular,
        \begin{align*}
            {\|\mathcal A^{h\ell}u\|}_\Y\le C\left(\sigma_2\,{\|(\mathcal S^{h\ell}u)(T)\|}_H+\sigma_1\,{\|\mathcal S^{h\ell}u\|}_{L^2(0,T;H)}\right)
        \end{align*}
        for a constant $C\ge0$ which is independent of $\ell$.
        \item [\rm 3)] The solution to \eqref{LQR:Dual-1} can be expressed as $p^{h\ell}=\hat p^{h\ell}+\mathcal A^{h\ell} u$.
    \end{enumerate}
\end{corollary}

The next result is proved in Section~\ref{SIAM-Book:Section4.7.3}.

\begin{corollary}
    \label{Cor:LQR-1}
    Let all assumptions of Proposition~{\rm\ref{Prop:LQR-1}} hold. Then
    \begin{align*}
        {\|p^{h\ell}\|}_\Y\le C\,\big({\|\ydT\|}_H+{\|\ydQ\|}_{L^2(0,T;H)}+{\|\mathcal P^{h\ell} y_\circ\|}_H+{\|\mathcal F\|}_{L^2(0,T;V')}+{\|u\|}_\U\big)
    \end{align*}
    for a constant $C\ge0$ which is independent of $\ell$.
\end{corollary}

\subsection{POD a-priori error analysis for the dual equation}
\label{SIAM-Book:Section4.3.4}

We suppose Assumption~\ref{A11} and proceed similarly to Section~\ref{SIAM-Book:Section4.2.2}. Then $\mathcal P^{h\ell}$ is an $H$-orthonormal projection for $X=H$ and $X=V$, respectively. Our next goal is to derive a-priori error estimates for the term
\begin{align*}
    \int_0^T {\|p^h(t)-p^{h\ell}(t)\|}_V^2\,\mathrm dt,
\end{align*}
where $p^h$ and $p^{h\ell}$ are the solutions to \eqref{MD-2-2} and \eqref{LQR:Dual-1}, respectively. Following the arguments of the proof of Theorem~\ref{Th:DualA-PrioriError-2} we derive the following result.

\begin{theorem}
    \label{Th:LQR-1}
    Assume \index{Error estimate!a-priori!dual variable}that Assumption~{\rm\ref{A11}} hold. Furthermore, $p^\ell$ and $p^{h\ell}$ denote the unique solutions to \eqref{MD-2-2} and \eqref{LQR:Dual-1}, respectively. Then  there exists a constant $C>0$ satisfying the \index{Error estimate!a-priori!dual variable}a-priori error estimate
    \begin{align*}
        &\int_0^T\big\|p^h(t)-p^{h\ell}(t)\big\|_V^2\,\mathrm dt\le C\cdot\left\{
        \begin{aligned}
            &\sum_{i=\ell+1}^d\lambda_i^H\,\big\|\psi_i^H\big\|_V^2&&\text{for }X=H,\\
            &\sum_{i=\ell+1}^d\lambda_i^V\,\big\|\psi_i^V-\mathcal Q_H^{h\ell}\psi_i^V\big\|_V^2&&\text{for }X=V.
        \end{aligned}
        \right.
    \end{align*}
\end{theorem}

\begin{example} \label{Example:Ch4_aprioriestimates_adjointequation}
    \rm Let us now verify the statement of Theorem~\ref{Th:LQR-1} numerically. To do so we use the guiding example for this section introduced in Example~\ref{Example:ch4_IntroducingOptimalControlProblem}. For a random control $u$ we solve the associated state and dual equations \eqref{MD-3} and \eqref{MD-2-2}, respectively, to obtain the state $y^h(u)$ and the dual state $p^h(u)$. Following the procedure presented in Section~\ref{SIAM-Book:Section4.3.2} a POD basis $\left\{ \Psi_i \right\}_{i=1}^\ell \subset V^h$ is computed using the snapshots $y^h(u)$ and $p^h(u)$. For this POD basis we can then compute the solution $p^{h\ell}(u)$ of the POD dual equation \eqref{LQR:Dual-1}. 
    \begin{figure}
        \begin{center}
            \includegraphics[height=50mm]{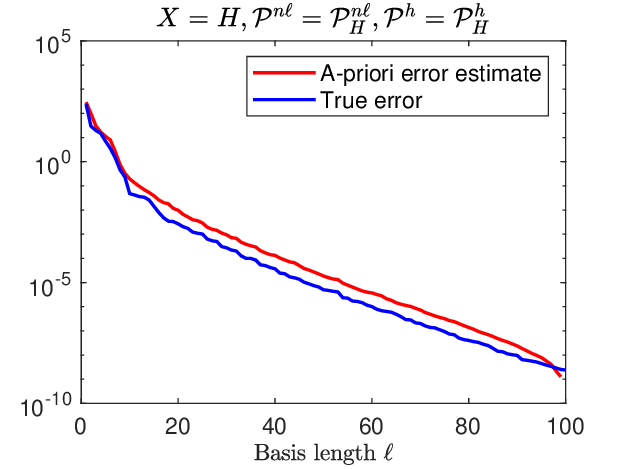}\hspace{10mm}
            \includegraphics[height=50mm]{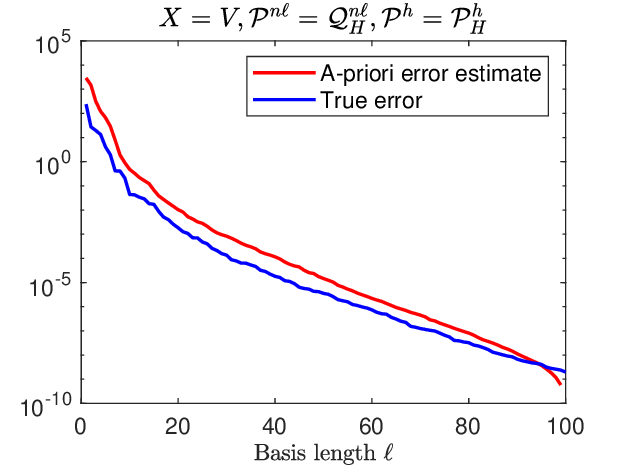}
        \end{center}
        \caption{Example~\ref{Example:Ch4_aprioriestimates_adjointequation}. Plot of the a-priori error estimate and the approximation error of the dual states for different choices of spaces and projections (y-axis: logarithmic scale).}
        \label{fig:Ch4_aprioriestimates_adjointequation}
    \end{figure}
    In Figure \ref{fig:Ch4_aprioriestimates_adjointequation} the error between $p^h(u)$ and $p^{h\ell}(u)$ is depicted for a varying number of POD basis elements. As expected the error decays when increasing the number of basis elements. At the same time it can be observed that the a-priori error estimates for both choices of the space and the projection predict the qualitative behaviour of the error nicely. Let us note again that we plot the a-priori estimates with the constant $C = 1$, since we are only interested in obtaining a qualitative but not a quantitative estimate of the error.\hfill$\blacklozenge$
\end{example}

\subsection{POD Galerkin scheme for the optimality system}
\label{SIAM-Book:Section4.3.5}

The POD Galerkin approximation for \eqref{PhatFE} is as follows:
\begin{equation}
    \label{PhatFEell}
    \tag{$\mathbf{\hat P}^{h\ell}$}
    \min\hat J^{h\ell}(u)\quad\text{s.t.}\quad u\in\Uad,
\end{equation}
where the objective is defined as $\hat J^{h\ell}(u)=J(\hat y^{h\ell}+\mathcal S^{h\ell}u,u)$ for $u\in\U$. Note that \eqref{PhatFEell} is a convex quadratic programming problem. An optimal solution $\bar u^{h\ell}\in\Uad$ to \eqref{PhatFEell} is characterized by the first-order sufficient optimality condition
\begin{equation}
    \label{LQR:Dual-7}
    {\langle\sigma(\bar u^{h\ell}-\un)-\mathcal B'\bar p^{h\ell},u-\bar u^{h\ell}\rangle}_\U\ge0\quad\text{for all }u\in\Uad
\end{equation}
with $\bar p^{h\ell}=\hat p^{h\ell}+\mathcal A^{h\ell}\bar u^{h\ell}$. In \eqref{LQR:Dual-7} we have utilized that the gradient $\nabla\hat J^{h\ell}$ of the reduced cost functional $\hat J^{h\ell}$ is given by the formula
\begin{equation}
    \label{GradientJhat-hell}
    \nabla\hat J^{h\ell}(u)=\sigma(u-\un)-\mathcal B'(\hat p^{h\ell}+\mathcal A^{h\ell}u)\quad\text{for }u\in\Uad.
\end{equation}
Since Theorem~\ref{GVLuminy:Theorem4.4.2} can be extended easily to the semidiscrete case, we derive the following result.

\begin{theorem}
    \label{Th:LQR-2}
    Assume that Assumption~{\rm\ref{A11}} \index{Error estimate!a-priori!control variable}hold. Let $\bar u^h\in\Uad$ be the optimal solution to \eqref{GVLuminy:Eq4.4.13}. Moreover, $\bar u^{h\ell}$ denotes the unique solution to \eqref{LQR:Dual-7}. Then there is a constant $C>0$ satisfying the \index{Error estimate!a-priori!control variable}a-priori error estimate
    \begin{equation}
        \label{GVLuminy:Eq4.4.15-2}
        \big\|\bar u^h-\bar u^{h\ell}\big\|_\U^2\le C\cdot\left\{
        \begin{aligned}
            &\sum_{i=\ell+1}^d\lambda_i^H\,\big\|\psi_i^H\big\|_V^2&&\text{for }X=H,\\
            &\sum_{i=\ell+1}^d\lambda_i^V\,\big\|\psi_i^V-\mathcal Q_H^{h\ell}\psi_i^V\big\|_V^2&&\text{for }X=V,
        \end{aligned}\right.
    \end{equation}
    where the constant $C$ depends on $c_V$, $\gamma$, $\gamma_1$, $T$, $\sigma_1$, $\sigma_2$ and $\|\mathcal B'\|_{\mathscr L(L^2(0,T;V),\U)}$.
\end{theorem}

\begin{example} \label{Example:ch4_AprioriEstimate_OptimalitySystem}
    \rm For the optimal control problem from Example~\ref{Example:ch4_IntroducingOptimalControlProblem} we investigate the statement of Theorem~\ref{Th:LQR-2} numerically. Therefore, we compute the solution $\bar{u}^h$ to \eqref{GVLuminy:Eq4.4.13}. For this control we compute the associated state $y^h(\bar{u})$ and dual state $p^h(\bar{u})$ by solving the state equation \eqref{MD-3} and the dual state equation \eqref{MD-2-2}, respectively, and use these as snapshots for the computation of a POD basis. In Figure~\ref{fig:Ch4_aprioriestimates_controlerror} the error between $u^h$ and $u^{h\ell}$ together with the a-priori error estimates for the two choices of the space and the projection can be seen.
    \begin{figure}
        \begin{center}
            \includegraphics[height=50mm]{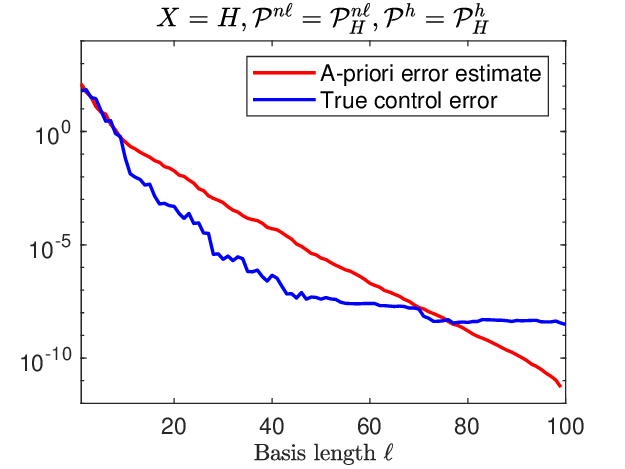}\hspace{10mm}
            \includegraphics[height=50mm]{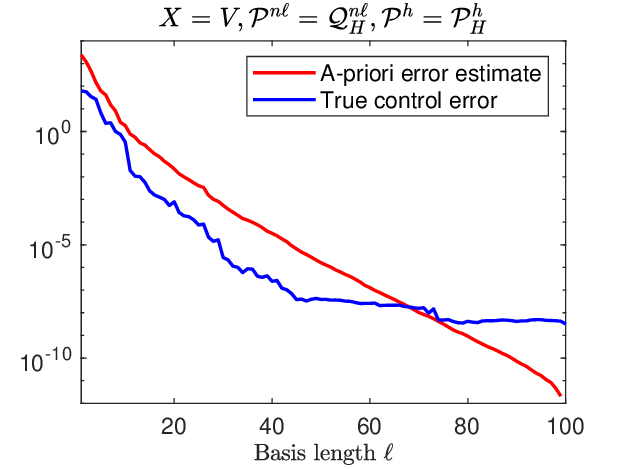}
        \end{center}
        \caption{Example~\ref{Example:ch4_AprioriEstimate_OptimalitySystem}. Plot of the a-priori error estimate and the approximation error $\|\bar u^h-\bar u^{h\ell}\|_\U^2$ between the optimal controls for different choices of spaces and projections (y-axis: logarithmic scale).}
        \label{fig:Ch4_aprioriestimates_controlerror}
    \end{figure} 
    In both cases the a-priori estimates describe the qualitative behaviour of the error quite well up to a basis size of around 50. When further enlarging the POD basis a stagnation of the error in the controls can be observed, which can be explained with the predefined accuracy with which we solve the two optimal control problems. Hence, for the given stopping tolerances for the optimal control problems it is not meaningful to use more than 50 basis functions, since the error in the controls already reaches the computational accuracy at this point. To underline the importance of a good choice of snapshots for the POD basis computation, we conduct the same experiment again with two different choices of snapshots. First, we take the snapshots $y^h(u)$ and $p^h(u)$ which are associated to a randomly chosen control $u \in \Uad$. Secondly, we use the snapshots $y^h(\tilde{u})$ and $p^h(\tilde{u})$, where $\tilde{u} \in \Uad$ is a suboptimal control, which is obtained by only computing a few iterations in the optimization process. The results of this experiment can be seen in Figure~\ref{fig:Ch4_controlerror_differentsnapshots}.
    \begin{figure}
        \begin{center}
            \includegraphics[height=50mm]{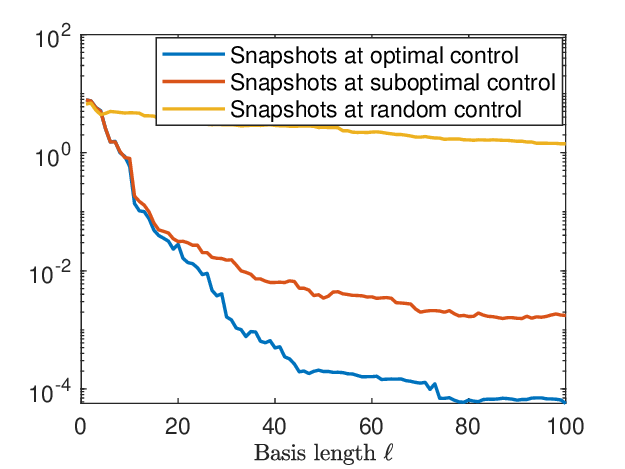}
        \end{center}
        \caption{Example~\ref{Example:ch4_AprioriEstimate_OptimalitySystem}. Plot of the approximation error $\|\bar u^h-\bar u^{h\ell}\|_\U$ between the optimal controls for different choices of snapshots (y-axis: logarithmic scale).}
        \label{fig:Ch4_controlerror_differentsnapshots}
    \end{figure} 
    As expected the snapshots taken from a randomly chosen control do not contain the dynamics close to the optimal control, so that the approximation error is large even for a large number of POD basis functions. Taking a suboptimal control for generating the snapshots yields already better results. In fact, for a basis size of up to 20 the approximation error is of the same order as for the optimal snapshots. Only when further increasing the number of basis function we can see a clear difference between the optimal and the suboptimal snapshots.\hfill$\blacklozenge$
\end{example}

\subsection{POD a-posteriori error analysis for the reduced gradient}
\label{SIAM-Book:Section4.3.7}

Now we want to control the error in the gradient of the reduced cost functional by an \index{Error estimate!a-posteriori!gradient}a-posteriori error estiomate. Let $u\in\Uad$ be chosen arbitrarily. Then we have from \eqref{GradientJhat-h} and \eqref{GradientJhat-hell} that
\begin{equation}
    \label{GradCostApo-0}
    \nabla J^h(u)-\nabla\hat J^{h\ell}(u)=\mathcal B'\big(p^{h\ell}(u)-p^h(u)\big)
\end{equation}
with the dual variables $p^h(u)=\hat p^h+\mathcal A^hu$ and $p^{h\ell}(u)=\hat p^{h\ell}+\mathcal A^{h\ell}u$. Thus, an a-posteriori error estimate for the difference $\nabla J^h(u)-\nabla\hat J^{h\ell}(u)$ is based on the availability of an \index{Error estimate!a-posteriori!dual variable}a-posteriori error estimate for the difference $p^{h\ell}(u)-p^h(u)$. This is stated in the next proposition.

\begin{proposition}
    \label{Prop1:GradCostApo}
    Suppose that Assumptions~{\rm\ref{A1}} and {\rm\ref{A9}} hold. For an arbitrarily given $u\in\U$ let $p^h(u)=\hat p^h+\mathcal A^hu$ and $p^{h\ell}(u)=\hat p^{h\ell}+\mathcal A^{h\ell}u$. Then the \index{Error estimate!a-posteriori!dual variable}a-posteriori error estimate
    \begin{equation}
        \label{GradCostApo-1}
        \begin{aligned}
            {\|p^h(t;u)-p^{h\ell}(t;u)\|}_H^2&\le\mathsf R_T^{h\ell}+\mathsf C_2^{h\ell}\bigg(\mathsf R_\circ^{h\ell}+\int_0^T\mathsf R_1^{h\ell}(s)\,\mathrm ds\bigg)+\frac{2}{\gamma_1}\int_t^T{\|\mathsf r^{h\ell}(s;u)\|}_{(V^h)'}^2\,\mathrm ds\\
            &\quad+\int_t^T\mathsf C_3^{h\ell}(s)\bigg(\mathsf R_\circ^{h\ell}+\int_0^s\mathsf R_1^{h\ell}(\tau)\,\mathrm d\tau\bigg)\,\mathrm ds
        \end{aligned}
    \end{equation}
    is satisfied for almost all $t\in[0,T]$, where
    \begin{equation}
        \label{GradCostApo-2}
        \mathsf R_T^{h\ell}=2\sigma_2\,{\|(\mathcal P^h-\mathcal P^{h\ell})\ydT\|}_H^2,~\mathsf C_2^{h\ell}=2\sigma_2\,\mathsf C_1^{h\ell}(T),~\mathsf C_3^{h\ell}(s)=\frac{2\sigma_1^2c_V^2}{\gamma_1}\,\mathsf C_1^{h\ell}(s)
    \end{equation}
    holds, the (dual) residual $\mathsf r^{h\ell}(t)\in(V^h)'$, $t\in[0,T]$ a.e., is defined as
    \begin{equation}
        \label{GradCostApo-2a}
        \mathsf r^{h\ell}(t)=-{\langle p_t^{h\ell}(t;u),\cdot\rangle}_{V',V}+a(t;\cdot\,,p^{h\ell}(t;u))-\sigma_1\,{\langle\ydQ(t)-y^{h\ell}(t;u)\rangle}_H
    \end{equation}
    and the constants $\mathsf C_1^{h\ell}(t)$, $\mathsf R_\circ^{h\ell}$ and $\mathsf R_1^{h\ell}(s)$ are given by \eqref{Alb-4a}. Moreover, we have
    \begin{equation}
        \label{GradCostApo-3}
        \begin{aligned}
            &\int_0^T{\|p^h(t;u)-p^{h\ell}(t;u)\|}_V^2\,\mathrm dt\\
            &\quad\le\frac{1}{\gamma_1}\,\bigg(\mathsf R_T^{h\ell}+\mathsf C_2^{h\ell}\Big(\mathsf R_\circ^{h\ell}+\int_0^T\mathsf R_1^{h\ell}(t)\,\mathrm dt\Big)\bigg)+\frac{2}{\gamma_1}\,{\|\mathsf r^{h\ell}(t;u)\|}_{(V^h)'}^2\,\mathrm dt\\
            &\qquad+\frac{1}{\gamma_1}\int_0^T\mathsf C_3^{h\ell}(s)\bigg(\mathsf R_\circ^{h\ell}+\int_0^s\mathsf R_1^{h\ell}(\tau)\,\mathrm d\tau\bigg).
        \end{aligned}
    \end{equation}
\end{proposition}

From \eqref{GradCostApo-0} and \eqref{GradCostApo-3} we derive the following estimate for the gradients of the reduced cost functionals for a given control $u\in\Uad$:
\begin{equation}
    \label{GradCostApo-10}
    {\|\nabla J^h(u)-\nabla\hat J^{h\ell}(u)\|}_\U\le{\|\mathcal B'\|}_{\mathscr L(L^2(0,T;V),\U)}{\|p^{h\ell}(u)-p^h(u)\|}_{L^2(0,T;V)}\le \Delta_\mathsf{apo}^{h\ell}(u)
\end{equation}
with the a-posteriori error bound
\begin{align*}
    \Delta_\mathsf{apo}^{h\ell}(u)&={\|\mathcal B'\|}_{\mathscr L(L^2(0,T;V),\U)}\Bigg(\frac{1}{\gamma_1}\,\bigg(\mathsf R_T^{h\ell}+\mathsf C_2^{h\ell}\Big(\mathsf R_\circ^{h\ell}+\int_0^T\mathsf R_1^{h\ell}(t)\,\mathrm dt\Big)\bigg)\\
    &\hspace{36mm}+\frac{1}{\gamma_1}\int_0^T\mathsf C_3^{h\ell}(s)\bigg(\mathsf R_\circ^{h\ell}+\int_0^s\mathsf R_1^{h\ell}(\tau)\,\mathrm d\tau\bigg)+\frac{2}{\gamma_1}\,{\|\mathsf r^{h\ell}(t;u)\|}_{(V^h)'}^2\,\mathrm dt\Bigg)^{1/2}.
\end{align*}
Estimate \eqref{GradCostApo-10} is very useful in inexact optimization methods. It can be utilized to control the inexactness of the gradient of reduced-order cost functional. We refer, e.g., to \index{POD method!trust-region}trust-region POD methods, where a so-called {\em Carter condition} is utilized to ensure global convergence of the iterative trust-region method; cf. \cite{AFS00,RTV17,Sch12}.

\subsection{POD a-posteriori error analysis for the optimality system}
\label{SIAM-Book:Section4.3.6}

Proceeding as in Section~\ref{SIAM-Book:Section4.2.4} we can derive an \index{Error estimate!a-posteriori!control variable}a-posteriori error estimate for the difference $\|\bar u^h-\bar u^{h\ell}\|_\U$, where $\bar u^h$ and $\bar u^{h\ell}$ are the solutions to \eqref{PhatFE} and \eqref{PhatFEell}, respectively.

\begin{theorem}
    \label{Th:LQR-3}
    Suppose that Assumptions~{\rm\ref{A1}} and {\rm\ref{A9}} hold. Let $u\in\U$ be arbitrarily given so that $\mathcal S^{h\ell}u\neq0$ and $\mathcal A^{h\ell}u\neq0$. To compute a POD basis $\{\psi_i^h\}_{i=1}^\ell$ of rank $\ell$ we choose $K=2$, $y^1=\hat y^{h\ell}+\mathcal S^{h\ell}u$ and $y^2=\hat p^{h\ell}+\mathcal A^{h\ell}u$. Define the function $\zeta^{h\ell}\in\U$ by
    \begin{align*}
        \zeta_i^{h\ell}(t)=\left\{
        \begin{aligned}
            &-\min(0,\xi_i^{h\ell}(s))&&\text{a.e. in }\mathscr A^{h\ell}_{\mathsf ai}=\big\{s\in\mathscr D\,|\bar u^{h\ell}_i(s)=u_{\mathsf ai}(s)\big\},\\
            &-\max(0,\xi^{h\ell}_i(s))&&\text{a.e. in }\mathscr A^{h\ell}_{\mathsf bi}=\big\{s\in\mathscr D\,|\bar u^{h\ell}_i(s)=u_{\mathsf bi}(s)\big\},\\
            &-\xi_i^{h\ell}(t) && \text{a.e. in }\mathscr D\setminus\big(\mathscr A^{h\ell}_{\mathsf ai}\cup\mathscr A_{\mathsf bi}^{h\ell}\big)
        \end{aligned}
        \right.
    \end{align*}
    for $i=1,\ldots,\mathsf m$, where $\xi^{h\ell}=\sigma(\bar u^{h\ell}-\un)-\mathcal B' (\hat p^{h\ell}+\mathcal A\bar u^{h\ell})$ in $\U$. Then the a-posteriori error estimate
    \begin{equation}
        \label{GVLuminy:Eq4.5.7}
        {\|\bar u^h-\bar u^{h\ell}\|}_\U\le\frac{1}{\sigma}\,{\|\zeta^{h\ell}\|}_\U
    \end{equation}
    holds. In particular, $\lim \limits_{\ell \to \infty}\|\zeta^{h\ell}\|_\U=0$.
\end{theorem}

\begin{remark}
    \rm From a-priori error analysis it is known that
    \begin{align*}
        {\|\bar u-\bar u^h\|}_\U\le C h^\alpha
    \end{align*}
    where $\bar u$ and $\bar u^h$ are the solutions to \eqref{GVLuminy:Eq4.1.9} and \eqref{PhatFE}, respectively. Then  we find
    \begin{align*}
        {\|\bar u-\bar u^{h\ell}\|}_\U\le{\|\bar u-\bar u^h\|}_\U+{\|\bar u^h-\bar u^{h\ell}\|}_\U
    \end{align*}
    which is studied in \cite{GNV17}, for instance.\hfill$\Diamond$
\end{remark}

\begin{example}
    \label{Example:ch4_Aposteriori_OptimalitySystem}
    \rm Now we want to verify the a-posteriori error estimate from Theorem~\ref{Th:LQR-3} numerically. Again all parameters of the optimal control problem are taken from Example~\ref{Example:ch4_IntroducingOptimalControlProblem}. As we have seen in Figure~\ref{fig:Ch4_controlerror_differentsnapshots} the choice of good snapshots is essential for the approximation quality. Therefore, we take a framework that might very well appear in practice: Given a suboptimal control $\tilde{u} \in \Uad$ of the optimal control problem, we compute the state $y^h(\tilde{u})$ as well as the dual state $p^h(\tilde{u})$, which are then taken as snapshots for computing the POD basis. With this basis we solve \eqref{PhatFEell} and then compare the error $\|\bar u^h-\bar u^{h\ell}\|_\U$ with the a-posteriori error estimate from Theorem~\ref{Th:LQR-3}. The results of this test can be seen in Figure~\ref{fig:Ch4_aposterioriestimate_controlerror}. 
    \begin{figure}
        \begin{center}
            \includegraphics[height=50mm]{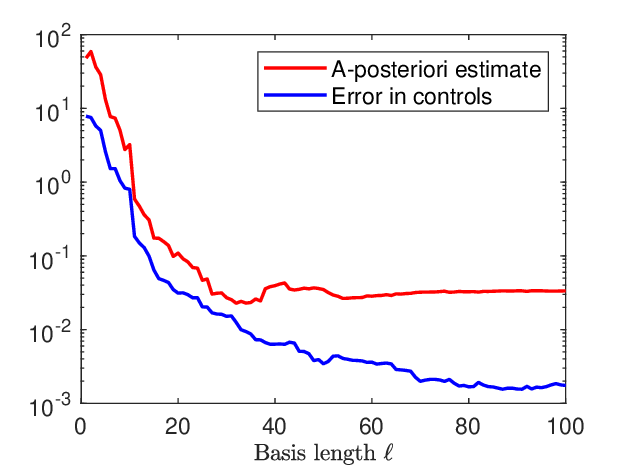}\hspace{10mm}
            \includegraphics[height=50mm]{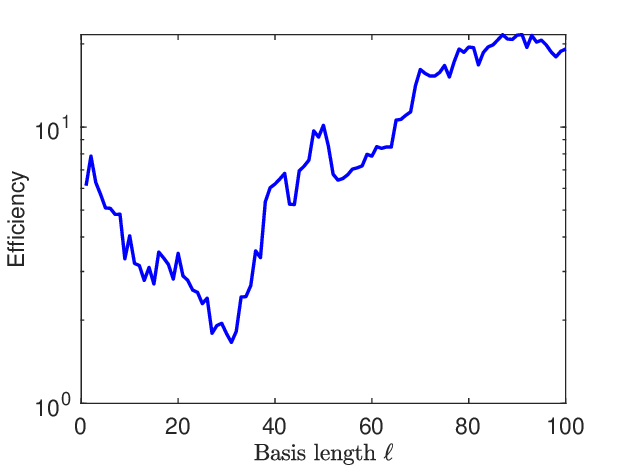}
        \end{center}
        \caption{Example~\ref{Example:ch4_Aposteriori_OptimalitySystem}. Left: Plot of the a-posteriori error estimate and the approximation error $\|\bar u^h-\bar u^{h\ell}\|_\U^2$ between the optimal controls(y-axis: logarithmic scale). Right: Efficiency of the a-posteriori error estimate.}
        \label{fig:Ch4_aposterioriestimate_controlerror}
    \end{figure} 
    Most importantly we observe that the a-posteriori error estimate is a rigorous error bound, as was shown in Theorem~\ref{Th:LQR-3}. The second crucial issue of an a-posteriori error estimate is its efficiency, i.e., the factor by which the estimate overestimates the true error. The reason for this is that one might want to apply the a-posteriori error estimate in an iterative POD basis update scheme as for example shown in Algorithm~\ref{GVLuminy:Algorithm4.5.1}. In this case it is important that the overestimation is moderate in order to avoid unnecessary basis updates. For this test case the a-posteriori error estimate turns out to be quite efficient, with the efficiency being between 2 and 20. In general, it can be seen that the efficiency gets worse if the parameter $\sigma$ in front of the control cost term $\| u^h - \un \|_\U^2$ is decreased, see Figure~\ref{fig:Ch4_aposterioriestimate_efficiencydifferentsimga}.
    \begin{figure}
        \begin{center}
            \includegraphics[height=50mm]{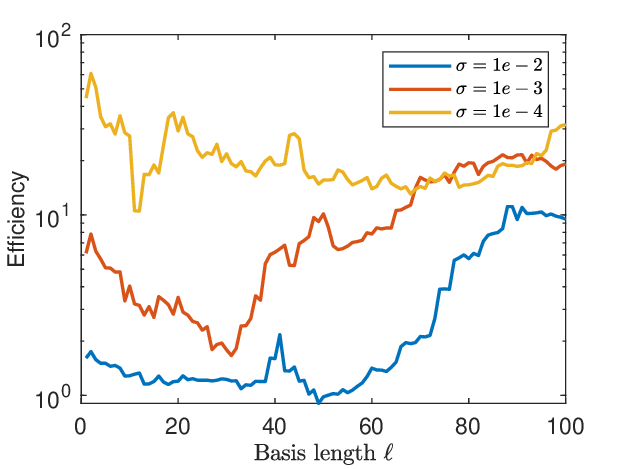}
        \end{center}
        \caption{Example~\ref{Example:ch4_Aposteriori_OptimalitySystem}. Efficiency of the a-posteriori error estimate for different values of $\sigma$ (y-axis: logarithmic scale)}
        \label{fig:Ch4_aposterioriestimate_efficiencydifferentsimga}
    \end{figure} 
    The reason for this lies in the proof of the a-posteriori estimate (cf. the lines before Theorem~\ref{GVLuminy:Theorem4.5.1}). There we made the estimate  
    \begin{align*}
        0\le-\sigma{\|\bar u-\bar u^\ell\|}^2_\U-{\langle\Theta(\bar y-\tilde y^\ell),\bar y-\tilde y^\ell\rangle}_{W_0(0,T)',W_0(0,T)}+{\langle\zeta^\ell,\bar u-\bar u^\ell\rangle}_\U\le-\sigma{\|\bar u-\bar u^\ell\|}^2_\U+{\langle\zeta^\ell,\bar u^\ell-\bar u^\ell\rangle}_\U.
    \end{align*}
    In our case it holds
    \begin{align*}
        {\langle\Theta(\bar y-\tilde y^\ell),\bar y-\tilde y^\ell\rangle}_{W_0(0,T)',W_0(0,T)} =\sigma_1 \int_0^T{\|y^h(\bar{u})(t) - y^h(\bar{u}^\ell)(t) \|}_H^2\,\mathrm dt+\sigma_2{\|y^h(\bar{u})(T) - y^h(\bar{u}^\ell)(T) \|}_H^2,
    \end{align*} 
    which is independent of $\sigma$. So especially when decreasing $\sigma$ while keeping $\sigma_1$ and $\sigma_2$ constant using the estimate ${\langle\Theta(\bar y-\tilde y^\ell),\bar y-\tilde y^\ell\rangle} \geq 0$ yields a large overestimation of the approximation error between the controls. Note that there are two values for the basis length $\ell$, for which the efficiency of the a-posteriori error estimate is smaller than one, i.e., it is no upper bound of the error. This is due to the fact that we cannot solve the optimal control problem \eqref{PhatFE} exactly, but only with respect to some tolerance. However, the a-posteriori estimate gives an upper bound of the error with respect to the exact optimal control of \eqref{PhatFE}. Therefore, it might happen that the error between $u^{h\ell}$ and the computed (inexact) optimal control $\bar{u}^h$ of \eqref{PhatFE} is larger than the a-posteriori error estimate.\hfill$\blacklozenge$
\end{example}

\section{Optimality system POD}
\label{SIAM-Book:Section4.4}
\setcounter{equation}{0}
\setcounter{theorem}{0}

The accuracy of the reduced-order model can be controlled by the a-posteriori error analysis presented, e.g., in Section~\ref{SIAM-Book:Section4.2.3}. However, if the POD basis is created from a reference trajectory containing features which are quite different from those of the optimally controlled trajectory, a rather huge number of POD basis functions has to be included in the reduced-order model. This fact may lead to non-efficient reduced-order models and numerical instabilities. To avoid these problems the POD basis is generated utilizing \index{POD method!optimality system, OS-POD}\index{Method!optimality system POD, OS-POD}{\em optimality system POD} introduced in \cite{KV08} and extended in \cite{Gri13,GGV15,Gub16,Vol11}.

Let us consider
\begin{equation}
    \label{OSPOD:Phatell}
    \tag{$\mathbf{\hat P}^\ell$}
    \min\hat J^\ell(u)\quad\text{s.t.}\quad u\in\Uad,
\end{equation}
where $\hat J^\ell(u)=J(\hat y+\mathcal S^\ell u)$ holds and $j$ has been introduced in \eqref{GVLuminy:Eq4.1.5}. Suppose that Assumption~\ref{A10} is satisfied and that for an appropriately chosen $\ell\in\mathbb N$ the POD basis $\{\psi_i\}_{i=1}^\ell\subset V$ solves the eigenvalue problem
\begin{equation}
    \label{OSPOD:EigPro}
    \mathcal R\psi_i=\lambda_i\psi_i\text{ for }i=1,\ldots,\ell,\quad\lambda_1\ge\ldots\ge\lambda_\ell>0
\end{equation}
with
\begin{align*}
    \mathcal R\psi=\int_0^T{\langle y(t),\psi\rangle}_X\,y(t)+{\langle p(t),\psi\rangle}_X\,p(t)\,\mathrm dt\quad\text{for }\psi\in X,
\end{align*}
where $y=\hat y+\mathcal Su$ and $p=\hat p+\mathcal Au$ solve \eqref{GVLuminy:Eq4.1.1} and \eqref{GVLuminy:Eq4.3.16}, respectively. Thus, the POD basis $\{\psi_i\}_{i=1}^\ell$ depends on $y$, $p$ and thus on the control $u$. This may deter from one of the main advantages of the approach for model reduction which consists in the fact that unlike typical finite element basis functions the elements of the basis functions reflect the dynamics of the system. In optimal control this feature gets lost if the dynamics of the state and the adjoint corresponding to the reference control are significantly different from that of the trajectories corresponding to the optimal control solving
\begin{equation}
    \label{OSPOD:Phat}
    \tag{$\mathbf{\hat P}$}
    \min\hat J(u)\quad\text{s.t.}\quad u\in\Uad.
\end{equation}
For the latter case we have the a-priori error analysis presented in Theorem~\ref{GVLuminy:Theorem4.4.2}. To eliminate this weakness of the conventional approach we propose to consider the following enlarged minimization problem
\begin{equation}
    \label{OSPOD:Pospod}
    \tag{$\mathbf P^\ell_\mathsf{os}$}
    \begin{aligned}
        &\min J(y^\ell_{[u]},u)\\
        &\hspace{0.5mm}\text{s.t. }y^\ell_{[u]}=\hat y^\ell+\mathcal S^\ell u,~u\in\Uad,\\
        &\hspace{0.5mm}\text{\phantom{s.t. }}y_{[u]}=\hat y+\mathcal Su,~p_{[u]}=\hat p+\mathcal Au,\\
        &\hspace{0.5mm}\text{\phantom{s.t. }}\mathcal R_{[u]}\psi_{i[u]}=\lambda_{i[u]}\psi_{i[u]}\text{ and }{\langle\psi_{i[u]},\psi_{j[u]}\rangle}_X=\delta_{ij}\text{ for }i,j=1,\ldots,\ell.
    \end{aligned}
\end{equation}
In \eqref{OSPOD:Pospod} we indicate by the subindex $[u]$ that the POD basis and that all reduced quantities depends on the control $u$. The first three lines coincide with the conventional POD problem \eqref{OSPOD:Phatell}. The next line stands for the non-reduced state and adjoint equations and the last line represents the eigenvalue problem characterizing the POD basis associated with the current control variable $u$. The optimal solution $(\bar y^\ell,\bar y,\bar p,\bar\psi,\bar\lambda,\bar u)$ to \eqref{OSPOD:Pospod} has the property that the associated POD reduced system is computed from the trajectory corresponding to the optimal control $\bar u$ and thus, differently from \eqref{OSPOD:Phatell}, the problem of unmodelled dynamics is removed. Of course, \eqref{OSPOD:Pospod} is more complicated than the original problem \eqref{OSPOD:Phat} and we thus need to justify our approach. 

Next we formulate \eqref{OSPOD:Pospod} as a constrained optimization problem. Recall that the \index{Space!state,$\Y$}{\em state space} is given as $\Y=W(0,T)$. We set
\begin{align*}
    \mathbb X^\ell=\underbrace{X\times\ldots\times X}_{\ell\text{-times}}.
\end{align*}
Let us introduce the solution space for the minimization variables as
\begin{align*}
    \X=H^1(0,T;\mathbb R^\ell)\times \Y\times \Y\times \mathbb X^\ell\times\mathbb R^\ell\times\U
\end{align*}
which is a product Hilbert space endowed with the common product topology. Next we turn to the equality constraints. Utilizing \eqref{SIAM:PODStateDis} we express the relationship $y^\ell_{[u]}=\hat y^\ell+\mathcal S^\ell u$ as an equality constraint of the form
\begin{align*}
    \big(e_1,e_2\big)(x)=0\quad\text{in }L^2(0,T;\mathbb R^\ell)'\times\mathbb R^\ell,\quad x\in\X,
\end{align*}
with the operators $e_1:\X\to\mP_1'=L^2(0,T;\mathbb R^\ell)'$, $e_2:\X\to\mP_2=\mathbb R^\ell$ given by
\begin{align*}
    {\langle e_1(x),\mathrm q^\ell\rangle}_{\mP_1',\mP_1}&=\int_0^T\big(\bM^\ell(\psi)\dot{\mathrm y}^\ell(t)+\bA^\ell(t;\psi)\mathrm y^\ell(t)-\mathrm g^\ell(t;u,\psi)\Big)^\top\mathrm q^\ell(t)\,\mathrm dt,\\
    {\langle e_2(x),\mathrm q_\circ\rangle}_{\mP_2}&=\big(\mathrm y^\ell(0)-\mathrm y_\circ^\ell(\psi)\big)^\top\mathrm q_\circ^\ell
\end{align*}
for every $(\mathrm q^\ell,q_\circ)\in \mP_1\times\mP_2$. In contrast to \eqref{SIAM:PODStateDis} -- we also indicate the dependence on the state variable $\psi$. The expression $y_{[u]}=\hat y+\mathcal Su$ is based on \eqref{GVLuminy:Eq4.1.1}. Thus, we define the operators $e_3:\X\to\mP_3'=L^2(0,T;V)'$, $e_4:\X\to\mP_4=H$ by
\begin{align*}
    {\langle e_3(x),q\rangle}_{\mP_3',\mP_3}&=\int_0^T{\langle(y_t-\mathcal F-\mathcal Bu)(t),q(t)\rangle}_{V',V}+a(t;y(t),q(t))\,\mathrm dt,\\
    {\langle e_4(x),q_\circ\rangle}_{\mP_4}&={\langle y(0)-y_\circ,q_\circ\rangle}_H
\end{align*}
for every $(q,q_\circ)\in\mP_3\times\mP_4$. The adjoint variable $p_{[u]}=\hat p+\mathcal Au$ satisfies \eqref{GVLuminy:Eq4.3.16}. Therefore, we introduce the operators $e_3:\X\to\mP_5'=L^2(0,T;V)'$, $e_4:\X\to\mP_6=H$ by
\begin{align*}
    {\langle e_5(x),\xi\rangle}_{\mP_5',\mP_5}&=\int_0^T-{\langle p_t(t),\xi(t)\rangle}_{V',V}+a(t;\xi(t),p(t))-\sigma_1\,{\langle(\ydQ-y)(t),\xi(t)\rangle}_H\,\mathrm dt,\\
    {\langle e_6(x),\xi_\circ\rangle}_{\mP_6}&={\langle p(T)-\sigma_2\,(\ydT-y(T)),\xi_\circ\rangle}_H
\end{align*}
for all $(\xi,\xi_\circ)\in\mP_5\times\mP_6$. To incorporate \eqref{OSPOD:EigPro} we define the operator $e_7:\X\to\mP_7=\mathbb X^\ell$ as
\begin{align*}
    {\langle e_7(z),\mu\rangle}_{\mP_7}=\sum_{i=1}^\ell{\langle\mathcal R(y,p)\psi_i-\lambda_i\psi_i,\mu_i\rangle}_X
\end{align*}
for every $\mu\in\mP_7$. Notice, that we indicate the dependence of the operator $\mathcal R$ on the state $y$ and the dual $p$. Finally, we have to ensure that the POD basis is orthonormal. Here we make use of the following assumption.

\begin{assumption}
    \label{Ass:OSPOS-1}
    The first largest eigenvalues of $\mathcal R(y,p)$ satisfy $\lambda_1>\ldots>\lambda_\ell$ and
    \begin{align*}
        \min\big\{\lambda_\ell(\mathcal R(y,p)\,\big|\,y=\hat y+\mathcal Su,~p=\hat p+\mathcal Au\text{ for }u\in\Uad\big\}>0.
    \end{align*}
\end{assumption}

It follows from Assumption~\ref{Ass:OSPOS-1} and Corollary~\ref{Lemma2.2.1}-1) that $\langle\psi_j,\psi_i\rangle_X=\delta_{ij}$ holds for $1\le i,j\le\ell$. Hence, we only have to ensure $\|\psi_i\|_X=1$ for $1\le i\le\ell$. Thus, we define the operator $e_8:\X\to\mP_8=\mathbb X^\ell$ by
\begin{align*}
    {\langle e_8(z),\eta\rangle}_{\mP_8}=\sum_{i=1}^\ell\big({\|\psi_i\|}_X^2-1\big)\eta_i
\end{align*}
for $\eta\in\mP_8$. Setting $\mP=\mP_1\times\ldots\times\mP_8$ and $e=(e_1,\ldots,e_8):\X\to\mP'$ the equality constraints of \eqref{OSPOD:Pospod} can be written as
\begin{align*}
    e(x)=0\quad\text{in }\mP',\quad x\in\X.
\end{align*}

The inequality constraints for the control variable $u$ are taken into account by defining the admissible set
\begin{align*}
    \Xad=\big\{x=(\mathrm y^\ell,y,p,\psi,\lambda,u)\in\X\,\big|\,u\in\Uad\text{ and }e(x)=0\text{ in }\mP'\big\},
\end{align*}
where we identify the dual space $\mP'$ with the space
\begin{align*}
    \mP'\simeq \mP_1'\times\mP_2\times\mP_3'\times\mP_4\times\mP_5'\times\mP_6\times\mP_7\times\mP_8.
\end{align*}
It follows that $\Xad$ is weakly closed. Now we can express \eqref{OSPOD:Pospod} equivalently as the following infinite-dimensional constrained minimization problem:
\begin{equation}
    \tag{$\mathbf P^\ell_\mathsf{os}$}
    \min J(y^\ell,u)\quad\text{s.t.}\quad x=(\mathrm y^\ell,y,p,\psi,\lambda,u)\in\Xad.
\end{equation}

\subsection{Existence of optimal solutions}
\label{Section:4.4.1}

We make use of the following regularity assumption.

\begin{assumption}
    \label{Ass:OSPOS-2}
    The subspace
    \begin{align*}
        \big\{(y,p)\in W(0,T)\times W(0,T)\,\big|\,y=\hat y+\mathcal Su,~p=\hat p+\mathcal Au\text{ for }u\in\Uad\big\}
    \end{align*}
    is compactly embedded in $L^2(0,T;V)\times L^2(0,T;V)$. Aubin's lemma implies that
\end{assumption}

\begin{remark}
    \em
    In Example~\ref{Example:RegResult} we have presented a setting where the state space
    \begin{align*}
        \big\{y\in W(0,T)\,\big|\,y=\hat y+\mathcal Su\text{ for }u\in\Uad\big\}
    \end{align*}
    is compactly embedded into $L^2(0,T;V)$. Let Assumption~\ref{A9} hold. An analogous result for the dual variable $p$ holds provided we even have $\ydT\in V$.\hfill$\blacksquare$
\end{remark}

In the next theorem we state existence of an optimal solution to \eqref{OSPOD:Pospod}. For a proof we refer the reader to Section~\ref{SIAM-Book:Section4.7.4}.

\begin{theorem}
    \label{Theorem:OSPOD-1}
    Let Assumptions~{\rm\ref{A1}}, {\rm\ref{A9}}, {\rm\ref{Ass:OSPOS-1}} and {\rm\ref{Ass:OSPOS-2}} hold. Then \eqref{OSPOD:Pospod} admits a globally optimal solution $\bar x$.
\end{theorem}

\subsection{Practical realization of optimality system POD}
\label{Section:4.4.2}

As we have already stated, \eqref{OSPOD:Pospod} is more complicated than the original problem \eqref{OSPOD:Phat}. Thus we have to explain how optimality system can be used efficiently in a computational realization. In Sections~\ref{SIAM-Book:Section4.2.4} and \ref{SIAM-Book:Section4.3.7} we have discussed a-posteriori error estimates and presented a strategy in Algorithm~4.2.1. The goal is to improve the quality of the POD approximation by including more POD basis functions in the POD Galerkin discretization. Due to the convergence results (cf.  Theorems~\ref{GVLuminy:Theorem4.5.1} and \ref{Th:LQR-3}) we know that -- at least theoretically -- for increasing $\ell$ the error between the exact (unknown) control and the associated suboptimal POD control should tend to zero. However, since all POD matrices are usually dense, the POD approximations become more and more ill-conditioned as $\ell$ gets bigger. Moreover, we are often not interested to involve too many POD basis functions on the POD Galerkin approximation. To overcome this problem while still ensuring a sufficient accurate approximation we apply updates of the whole POD basis. Optimality system POD offers the possibility to update the POD basis with respect to the minimization goal. In Algorithm~4.4.2 we propose a strategy which combines optimality system POD with the a-posteriori error analysis. In a first phase of the algorithm a few gradient steps for \eqref{OSPOD:Pospod} are applied in order to get a ``good'' POD basis and then fix this POD basis and enlarge $\ell$ if necessary.

\bigskip
\hrule
\vspace{-7mm}
\begin{algorithm}[(OS-POD combined with a-posteriori error analysis)]
    \label{Algorithm:OSPOD}
    \vspace{-3mm}
    \hrule
    \vspace{0.5mm}
    \begin{algorithmic}[1]
        \REQUIRE Maximal number $\ell_\mathsf{max}$ of POD basis elements, $\ell<\ell_\mathsf{max}$, initial control $u^0$, and a-posteriori error tolerance $\varepsilon_\mathsf{apo}>0$;
        \STATE Determine the state $y^0=\hat y+\mathcal Su^0$ and adjoint $p^0=\hat p+\mathcal Au^0$;
        \STATE Compute a POD basis $\{\psi_i (u^0)\}_{i=1}^\ell$ utilizing the two snapshots $y^1=y$ and $y^2=p$;
        \STATE Perform $k\ge 0$ (projected) gradient steps with an Armijo line search for \eqref{OSPOD:Pospod} in order to get $u^k$ and associated POD basis $\{\psi_i(u^k)\}_{i=1}^\ell$;
        \STATE Solve \eqref{OSPOD:Phatell} for $\bar u^\ell$ utilizing the fixed POD basis $\{\psi_i(u^k)\}_{i=1}^\ell$;
        \STATE Compute the perturbation $\zeta=\zeta(\bar u^\ell)$ as explained in Section~\ref{SIAM-Book:Section4.2.4};
        \IF{$\|\zeta\|_\U/\sigma >\varepsilon_\mathsf{apo}$ {\bf and} $\ell<\ell_\mathrm{max}$}
            \STATE Enlarge $\ell$ and go back to step 4;\hfill\\
        \ENDIF
    \end{algorithmic}
    \hrule
\end{algorithm}

\section{Proofs of Section~\ref{SIAM-Book:Section4}}
\label{SIAM-Book:Section4.7}
\setcounter{equation}{0}
\setcounter{theorem}{0}

\subsection{Proofs of Section~\ref{SIAM-Book:Section4.1}}
\label{SIAM-Book:Section4.7.1}

\noindent{\bf\em Proof of Lemma~{\em\ref{Lemma:HI-100}}.}
\begin{enumerate}
    \item [\rm 1)] The cost functional $J:\X\to\mathbb R$ is convex if we have
    \begin{equation}
    \label{HI-101a}
        J(\tau x+(1-\tau)\tilde x)\le\tau J(x)+(1-\tau)J(\tilde x)
    \end{equation}
    for every $\tau\in(0,1)$, $x=(y,u)\in\X$ and $\tilde x=(\tilde y,\tilde u)\in\X$. From \eqref{GVLuminy:Eq4.1.5} we infer that
    \begin{align*}
        J(\tau x+(1-\tau)\tilde x)&= \frac{\sigma_1}{2}\int_0^T\tau^2\,{\|y(t)-\ydQ(t)\|}_H^2\,\mathrm dt+\frac{\sigma_2\tau^2}{2}\,{\|y(T)-\ydT\|}_H^2+\frac{\sigma\tau^2}{2}\,{\|u-\un\|}_\U^2\\
        &\quad+\frac{\sigma_1}{2}\int_0^T(1-\tau)^2\,{\|\tilde y(t)-\ydQ(t)\|}_H^2\,\mathrm dt+\frac{\sigma_2(1-\tau)^2}{2}\,{\|\tilde y(T)-\ydT\|}_H^2\\
        &\quad+\frac{\sigma(1-\tau)^2}{2}\,{\|\tilde u-\un\|}_\U^2\\
        &\quad+\sigma_1\tau(1-\tau)\int_0^T{\langle y(t)-\ydQ(t),\tilde y(t)-\ydQ(t)\rangle}_H\,\mathrm dt\\
        &\quad+\sigma_2\tau(1-\tau)\,{\langle y(T)-\ydT,\tilde y(T)-\ydT\rangle}_H+\sigma\tau(1-\tau)\,{\langle u-\un,\tilde u-\un\rangle}_\U.
    \end{align*}
    From
    \begin{align*}
        {\langle y(t)-\ydQ(t),\tilde y(t)-\ydQ(t)\rangle}_H\,\mathrm dt&\le\frac{1}{2}\left({\|y(t)-\ydQ(t)\|}_H^2+{\|\tilde y(t)-\ydQ(t)\|}_H^2\right)\\
        {\langle y(T)-\ydT,\tilde y(T)-\ydT\rangle}_H&\le\frac{1}{2}\left({\|y(T)-\ydT\|}_H^2+{\|\tilde y(T)-\ydT\|}_H^2\right)\\
        {\langle u-\un,\tilde u-\un\rangle}_\U&\le\frac{1}{2}\left({\|u-\un\|}_\U^2+{\|\tilde u-\un\|}_\U^2\right)
    \end{align*}
    and from $\tau^2+\tau(1-\tau)=\tau$, $(1-\tau)^2+\tau(1-\tau)=1-\tau$ we infer that
    \begin{align*}
        J(\tau x+(1-\tau)\tilde x)&\le \frac{\sigma_1}{2}\int_0^T\tau\,{\|y(t)-\ydQ(t)\|}_H^2\,\mathrm dt+\frac{\sigma_2\tau}{2}\,{\|y(T)-\ydT\|}_H^2+\frac{\sigma\tau}{2}\,{\|u-\un\|}_\U^2\\
        &\quad+\frac{\sigma_1}{2}\int_0^T(1-\tau)\,{\|\tilde y(t)-\ydQ(t)\|}_H^2\,\mathrm dt+\frac{\sigma_2(1-\tau)}{2}\,{\|\tilde y(T)-\ydT\|}_H^2\\
        &\quad+\frac{\sigma(1-\tau)}{2}\,{\|\tilde u-\un\|}_\U^2=\tau J(x)+(1-\tau)J(\tilde x)
    \end{align*}
    for every $\tau\in(0,1)$, $x=(y,u)\in\X$ and $\tilde x=(\tilde y,\tilde u)\in\X$, which is \eqref{HI-101a}.
    \item [\rm 2)] $J$ is convex by part 1) and $\Xad$ is convex by assumption. Let $\bar x \in \Xad$ be a local solution to \eqref{GVLuminy:Eq4.1.6} and assume that there is another admissible point $\hat x \in \Xad$ with $J(\hat x)<J(\bar x)$. Since $\Xad$ is convex, we have $x_\tau := \tau \hat x + (1-\tau) \bar x \in \Xad$ for every $\tau \in (0,1)$ as well. By convexity of $J$, we get
    \begin{align*}
	   J(x_\tau) \le \tau J(\hat x) + (1-\tau) J(\bar x) < J(\bar x)
    \end{align*}
    Also, $x_\tau \to \bar x$ as $\tau \searrow 0$ which contradicts the fact that $\bar x$ is a local minimizer of $J$. Therefore, it holds $J(\bar x) \le J(\hat x)$ for all $\hat x \in \Xad$, so $\bar x$ is indeed a global solution of \eqref{GVLuminy:Eq4.1.6}.\hfill$\Box$
\end{enumerate}

\noindent{\bf\em Proof of Theorem~{\em\ref{GVLuminy:Theorem4.2.2}}.}
Let us choose the Hilbert spaces $\Ha=L^2(0,T;H)\times H$ and $\tU=\U$. Moreover, $\mathcal C_1:\Y\to L^2(0,T;H)$ is the canonical embedding operator, which is linear and bounded. We define the linear operator $\mathcal C_2:\Y\to H$ by $\mathcal C_2\varphi=\varphi(\te)$ for $\varphi\in \Y$. Since $\Y$ is continuously embedded into $C([0,T];H)$, the linear operator $\mathcal C_2$ is continuous. Finally, we set
\begin{equation}
    \label{GVLuminy:Eq4.2.2}
    \mathcal G=\left(
    \begin{array}{c}
        \sqrt{\sigma_1}\,\mathcal C_1\mathcal S\\begin{align*}1mm]
        \sqrt{\sigma_2}\,\mathcal C_2\mathcal S
    \end{array}
    \right)\in \mathscr L(\tU,\Ha),\quad\yd=\left(
    \begin{array}{c}
        \sqrt{\sigma_1}\,(\ydQ-\hat y)\\[1mm]
        \sqrt{\sigma_2}\,\big(\ydT-\hat y(\te)\big)
    \end{array}
    \right)\in\Ha,\quad\tilde u^\mathsf n=\un\in\tU
\end{equation}
and $\tUad=\Uad$. Then \eqref{GVLuminy:Eq4.1.9} and \eqref{GVLuminy:Eq4.2.1} coincide. Consequently, by Theorem~\ref{GVLuminy:Theorem4.2.1} the problem \eqref{GVLuminy:Eq4.1.9} possesses a solution. If $\sigma > 0$ we can directly conclude from Theorem~\ref{GVLuminy:Theorem4.2.1} that the solution is unique. If $\sigma_1 > 0$, we can conclude from the fact that the operator $\mathcal S$ is injective (cf. Remark~\ref{Remark:HI-33}-2)) that $\mathcal{G}$ is injective as well. Thus, the solution is also unique in this case by Theorem~\ref{GVLuminy:Theorem4.2.1}.\hfill$\Box$

\medskip\noindent{\bf\em Proof of Lemma~{\em\ref{Lemma:GradRepr}}.} First we compute the directional derivative of $\mathcal J$ for a given $u\in\tU$. Let $u^\delta\in\tU$ be an arbitrary direction and $\varepsilon>0$. 
\begin{align*}
    \mathcal J(u+\varepsilon u^\delta)-\mathcal J(u)&=\frac{1}{2}\,{\|\mathcal Gu-\yd+\varepsilon\mathcal Gu^\delta\|}_\Ha^2+\frac{\sigma}{2}\,{\|u-\tilde u^\mathsf n+\varepsilon u^\delta\|}_\tU^2-\frac{1}{2}\,{\|\mathcal Gu-\yd\|}_\Ha^2+\frac{\sigma}{2}\,{\|u-\tilde u^\mathsf n\|}_\tU^2\\
    &=\frac{\varepsilon^2}{2}\,{\|\mathcal Gu^\delta\|}_\Ha^2+\varepsilon\,{\langle\mathcal Gu-\yd,\mathcal Gu^\delta\rangle}_\Ha+\frac{\varepsilon^2 \sigma}{2}\,{\|u^\delta\|}_\tU^2+\varepsilon\sigma\,{\langle u-\tilde u^\mathsf n,u^\delta\rangle}_\tU.
\end{align*}
Thus, we have
\begin{align*}
    D\mathcal J(u;u^\delta)&=\lim_{\varepsilon\searrow 0}\bigg(\frac{\varepsilon}{2}\,{\|\mathcal Gu^\delta\|}_\Ha^2+\,{\langle\mathcal Gu-\yd,\mathcal Gu^\delta \rangle}_\Ha+\frac{\varepsilon \sigma}{2}\,{\|u^\delta\|}_\tU^2+\sigma\,{\langle u-\tilde u^\mathsf n,u^\delta\rangle}_\tU\bigg)\\
    &={\langle\mathcal Gu-\yd,\mathcal Gu^\delta\rangle}_\Ha+\sigma\,{\langle u-\tilde u^\mathsf n,u^\delta\rangle}_\tU.
\end{align*}
Next we show that $D\mathcal J(u;\cdot):\tU\to\mathbb R$ is already the G\^ateaux derivative of $\mathcal J$ at $u\in\tU$. Since $\mathcal G$ is continuous on $\tU$, the mapping $D\mathcal J(u;\cdot):\tU\to\mathbb R$ is well-defined for every $u\in\tU$. Moreover, we have
\begin{align*}
    {\|D\mathcal J(u;\cdot)\|}_{\tU'}&=\sup_{\|u^\delta\|_\tU=1}D\mathcal J(u;u^\delta)=\sup_{\|u^\delta\|_\tU=1}\big({\langle\mathcal Gu-\yd,\mathcal Gu^\delta\rangle}_\Ha+\sigma\,{\langle u-\tilde u^\mathsf n,u^\delta\rangle}_\tU\big)\\
    &\le\sup_{\|u^\delta\|_\tU=1}\big({\|\mathcal Gu-\yd\|}_\Ha{\|\mathcal Gu^\delta\|}_\Ha+\sigma\,{\|u-\tilde u^\mathsf n\|}_\tU{\|u^\delta\|}_\tU\big)\\
    &\le{\|\mathcal Gu-\yd\|}_\Ha{\|\mathcal G\|}_{\mathscr L(\tU,\Ha)}+\sigma\,{\|u-\tilde u^\mathsf n\|}_\tU<\infty
\end{align*}
by assumption. Thus, $D\mathcal J(u;\cdot)\in\tU'=\mathscr L(\tU,\mathbb R)$ holds. We conclude that the Gateaux derivative of $\mathcal J$ at $u$ is given as $\mathcal J'(u)=D\mathcal J(u;\cdot)$.\hfill$\Box$

\medskip\noindent{\bf\em Proof of Theorem~{\em\ref{GVLuminy:Theorem4.3.1}}.}
Let $\bar u\in\tUad$ be a solution to \eqref{GVLuminy:Eq4.2.1} and $u\in\tUad$ be chosen arbitrarily. Since $\tUad$ is convex, the element $\bar u+\varepsilon(u-\bar u)=\varepsilon u+(1-\varepsilon)\bar u$ belongs to $\tUad$ for any $\varepsilon\in(0,1]$. Since $\bar u$ is optimal, we also have
\begin{equation}
    \label{HI-3}
    \mathcal J(\bar u+\varepsilon(u-\bar u))-\mathcal J(\bar u)\ge0\quad\text{for any }\varepsilon\in(0,1].
\end{equation}
Utilizing \eqref{HI-3} and the fact that $\mathcal J$ is G\^{a}teaux-differentiable at $\bar u$, we derive
\begin{align*}
    {\langle\nabla\mathcal J(\bar u),u-\bar u\rangle}_\tU=D\mathcal J(\bar u,u-\bar u)=\lim_{\varepsilon\searrow 0}\frac{1}{\varepsilon}\,\big(\mathcal J(\bar u+\varepsilon(u-\bar u))-\mathcal J(\bar u)\big)\ge0
\end{align*}
which is \eqref{GVLuminy:Eq4.3.6}. Let now reversely $\bar u\in\tUad$ be a control that satisfies \eqref{GVLuminy:Eq4.3.6}. As the function $\mathcal{J}$ is convex it holds for an arbitrary $u \in \tUad$
\begin{align*}
    \mathcal{J}(u) \geq \mathcal{J}(\bar{u}) + {\langle\nabla\mathcal J(\bar u),u-\bar u\rangle}_\tU \geq \mathcal{J}(\bar{u}),
\end{align*}
where the first inequality is a well-known property of convex functions and we used \eqref{GVLuminy:Eq4.3.6} for the second inequality. Thus, $\bar{u}$ is a solution to \eqref{GVLuminy:Eq4.2.1}.\hfill$\Box$

\medskip\noindent{\bf \em Proof of Lemma~{\em\ref{Lemma:HI-31}}.}
\begin{enumerate}
    \item [1)] For $t\in[0,T]$ we set
    \begin{align*}
        \tilde a(t;\varphi,\phi)=a(T-t;\phi,\varphi)\quad\text{for all }\varphi,\phi\in V.
    \end{align*}
    Clearly, for almost all $t\in[0,T]$ the bilinear form $\tilde a(t;\cdot\,,\cdot):V\times V\to\mathbb R$ satisfies
    \begin{subequations}
        \label{tildeA}
        \begin{align}
            \label{tildeA-1}
            \big|\tilde a(t;\varphi,\psi)\big|&\le\gamma\,{\|\varphi\|}_V{\|\psi\|}_V&&\text{for all }\varphi,\psi\in V,\\
            \label{tildeA-2}
            \tilde a(t;\varphi,\varphi)&\ge\gamma_1\,{\|\varphi\|}_V^2&&\text{for all }\varphi\in V
        \end{align}
    \end{subequations}
    for the constants $\gamma\ge0$ and $\gamma_1>0$ introduced in \eqref{HI-101}. Utilizing $\mathcal Su\in\Y\hookrightarrow L^2(0,T;V')$ we can argue as in the proof of Theorem~\ref{SIAM:Theorem3.1.1} and conclude that there exists a unique solution $\tilde p \in \Y$ to
    \begin{equation}
        \label{tildeA2}
        \begin{aligned}
            \frac{\mathrm d}{\mathrm dt}\,{\langle\tilde p(t),\varphi\rangle}_H+\tilde a(t;\tilde p(t),\varphi)&=-\sigma_1\,{\langle (\mathcal Su)(T-t),\varphi\rangle}_H&&\text{for all }\varphi\in V\text{ a.e. in}[0,T),\\
            \tilde p(T)&=-\sigma_2\,(\mathcal Su)(T)&&\text{in } H.
        \end{aligned}
    \end{equation}
    Moreover, a constant $C>0$ exists with
    \begin{equation}
        \label{tildeA3}
        {\|\tilde p\|}_\Y\le C\,{\|u\|}_\U.
    \end{equation}
    Clearly, $\tilde{p}$ is the solution of \eqref{tildeA2} if and only if $p$ defined by $p(t) := \tilde{p}(T-t)$ for $t\in[0,T]$ is the solution to \eqref{GVLuminy:Eq4.3.12}, i.e. \eqref{GVLuminy:Eq4.3.12} is also uniquely solvable. From \eqref{tildeA3} and $\|\tilde p\|_\Y=\|p\|_\Y$ the claim follows.
    \item [2)] We derive from \eqref{GVLuminy:Eq4.3.12}, $y=\mathcal Su$, $w=\mathcal Sv$, $p=\mathcal Av$, $y\in W_0(0,T)$ and integration by parts \eqref{eq:partialIntegration}:
    \begin{align*}
        \int_0^T {\langle (\mathcal Bu)(t),p(t)\rangle}_{V',V}\,\mathrm dt&=\int_0^T {\langle y_t(t),p(t)\rangle}_{V',V}+a(t;y(t),p(t))\,\mathrm dt\\
        &=\int_0^T-{\langle p_t(t),y(t)\rangle}_{V',V}+a(t;y(t),p(t))\,\mathrm dt+{\langle p(\te),y(\te)\rangle}_H\\
        &=-\int_0^T\sigma_1\,{\langle w(t),y(t)\rangle}_H\,\mathrm dt-\sigma_2\,{\langle w(T),y(T)\rangle}_H
    \end{align*}
    which is the claim.\hfill$\Box$
\end{enumerate}

\noindent{\bf \em Proof of Lemma~{\em\ref{Lemma:HI-30}}.}
\begin{enumerate}
    \item [1)] Utilizing the definition of $\mathcal G$, we obtain for all $h \in \tilde \U$:
    \begin{align*}
        {\langle\mathcal G u-\yd,\mathcal Gh\rangle}_\Ha &=\sigma_1\,{\langle\mathcal Su-(\ydQ-\hat y),\mathcal Sh\rangle}_{L^2(0,T;H)}+\sigma_2\,{\langle(\mathcal Su)(T)-(\ydT-\hat y(T)),(\mathcal Sh)(T)\rangle}_H\\
        &=\sigma_1\,{\langle\mathcal Su,\mathcal Sh\rangle}_{L^2(0,T;H)}+\sigma_2\,{\langle(\mathcal Su)(T),(\mathcal Sh)(T)\rangle}_H-\sigma_1\,{\langle \ydQ-\hat y,\mathcal Sh\rangle}_{L^2(0,T;H)}\\
        &\quad-\sigma_2\,{\langle \ydT-\hat y(T),(\mathcal Sh)(T)\rangle}_H \\
        &={\langle\Theta\mathcal Su,\mathcal Sh\rangle}_{W_0(0,T)',W_0(0,T)} -{\langle\Xi(y_1^\mathsf d-\hat y,y_2^\mathsf d-\hat y(T)),\mathcal Sh\rangle}_\Ha
    \end{align*}
    Note that in the first summand, we understand $\mathcal S$ as a mapping from $\U$ to $W_0(0,T)$ whereas in the second, it is considered to map to $\Ha$. In both cases, the dual operator $\mathcal S'$ maps to $\tilde U' \simeq \tilde U$ which directly gives \eqref{GVLuminy:Eq4.3.11}.
    \item [2)] Let $u,v\in\mathscr U$ be chosen arbitrarily. We set $y=\mathcal Su\in W_0(0,T)$ and $w=\mathcal Sv \in W_0(0,T)$. Recall that we identify $\U$ with its dual space $\U'$. From the integration by parts formula and Lemma~\ref{Lemma:HI-31}-2 we infer that
    \begin{align*}
        {\langle\mathcal S'\Theta\mathcal Sv,u\rangle}_\U&={\langle\Theta\mathcal Sv,\mathcal Su\rangle}_{W_0(0,T)',W_0(0,T)}={\langle\Theta w,y\rangle}_{W_0(0,T)',W_0(0,T)}\\
        &=\int_0^T \sigma_1\, {\langle w(t),y(t)\rangle}_H\,\mathrm dt+\sigma_2\,{\langle w(T),y(T)\rangle}_H\\
    &=-{\langle \mathcal Bu,p\rangle}_{L^2(0,T;V'),L^2(0,T;V)}=-{\langle u,\mathcal B'p\rangle}_\U=-{\langle \mathcal B'\mathcal Av,u\rangle}_\U.
    \end{align*}
    Since $u,v\in\U$ are chosen arbitrarily, we have $\mathcal B'\mathcal A=-\mathcal S'\Theta\mathcal S$. Further, we find utilizing \eqref{GVLuminy:Eq4.3.15}
    \begin{align*}
        &{\langle\mathcal S'\Xi(\ydQ-\hat y,\ydT-\hat y(T)),u\rangle}_\U={\langle\Xi(\ydQ-\hat y,\ydT-\hat y(T)),\mathcal Su\rangle}_{W_0(0,T)',W_0(0,T)}\\
        &\quad=\int_0^T \sigma_1\, {\langle\ydQ-\hat y(t),y(t)\rangle}_H\,\mathrm dt+\sigma_2\,{\langle\ydT-\hat y(T),y(T)\rangle}_H\\
        &\quad=\int_0^T -{\langle \hat p_t(t),y(t)\rangle}_H+a(t;y(t),\hat p(t))\,\mathrm dt+{\langle \hat p(T),y(T)\rangle}_H\\
        &\quad=\int_0^T {\langle y_t(t),\hat p(t)\rangle}_H+a(t;y(t),\hat p(t))\,\mathrm dt=\int_0^T{\langle(\mathcal Bu)(t),\hat p(t)\rangle}_{V',V}\,\mathrm dt={\langle\mathcal B'\hat p,u\rangle}_\U,
    \end{align*}
    which gives the claim since $u \in \U$ was chosen arbitrarily.\hfill$\Box$
\end{enumerate}

\noindent{\bf\em Proof of Lemma~{\em\ref{Lemma:Jprime}}.} Let $x=(y,u)\in\X$ and $x^\delta=(y^\delta,u^\delta)\in\X$ be given arbitrarily. Then
\begin{align*}
    &J(x+x^\delta)-J(x)-J'(x)x^\delta\\
    &=\frac{\sigma_1}{2}\int_0^T{\|y(t)+y^\delta(t)-\ydQ(t)\|}_H^2\,\mathrm dt+\frac{\sigma_2}{2}\,{\|y(T)+y^\delta(T)-\ydT\|}_H^2+\frac{\sigma}{2}\,{\|u+u^\delta-\un\|}_\U^2\\
    &\quad-\frac{\sigma_1}{2}\int_0^T{\|y(t)-\ydQ(t)\|}_H^2\,\mathrm dt-\frac{\sigma_2}{2}\,{\|y(T)-\ydT\|}_H^2-\frac{\sigma}{2}\,{\|u-\un\|}_\U^2\\
    &\quad-\sigma_1\int_0^T{\langle y(t)-\ydQ(t),y^\delta(t)\rangle}_H\,\mathrm dt-\sigma_2\,{\langle y(T)-\ydT,y^\delta(T)\rangle}_H-\sigma\,{\langle u-\un,u^\delta\rangle}_\U\\
    &=\frac{\sigma_1}{2}\,{\|y^\delta\|}_{L^2(0,T;H)}^2+\frac{\sigma_2}{2}\,{\|y^\delta(T)\|}_H^2+\frac{\sigma}{2}\,{\|u^\delta\|}_\U^2.
\end{align*}
Since $\Y$ is continuously embedded into $L^2(0,T;H)$ and $C([0,T];H)$, respectively, there exists a constant $c_1>0$ independent of $y^\delta$, satisfying
\begin{align*}
    \max\big({\|y^\delta\|}_{L^2(0,T;H)},{\|y^\delta(T)\|}_H\big)\le c_1\,{\|y^\delta\|}_\Y.
\end{align*}
Therefore, we have
\begin{align*}
    \big|J(x+x^\delta)-J(x)-J'(x)x^\delta\big|\le c_2\big({\|y^\delta\|}_\Y^2+{\|u^\delta\|}_\U^2\big)=c_2\,{\|x^\delta\|}_\X^2
\end{align*}
with $c_2=\max(c_1^2\sigma_1,c_1^2\sigma_2,\sigma)/2$. Consequently,
\begin{align*}
    \frac{1}{\|x^\delta\|_\X}\,\big|J(x+x^\delta)-J(x)-J'(x)x^\delta\big|\le c_2 \|x^\delta\|_\X \to 0 \quad \text{as } x^\delta \to 0
\end{align*}
which implies that $J$ is Fr\'echet-differentiable in $x$ and its Fr\'echet-derivative is given by \eqref{Custer-1}. For the continuity of the mapping $\X\ni x\mapsto J'(x)\in\X'$, consider that an alternative way to represent the derivative $J'(x) \in \mathscr X'$ is
\begin{align*}
    J'(x)x^\delta=\sigma_1\,{\langle y-y_1^\mathsf d,y^\delta\rangle}_{L^2(0,T;H)} + \sigma_2\,{\langle y(T)-y_2^\mathsf d, y^\delta_2(T)\rangle}_H+\sigma\,{\langle u- \un, u^\delta\rangle}_\U
\end{align*}
From this representation, the continuity can be easily shown if we consider that the embeddings $\Y\hookrightarrow L^2(0,T;H)$, $\Y\hookrightarrow C([0,T];H)$ as well as the inner products are continuous.\hfill$\Box$

\medskip\noindent{\bf\em Proof of Lemma~{\em\ref{Lemma:EGprime}}.} The assertion follows from the fact that $\X\ni x\mapsto e(x)\in\Lambda'$ is a bounded, affine-linear operator. More precisely, using \eqref{ST-1} and \eqref{ST-2} we infer that
\begin{align*}
    {\langle e_1(x+x^\delta)-e_1(x)-e_1'(x)x^\delta,p^1\rangle}_{L^2(0,T;V'),L^2(0,T;V)}=0
\end{align*}
and
\begin{align*}
    {\langle e_2(x+x^\delta)-e_2(x)-e_2'(x)x^\delta,p^2\rangle}_H=0
\end{align*}
for all $x=(y,u)\in\X$, $x^\delta=(y^\delta,u^\delta)\in X$ and $p=(p^1,p^2)\in \Lambda$. Thus,
\begin{align*}
    {\|e(x+x^\delta)-e(x)-e'(x)x^\delta\|}_{\Lambda'}&=\sup_{\|p\|_\Lambda=1}\Big({\langle e_1(x+x^\delta)-e_1(x)-e_1'(x)x^\delta,p^1\rangle}_{L^2(0,T;V'),L^2(0,T;V)}\\
    &\hspace{18mm}+{\langle e_2(x+x^\delta)-e_2(x)-e_2'(x)x^\delta,p^2\rangle}_H\Big)=0.
\end{align*}
Consequently,
\begin{align*}
\lim_{\|x^\delta\|_\X\to 0}\frac{\|e(x+x^\delta)-e(x)-e'(x)x^\delta\|_{\Lambda'}}{\|x^\delta\|_\X}=0
\end{align*}
holds. The continuity of the mapping $\X\ni x\mapsto e'(x)$ is trivial because $e'(x)$ does not depend on $x$.\hfill$\Box$

\medskip\noindent{\bf\em Proof of Proposition~{\em\ref{Prop:HI-30}}.}
\begin{enumerate}
    \item [$\Rightarrow$:] 
    Suppose that \eqref{Eq:HI-35} holds, i.e. for every $y^\delta \in \Y$, we have
    \begin{align*}
        0=\mathcal L_y(x(u),p(u))y^\delta = J_y(x(u))y^\delta + \langle e_y(x(u))y^\delta, p(u) \rangle_{\Lambda',\Lambda}={\langle J_y(x(u))+e_y(x(u))'p(u),y^\delta\rangle}_{\Y',\Y}.
    \end{align*}
    Then we infer that
    \begin{align*}
        e_y(x(u))'p(u)=-J_y(x(u))\quad\text{in }\Y'.
    \end{align*}
    Since $e_y(x)'$ is invertible by Remark \ref{Remark:HI-35}, we find that $p(u)$ satisfies \eqref{HI-7a}.
    \item [$\Leftarrow$:] Let the adjoint variable $p(u)=(p^1(u),p^2(u))$ be given by \eqref{HI-7a}. Then 
    \begin{align*}
        \mathcal L_y(x(u),p(u))y^\delta&={\langle J_y(x(u))-e_y(x(u))'(e_y(x(u))')^{-1}J_y(x(u)),y^\delta\rangle}_{\Y',\Y}\\
        &={\langle J_y(x(u))-J_y(x(u)),y^\delta\rangle}_{\Y',\Y}=0.
    \end{align*}
    which is \eqref{Eq:HI-35}.
\end{enumerate}
Using $\mathcal L_y(y,u,p)=0$ in $\Y'$, Lemmas~\ref{Lemma:Jprime} and \ref{Lemma:EGprime} we find that
\begin{equation}
    \label{ST-7}
    \begin{aligned}
        0&=\sigma_1\int_0^T{\langle y(t)-\ydQ(t),y^\delta(t)\rangle}_H\,\mathrm dt+\sigma_2\,{\langle y(\te)-\ydT,y^\delta(T)\rangle}_H\\
        &\quad+\int_0^T{\langle y^\delta_t(t),p^1(t)\rangle}_{V',V}+a(t;y^\delta(t),p^1(t))\,\mathrm dt+{\langle y^\delta(0), p^2\rangle}_H
    \end{aligned}
\end{equation}
for all $y^\delta\in\Y$. Suppose that $y^\delta(t)=\chi(t)\varphi$ with $\chi\in C^\infty_0(0,T)$ and $\varphi\in V$. Then $y^\delta\in\Y$ and $y^\delta(0)=y^\delta(T)=0$ hold. Moreover, we have
\begin{equation}
    \label{ST-8}
    \int_0^T{\langle y^\delta_t(t),p^1(t)\rangle}_{V',V}\,\mathrm dt=\bigg\langle\int_0^T p^1(t)\dot\chi(t)\,\mathrm dt,\varphi\bigg\rangle_H=-\bigg\langle\int_0^T p^1_t(t)\chi(t)\,\mathrm dt,\varphi\bigg\rangle_{V',V},
\end{equation}
where $p^1_t$ stands for the distributional derivative of $p^1$. Inserting \eqref{ST-8} into \eqref{ST-7} leads to
\begin{align*}
    \int_0^T\Big(-{\langle p_t^1(t),\varphi\rangle}_{V',V}+a(t;\varphi,p^1(t))+\sigma_1\,{\langle y(t)-\ydQ(t),\varphi\rangle}_H\Big)\chi(t)\,\mathrm dt=0
\end{align*}
for all $\varphi\in V$ and $\chi\in C^\infty_0(0,T)$. Since $C_0^\infty(0,T)$ is dense in $L^2(0,T)$ we get
\begin{align*}
    -{\langle p^1_t(t),\varphi\rangle}_{V',V}+a(t;\varphi,p^1(t))+\sigma_1\,{\langle y(t)-\ydQ(t),\varphi\rangle}_H=0\quad\text{for all }\varphi\in V\text{ a.e. in }[0,T).
\end{align*}
Note that for almost all $t\in[0,T]$ we have
\begin{align*}
    a(t;\cdot\,,p^1(t))+\sigma_1\,{\langle y^1(t)-\ydQ(t),\cdot\rangle}_H\in V'.
\end{align*}
Hence, $p^1_t\in L^2(0,T;V')$, which implies $p^1\in\Y$ and
\begin{align*}
    {\langle p^1_t(t),\varphi\rangle}_{V',V}=\frac{\mathrm d}{\mathrm dt}\,{\langle p^1(t),\varphi\rangle}_H\quad\forall\varphi\in V\text{ a.e. in }[0,T].
\end{align*}
Thus, we have derived
\begin{align*}
    -\frac{\mathrm d}{\mathrm dt}\,{\langle p^1(t),\varphi\rangle}_H+a(t;\varphi,p^1(t))=\sigma_1\,{\langle\ydQ(t)-y(t),\varphi\rangle}_H\quad\text{for all }\varphi\in V\text{ a.e. in }[0,T)
\end{align*}
which coincides with \eqref{GVLuminy:Eq4.3.16a}. Combining this and \eqref{ST-7} yields:
\begin{align*}
    0&=\sigma_1\int_0^T{\langle y(t)-\ydQ(t),y^\delta(t)\rangle}_H\,\mathrm dt+\sigma_2\,{\langle y(T)-\ydT,y^\delta(T)\rangle}_H\\
    &\quad- \sigma_1\int_0^T \langle y_1^\mathsf d(t)-y(t),y^\delta(t) \rangle_H\,\mathrm dt+{\langle y^\delta(T),p^1(T)\rangle}_H - \langle y^\delta(0),p^1(0)\rangle_H+{\langle y^\delta(0), p^2\rangle}_H \\
    &= {\langle p^1(T)-\sigma_2(\ydT-y(T)),y^\delta(T)\rangle}_H-{\langle p^1(0)-p^2,y^\delta(0)\rangle}_H
\end{align*}
where we have used partial integration for the second interval in \eqref{ST-7}. Choosing $y^\delta$ with $y^\delta(0)=0$ implies $p^1(T)=\sigma_2(y_2^\mathsf d-y(T))$ in $H$ which is \eqref{GVLuminy:Eq4.3.16b}. Furthermore, we get $p^2=p^1(0)$ in $H$ if we choose $y^\delta$ with $y^\delta(0)=0$.\hfill$\Box$

\medskip\noindent{\bf\em Proof of Lemma~{\em\ref{Lemma:HI-106}}.} Let $x=(y,u)\in\X$ and $x^\delta=(y^\delta,u^\delta)\in\X$ and $\tilde x^\delta=(\tilde y^\delta,\tilde u^\delta)\in\X$. Then we infer from \eqref{Custer-1} and \eqref{HI-106} that
\begin{align*}
    {\|J'(x+x^\delta)-J'(x)-J''(x)x^\delta\|}_{\X'}=\sup_{\|\tilde x^\delta\|_\X=1}{\langle J'(x+x^\delta)-J'(x)-J''(x)x^\delta,\tilde x^\delta\rangle}_{\X',\X}=0
\end{align*}
so that the second Fr\'echet derivative is given by \eqref{HI-106}. Since $J''(x)$ is independent of $x$, we conclude that $x\mapsto J''(x)$ is continuous, of course. Now we turn to the mapping $e$. It follows from \eqref{ST-2} that for every $x\in\X$ the linear and bounded operator $e'(x)$ is independent of $x$. Thus, $e$ is twice continuously Fr\'echet-differentiable and its second derivative is zero.\hfill$\Box$

\medskip\noindent{\bf\em Proof of Lemma~{\em\ref{Lemma:HI-200}}.} Let $x^\delta=(y^\delta,u^\delta)\in\mathrm{ker}\,e'(x)$ be given. Then we have
\begin{equation}
    \label{HI-109}
    e'(x)x^\delta=0\quad\text{in }\mP'.
\end{equation}
Due to Lemma~\ref{Lemma:EGprime} we deduce that \eqref{HI-109} is equivalent with
\begin{equation}
    \label{HI-110}
    \begin{aligned}
        \int_0^T{\langle y^\delta_t(t),\varphi(t)\rangle}_{V',V}+a(t;y^\delta(t),\varphi(t))\,\mathrm dt&= \int_0^T{\langle (\mathcal{B}u^{\delta})(t),\varphi(t)\rangle}_{V',V}\,\mathrm dt\quad\forall \varphi\in L^2(0,T;V),\\
        y^\delta(\tc)&=0.
    \end{aligned}
\end{equation}
System \eqref{HI-110} is equivalent with
\begin{equation}
    \label{HI-110a}
    \begin{aligned}
        \frac{\mathrm d}{\mathrm dt}\,{\langle y^\delta(t),\varphi\rangle}_H+a(t;y^\delta(t),\varphi) &={\langle (\mathcal{B}u^{\delta})(t),\varphi\rangle}_{V',V}\quad \text{for all } \varphi\in V\text{ a.e. in }(0,T],\\
        y^\delta(\tc)&=0.
    \end{aligned}
\end{equation}
Applying Corrollary~\ref{Corollary:HI-20}-2) we find that $y^\delta$ is uniquely determined as $y^\delta=\mathcal Su^\delta$. The estimate follows directly from \eqref{HI-100}.\hfill$\Box$

\medskip\noindent{\em\bf Proof of Proposition~{\em\ref{Proposition:HI-110}}.} Applying \eqref{HI-200} and Lemma~\ref{Lemma:HI-200} we find
\begin{align*}
    {\langle\mathcal L_{xx}(x,p)x^\delta,x^\delta\rangle}_{\X',\X}\ge\sigma\,{\|u^\delta\|}_\U^2\ge\frac{\sigma}{2}\,{\|u^\delta\|}_\U^2+\frac{\sigma}{2c_\mathsf{ker}}\,{\|y^\delta\|}_\Y^2\ge\kappa\,{\|x^\delta\|}_\X^2
\end{align*}
for all $x^\delta\in\mathrm{ker}\,e'(x)$ with $\kappa=\tfrac{\sigma}{2}\min\{1,1/c_\mathsf{ker}\}>0$.\hfill\endproof

\subsection{Proofs of Section~\ref{SIAM-Book:Section4.2}}
\label{SIAM-Book:Section4.7.2}

\noindent{\bf\em Proof of Corollary~{\em\ref{Corollary:HI-40}}.} For part 1), it suffices to choose $u = 0$ in Proposition \ref{Prop:DualPODApriori}, resulting in $p^\ell = \hat p^\ell$. For part 2), assume that $f=0$ and $y_\circ=0$ holds in \eqref{GVLuminy:Eq4.1.1}, and therefore $\hat y=0$. Furthermore, let $y^\mathsf d_2=0$ hold in $H$ and $y_1^\mathsf d=0$ in $L^2(0,T;H)$. Altogether, we then have $p^\ell = \mathcal A^\ell u$. Applying Proposition \ref{Prop:DualPODApriori} yields \eqref{eq:AellContinuity1}. Since $\mathcal S^\ell: \U \to \Y$ is continuous and $\Y$ embeds continuously into $C([0,T];H)$ as well as $L^2(0,T;H)$, there also exists $C_1 > 0$ with 
\begin{align*}
    \max \left ( \lVert (\mathcal S^\ell u)(T) \rVert_H, \lVert \mathcal S^\ell u \rVert_{L^2(0,T;H)} \right) \le C_1 \|u\|_\U
\end{align*}
A further estimation of \eqref{eq:AellContinuity1} then results in 
\begin{align*}
    \| \mathcal A^\ell u \|_\Y\le C C_1 (\sigma_2 + \sigma_1) \| u \|_\U
\end{align*}
which proves the boundedness of $\mathcal A^\ell$.\hfill$\Box$

\medskip\noindent{\bf\em Proof of Corollary~{\em\ref{Cor:pell}}.} Let $y^\ell=\hat y^\ell+\mathcal S^\ell u$. Due to \eqref{Wemb} we have 
\begin{equation}
    \label{Ida-2}
    {\|y^\ell(T)\|}_H\le {\|y^\ell\|}_{C([0,T];H)}\le c_\Y\,{\|y^\ell\|}_\Y.
\end{equation}
Moreover, we infer from the embedding inequality \eqref{Poincare} that
\begin{equation}
    \label{Ida-3}
    {\|y^\ell\|}_{L^2(0,T;H)}^2=\int_0^T{\|y^\ell(t)\|}_H^2\,\mathrm dt\le \int_0^T c_V^2\,{\|y^\ell(t)\|}_V^2\,\mathrm dt=c_V^2\,{\|y^\ell\|}_{L^2(0,T;V)}^2.
\end{equation}
Thus, using Proposition~\ref{Prop:DualPODApriori}, triangle inequality, \eqref{Ida-1}, \eqref{Ida-2} and \eqref{Ida-3} we deduce
\begin{align*}
    {\|p^\ell\|}_\Y&\le C\left(\sigma_2\big({\|\ydT\|}_H+{\|y^\ell(T)\|}_H\big)+\sigma_1\big({\|\ydQ\|}_{L^2(0,T;H)}+{\|y^\ell\|}_{L^2(0,T;H)}\big)\right)\\
    &\le c_1\big({\|\ydT\|}_H+{\|\ydQ\|}_{L^2(0,T;H)}\big)+c_2\,{\|y^\ell\|}_\Y\\
    &\le c_3\left({\|\ydT\|}_H+{\|\ydQ\|}_{L^2(0,T;H)}+{\|y_\circ\|}_H+{\|f\|}_{L^2(0,T;V')}+{\|u\|}_\U\right)
\end{align*}
with $c_1=C\max(\sigma_2,\sigma_1)$, $c_2=C(\sigma_2 C_\Y+\sigma_1c_V)$ and $c_3=\max(c_1,c_2C_y)$.\hfill$\Box$

\medskip\noindent{\bf\em Proof of Lemma~{\em\ref{Lemma:HI-50}}.} From $p,p^\ell\in H^1(0,T;V)$ we infer that $\theta^\ell\in H^1(0,T;V)$ and
\begin{equation}
    \label{Dual:APriori-Est-1b}
    {\langle \theta^\ell_t(t),\psi\rangle}_{V',V}={\langle \theta^\ell_t(t),\psi\rangle}_H\quad \text{for all }\psi\in V\text{ in }[0,T]\text{ a.e.}
\end{equation}
Moreover, we have (cf. \eqref{Eq:SnapOhneDQ-2})
\begin{equation}
    \label{Eq:SnapOhneDQ-20}
    {\langle \mathcal P^\ell p_t(t),\psi\rangle}_H={\langle p_t(t),\psi\rangle}_H\quad\text{for all }\psi\in X^\ell\text{ in }[0,T]\text{ a.e.}
\end{equation}
From \eqref{GVLuminy:Eq4.3.16}, \eqref{SIAM:DualPOD}, \eqref{Eq:SnapOhneDQ-20}, \eqref{Dual:APriori-Est-1b} and \eqref{HI-101-1} we obtain that
\begin{align*}
    -{\langle\theta^\ell_t(t),\psi\rangle}_H+a(t;\psi,\theta^\ell(t))
    &=-{\langle p_t(t),\psi\rangle}_H+a(t;\psi,\mathcal P^\ell p(t))+{\langle p_t^\ell(t),\psi\rangle}_H-a(t;\psi,p^\ell(t))\\
    &=a(t;\psi,(\mathcal P^\ell p-p)(t))+\sigma_1\,{\langle (y^\ell-y)(t),\psi\rangle}_H\\
    &\le\gamma\,{\|\psi\|}_V{\|(\mathcal P^\ell p-p)(t)\|}_V+\sigma_1\,{\|(y^\ell-y)(t)\|}_H{\|\psi\|}_H
\end{align*}
in $[0,T]$ a.e. Choosing $\psi=\theta^\ell(t)\in X^\ell$ and using \eqref{HI-101-2} we find
\begin{align*}
    -\frac{1}{2}\frac{\mathrm d}{\mathrm dt}\,{\|\theta^\ell(t)\|}_H^2+\gamma_1\,{\|\theta^\ell(t)\|}_V^2\le\gamma\,{\|\theta^\ell(t)\|}_V{\|(\mathcal P^\ell p-p)(t)\|}_V+\sigma_1\,{\|(y^\ell-y)(t)\|}_H{\|\theta^\ell(t)\|}_H
\end{align*}
in $[0,T]$ a.e. We apply Young's inequality (cf. Lemma~\ref{lem:youngsInequality} with $p=q=2$ and $\varepsilon=\gamma_1^2/\gamma^2$ for the first, $\varepsilon=\sigma_1$ for the second product):
\begin{equation}
    \label{Albi-1}
    -\frac{\mathrm d}{\mathrm dt}\,{\|\theta^\ell(t)\|}_H^2+\gamma_1\,{\|\theta^\ell(t)\|}_V^2\le{\|\theta^\ell(t)\|}_H^2+\frac{\gamma^2}{\gamma_1^2}\,{\|(\mathcal P^\ell p-p)(t)\|}_V^2+\sigma_1^2\,{\|(y^\ell-y)(t)\|}_H^2
\end{equation}
in $[0,T]$ a.e. Recall that the embedding constant $c_V$ has been introduced in \eqref{Poincare}. We set
\begin{align*}
    c_1=\left\{
    \begin{array}{ll}
        0&\text{if }1-\frac{\gamma_1}{c_V^2}\le0,\\[2ex]
        1-\frac{\gamma_1}{c_V^2}&\text{otherwise.}
    \end{array}
    \right.
\end{align*}
and infer
\begin{align*}
    -\frac{\mathrm d}{\mathrm dt}\,{\|\theta^\ell(t)\|}_H^2\le c_1\,{\|\theta^\ell(t)\|}_H^2+\frac{\gamma^2}{\gamma_1^2}\,{\| (\mathcal P^\ell p-p)(t)\|}_V^2 + \sigma_1^2\,{\|(y^\ell-y)(t)\|}_H^2.
\end{align*}
Using the Gronwall Lemma (cf. Corollary \ref{DualGronwall}), it follows that
\begin{align*}
    {\|\theta^\ell(t)\|}_H^2\le e^{c_1(T-t)}\bigg({\|\theta^\ell(T)\|}_H^2 + \int_t^T \frac{\gamma^2}{\gamma_1^2}\,{\|(\mathcal P^\ell p-p)(s)\|}_V^2 + \sigma_1^2\,{\|(y^\ell-y)(s)\|}_H^2 ~\mathrm ds \bigg).
\end{align*}
For the terminal value, we use the definitions of $p$ and $p^\ell$ in \eqref{GVLuminy:Eq4.3.16b} and \eqref{SIAM:DualPODb}:
\begin{equation}
    \label{Albi-2}
    \theta^\ell(T)=\mathcal P^\ell p(T)-p^\ell(T)=\sigma_2\big(y^\ell(T)-\mathcal P^\ell y(T)\big).
\end{equation}
Thus we arrive at
\begin{equation}
    \label{eq:thetaEstimate1}
    \begin{aligned}
        {\|\theta^\ell(t)\|}_H^2&\le e^{c_1(T-t)} \Big( \sigma_2^2 \|(y^\ell-\mathcal P^\ell y)(T)\|_H^2 + \frac{\gamma^2}{\gamma_1^2}\,{\| \mathcal P^\ell p - p \|}_{L^2(t,T;V)}^2+\sigma_1^2\,{\| y^\ell-y \|}_{L^2(t,T;H)}^2 \Big)\\
        &\le c_3\Big({\|(y^\ell-\mathcal P^\ell y)(T)\|}_H^2+{\|\mathcal P^\ell p-p\|}_{L^2(0,T;V)}^2+{\|y^\ell-y\|}_{L^2(0,T;H)}^2\Big) 
    \end{aligned}
\end{equation}
with $c_3=e^{c_1T}\max(\sigma_2^2,\gamma^2/\gamma_1,\sigma_1^2)>0$. Next, we integrate \eqref{Albi-1} over $[0,T]$ from which it follows:
\begin{align*}
    \gamma_1\,{\|\theta^\ell\|}_{L^2(0,T;V)}^2 &\le {\|\theta^\ell(T)\|}_H^2 + \int_0^T{\|\theta^\ell(t)\|}_H^2\,\mathrm dt+ \frac{\gamma^2}{\gamma_1^2} \|\mathcal P^\ell-p\|_{L^2(0,T;V)}^2+\sigma_1^2\,{\|y^\ell-y\|}_{L^2(0,T;H)}^2.
 \end{align*}
Using \eqref{Albi-2} for $\theta^\ell(T)$ and applying \eqref{eq:thetaEstimate1}, we find that
\begin{align*}
    {\|\theta^\ell\|}_{L^2(0,T;V)}^2\le c_4\Big({\|(y^\ell-\mathcal P^\ell y)(T)\|}^2_H+{\|\mathcal P^\ell p-p\|}_{L^2(0,T;V)}^2+{\|y^\ell-y\|}_{L^2(0,T;H)}^2\Big)
\end{align*}
with $c_4=\nicefrac{1}{\gamma}\max(\sigma_2^2+c_3T,\gamma^2/\gamma_1^2+c_3T,\sigma_1^2+c_3T)$. Using \eqref{Ida-3} and \eqref{Wemb} we get \eqref{Albi-4a} with $c_\theta=c_4\max(c_\Y^2,1,c_V^2)$. \hfill$\Box$

\medskip\noindent{\bf\em Proof of Theorem~{\em\ref{Th:DualA-PrioriError-2}}.} Combining \eqref{Dual:APriori-Est-1} and \eqref{Albi-4a} we have
\begin{equation}
    \label{Albi-4}
    \begin{aligned}
        \int_0^T{\|p(t)-p^\ell(t)\|}_V^2\,\mathrm dt&\le C\Big({\|\mathcal P^\ell p-p\|}_{L^2(0,T;V)}^2+{\|\mathcal P^\ell y-y\|}_{L^2(0,T;V)}^2\Big)\\
        &\le C\sum_{k=1}^\wp\omega_k^\wp\int_0^T\big\|y^k(t)-\mathcal P^\ell y^k(t)\big\|_V^2\,\mathrm dt
    \end{aligned}
\end{equation}
with $C=2c_\theta+2$. It follows from \eqref{RatePH-V} that for $X=H$ and $\mathcal P^\ell=\mathcal P^\ell_H$ we have
\begin{align*}
    \sum_{k=1}^\wp\omega_k^\wp\int_0^T\big\|y^k(t)-\mathcal P^\ell y^k(t)\big\|_V^2\,\mathrm dt=\sum_{i>\ell}\lambda_i^H\big\|\psi_i^H\big\|_V^2<\infty.
\end{align*}
Thus, \eqref{Albi-4} yields
\begin{align*}
    \int_0^T{\|p(t)-p^\ell(t)\|}_V^2\,\mathrm dt\le C\sum_{i>\ell}\lambda_i^H\big\|\psi_i^H\big\|_V^2. 
\end{align*}
For the choices $X=V$ and $\mathcal P^\ell=\mathcal Q^\ell_H$ the claim follows analogously from \eqref{Albi-4} and \eqref{RateQH-V}.\hfill$\Box$

\medskip\noindent{\bf\em Proof of Theorem~{\em\ref{GVLuminy:Theorem4.4.2}}.} Choosing $u=\bar u^\ell$ in \eqref{GVLuminy:Eq4.3.20-2} and $u=\bar u$ in \eqref{GVLuminy:Eq4.4.13} we get the variational inequality
\begin{equation}
    \label{GVLuminy:Eq4.4.16}
    0 \le {\langle \sigma (\bar u-\bar u^\ell)-\mathcal B' (\bar p-\bar p^\ell),\bar u^\ell-\bar u \rangle}_\U.
\end{equation}
Utilizing Lemma~\ref{Lemma:HI-30} and $\langle\Theta \varphi,\varphi\rangle_{W_0(0,T)',W_0(0,T)}\ge 0$ for all $\varphi\in W_0(0,T)$ we infer from \eqref{GVLuminy:Eq4.4.16} that
\begin{align*}
    0 \le & {\langle \mathcal B' \mathcal A^\ell\bar u^\ell-\mathcal B'\mathcal A\bar u+\mathcal B'(\hat p^\ell-\hat p),\bar u^\ell-\bar u \rangle}_\U-\sigma\,{\|\bar u-\bar u^\ell\|}_\U^2\\
    = & {\langle \mathcal B' \mathcal A^\ell(\bar u^\ell-\bar u)+\mathcal B' (\mathcal A^\ell-\mathcal A)\bar u+\mathcal B'(\hat p^\ell-\hat p),\bar u^\ell-\bar u \rangle}_\U-\sigma\,{\|\bar u-\bar u^\ell\|}_\U^2\\
    = & {\langle \mathcal B' \mathcal A^\ell(\bar u^\ell-\bar u)+\mathcal B' (\hat p^\ell+\mathcal A^\ell\bar u-\hat p-\mathcal A\bar u),\bar u^\ell-\bar u \rangle}_\U-\sigma\,{\|\bar u-\bar u^\ell\|}_\U^2\\
    \le & {\langle \Theta\mathcal S^\ell(\bar u-\bar u^\ell),\mathcal S^\ell(\bar u^\ell-\bar u)\rangle}_{W_0(0,T)',W_0(0,T)} +{\| \mathcal B' (\bar p^\ell(\bar u)-\bar p(\bar u))\|}_\U{\|\bar u^\ell-\bar u\|}_\U-\sigma\,{\|\bar u-\bar u^\ell\|}_\U^2 \\
    \le & {\|\mathcal B'\|}_{\mathscr L(L^2(0,T;V),\U)}{\| \bar p^\ell(\bar u)-\bar p(\bar u)\|}_{L^2(0,T;V)}{\|\bar u^\ell-\bar u\|}_\U-\sigma\,{\|\bar u-\bar u^\ell\|}_\U^2,
\end{align*}
where we have set $\bar p^\ell(\bar u)=\hat p^\ell+\mathcal A^\ell\bar u$ and $\bar p(\bar u)=\hat p+\mathcal A\bar u$. Consequently,
\begin{align*}
    {\|\bar u-\bar u^\ell\|}_\U\le\frac{1}{\sigma}\,{\|\mathcal B'\|}_{\mathscr L(L^2(0,T;V),\U)}{\| \bar p^\ell(\bar u)-\bar p(\bar u)\|}_{L^2(0,T;V)}.
\end{align*}
Now \eqref{GVLuminy:Eq4.4.15-2} follows from Theorem~\ref{Th:DualA-PrioriError-2}.\hfill$\Box$

\noindent{\bf\em Proof of Theorem~{\em\ref{GVLuminy:Theorem4.5.1}}.} Estimate \eqref{GVLuminy:Eq4.5.7a} has already be shown. We proceed by constructing the function $\zeta^\ell$. Here we adapt the lines of the proof of Proposition~3.2 in \cite{TV09} to our optimal control problem. Suppose that we know $\bar u^\ell$ and $\tilde p^\ell=\hat p+\mathcal A\bar u^\ell$. The goal is to determine $\zeta^\ell \in\U$ satisfying \eqref{GVLuminy:Eq4.5.1}. We distinguish three different cases.
\begin{itemize}
    \item  Case $\bar u^\ell_i(s)=u_{\mathsf ai}(s)$ for fixed $s\in\mathscr D$ and $i\in\{1,\ldots,\mathsf m\}$: Then $u_i(s)-\bar u_i^\ell(s)=u_i(s)-u_{\mathsf ai}(s) \ge 0$ for all $u \in\Uad$. Hence, $\zeta^\ell_i(s)$ has to satisfy
    \begin{equation}
        \label{GVLuminy:Eq4.5.8}
        \big(\sigma (\bar u^\ell-\un)-\mathcal B'\tilde p^\ell\big)_i(s)+\zeta_i^\ell(s)\ge 0.
    \end{equation}
    Setting $\zeta_i^\ell(s)= -\min(0,(\sigma (\bar u^\ell-\un)-\mathcal B'\tilde p^\ell)_i(s))$ the value $\zeta_i^\ell(s)$ satisfies \eqref{GVLuminy:Eq4.5.8}. 
    \item  Case $\bar u^\ell_i(s)=u_{\mathsf bi}(s)$ for fixed $s\in\mathscr D$ and $i\in \{1,\ldots,\mathsf m\}$: Now, $u_i(s)-\bar u_i^\ell(s)=u(s)-u_{\mathsf bi}(s) \le 0$ for all $u \in\Uad$. Analogously to the first case we define $\zeta_i^\ell(s)=-\max(0,(\sigma (\bar u^\ell-\un)-\mathcal B'\tilde p^\ell)_i(s))$ to ensure 
    \begin{equation}
        \label{GVLuminy:Eq4.5.8b}
        \big(\sigma (\bar u^\ell-\un)-\mathcal B'\tilde p^\ell\big)_i(s)+\zeta_i^\ell(s)\le 0.
    \end{equation}
    \item Case $u_{\mathsf ai}(s)<\bar u_i^\ell(s)<u_{\mathsf bi}(s)$ for fixed $s\in\mathscr D$ and $i\in\{1,\ldots,\mathsf m\}$: Consequently, %
    \begin{equation}
        \label{GVLuminy:Eq4.5.8c}
        (\sigma (\bar u^\ell-\un)-\mathcal B'\tilde p^\ell)_i(s)+\zeta^\ell_i(s)=0
    \end{equation} 
    has to hold, which is guaranteed by setting $\zeta_i^\ell(s)=-(\sigma (\bar u^\ell-\un)-\mathcal B'\tilde p^\ell)_i(s)$.
\end{itemize}
It remains to prove that $\zeta^\ell$ tends to zero for $\ell\to\infty$. Here we follow the proof of Theorem~4.11 in \cite{TV09}. By Theorem~\ref{GVLuminy:Theorem4.4.2}-1), the sequence $\{\bar u^\ell\}_{\ell\in\mathbb N}$ converges to $\bar u$ in $U$. Since the linear operator $\mathcal B'\mathcal A$ is bounded and $\tilde p^\ell=\hat p+\mathcal A\bar u^\ell$ holds, $\{\mathcal B' \tilde p^\ell\}_{\ell\in\mathbb N}$ tends to $\mathcal B'\bar p=\hat p+\mathcal B'\mathcal A\bar u$ as well. Furthermore, for $\bar p^\ell = \hat{p} + \mathcal{A}^{\ell} \bar{u}^{\ell}$ it holds 
\begin{align*}
    \mathcal B'\bar p^\ell - \mathcal B'\bar{p} = \mathcal{B}'(\mathcal{A}^\ell \bar{u}^\ell - \mathcal{A}\bar{u}) = \mathcal{B}'\mathcal{A}^\ell (\bar{u}^\ell - \bar{u}) + \mathcal{B}'(\mathcal{A}^\ell - \mathcal{A}) \bar{u}.
\end{align*}
Since the operator $\mathcal A^\ell$ is bounded independently of $\ell$, the sequence $\{\bar u^\ell\}_{\ell\in\mathbb N}$ converges to $\bar u$ and $\lim_{\ell \to \infty} (\mathcal A^\ell - \mathcal A) \bar{u} = 0$ (due to \eqref{ConvDualVar}), we can conclude that the sequence $\{ \mathcal B'\bar p^\ell \}_{\ell\in\mathbb N}$ converges to $\mathcal B'\bar p$. \\
Hence, there exist subsequences $\{\bar u^{\ell_k}\}_{k \in \mathbb N}$, $\{\mathcal B' \bar p^{\ell_k}\}_{k \in \mathbb N}$ and $\{\mathcal B' \tilde p^{\ell_k}\}_{k \in \mathbb N}$ satisfying
\begin{subequations}
    \begin{align}
        \lim_{k \to \infty}\bar u_i^{\ell_k}(s) & =\bar u_i(s) \label{GVLuminy:Eq4.55.1}\\
        \lim_{k \to \infty}(\mathcal B'\tilde p^{\ell_k})_i(s) & =(\mathcal B' \bar p)_i(s) \label{GVLuminy:Eq4.55.2} \\
        \lim_{k \to \infty}(\mathcal B'\bar p^{\ell_k})_i(s) & =(\mathcal B' \bar p)_i(s) \label{GVLuminy:Eq4.55.3}
    \end{align}
\end{subequations}
f.a.a. $s\in\mathscr D$ and for $1\le i\le\mathsf m$. Next we consider the active and inactive sets for $\bar u$.
\begin{itemize}
    \item Let $s\in\mathscr J_i=\{s\in\mathscr D\,|\, u_{\mathsf ai}(s)<\bar u_i(s)<u_{\mathsf bi}(s)\}$ for $i\in\{1,\ldots,\mathsf m\}$. By using statement 2) from Proposition~\ref{Pro:HI-1} this immediately implies that $(\sigma(\bar u-\un)-\mathcal B'\bar p)_i(s) = 0$. Additionally, from \eqref{GVLuminy:Eq4.55.1} we can deduce that there exists $k_\circ=k_\circ(s)\in\mathbb N$ such that $\bar u^{\ell_k}_i(s) \in (u_{\mathsf ai}(s),u_{\mathsf bi}(s))$ holds for all $k\ge k_\circ$. Thus, $\zeta^{\ell_k}_i(s) =  -(\sigma (\bar u^{\ell_k}-\un) - \mathcal B'\tilde p^{\ell_k})_i(s)$ for all $k \ge k_\circ$. Now \eqref{GVLuminy:Eq4.55.1} and \eqref{GVLuminy:Eq4.55.2} imply
    \begin{equation}
        \label{GVLuminy:Eq4.5.10}
        \lim_{k \to \infty} \zeta^{\ell_k}_i(s) = (\sigma(\bar u-\un)-\mathcal B'\bar p)_i(s) = 0.
    \end{equation}
    \item Suppose that $s\in\mathscr A_{\mathsf ai}=\{s\in\mathscr D\,|\,u_{\mathsf ai}(s)=\bar u_i(s)\}$ for $i\in\{1,\ldots,\mathsf m\}$. Then it holds $(\sigma (\bar u-\un)-\mathcal B' \bar p)_i(s)\ge 0$. If $(\sigma (\bar u-\un)-\mathcal B' \bar p)_i(s) = 0$ we can proceed similarly to the first case and conclude that $\lim_{k \to \infty} \zeta^{\ell_k}_i(s) = 0$. If $(\sigma (\bar u-\un)-\mathcal B' \bar p)_i(s) > 0$ holds, there is a $k_\circ=k_\circ(s)\in\mathbb N$ such that $(\sigma (\bar u^{\ell_k}-\un) - \mathcal B'\bar p^{\ell_k})_i(s) > 0$ for all $k \geq k_\circ$ by \eqref{GVLuminy:Eq4.55.3}. Thus, Proposition~\ref{Pro:HI-1}-2) yields $\bar u^{\ell_k}_i(s)=u_{\mathsf ai}(s)$ for all $k \geq k_\circ$. As previously shown, $\zeta^{\ell_k}_i(s)$ is then given by $\zeta^{\ell_k}_i(s) = -\min(0,(\sigma (\bar u^\ell-\un)-\mathcal B'\tilde p^\ell)_i(s))$. Because of \eqref{GVLuminy:Eq4.55.2} we get
    \begin{equation}
        \label{GVLuminy:Eq4.5.10.1}
        \lim_{k \to \infty} \zeta^{\ell_k}_i(s) = -\min(0,(\sigma (\bar u-\un)-\mathcal B'\bar p)_i(s)) = 0.
    \end{equation} 
    \item Suppose that $s\in\mathscr A_{\mathsf bi}=\{s\in\mathscr D\,|\,u_{\mathsf bi}(s)=\bar u_i(s)\}$. Analogously to the second case we find
    \begin{equation}
        \label{GVLuminy:Eq4.5.10.2}
        \lim_{k \to \infty}\zeta^{\ell_k}_i(s)=0. 
    \end{equation}
\end{itemize}
Combining \eqref{GVLuminy:Eq4.5.10}-\eqref{GVLuminy:Eq4.5.10.2} we conclude that $\lim_{k \to \infty}\zeta^{\ell_k}_i = 0$ a.e. in $\mathscr D$ and for $1\le i\le\mathsf m$. Utilizing the dominated convergence theorem \cite[p.~24]{RS80} we have
\begin{align*}
    \lim_{k \to \infty} \big\| \zeta^{\ell_k}\big\|_\U=0.
\end{align*}
Since all subsequences contain a subsequence converging to zero, the claim follows from a standard argument.\hfill$\Box$

\subsection{Proofs of Section~\ref{SIAM-Book:Section4.3}}
\label{SIAM-Book:Section4.7.3}

\noindent{\bf\em Proof of Corollary~{\em\ref{Cor:LQR-1}}.} From \eqref{HI-220} we infer that
\begin{equation}
    \label{LQR:Dual-4}
    {\|y^{h\ell}\|}_\Y\le c_y\left({\|\mathcal P^{h\ell} y_\circ\|}_H+{\|f\|}_{L^2(0,T;V')}+{\|u\|}_\U\right)
\end{equation}
for a constant $c_y>0$ which is independent of $\ell$. Utilizing \eqref{Wemb} we have
\begin{equation}
    \label{LQR:Dual-5}
    {\|y^{h\ell}(T)\|}_H\le {\|y^{h\ell}\|}_{C([0,T];H)}\le c_\Y\,{\|y^{h\ell}\|}_\Y.
\end{equation}
Moreover, we infer from \eqref{Poincare} that
\begin{equation}
    \label{LQR:Dual-6}
    {\|y^{h\ell}\|}_{L^2(0,T;H)}\le c_V\,{\|y^{h\ell}\|}_{L^2(0,T;V)}\le c_V\,{\|y^{h\ell}\|}_\Y.
\end{equation}
Thus, using Proposition~\ref{Prop:LQR-1}, and \eqref{LQR:Dual-4}-\eqref{LQR:Dual-6} we deduce
\begin{align*}
    {\|p^{h\ell}\|}_\Y&\le C\left(\sigma_2\big({\|\ydT\|}_H+{\|y^{h\ell}(T)\|}_H\big)+\sigma_1\big({\|\ydQ\|}_{L^2(0,T;H)}+{\|y^{h\ell}\|}_{L^2(0,T;H)}\big)\right)\\
    &\le c_1\big({\|\ydT\|}_H+{\|\ydQ\|}_{L^2(0,T;H)}\big)+c_2\,{\|y^{h\ell}\|}_\Y\\
    &\le c_3\left({\|\ydT\|}_H+{\|\ydQ\|}_{L^2(0,T;H)}+{\|\mathcal P^{h\ell} y_\circ\|}_H+{\|f\|}_{L^2(0,T;V')}+{\|u\|}_\U\right)
\end{align*}
with $c_1=C\max(\sigma_2,\sigma_1)$, $c_2=C(\sigma_2 c_\Y+\sigma_1c_V)$ and $c_3=\max(c_1,c_2c_\Y)$.\hfill$\Box$

\medskip\noindent{\bf\em Proof of Proposition~{\em\ref{Prop1:GradCostApo}}.} Let the associated state variables are given by $y^h(u)=\hat y^h+\mathcal S^hu$ and $y^{h\ell}(u)=\hat y^{h\ell}+\mathcal S^{h\ell}u$, respectively. Then, we have proved in Proposition~\ref{Prop:NSAposti} the following a-posteriori error estimate for the difference $e^{h\ell}(t;u)=y^h(t;u)-y^{h\ell}(t;u)\in V^h$
\begin{equation}
    \label{GradCostApo-4}
    {\|e^{h\ell}(t;u)\|}_H^2\le \mathsf C_1^{h\ell}(t)\bigg(\mathsf R_\circ^{h\ell}+\int_0^t\mathsf R_1^{h\ell}(s)\,\mathrm ds\bigg),
\end{equation}
where the constants $\mathsf C_1^{h\ell}(t)$, $\mathsf R_\circ^{h\ell}$ and $\mathsf R_1^{h\ell}(s)$ are given by \eqref{Alb-4a} for almost all $s\in[0,T]$. Now we set $\mathsf e^{h\ell}(t;u)=p^h(t;u)-p^{h\ell}(t;u)$ for almost all $t\in[0,T]$. Then we derive from \eqref{MD-2-2} and \eqref{LQR:Dual-1} for all $\varphi^h\in V^h$ and $t\in[0.T]$ a.e.
\begin{align*}
    &-{\langle\mathsf e^{h\ell}(t;u),\varphi^h\rangle}_{V',V}+a(t;\varphi^h,\mathsf e^{h\ell}(t;u))\\
    &\quad=\sigma_1\,{\langle\ydQ(t)-y^h(t;u)\rangle}_H+{\langle p_t^{h\ell}(t;u),\varphi^h\rangle}_{V',V}-a(t;\varphi^h,p^{h\ell}(t;u))\\
    &\quad=-\sigma_1\,{\langle e^{h\ell}(t;u),\varphi^h\rangle}_H-{\langle\mathsf r^{h\ell}(t;u),\varphi^h\rangle}_{V',V}
\end{align*}
with the (dual) residual defined in \eqref{GradCostApo-2a}. Choosing $\varphi=\mathsf e^{h\ell}(t;u)\in V^h$ we derive from \eqref{HI-101-2} that
\begin{align*}
    -\frac{1}{2}\frac{\mathrm d}{\mathrm dt}\,{\|\mathsf e^{h\ell}(t;u)\|}_H^2+\gamma_1\,{\|\mathsf e^{h\ell}(t;u)\|}_V^2&\le-{\langle\mathsf e^{h\ell}(t;u),\varphi^h\rangle}_{V',V}+a(t;\varphi^h,\mathsf e^{h\ell}(t;u))\\
    &\le\sigma_1\,{\|e^{h\ell}(t;u)\|}_H{\|\mathsf e^{h\ell}(t;u)\|}_H+{\|\mathsf r^{h\ell}(t;u)\|}_{(V^h)'}{\|\mathsf e^{h\ell}(t;u)\|}_V.
\end{align*}
Next we apply the Young inequality and \eqref{Poincare}. We get
\begin{align*}
    &-\frac{1}{2}\frac{\mathrm d}{\mathrm dt}\,{\|\mathsf e^{h\ell}(t;u)\|}_H^2+\gamma_1\,{\|\mathsf e^{h\ell}(t;u)\|}_V^2\\
    &\quad\le\sigma_1c_V\,{\|e^{h\ell}(t;u)\|}_H{\|\mathsf e^{h\ell}(t;u)\|}_V+\frac{1}{\gamma_1}\,{\|\mathsf r^{h\ell}(t;u)\|}_{(V^h)'}^2+\frac{\gamma_1}{4}\,{\|\mathsf e^{h\ell}(t;u)\|}_V^2\\
    &\quad\le\frac{1}{\gamma_1}\big(\sigma_1^2c_V^2\,{\|e^{h\ell}(t;u)\|}_H^2+{\|\mathsf r^{h\ell}(t;u)\|}_{(V^h)'}^2\big)+\frac{\gamma_1}{2}\,{\|\mathsf e^{h\ell}(t;u)\|}_V^2
\end{align*}
which implies that
\begin{equation}
    \label{GradCostApo-5}
    -\frac{\mathrm d}{\mathrm dt}\,{\|\mathsf e^{h\ell}(t;u)\|}_H^2+\gamma_1\,{\|\mathsf e^{h\ell}(t;u)\|}_V^2\le\frac{2}{\gamma_1}\big(\sigma_1^2c_V^2\,{\|e^{h\ell}(t;u)\|}_H^2+{\|\mathsf r^{h\ell}(t;u)\|}_{(V^h)'}^2\big).
\end{equation}
After integration on the interval $[t,T]$ with $t\in [0,T)$ and using
\begin{equation}
    \label{GradCostApo-6}
    \begin{aligned}
    \mathsf e^{h\ell}(T;u)&=p^h(T;u)-p^{h\ell}(T;u)=\sigma_2\big(\mathcal P^h\ydT-y^h(T;u)\big)-\sigma_2\big(\mathcal P^{h\ell}\ydT-y^{h\ell}(T;u)\big)\\
    &=\sigma_2\big(\mathcal P^h-\mathcal P^{h\ell}\big)\ydT-\sigma_2e^{h\ell}(T;u)
    \end{aligned}
\end{equation}
(cf. \eqref{MD-2-2b} and \eqref{LQR:Dual-1b}) we find that
\begin{align*}
    {\|\mathsf e^{h\ell}(t;u)\|}_H^2&\le{\|\mathsf e^{h\ell}(T;u)\|}_H^2+\int_t^T\frac{2}{\gamma_1}\big(\sigma_1^2c_V^2\,{\|e^{h\ell}(s;u)\|}_H^2+{\|\mathsf r^{h\ell}(s;u)\|}_{(V^h)'}^2\big)\,\mathrm ds\\
    &\le2\sigma_2\,{\|(\mathcal P^h-\mathcal P^{h\ell})\ydT\|}_H^2+2\sigma_2\,{\|e^{h\ell}(T;u)\|}_H^2\\
    &\quad+\frac{2}{\gamma_1}\int_t^T\sigma_1^2c_V^2\,{\|e^{h\ell}(s;u)\|}_H^2+{\|\mathsf r^{h\ell}(s;u)\|}_{(V^h)'}^2\,\mathrm ds.
\end{align*}
We proceed by applying \eqref{GradCostApo-4}. Using the constants defined in \eqref{GradCostApo-2} we conclude that
\begin{align*}
    {\|\mathsf e^{h\ell}(t;u)\|}_H^2&\le\mathsf R_T^{h\ell}+\mathsf C_2^{h\ell}\bigg(\mathsf R_\circ^{h\ell}+\int_0^T\mathsf R_1^{h\ell}(s)\,\mathrm ds\bigg)\\
    &\quad+\int_t^T\mathsf C_3^{h\ell}(s)\bigg(\mathsf R_\circ^{h\ell}+\int_0^s\mathsf R_1^{h\ell}(\tau)\,\mathrm d\tau\bigg)+\frac{2}{\gamma_1}\,{\|\mathsf r^{h\ell}(s;u)\|}_{(V^h)'}^2\,\mathrm ds
\end{align*}
holds for almost all $t\in[0,T]$, where the constants $\mathsf R_\circ^{h\ell}$ and $\mathsf R_1^{h\ell}(s)$ are given by \eqref{Alb-4a}. Thus, we have shown estimate \eqref{GradCostApo-1}. Integrating \eqref{GradCostApo-5} over the interval $[0,T]$ and using \eqref{GradCostApo-4}, \eqref{GradCostApo-6} we find that
\begin{align*}
    \int_0^T{\|\mathsf e^{h\ell}(t;u)\|}_V^2\,\mathrm dt&\le\frac{1}{\gamma_1}\,{\|\mathsf e^{h\ell}(T;u)\|}_H^2+\frac{2}{\gamma_1^2}\int_0^T\sigma_1^2c_V^2\,{\|e^{h\ell}(t;u)\|}_H^2+{\|\mathsf r^{h\ell}(t;u)\|}_{(V^h)'}^2\,\mathrm dt\\
    &\le\frac{1}{\gamma_1}\,\bigg(\mathsf R_T^{h\ell}+\mathsf C_2^{h\ell}\Big(\mathsf R_\circ^{h\ell}+\int_0^T\mathsf R_1^{h\ell}(t)\,\mathrm dt\Big)\bigg)\\
    &+\frac{1}{\gamma_1}\int_0^T\mathsf C_3^{h\ell}(s)\bigg(\mathsf R_\circ^{h\ell}+\int_0^s\mathsf R_1^{h\ell}(\tau)\,\mathrm d\tau\bigg)+\frac{2}{\gamma_1}\,{\|\mathsf r^{h\ell}(t;u)\|}_{(V^h)'}^2\,\mathrm dt
\end{align*}
which gives estimate \eqref{GradCostApo-3}.\hfill$\Box$

\subsection{Proofs of Section~\ref{SIAM-Book:Section4.4}}
\label{SIAM-Book:Section4.7.4}

\noindent{\bf\em Proof of Theorem~{\em\ref{Theorem:OSPOD-1}}.}
\begin{enumerate}
    \item [1)] Boundedness of the minimizing sequence: Let
    \begin{align*}
        \bar J=\inf\big\{J(y^\ell,u)\,\big|\,x=(\mathrm y^\ell,y,p,\psi,\lambda,u)\in\Xad\big\}.
    \end{align*}
    Since $J$ is bounded from below by $0$, we have $\bar J\ge0$. Consequently, there exists a minimizing sequence $\{x_n\}_{n\in\mathbb N}\subset\Xad$ with $x_n=(\mathrm y^\ell_n,y_n,p_n,\psi_n,\lambda_n,u_n)$ and
    \begin{equation}
        \label{MinSeq}
        \lim_{n\to\infty}J(x_n)=\bar J.
    \end{equation}
    Moreover, $J$ is radially unbounded with respect to the control variabe so that $\{u_n\}_{n\in\mathbb N}\subset\Uad$ is bounded in $\U$. We deduce from Remark~\ref{RemH1WEq}-1) that $\mathrm y^\ell$ is bounded in $H^1(0,T;\mathbb R^\ell)$. The bound depends on $\ell$, but $\ell$ is fixed here. It follows from \eqref{SIAM:Eq3.1.7} and Lemma~\ref{Lemma:HI-31}-1) that $\{y_n\}_{n\in\mathbb N}$, $y_n=\hat y+\mathcal Su_n$, and $\{p_n\}_{n\in\mathbb N}$, $p_n=\hat p+\mathcal Au_n$, are bounded sequences in $W(0,T)$. Clearly, $\{\psi_n\}_{n\in\mathbb N}$ is bounded in $\mathbb X^\ell$. Furthermore, it follows from
    \begin{align*}
        {\|\mathcal R(y_n,p_n)\|}_{\mathscr L(X)}&=\sup_{\|\psi\|_X=1}{\|\mathcal R(y_n,p_n)\psi\|}_X\\
        &=\sup_{\|\psi\|_X=1}\bigg\|\int_0^T{\langle y_n(t),\psi\rangle}_X\,y_n(t)+{\langle p_n(t),\psi\rangle}_X\,p_n(t)\,\mathrm dt\bigg\|_X\\
        &\le\sup_{\|\psi\|_X=1}\int_0^T{\|y_n(t)\|}_X^2{\|\psi\|}_X+{\|p_n(t)\|}_X^2{\|\psi\|}_X\,\mathrm dt\\
        &\le{\|y_n\|}_{L^2(0,T;V)}^2+{\|p_n\|}_{L^2(0,T;V)}^2\\
        &\le C\left({\|y_\circ\|}_H+{\|f\|}_{L^2(0,T;V')}+{\|\ydT\|}_H+{\|\ydQ\|}_{L^2(0,T;H)}+{\|u_n\|}_\U\right)
    \end{align*}
    that the sequence $\{\mathcal R(y_n,p_n)\}_{n\in\mathbb N}$ is bounded in $\mathscr L(X)$. Consequently,
    \begin{align*}
        \lambda_{i,n}=\lambda_{i,n}\,{\|\psi_{i,n}\|}_X={\|\lambda_{i,n}\psi_{i,n}\|}_X={\|\mathcal R(y_n,p_n)\psi_{i,n}\|}_X\le C\quad\text{for all }n\in\mathbb N
    \end{align*}
    implies that $\{\lambda_{i,n}\}_{n\in\mathbb N}$ is bounded in $\mathbb R$ as well for $i=1,\ldots,\ell$. Therefore, we have proved that the sequence $\{x_n\}_{n\in\mathbb N}\subset\Xad$ is bounded in the Hilbert space $\X$. Since $\X$ is a Hilbert space, there exist a subsequence of $\{x_n\}_{n\in\mathbb N}$ which we again denote by $\{x_n\}_{n\in\mathbb N}$ and an element $\bar x=(\bar y^\ell,\bar y,\bar p,\bar\psi,\bar\lambda,\bar u)\in\X$ satisfying
    \begin{equation}
        \label{WeakConvZ}
        x_n\rightharpoonup\bar x\text{ in }\X\text{ for }n\to\infty\text{ in }\mathscr X.
    \end{equation}
    \item [2)] Feasibility of the limit point $\bar x$: Since $\Uad$ is weakly closed, we have $\bar u\in\Uad$. From $y_n\rightharpoonup \bar y$ in $W(0,T)$ and $u_n\rightharpoonup \bar u$ in $\U$ for $n\to\infty$ it follows that
    \begin{align*}
        {\langle y_{n,t}(t),\varphi\rangle}_{V',V}&\to{\langle \bar y_t(t),\varphi\rangle}_{V',V},\quad a(t;y_n(t),\varphi)\to a(t;\bar y(t),\varphi),\\
        {\langle (\mathcal Bu_n)(t),\varphi\rangle}_{V',V}&\to{\langle (\mathcal B\bar u)(t),\varphi\rangle}_{V',V}
    \end{align*}
    for every $\varphi\in V$ a.e. in $[0,T]$. Hence,
    \begin{align*}
        e_3(\bar x)=0\text{ in }\mP_3'\quad\text{and}\quad e_4(\bar x)=0\text{ in }\mP_4.
    \end{align*}
    Using $y_n\rightharpoonup \bar y$ and $p_n\rightharpoonup \bar p$ in $W(0,T)$ for $n\to\infty$ implies that
    \begin{align*}
        {\langle p_{n,t}(t),\varphi\rangle}_{V',V}&\to{\langle \bar p_t(t),\varphi\rangle}_{V',V},\quad a(t;\varphi,p_n(t))\to a(t;\varphi,\bar p(t)),\\
        \sigma_1\,{\langle y_{n,t}(t),\varphi\rangle}_H&\to\sigma_1\,{\langle \bar y(t),\varphi\rangle}_{V',V}
    \end{align*}
    for every $\varphi\in V$ a.e. in $[0,T]$. We conclude that
    \begin{align*}
        e_5(\bar x)=0\text{ in }\mP_5'\quad\text{and}\quad e_6(\bar x)=0\text{ in }\mP_6.
    \end{align*}
    Next we argue that $\bar\lambda_i$, $1\le i\le\ell$, is an eigenvalue of $\mathcal R(\bar y,\bar p)$ with eigenvector $\bar\psi_i$. By Assumption~\ref{Ass:OSPOS-2} we have
    \begin{equation}
        \label{OSPOD:StrConv}
        (y_n,p_n)\to(\bar y,\bar p)\quad\text{in }L^2(0,T;V)\times L^2(0,T;V)\text{ for }k\to\infty.
    \end{equation}
    Therefore, it follows from \eqref{WeakConvZ} and \eqref{OSPOD:StrConv} that
    \begin{equation}
        \label{OSPOD-A}
        \begin{aligned}
            \lambda_{i,n}\,{\langle\psi_{i,n},\psi\rangle}_X&={\langle\mathcal R(y_n,p_n)\psi_{i,n},\psi\rangle}_X\\
            &=\int_0^T{\langle y_n(t),\psi_{i,n}\rangle}_X\,{\langle y_n(t),\psi\rangle}_X+{\langle p_n(t),\psi_{i,n}\rangle}_X\,{\langle p_n(t),\psi\rangle}_X\,\mathrm dt\\
            &\stackrel{k\to\infty}{\longrightarrow}\int_0^T{\langle\bar y(t),\bar\psi_i\rangle}_X\,{\langle\bar y(t),\psi\rangle}_X+{\langle \bar p(t),\bar\psi_i\rangle}_X\,{\langle\bar p(t),\psi\rangle}_X\,\mathrm dt\\
            &\qquad={\langle\mathcal R(\bar y,\bar p)\bar\psi_i,\psi\rangle}_X\quad\text{for all }\psi\in X.
        \end{aligned}
    \end{equation}
    Since $\{\lambda_n\}_{k\in\mathbb N}$ is a sequence in $\mathbb R^\ell$, we also have
    \begin{equation}
        \label{OSPOD:StrConv-c}
        \lambda^n\to\bar\lambda\quad\text{in }\mathbb R^\ell\text{ for }k\to\infty.
    \end{equation}
    Thus, passing to the limit on the left in \eqref{OSPOD-A} yields
    \begin{align*}
        \mathcal R(\bar y,\bar p)\bar\psi_i=\bar\lambda_i\bar\psi_i\quad\text{for }1\le i\le\ell.
    \end{align*}
    Consequently, $(\bar\psi_i,\bar\lambda_i)$, $i=1,\ldots,\ell$, is an eigenvector-eigenvalue pair for the operator $\mathcal R(\bar y,\bar p)$ so that we have
    \begin{align*}
        e_7(\bar x)=0\text{ in }\mP_7.
    \end{align*}
    Furthermore, by \eqref{WeakConvZ} we find that
    \begin{align*}
        \bar\lambda_i&=\lim_{n\to\infty}\lambda_{i,n}=\lim_{n\to\infty}{\langle\lambda_{i,n}\psi_{i,n},\psi_{i,n}\rangle}_X=\lim_{k\to\infty}{\langle\mathcal R(y_n,p_n)\psi_{i,n},\psi_{i,n}\rangle}_X\\
        &=\lim_{n\to\infty}\int_0^T\big|{\langle y_n,\psi_{i,n}\rangle}_X\big|^2+\big|{\langle p_n,\psi_{i,n}\rangle}_X\big|^2\,\mathrm dt\\
        &=\int_0^T\big|{\langle\bar y,\bar\psi_i\rangle}_X\big|^2+\big|{\langle\bar p,\bar \psi_i\rangle}_X\big|^2\,\mathrm dt={\langle\mathcal R(\bar y,\bar p)\bar \psi_i,\bar \psi_i\rangle}_X={\langle\bar \lambda_i\bar \psi_i,\bar \psi_i\rangle}_X=\bar \lambda_i\,{\|\bar \psi_i\|}_X^2
    \end{align*}
    for $1\le i\le\ell$. By assumption, we have $\bar\lambda_i>0$ for $1\le i\le\ell$. Thus, we have $\|\bar\psi_i\|_X=1$ and
    \begin{align*}
        e_8(\bar x)=0\text{ in }\mP_8.
    \end{align*}
    Since the norm is weakly lower semicontinuous, we deduce
    \begin{align*}
        1={\|\bar\psi_i\|}_X\le\liminf_{n\to\infty}{\|\psi_{i,n}\|}_X=1\quad\text{for }1\le i\le\ell
    \end{align*}
    which implies that
    \begin{equation}
        \label{OSPOD:NormConv}
        \lim_{n\to\infty}{\|\psi_{i,n}\|}_X={\|\bar\psi_i\|}_X\quad\text{for }1\le i\le\ell.
    \end{equation}
    Thus,
    \begin{equation}
        \label{OSPOD:PsiConv}
        \begin{aligned}
            \lim_{n\to\infty}{\|\psi_{i,n}-\bar\psi_i\|}_X^2&=\lim_{n\to\infty}{\langle\psi_{i,n}-\bar\psi_i,\psi_{i,n}-\bar\psi_i\rangle}_X\\
            &=\lim_{n\to\infty}{\langle\psi_{i,n}-\bar\psi_i,\psi_{i,n}\rangle}_X+\lim_{n\to\infty}{\langle\bar\psi_i-\psi_{i,n},\bar\psi_i\rangle}_X\\
            &=\lim_{n\to\infty}{\|\psi_{i,n}\|}_X^2-\lim_{n\to\infty}{\langle\bar\psi_i,\psi_{i,n}\rangle}_X+\lim_{k\to\infty}{\langle\bar\psi_i-\psi_{i,n},\bar\psi_i\rangle}_X\\
            &={\|\bar\psi_i\|}_X^2-{\|\bar\psi_i\|}_X^2+0=0
        \end{aligned}
    \end{equation}
    because of \eqref{OSPOD:NormConv} and \eqref{WeakConvZ}. This also implies that
    \begin{equation}
        \label{OSPOD:PsiOrtho}
        {\langle\bar\psi_i,\bar\psi_j\rangle}_X=\delta_{ij}\quad\text{for }1\le i,j\le\ell.
    \end{equation}
    Moreover, we have  $\psi_{i,n}\to\bar\psi$ in $V$ for $1\le i\le\ell$: In fact, for the choice $X=V$ we have already proved this claim. Let u consider the case $X=H$. We define the four functions
    \begin{align*}
        r_n&=\frac{1}{\lambda_{i,n}}\,{\langle y_n(\cdot),\psi_{i^n}\rangle}_H\in L^2(0,T),&\bar r&=\frac{1}{\bar\lambda_i}\,{\langle \bar y(\cdot),\bar\psi_i\rangle}_H\in L^2(0,T),\\
        s_n&=\frac{1}{\lambda_{i,n}}\,{\langle p_n(\cdot),\psi_{i^n}\rangle}_H\in L^2(0,T),&\bar s&=\frac{1}{\bar\lambda_i}\,{\langle \bar p(\cdot),\bar\psi_i\rangle}_H\in L^2(0,T).
    \end{align*}
    It follows from \eqref{WeakConvZ} and \eqref{OSPOD:StrConv} that
    \begin{equation}
        \label{OSPOD:srConv}
        \lim_{n\to\infty}r_n=\bar r\text{ in }L^2(0,T)\quad\text{and}\quad\lim_{n\to\infty}s_n=\bar s\text{ in }L^2(0,T).
    \end{equation}
    Therefore, we derive from \eqref{OSPOD:StrConv} and \eqref{OSPOD:PsiConv} and \eqref{OSPOD:srConv}
    \begin{align*}
        {\|\psi_{i,n}-\bar\psi_i\|}_V&=\bigg\|\frac{1}{\lambda_{i,n}}\,\mathcal R(y_n,p_n)\psi_{i,n}-\frac{1}{\bar\lambda_i}\,\mathcal R(\bar y,\bar p)\bar\psi_i\bigg\|_V\\
        &=\bigg\|\frac{1}{\lambda_{i,n}}\int_0^T{\langle y_n(t),\psi_{i,n}\rangle}_X\,y_n(t)+{\langle p_n(t),\psi_{i,n}\rangle}_X\,p_n(t)\,\mathrm dt\\
        &\qquad-\frac{1}{\bar\lambda_i}\int_0^T{\langle \bar y(t),\bar\psi_i\rangle}_X\,\bar y(t)+{\langle \bar p(t),\bar\psi_i\rangle}_X\,\bar p(t)\,\mathrm dt\bigg\|_V\\
        &\le\int_0^T\big\|r_n(t)\,y_n(t)+s_n\,p_n(t)-\bar r(t)\,\bar y(t)-\bar s(t)\,\bar p(t)\big\|_V\,\mathrm dt\\
        &\le\int_0^T\big|r_n(t)-\bar r(t)\big|{\|y_n(t)\|}_V+|\bar r(t)|{\|y_n(t)-\bar y(t)\|}_V\,\mathrm dt\\
        &\quad+\int_0^T\big|s_n-\bar s(t)\big|{\|p_n(t)\|}_V+|\bar s(t)|{\|p_n(t)-\bar p(t)\|}_V\,\mathrm dt\\
        &\le{\|r_n-\bar r\|}_{L^2(0,T)}{\|y_n\|}_{L^2(0,T,V)}+{\|\bar r\|}_{L^2(0,T)}\,{\|y_n-\bar y\|}_{L^2(0,T;V)}\\
        &\quad+{\|s^{n_k}-\bar s\|}_{L^2(0,T)}{\|p_n\|}_{L^2(0,T,V)}+{\|\bar s\|}_{L^2(0,T)}\,{\|p_n-\bar p\|}_{L^2(0,T;V)}\stackrel{n\to\infty}{\longrightarrow}0
    \end{align*}
    for $1\le i\le\ell$. Using the weak convergence of the sequence $\{\mathrm y_n^\ell\}_{n\in\mathbb N}$ in $H^1(0,T;\mathbb R^\ell)$ and the strong convergence of $\{\phi_{i,n}\}_{n\in\mathbb N}$, $1\le i\le \ell$, we derive that
    \begin{align*}
        e_1(\bar x)=0\text{ in }\mP_1',\quad e_2(\bar x)=0\text{ in }\mP_2.
    \end{align*}
    holds. Summarizing we have proved that $\bar u\in\Uad$ and $e(\bar x)=0$ in $\mP'$ is true which gives $\bar x\in\Xad$.
    \item [3)] Optimality of $\bar z$: Since the cost functional is weakly lower semicontinuous, we have
    \begin{align*}
        0\ge \bar J=\inf\big\{J(y^\ell,u)\,\big|\,x\in\Zad\big\}=\lim_{n\to\infty}J(x_n)=\liminf_{n\to\infty}J(x_n)\ge J(\bar x)
    \end{align*}
    which gives the claim.\hfill$\Box$
\end{enumerate}


\chapter{Advanced Topics in POD Suboptimal Control}
\label{Advanced Topics in POD Suboptimal control}

\section{Optimal control of nonlinear evolution problems}
\label{Section:5.1}

It has been shown in the previous chapters that POD can be easily and efficiently applied to optimal control problems with linear state equations. This linearity made it straightforward to define a POD Galerkin scheme in \eqref{SIAM:Eq3.1.6POD} which could be evaluated exclusively in the low-order subspace $V^\ell$ by the system of ordinary differential equations \eqref{SIAM:PODStateDis} for the basis representation. In this chapter, we will elaborate on the difficulties that arise with a nonlinear state equation and how these can be overcome. For the sake of simplicity, we work with the particular semilinear example below, even though the presented strategies extend to the general case. 
	
\subsection{State equation} 
	
Suppose that $\Omega \subset \mathbb R^d$ is a bounded Lipschitz domain with boundary $\Gamma = \partial \Omega$ and $T>0$ is a final time. As before, we denote the space-time cylinder by $Q = (0,T) \times \Omega$ and the boundary-time cylinder by $\Sigma = (0,T) \times \Gamma$. We consider the semilinear heat equation\index{Problem!nonlinear heat equation}\index{Equation!nonlinear heat}
\begin{align}
	\label{eq:semilin_stateEq1}
	\begin{aligned}
		y_t(t,\bx) - c \Delta y(t,\bx) + y(t,x)^3 &= \sum_{i=1}^\md u_i(t) \chi_i(\bx), &&(t,\bx) \in Q, \\
		c \frac{\partial y}{\partial \bn}(t,\bs) &= 0, &&(t,\bs) \in \Sigma, \\
		y(0,\bx) &= y_\circ(\bx), &&x \in \Omega.
	\end{aligned}
\end{align}	
In contrast to the general example in Section \ref{SIAM-Book:Section1.2}, there is now a constant diffusivity $c > 0$ and a nonlinear cubic term that is evaluated at the state value $y(t,\bx)$. Because there are no nonlinear evaluations of $y$-derivatives, the above system is referred to as a \textit{semilinear} equation. All remaining terms are identical to Section \ref{SIAM-Book:Section1.2}, to which we refer for further details. We use the typical Gelfand triple given by the spaces $V = H^1(\Omega)$, $H = L^2(\Omega)$ and define the solution space $\Y = W(0,T) \cap L^\infty(\overline Q)$, where the space $W(0,T) = H^1(0,T;V') \cap L^2(0,T;V)$ is introduced in Section \ref{App:A.4}. Furthermore, the control space is given by $\U = L^2(0,T;\mathbb R^\md)$, with an inner product on $\mathbb R^m$ that is still left to be specified. Let us assume that it holds $\chi_1,...,\chi_\md \in L^\infty(\Omega)$ and $y_\circ \in C(\overline \Omega)$. 

These conditions allow us to define the symmetric bilinear form $a: V \times V \to \mathbb R$ as
\begin{subequations}
    \begin{align}
    	\label{eq:semilin_maps_1}
    	a(\varphi,\psi) = c \int_\Omega \nabla \varphi(t,\bx) \cdot \nabla \psi(\bx) ~\mathrm d\bx \qquad \forall \varphi,\psi \in V.
    \end{align}
    For all $u \in \U$ and almost all $t \in (0,T)$, we also define the time-dependent linear form $(\mathcal B u)(t): V \to \mathbb R$ by 
    \begin{align}
    	\label{eq:semilin_maps_3}
    	{\langle (\mathcal B u)(t),\psi \rangle}_{V',V} = \sum_{i=1}^\md u_i(t) \int_\Omega \xi_i(\bx) \psi(\bx) ~\mathrm d\bx \quad \text{for all }\psi \in V.
    \end{align}
    Finally, let us denote for $z: Q \to \mathbb R$ the function $\mathcal N(z): Q \to \mathbb R$ as
    \begin{align}
    	\label{eq:semilin_maps_4}
    	\mathcal N(z)(t,\bx) = z(t,\bx)^3 \quad \text{f.a.a. } \bx \in Q.
	\end{align}
\end{subequations}
Similar to the function $y(t)$ standing for $y(t,\cdot): \Omega \to \mathbb R$ for $t \in (0,T)$, we also abbreviate the function $\mathcal N(y(t)) := \mathcal N(y(t,\cdot)): \Omega \to \mathbb R$. 
	
\begin{definition}
	A function $y \in \Y$ is called a weak solution of \eqref{eq:semilin_stateEq1} if it satisfies for all $\varphi \in V$:
    \begin{align}
		\label{eq:semilin_stateEq2}
		\begin{aligned}
			\frac{\mathrm d}{\mathrm dt} \langle y(t), \varphi \rangle_{H} + a(y(t),\varphi) + \langle \mathcal N(y(t)),\varphi \rangle_H &= \langle (\mathcal Bu)(t), \varphi \rangle_{V',V} &&\text{f.a.a. } t \in (0,T), \\
			\langle y(0), \varphi \rangle_{H} &= \langle y_\circ, \varphi \rangle_H. &&
		\end{aligned}
	\end{align}
\end{definition}
	
Notice that the above definition is well-defined since it follows from $y(t) \in L^\infty(\Omega)$ for almost all $t \in (0,T)$ that also $\mathcal N(y(t)) \in L^\infty(\Omega)$, which embeds into $L^2(\Omega)$. 
	
\begin{theorem}
	\label{thm:semilin_solOp}
	For every $u \in L^r(0,T;\mathbb R^\md)$ with $r>d/2 + 1$, the weak state equation \eqref{eq:semilin_stateEq2} admits a unique solution $y \in \Y \cap C(\overline Q)$ which satisfies the energy estimate
    \begin{align}
		\label{eq:energy_solOp_cubic}
		{\|y\|}_{W(0,T)} + {\|y\|}_{C(\overline Q)} \le c_\infty \big( \|u\|_{L^r(0,T;\mathbb R^\md)} + {\|y_\circ\|}_{C(\overline Q)} \big) 
	\end{align}
    with a constant $c_\infty>0$ that does not depend on $u$ or $y_\circ$.
\end{theorem}

\noindent{\bf\em Proof.}
Let us define the inhomogeneity function
\begin{align*}
	f: Q \to \mathbb R, \quad f(t,\bx) = \sum_{i=1}^\md u_i(t) \chi_i(\bx).
\end{align*}
We can utilize the fact that $\chi_i \in L^\infty(\Omega)$ for $i=1,...,\md$ and easily verify that $f \in L^r(Q)$ with
\begin{align*}
	{\|f\|}_{L^r(Q)} \le c\,{\|u\|}_{L^r(0,T;\mathbb R^\md)}
\end{align*}
with a constant $c>0$ that does not depend on $u$. Inserting these two properties into \cite[Theorem 5.5]{Tro10} gives the claim.\hfill$\Box$

From Theorem \ref{thm:semilin_solOp}, we get a well-defined solution operator $\mathcal S: \U \to \Y \cap C(\overline \Omega)$ which is defined on $\mathcal D(\mathcal S) := L^\infty(0,T;\mathbb R^{\md})$.

\subsection{Semidiscretization by FE and POD}
\label{sec:cubicSemidiscretization}
	
As it was extensively studied in Section \ref{SIAM-Book:Section3.4} for the linear problem class, the standard way to proceed is to introduce a space discretization, for which we choose linear Finite Elements (FE). Let therefore 
\begin{align*}
	V^h = \text{span } \{ \varphi_1^h,\hdots,\varphi_m^h \} \subset V
\end{align*}
denote the ansatz space. The semidiscrete FE Galerkin discretization of \eqref{eq:semilin_stateEq2} is given by a function $y^h \in H^1(0,T;V^h) \cap L^\infty(\overline Q)$ which solves for all $\varphi^h \in V^h$
\begin{align}
	\label{eq:semilin_stateEq_FE}
	\begin{aligned}
		\frac{\mathrm d}{\mathrm dt}\,{\langle y^h(t),\varphi^h \rangle}_{H} + a(y^h(t),\varphi^h) + {\langle \mathcal N(y^h(t)),\varphi^h \rangle}_H&= {\langle (\mathcal Bu)(t), \varphi^h \rangle}_{V',V} \quad \text{f.a.a. } t \in (0,T), \\
		{\langle y(0), \varphi^h \rangle}_{H} &= {\langle y_\circ, \varphi^h \rangle}_H.
	\end{aligned}
\end{align}
Along the lines of Section \ref{sec:nonlinearGalerkinDiscretiation}, we can show that there exists an equivalent algebraic ODE system to \eqref{eq:semilin_stateEq_FE}. If $\mathrm y^h: (0,T) \to \mathbb R^m$ is the coordinate function to $y^h$, it satisfies (cf. \eqref{NonlFineModel})
\begin{equation}
	\label{eq:semilin_stateEq_FE2}
	\begin{aligned}
		\bM^h\dot{\mathrm y}^h(t)+\bA^h\mathrm y^h(t)+\mathrm n^h\big(\mathrm y^h(t)\big)&=\bB^hu(t)\quad\text{for }t\in(0,T],\\
		\mathrm y^h(0)&=\mathrm y_\circ^h.
	\end{aligned}
\end{equation}
Here, the matrices $\bM^h$ and $\bA^h \in \mathbb R^{m \times m}$ are given by 
\begin{align*}
	\bM^h_{ij} = \int_\Omega \varphi_i^h(\bx) \varphi_j^h(\bx) ~\mathrm d\bx, \quad \bA_{ij}^h = a(\varphi^h_i,\varphi^h_j), \quad \text{for all } i,j=1,...,m.
\end{align*}
The matrix $\bB^h \in \mathbb R^{m \times \md}$ is given by the elements
\begin{align*}
	\bB^h_{ij} = \int_\Omega \xi_j(\bx) \varphi^h_i(\bx) ~\mathrm d\bx \quad \text{for all } i=1,\hdots,m, ~j=1,\hdots,\md.
\end{align*}
The initial condition is the coordinate representation of the $H$-orthogonal projection $y_\circ^h = \mathcal P^h y_\circ$ onto $V^h$. Finally, the nonlinearity $\mathrm n^h: \mathbb R^m \to \mathbb R^m$ reads
\begin{align*}
	\mathrm n^h(\mathrm v) = \bigg( \int_\Omega \Big( \sum_{j=1}^m \mathrm v_j \varphi_j^h(\bx) \Big)^3 \varphi_i^h(\bx) ~\mathrm d\bx \bigg)_{1 \le i \le m}
\end{align*}
Since a correct evaluation of this function would force us to leave the linear FE space, we make the nodal approximation
\begin{align*}
	\mathrm n^h(\mathrm v) \approx \bigg( \sum_{j=1}^m \int_\Omega \mathrm v_j^3 \varphi_j^h(\bx) \varphi_i^h(\bx) ~\mathrm d\bx \bigg)_{1 \le i \le m} = \bM^h \mathrm N^h(\mathrm v),
\end{align*}
where we have defined the coordinate nonlinearity $\mathrm N^h: \mathbb R^m \to \mathbb R^m$ with $\mathrm N^h_i(\mathrm v) = \mathrm v_i^3$ for $1 \le i \le m$. In replacement of \eqref{eq:semilin_stateEq_FE2}, we now get the system
\begin{equation}
	\label{eq:semilin_stateEq_FE3}
	\begin{aligned}
		\bM^h\dot{\mathrm y}^h(t)+\bA^h\mathrm y^h(t)+ \bM^h \mathrm N^h(y^h(t))&=\bB^hu(t)\quad\text{for }t\in(0,T],\\
		\mathrm y^h(0)&=\mathrm y_\circ^h.
	\end{aligned}
\end{equation}
In a next step, let us suppose that we have computed a POD subspace
\begin{align*}
	V^{h\ell} = \text{span } \{\varphi^{h\ell}_1, \hdots, \varphi^{h\ell}_\ell \} \subset V^h \subset V
\end{align*}
of dimension $\ell \ll m$. The POD Galerkin discretization of \eqref{eq:semilin_stateEq_FE} is given by a function $y^{h\ell} \in H^1(0,T;V^{h\ell}) \cap L^\infty(Q)$ which solves for all $\varphi^{h\ell} \in V^{h\ell}$
\begin{align}
	\label{eq:semilin_stateEq_POD}
	\begin{aligned}
		\frac{\mathrm d}{\mathrm dt}\,{\langle y^{h\ell}(t), \varphi^{h\ell} \rangle}_{H} + &a(y^{h\ell}(t),\varphi^{h\ell}) + {\langle \mathcal N(y^{h\ell}(t)),\varphi^{h\ell} \rangle}_H&= {\langle (\mathcal Bu)(t), \varphi^{h\ell} \rangle}_{V',V} \qquad \text{f.a.a. } t \in (0,T), \\
		{\langle y(0), \varphi^{h\ell} \rangle}_{H} &= \langle y_\circ, \varphi^{h\ell} \rangle_H.
	\end{aligned}
\end{align}
With the same nodal appoximation for the nonlinearity, we get the ODE system for the POD coordinate solution $\mathrm y^\ell: (0,T) \to \mathbb R^\ell$ 
\begin{align}
	\label{eq:semilin_stateEq_POD2}
	\begin{aligned}
		\bM^{h\ell} \dot{\mathrm y}^{h\ell}(t) + \bA^{h\ell} \mathrm y^{h\ell}(t) + (\bPsi^{h\ell})^\top \mathrm N^h \big(\bPsi^{h\ell} \mathrm y^\ell(t)\big) &= \bB^{h\ell} u(t) \quad &\text{for } t \in (0,T], \\
		\bM^{h\ell} \mathrm y^{h\ell}(0) &= \mathrm y_\circ^{h\ell}.
	\end{aligned}
\end{align}
The matrix $\bPsi^{h\ell} \in \mathbb R^{m \times \ell}$ holds the coordinates of the POD basis vectors such that
\begin{align*}
	\varphi_j^{h\ell} = \sum_{i=1}^m \bPsi^{h\ell}_{ij} \varphi^h_i, \qquad \text{for } j=1,\hdots,\ell.
\end{align*}
The matrices $\bM^{h\ell}$ and $\bA^{h\ell} \in \mathbb R^{\ell \times \ell}$ indicate the mass and stiffness matrices of the POD system and are given by
\begin{align*}
	\bM^{h\ell} = \Big( \int_\Omega \varphi^{h\ell}_i(\bx) \varphi^{h\ell}_j(\bx) ~\mathrm d\bx \Big)_{1 \le i,j \le \ell} =(\bPsi^{h\ell})^\top \bM^h \bPsi^{h\ell},
	~\bA^{h\ell} = \Big( a(\varphi^{h\ell}_i, \varphi^{h\ell}_j) \Big)_{1 \le i,j \le \ell}=(\bPsi^{h\ell})^\top \bA^h \bPsi^{h\ell}.
\end{align*}
The same holds true for the control matrix $\bB^{h\ell} \in \mathbb R^{\ell \times \md}$, which satisfies
\begin{align*}
	\bB^{h\ell} = \Big( \int_\Omega \xi_j(\bx) \varphi_j^{h\ell}(\bx) ~\mathrm d\bx \Big)_{\substack{1 \le i \le \ell \\ 1 \le j \le \md}} = (\bPsi^{h\ell})^\top \bB^h.
\end{align*}
The initial vector $\mathrm y^{h\ell}_\circ \in \mathbb R^\ell$ is given by the coordinate representatin of the $H$-orthogonal projection $y_\circ^{h\ell} = \mathcal P^{h\ell} y_\circ^h$ onto $V^{h\ell}$. It can be easily shown that it is implicitly given be the linear system $\bM^{h\ell} \mathrm y_\circ^{h\ell} = (\bPsi^{h\ell})^\top \bM^h \mathrm y_\circ^h$. 
		
The efficiency-related problem in \eqref{eq:semilin_stateEq_POD2} concerning the nonlinearity was already mentioned in Section \ref{SIAM-Book:Section3.6.5}. It lies in the fact that it enforces us to enhance the low-dimensional vector $\mathrm y^\ell(t)$ to the FE size $m$ by multiplication with $\bPsi^{h\ell}$, perform a full-order evaluation of the nonlinearity $\mathrm N$ and then re-compress the result to $\mathbb R^\ell$ by multiplication with $(\bPsi^{h\ell})^\top$. As a workaround, both the Empirical Interpolation method (EIM, Algorithm \ref{algo:EIM} and the Discrete Empirical Interpolation Method (DEIM, Algorithm \ref{algo:DEIM}) are available. After the computation of a matrix $\bPhi^{\mathsf p} \in \mathbb R^{m \times \mathsf p}$ and a projection matrix $\bP \in \mathbb R^{m \times \mathsf p}$, we have the interpolated nonlinearity
\begin{align*}
	\mathrm N^{\mathsf p \ell}: \mathbb R^{\ell} \to \mathbb R^{\mathsf p}, \quad \mathrm N^{\mathsf p \ell}(\mathrm y^{h\ell}) = (\bP^\top \bPhi^{\mathsf p})^{-1} \bP^\top \mathrm N^h(\bPsi^{h\ell} \mathrm y^{h\ell}),
\end{align*}

which was indicated by the vector $\mathrm c$ in Section \ref{SIAM-Book:Section3.6.5}. Since multiplication with the projection matrix $\bP^\top$ only selects $\mathsf p$ components of a vector in $\mathbb R^m$ and the nonlinearity $\mathrm N^h$ works componentwise, the above function satisfies
\begin{align*}
	\mathrm N^{\mathsf p \ell}(\mathrm y^{h\ell}) = (\bP^\top \bPhi^{\mathsf p})^{-1} \mathrm N^{\mathsf p}( \bP^\top \bPsi^{h\ell} \mathrm y^{h\ell}),
\end{align*}
where the nonlinearity $\mathrm N^{\mathsf p}: \mathbb{R}^{\mathsf p} \to \mathbb R^{\mathsf p}$ is given by $\mathrm N_i^\ell \mathrm v = \mathrm v_i^3$ for $i=1,\hdots,p$. Replacing the system \eqref{eq:semilin_stateEq_POD2} with the inclusion of this interpolated nonlinearity, we end up with
\begin{align}
	\label{eq:semilin_stateEq_POD3}
	\begin{aligned}
		\bM^{h\ell} \dot{\mathrm y}^{h\ell}(t) + \bA^{h\ell} \mathrm y^{h\ell}(t) + (\bPsi^{h\ell})^\top \bM^{h\ell} \bPhi^{\mathsf p} (\bP^\top \bPhi^{\mathsf p})^{-1} \mathrm N^{\mathsf p} \big(\bP^\top \bPsi^{h\ell} \mathrm y^\ell(t)\big)&= \bB^{h\ell} u(t), \\
		\bM^{h\ell} \mathrm y^{h\ell}(0) &= \mathrm y_\circ^{h\ell}.
	\end{aligned}
\end{align}
for $t \in (0,T]$. Even though this may look more complex than the original version \eqref{eq:semilin_stateEq_POD2}, the matrices $(\bPsi^{h\ell})^\top \bM^{h\ell} \bPhi^{\mathsf p} \in \mathbb R^{\ell \times \mathsf p}$, $\bP^\top \bPhi^{\mathsf p} \in \mathbb R^{\mathsf p \times \mathsf p}$ and $\bP^\top \bPsi^{h \ell} \in \mathbb R^{\mathsf p \times \ell}$ are of reduced dimension and can be pre-computed to save on efficiency. The remaining nonlinearity $\mathrm N^{\mathsf p}$ only operates in the (D)EIM dimension $\mathsf p \ll m$, so the complexity is massively reduced compared to \eqref{eq:semilin_stateEq_POD2}. 
		
\subsection{Cost function and optimal control problem}
	
For the cost function, we choose a desired final-time state $\yd \in C(\overline Q)$ and set 
\begin{align*}
	J: \Y \times \U \to \mathbb R, \quad J(y,u) = \frac{1}{2} \int_\Omega \big| y(T,\bx) - \yd(\bx) \big|^2 ~\mathrm d\bx + \frac{\sigma}{2} \lVert u \rVert_\U^2.
\end{align*}

As usual, the reduced cost function is defined as $\Jhat(u) = J(\mathcal S(u),u)$ for all $u \in L^\infty(0,T;\mathbb R^\md)$, cf. Theorem \ref{thm:semilin_solOp}. On the control space $\U$, we impose the partial order relation $u \le \tilde u$ as in Remark \ref{Remark:HI-30}. For two bilateral constraints given by bounding variables $\ua, \ub \in L^\infty(0,T;\mathbb R^\md)$ with $\ua \le \ub$, the admissible set is given as in Definition \ref{Definition:HI-1001} by 
\begin{align*}
	\Uad = \big \{ u \in \U ~\big|~ \ua \le u \le \ub \big \}
\end{align*}
Notice that $\Uad$ is a bounded subset of $L^\infty(0,T;\mathbb R^\md)$. Therefore, the solution operator $\mathcal S$ and with it the reduced cost function $\Jhat$ are well-defined on $\Uad$. Within this admissible set, the continuous optimal control problem is now given by\index{Problem!nonlinear programming}
\begin{align}
	\label{eq:ocProb_cubic}
	\min\limits_{u \in \U} \quad \Jhat(u) \quad \text{s.t.} \quad u \in \Uad.
\end{align}
Notice that unlike the main linear-quadratic problem \eqref{GVLuminy:Eq4.1.9} considered in Chapter \ref{SIAM-Book:Section4}, this problem is not convex because of the nonlinearity of the solution operator $\mathcal S$ contained in the function $\Jhat$. Nevertheless, a result about solvability of this problem can be derived.
	
\begin{theorem}
	The constrained optimal control problem \eqref{eq:ocProb_cubic} admits at least one globally optimal solution $\bar u \in \Uad$.
\end{theorem}

\noindent{\bf\em Proof.}
We follow along the lines of \cite[Theorem 5.7]{Tro10} and take $r > d/2+1$. Since $\Uad$ is bounded in $L^r(0,T;\mathbb R^\md)$, \eqref{eq:energy_solOp_cubic} gives us the existence of a constant $M>0$ with
\begin{align*}
	{\| \mathcal S(u) \|}_{C(\overline Q)} \le M \quad \text{for all } u \in \Uad.
\end{align*}
Next, we take a minimizing sequence $\{u_n\}_{n \in \mathbb N} \subset \Uad$ with 
\begin{align*}
	\Jhat(u_n) \to J^* := \inf\limits_{u \in \Uad} \Jhat(u) \quad \text{as } n \to \infty.
\end{align*}

This is possible since the cost function $\Jhat$ is nonnegative. We also define $y_n := \mathcal S(u_n)$. Recall that $\Uad$ is a closed, convex, and bounded subset of the reflexive Banach space $L^r(0,T;\mathbb R^\md)$. As such, it is weakly sequentially compact, cf. \cite[Theorem 2.11]{Tro10}. Therefore we have a weakly convergent subsequent of $\{u_n\}_{n \in \mathbb N}$, which without loss of generality we write again as $\{u_n\}_{n \in \mathbb N}$. Let the weak limit be denoted by $\bar u \in \Uad$. In a next step, we write the nonlinearity and the control input as
\begin{align*}
	z_n: Q &\to \mathbb R, & z_n(t,\bx) &= y(t,\bx)^3 &&\text{for } n \in \mathbb N, \\
	\mathcal B: L^\infty(0,T;\mathbb R^\md) &\to L^\infty(Q), &(\mathcal Bu)(t,\bx) &= \sum_{i=1}^\md u_i(t) \chi_i(\bx)
\end{align*}
With it, the state equation \eqref{eq:semilin_stateEq1} for the sequence states can be equivalently formulated as 
\begin{align}
	\label{eq:rewrittenCubicStateEq}
	\begin{aligned}
		y^n_t(t,\bx) - c \Delta y^n(t,\bx) &= -z_n(t,\bx) + (\mathcal B u_n)(t,\bx), &&(t,\bx) \in Q, \\
		c\,\frac{\partial y^n}{\partial \bn}(t,\bs) &= 0, &&(t,\bs) \in \Sigma, \\
		y^n(0,\bx) &= y_\circ(\bx), &&x \in \Omega.
	\end{aligned}
\end{align}
This linear equation possesses a well-defined, affine linear and continuous solution operator mapping from $L^2(Q) \to W(0,T)$, cf. Corollary \ref{cor:numericalExample_solvability}. As such, this solution operator is also weakly continuous (cf. \cite[p.45]{Tro10}). The same reasoning applies to $\mathcal B$, which gives $\mathcal B u_n \to \mathcal B \bar u$ in $L^\infty(Q)$ and therefore also in $L^2(Q)$. Furthermore, the sequence $\{z_n\}_{n \in \mathbb N} \subset L^\infty(Q)$ is bounded in $L^\infty(Q)$ and theefore has a weakly convergent subsequence. Without loss of generality, we again write $z_n \rightharpoonup \bar z \in L^\infty(Q) \hookrightarrow L^2(Q)$. Together, this gives us 
\begin{align*}
	-z_n+\mathcal B u_n \rightharpoonup -\bar z + \mathcal B \bar u \quad \text{in } L^2(Q) \quad \text{as } n \to \infty.
\end{align*}
With the already mentioned weak continuity of the solution operator of \eqref{eq:rewrittenCubicStateEq}, this gives us the weak convergence of $y_n$ in $W(0,T)$ to an element $\bar y \in W(0,T)$. The rest of the proof is identical to that of \cite[Theorem 5.7]{Tro10} and we will only outline the last steps: The fact that $y_\circ \in C(\overline \Omega)$ allows to conclude that $y_n$ even convergences strongly in $W(0,T)$ to $\bar y$. Further, we find by local Lipschitz continuity of the cubic nonlinearity that it holds $y_n(\cdot)^3 \to \bar y(\cdot)^3$ in $L^\infty(Q)$. This can be used to show $\bar y$ is the solution of the state equation \eqref{eq:semilin_stateEq2} for $u = \bar u$. Finally, this allows to conclude:
\begin{align*}
	J^* = \liminf_{n \to \infty} \Jhat(u^n) &= \frac{1}{2} \Big( \lim\limits_{n \to \infty} \int_\Omega \big| y^n(T,\bx)-\yd(\bx)\big|^2 ~\mathrm d\bx + \sigma \liminf_{n \to \infty} \|u^n\|_\U^2 \Big) \\
	&\ge \frac{1}{2} \Big( \int_\Omega \big| \bar y(T,\bx)-\yd(\bx)\big|^2 ~\mathrm d\bx + \sigma \|\bar u\|_\U^2 \Big) = \Jhat(\bar u).
\end{align*}
Notice that the strong convergence $y^n \to \bar y$ in $W(0,T)$ extends to strong convergence $y^n(T) \to \bar y(T)$ in $H = L^2(\Omega)$ due to the embedding $W(0,T) \hookrightarrow C([0,T];H)$, such that the first term converges exactly. Furthermore, the weak convergence $u^n \rightharpoonup \bar u$ together with the lower semicontinuity of the square norm as a convex function guarantees the estimate for the second term. Since $J^*$ is the infimum over all cost function values in $\Uad$, $\bar u$ is a global minimizer of $\Jhat$ in $\Uad$. \hfill$\Box$
			
\medskip
In a next step, we will comment on the differentiability of the solution operator $\mathcal S: \U \to \Y$ and the reduced cost function $\Jhat: \U \to \mathbb R$. Proofs of the Lemmas  \ref{lem:cubicSolutionOperatorDerivatives} and \ref{lem:cubicReducedCostDerivatives} can be found in \cite{Tre17} or \cite[Section 1.6]{HPUU09}.
	
\begin{lemma}
	\label{lem:cubicSolutionOperatorDerivatives}
	The solution operator $\mathcal S$ is twice continuously Fréchet differentiable in every $u \in L^\infty(0,T;\mathbb R^\md)$ and every direction $h \in L^\infty(0,T;\mathbb R^\md)$. Let $y = \mathcal S(u)$ and $y^h = \mathcal S'(u)h$. Then $y^h$ is the solution to the \textit{linearized state equation}
	\begin{align}
		\label{eq:cubicLinearizedStateEq}
		\begin{aligned}
			\frac{\mathrm d}{\mathrm dt}\,{\langle y^h(t), \varphi \rangle}_H + a(y^h(t),\varphi) + {\langle \mathcal N'(y(t)) y^h(t), \varphi \rangle}_H&={\langle (\mathcal B h)(t),\varphi  \rangle}_{V',V},\\
			y^h(0) &= 0.
		\end{aligned}
	\end{align}
    for all $\varphi \in V$ and f.a.a. $t \in (0,T]$.
\end{lemma}
	
\begin{remark}
	\rm
	The nonlinearity in \eqref{eq:cubicLinearizedStateEq} is given as
    \begin{align*}
		{\langle \mathcal N'(y(t)) y^{h}(t), \varphi \rangle}_H = 3\int_\Omega y(t,\bx)^2 y^{h}(t,\bx) \varphi(\bx) ~\mathrm d\bx.
	\end{align*}
	\hfill$\blacksquare$
\end{remark}
	
\begin{lemma}
	\label{lem:cubicReducedCostDerivatives}
	The reduced cost function $\Jhat$ is twice continuously Fréchet differentiable in every $u \in L^\infty(0,T;\mathbb R^\md)$. Its $\U$-gradient $\nabla \Jhat(u) \in \U$ is given by
	\begin{align*}
		[\nabla \Jhat(u)]_i(t) = \int_\Omega \chi_i(\bx) p(t,\bx) ~\mathrm d\bx + \sigma u_i, \quad \text{for } i=1,...,\md, ~\text{f.a.a. } t \in (0,T).
	\end{align*}
	Here, $p \in \Y$ is the solution to the \textit{adjoint equation} for all $\varphi \in V$.
	\begin{align}
		\begin{aligned}
			-\frac{\mathrm d}{\mathrm dt}\,{\langle p(t), \varphi \rangle}_H + a(\varphi, p(t)) + {\langle N'(y(t)) p(t), \varphi \rangle}_H &= 0, &&\text{f.a.a. } t\in (0,T), \\
			p(T) &= y(T) - \yd &&\text{in } H.
		\end{aligned}
	\end{align}
\end{lemma}
	
Of course, there exist similar representations for the second derivatives of $\mathcal S$ and $\Jhat$, cf. \cite{Tre17} or \cite[Section 1.6]{HPUU09}.
	
\subsection{POD error estimation}
	
To find a solution to \eqref{eq:ocProb_cubic}, we want to employ a descent method that is based on first-order gradient information and incorporates the POD-reduced model from Section \ref{sec:cubicSemidiscretization}. After such a model is built, the objective is to avoid any superfluous evaluations of the expensive FE model, while at the same time keeping track of how well the reduced model works in context of the optimization. In order to do this, we have the following a-posteriori error estimator from \cite{KTV13} available to us.
	
\begin{theorem}
	\label{thm:cubicErrEst}
	Let $\bar u \in \Uad$ be a locally optimal control for \eqref{eq:ocProb_cubic} that satisfies the following coercivity condition for some $\kappa>0$.
    \begin{align}
		\label{eq:cubicCoercivityAss}
		\Jhat''(\bar u) (v,v) \ge \kappa \|v\|_\U^2 \quad \text{for all } v \in \U.
	\end{align}
    Then there exists a radius $r>0$ such that for all $u \in \U$ with $\|\bar u-u\| \le r$,  it holds
    \begin{align}
		\label{eq:cubicErrEst}
		{\|\bar u - u\|}_\U \le \frac{2}{\kappa}\,{\| \zeta_u \|}_\U.
	\end{align}
    Here, $\zeta_u \in \U$ is given by the \textit{perturbation function}
    \begin{align}
		[\zeta_u]_i(t) = \left \{
        \begin{array}{r l}
			\max \big(0,-[\nabla \Jhat(u)]_i(t) \big), &\text{if } [u]_i(t) = [\ua]_i(t), \\
			\min \big(0,-[\nabla \Jhat(u)]_i(t) \big), &\text{if } [u]_i(t) = [\ub]_i(t), \\
			-[\nabla \Jhat(u)]_i(t), &\text{otherwise}, \\
        \end{array}
        \right .
	\end{align}
\end{theorem}
	
The main problem in Theorem \ref{thm:cubicErrEst} when compared to Theorem \ref{GVLuminy:Theorem4.5.1} for the linear-quadratic case lies in the coercivity constant $\kappa>0$ and the unspecified radius $r>0$ in which \eqref{eq:cubicErrEst} holds. From a numerical reality, actually evaluating this constant is a quite expensive and complex task. We refer to \cite{Tre17}, where this manner has been thoroughly investigated. In some cases, the concrete problem structure can be exploited and $\kappa$ is analytically known. An example for such a case will come up in Section \ref{Section:MOP} where the reference point function - a quartic nonlinearity - represents a nonquadratic cost function whose coercivity constant can nevertheless be estimated in the linear-quadratic multiobjective case. 
	
\section{State-constrained model predictive control}
\label{Section:5.2}
The basic problem \eqref{GVLuminy:Eq4.1.6} has only bilateral constraints on the control variable $u$. In this section, we will examine the case when the state variable $y$ is likewise constrained by pointwise lower and upper bounds. Further, we present a Model Predictive Control algorithm to handle long-time horizon and parameter updates. In the last part, we show how to combine the previous methodology with the POD method in a successful way, utilizing an a-posteriori error estimate.   

\subsection{Problem formulation}
\label{Section:5.2.1}

We introduce the space $\W=L^2(Q)$. Let $\mathcal E:L^2(0,T;H)\to\mathscr W$ be a linear, bounded operator mapping to this space. For a given state $y\in L^2(0,T;H)$, we deal with pointwise constraints of the following type
\begin{equation}
    \label{Section5.2:StateConstraints}
    \ya(t,\bx)\le (\mathcal Ey)(t,\bx)\le \yb(t,\bx)\quad\text{a.e. in } Q,
\end{equation}
where $\ya$, $\yb\in\W$ are given lower and upper bounds with $\ya \leq \yb$ a.e. in $Q$. To obtain regular Lagrange multipliers we utilize a virtual control approach \cite{KR09,Mec19,MV18}. Let $\varepsilon>0$ be chosen and $w\in\W$ an additional (artificial) control. Then \eqref{Section5.2:StateConstraints} is replaced by the mixed constraints
\begin{align*}
    \ya(t,\bx)\le\mathcal Ey(t,\bx)+\varepsilon w(t,\bx)\le \yb(t,\bx)\quad\text{ a.e in } Q.
\end{align*}
We introduce the Hilbert space
\begin{align*}
    \mathscr X=W(0,T)\times\U\times\W
\end{align*}
endowed with the common product topology. The quadratic cost functional $J:\mathscr X\to\mathbb R$ is given by
\begin{align*}
    J(x)=\frac{\sigma_1}{2}\int_0^T{\|y(t)-\ydQ(t)\|}_H^2\,\mathrm dt+\frac{\sigma_2}{2}\,{\|y(T)-\ydT\|}_H^2+\frac{\sigma}{2}\,{\|u-\un\|}_\U^2+\frac{\sigma_w}{2}\,{\|w\|}_\W^2
\end{align*}
for $x=(y,u,w)\in\mathscr X$. Furthermore, $\sigma_1,\,\sigma_2\geq0$, and $\sigma,\sigma_w>0$. The set of admissible solutions is given by
\begin{align*}
    \Xadeps=\big\{x=(y,u,w)\in\X\,\big|\,y=\hat y+\mathcal Su,~\ya\le\mathcal Ey+\varepsilon w\le \yb\text{ and }u\in\Uad\big\}.
\end{align*}
Then the optimal control problem reads
\begin{equation}
    \label{Section5.2:Pw}
    \tag{$\mathbf P_\varepsilon$}
    \min J(x)\quad\text{s.t.}\quad x\in \Xadeps.
\end{equation}
Problem \eqref{Section5.2:Pw} can also be formulated as a pure control constrained problem. We set $\hya=\ya-\mathcal E\hat y\in\W$ and $\hyb=\yb-\mathcal E\hat y\in\W$. Then \eqref{Section5.2:StateConstraints} can be formulated equivalently in the variables $u$ and $w$ as follows:
\begin{align}
    \label{Section5.2:eq:stateConstraints}
    \hya(t,\bx)\le(\mathcal E\mathcal Su)(t,\bx)+\varepsilon w(t,\bx)\le\hyb(t,\bx)\quad\text{ a.e. in } Q.
\end{align}
In light of this, let us define $\Z=\U\times\W$ and introduce the according bounded and linear mapping $\mathcal T:\Z\to\Z$
\begin{equation}
    \label{Section5.2:Teps}
    z=(u,w)\mapsto\mathcal T(z)=\left(
    \begin{array}{c}
        u\\
        \varepsilon w+\mathcal E\mathcal Su
    \end{array}
    \right)=\left(
    \begin{array}{cc}
        \mathcal I_\U&0\\
        \mathcal E\mathcal S&\varepsilon\mathcal I_\W
    \end{array}
    \right)\left(
    \begin{array}{c}
        u\\
        w
    \end{array}
    \right),
\end{equation}
where $\mathcal I_\U:\U\to\U$ and $\mathcal I_\W:\W\to\W$ stand for the identity operators in $\U$ and $\W$, respectively. Using this operator and $\hat z_a=(u_a,\hya)$, $\hat z_b=(u_b,\hyb)\in\Z$, \eqref{Section5.2:eq:stateConstraints} can equivalently be formulated as 
\begin{align*}
    z=(u,w) \in \Zad=\big\{z=(u,w)\in\Z\,\big|\,\hat z_a\le\mathcal T(z)\le\hat z_b\big\}.
\end{align*}
where $\Zad$ is closed and convex; cf. \cite{HKV09}. Notice that $\mathcal T$ is invertible and $\mathcal T^{-1}$ is explicitly given as
\begin{align*}
    \mathcal T^{-1}(\mathfrak u,\mathfrak w)=\left(
    \begin{array}{cc}
        \mathcal I_\U&0\\
        -\varepsilon^{-1}\mathcal E\mathcal S&\varepsilon^{-1}\mathcal I_\W
    \end{array}
    \right)\left(
    \begin{array}{c}
        \mathfrak u\\
        \mathfrak w
    \end{array}
    \right)=\left(
    \begin{array}{c}
        \mathfrak u\\
        (\mathfrak w-\mathcal E\mathcal S \mathfrak u)/\varepsilon
    \end{array}
    \right)
\end{align*}
for all $\mathfrak z=(\mathfrak u,\mathfrak w)\in\Z$. In fact, we have
\begin{align*}
    \mathcal T^{-1}(\mathcal T(u,w))=\left(
    \begin{array}{cc}
        \mathcal I_\U&0\\
        -\varepsilon^{-1}\mathcal E\mathcal S&\varepsilon^{-1}\mathcal I_\W
    \end{array}
    \right)\left(
    \begin{array}{c}
        u\\
        \varepsilon  w+\mathcal E\mathcal S  u
    \end{array}
    \right)=\left(
    \begin{array}{c}
        u\\
        w
    \end{array}
    \right)
\end{align*}
for all $z=(u,w)\in\Z$ and viceversa. Furthermore, we find
\begin{equation}
    \label{Section5.2:Tinvstar}
    \mathcal T^{-\star}=\left(
    \begin{array}{cc}
        \mathcal I_\U& -\varepsilon^{-1}\mathcal S^\star\mathcal E^\star \\
        0&\varepsilon^{-1}\mathcal I_\W
    \end{array}
    \right):\Z\to\Z,
\end{equation}
with $\mathcal T^{-\star}=(\mathcal T^\star)^{-1}$. Let $\hydQ=\ydQ-\hat y\in L^2(0,T;H)$ and $\hydT=\ydT-\hat y(T)\in H$. With this we introduce the reduced cost functional $\hat J: \Z \to \mathbb{R}$ by
\begin{align*}
    \hat J(z)=J(\hat y+\mathcal Su,u,w)&=\frac{\sigma_1}{2}\int_0^T{\|(\mathcal Su)(t)-\hydQ(t)\|}_H^2\,\mathrm dt+\frac{\sigma_2}{2}\,{\|(\mathcal Su)(T)-\hydT\|}_H^2+\frac{\sigma}{2}\,{\|u-\un\|}_\U^2\\
    &\quad +\frac{\sigma_w}{2}\,{\|w\|}_\W^2
\end{align*}
for all $z = (u,w) \in \Z$. Then  \eqref{Section5.2:Pw} is equivalent to the following reduced problem
\begin{equation}
    \label{Section5.2:Pwhat}
    \tag{$\mathbf{\hat P}_\varepsilon$}
    \min\hat J(z)\quad\text{s.t.}\quad z\in\Zad.
\end{equation}
Applying standard arguments \cite{Lio71} one can prove that there exists a unique optimal solution $\bar z=(\bar u,\bar w)\in\Zad$ to \eqref{Section5.2:Pwhat}, assuming that the admissible set $\Zad$ is non-empty. The uniqueness follows from the strict convexity properties of the reduced cost functional. 

\subsection{Primal-dual active set strategy}
\label{Section:5.2.2}

Optimality conditions for problem \eqref{Section5.2:Pwhat} can be derived following the optimization approach presented in Section~\ref{SIAM-Book:Section4.1.3}.
\begin{theorem}
\label{Section5.2:TheoremFOC}
Let Assumptions~{\rm\ref{A1}} and {\rm\ref{A9}} hold. Suppose that the feasible set $\Zad$ is nonempty and that $\bar z=(\bar u,\bar w)\in\Zad$ is the solution to \eqref{Section5.2:Pwhat} with associated optimal state $\bar y=\hat y+\mathcal S\bar u$. Then there exist unique Lagrange multipliers $\bar p\in\Y$ and $\bar\nu\in\W$, $\bar\mu=(\bar\mu_i)_{1\le i\le\mathsf m}\in\U$ satisfying the dual equation
\begin{equation}
    \label{Section5.2:DualEquation}
    \begin{aligned}
        -\frac{\mathrm d}{\mathrm dt}\,{\langle \bar p(t),\varphi\rangle}_H+a(t;\varphi,\bar p(t))+{\langle(\mathcal E^\star\bar\nu)(t),\varphi\rangle}_H&=\sigma_1\,{\langle(\ydQ-\bar y)(t),\varphi\rangle}_H~\forall\varphi\in V\text{ a.e. in }[0,T),\\
        \bar p(T)&=\sigma_2\big(\ydT-\bar y(T)\big)
    \end{aligned}
\end{equation}
and the optimality system
\begin{equation}
    \label{Section5.2:OptCond}
    \sigma(\bar u-\un)-\mathcal B'\bar p+\bar\mu=0\text{ in }\U,\quad\sigma_w \bar w+\varepsilon\bar\nu=0\text{ in }\W.
\end{equation}
Moreover,
\begin{subequations}
    \label{Section5.2:NCP}
    \begin{align}
        \label{Section5.2:NCP-1}
        \bar\nu&=\max\big\{0,\bar\nu+\eta(\mathcal E\bar y+\varepsilon\bar w-\yb)\big\}+\min\big\{0,\bar\nu+\eta(\mathcal E\bar y+\varepsilon\bar w-\ya)\big\},\\
        \label{Section5.2:NCP-2}
        \bar\mu_i&=\max\big\{0,\bar\mu_i+\zeta_i(\bar u-u_{bi})\big\}+\min\big\{0,\bar\mu_i+\zeta_i(\bar u_i-u_{ai})\big\}
    \end{align}
\end{subequations}
for $i=1\ldots,\mathsf m$, and for arbitrarily chosen $\eta,\zeta_i>0$, where the max- and min-operations are interpreted componentwise in the pointwise everywhere sense.
\end{theorem}

\begin{remark}
    \label{Section5.2:Remark:phat}
    \rm
    \begin{enumerate}
        \item [1)] Note that \eqref{Section5.2:NCP-1} is a nonlinear complementarity problem (NCP) function based reformulation of the complementarity system
        \begin{align*}
            \bar\nu_a&\ge0,&\ya-\mathcal E\bar y-\varepsilon\bar w&\le0,&{\langle\bar\nu_\mathsf a,\ya-\mathcal E\bar y-\varepsilon\bar w\rangle}_\W&=0,\\
            \bar\nu_b&\ge0,&\mathcal E\bar y+\varepsilon\bar w-\yb&\le0,&{\langle\bar\nu_\mathsf b,\mathcal E\bar y+\varepsilon\bar w-\yb\rangle}_\W&=0
        \end{align*}
        with $\bar\nu=\bar\nu_b-\bar\nu_a\in\W$. Analogously, \eqref{Section5.2:NCP-2} is a NCP function based reformulation of the complementarity system
        \begin{align*}
            \bar\mu_a&\ge0,&u_a-\bar u&\le0,&{\langle \bar\mu_a,u_a-\bar u\rangle}_\U&=0,\\
            \bar\mu_b&\ge0,&\bar u-u_b&\le0,&{\langle \bar\mu_b,\bar u-u_b\rangle}_\U&=0
        \end{align*}
        with $\bar\mu=\bar\mu_b-\bar\mu_a\in\U$.
        \item [2)] Analogously to Lemma~\ref{Lemma:HI-31} and Remark~\ref{Remark:HI-105} we split the adjoint variable $p$ into one part depending on the fixed desired states and into two other parts, which depend linearly on the control variable and on the multiplier $\nu$. Let $\hat p\in\Y$ be the solution to \eqref{GVLuminy:Eq4.3.15}. Further, we define the linear, bounded operators $\mathcal A_1:\U\to\Y$ and $\mathcal A_2:\W\to\Y$ as follows: for given $u\in\U$ the function $p=\mathcal A_1u$ is the unique solution to \eqref{GVLuminy:Eq4.3.12} and for given $\nu\in\W$ the function $p=\mathcal A_2\nu$ uniquely solves
        \begin{align*}
            -\frac{\mathrm d}{\mathrm dt}\,{\langle p(t),\varphi\rangle}_H+a(t;\varphi,p(t))&=-{\langle(\mathcal E^\star\nu)(t),\varphi\rangle}_H&&\text{for all }\varphi\in V\text{ a.e. in }[0,T),\\
            p(T)&=0&&\text{in }H.
        \end{align*}
        In particular, the solution $\bar p$ to \eqref{Section5.2:DualEquation} is given by $\bar p=\hat p+\mathcal A_1\bar u+\mathcal A_2\bar\nu$.\hfill$\blacksquare$
    \end{enumerate}
\end{remark}

It follows from Theorem~\ref{Section5.2:TheoremFOC} that the first-order conditions for \eqref{Section5.2:Pwhat} can be equivalently written as the nonsmooth nonlinear system
\begin{subequations}
    \label{Section5.2:SemiSmoothSystem}
    \begin{align}
        \label{Section5.2:SemiSmoothSystem-a}
        &\sigma\big(\bar u-\un\big)-\mathcal B'\bar p+\bar\mu=0&&\text{in }\U,\\
        \label{Section5.2:SemiSmoothSystem-b}
        &\sigma_w \bar w+\varepsilon\bar\nu=0&&\text{in }\W,\\
        \label{Section5.2:SemiSmoothSystem-c}
        &\bar\nu=\max\big\{0,\bar\nu+\eta(\mathcal E\bar y+\varepsilon\bar w-\yb)\big\}+\min\big\{0,\bar\nu+\eta(\mathcal E\bar y+\varepsilon\bar w-\ya)\big\} \\
        \label{Section5.2:SemiSmoothSystem-d}
        &\bar\mu_i=\max\big\{0,\bar\mu_i+\zeta_i(\bar u_i-u_{bi})\big\}+\min\big\{0,\bar\mu_i+\zeta_i(\bar u_i-u_{ai})\big\}&&\text{for }i=1,\ldots,\mathsf m
    \end{align}
\end{subequations}
with the unknowns $\bar u$, $\bar w$, $\bar\nu$ and $\bar\mu$. As presented for example in \cite{MV18}, one can apply the primal-dual active set strategy (PDASS) in order to solve \eqref{Section5.2:SemiSmoothSystem}, which is equivalent to a semismooth Newton method \cite{HIK03} and therefore super-linearly convergent. We choose $\eta=\sigma_w/\varepsilon^2>0$ and $\zeta_i=\sigma>0$ for all $i=1,\ldots,\mathsf m$. We indicate with $z^k=(u^k,w^k)\in\Z$ the current semi-smooth Newton iterate for $k\in\mathbb N_0$ and we set $y^0=\hat y+\mathcal Su^0$, $p^0=\hat p+\mathcal A_1u^0-\sigma_w\mathcal A_2w^0/\varepsilon$. Moreover, from \eqref{Section5.2:SemiSmoothSystem-a}-\eqref{Section5.2:SemiSmoothSystem-b} we have
\begin{align*}
    y^k&=\hat y+\mathcal Su^k, &\nu^k&=-\frac{\sigma_w}{\varepsilon}\,w^k,\\
    p^k&=\hat p+\mathcal A_1u^k+\mathcal A_2\nu^k,& \mu^k&=\mathcal B'p^k-\sigma \big(u^k-\un\big).
\end{align*}
Defining the associated active sets
\begin{subequations}
    \label{Section5.2:AS}
    \begin{equation}
        \label{Section5.2:ActiveSets}
        \begin{aligned}
            \AaiU(z^k)&=\big\{t\in[0,T]\,\big|\,\mu_i^k+\sigma(u_i^k-u_a)<0\text{ a.e.}\big\},\quad i=1,\ldots,\mathsf m,\\
            \AbiU(z^k)&=\big\{t\in[0,T]\,\big|\,\mu_i^k+\sigma(u_i^k-u_b)>0\text{ a.e.}\big\},\quad i=1,\ldots,\mathsf m,\\
            \AaW(z^k)&=\bigg\{(t,\bx)\in Q\,\big|\,\nu^k+\frac{\sigma_w}{\varepsilon^2}\big(y^k+\varepsilon w^k-\ya\big)<0\text{ a.e.}\bigg\},\\
            \AbW(z^k)&=\bigg\{(t,\bx)\in Q\,\big|\,\nu^k+\frac{\sigma_w}{\varepsilon^2}\big(y^k+\varepsilon w^k-\yb\big)>0\text{ a.e.}\bigg\}
        \end{aligned}
    \end{equation}
    and associated inactive sets
    \begin{equation}
        \label{Section5.2:InactiveSets}
        \begin{aligned}
            \IiU(z^k)&=[0,T]\setminus\big(\AaiU(z^k)\cup\AbiU(z^k)\big)\quad\text{for }i=1,\ldots,\mathsf m,\\
            \IW(z^k)&=Q\setminus\big(\AaW(z^k)\cup\AbW(z^k)\big),
        \end{aligned}
    \end{equation}
\end{subequations}
it holds that the new state $y^{k+1}$ and the new adjoint $p^{k+1}$ are given by the two coupled problems
\begin{align*}
    &\frac{\mathrm d}{\mathrm dt} {\langle y^{k+1}(t),\varphi \rangle}_H+a(t;y^{k+1}(t),\varphi)-\sum_{i=1}^{\mathsf m}\chi_{\IiU(z^k)}(t)\frac{1}{\sigma}\int_\Gamma b_ip^{k+1}(t)\mathrm d\boldsymbol{\tilde s}\,\int_\Gamma b_i\varphi\mathrm d\bs\\
    &={\langle\mathcal F(t),\varphi\rangle}_{V',V}+\sum_{i=1}^{\mathsf m}\big(\chi_{\AaiU(z^k)}(t)u_{ai}(t)+\chi_{\AbiU(z^k)}(t)u_{bi}(t)+\chi_{\IiU(z^k)}(t)\uni \big)\int_\Gamma b_i\varphi\mathrm d\bs\\
    &\hspace{90mm}\forall\varphi\in V\text{ a.e. in }(0,T],\\
    &y^{k+1}(0)=y_\circ.
\end{align*}
and
\begin{align*}
    &-\frac{\mathrm d}{\mathrm dt}\,{\langle p^{k+1}(t),\varphi\rangle}_H+a(t;\varphi,p^{k+1}(t))+\sigma_1\,{\langle y^{k+1}(t),\varphi\rangle}_H&\\
    &\quad\qquad+\frac{\sigma_w}{\varepsilon^2}\,\left\langle y^{k+1}(t)\big(\chi_{\AaW(z^k)}(t)+\chi_{\AbW(z^k)}(t)\big),\varphi\right\rangle_H\\
    &\qquad=\sigma_1\,{\langle \ydQ(t),\varphi\rangle}_H+\frac{\sigma_w}{\varepsilon^2}\,\left\langle \ya(t)\chi_{\AaW(z^k)}(t)+\yb(t)\chi_{\AbW(z^k)}(t),\varphi\right\rangle_H ,\\
    &\hspace{60mm}\forall\varphi\in V\text{ a.e. in }[0,T),\\
    &\quad p^{k+1}(T)=\sigma_2\big(\ydT-y^{k+1}(T)\big),
\end{align*}
respectively, which can be expressed as
\begin{equation}
    \label{Section5.2:OpSystem}
    \left(
    \begin{array}{cc}
        \mathcal A_{11}^k&\mathcal A_{12}^k\\[1ex]
        \mathcal A_{21}^k&\mathcal A_{22}^k
    \end{array}
    \right)\left(
    \begin{array}{c}
        y^{k+1}\\[1ex]
        p^{k+1}
    \end{array}
    \right)=\left(
    \begin{array}{c}
        \mathcal Q_1(z^k;y_\circ,u_a,u_b,\un,b_i,\sigma)\\[1ex]
        \mathcal Q_2(z^k;\ya,\yb,\ydQ,\ydT,\varepsilon,\sigma_w,\sigma_1,\sigma_2)
    \end{array}
    \right).
\end{equation}
We have $\mathcal A_{11}^k=\mathcal H+\tilde{\mathcal A}_{11}^k$ and $\mathcal A_{22}^k=\mathcal H^\star+\tilde{\mathcal A}_{22}^k$, where the $k$-independent operator $\mathcal H:W(0,T)\to L^2(0,T,V')$ is defined as
\begin{align*}
{\langle\mathcal Hy,\varphi\rangle}_{L^2(0,T;V'),L^2(0,T;V)}=\int_0^T{\langle y_t(t),\varphi(t)\rangle}_{V',V}+a(t;y(t),\varphi(t))\,\mathrm dt
\end{align*}
for $y\in W(0,T)$ and $\varphi\in L^2(0,T;V)$. Furthermore, the new control variable $z^{k+1}=(u^{k+1},w^{k+1})$ is given by the linear system
\begin{equation}
    \label{Section5.2:PDASS-System}
    \begin{aligned}
        \int_\Gamma b_ip^{k+1}\mathrm d \bs-\sigma (u_i^{k+1}-\uni)&=0 &&\text{in }\IiU(z^k),&&i=1,\ldots,\mathsf m,\\
        u_i^{k+1}&=u_a&&\text{in }\AaiU(z^k),&&i=1,\ldots,\mathsf m,\\
        u_i^{k+1}&=u_b&&\text{in }\AbiU(z^k),&&i=1,\ldots,\mathsf m,\\
        w^{k+1}&=0&&\text{in }\IW(z^k),\\
        y^{k+1}+\varepsilon\,w^{k+1}&=\ya&&\text{in }\AaW(z^k),\\
        y^{k+1}+\varepsilon\,w^{k+1}&=\yb&&\text{in }\AbW(z^k).
    \end{aligned}
\end{equation}
We resume the previous strategy in Algorithm~\ref{Section5.2:Alg:PDASS}.

\bigskip
\hrule
\vspace{-3.5mm}
\begin{algorithm}[(PDASS method for \eqref{Section5.2:Pw})]
    \label{Section5.2:Alg:PDASS}
    \vspace{-3mm}
    \hrule
    \vspace{0.5mm}
    \begin{algorithmic}[1]
        \STATE Choose starting value $z^0=(u^0,w^0)\in\Z$; set $k=0$ and {\tt flag}~=~{\tt false};
        \STATE Determine $y^0=\hat y+\mathcal Su^0$ and $p^0=\hat p+\mathcal A_1u^0-\sigma_w\mathcal A_2w^0/\varepsilon$;
        \STATE Get $\AaiU(z^0)$, $\AbiU(z^0)$, $\IiU(z^0)$, $i=1,\ldots,\mathsf m$, and $\AaW(z^0)$, $\AbW(z^0)$, $\IW(z^0)$ from \eqref{Section5.2:AS};
        \REPEAT
            \STATE Compute the solution $(y^{k+1},p^{k+1})$ by solving \eqref{Section5.2:OpSystem};
            \STATE Compute $z^{k+1}=(u^{k+1},w^{k+1})\in\Z$ from \eqref{Section5.2:PDASS-System}; 
            \STATE Set $k=k+1$;
            \STATE Get $\AaiU(z^k)$, $\AbiU(z^k)$, $\IiU(z^k)$, $i=1,\ldots,\mathsf m$, and $\AaW(z^k)$, $\AbW(z^k)$, $\IW(z^k)$ from \eqref{Section5.2:AS};
            \IF{$\mathscr A_{a1}^\U(z^k)=\mathscr A_{a1}^\U(z^{k-1})$ {\bf and \rm}~\ldots~{\bf and \rm} $\mathscr A_{a\mathsf m}^\U(z^k)=\mathscr A_{a\mathsf m}^\U(z^{k-1})$}
                \IF{$\mathscr A_{b1}^\U(z^k)=\mathscr A_{b1}^\U(z^{k-1})$ {\bf and \rm}~\ldots~{\bf and \rm } $\mathscr A_{b\mathsf m}^\U(z^k)=\mathscr A_{b\mathsf m}^\U(z^{k-1})$}
                    \IF{$\AaW(z^k)=\AaW(z^{k-1})$ {\bf and \rm} $\AbW(z^k)=\AbW(z^{k-1})$}
                        \STATE {\tt flag}~=~{\tt true};
                    \ENDIF
                \ENDIF
            \ENDIF
        \UNTIL{{\tt flag}~=~{\tt true};}
    \end{algorithmic}
    \hrule
\end{algorithm}

\begin{remark}
    \label{Section5.2:remarkonPDASS}
    \rm 
    For the numerical realization of system \eqref{Section5.2:PDASS-System} one can use for example piecewise linear finite element in space and implicit Euler scheme in time. Unfortunately, the discretized linear system matrix $A\in\mathbb{R}^{2N_xN_t\times 2N_xN_t}$, where $N_x$ corresponds to the number of finite element nodes and $N_t$ to the number of time steps. It is clear that for a fine triangulation of the domain $\Omega$ store the matrix $A$ and solve the linear system may cause problem. The first problem can be avoided thanks to the sparse block structure of the matrix $A$ that allows to use a matrix-free GMRES method (even if we still need to store the finite element blocks). The second problem can be solved applying POD. Moreover, it is still clear that the same problem may occur if the horizon $T$ is large (or even infinite!). In this case, a Model Predictive Control (MPC)\index{Model predictive control, MPC} approach can be utilized; cf. \cite{GP16,RM09}. To be precise, in this context the presence of a terminal target $\ydT$ may lead to ill-posed algorithm, depending on the target itself. Therefore, in the next section we set $\sigma_2=0$ and we briefly introduce MPC, showing also how to combine it with POD in a successful way. Let us mention that there are several contributions in the literature, where MPC and POD are combined; cf., e.g., \cite{AV15,GU14,Mec19,MV19}.\hfill$\blacksquare$
\end{remark}

\subsection{Model Predictive Control}

The basic idea of MPC is to predict, stabilize and optimize a given dynamical system by reconstructing the optimal control $u(t)=\Phi(t,y(t))$ in a feedback form. In order to do that, one has to solve repetitively many open-loop optimal control problems on smaller time horizons $N\Delta T < T$, $N\in\mathbb N$. Then, the open-loop control on the first time step is stored and used to compute the dynamical system trajectory, before solving the next open-loop problem on a shifted time horizon. More details on MPC and its properties can be found in \cite{GP16,RM09}, for instance. Now, how can MPC be applied to \eqref{Section5.2:Pw}? As first, for chosen $0\le t_n<t_n^N\le T$ with $t_n^N=t_n+N\Delta T$ and $y_n\in H$ we consider the following dynamical system on the time horizon $[t_n,t_n^N]$:
\begin{equation}
    \label{Section5.2:WeakDynamicalSystem-n}
    y_t(t)=\mathcal D(t; y(t),u(t))\in V'\text{ a.e. in }(t_n,t_n^N],\quad y(t_n)=y_n\text{ in }H.
\end{equation}
where the time-dependent mapping $\mathcal D(t; \cdot,\cdot):V\times\mathbb{R}^m\to V'$ is defined as
\begin{align*}
    {\langle\mathcal D(t;\phi,\mathfrak u),\varphi\rangle}_{V',V}&=-a(t;\phi,\varphi)+{\langle f(t),\varphi\rangle}_{V',V}+ \sum_{i=1}^\mathsf m \mathfrak u_i(t)\int_\Gamma b_i(\bs)\varphi(\bs)\,\mathrm d\bs
\end{align*}
for $\phi,\varphi\in V$, $\mathfrak u=(\mathfrak u_i)\in\mathbb R^m$ a.e. in $[0,T]$. Furthermore, we consider the inequality constraints
\begin{equation}
    \label{Section5.2:Eq:MPC-10}
    \begin{aligned}
        u_{ai}(t)&\le u_i(t)\le u_{bi}(t),&&i=1,\ldots,\mathsf m\text{ a.e. in }[t_n,t_n^N],\\
        \ya(t,\bx)&\le y(t,\bx)+\varepsilon w(t,\bx)\le \yb(t,\bx)&&\text{a.e. in }Q_n=(t_n,t_n^N)\times\Omega.
    \end{aligned}
\end{equation}
Next we define the function spaces related to $[t_n,t_n^N]$
\begin{align*}
    \U_n=L^2(t_n,t_n^N;\mathbb R^m),\quad\W_n=L^2(t_n,t_n^N;H),\quad\X_n=W(t_n,t_n^N)\times\U_n\times\W_n. 
\end{align*}
Moreover, let the set of admissible solutions be given as
\begin{align*}
    \Xadn=\big\{x=(y,u,w)\in\X\,\big|\,~x\text{ satisfies \eqref{Section5.2:WeakDynamicalSystem-n} and \eqref{Section5.2:Eq:MPC-10}}\big\}.
\end{align*}
Now, the open-loop problem can be adapted by choosing the following cost:
\begin{align*}
    J_n(x)=\frac{\sigma_1}{2}\int_{t_n}^{t_n^N}{\|y(t)-\ydQ(t)\|}_H^2\,\mathrm dt+\frac{\sigma}{2}\,{\|u-\un\|}_{\U_n}^2+\frac{\sigma_w}{2}\,{\|w\|}_{\W_n}^2 
\end{align*}
for $x=(y,u,w)\in\X_n$. Note that we choose $\sigma_2=0$ as briefly justified in Remark~\ref{Section5.2:remarkonPDASS}. The MPC method is summarized in Algorithm~\ref{Section5.2:Alg:MPC}.
\bigskip
\hrule
\vspace{-3.5mm}
\begin{algorithm}[(MPC method)]
    \label{Section5.2:Alg:MPC}
    \vspace{-3mm}
    \hrule
    \vspace{0.5mm}
    \begin{algorithmic}[1]
        \REQUIRE Initial state $y_\circ$, time horizon $N \Delta t$ and regularization parameter $\varepsilon>0$;
        \STATE Set $t_0=0$ and $y_0(0)=y_\circ$;
        \FOR{ $n=0,1,2,\dots$ }
            \STATE Set $t_n^N=t_n+N\Delta t$;
            \STATE Compute the optimal solution $\bar x_n^{\varepsilon}=(\bar y_n^{\varepsilon},\bar u_n^{\varepsilon},\bar w_n^{\varepsilon})$ to
            \begin{equation}
                \label{Section5.2:PhatNt0}
                \tag{${\mathbf P}_n^\varepsilon$}
                \min J_n(x)\quad\text{subject to}\quad x\in\Xadn;
            \end{equation}
            \STATE Define the MPC feedback law $\Phi^N(t;y(t))=\bar u_n^{\varepsilon}(t)$ for $t\in(t_n,t_n+\Delta t]$;
            \STATE Get the MPC state $y_{n+1}$ solving
            \begin{align*}
                (y)_t(t)=\mathcal D(t; y(t),\Phi^N(t;y(t)))\in V'\text{ a.e. in }(t_n,t_n+\Delta t],\quad y(t_n)=y_n\text{ in }H;
            \end{align*}
            \STATE Set $t_{n+1}=t_n+\Delta t$.
        \ENDFOR
    \end{algorithmic}
    \hrule
\end{algorithm}

\bigskip
If $\Phi^N$ is computed by the MPC algorithm, then the state $\bar y^N$ solves \eqref{Section5.2:WeakDynamicalSystem-n} for the closed-loop control $\bar u^N=\Phi^N(\cdot\,;y(\cdot))$ with a given initial condition $y_\circ$. Another advantage of MPC is that we can update the data during the for-loop. For example, suppose that we have a good forecast for the parameters involved in the problem until a certain time $\tilde t\in(0,T)$. Then, we can incorporate a new forecast at that time and update the data for the parameters in the next open-loop solves. This updating strategy can not be done when we solve \eqref{Section5.2:Pw} directly with PDASS: we can not ensure the convergence of PDASS while changing the data during the iteration. Moreover, as already mentioned, when the horizon $T$ is sufficiently large, it is not convenient, or even impossible, to use the PDASS directly, because of large memory consumption and high computational time. With MPC, instead, the horizon length depends on $N$, which is independent from $T$, although the larger $N$ is chosen the better the MPC algorithm converges to the optimal solution in $[0,T]$; cf. \cite{GP16}. 
\subsection{MPC-POD and a-posteriori error estimator}
As already mentioned, one can think to apply POD to \eqref{Section5.2:PDASS-System} improving the computational time. This can be essentially done as shown in Section~\ref{SIAM-Book:Section4.2} and Section~\ref{SIAM-Book:Section4.3}, therefore we do not repeat the procedure. The purpose of this section is to give an idea on how MPC and POD can be successfully combined. As first, note that in Algorithm~\ref{Section5.2:Alg:MPC} we solve \eqref{Section5.2:PhatNt0} in $[t_n,t^N_n]$, therefore we compute an open-loop optimal control that can be used as an indicator of the future closed-loop control, which will be computed as $n$ increases. Therefore, we expect that the snapshots computed through the open-loop optimal control can be close enough to the optimal solution, especially when the dynamics involved are not so complicated. Then, we propose the following method to generate the POD-basis: we solve \eqref{Section5.2:PhatNt0} by using the FE Galerkin scheme for $n=0$. Then, we take the state $\bar y_0^\varepsilon$ and the associated adjoint variable $\bar p^\varepsilon_0$ to build a POD basis of rank $\ell$ as shown in Section~\ref{SIAM-Book:Section4.2} and Section~\ref{SIAM-Book:Section4.3}. Now, \eqref{Section5.2:PhatNt0} is solved for all $n>0$ applying its (from now on) fixed POD Galerkin approximation. As said, we can expect that the POD Galerkin approximation is good enough for some subsequent steps $n$, but clearly at a certain point we need to update the basis. The quality of the approximation can be estimated using an a-posteriori error estimator, which is based on an perturbation argument \cite{DHPY95} and has been already utilized in \cite{TV09}, as the one in Section~\ref{SIAM-Book:Section4.2}. If such a-posteriori error is too big, we solve the current problem \eqref{Section5.2:PhatNt0} by using the FE Galerkin scheme and recompute the POD basis using the obtained optimal FE state and associated adjoint. As done in \cite{GGV15}, this estimate can be generalized for the mixed control-state constraints case. At first, suppose that Assumptions~{\rm\ref{A1}} and {\rm\ref{A9}} hold. Recall that the linear, invertible operator $\mathcal T_\varepsilon$ has been introduced in \eqref{Section5.2:Teps}. In particular, $z=(u,w)$ belongs to $\Zad$ if $\mathfrak z=(\mathfrak u,\mathfrak w)=\mathcal T(z)\in\ZAD$ holds with the closed, bounded and convex subset
\begin{align*}
    \ZAD=\big\{\mathfrak z=(\mathfrak u,\mathfrak w)\in\Z\,\big|\,u_a\le\mathfrak u\le u_b\text{ in }\U\text{ and }\hya \le\mathfrak w\le \hyb \text{ in }\W\big\}\subset\Z.
\end{align*}
Note that -- compared to the definition of the admissible set $\Zad$ -- the set $\ZAD$ does not depend on the solution operator $\mathcal S$ and on the regularization parameter $\varepsilon$. Now, we consider instead of \eqref{Section5.2:Pwhat} the following optimal control problem
\begin{equation}
    \label{Section5.2:PhatT}
    \tag{$\mathbf{\hat{\mathtt P}}^\varepsilon$}
    \min \hat J\big(\mathcal T_\varepsilon^{-1}\mathfrak z\big)\quad\text{s.t.}\quad\mathfrak z=(\mathfrak u,\mathfrak w)\in\ZAD.
\end{equation}
If $\bar z=(\bar u,\bar w)$ solves \eqref{Section5.2:Pwhat}, then $\bar{\mathfrak z}=\mathcal T_\varepsilon(\bar z)$ is the solution to \eqref{Section5.2:PhatT}. Conversely, if $\bar{\mathfrak z}$ solves \eqref{Section5.2:PhatT}, then $\bar z=\mathcal T_\varepsilon^{-1}(\bar{\mathfrak z})$ is the solution to \eqref{Section5.2:Pwhat}. According to \cite{GGV15} we have the following result.

\begin{theorem}
    \label{Section5.2:Result:Aposti}
    Suppose that Assumptions~{\rm \ref{A1}} and {\rm\ref{A9}} hold. Let $\bar z=(\bar u,\bar w)$ be the optimal solution to \eqref{Section5.2:Pwhat}. Then, $\bar{\mathfrak z}=\mathcal T_\varepsilon(\bar z)$ is the solution to \eqref{Section5.2:PhatT}. Moreover, suppose that a point $\mathfrak z^\mathsf{ap}=(\mathfrak u^\mathsf{ap},\mathfrak w^\mathsf{ap})\in\ZAD$ is computed. We set $z^\mathsf{ap}=\mathcal T^{-1}_\varepsilon(\mathfrak z^\mathsf{ap})$, i.e., $z^\mathsf{ap}=(u^\mathsf{ap},w^\mathsf{ap})$ fulfils $u^\mathsf{ap}=\mathfrak u^\mathsf{ap}$ and $w^\mathsf{ap}=\varepsilon^{-1}\,(\mathfrak w^\mathsf{ap}-\mathcal E\mathcal S\mathfrak u^\mathsf{ap})$. Then, there exists a perturbation $\zeta=(\zeta^u,\zeta^w)\in\Z$, which is independent of $\bar z$, so that
    \begin{equation}
        \label{Section5.2:APostError}
        {\|\bar z-z^\mathsf{ap}\|}_\Z\le\frac{1}{\sigma_z}\,{\|\mathcal T_\varepsilon^\star\zeta\|}_\Z\quad\text{with }\sigma_z=\min\{\sigma,\sigma_w\}>0.
    \end{equation}
    where $\mathcal T_\varepsilon^\star$ denotes the adjoint of the operator $\mathcal T_\varepsilon$. 
\end{theorem}

\noindent{\bf\em Proof.}
Since $\mathcal T_\varepsilon$ has a bounded inverse, the first claim follows. The second one can be shown by adapting the proof of Proposition~1 in \cite{GGV15}.\hfill$\Box$

\begin{remark}
    \rm
    The perturbation $\zeta$ can be computed as follows: let $\xi=(\xi^u,\xi^w)\in\Z$ be given as $\xi=\mathcal T^{-\star}\nabla\hat J(\bar z^\mathsf{ap})\in\Z$. Then, $\xi$ solves the linear system $\mathcal T^\star\xi=\nabla\hat J(z^\mathsf{ap})$, i.e.,
    \begin{subequations}
        \label{Section5.2:PertZeta}
        \begin{equation}
            \label{Section5.2:SystemAPostError}
            \left(
            \begin{array}{cc}
                \mathcal I_\U& \mathcal S^\star \mathcal E^\star \\
                0&\varepsilon \mathcal I_\W
            \end{array}
            \right)\left(
            \begin{array}{c}
                \xi^u\\
                \xi^w
            \end{array}\right)
            =\left(
            \begin{array}{c}
                \sigma (u^\mathsf{ap}-\un) - \mathcal{B} p^\mathsf{ap} \\
                \sigma_w w^\mathsf{ap}
                \end{array}
            \right)
        \end{equation}
    where $p^\mathsf{ap}=\hat p+\mathcal A_1 u^\mathsf{ap}$. Note that the first-order optimality conditions for problem \eqref{Section5.2:PhatT} can be written as
    \begin{align*}
        {\langle \xi+\zeta,\mathfrak z-\mathfrak z^\mathsf{ap}\rangle}_\Z\ge 0\quad\text{for all }\mathfrak z\in\ZAD.
    \end{align*}
    Hence, we find
    \begin{equation}
        \label{Section5.2:PertZeta-a}
        \zeta_i^u(t)=\left\{
        \begin{array}{ll}
            -\min\{0,\xi_i^u(t)\}&\text{for }t\in \AaiU(z^\mathsf{ap}),\\[1mm]
            -\max\{0,\xi_i^u(t)\}&\text{for }t\in\AbiU(z^\mathsf{ap}),\\[1mm]
            -\xi^u_i(t)&\text{for }t\in\IiU(z^\mathsf{ap})
        \end{array}
        \right.
    \end{equation}
    for $i=1,\ldots,m$ and
    \begin{equation}
        \label{Section5.2:PertZeta-b}
        \zeta^w(t,\bx)=\left\{
        \begin{array}{ll}
            -\min\{0,\xi^w(t,\bx)\}&\text{for }(t,\bx)\in\AaW(z^\mathsf{ap}),\\[1mm]
            -\max\{0,\xi^w(t,\bx)\}&\text{for }(t,\bx)\in\AbW(z^\mathsf{ap}),\\[1mm]
            -\xi^w(t,\bx)&\text{for }(t,\bx)\in\IW(z^\mathsf{ap}).
            \end{array}\right.
        \end{equation}
    \end{subequations}
\hfill$\blacksquare$
\end{remark}

Finally, the MPC-POD algorithm is resumed in Algorithm~\ref{Section5.2:Alg:MPCPOD}.

\bigskip
\vspace{-1mm}
\hrule
\vspace{-3.5mm}
\begin{algorithm}[(MPC-POD method)]
    \label{Section5.2:Alg:MPCPOD}
    \vspace{-3mm}
    \hrule
    \vspace{0.5mm}
    \begin{algorithmic}[1]
        \REQUIRE Initial state $y_\circ$, time horizon $N \Delta t$, regularization parameter $\varepsilon>0$, rank of POD basis $\ell$ and tolerance $\tau$ for the a-posteriori estimate;
        \STATE Set $t_0=0$, $y_0(0)=y_\circ$, and {\tt flag} = {\tt true};
        \FOR{ $n=0,1,2,\dots$ }
            \STATE Set $t_n^N=t_n+N\Delta t$;
            \IF { {\tt flag} = {\tt true} }
                \STATE Compute the optimal solution $\bar x_n^{\varepsilon}=(\bar y_n^{\varepsilon},\bar u_n^{\varepsilon},\bar w_n^{\varepsilon})$ to \eqref{Section5.2:PhatNt0};
                \STATE Compute a POD basis $\{\psi_i\}_{i=1}^\ell$ of rank $\ell$ using as snapshots $\bar y_n^\varepsilon$ and $\bar p_n^\varepsilon$;
                \STATE Define the MPC feedback law $\Phi^N(t;y(t))=\bar u_n^{\varepsilon}(t)$ for $t\in(t_n,t_n+\Delta t]$;
                \STATE Set {\tt flag} = {\tt false};
            \ELSE
                \STATE Compute the POD suboptimal solution $\bar x_n^{\varepsilon,\ell}=(\bar y_n^{\varepsilon,\ell},\bar u_n^{\varepsilon,\ell},\bar w_n^{\varepsilon,\ell})$ to \eqref{Section5.2:PhatNt0};
                \STATE Define the MPC feedback law $\Phi^N(t;y(t))=\bar u_n^{\varepsilon,\ell}(t)$ for $t\in(t_n,t_n+\Delta t]$;
                \STATE Compute the a-posteriori error estimate $e:= {\|\mathcal T_\varepsilon^\star\zeta\|}_\Z$ from \eqref{Section5.2:PertZeta};
                \IF { $e> \tau$}
                    \STATE Set {\tt flag} = {\tt true};
                \ENDIF
            \ENDIF
            \STATE Set $t_{n+1}=t_n+\Delta t$ and get the MPC state $y_{n+1}$ solving
            \begin{align*}
            (y)_t(t)=\mathcal D(t; y(t),\Phi^N(t;y(t)))\in V'\text{ a.e. in }(t_n,t_{n+1}],\quad y(t_n)=y_n\text{ in }H;
            \end{align*}
        \ENDFOR
    \end{algorithmic}
    \hrule
    \end{algorithm}
\begin{example}
    \rm For this example we want to solve the guiding model problem with slight modification on its parameters and few additional requirements. In this case, we do not consider distributed control, assuming that there exists only a control on the boundary region $\Gamma_1= \{0\}\times[0.75,1.0]$, which is also where we have the inflow for the convection. We assume a different outside temperature 
    \begin{align*}
        \yout= 15+\left(\frac{1}{2}+\cos(\pi t)\right)\left(\frac{1}{2}+\frac{\sin(\pi t)}{2}\right),
    \end{align*}
    see Figure~\ref{Fig:outside_MPC}. 
    \begin{figure}
        \centering
        \includegraphics[height=50mm]{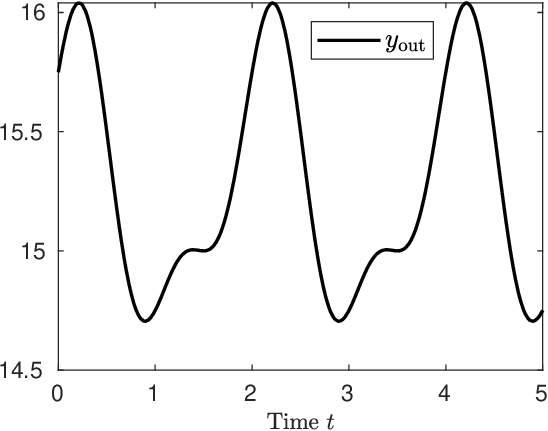}
        \caption{Outside temperature $\yout$.}
        \label{Fig:outside_MPC}
    \end{figure}
    We define also a modified advection field
    \begin{align*}
        \widetilde{\bv}(t,\bx) = (0.1 + 0.1t)\bv(\bx),
    \end{align*}
    where $\bv$ is the advection field introduced in \ref{SIAM-Book:Section1.3}. Furthermore, we choose $\ua = 0$, $\ub = 10^7$ as control constraints and
    \begin{align*}
        \ya(t,\bx) = \min\left(16+\frac{t}{4},18\right), \quad \yb(t,\bx)= 32
    \end{align*}
    as state constraints. The lower bound imposes a progressive rise of the minimum temperature in the room, instead the upper bound fix the value for the maximum temperature. Lastly we treat an economic MPC problem \cite{GP16}, i.e., we choose $\sigma_1=\sigma_2=0$ and $\sigma= \sigma_w= 1.0$. Therefore, the optimization is focused on keeping the temperature $y$ inside the bounds using less control energy as possible. Note that, since we do not have a target, the existence of a non-trivial solution for the optimal control problem is guaranteed only by the presence of the state constraints.
    
    We focus on three related solution approaches: MPC using the full-order model (MPC-FE), MPC using the POD model without updating the basis functions (MPC-POD-NU) and Algorithm~\ref{Section5.2:Alg:MPCPOD} (MPC-POD) with the tolerance $\tau=3.0$. We point out that there is no general recipe to choose the tolerance $\tau$. Heuristically, one should consider the expected range for the optimal controls, according also to the control constraints. In the MPC framework, the possible determination of the turnpike of the closed-loop trajectory will give access to this expected range; cf. \cite{GP16}. This expected range can be possibly determined also by solving the open-loop optimal control problem for a sufficiently large horizon. Of course, the tolerance $\tau$ can be adjusted during the MPC algorithm. In our simulations, the magnitude of the controls for the open-loop problems was always in between 50 and 150, leading to the decision of a fixed $\tau=3.0$, which indicates our will to have a relative error of approximately the $3\%$ between the full-order and reduced-order controls. The MPC open-loop time horizon is set to $N=120$ and $\Delta t = 5/199$.
     
    The computed trajectory with the three MPC algorithms are reported in Figure~\ref{Section5:Fig:MPCtraj}.
    \begin{figure}
    	\begin{subfigure}[b]{0.32\textwidth}
    		\centering
    		\includegraphics[height=40mm]{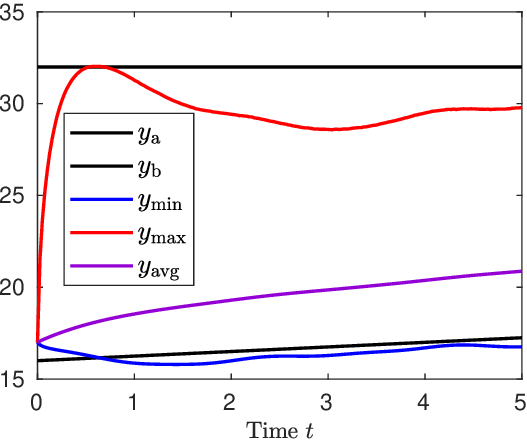}
    		\caption{MPC-FE}
    	\end{subfigure}
    	\begin{subfigure}[b]{0.32\textwidth}
    		\centering
    		\includegraphics[height=40mm]{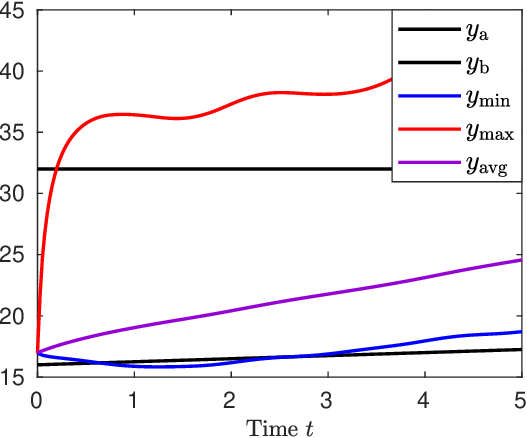}
    		\caption{MPC-POD-NU}
    	\end{subfigure}
    	\begin{subfigure}[b]{0.32\textwidth}
    		\centering
    		\includegraphics[height=40mm]{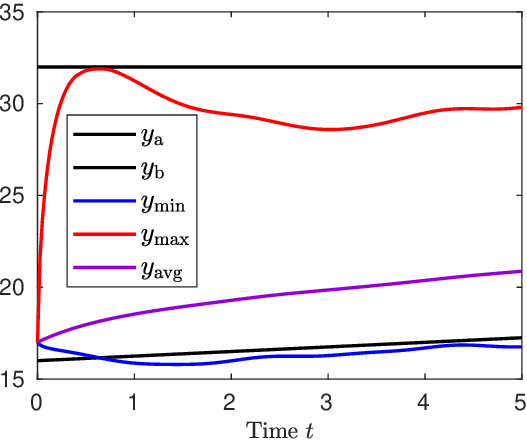}
    		\caption{MPC-POD}
    	\end{subfigure}
    	\caption{MPC trajectory up to $t_{199} = 5$.}
        \label{Section5:Fig:MPCtraj}
    \end{figure}
    As one can see, the minimum temperature in the domain $\Omega$ at each time step (in blue) is not always greater than the constraint $\ya$. This is caused by the virtual control regularization. For this example, the parameter $\varepsilon$ is fixed in fact to 0.001 and therefore allows small violations of the constraints. The behavior of the virtual control approach diminishing the parameter $\varepsilon$ is studied in \cite{KR09,Mec19}. The average temperature (in violet) is anyway maintained inside the constraints' range. Furthermore, Figure~\ref{Section5:Fig:MPCtraj} shows that the MPC-POD-NU method (i.e. without basis updates) is not able to reconstruct the MPC-FE trajectory and even keeping the temperature inside the relaxed virtual constraints. This is confirmed by Table~\ref{Section5:Tab:MPCerrors}, where the relative errors between the MPC-FE trajectory and the MPC-FE feedback control are reported. These errors are respectively defined as
    \begin{align*}
        \text{Rel.Err.}y = \frac{\|y^\text{MPC-FE}-y\|_\W}{\|y^\text{MPC-FE}\|_\W},\quad\text{Rel.Err.}u = \frac{\|u^\text{MPC-FE}-u\|_\U}{\|u^\text{MPC-FE}\|_\U}.
    \end{align*}
    In Table~\ref{Section5:Tab:MPCerrors} are reported also the computational time speed-up and the number of POD basis functions (also after each update for MPC-POD).
    \begin{table}
    	\centering
    	\begin{tabular}{l|cccc}\toprule
    		Method & Rel. Err. $y$ & Rel. Err. $u$ & POD basis $\ell$ & Speed-up \\ \midrule
    		MPC-POD-NU & 0.1178 & 0.6131 & 14 & 28.4 \\
    		MPC-POD & 0.0053 & 0.0004 & $\left\{14,21,33,48,64\right\}$ & 10.9\\\bottomrule
    	\end{tabular}
    	\caption{Relative errors, $\ell$ and speed-up in approximating the MPC-FE solution.} \label{Section5:Tab:MPCerrors}
    \end{table}
    As one can notice, the risk is that $\ell$ increases at each update, so that the MPC-POD method becomes slower. Furthermore, as Algorithm~\ref{Section5.2:Alg:MPCPOD} proceeds, $\ell$ can potentially become greater than the number of degrees of freedom of the full-order model. To avoid this, a more elaborated basis update strategy can be developed \cite{Mec19}. 
    
    The computed feedback controls are shown in Figure~\ref{Section5:Fig:MPCoptcont} together with the MPC-POD-NU and MPC-POD error approximation of the MPC-FE control. Clearly, the MPC-POD-NU method can not reconstruct the full-order feedback control without updating the POD basis, due to the fact that the involved dynamics change according to the time-dependent parameters. A change of these parameters after the training horizon used for the POD method will cause an important lack in approximation. This lack in approximation can even be present since the begin, making necessary an immediate correction, as in this case ($\ell = 21$). The subsequent basis updates happening respectively at time steps $n=\left\{4,17,42\right\}$ ($\ell = \left\{33,48,64\right\}$) show anyway that the update of the POD model is important also during the MPC algorithm.\hfill$\blacklozenge$
    \begin{figure}
        \centering
        \begin{subfigure}[b]{0.4\textwidth}
    	   \centering
    	   \includegraphics[height=40mm]{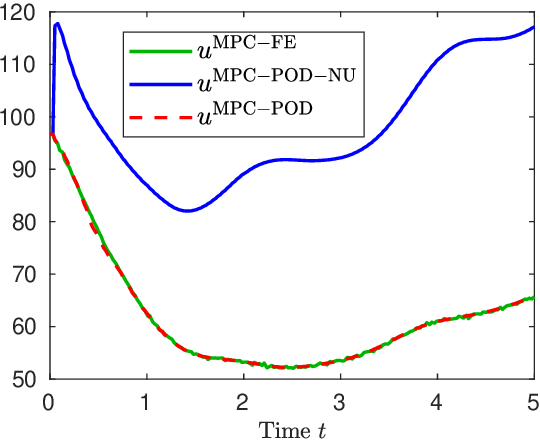}
    	   \caption{MPC feedback controls}
        \end{subfigure}
        \begin{subfigure}[b]{0.4\textwidth}
    	   \centering
    	   \includegraphics[height=40mm]{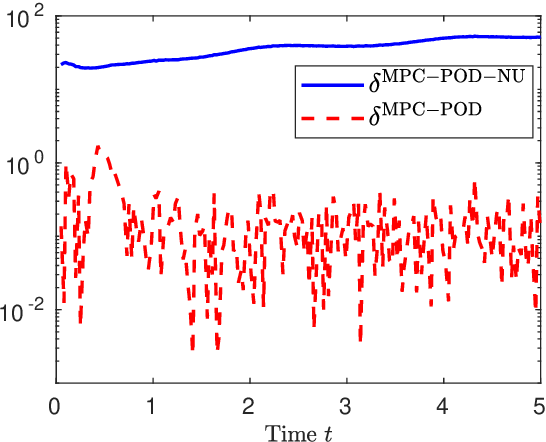}
    	   \caption{$\delta:= \|u^\text{MPC-FE}-u\|_\U$}
        \end{subfigure}
        \caption{Feedback control and POD approximation error up to $t_{199}=5$.}
        \label{Section5:Fig:MPCoptcont}
    \end{figure}
\end{example}

\section{Multiobjective PDE-constrained optimization}
\label{Section:MOP}

In real applications, optimization problems are often described by introducing several objective functions conflicting with each other. This leads to {\em multiobjective}\index{Problem!multiobjective optimization} or {\em multicriterial} optimization problems\index{Problem!multicriterial optimization}; cf. \cite{Ehr05,Mie98,Sta88}. One prominent example is given by an energy efficient heating, ventilation and air-conditioning (HVAC) operation of a building with conflicting objectives such as minimal energy consumption and maximal comfort; cf. \cite{FHC06,KTX11}. Finding the optimal control that represents a good compromise is the main issue in these problems. For that reason the concept of Pareto optimal or efficient points is developed. In contrast to scalar-valued optimization problems, the computation of a set of Pareto optimal points is required. Consequently, many scalar-valued constrained optimization problems have to be solved.

In this section we follow \cite{Ban17,BBV16,BBV17} and apply the reference point method \cite{Wie80,Wie86,RBWSZ09} in order to transform a bicriterial optimal control problem into a sequence of scalar-valued optimal control problems and solve them using well-known optimal control techniques; see \cite{Tro10}. Let us mention that preliminary results combining reduced-order modeling and multiobjective PDE-constrained optimization have recently been derived; cf. \cite{ITV16,IUV17,POD16}. 

In this section we will utilize the following notation: If $x^1,x^2 \in \mathbb R^2$ are two vectors, we write $x^1 \le x^2$ if $x^1_i \le x^2_i$ for $i=1,2$, and analogously for $x^1\ge x^2$. For $x\in\mathbb R^2$ we set $\mathbb R^2_{\ge x}=\{y\in\mathbb R^2\,|\,y\ge x\}$ and $\mathbb R^2_{\le x}=\{y\in\mathbb R^2\,|\,y\le x\}$. For brevity, we define $\mathbb R^2_{\le} := \mathbb R^2_{\le 0}$ and $\mathbb R^2_{\ge} := \mathbb R^2_{\ge 0}$. 

\subsection{Problem formulation}
\label{Section:MOP-1}

The Hilbert space $\U$ is either $\mathbb R^\mathsf m$ or $L^2(\mathscr D;\mathbb R^\mathsf m)$ with an open subset $\mathscr D\subset\mathbb R^{\mathsf m_\circ}$. Let Assumption~\ref{A9} hold. Suppose that $\Ha$ is a real Hilbert space and $\mathcal C:\Y\to\Ha$ a linear, bounded operator. Then we consider the following bicriterial optimal control problem:
\begin{subequations}
    \label{MOP-1}
    \begin{equation}
        \label{MOP-1a}
        \min J(y,u)=\left(
        \begin{array}{c}
            J_1(y,u)\\
            J_2(y,u)
        \end{array}
        \right)=\frac{1}{2}\left(
        \begin{array}{c}
            {\|\mathcal Cy-\yd\|}_\Ha^2\\[1mm]
            {\|u-\un\|}_\U^2
        \end{array}
    \right),
\end{equation}
subject to $(y,u)\in\X=\Y\times\U$ solves
\begin{equation}
    \label{MOP-1b}
    \begin{aligned}
    \frac{\mathrm d}{\mathrm dt} {\langle y(t),\varphi \rangle}_H+a(t;y(t),\varphi)&={\langle (f+\mathcal Bu)(t),\varphi\rangle}_{V',V}&&\text{for all }\varphi\in V\text{ a.e. in }(\tc,\te],\\
    {\langle y(\tc),\varphi\rangle}_H&={\langle y_\circ,\varphi\rangle}_H&&\text{for all }\varphi\in H,
    \end{aligned}
\end{equation}
and
\begin{equation}
\label{MOP-1c}
u\in\Uad=\big\{\tilde u\in\U\,\big|\,u_a\le\tilde u\le u_b\text{ in }\U\big\},
\end{equation}
\end{subequations}
where we additionally assume $u_a, u_b \in L^{\infty}(\mathscr D;\mathbb R^\mathsf m)$ in the case of $U = L^2(\mathscr D;\mathbb R^\mathsf m)$. \\
Recall that the solution to \eqref{MOP-1b} can be expressed as $y=\hat y+\mathcal Su$. We set $\mathcal G=\mathcal C\mathcal S$, which implies that $\mathcal G\in\mathscr L(\U,\Ha)$ holds.

Then we define the reduced objective $\hat J:\U\to\mathbb R^2$ by
\begin{align*}
    \hat J(u)=\left(
    \begin{array}{c}
        \hat J_1(u)\\
        \hat J_2(u)
    \end{array}
    \right)=J(\hat y+\mathcal Su,u)\quad \text{for }u\in\U
\end{align*}
and consider the reduced bicriterial optimal control problem\index{Problem!bicriterial optimization}
\begin{equation}
    \label{BP-1}
    \tag{$\mathbf{\widehat{BP}}$}
    \min \hat J(u)\quad\text{s.t.}\quad u\in\Uad.
\end{equation}

Problem \eqref{BP-1} involves the minimization of a vector-valued objective. This is done by using the concept of order relation and Pareto optimality; cf. \cite{Ehr05}, for instance.

\begin{definition}
    \label{Definition:BP-1}
    The point $\bar u\in\Uad$ is called {\em Pareto optimal} for \eqref{BP-1} provided there is no other control $u\in\Uad\setminus\{\bar u\}$ with $\hat J_i(u)\le\hat J_i(\bar u)$, $i=1,2$, and $J_j(u)<\hat J_j(\bar u)$ for at least one index $j\in\{1,2\}$.
\end{definition} 

The goal of multiobjective optimization is to compute (a numerical approximation of) the {\em Pareto set}\index{Pareto set}\index{Set!Pareto}
\begin{align*}
    \Ps=\big\{u\in\Uad\,\big|\,u\text{ is Pareto optimal}\big\}\subset\U
\end{align*}
and the {\em Pareto front}\index{Pareto front}
\begin{align*}
    \Pf=\big\{\hat J(u)\,\big|\,u\in\Ps\big\}\subset\mathbb R^2.
\end{align*}
In the next section we will introduce the Euclidean reference point method as a tool to achieve this goal. For the theoretical investigation we need the following properties of the functions $\hat J_1$ and $\hat J_2$, which are proved in \cite{BBV17}.

\begin{lemma}
    \label{Lemma:BP-1}
    \begin{enumerate}
        \item [\rm 1)] Both $\hat J_1$ and $\hat J_2$ are bounded from below.
        \item [\rm 2)] Both $\hat J_1$ and $\hat J_2$ are twice continuously Fr\'{e}chet-differentiable.
        \item [\rm 3)] $\hat J_1$ is convex and $\hat J_2$ is strictly convex.
        \item [\rm 4)] $\hat J_1$ is even strictly convex provided $\mathcal G=\mathcal C\mathcal S$ is injective.
    \end{enumerate}
\end{lemma}

\begin{remark}
    \label{Remark:BP-1}
    \rm Since the functions $\hat J_1$ and $\hat J_2$ are bounded from below, convex and continuous, and $\Uad$ is bounded, the ideal vector $\yid=(y^\mathsf{id}_1,y^\mathsf{id}_2)\in\mathbb R^2$ given as
    \begin{align*}
        y^\mathsf{id}_i=\min\big\{\hat J_i(u)\,\big|\,u\in\Uad\big\}
    \end{align*}
    is well-defined.\hfill$\blacksquare$
\end{remark}

\subsection{The Euclidean reference point method}
\label{Section:MOP-2}

One popular methodology to solve multiobjective optimization problems is the so-called scalarization method, where the multiobjective problem is transformed into a series of scalar optimization problems. These can then be solved by using well-known techniques for scalar optimization. Here we follow \cite{BBV17,POD16,RBWSZ09}: Given a {\em reference point}\index{Reference point}
\begin{align*}
    z=(z_1,z_2)\in\Pf+\mathbb R^2_\le=\big\{z+x\,\big|\,z\in\Pf\text{ and }x\in\mathbb R^2_\le\big\}
\end{align*}
we introduce the (Euclidean) {\em distance function}\index{Distance function} $F_z:\U\to\mathbb R^2$ by
\begin{align*}
    F_z(u)=\frac{1}{2}\,{\|\hat J(u)-z\|}_2^2 = \frac{1}{2}\,\big(\hat J_1(u)-z_1\big)^2+\frac{1}{2}\,\big(\hat J_2(u)-z_2\big)^2\quad \text{for }u\in\U.
\end{align*}

\begin{lemma}
    \label{Lemma:BP-2}
    If $z\le\yid$ holds, $F_z$ is strictly convex.
\end{lemma}

\noindent{\bf\em Proof.} The distance functions has the form 
\begin{align*}
    F_z(u) = G_1(\hat J_1(u)) + G_2(\hat J_2(u))
\end{align*}
with $G_i(\xi) = \tfrac{1}{2}(\xi-z_i)^2$ being strictly monotonically increasing on $[z_i,\infty)$. Since $\hat J_1$ is convex and $\hat J_2$ strictly convex by Lemma \ref{Lemma:BP-1}, $F_z$ itself is strictly convex because $\hat J(u) \ge y^\mathsf{id} \ge z$ for all $u \in \Uad$. \hfill$\Box$

\medskip
To find a point that approximates $z$ as good as possible we consider\index{Problem!reference point}
\begin{equation}
    \tag{$\mathbf{\widehat{BP}}_z$}
    \label{BP-2}
    \min F_z(u)\quad\text{s.t.}\quad u\in\Uad.
\end{equation}
Notice that \eqref{BP-2} is a scalar-valued optimization problem which is convex provided $z\le\yid$ holds true. The next result is proved in \cite{Ban17}.

\begin{proposition}
    \label{Proposition:BP-1}
    Let $z\in\Pf+\mathbb R^2_\le$. Then \eqref{BP-2} has a unique solution $\bar u_z=(\bar u_{zj})_{1\le j\le\mathsf m}$, which is Pareto optimal for \eqref{BP-1}.
\end{proposition}

The solution $\bar u_z$ can be characterized by the following optimality conditions; cf. \cite{BBV17}.

\begin{proposition}
    \label{Proposition:BP-2}
    Let $z\in\Pf+\mathbb R^2_\le$ and $\mathcal G$ be injective. If the point $\bar u_z$ is optimal for \eqref{BP-2}, then the first-order necessary optimality condition
    \begin{equation}
        \label{BP-3}
        \big\langle(\hat J_1(\bar u_z)-z_1)\mathcal G^\star(\mathcal G\bar u_z + \hat{y} - \yd)+(\hat J_2(\bar u_z)-z_2)(\bar u_z-\un),u-\bar u_z\big\rangle_\U\ge 0
    \end{equation}
    holds for all $u\in\Uad$. Moreover, if $z_1\le\hat J_1(\bar u_z)$ and $z_2<\hat J_2(\bar u_z)$ we have that
    \begin{align*}
        {\langle F_z''(\bar u_z)u^\delta,u^\delta\rangle}_{\U',\U}\ge\kappa_z\,{\|u^\delta\|}_\U^2\quad\forall u^\delta\in\U,
    \end{align*}
    where the hessian of $F_z''(u):\U\to\U'$ at a point $u\in\U$ is given as
    \begin{align*}
        F_z''(u)u^\delta=\left\langle\big((\hat J_1(u)-z_1)\mathcal G^\star\mathcal G+(\hat J_2(u)-z_2)\mathcal I_\U\big)u^\delta,\cdot\right\rangle_\U\quad \text{for }u^\delta\in\U
    \end{align*}
    and $\kappa_z=\hat J_2(\bar u_z)-z_2>0$. If even $z\le\yid$ is valid, \eqref{BP-3} is a sufficient optimality condition for \eqref{BP-2}.
\end{proposition}

The idea is to solve \eqref{BP-2} for different appropriately chosen values of $z$ to get a uniform approximation of the Pareto front. To be more precise, we want to construct a sequence $\{z^{(k)}\}_{k\in\mathbb N}\subset \mathbb R^2$ of reference points along with optimal controls $\{u^{(k)}\}_{k\in\mathbb N}\subset\Uad$ that solve \eqref{BP-2} with $z = z^{(k)}$ as well as $\{\hat J^{(k)}\}_{k\in\mathbb N} \subset \mathbb{R}^2$ with $\hat J^{(k)}=\hat J(u^{(k)})$, such that a predefined quality criteria for the approximation of the Pareto front is fulfilled. That this is possible in theory is shown in the next proposition, which was proved in \cite[Theorem 3.43, Theorem 3.47]{Ban17}.

\begin{proposition}
    \label{Proposition:BP-AllPoints}
    \begin{enumerate}
        \item [\rm 1)] For any $u \in \Ps$ there is $z \in \Pf+\mathbb R^2_\le$ such that $u$ is the unique solution of \eqref{BP-2}.
        \item [\rm 2)] Let $z^{(1)},z^{(2)} \in \Pf+\mathbb R^2_\le$ be arbitrary and denote by $\bar{u}_{z^{(1)}},\bar{u}_{z^{(2)}} \in \Ps$ the unique solutions of the respective Euclidean reference point problems \eqref{BP-2}. Then it holds
        \begin{equation}
        \label{Equation:UniformDistance}
        \left\Vert J(\bar{u}_{z^{(1)}}) - J(\bar{u}_{z^{(2)}}) \right\Vert \leq \left\Vert z^{(1)} - z^{(2)} \right\Vert.
        \end{equation}
    \end{enumerate}
\end{proposition}

\begin{figure}
    \begin{center}
        \includegraphics[height=50mm]{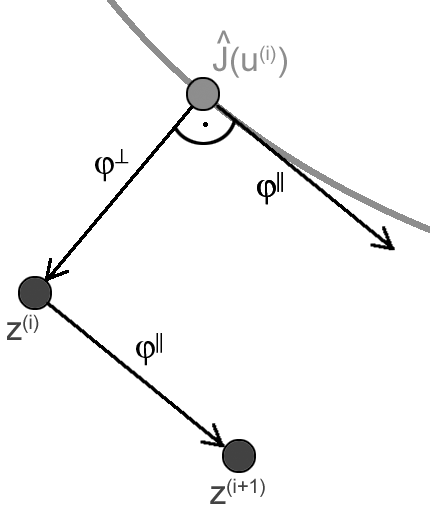}\hspace{10mm}
        \includegraphics[height=50mm]{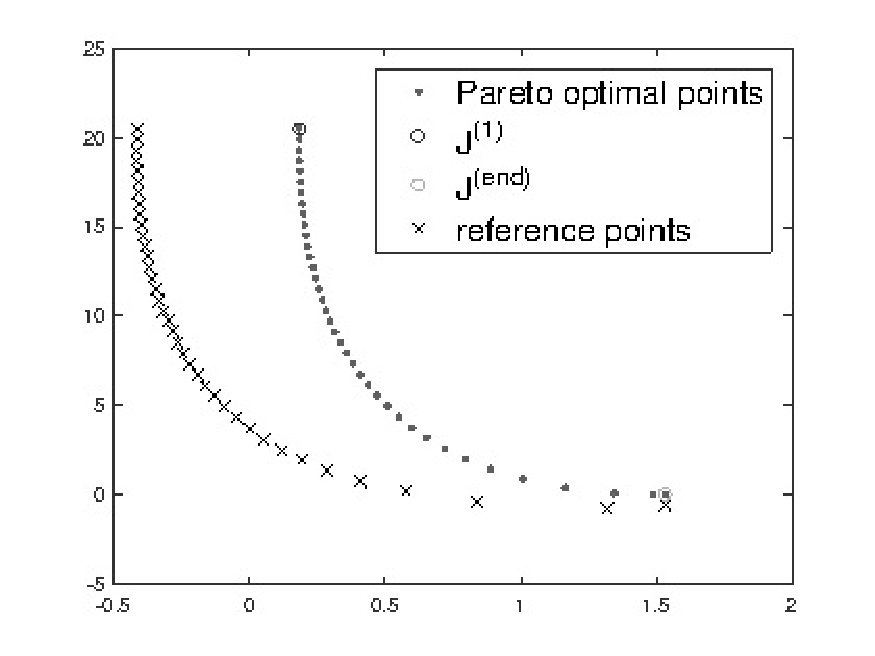}
    \end{center}
    \caption{
    Graphical visualization of the iterative reference point generation (left) and a Pareto front which was computed by Algorithm \ref{alg:rpm_full} (right).}
    \label{fig:RefPointGeneration}
\end{figure}

The first statement of the previous proposition tells us that we can indeed compute every Pareto optimal control $u \in \Ps$ by solving a Euclidean reference point problem \eqref{BP-2}. Moreover, the second statement implies that the distance between two reference points $z^{(1)},z^{(2)}$ is an upper bound of the distance between the corresponding Pareto optimal points $J(\bar{u}_{z^{(1)}}),J(\bar{u}_{z^{(2)}}) \in \Pf$. This enables us to define an iterative sequence of reference points with a guaranteed upper bound on the maximal distance between two neighbouring Pareto optimal points. The exact procedure is the following: We first compute the minimizers of $\hat{J_1}$ and $\hat{J}_2$, which are denoted by $\bar{u}^{(1)}$ and $\bar{u}^{(\text{end})}$. Then the corresponding Pareto optimal points $\hat{J}^{(1)} := \hat{J}(\bar{u}^{(1)})$ and $\hat{J}^{(\text{end})} := \hat{J}(\bar{u}^{(\text{end})})$ mark the start and the end point of our iteration. Note that in our example computing the minimizer of $\hat J_1$ numerically is quite challenging. Therefore, instead of the minimizer of $\hat J_1$, we compute the minimizer of the weighted-sum function $\hat J_1 + \alpha \hat J_2$ for a weight $0 < \alpha \ll 1$, which is on the one hand a good approximation of the minimizer of $\hat J_1$ and on the other hand regular enough to solve it numerically. Starting at the point $\hat{J}^{(1)}$, we generate the first reference point by the formula
\begin{equation}
    \label{Equation:FirstReferencePoint}
    z^{(1)} = \hat{J}^{(1)} + h^{||} \cdot \frac{(\alpha,-1)^T}{\left\Vert (\alpha,-1)^\top \right\Vert} + h^\perp \cdot \frac{(-1,-\alpha)^\top}{\left\Vert (-1,-\alpha)^T \right\Vert}.
\end{equation}
Here, we make use of the fact that the Pareto front at the point $\hat{J}^{(1)}$ is tangential to the vector $\left( \alpha,-1 \right)^T$. Thus, the parameter $h^{||} > 0$ determines how far along the Pareto front the reference point will be chosen. By $h^\perp \geq 0$ we can choose how far the new reference point should be located from the Pareto front. Analogously, one can show that after solving a Euclidean reference point problem \eqref{BP-2} to a reference point $z$, the vectors $\varphi^\perp = z - \hat{J}(\bar{u}^z)$ and $\varphi^{||} = ( -\varphi^\perp_2,\varphi^\perp_1 )^\top$ are perpendicular and tangential to the Pareto front at the point $\hat{J}(\bar{u}^z)$, respectively. This motivates the iterative formula for the reference points
\begin{equation}
    \label{Equation:IterationReferencePoints}
    z^{(i+1)} = z^{(i)} + h^{||} \frac{\varphi^{||}}{\left\Vert \varphi^{||} \right\Vert} + h^\perp \frac{\varphi^\perp}{\left\Vert \varphi^\perp \right\Vert}
\end{equation}
for $i = 1,\ldots$, see Figure \ref{fig:RefPointGeneration}.
Using this formula it is guaranteed that 
\begin{equation*}
    \left\Vert J^{(i+1)} - J^{(i)} \right\Vert \leq h^{||}
\end{equation*}
holds for all $i = 1,\ldots$, which is why we can control the approximation quality by adapting the parameter $h^{||}$. We stop this iteration when the first component of the new reference point is larger than the first component of our end point $\hat{J}^{(\text{end})}$. The complete methodology is summarized in Algorithm \ref{alg:rpm_full}.\index{Method!Euclidean reference point} 

\bigskip
\hrule
\vspace{-3.5mm}
\begin{algorithm}[(Euclidean Reference Point Method)]
    \label{alg:rpm_full}
    \vspace{-3mm}
    \hrule
    \vspace{0.5mm}
	\begin{algorithmic}[1]
		\REQUIRE Parameters $h^{||} > 0$, $h^\perp \geq 0$, weighted-sum parameter $\alpha > 0$.
		\STATE Minimize $\hat{J}_1 + \alpha \hat J_2$ and $\hat J_2$ and save the solutions as $\bar{u}^1$ and $\bar{u}^{(\text{end})}$. Set $\hat{J}^{(1)} := \hat{J}(\bar{u}^1)$ and $\hat{J}^{(\text{end})} := \hat{J}(\bar{u}^{(\text{end})})$.
		\STATE Compute $z^{(1)}$ by \eqref{Equation:FirstReferencePoint}. 
		\FOR{$n=1,\ldots$}
			\STATE Solve \eqref{BP-2} with reference point $z^{(n)}$ and save the solution as $\bar{u}^{(n+1)}$ and $\hat{J}^{(n+1)} := \hat{J}(\bar{u}^{(n+1)})$.
			\STATE Compute $z^{(n+1)}$ by \eqref{Equation:IterationReferencePoints}.
			\IF{$z^{(n+1)}_1 > \hat{J}^{(\text{end})}_1$}
				\RETURN $(\bar{u}^{(1)},\ldots,\bar{u}^{(n+1)},\bar{u}^{(\text{end})})$ and $(\hat{J}^{(1)},\ldots,\hat{J}^{(n+1)},\hat{J}^{(\text{end})})$
			\ENDIF
		\ENDFOR
	\end{algorithmic}
    \hrule
\end{algorithm}

\bigskip\noindent
For more details we refer to \cite[Section 6.1]{BBV17} and \cite[Section 3.4]{Ban17}.

\subsection{The POD approximation}
\label{Section:MOP-3}

In this section we want to show how the POD method can be successfully applied to the bicriterial optimal control problem \eqref{MOP-1}. In particular, we will follow the outline of Section \ref{SIAM-Book:Section4.2} and show a-priori convergence results and a-posteriori estimates. The analysis is only done for the continuous POD method, but can also be transferred to the POD method for the semidiscrete and the fully discrete approximation. \\

Suppose that we have computed a POD basis $\{\psi_i\}_{i=1}^\ell$ of rank $\ell$. Then the reduced-order cost function is given as
\begin{align*}
    \hat J^\ell(u)=J(\hat y^\ell+\mathcal S^\ell u,u)\quad \text{for }u\in\Uad.
\end{align*}

\begin{remark}
\rm
From the definition it follows that $\hat{J}_2 = \hat{J}_2^\ell$. However, to avoid misunderstandings we will still write $\hat{J}_2^\ell$ if we want to refer to the second cost function in the POD case.\hfill$\blacksquare$
\end{remark}
The reduced-order bicriterial optimal control problem then reads
\begin{equation}
    \label{BP-1-POD}
    \tag{$\mathbf{\widehat{BP}}^\ell$}
    \min \hat J^\ell(u)\quad\text{s.t.}\quad u\in\Uad.
\end{equation}
Its Pareto set and Paret front are denoted by $\Psell$ and $\Pfell$, respectively. For a given reference point $z \in \Pfell+\mathbb R^2_\le$ the Euclidean reference point problem \eqref{BP-2} is replaced by its reduced-order equivalent
\begin{equation}
    \tag{$\mathbf{\widehat{BP}}^\ell_z$}
    \label{BP-2-POD}
    \min F^\ell_z(u)\quad\text{s.t.}\quad u\in\Uad,
\end{equation}
where the reduced-order distance function $F^\ell_z:\U\to\mathbb R^2$ is given by
\begin{align*}
    F^\ell_z(u)=\frac{1}{2}\,{\|\hat J^\ell(u)-z\|}_2^2 = \frac{1}{2}\,\big(\hat J^\ell_1(u)-z_1\big)^2+\frac{1}{2}\,\big(\hat J^\ell_2(u)-z_2\big)^2\quad \text{for }u\in\U.
\end{align*}

\begin{remark}
    \label{Remark:BP-ResultsPOD}
    \rm The results of the Lemmas \ref{Lemma:BP-1} and \ref{Lemma:BP-2}, and the Propositions \ref{Proposition:BP-1} and \ref{Proposition:BP-AllPoints} still hold true for the reduced-order problem. Moreover, the first-order necessary optimality condition for the optimal solution $\bar u_z^\ell$ of \eqref{BP-2-POD} is given by
    \begin{equation}
        \label{BP-3-A}
        \big\langle(\hat J_1^\ell(\bar u_z^\ell)-z_1)(\mathcal G^\ell)^\star(\mathcal G^\ell\bar u_z^\ell + \hat{y}^\ell - \yd)+(\hat J_2^\ell(\bar u_z^\ell)-z_2)(\bar u_z^\ell-\un),u-\bar u_z^\ell\big\rangle_\U\ge 0
    \end{equation}
    for all $u\in\Uad$, where the operator $\mathcal G^\ell$ is defined as $\mathcal G^\ell = \mathcal{C} \mathcal{S}^\ell$.\hfill$\blacksquare$
\end{remark}

\subsubsection{POD a-priori analysis}
\label{SIAM-Book:Section5.3.3.1}

In this section we deal with the a-priori analysis for the optimality system of the Euclidean reference point problem \eqref{BP-2}, which is no linear-quadratic optimal control problem anymore due to the concatenation of the linear-quadratic cost functions with the Euclidean norm, for the continuous POD method. However, the idea and the procedure of our analysis strongly follows Section \ref{SIAM-Book:Section4.2.3}. A-priori estimates for both the state and adjoint equation for the continuous POD method can be found in Section \ref{SIAM-Book:Section3.3.2} and Section \ref{SIAM-Book:Section4.2.2}. In particular, for a given reference point $z$ we are interested in a-priori convergence results for the error between the solutions $\bar{u}$ and $\bar{u}^\ell$ to \eqref{BP-2} and \eqref{BP-2-POD}, respectively. To be able to prove an a-priori convergence result for the optimal controls of \eqref{BP-2-POD} to the one of \eqref{BP-2}, we first need to show the following result involving the first cost function $\hat{J}_1$ and its gradient $\nabla \hat{J}_1$.

\begin{lemma}
    \label{Lemma:CF-AprioriConvergence}
    Assume that Assumptions~{\rm\ref{A8}}, {\rm\ref{A9}} and {\rm\ref{A10}} for a $u \in \U$ are satisfied. Then it holds
    \begin{equation}
        \label{Equation:AprioriEstimateCostFunction}
        \left|\hat{J}_1(u)-\hat{J}_1^\ell(u) \right|^2 \le C\cdot\left\{
        \begin{aligned}
            &\sum_{i>\ell}\lambda_i^H\big({\|\psi_i^H\|}_V^2+1\big)&&\text{for }X=H,~\mathcal P^\ell=\mathcal P^\ell_H,\\
            &\sum_{i>\ell}\lambda_i^H\big\|\psi_i^H-\mathcal P^\ell_V\psi_i^H\big\|_V^2&&\text{for }X=H,~\mathcal P^\ell=\mathcal P^\ell_V,\\
            &\sum_{i>\ell}\lambda_i^V{\|\psi_i^V-\mathcal Q_H^\ell\psi_i^V\|}_V^2&&\text{for }X=V,~\mathcal P^\ell=\mathcal Q^\ell_H,\\
            &\sum_{i>\ell}\lambda_i^V&&\text{for }X=V,~\mathcal P^\ell=\mathcal Q^\ell_V,
        \end{aligned}
        \right.
    \end{equation}
    for a constant $C > 0$ independent of $\ell$ and
    \begin{equation}
        \label{Equation:AprioriEstimateGradient}
        \left\Vert \nabla \hat{J}_1(u) - \nabla \hat{J}_1^\ell(u) \right\Vert_\U^2 \le C \cdot \left\{
        \begin{aligned}
            &\sum_{i>\ell}\lambda_i^H\,\big\|\psi_i^H\big\|_V^2,&&X=H,~\mathcal P^\ell=\mathcal P^\ell_H,\\
            &\sum_{i>\ell}\lambda_i^V\,\big\|\psi_i^V-\mathcal P^\ell\psi_i^V\big\|_V^2,&&X=V,~\mathcal P^\ell=\mathcal Q^\ell_H,
        \end{aligned}
        \right.
    \end{equation}
    for a constant $C > 0$ independent of $\ell$. In particular, it holds
    \begin{equation}
        \lim_{\ell \to \infty} \left|\hat{J}_1(u)-\hat{J}_1^\ell(u) \right| = 0, \quad\lim_{\ell \to \infty} \left\Vert \nabla \hat{J}_1(u) - \nabla \hat{J}_1^\ell(u) \right\Vert_\U = 0.
    \end{equation}
\end{lemma}

\noindent{\bf\em Proof.}
To prove \eqref{Equation:AprioriEstimateCostFunction} we note that it is possible to derive the following estimate for an arbitrary control $u \in \U$, (cf. \cite[Proof of Theorem 5.31]{Ban17}, where a similar estimate was shown for the case that $\hat{y}^\ell = \hat{y}$),
\begin{equation}
    \begin{aligned}
        &\left| \hat{J}_1(u)-\hat{J}_1^\ell(u) \right| \\
        &\quad\leq \frac{1}{2} \left| \big\| \mathcal Su + \hat{y} \big\|_{L^2(0,T;H)} + \big\| \mathcal S^{\ell} u + \hat{y}^\ell \big\|_{L^2(0,T;H)} \right| \big\|\mathcal Su + \hat{y} - (\mathcal S^{\ell}u+\hat{y}^\ell) \big\|_{L^2(0,T;H)} \\
        &\qquad+ \big\| \mathcal Su + \hat{y} -(\mathcal S^{\ell}u + \hat{y}^\ell) \big\|_{L^2(0,T;H)} \big\| y_d \big\|_{L^2(0,T;H)}.
    \end{aligned}
\end{equation}
Now \eqref{Equation:AprioriEstimateCostFunction} follows from Theorem \ref{Th:A-PrioriError} and Corollary \ref{Corollary:HI-21}. For \eqref{Equation:AprioriEstimateGradient} we can use the formula for the gradient to conclude
\begin{equation}
    \begin{aligned}
        \big\|\nabla \hat{J}_1(u) - \nabla \hat{J}_1^\ell(u) \big\|_\U \leq \|\mathcal B'\|_{\mathscr L(L^2(0,T;V),\U)} \big\|\mathcal{A}u + \hat{p} - (\mathcal{A}^\ell u + \hat{p}^\ell) \big\|_{L^2(0,T;V)}.
    \end{aligned}
\end{equation}
Together with Theorem \ref{Th:DualA-PrioriError-2} this implies \eqref{Equation:AprioriEstimateGradient}.
\hfill$\Box$

\medskip
Now it is possible to show a-priori convergence with a convergence rate.

\begin{theorem}
    \label{Theorem:BP-AprioriConvergence1}
    Let $z \in \bigcap_{\ell \in \mathbb{N}} ( \Pfell + \mathbb{R}_\le^2 ) \bigcap (\Pf + \mathbb{R}_\le^2)$ be a reference point and denote by $\bar{u}$ and $\bar{u}^\ell$ the solutions to \eqref{BP-2} and \eqref{BP-2-POD}, respectively. Assume that Assumptions~{\rm\ref{A8}}, {\rm\ref{A9}} and {\rm\ref{A10}} with $u = \bar{u}$ hold. Let $\|\mathcal{Q_H^\ell} \|_{\mathcal{L}(V)}$ be bounded independently of $\ell$. If it holds $\hat{J}_2(\bar{u}) > z_2$, then
    \begin{equation}
        \big\|\hat{J}(\bar u)-\hat{J}^\ell(\bar u^\ell) \big\|^2 + \big\|\bar u-\bar u^\ell\big\|_\U^2 \le C\cdot\left\{
        \begin{aligned}
            &\sum_{i>\ell}\lambda_i^H\,\big\|\psi_i^H\big\|_V^2,&&X=H,~\mathcal P^\ell=\mathcal P^\ell_H,\\
            &\sum_{i>\ell}\lambda_i^V\,\big\|\psi_i^V-\mathcal P^\ell\psi_i^V\big\|_V^2,&&X=V,~\mathcal P^\ell=\mathcal Q^\ell_H.
        \end{aligned}\right.
    \end{equation}
    In particular, it holds
    \begin{equation}
        \lim_{\ell \to \infty} \big\|\hat{J}(\bar u)-\hat{J}^\ell(\bar u^\ell) \big\| + \big\| \bar{u} - \bar{u}^\ell \big\|_{\U} = 0.
    \end{equation}
\end{theorem}

\noindent{\bf\em Proof.}
We follow the lines of the proof of \cite[Theorem 5.41]{Ban17} and get the estimate
\begin{align*}
    & \big\|\hat{J}^{\ell}(\bar{u}^{\ell}) - \hat{J}(\bar{u}) \big\|_{\mathbb{R}^2}^2 +  \left( \frac{\hat{J}_1^{\ell}(\bar{u}^{\ell}) + \hat{J}_1(\bar{u})}{2} - z_1 \right) \big\| \mathcal S^{\ell}(\bar{u}^{\ell} - \bar{u}) \big\|_{L^2(0,T;H)}^2 \\
    & + \left( \frac{\hat{J}_2(\bar{u}) + \hat{J}_2^\ell(\bar{u}^{\ell})}{2} - z_2 \right) \big\|\bar{u}^{\ell} - \bar{u} \big\|_\U^2\\
    & \leq \left( \hat{J}_1(\bar{u}) - z_1 \right) \langle \nabla \hat{J}_1(\bar{u}) - \nabla \hat{J}_1^{\ell}(\bar{u}) , \bar{u}^{\ell} - \bar{u} \rangle_U + \left( \hat{J}_1(\bar{u}) - \hat{J}_1^{\ell}(\bar{u}^{\ell}) \right) \left( \hat{J}_1(\bar{u}) - \hat{J}_1^{\ell}(\bar{u}) \right).
\end{align*}
Using $\hat{J}_1^{\ell}(\bar{u}^{\ell}), \hat{J}_1(\bar{u}) \geq z_1$, $\hat{J}_2^\ell(\bar{u}^{\ell}) \geq z_2$ (which follow from \cite[Corollary 3.41]{Ban17}), the Cauchy-Schwarz inequality and Young's inequality we get
\begin{align*}
    \big\|\hat{J}^{\ell}(\bar{u}^{\ell}) - \hat{J}(\bar{u}) \big\|_{\mathbb{R}^2}^2 + \frac{\hat{J}_2(\bar{u}) - z_2}{2} \big\|\bar{u}^{\ell} - \bar{u} \big\|_\U^2 &\leq C(\epsilon) \left( \hat{J}_1(\bar{u}) - z_1 \right)^2\big\|\nabla \hat{J}_1(\bar{u}) - \nabla \hat{J}_1^{\ell}(\bar{u})\big\|_\U^2 + \epsilon\,\big\|\bar{u}^{\ell} - \bar{u}\big\|_\U^2 \\
    &\quad + \frac{1}{2} \left( \hat{J}_1(\bar{u}) - \hat{J}_1^{\ell}(\bar{u}^{\ell}) \right)^2 + \frac{1}{2} \left( \hat{J}_1(\bar{u}) - \hat{J}_1^{\ell}(\bar{u}) \right)^2,
\end{align*}
which yields
\begin{align*}
    \big\|\hat{J}^{\ell}(\bar{u}^{\ell}) - \hat{J}(\bar{u}) \big\|_{\mathbb{R}^2}^2 + \big\|\bar{u}^{\ell} - \bar{u} \big\|_\U^2 \leq C \left( \big\|\nabla \hat{J}_1(\bar{u}) - \nabla \hat{J}_1^{\ell}(\bar{u}) \big\|_\U^2 + \left( \hat{J}_1(\bar{u}) - \hat{J}_1^{\ell}(\bar{u}) \right)^2 \right)
\end{align*}
for a constant $C > 0$ independent of $\ell$. Now the claim follows from Lemma \ref{Lemma:CF-AprioriConvergence}.\hfill\endproof

Theorem \ref{Theorem:BP-AprioriConvergence1} together with the continuity of $\hat{J}_1$ and $\nabla \hat{J}_1$ and a-priori convergence statements for the state and adjoint equation also yields a-priori convergence results for the cost functions and its gradients, cf. \cite[Corollary 5.42]{Ban17}.

\begin{theorem}
    \label{Theorem:BP-AprioriConvergenceCostFunc}
    Let the assumptions of Theorem \ref{Theorem:BP-AprioriConvergence1} be satisfied. Then it holds
    \begin{enumerate}
        \item [\rm 1)] $\hat{J}_1(\bar{u}^\ell) \to \hat{J}_1(\bar{u})$ as $\ell \to \infty$.
        \item [\rm 2)] $\nabla \hat{J}_1(\bar{u}^\ell) \to \nabla \hat{J}_1(\bar{u})$ and $\nabla \hat{J}_1^\ell(\bar{u}^\ell) \to \nabla \hat{J}_1(\bar{u})$ as $\ell \to \infty$.
        \item [\rm 3)] $F_z(\bar{u}^\ell) \to F_z(\bar{u})$ and $F_z^\ell(\bar{u}^\ell) \to F_z(\bar{u})$ as $\ell \to \infty$.
        \item [\rm 4)] $\nabla F_z(\bar{u}^\ell) \to \nabla F_z(\bar{u})$ and $\nabla F_z^\ell(\bar{u}^\ell) \to \nabla F_z(\bar{u})$ as $\ell \to \infty$.
    \end{enumerate}
\end{theorem}

\subsubsection{POD a-posteriori analysis}
\label{SIAM-Book:Section5.3.3.2}
 
The focus of this section is to discuss an a-posteriori estimate for the error $\|\bar{u} - \bar{u}^\ell\|_{\U}$ between the solution $\bar{u}$ of \eqref{BP-2} and a suboptimal control $\bar{u}^\ell$, which is possibly -- but not necessarily -- the solution of \eqref{BP-2-POD}. As in Section \ref{SIAM-Book:Section4.2.4} we make use of a perturbation method (cf. \cite{DHPY95}) and the fact that the cost functions $\hat{J}_1$ and $\hat{J}_2$ are still linear-quadratic. Thus, for details of the idea of the perturbation method we directly refer to Section \ref{SIAM-Book:Section4.2.4}. Here, we only state the following theorem, which was proved in \cite[Corollary 5.43 and Theorem 5.45]{Ban17}.

\begin{theorem}
    \label{Theorem:BP-Aposteriori1}
    Suppose that Assumption~{\rm\ref{A9}} holds. Let $z \in \Pf + \mathbb{R}_\le^2$ be a reference point and $\{\psi_i\}_{i=1}^\ell$ be a POD basis of rank $\ell$. Define the function $\zeta^\ell\in\U$ by
    \begin{align*}
        \zeta_i^\ell(t)=\left\{
        \begin{aligned}
            &-\min(0,(\nabla F_z^\ell(\bar{u}^\ell))_i(s))&&\text{a.e. in }\mathscr A^\ell_{\mathsf ai}=\big\{s\in\mathscr D\,|\bar u^\ell_i(s)=u_{\mathsf ai}(s)\big\},\\
            &-\max(0,(\nabla F_z^\ell(\bar{u}^\ell))_i(s))&&\text{a.e. in }\mathscr A^\ell_{\mathsf bi}=\big\{s\in\mathscr D\,|\bar u^\ell_i(s)=u_{\mathsf bi}(s)\big\},\\
            &-(\nabla F_z^\ell(\bar{u}^\ell))_i(t) && \text{a.e. in }\mathscr D\setminus\big(\mathscr A^\ell_{\mathsf ai}\cup\mathscr A_{\mathsf bi}^\ell\big).
        \end{aligned}
        \right.
    \end{align*}
    If $\hat{J}(\bar{u}^\ell) \geq z$ and $\hat{J}_2(\bar{u}^\ell) > z_2$ are satisfied, then the a-posteriori error estimates
    \begin{equation}
        \label{Equation:AposterioriFormula}
        {\|\bar u-\bar u^\ell\|}_\U \le \left( \frac{\hat{J}_2(\bar{u}^\ell) - z_2}{2} \right)^{-1} \,{\|\zeta^\ell\|}_\U,\quad{\|\hat{J}(\bar u)-\hat{J}(\bar u^\ell)\|}_\U  \le \frac{1}{2} \left( \frac{\hat{J}_2(\bar{u}^\ell) - z_2}{2} \right)^{-1/2} \,{\|\zeta^\ell\|}_\U
    \end{equation}
    hold. In particular, it also holds $\lim \limits_{\ell \to \infty}\|\zeta^\ell\|_\U=0$.
\end{theorem}

\begin{remark}
\label{Remark:ERPM-WSM}
    \rm In contrast to the a-posteriori estimate in Theorem \ref{GVLuminy:Theorem4.5.1}, we also get an estimate for the error between the optimal and the suboptimal cost function values $\hat{J}(\bar{u})$ and $\hat{J}(\bar{u}^\ell)$. Also for the a-priori convergence statement in Theorem \ref{Theorem:BP-AprioriConvergence1} we get a convergence rate for the error $\left\Vert \hat{J}(\bar{u}) - \hat{J}^\ell(\bar{u}^\ell) \right\Vert$ -- again in contrast to the a-priori convergence result in Theorem \ref{GVLuminy:Theorem4.4.2}. Together with the statement from Proposition \ref{Proposition:BP-AllPoints} 2), this is the main reason to use the Euclidean reference point method instead of the weighted-sum method (see e.g. \cite[Chapter 3]{Ehr05}) to solve the bicriterial optimal control problem \eqref{BP-1}. By concatenating the cost functions with the Euclidean norm the Euclidean reference point problem gives us more geometrical information about the Pareto front, which can then be used to generate a uniform approximation of it, even in the case when POD is applied. This makes up for the fact that Euclidean reference point problems are not linear-quadratic anymore and thus harder to solve than the weighted-sum problems in general.\hfill$\blacksquare$
\end{remark}

Using the a-posteriori estimate from Theorem \ref{Theorem:BP-Aposteriori1} we set up an algorithm using POD for model-order reduction and controlling the maximal error between the POD and the finite element solutions; see Algorithm~\ref{alg:rpm_full-a}.

\bigskip
\hrule
\vspace{-3.5mm}
\begin{algorithm}[(Bicriterial Optimal Control using POD)]
    \label{alg:rpm_full-a}
    \vspace{-3mm}
    \hrule
    \vspace{0.5mm}
    \begin{algorithmic}[1]
		\REQUIRE Recursive parameters $h^{||} > 0$, $h^\perp \geq 0$, weighted-sum parameter $\alpha > 0$, initial length of POD basis $\ell_0$, number of incremented basis elements $\ell_{\text{incr}}$, maximal error in control space $\varepsilon_{\text{max}}$.
		\STATE Minimize $\hat{J}_1 + \alpha \hat J_2$ and $\hat J_2$ and save the solutions as $\bar{u}^{(1),\ell}$ and $\bar{u}^{(\text{end}),\ell}$. Set $\hat{J}^{(1),\ell} := \hat{J}(\bar{u}^{(1),\ell})$ and $\hat{J}^{(\text{end}),\ell} := \hat{J}(\bar{u}^{(\text{end}),\ell})$.
		\STATE Compute a POD basis of length $\ell_0$ using the snapshots $\mathcal{S} \bar{u}^{(1),\ell} + \hat{y}$ and $\mathcal{A} \bar{u}^{(1),\ell} + \hat{p}$.
		\STATE Compute $z^{(1),\ell}$ by \eqref{Equation:FirstReferencePoint}. 
		\FOR{$n=1,\ldots$}
			\STATE Set bool $\leftarrow$ FALSE.
			\WHILE{bool is FALSE}
				\STATE Solve \eqref{BP-2-POD} with reference point $z^{(n),\ell}$.
				\STATE Compute the a-posteriori estimate $\eta$ from \eqref{Equation:AposterioriFormula}
				\IF{$\eta > \varepsilon_{\text{max}}$}
					\STATE Set size $\ell$ of POD basis to $\ell \leftarrow \ell + \ell_{\text{incr}}$.
				\ELSE
					\STATE Save the solution to \eqref{BP-2-POD} as $\bar{u}^{(n+1),\ell}$ and $\hat{J}^{(n+1),\ell} := \hat{J}(\bar{u}^{(n+1),\ell})$.
					\STATE Set bool $\leftarrow$ TRUE.
				\ENDIF
			\ENDWHILE
			\STATE Compute $z^{(n+1),\ell}$ by \eqref{Equation:IterationReferencePoints}.
			\IF{$z^{(n+1),\ell}_1 > \hat{J}^{(\text{end}),\ell}_1$}
				\RETURN $(\bar{u}^{(1),\ell},\ldots,\bar{u}^{(n+1),\ell},\bar{u}^{(\text{end}),\ell})$ and $(\hat{J}^{(1),\ell},\ldots,\hat{J}^{(n+1),\ell},\hat{J}^{(\text{end}),\ell})$.
			\ENDIF
		\ENDFOR
	\end{algorithmic}
    \hrule
\end{algorithm}

\bigskip\noindent
Of course, this simple POD basis extension strategy might not work for complex problems. More elaborate strategies to adapt the POD basis when necessary can be found in \cite{Mak18,Spu19,BMV19}.


\renewcommand{\thechapter}{\Alph{chapter}}
\setcounter{chapter}{0}
\chapter{Appendix}
\label{Chapter:Appendix}

\renewcommand{\theequation}{\Alph{chapter}.\arabic{section}.\arabic{equation}}
\renewcommand{\thecorollary}{\Alph{chapter}.\arabic{section}.\arabic{corollary}}
\renewcommand{\thetheorem}{\Alph{chapter}.\arabic{section}.\arabic{theorem}}
\renewcommand{\theproposition}{\Alph{chapter}.\arabic{section}.\arabic{proposition}}
\renewcommand{\theremark}{\Alph{chapter}.\arabic{section}.\arabic{remark}}
\renewcommand{\thedefinition}{\Alph{chapter}.\arabic{section}.\arabic{definition}}
\renewcommand{\thefigure}{\Alph{chapter}.\arabic{section}.\arabic{figure}}
\renewcommand{\thetable}{\Alph{chapter}.\arabic{section}.\arabic{figure}}
\renewcommand{\theexample}{\Alph{chapter}.\arabic{section}.\arabic{example}}
\renewcommand{\theassumption}{\Alph{chapter}.\arabic{section}.\arabic{assumption}}
\renewcommand{\theproblem}{\Alph{chapter}.\arabic{section}.\arabic{probleme}}

\section{Basics}
\label{App:Basics}

\begin{definition}[Complex numbers]
	\label{def:complexNumbers}
	\begin{enumerate}
	    \item [\rm 1)] For a complex number $\alpha \in \mathbb C$, let $\Re(\alpha)$ and $\Im(\alpha) \in \mathbb R$ denote its real and imaginary part, so we have $\alpha = \Re(\alpha) + \jmath \cdot \Im(\alpha)$ where $\jmath$ is the complex unity. The \emph{complex conjugate} of $\alpha$ is given by $\overline \alpha = \Re(\alpha) - \jmath \cdot \Im(\alpha)$\index{a@$\overline{\alpha}$}.
	    \item [\rm 2)] We write a vector $z \in \mathbb C^n$ as $z = (z_i)_i$ with $z_1,\hdots,z_n \in \mathbb C$ and a matrix $Z \in \mathbb C^{m \times n}$ as $Z = (Z_{ij})_{ij}$ with $Z_{ij} \in \mathbb C$ for $i=1,...,m$, $j=1,...,n$.\index{z@$(z_i)_i,~(Z_{ij})_{ij}$} The vector $z$ can then also be understood as a matrix $z \in \mathbb C^{n \times 1}$. For the matrix, we define its \emph{transpose} as $Z^\top = (Z_{ij})_{ji} \in \mathbb C^{n \times m}$ and its \emph{Hermitian} as $Z^\mathsf H = (\overline{Z}_{ij})_{ji} \in \mathbb C^{n \times m}$.\index{z@$Z^\top,~Z^\mathsf H$}.
	    \item [\rm 3)] The canonical inner product in $\mathbb C^n$ is given by $\langle y, z \rangle_{\mathbb C^n} = y^\top \overline z$ for $y,z \in \mathbb C^n$. The induced norm is given by $\|z\|_{\mathbb C^n} = \sqrt{ \langle z, z \rangle_{\mathbb C^n} }$ for $z \in \mathbb C^n$.
	    \item [\rm 4)] We write $|\cdot|_2$ for the Euclidean norm, and $\langle\cdot\,,\cdot\rangle_2$ stands for the canonical inner product in $\mathbb R^m$.
	\end{enumerate}
\end{definition}

Recall that the Frobenius norm is given by
\begin{align*}
    {\|\bA \|}_F=\sqrt{\sum_{i=1}^m \sum_{j=1}^n \big|A_{ij}\big|^2} \quad \text{for }\bA \in \mathbb C^{m \times n}.
\end{align*}

\begin{definition}[Tensor notation]
    Let $X$ be a set, $n,m \in \mathbb{N}$ and $b_i \in X$ for $1 \leq i \leq n$ as well as $a_{ij} \in X$ for $1 \leq i \leq n$ and $1 \leq j \leq m$. Then we write $(b_i)$ for the vector $(b_1,\ldots,b_n) \in X^n$ and $((a_{ij}))$ for the matrix $(a_{ij})_{1 \leq i \leq n, \, 1 \leq j \leq m} \in X^{n \times m}$.
\end{definition}

\begin{lemma}[Young's inequality, \cite{Alt92}]
	\label{lem:youngsInequality}
	For \index{Inequality!Young}any real numbers $a,b,\varepsilon>0$ and for all $p,q \in (1,\infty)$ with $\tfrac{1}{p} + \tfrac{1}{q} = 1$, it holds
	\begin{align*}
		ab \le \frac{\varepsilon a^p}{p} + \frac{b^q}{q \varepsilon^{q/p}}
	\end{align*}
\end{lemma}

\begin{definition}[Gradient and Hessian]
	\begin{enumerate}
	    \item [\rm 1)] If $f: \mathbb R^n \to \mathbb R^m$ is a differentiable function with $m>1$, we write its \emph{Jacobian matrix} as $\nabla f(x) = (\tfrac{\partial f_i}{\partial x_j}(x))_{ij} \in \mathbb R^{m \times n}$ for $x \in \mathbb R^n$. If $m=1$, we use this symbol for the \emph{gradient} $\nabla f(x) = (\tfrac{\partial f}{\partial x_1}(x),\hdots,\tfrac{\partial f}{\partial x_n}(x))^\top \in \mathbb R^n$. 
	   \item [\rm 2)] If $f$ is twice diffentiable with $m=1$, we denote the \emph{Hessian matrix} by
    	\begin{align*}
    		\nabla^2f(x) = \nabla_{xx} f(x) = \big( \tfrac{\partial^2 f}{\partial x_i \partial x_j}(x) \big)_{ij} \in \mathbb R^{n \times n} \qquad \text{for } x \in \mathbb R^n
    	\end{align*}
    \end{enumerate}
\end{definition}

Recall that if $f$ is twice continuously differentiable, the Schwarz theorem states that $(\nabla^2 f(x))^\top = \nabla^2 f(x)$ holds for all $x \in \mathbb R^n$. 
 
\begin{proposition}[Gronwall inequality, {\cite[p.~708]{Eva08}}]
    \label{Gronwall}
    Let \index{Inequality!Gronwall}$\eta:[\tc,\te]\to\mathbb R$ be a nonnegative, differentiable function satisfying
    \begin{align*}
        \eta'(t)\le\varphi(t)\eta(t)+\chi(t)\quad\text{a.e. in }(\tc,\te],
    \end{align*}
    where $\varphi$ and $\chi$ are real-valued, nonnegative, integrable functions on $[\tc,\te]$. Then,
    \begin{align*}
        \eta(t)\le\exp\bigg(\int_\tc^t\varphi(s)\,\mathrm ds\bigg)\bigg(\eta(\tc)+\int_\tc^t\chi(s)\,\mathrm ds\bigg)
    \end{align*}
    In particular, if $\chi=0$ almost everywhere in $[0,T]$ and $\eta(0)=0$ hold true, we have $\eta=0$ almost everywhere in $[0,T]$.
\end{proposition}

Setting $\tilde\eta(t)=\eta(\tc+\te-t)$ for $t\in[0,T]$ we easliy derive the following result from Proposition~\ref{Gronwall}.

\begin{corollary}
    \label{DualGronwall}
    Let $\eta:[0,T]\to\mathbb R$ be a nonnegative, differentiable function satisfying
    \begin{align*}
        -\eta'(t)\le\varphi(t)\eta(t)+\chi(t)\quad\text{a.e. in }(0,T],
    \end{align*}
    where $\varphi$ and $\chi$ are real-valued, nonnegative, integrable functions on $[\tc,\te]$. Then,
    \begin{align*}
        \eta(t)\le\exp\bigg(\int_\tc^t\varphi(s)\,\mathrm ds\bigg)\bigg(\eta(\te)+\int_\tc^t\chi(s)\,\mathrm ds\bigg)
    \end{align*}
    In particular, if $\chi=0$ almost everywhere in $[0,T]$ and $\eta(\te)=0$ hold true, we have $\eta=0$ almost everywhere in $[0,T]$.
\end{corollary}

Recall the following wel-known inequality.

\begin{lemma}[Bernoulli Inequality]
    \label{app_bernIneq}
    \index{Inequality!Bernoulli}
    For $x \ge -1$ and $n \in \mathbb N$, we have \index{Inequality!Bernoulli}
    \begin{align}
    	(1+x)^n \ge 1+nx
    \end{align}
\end{lemma}

\section{Linear operators in Banach spaces}

\begin{definition}[Banach space]
    A complex vector space $X$ together with a map $\| \cdot \|_X: X \to \mathbb R$ is called a \index{Space!normed}\emph{normed space} if
	\begin{enumerate}
		\setlength{\itemsep}{0.2ex}
		\item [\rm 1)] $\| \varphi \|_X \ge 0$ for all $\varphi \in X$ and $\|\varphi\|_X=0$ if and only if $\varphi = 0$. 
		\item [\rm 2)] $\| \alpha \varphi \|_X = |\alpha| \cdot \|\varphi\|_X$ for all $\varphi \in X$, $\alpha \in \mathbb C$. 
		\item [\rm 3)] $\| \varphi + \phi \|_X \le \| \varphi \|_X + \| \phi \|_X$ for all $\varphi, \phi \in X$. 
	\end{enumerate}
	The space $X$ is called a \index{Space!Banach}\emph{Banach space} if it is additionally complete, meaning that every Cauchy sequence converges in $X$. 
\end{definition}

\begin{definition}[Hilbert space]
    \label{App:DefHSpace}
	A complex vector space $X$ together with a map $\langle \cdot, \cdot \rangle_X: X \times X \to \mathbb C$ is called an \index{Space!Hilbert!inner product}\emph{inner product space} if
    \begin{enumerate}
		\item [\rm 1)] $\langle \varphi, \varphi \rangle_X \ge 0$ for all $\varphi \in X$ and $\langle \varphi, \varphi \rangle_X = 0$ if and only if $\varphi = 0$.
		\item [\rm 2)] The map $\langle \cdot, \varphi \rangle_X: X \to \mathbb C$ is linear for all $\varphi \in X$.
		\item [\rm 3)] $\langle \varphi, \phi \rangle_X = \overline{ \langle \phi, \varphi \rangle }_X$ for all $\varphi, \phi \in X$. 
	\end{enumerate}
    The space $X$ is then also a normed space with the induced norm $\|\varphi\|_X = \sqrt{ \langle \varphi, \varphi \rangle_X }$ for all $x \in X$. $X$ is called a \index{Space!Hilbert}\emph{Hilbert space} if it is additionally complete. A subset $\mathscr S \subset X$ is called \index{Space!Hilbert!orthoogonal}\emph{orthogonal} if it holds $\langle \phi, \psi \rangle_X = 0$ for all $\phi,\psi \in\mathscr S$ with $\phi \neq \psi$. It is called \index{Space!Hilbert!orthonormal}\emph{orthonormal} if we additionally have $\|\phi\|_X = 1$ for all $\phi \in\mathscr S$. Lastly, $\mathscr S$ is called an \index{Space!Hilbert!orthonormal basis}\emph{orthonormal basis} of $X$ if it is orthonormal and its linear hull is dense in $X$. The \index{Space!Hilbert!orthogonal complement}{\em orthogonal complement} of $\mathscr S$ is defined as
	\begin{align*}
	   \mathscr S^\bot=\big\{\varphi\in X\,\big|\,{\langle\varphi,\phi\rangle}_X=0\text{ for all }\phi\in\mathscr S\big\}.
	\end{align*}
	The space $X$ is called \index{Space!Hilbert!separable}\emph{separable} if it has a countable orthonormal basis $\{\psi_i\}_{i \in \mathbb I} \subset X$. Here we have $\mathbb I = \{1,...,n\}$ (if $X$ has the finite dimension $n \in \mathbb N$) or $\mathbb I = \mathbb N$ (if $X$ has infinite dimension). 
\end{definition}

\begin{example}
    \rm For any $n \in \mathbb{N}$ the spaces $\mathbb{R}^n$ and $\mathbb{C}^n$ endowed with the respective inner products
    \begin{align*}
		{\langle x , y \rangle}_{\mathbb{R}^n} & = x^\top y = \sum_{i=1}^n x_i y_i \quad \text{ for all } x,y \in \mathbb{R}^n \\
		\langle x , y \rangle_{\mathbb{C}^n} & := x^\top \bar{y} = \sum_{i=1}^n x_i \bar{y}_i \quad \text{ for all } x,y \in \mathbb{C}^n
	\end{align*}
    are Hilbert spaces. In both cases the induced norm is just given by the Euclidean norm. We refer to these inner products as the canonical inner product in $\mathbb{R}^n$ and $\mathbb{C}^n$, respectively.\hfill$\blacklozenge$
\end{example}

If $X$ is separable with orthonormal basis $\{\psi_i\}_{i \in \mathbb I}$, it can be shown \cite{RS80} that
\begin{align*}
	\varphi = \sum_{i \in \mathbb I} \langle \varphi, \psi_i \rangle_X \psi_i \quad \text{and} \quad {\| \varphi \|}_X^2 = \sum_{i \in \mathbb I} | \langle \varphi, \psi_i \rangle_X |^2 \qquad \text{for all } \varphi \in X.
\end{align*}

The following result is proven in \cite[p.~43]{RS80}.

\begin{theorem}[Riesz isomorphism]
    \label{Theorem:Riesz_isomorphism}
	Let $X$ be a Hilbert space with inner product $\langle \cdot,\cdot \rangle_X$. Then the mapping
    \begin{align}
		\label{app_riesz}
		\mathcal J_X: X \to X', \quad {\langle \mathcal J_X x_1, x_2 \rangle}_{X',X} = {\langle x_1,x_2 \rangle}_X
	\end{align}
    is an isometric isomorphism called the \index{Linear operator!Riesz isomorphism}\emph{Riesz isomorphism}.
\end{theorem}

\begin{theorem}[Projection theorem]
	\label{thm:projectionTheorem}
	Let $X$ \index{Theorem!projection}be a Hilbert space and $M \subset X$ a nonempty and closed subspace. Then the space $X$ can be written as the direct sum $X = M \oplus M^\perp$, i.e., for every $\varphi \in X$, there is a unique decomposition $\varphi = \varphi_M + \varphi_\perp$ with $\varphi_M \in M$ and $\varphi_\perp \in M^\perp$. Furthermore, $\varphi_M$ is the unique solution to the approximation problem
    \begin{align}
		\label{eq:projectionTheorem}
		\min\limits_{\phi \in M} ~~ \| \varphi - \phi \|_X.
	\end{align}
\end{theorem}

\noindent{\bf\em Proof}. 
\hfill The fact that $X = M \oplus M^\perp$ is a standard result and is proven, e.g., in \cite[Theorem~II.3]{RS80}. Let now $\varphi = \varphi_M + \varphi_\perp \in X$ with $\varphi_M \in M$ and $\varphi_\perp \in M^\perp$ and $\phi \in M$ be arbitrarily given. It follows that
\begin{align*}
	\| \varphi - \phi \|_X^2 = \| \varphi_M - \phi \|_X^2 + \| \varphi_\perp \|_X^2 \ge \| \varphi_\perp \|_X^2 = \| \varphi - \varphi_M \|_X^2.
\end{align*}
Therefore, $\varphi_M$ solves \eqref{eq:projectionTheorem}; cf. Figure~\ref{fig:projectionTheorem}.\hfill$\Box$

\begin{figure}
    \centering
    \begin{tikzpicture}
        \coordinate (Origin) at (-5,1.5);
        \coordinate (Plane) at (0,2.25);
        \coordinate (Vector) at (0,4.75);
        \filldraw [thick, fill=gray!60, even odd rule] (Origin)  coordinate (GeneralStart) -- ++(8,0) -- ++(2,1.5) -- ++(-8,0) -- ++(-2,-1.5) -- cycle;
        \node at ($(Plane)+(2,0)$) {\LARGE $M$};
        \draw[->,line width=0.6mm] (Origin) -- ($(Vector)-(0.1,0.1)$);
        \node at ($(Origin) + (2.5,2)$) {\LARGE $\varphi$};
        \draw[->,line width=0.6mm] (Origin) -- (Plane);
        \node at ($(Origin) + (4,0.25)$) {\LARGE $\varphi_M$};
        \draw[->,line width=0.6mm] (Plane) -- (Vector);
        \node at ($(Plane) + (0.5,1.3)$) {\LARGE $\varphi_\perp$};
    \end{tikzpicture}
    \caption{Visualization of the Theorem \ref{thm:projectionTheorem}.}
    \label{fig:projectionTheorem}
\end{figure}

\begin{lemma}[Cauchy-Schwarz inequality, \cite{RS80}]
\label{App:CSIneq}
\index{Inequality!Cauchy-Schwarz}
	Let $X$ be a Hilbert space with inner product $\langle \cdot, \cdot \rangle_X$. Then it holds
	\begin{align*}
		\big | {\langle \varphi, \phi \rangle}_X \big | \le {\| \varphi \|}_X{\| \phi \|}_X \qquad \text{for all } \varphi, \phi \in X.
	\end{align*}
\end{lemma}

\begin{lemma}
	\label{app_lem_innerProd1}
	Let $X$ be an inner product space with inner product $\langle \cdot, \cdot \rangle_X$. Then it holds
	\begin{align}
		\label{app_innerProd1}
		{\langle \varphi - \phi, \phi \rangle}_X = -\frac{1}{2}\,{\| \varphi - \phi \|}_X^2 + \frac{1}{2}\,{\|\varphi\|}_X^2 - \frac{1}{2}\,{\|\phi\|}_X^2 \quad \text{for all } \varphi, \phi \in X
	\end{align}
\end{lemma}

\noindent {\bf \em Proof}.
We have
\begin{align*}
    {\langle \varphi - \phi, \phi \rangle}_X + \frac{1}{2}\,{\| \varphi - \phi \|}_X^2 + \frac{1}{2}\,{\|\phi\|}_X^2 &= \frac{1}{2}\,{\langle \varphi - \phi, \varphi + \phi \rangle}_X + \frac{1}{2}\,{\|\phi\|}_X^2 = \frac{1}{2}\, {\| \varphi \|}_X^2
\end{align*}
where we have utilized the fact that $\langle \cdot\,, \cdot \rangle_X$ is bilinear.\hfill$\Box$

\begin{definition}[Linear operators in Banach spaces]
	\label{SIAM:Definition-I.1.1.1}
	Let $X$, $Y$ be two complex Banach spaces. The operator $\mathcal T: X\to X$ is called a \index{Linear Operator}\index{Linear operator!bounded}{\em linear, bounded operator} if these conditions are satisfied:
    \begin{enumerate}
		\item [\rm 1)] $\mathcal T(\alpha_1 \varphi + \beta \psi)= \alpha \mathcal T \varphi +\beta\mathcal T \phi$ for all $\alpha,\beta \in\mathbb C$ and $\varphi,\,\phi\in X$.
	   \item [\rm 2)] There exists a constant $c>0$ such that $\|\mathcal T \varphi\|_{Y}\le c\,\|\varphi\|_{X}$ for all $\varphi\in X$.
	\end{enumerate}
    The set of all linear, bounded operators from $X$ to $Y$ is denoted by $\mathscr L(X,Y)$ which is a Banach space equipped with the \index{Linear operator!norm}{\em operator norm} {\rm \cite[pp.~69-70]{RS80}}
    \begin{align*}
	{\|\mathcal T\|}_{\mathscr L(X,Y)}=\sup_{\|\varphi\|_{X}=1}{\|\mathcal T \varphi\|}_{Y}\quad\text{for }\mathcal T\in\mathscr L(X,Y).
	\end{align*}
    If $X=Y$ holds, we briefly write $\mathscr L(X)$ instead of $\mathscr L(X,X)$. For $\mathcal T\in\mathscr L(X,Y)$ the subspace
    \begin{align*}
		\mathrm{ran}\,\mathcal T=\mathcal T(X)\subset Y
	\end{align*}
    is called the \index{Linear operator!range, $\mathrm{ran}\,\mathcal T$}{\em range of $\mathcal T$}. The subspace
    \begin{align*}
		\mathrm{ker}\,\mathcal T=\{\varphi\in X\,|\,\mathcal T \varphi=0\}\subset X
	\end{align*}
    is said to be the \index{Linear operator!kernel, $\mathrm{ker}\,\mathcal T$}{\em kernel} or \index{Linear operator!null space, $\mathrm{ker}\,\mathcal T$}{\em null space of $\mathcal T$}.
\end{definition}

\begin{definition}[Dual operators in Banach spaces]
	\label{SIAM:Definition-Dual}
	Let $X,Y$ be two complex Banach spaces and $\mathcal T\in\mathscr L(X,Y)$. The \index{Linear operator!dual, $\mathcal T'$}{\em dual operator} $\mathcal T':Y'\to X'$ of a linear operator $\mathcal T\in\mathscr L(X,Y)$ is defined as
	\begin{align*}
		{\langle \mathcal T'f,\varphi \rangle}_{X',X}={\langle f,\mathcal T \varphi \rangle}_{Y',Y}\quad\text{for all }(\varphi,f)\in X\times Y',
	\end{align*}
	where, for instance, $\langle\cdot\,,\cdot\rangle_{X',X}$ denotes the dual pairing of the space $X$ with its dual space $X'=\mathscr L(X,\mathbb C)$.
\end{definition}

\begin{definition}[Adjoint operators in Hilbert spaces]
	\label{SIAM:Definition-Adjoint}
	Let $X$ and $Y$ denote two complex Hilbert spaces. For a given $\mathcal T\in\mathscr L(X,Y)$ the \index{Linear operator!adjoint, $\mathcal T^\star$}{\em adjoint operator} $\mathcal T^\star: Y\to X$ is uniquely defined by
    \begin{align*}
		{\langle\mathcal T^\star \phi,\varphi\rangle}_{X}={\langle \phi,\mathcal T \varphi\rangle}_{Y}=\overline{\langle \mathcal T \varphi,\phi\rangle}_{Y}\quad\text{for all }(\varphi,\phi)\in X\times Y.
	\end{align*}
\end{definition}

With the Riesz isomorphism, we can draw a connection between the dual operator from Definition \ref{SIAM:Definition-Dual} and the adjoint operator from Definition \ref{SIAM:Definition-Adjoint}. This result is proven e.g. in \cite[p.~186]{Tro10}.

\begin{theorem}
	\label{Th:Riesz}
	\rm
	Let\index{Theorem!Riesz} $X$ and $Y$ be two complex Hilbert spaces and $\mathcal J_{X}:X\to X'$, $\mathcal J_Y:Y\to Y'$ denote the Riesz isomorphisms from Theorem \ref{Theorem:Riesz_isomorphism}. For an operator $\mathcal T \in \mathscr L(X,Y)$, it then holds
    \begin{align*}
		\mathcal T^\star = \mathcal J_{X}^{-1} \mathcal T' \mathcal J_{Y}.
	\end{align*}
    Moreover, $(\mathcal T^\star)^\star = \mathcal T$.
\end{theorem}

\begin{definition}
	Let $X,Y$ be two complex Hilbert spaces and $\mathcal T \in \mathscr L(X,Y)$.
	\begin{enumerate}
		\item [\rm 1)] If $X = Y$, the operator $\mathcal T$ is called \emph{self-adjoint} if $\mathcal T^\star = \mathcal T$.
		\item [\rm 2)] If $\mathcal T$ is self-adjoint, we call it \emph{nonnegative} if $\langle \mathcal T \varphi, \varphi \rangle_{X} \ge 0$ for all $\varphi \in X$.  
		\item [\rm 3)] $\mathcal T$ is called \emph{compact} if for every bounded sequence $\{\varphi_n\}_{n \in \mathbb N} \subset X$, the sequence $\{\mathcal T \varphi_n\}_{n \in \mathbb N} \subset Y$ contains a convergent subsequence. 
	\end{enumerate}
\end{definition}

\begin{remark}
	\label{SIAM:Remark-RealNonnegativeEigenvalue}
	\rm Let $X$ be a complex Hilbert space. If $\mathcal T \in \mathscr L(X)$ is self-adjoint it is easy to show that each eigenvalue $\lambda$ of $\mathcal T$ fulfils $\lambda \in \mathbb{R}$. Additionally, if $\mathcal T$ is even nonnegative, it holds $\lambda \geq 0$ for each eigenvalue $\lambda$ of $\mathcal{T}$. \hfill$\blacksquare$
\end{remark}

\section{Spectral theory of linear operators}
\label{SIAM:Section_SpectralTheory}

\begin{definition}[Spectrum and resolvent]
	\rm
	Let $X$ be a complex Hilbert space and $\mathcal T \in \mathscr L(X)$.
	\begin{enumerate}
		\setlength{\itemsep}{0.2ex}
		\item [\rm 1)] A complex number $\lambda \in \mathbb C$ belongs to the \emph{resolvent set} $\rho(\mathcal T)$ if $\lambda \mathcal I - \mathcal T$ is a bijection with a bounded inverse. Here, $\mathcal I \in \mathscr L(X)$ stands for the identity operator in $X$. We define the \emph{spectrum} of $\mathcal T$ as $\sigma(\mathcal T) = \mathbb C \setminus \rho(\mathcal T)$.
		\item [\rm 2)] For a spectral point $\lambda \in \mathbb C$, if $\lambda \mathcal I - \mathcal T$ is not injective, $\lambda$ is said that to belong to the \emph{point spectrum} of $\mathcal T$. In this case, a vector $\varphi \in X \setminus \{0\}$ with $\mathcal T \varphi = \lambda \varphi$ is called an \emph{eigenvector} of $\mathcal T$. 
	\end{enumerate}
	To summarize, the complex numbers are separated into $\mathbb C = \sigma(\mathcal T) \cup \rho(\mathcal T)$ and part of the spectrum $\sigma(\mathcal T)$ is made up of the point spectrum.
\end{definition}

The following two theorems are essential for the spectral analysis of compact operators and can be found in \cite[p.~203]{RS80}.

\begin{theorem}[Riesz-Schauder]
	\label{SIAM:Theorem-I.1.1.1}
	Let\index{Theorem!Riesz-Schauder} $X$ be a complex Hilbert space and $\mathcal T: X \to X$ a linear, compact operator. Then the spectrum $\sigma(\mathcal T)$ is a discrete set having no limit points except perhaps $0$. Furthermore, every nonzero $\lambda \in \sigma(\mathcal T)$ belongs to the point spectrum of $\mathcal T$ and the space of eigenvectors to $\lambda$ is finite dimensional.
\end{theorem}

\begin{theorem}[Hilbert-Schmidt]
	\label{SIAM:Theorem-I.1.1.2}
	Let\index{Theorem!Hilbert-Schmidt} $X$ be a complex separable Hilbert space and $\mathcal T \in \mathscr L(X)$ a compact and selfadjoint operator. Then there is a sequence of real eigenvalues $\{\lambda_i\}_{i\in\mathbb I}$ and corresponding eigenvectors $\{\psi_i\}_{i \in \mathbb I} \subset X$ that form a complete orthonormal basis of $X$. If $X$ has infinite dimension, we further get $\lambda_i \to 0$ as $i \to \infty$. 
\end{theorem}	

The following result comes from \cite[III-\S 6.4]{Kat80}.

\begin{theorem}
	\label{thm:spectralDecomposition}
	Let $X$ be a Banach space, $\mathcal T \in \mathscr L(X)$ and $\Sigma \subset \sigma(\mathcal T)$ a part of the spectrum such that there exists a closed curve $\Gamma \subset \rho(\mathcal T)$ with $\Sigma$ on the inside and $\Sigma' := \sigma(\mathbb C) \setminus \Sigma$ on the outside. Then there exists a decomposition $X = M \oplus M'$ which satisfies
	\begin{enumerate}
		\item [\rm 1)] $\mathcal T(M) \subset M$ and $\mathcal T(M') \subset M'$.
		\item [\rm 2)] The spectra of the restricted operators $\mathcal T |_M \in \mathscr L(M)$ and $\mathcal T |_{M'} \in \mathscr L(M')$ coincide with $\Sigma$ respectively $\Sigma'$. 
	\end{enumerate}
\end{theorem}

\begin{corollary}
	\label{cor:eigenvalueDecomposition}
	Let $X$ be a Banach space, $\mathcal T \in \mathscr L(X)$ be compact and self-adjoint and $\Sigma := \{\lambda\}$ with $\lambda$ being a nonzero eigenvalue of $\mathcal T$. Then in the decomposition $X = M \oplus M'$ from Theorem \ref{thm:spectralDecomposition}, $M$ is given by the space of eigenvectors to $\lambda$. 
\end{corollary}

\noindent {\bf \em Proof.}
First of all, the closed curve $\Gamma \subset \rho(\mathcal T)$ may be chosen as a circle around $\lambda$. Since $\lambda$ is nonzero, Theorem \ref{SIAM:Theorem-I.1.1.1} tells us that it will be the only spectral point inside $\Gamma$ if the radius is chosen small enough. Therefore, Theorem \ref{thm:spectralDecomposition} is applicable. Furthermore, it was shown in Section III-\S 6.5 of \cite{Kat80} that the operator $D_\lambda := \lambda \mathcal I_M - \mathcal T |_M$ is quasi-nilpotent. Because $\mathcal T$ is compact, Theorem \ref{SIAM:Theorem-I.1.1.2} tells us that $M$ is finite-dimensional, so $D_\lambda$ is even nilpotent. However, the only self-adjoint and nilpotent operator is zero which yields $\lambda \mathcal I_M = \mathcal T|_M$. In other words, every nonzero vector in $M$ is an eigenvector to $\lambda$. On the other hand, if $x \in X$ is an arbitrary eigenvector of $\mathcal T$, it can be written as $x = x_1 + x_2$ with $x_1 \in M$, $x_2 \in M'$. Since $M, M'$ are $\mathcal T$-invariant by Theorem \ref{thm:spectralDecomposition}, it follows that $\mathcal T x_i = \lambda x_i$ for $i=1,2$. Because the spectrum of $\mathcal T|_{M'}$ does not contain $\lambda$, we infer $x_2 = 0$, so $x=x_1 \in M$.\hfill$\Box$

\begin{remark}
	\label{rem:eigenvalueDecomposition}
	\rm In \cite[IV-\S 3.5]{Kat80} it is stated that the previous claim also holds true for a finite system of eigenvalues, i.e., if $\mathcal T \in \mathscr L(X)$ is compact and self-adjoint and $\Sigma = \{\lambda_1,\ldots,\lambda_n\}$ is a set of $n$ distinct nonzero eigenvalues of $\mathcal T$, then in the decomposition $X = M \oplus M'$ from Theorem \ref{thm:spectralDecomposition}, $M$ is given by the space of eigenvectors to $\lambda_1,\ldots,\lambda_n$.\hfill$\blacksquare$
\end{remark}

The following result is taken from \cite[IV-\S 3.4]{Kat80}.

\begin{theorem}
	\label{thm:continuityOfEigenvalues}
	Let $X$ be a Banach space, $\mathcal T \in \mathscr L(X)$ and $\Sigma, \Sigma' \subset \sigma(\mathcal T)$ be parts of the spectrum separated by a closed curve $\Gamma \in \mathbb C \setminus \sigma(\mathcal T)$ that encloses $\Sigma$. Let further $X = M(\mathcal T) \oplus M'(\mathcal T)$ be the decomposition from Theorem \ref{thm:spectralDecomposition}. Then there exists a $\delta > 0$ depending on $\mathcal T$ and $\Gamma$ with the following properties:
	\begin{enumerate}
		\item [\rm 1)] Any $\mathcal S \in \mathscr L(X)$ with $\| \mathcal S - \mathcal T \|_{\mathscr L(X)} < \delta$ has its spectrum $\sigma(\mathcal S)$ likewise split into two parts $\Sigma(\mathcal S)$, $\Sigma'(\mathcal S)$ separated by $\Gamma$. It holds $\Gamma \subset \rho(\mathcal S)$. 
		\item [\rm 2)] In the associated decomposition $X = M(\mathcal S) \oplus M'(\mathcal S)$ from Theorem \ref{thm:spectralDecomposition}, $M(\mathcal S)$ and $M'(\mathcal S)$ are isomorphic to $M(\mathcal T)$ and $M'(\mathcal T)$, respectively.
	   \item [\rm 3)] Let $P(\mathcal T)$, $P(\mathcal S) \in \mathscr L(X)$ denote the projections from $X$ to $M(\mathcal T)$ and $M(\mathcal S)$, respectively. Then it holds $P(\mathcal S) \to P(\mathcal T)$ in $\mathscr L(X)$ as $\mathcal S \to \mathcal T$ in $\mathscr L(X)$.
	\end{enumerate}
\end{theorem}

\section{Integration and differentiation in Banach spaces}
\label{App:A.4}

\begin{definition}
\label{Definition_L2_BanachValued}
	Let $X$ be a Banach space and $\mathcal D \subset \mathbb R^{\mathfrak p}$. A function $f: \mathcal D \to X$ is called \emph{measurable} if there exists a sequence $\{f_n\}_{n \in \mathbb N}$ of step functions $f_n: \mathcal D \to X$ such that $f_n \to f$ almost everywhere. 
	For $p \in [1,\infty)$, the $p$-th \emph{Bochner space} is defined as
	\begin{align*}
		L^p(\mathcal D;X) = \big \{ f: \mathcal D \to X ~\big|~ f \text{ is measurable and } \int_{\mathcal D} \|f(\bx) \|_X^p ~\mathrm d \bx < \infty \big \}
	\end{align*}
	Endowed with the norm
	\begin{align*}
	\| f \|_p := \| f \|_{L^p(\mathcal D;X)} := \bigg( \int_{\mathcal D} \|f(\bx) \|_X^p ~\mathrm d \bx \bigg)^{1/p} \quad \text{ for all } f \in L^p(\mathcal D;X),
	\end{align*}
	The space $L^p(\mathcal D;X)$ is a Banach space. If $X=\mathbb R$ holds, we write briefly $L^p(\mathcal D)$ instead of $L^p(\mathcal D;\mathbb R)$. If $X$ together with the inner product $\langle \cdot , \cdot \rangle_X$ is even a Hilbert space, $L^2(\mathcal D;X)$ is also a Hilbert space \emph{\cite{AE01}} endowed with the inner product
	\begin{align*}
		{\langle f, g \rangle}_{L^2(\mathcal D;X)} = \int_{\mathcal D} {\langle f(\bx), g(\bx) \rangle}_X ~\mathrm d\bx \qquad \text{for all } f,g \in L^2(\mathcal D;X).
	\end{align*}
	Furthermore, the {\em essential Bochner supremum} is defined for a measurable funnction $f: \mathcal D \to X$ as
	\begin{align*}
		{\| f \|}_\infty = {\| f \|}_{L^\infty(\mathcal D;X)} := \inf \big \{ C \ge 0 ~\big|~ {\|f(\bx)\|}_X \le C \text{ for almost all } \bx \in \mathcal D \big \}.
	\end{align*}
	Finally, the {\em Bochner space of essentially bounded fuctions} consists of all measurable functions $f: \mathcal D \to X$ with $\|f\|_\infty < \infty$. 
\end{definition}

\begin{theorem}[Hölder inequality, \cite{AE01}]
	\index{Inequality!Hölder}
	Let $X$ be a Banach space and $\mathcal D \subset \mathbb R^{\mathfrak p}$. Then it holds
	\begin{align*}
		\int_{\mathcal D} \|f(\bx)\|_X \cdot \|g(\bx)\|_X ~\mathrm d \bx ~\le~ \|f\|_{L^2(\mathcal D; X)} \cdot \|g\|_{L^2(\mathcal D;X)} \quad \text{for all } f,g \in L^2(\mathcal D;X).
	\end{align*}
\end{theorem}

In the following definition, we use the symbol \index{n@$\mathbb N_0$} $\mathbb N_0 = \mathbb N \cup \{0\}$. An element $\alpha \in \mathbb N_0^{\mathfrak p}$ is called a \emph{multi-index} of order $|\alpha| = \sum_{i=1}^{\mathfrak p} \alpha_i$. For a function $\varphi: \mathcal D \to \mathbb R$ that is $|\alpha|$-times classically differentiable, we use the short notation \index{$\partial^\alpha$}
\begin{align*}
	\partial^\alpha \varphi(\bx) = \left( \frac{\partial}{\partial x_1} \right)^{\alpha_1} \hdots \left( \frac{\partial}{\partial x_{\mathfrak p}} \right)^{\alpha_{\mathfrak p}} \varphi(\bx) \quad \text{for } x \in \mathcal D.
\end{align*}

\begin{definition}[Sobolev space]
    \label{Definition_H1_BanachValued}
	Let $X$ be a Banach space, $\mathcal D \subset \mathbb R^{\mathfrak p}$ open and $f \in L^2(\mathcal D;X)$. For a multi-index $\alpha \in \mathbb N^{\mathfrak p}_0$, a function $f_{\alpha} \in L^2(\mathcal D;X)$ is called the \emph{weak $\alpha$-derivative} of $f$ if 
	\begin{align*}
		\int_{\mathcal D} f(\bx) \cdot \partial^\alpha \varphi(\bx) ~\mathrm d\bx = (-1)^{|\alpha|} \int_{\mathcal D} f_\alpha(\bx) \cdot \varphi(\bx) ~\mathrm d\bx \quad \text{for all } \varphi \in C_0^\infty(\mathcal D)
	\end{align*}
	Here, the space $C_0^\infty(\mathcal D)$ is the space of all infinitely differentiable functions $\varphi: \mathcal D \to \mathbb R$ such that the support $\text{supp } \varphi = \overline{ \{ \bx \in \mathcal D ~|~ \varphi(\bx) \neq 0 \} }$ is compact in $\mathbb R^{\mathfrak p}$. If $f$ has such a weak $\alpha$-derivative, we indicate this by writing $\partial^\alpha f \in L^2(\mathcal D;X)$ and identify $\partial^\alpha f$ with $f_\alpha$. For $k \in \mathbb N$, we define the \emph{Sobolev space}
	\begin{align*}
		H^k(\mathcal D;X) = \big \{ f \in L^2(\mathcal D;X) ~\big|~ \partial^\alpha f \in L^2(\mathcal D;X) \text{ for } |\alpha| \le k \big \}.
	\end{align*}
	If $X$ is is a Hilbert space, $H^1(\mathcal D;X)$ is also a Hilbert space \cite{DL00,Eva08} endowed with the inner product
    \begin{align*}
		{\langle f, g \rangle}_{H^1(\mathcal D;X)} = {\langle f, g \rangle}_{L^2(\mathcal D;X)} + \sum_{i=1}^{\mathfrak p} {\langle f_{x_i}, g_{x_i} \rangle}_{L^2(\mathcal D;X)} \quad \text{for } f,g \in H^1(\mathcal D;X).
	\end{align*}
	Here, we have used the notion of weak partial derivatives $f_{x_i} = \tfrac{\partial}{\partial x_i} f = \partial^\alpha f$ with $\alpha_i = 1$ and $\alpha_j = 0$ for $j \neq i$. 
\end{definition}

\begin{definition}[Gelfand triple]
	Let $V$ and $H$ be real, separable Hilbert spaces such that $V$ lies dense in $H$ with compact embedding. We identify $H$ with its dual space $H'$ so that we have the Gelfand triple\index{Gelfand triple}
	\begin{align*}
		V \hookrightarrow H \simeq H' \hookrightarrow V'
	\end{align*}
	where the last embedding is understood as follows: For every element $h' \in H'$ and $v \in V$, we also have $v \in H$ by the first embedding, so we can define $\langle h',v \rangle_{V',V} = \langle h',v \rangle_{H',H}$. 
\end{definition}

Let now $T>0$ and the Hilbert spaces $V$, $H$ form a Gelfand triple. Then we define the space
\begin{equation}
\label{EqW(0,T)}
W(0,T)=L^2(0,T;V)\cap H^1(0,T;V')=\big\{\varphi\in L^2(0,T;V)\,\big|\,\varphi_t \in L^2(0,T;V')\big\}
\end{equation}
	
The following properties of $W(0T)$ are proven e.g. in \cite[pp.~472-479]{DL00} and \cite[pp.~146-148]{Tro10}.
	
\begin{lemma}
	\label{Lemma:HI-20}
	Let $V, H$ form a Gelfand triple and $T>0$.
    \begin{enumerate}
		\item [\rm 1)] $W(0,T)$ is a Hilbert space endowed with the inner product
		\begin{align*}
			{\langle \varphi,\phi\rangle}_{W(0,T)}=\int_0^T{\langle\varphi(t),\phi(t)\rangle}_V+{\langle\varphi_t(t),\phi_t(t)\rangle}_{V'}\,\mathrm dt\quad\text{for } \varphi,\phi\in W(0,T)
		\end{align*}
		and the induced norm $\|\varphi\|_{W(0,T)}=\langle\varphi,\varphi\rangle_{W(0,T)}^{1/2}$. 
		\item [\rm 2)] The space $W(0,T)$ is continuously embedded into the space $C([0,T];H)$, the space of all continuous mappings from $[0,T]$ to $H$. Thus, there exists a constant $C_W>0$ such that
		\begin{equation}
    		\label{Wemb}
    		{\|\varphi\|}_{C([0,T];H)}\le C_W\,{\|\varphi\|}_{W(0,T)}\quad\text{for all }\varphi\in W(0,T)
		\end{equation}
		In particular, $\varphi(t)$, $t\in[0,T]$, is meaningful in $H$ for an element $\varphi\in W(0,T)$.
		\item [\rm 3)] The integration by parts formula reads for all $\varphi, \phi \in W(0,T)$
		\begin{align}
			&\int_0^T{\langle \varphi_t(t),\phi(t) \rangle}_{V',V}\,\mathrm dt+\int_0^T{\langle \phi_t(t),\varphi(t)\rangle}_{V',V}\,\mathrm dt \label{eq:partialIntegration}\\
			&\quad= \int_0^T\frac{\mathrm d}{\mathrm dt}{\langle \varphi(t),\phi(t)\rangle}_H\,\mathrm dt= {\langle\varphi(\te),\phi(\te)\rangle}_H-{\langle\varphi(\tc),\phi(\tc)\rangle}_H\notag
		\end{align}
		\item [\rm 4)] We have the formula
		\begin{align*}
			{\langle \varphi_t(t),\phi \rangle}_{V',V} = \frac{\mathrm d}{\mathrm dt}\,{\langle \varphi(t),\phi\rangle}_H\quad\text{for }(\varphi,\phi) \in W(0,T)\times V, \text{ f.a.a. } t \in [0,T].
		\end{align*}
	\end{enumerate}
\end{lemma}

\begin{definition}[Continuous functions]
	Let $X,Y$ be Banach spaces and $\mathcal O \subset X$ an open subset. Then we denote by 
    \begin{align*}
	   C(\mathcal O;Y) := \{ f: \mathcal{O} \subset X \to Y \mid f \text{ is continuous} \}
	\end{align*}
	the set of all continuous functions from $\mathcal{O}$ to $Y$.
\end{definition}

\begin{definition}[G\^ateaux and Fr\'echet derivatives]
	\label{Definition:GFDer}
	Let $X,Y$ be real Banach spaces, $\mathcal O \subset X$ an open subset and $f: \mathcal O \to Y$ a given mapping.
    \begin{enumerate}
        \item [\rm 1)] For a point $x \in \mathcal O$ and a given direction $x^\delta \in X$, the \index{Derivative!directional}\emph{directional derivative} of $x$ towards $x^\delta$ is defined by
    	\begin{align*}
    		Df(x,x^\delta) := \lim_{\varepsilon \searrow 0} \frac{1}{\varepsilon} \big( f(x + \varepsilon x^\delta) - f(x) \big)
    	\end{align*}
    	provided the limit exists in $Y$.
        \item [\rm 2)] If the directional derivative exists for every direction $x^\delta \in X$ and we have $Df(x,\cdot) \in \mathscr L(X,Y)$, then $f$ is said to be \emph{G\^ateaux differentiable} at $x$ with \index{Derivative!G\^ateaux}\emph{G\^ateaux derivative} $f'(x) = Df(x,\cdot)$. If $f$ is G\^ateaux differentiable in every $x \in \mathcal O$, we say that $f$ is G\^ateaux differentiable (in $\mathcal O$).
        \item [\rm 3)] If $f$ is G\^ateaux differentiable in $x \in \mathcal O$ and it holds
    	\begin{align*}
    		\frac{\|f(x+x^\delta)-f(x)-f'(x)x^\delta\|_Y}{\|x^\delta\|_X} \to 0 \quad \text{as} \quad x^\delta \to 0 \text{ in } X,
    	\end{align*}
    	then $f$ is said to be \index{Derivative!Fr\'echet}\emph{Fr\'echet differentiable} in $x$. If $f$ is Fr\'echet differentiable in every $x \in \mathcal O$, we call $f$ Fr\'echet differentiable (in $\mathcal O$). In that case, if the mapping
    	\begin{align}
    		\label{app_derivativeFunction}
    		f': \mathcal O \to L(X,Y), \quad x \mapsto f'(x)
    	\end{align}		 
    	is continuous, $f$ is called \emph{continuously Fr\'echet differentiable} (in $\mathcal O$).
    \end{enumerate} 
\end{definition}

\begin{definition}[Second-order Fr\'echet derivatives]
	\label{def:secondOrderDerivatives}
	Let $X, Y$ be real Banach spaces, $\mathcal O \subset X$ an open subset and $f: \mathcal O \to Y$ Fr\'echet differentiable. Then we call $f$ \emph{twice Fr\'echet differentiable} in $x \in \mathcal O$ if the function $f'$ from \eqref{app_derivativeFunction} is Fr\' echet differentiable in $x$, i.e. if a mapping $f''(x) \in L(X,L(X,Y))$ exists with
	\begin{align*}
		\frac{\|f'(x+x^\delta)-f'(x)-f''(x)x^\delta \|_{L(X,Y)}}{\|x\|_X} \to 0 \quad \text{as} \quad x^\delta \to 0 \text{ in } X.
	\end{align*}
	In that case we call $f''(x)$ the \index{Derivative!second Fr\'echet}\emph{second Fr\'echet derivative} of $f$ in $x$. If $f$ is twice Fr\'echet differentiable in every $x \in \mathcal O$, $f$ is called twice Fr\'echet differentiable (in $\mathcal O$). If the mapping 
	\begin{align*}
		f'': \mathcal O \to L(X,L(X,Y)), \quad x \mapsto f''(x)
	\end{align*}
	is continuous, $f$ is said to be \emph{twice continuously Fr\'echet differentiable} (in $\mathcal O$).
\end{definition}

\begin{remark}
	\label{app_secondDerAsBilinear}
	\rm For two real Banach spaces $X$ and $Y$, a function $A: X \times X \to Y$ is said to be bilinear if for every $x \in X$, the functions $A(x,\cdot)$ and $A(\cdot,x)$ are linear mappings from $X$ to $Y$. It is said to be bounded if there exists a constant $c > 0$ such that
	\begin{align*}		
		{\|A(x_1,x_2)\|}_Y\le c\,{\|x_1\|}_X{\|x_2\|}_X\quad \text{for all } x_1,x_2 \in X
	\end{align*}
	If we denote the space of bilinear and bounded functions with $\mathscr L_2(X,Y)$ and endow it with the norm
	\begin{align*}
		{\| A \|}_{\mathscr L_2(X,Y)} = \sup \big \{ \|A(x_1,x_2)\|_Y ~\big|~{\|x_1\|}_X ={\|x_2\|}_X = 1 \big \},
	\end{align*}
	it can be shown (see, e.g., \cite{Ger06}) that it is a Banach space and isometrically isomorph to the space $\mathscr L(X,\mathscr L(X,Y))$. Therefore, we can identify the second Fréchet derivative with a bilinear, continuous function $f''(x) \in \mathscr L_2(X,Y)$. Schwarz's theorem \cite{Ger06} states that $f''(x)$ is always symmetric, so we have $f''(x)(x^\delta_1,x^\delta_2) = f''(x)(x^\delta_2,x^\delta_1)$ for all directions $x^\delta_1, x^\delta_2 \in X$.\hfill$\blacksquare$ 
\end{remark}

\begin{definition}[Gradient and Hessian]
	Let $X$ be a real Hilbert space, $\mathcal O \subset X$ open and $f: \mathcal O \to \mathbb R$ Fréchet differentiable in $x \in \mathcal O$. With the Riesz isomorphism $\mathcal J_X: X \to X'$, the \emph{Gradient} of $f$ at $x$ is defined as $\nabla f(x) = \mathcal J_X^{-1} f'(x) \in X$ and satisfies
	\begin{align*}
		{\langle \nabla f(x), x^\delta \rangle}_X = f'(x) x^\delta \quad \text{for all } x^\delta \in X.
	\end{align*}
	If $f$ is twice differentiable in $x$, the \emph{Hessian} of $f$ in $x$ is defined as a mapping
	\begin{align*}
		\nabla^2 f(x): X \to X, \quad \nabla^2 f(x)x^\delta = \mathcal J_X^{-1} f''(x) x^\delta.
	\end{align*}
	If we understand $f''(x)$ as a bilinear function as mentioned in Remark \ref{app_secondDerAsBilinear}, it holds by definition of the Riesz isomorphism
	\begin{align*}
		{\langle \nabla^2 f(x) x^\delta_1, x^\delta_2 \rangle}_X = f''(x)(x^\delta_1, x^\delta_2) \qquad \text{for all } x^\delta_1, x^\delta_2 \in X.
	\end{align*}
	Since the Riesz isomorphism is an isometry, it holds $\|\nabla f(x)\|_X = \|f'(x)\|_{X'}$ and $\nabla^2 f(x) \in \mathscr L(X)$ with $\|\nabla^2 f(x)\|_{\mathscr L(X)} = \|f''(x)\|_{\mathscr L_2(X)}$. 
\end{definition}


\bibliographystyle{alpha}
\bibliography{Ref.bib}

\backmatter

\printindex

\end{document}